\documentclass[12pt,centertags,oneside]{amsart}

\usepackage{tikz}
\usepackage[mathscr]{eucal}
\usepackage{asymptote}
\usepackage{caption}
\usetikzlibrary{arrows}
\usepackage{import}
\usepackage{graphicx,graphics,array,booktabs}

\usepackage{amsmath,amstext,amsthm,amscd,typearea,hyperref}
\usepackage{amssymb}
\usepackage{a4wide}
\usepackage[mathscr]{eucal}
\usepackage{mathrsfs}
\usepackage{typearea}
\usepackage{charter}
\usepackage{pdfsync}
\usepackage{multirow}
\usepackage{array}
\usepackage{amsmath, amssymb, latexsym, amscd, amsthm,amsfonts,amstext}
\usepackage[mathscr]{eucal}

\usepackage{amsmath,amstext,amsthm,amscd,typearea,stmaryrd,mathabx}
\usepackage{amssymb}
\usepackage{a4wide}
\usepackage[mathscr]{eucal}
\usepackage{mathrsfs}
\usepackage{typearea}
\usepackage{charter}
\usepackage{pdfsync}
\usepackage[a4paper,width=16.2cm,top=3cm,bottom=3cm,lines=50]{geometry}

\numberwithin{equation}{section}

\newtheorem{theorem}{Theorem}[section]
\newtheorem{definition}[theorem]{Definition}
\newtheorem{proposition}[theorem]{Proposition}
\newtheorem{corollary}[theorem]{Corollary}
\newtheorem{lemma}[theorem]{Lemma}
\newtheorem{remark}[theorem]{Remark}

\newtheorem{example}[theorem]{Example}

\newtheorem{notation}[theorem]{Notation}

\newcommand{\cali}[1]{\mathscr{#1}}
\newcommand{\Tan}{\mathop{\mathrm{Tan}}\nolimits}

\newcommand{\Cotan}{\mathop{\mathrm{Cotan}}\nolimits}

\newcommand{\comp}{{\mathop{\mathrm{comp}}\nolimits}}
\usepackage{graphicx} 
\newcommand{\intprod}{\mathbin{\raisebox{\depth}{\scalebox{1}[-2]{$\lnot$}}}}

\newcommand{\AB}{{\rm AB}}

\newcommand{\DO}{{\rm DO}}
\newcommand{\Leb}{{\rm Leb}}

\newcommand{\Tube}{{\rm Tube}}

\newcommand{\supp}{{\rm supp}}

\newcommand{\dist}{{\rm dist}}

\newcommand{\ver}{{\rm ver}}
\newcommand{\hor}{{\rm hor}}

\newcommand{\loc}{{loc}}
\renewcommand{\top}{{top}}
\newcommand{\ddc}{{dd^c}}

\newcommand{\ddcw}{{dd^c_w}}
\newcommand{\ddcwz}{{dd^c_{w,z}}}
\newcommand{\ddct}{{dd^c_t}}
\newcommand{\ddczeta}{{dd^c_\zeta}}
\newcommand{\ddczetap}{{dd^c_{\zeta'}}}
\newcommand{\ddczetat}{{dd^c_{\zeta',t}}}

\newcommand{\ddcwzeta}{{dd^c_{w,\zeta}}}
\newcommand{\ddcwzetap}{{dd^c_{w,\zeta'}}}
\newcommand{\ddcv}{{dd^c_\ver}}
\newcommand{\dc}{{d^c}}
\newcommand{\dbar}{{\overline\partial}}
\newcommand{\ddbar}{{\partial\overline\partial}}

\newcommand{\SH}{{\rm SH}}
\newcommand{\PH}{{\rm PH}}
\newcommand{\CL}{{\rm CL}}
\newcommand{\GL}{{\rm GL}}

\newcommand{\tot}{{\rm tot}}
\newcommand{\ind}{{\bf 1}}
\newcommand{\bfm}{{\bf m}}
\newcommand{\bfn}{{\bf n}}
\newcommand{\id}{{\rm id}}

\newcommand{\exponent}{{\rm exponent}}

\newcommand{\codim}{{\rm codim\ \!}}

\renewcommand{\Re}{{\rm Re}}
\renewcommand{\Im}{{\rm Im}}
\newcommand{\bfc}{{\rm \mathbf{c}}}

\newcommand{\bfr}{{\rm \mathbf{r}}}
\newcommand{\bfx}{{\rm \mathbf{x}}}
\newcommand{\bfi}{{\rm \mathbf{i}}}
\newcommand{\bfI}{{\rm \mathbf{I}}}
\newcommand{\bfj}{{\rm \mathbf{j}}}
\newcommand{\bfU}{{\rm \mathbf{U}}}
\newcommand{\bfW}{{\rm \mathbf{W}}}
\newcommand{\bfe}{{\rm \mathbf{e}}}
\newcommand{\upm}{{\mathrm{ \overline{m}}}}
\newcommand{\lowm}{{\mathrm{ \underline{m}}}}

\newcommand{\Cc}{\cali{C}}
\newcommand{\Dc}{\cali{D}}
\newcommand{\Ec}{\cali{E}}
\newcommand{\Fc}{\cali{F}}

\newcommand{\Hc}{\cali{H}}
\newcommand{\Kc}{\cali{K}}
\newcommand{\Lc}{\cali{L}}
\renewcommand{\Mc}{\cali{M}}

\newcommand{\Uc}{\cali{U}}
\newcommand{\Vc}{\cali{V}}

\newcommand{\Ic}{\cali{I}}
\newcommand{\Nc}{\cali{N}}
\newcommand{\FS}{{\rm FS}}

\newcommand{\C}{\mathbb{C}}
\newcommand{\D}{\mathbb{D}}
\newcommand{\E}{\mathbb{E}}
\renewcommand{\H}{\mathbb{H}}
\newcommand{\N}{\mathbb{N}}

\newcommand{\R}{\mathbb{R}}
\newcommand{\T}{\mathbb{T}}
\newcommand{\B}{\mathbb{B}}
\newcommand{\G}{\mathbb{G}}
\newcommand{\U}{\mathbb{U}}

\newcommand{\M}{\mathbb{M}}

\renewcommand{\P}{\mathbb{P}}
\newcommand{\X}{\mathbb{X}}

\title[]{Positive plurisubharmonic currents: Generalized Lelong numbers and Tangent  theorems}

\author{Vi{\^e}t-Anh Nguy{\^e}n}
\address{Universit\'e de Lille, 
Laboratoire de math\'ematiques Paul Painlev\'e, 
CNRS U.M.R. 8524,  
59655 Villeneuve d'Ascq Cedex, 
France. }

\address{and Vietnam Institute for Advanced Study in Mathematics (VIASM),  157 Chua Lang Street, Hanoi, Vietnam.
}

\email{Viet-Anh.Nguyen@univ-lille.fr, {\tt   https://pro.univ-lille.fr/viet-anh-nguyen/}}

\date{May 30, 2022}

\begin{document}


\begin{abstract}
   Dinh--Sibony  theory of tangent  and density currents is a recent but   powerful tool  to study  positive closed currents.  Over twenty years  ago, Alessandrini and Bassanelli initiated the theory
   of  the Lelong number of a positive  plurisubharmonic  current in $\C^k$ along a linear subspace. Although  the latter theory   is    intriguing,  it has not yet been  explored in-depth  since then.
Introducing    the   concept of  the generalized  Lelong numbers and  studying these new numerical values, we   extend both  theories   to  a more general  class of positive  plurisubharmonic  currents  and  in a more general context  of   ambient manifolds.

More specifically, in the first part  of our article, 
we consider  a positive plurisubharmonic current $T$ of bidegree $(p,p)$   on a    complex  manifold $X$ of dimension $k,$ 
and let $V\subset X$ be   a K\"ahler submanifold  of dimension $l$  and $B$ a  relatively compact piecewise $\Cc^2$-smooth  open subset of $V.$
We impose  a mild reasonable condition on $T$ and $B,$   namely,   $T$ is weakly  approximable  by $T_n^+-T_n^-$ on a neighborhood $U$ of $\overline B$ in $X,$   where 
$(T^\pm_n)_{n=1}^\infty$ are  some  positive plurisubharmonic $\Cc^3$-smooth forms of bidegree $(p,p)$  defined on $U$ such that  the masses $\|T^\pm_n\|$ on $U$   are uniformly bounded and 
that   the  $\Cc^3$-norms  of  $T^\pm_n$ are   uniformly bounded near  $\partial B$
if $\partial B\not=\varnothing.$ Note  that if $X$ is K\"ahler and  $T$  is of class  $\Cc^3$ near $\partial B,$  then the   above mild condition  is satisfied. In particular,  this   $\Cc^3$-smoothness near $\partial B$ is   automatically fulfilled  if either  $\partial B=\varnothing$ or $V\cap \supp( T)\subset B.$
\begin{itemize}
\item[$\bullet$] We define  the notion  of  the $j$-th Lelong number of $T$ along $B$   for every $j$  with $\max(0,l-p)\leq j\leq \min(l,k-p)$ and  prove  their existence as well as their basic properties.
We also  show that $T$  admits   tangent currents and that  all tangent  currents  are not only  positive plurisubharmonic, but also  partially $V$-conic and    partially pluriharmonic. 

\item[$\bullet$] When  the currents $T^\pm_n$  are  moreover pluriharmonic (resp.  closed), we show, under a less restrictive smoothness of $T^\pm_n$ near $\partial B,$ that  every tangent current is 
also  $V$-conic pluriharmonic (resp.  $V$-conic closed).

\item[$\bullet$] We also prove that the  generalized  Lelong numbers  are intrinsic. 

\item[$\bullet$] In fact, if we are only interested in    the top degree Lelong  number of $T$  along $B$ (that is, the $j$-th Lelong number
for the maximal value  $j=\min(l,k-p)$), then  under  a  suitable  holomorphic context, the above condition  on the uniform regularity  of $T^\pm_n$  near $\partial B$ can be removed. 
\end{itemize}
Our  method relies on some  Lelong-Jensen formulas  for
the normal  bundle to $V$ in $X,$ which are of independent interest.

The second part of our article is  devoted  to  geometric  characterizations of the  generalized  Lelong  numbers. 
As a consequence of this  study, we  show that   the top degree Lelong  number of $T$  along $B$   is  totally intrinsic.
This is  a generalization  of the fundamental
result  of Siu (for positive closed  currents) and of Alessandrini--Bassanelli (for positive  plurisubharmonic currents)
on the independence of Lelong numbers at a single point on the choice of  coordinates.  
 
\end{abstract}

\maketitle
\tableofcontents

\medskip\medskip


 \medskip



\section{Introduction} \label{S:Intro}

\subsection{Motivations}
 Let $X$ be a   complex   manifold of dimension $k.$ Let   $d,$ $\dc$ denote the  real  differential operators on $X$  defined by
$d:=\partial+\overline\partial,$  $\dc:= {1\over 2\pi i}(\partial -\overline\partial) $ so that 
$\ddc={i\over \pi} \partial\overline\partial.$  A  $(p,p)$-current  $T$ defined on $X$  is  said to be {\it closed}  (resp. {\it  pluriharmonic\footnote{Some  authors
uses the  terminology  {\it harmonic} instead of {\it pluriharmonic}.}}), (resp. {\it plurisubharmonic})
if $dT=0$  (resp. $\ddc T=0$), (resp. $\ddc T$  is a positive current).  Here are   
relations of strict inclusions between several  well-known classes of currents on $X:$
\begin{multline*}
\{ \text{currents of integrations on complex subvarieties}\}\subsetneq\{ \text{positive  closed currents}\}\\
\subsetneq \{ \text{positive  pluriharmonic currents}\}
\subsetneq\{ \text{positive  plurisubharmonic currents}\}.
\end{multline*}
Let $T$ be  a positive plurisubharmonic  $(p,p)$-current  defined on $X$  and $x\in X$   a point.
We first recall the  notion  of  Lelong number  $\nu(T,a)$ of $T$ at  $x.$
This  notion  was first introduced  by   Lelong in \cite{Lelong}  for the class of positive closed currents.
It was later  formulated  by Skoda in \cite{Skoda} for the wider class of positive  plurisubharmonic currents.  The  notion  plays a fundamental role
in  Complex  Analysis, Complex  Geometry, Algebraic Geometry and Complex  Dynamics. The readers could  find more systematic  developments on Lelong numbers  for  positive closed  currents in Siu \cite{Siu}  and Demailly
\cite{Demailly93,Demailly}  as well as   the references therein.  As for   positive  plurisubharmonic currents,  the reader could  consult Alessandrini--Bassanelli  \cite{AlessandriniBassanelli96}. 

Choose a local holomorphic coordinate system $z$ near $x$ such that $x = 0$ in
these coordinates. The Lelong number $\nu(T,x)$ of $T$ at $x$ is the limit of the normalized
mass of $\|T \|$  on the ball $\B(0, r)$ of center $0$ and radius $r$ when $r$ tends to $0.$ More
precisely, we have
\begin{equation}\label{e:Lelong-number-point}
\nu(T , x) := \lim\limits_{r\to 0}\nu(T,x,r),\quad\text{where}\quad \nu(T,x,r):=
{\sigma_T(\B(0,r))\over 
(2\pi)^{k-p} r^{ 2k-2p}}. 
\end{equation}
Here, $\sigma_T:={1\over (k-p)!}\, T\wedge ({i\over 2}\ddbar \|z\|^2)^{k-p}$ is  the trace measure of $T.$
Note that $(2\pi)^{k-p} r^{2k-2p}$  is the mass on $\B(0, r)$ of the $(p , p  )$-current of
integration on a linear subspace of dimension $k - p$ through $0.$  When $T$ is a positive closed current, Lelong establishes  in \cite{Lelong} (see also \cite{Lelong68})
that 
the  {\it  average mean} $\nu(T,x,r)$  is a  non-negative-valued  increasing function in the  radius $r.$
So  the limit \eqref{e:Lelong-number-point} always exists. 
Skoda \cite{Skoda} proves  the  same result  for  positive  plurisubharmonic currents.   Thie \cite{Thie} shows
that when $T$ is given by an
analytic set this number is the multiplicity of this set at $x.$ Siu proves  that when $T$ is  a positive  closed  current, the
limit \eqref{e:Lelong-number-point} does not depend on the choice of coordinates.

There is  another equivalent  logarithmic definition of the Lelong number of a positive closed  current that we want to discuss in this   work.
Namely,  we have 
\begin{equation}\label{e:Lelong-number-point-bis}
\nu(T , x) := \lim\limits_{r\to 0} \kappa(T,x,r),\quad\text{where}\quad \kappa(T,x,r):= \int_{\B(0,r)} T(z)\wedge (\ddc\log{(\|z\|^2}))^{k-p}\,\cdot
\end{equation}
The {\it  logarithmic mean}  $\kappa(T,x,r)$ is a  non-negative-valued  increasing function in the  radius $r.$  
Observe  that  in the expression of  $\kappa(T,x,r)$  in \eqref{e:Lelong-number-point-bis}, the  wedge-product of currents  is only well-defined  outside the origin $0$ because
the  second factor $ (\ddc\log{(\|z\|^2}))^{k-p}$ is  only smooth there. 
In this  article we   deal  with two simple interpretations of  \eqref{e:Lelong-number-point-bis}  which correspond to
regularizing either  the first  or the second factor of the  wedge-product of currents in the  expression of
$ \kappa(T,x,r).$
 The  first interpretation concerns the  notion of approximation of currents.  By a standard regularization (e.g.  a convolution),  we see that there is a sequence of
 positive smooth closed  $(p,p)$-form  on   $\B(0,r+\epsilon)$ for some $\epsilon>0$  such that $T_n$  converges weakly  to $T.$ 
 The first  interpretation   
 of the integral on the RHS  of  \eqref{e:Lelong-number-point-bis}  is formulated as follows:
\begin{equation}\label{e:Lelong-number-point-bisbis(1)}
 \int_{\B(0,r)} T(z)\wedge (\ddc\log{(\|z\|^2}))^{k-p}:=\lim\limits_{n\to \infty}   \int_{\B(0,r)} T_n(z)\wedge (\ddc\log{(\|z\|^2)})^{k-p} \,\cdot
\end{equation}
provided that the limit exists. In fact, this  is indeed the case.
The  second  interpretation  consists in  regularizing  the integral kernel $(\ddc\log{(\|z\|^2}))^{k-p}$ in a    canonical way:
\begin{equation}\label{e:Lelong-number-point-bisbis(2)}
 \int_{\B(0,r)} T(z)\wedge (\ddc\log{(\|z\|^2)})^{k-p}:=\lim\limits_{\epsilon\to 0+}   \int_{\B(0,r)} T(z)\wedge (\ddc\log{(\|z\|^2+\epsilon^2)})^{k-p} \,\cdot
\end{equation}
provided that the limit exists. In fact, this is always  the case.

\smallskip

Next,  we revisit  the Lelong number of positive  closed  currents from another geometric point of
view related to Harvey's exposition \cite{Harvey}. 
Let $X$ be an  open  neighborhood of $0$ in $\C^k.$ Let $A_\lambda : \C^k \to  \C^k$ be defined by
$A_\lambda (x) := \lambda x$ with $\lambda \in \C^* .$  When $\lambda$ goes to infinity, the domain of definition
of the current $T_\lambda := (A_\lambda)_* (T)$ converges to $\C^k .$ This family of currents
is relatively compact,  and any limit current $T_\infty$ for $\lambda \to\infty,$ is  
called a {\it tangent current} to $T.$
A tangent  current  is defined on the  whole $\C^k,$ and it is 
conic  in the  sense  that  it is  invariant under $(A_\lambda)_* .$  
Given a  tangent  current $T_\infty$  to $T,$  we can
extend it to $\P^k$ with zero mass on the hyperplane at infinity. Thus, there is a
positive closed current $\T_\infty$ on $\P^{k-1}$  such that $T_\infty=\pi^*_\infty
(\T_\infty).$ Here we identify the hyperplane at infinity with $\P^{k-1}$
and we denote by $\pi_\infty: \P^k \setminus \{0\} \to \P^{k-1}$
the canonical central projection. The class of $\T_\infty$ (resp. of $T_\infty$) in the de Rham
cohomology of $\P^{k-1}$ (resp., of $\P^k$) is equal to $\nu(T , x)$ times the class of a
linear subspace. So these cohomology classes do not depend on the choice of
$T_\infty.$ In  general, the tangent current $T_\infty$ is not unique, see  Kiselman   \cite{Kiselman}. 

   \smallskip
   
   {\bf  Notation.} Throughout  the article,  we denote by
   \begin{itemize}
   
 \item   $\D$  the unit disc in $\C;$
 \item  $\C^*$ the punctured complex plane $\C\setminus \{0\};$  
 \item  $\R^+:=[0,\infty)$ and  $\R^+_*:=(0,\infty);$
 \item  $\partial B$  the boundary  of an open  set $B$ in a manifold $Y.$

   \end{itemize}

   If $X$ is an oriented manifold, denote by $H^* (X, \C)$ the de Rham cohomology
group of $X$ and $H^*_\comp (X, \C)$ the de Rham cohomology group deﬁned by forms
or currents with compact support in $X.$ If $V$ is a submanifold of $X,$ denote
by $H^*_V (X, \C)$ the de Rham cohomology group deﬁned in the same way using
only forms or currents on $X$ whose supports intersect $V$ in a compact set.

If $T$ is either  a closed current on $X$  or a $\ddc$-closed current on a compact K\"ahler manifold $X,$ denote by $\{T \}$
its class in $H^* (X, \C).$ When
$T$ is supposed to have compact support, then $\{T \}$ denotes the class of $T$
in $H_\comp^* (X, \C).$ If we only assume that $\supp(T ) \cap V$ is compact, then $\{T \}$ denotes the class of $ T$ in $H_V^* (X, \C).$ The current of integration on an oriented
submanifold or a  complex  variety $Y$ is denoted by $[Y ].$ Its class is denoted by $\{Y \}.$

  For a differentiable map $\pi: X\to Y$ between  manifolds,  $\pi^*$ (resp. $\pi_*$) denotes the pull-back (resp. the push-forward) operator
  acting  on forms and currents  defined on  $Y$ (resp. on $X$). These operators induce  natural maps on cohomological levels:
  $\pi^*:\  H^*(Y,\C)\to  H^*(X,\C)$ and $\pi_*:\  H^*(X,\C)\to  H^*(Y,\C).$
   
In the next subsection  we  present  a short digression  to two theories which are  the main  sources  of inspirations for this  work.    
   
\subsection{Dinh-Sibony theory}

 A fundamental achievement  has recently been  attained  by   Dinh and  Sibony \cite{DinhSibony18}.  
These authors develop  a  satisfactory  theory of tangent  currents and  density currents   for  positive  closed currents  in  the
context  where  the single point $x$  is  replaced  by   a submanifold $V\subset X$ of positive dimension $l$  ($1\leq l< k$).

  Let $\E$ be the normal
vector bundle to $V$ in $X$ and $\pi:\ \E\to V$ be the  canonical projection.
Let $\pi_0:\ \overline \E:=\P(\E\oplus\C)\to V$ be its canonical compactification. Denote
by $A_\lambda :\ \E \to \E$ the map induced by the multiplication by $\lambda$ on fibers of $\E$
with $\lambda \in \C^*.$  We identify $V$ with the zero section of $\E.$    We  expect as  in  Harvey's exposition \cite{Harvey} that every  tangent  current  $T_\infty$ lives on $\E.$
However,  a basic difficulty arises.  When $V$ has
positive dimension, in general, no neighbourhood of $V$ in $X$ is biholomorphic
to a neighbourhood of $V$ in $\E.$

To  overcome  this  difficulty,  Dinh and  Sibony propose  a  softer  notion:  {\it the  admissible maps.}  More  precisely, let $\tau$  be a diffeomorphism between a neighbourhood of $V$ in $X$ and a
neighbourhood of $V$ in $\E$ whose restriction to $V$ is identity. We assume that
$\tau$ is admissible in the sense that the endomorphism of $\E$ induced by the
differential of $\tau$  when restricted to $V$ is the identity map from $\E$ to $\E.$

 Fix $0\leq p\leq  k$ and set
\begin{equation}\label{e:m}\upm:= \min(l,k-p)\qquad\text{and}\qquad
  \lowm:=\max(0,l-p).
  \end{equation}
Here is  the main result of  Dinh and  Sibony.
\begin{theorem}\label{T:Dinh-Sibony-first}{\rm  (Dinh-Sibony \cite[Theorems 1.1, 4.6 and Definition 4.8]{DinhSibony18})}  Let $X,$ $V,$ $\E,$  $\overline \E,$ $A_\lambda$ and $\tau$ be as above. 
Let $T$ be  a positive closed $(p,p)$-current on $X.$
Assume in addition  that $X$ is  K\"ahler and  $\supp(T)\cap V$ is  compact. Then:
\begin{enumerate} \item The family
of currents $T_\lambda:= (A_\lambda)_* \tau_* (T )$ is relatively compact and any limit current, for
$\lambda\to\infty,$ is a positive closed $(p, p)$-current on $\E$ whose trivial extension is a
positive closed $(p, p)$-current on $\overline \E.$   Such a limit current $S$ is  called a {\rm tangent current to $T$ along $V.$} 
\item If $S$ is a tangent current to $T$ along $V$, then it is
$V$-conic, i.e., invariant under $(A_\lambda)_* ,$  and its de Rham cohomology class $\{S\}$ in
$H^{2p}_\comp (\E, \C)$ does not depend on the choice of $\tau$ and $S.$

\item  Let $-h_{\overline \E}$ denote the  tautological class of the bundle $\pi_0:\  \overline \E\to  V.$ Then  we have the following   decomposition
of the  cohomology  class $\{S\}$:
$$
\{S\}=\sum_{j=\lowm}^\upm \pi_0^*(\kappa_j(T))\smile  h_{\overline \E}^{j-l+p},
$$
where  $\kappa_j(T)$  is a class in  $H^{2l-2j}_\comp(V,\C).$ Moreover, this decomposition is  unique.
\end{enumerate}
\end{theorem}

 When $V$  has positive dimension $l,$
 according   to  Dinh  and Sibony,  the    notion of Lelong number of the current $T$ at a single point should be replaced by
 the  family of  cohomology classes $\{ \kappa_j(T):\  \lowm \leq j\leq \upm\}$  given by Theorem \ref{T:Dinh-Sibony-first} (3) above. 
This is an important and original viewpoint  of  Dinh and Sibony.

Since  then, this theory  finds  numerous applications in Complex Analysis  and  Algebraic  Geometry, especially in  Complex  Dynamics and the theory of foliations. In particular, Dinh and  Sibony apply their theory in order to obtain  the  equidistribution  of  saddle periodic points for  regular H\'enon type  automorphisms of  $\C^k$
 (see \cite{DinhSibony16}). This is  an important progress in Complex  Dynamics since  previous results are  only known for the dimension 
 $k=2$ (see \cite{BedfordLyubichSmillie}). 
Moreover, the theory itself  has also been  developed in many directions.
While  studying positive harmonic  currents  directed  by a singular holomorphic  foliation on compact K\"ahler  surfaces, Dinh and  Sibony  \cite{DinhSibony18} introduce a theory of tangent currents
in the  following new context:  $X$ is a  compact K\"ahler surface, $V$ is a  curve and  $T$ is a positive  harmonic $(1,1)$-current. Next,
in order to  establish   the unique  ergodicity  in the theory of singular  holomorphic  foliations, Dinh, Sibony and   the author \cite{DinhNguyenSibony18}   develop  a theory
of tangent and  density  currents for  tensor  product of  positive $\ddc$-closed currents on  compact K\"ahler surfaces.  A  further development in this  direction has been pursued  in  \cite{DinhNguyen20}.
On the other hand,
by optimizing  the original approach of  Dinh--Sibony, Vu \cite{Vu19a}  weakens the  K\"ahler  assumption on  $X,$  see the discussion  after Theorem \ref{T:main_2'} below.
Dinh, Huynh, Kaufmann, the  author, Truong, Vu  and  several other authors    apply this  theory to  many interesting  problems,  see \cite{DinhNguyenTruong15, DinhNguyenTruong17,HuynhKaufmannVu,HuynhVu,Kaufmann,KaufmannVu,Vu19a,Vu19b,Vu20a}  etc.  
\subsection{Alessandrini--Bassanelli theory}


 In  \cite{AlessandriniBassanelli96} Alessandrini and Bassanelli introduce  a remarkable  notion of Lelong number  of a  positive plurisubharmonic current
 in  a special  setting  of manifolds
$(X,V).$

\begin{theorem}\label{T:AB-1}{\rm  (Alessandrini and Bassanelli \cite[Theorem I and Definition 2.2]{AlessandriniBassanelli96})}
Consider $X=\C^k$ and $V$ is a linear complex subspace of dimension $l\geq 0. $   
  We use the coordinates $(z,w)\in\C^{k-l}\times \C^l$ so  that  $V=\{z=0\}.$ Let $0\leq p<k-l$ and  
  let $T$  be    a positive plurisubharmonic $(p,p)$-current on an open  neighborhood $\Omega$ of $0$ in $\C^k.$ 
  \begin{enumerate}
   \item 
Then, for every open ball $B$ in $V,$ $B\Subset \Omega,$ the  following limit  exists and is  finite
$$
 \nu_\AB(T,B):=\lim_{r\to 0+} {1\over r^{2(k-l-p)}}\int_{\Tube(B,r)}  T(z,w)\wedge (\ddc \|z\|^2)^{k-l-p}\wedge (\ddc \|w\|^2)^{l},
 $$
 where the  tube  $\Tube(B,r)$ of radius $r$ over $B$ is  given by 
 \begin{equation} \label{e:Tube-AB}
\Tube(B,r):=\left\lbrace (z,w)\in \C^{k-l}\times\C^{l}:\ \|z\|<r,\ w\in B\right\rbrace .
\end{equation}
 $\nu_\AB(T,B)$ is  called    the {\it  Alessandrini--Bassanelli's  Lelong  number of $T$ along  $B$}.
\item  There exist an open  neighborhood  $W$ of $0$ in $L,$ $W\subset \Omega,$  and a nonnegative  plurisubharmonic  function $f$ on $W$  such that
$$
\nu_\AB(T,B)=\int_B f(w) (\ddc \|w\|^2)^{l}
$$
for every open ball $B$ in $V$ with $B\Subset  W.$
\end{enumerate}

\end{theorem}

The  important viewpoint of Alessandrini--Bassanelli is  that when $V$ is of positive dimension,   tubular neighborhoods $\Tube(B,r)$ of $B$ and a  mixed form $  (\ddc \|z\|^2)^{k-l-p}\wedge (\ddc \|w\|^2)^{l}$ should replace
the usual balls $\B(x,r)$ around a single point $x$ with  the usual form  $  (\ddc \|z\|^2)^{k-p}.$
When $V$ is  a single point $\{x\}$ and $B=\{x\},$   Alessandrini--Bassanelli's  Lelong  number $\nu_\AB(T,x)$  coincides with the classical Lelong number
$\nu(T,x).$  They also obtain equivalent formulations of their  Lelong number in the spirit of  \eqref{e:Lelong-number-point-bisbis(1)}-\eqref{e:Lelong-number-point-bisbis(2)}.

 Alessandrini--Bassanelli's method relies on some Lelong-Jensen formulas  which can be obtained  from the usual Lelong-Jensen formula (see  \cite{Demailly, Skoda})  by slicing. 
They  also  characterize  this Lelong number geometrically in the sense   of Siu \cite{Siu}.
Namely, they use   the total space of the tautological vector bundles over suitable  Grassmannian manifolds  and  pull-back the  given current to this space,
and then  study  the cut-off  of this current on  the exceptional fibers.  In order to  state  a brief version of their result, we need to introduce
some more notations.

For every $1\leq j\leq k-l,$ let
$$
\X_j:=\left\lbrace (z,w,H)\in \C^{k-l}\times V\times \G_j(\C^{k-l}):\  z\in H   \right\rbrace ,
$$
where  $\G_j(\C^{k-l})$ is  the Grassmannian of all $j$-dimensional linear subspaces of $\C^{k-l}.$
Denote by $\Pi_j:\ \X_j\to\C^k=\C^{k-l}\times V$  the canonical projection. 
\begin{theorem}\label{T:AB-2}{\rm  (Siu \cite[Section 11]{Siu} for positive closed currents near a single point,   Alessandrini--Bassanelli \cite[Corollary 3.6 and Remark 3.7]{AlessandriniBassanelli96} for
positive  pluriharmonic and  positive plurisubharmonic  currents  near  a ball in a linear subspace)}
We keep the hypothesis of Theorem \ref{T:AB-1}.
Let $(T_n)_{n=1}^\infty$ be  a sequence of smooth  positive plurisubharmonic  forms on a neighborhood of $\overline B$ in $\Omega$
with  uniformly bounded masses such that  $T_n$ converge  to $T$ weakly as $n\to\infty$\footnote{Such a sequence can be obtained by a standard convolution with $T$.}. 
\begin{enumerate} \item Then, for a  suitable  subsequence  $(T_{N_n})_{n=1}^\infty,$ the  following  weak limit  exists
$$
\widetilde T:= \lim_{n\to\infty} \Pi_p^* (T_{N_n})
$$
and $\widetilde T$ is a positive plurisubharmonic  $(p,p)$-current on $\X_p.$
\item   The following  current
$$
\widetilde S_n:=  (\Pi_{p+1}|_{\Omega\setminus V}^*)\big((-\log{\|z\|})\ddc T_{N_n}\big).
$$
exists and has trivial  extension $(\widetilde S_n)_\bullet$ through $ \Pi^{-1}_{p+1}(V).$
By passing  to a subsequence  if necessary, we can define the current 
$$\widetilde S:=\lim_{n\to\infty}(\widetilde S_n)_\bullet$$
which is a positive $(p+1,p+1)$-current on $\X_{p+1}.$
\item The  following  identity holds
 $$
 \nu(T,B)=\| \widetilde T \|(\Pi^{-1}_p(B))+\|\widetilde S\| (\Pi_{p+1}^{-1}(B)),
 $$
 where, for a positive current $R$ defined on a complex manifold $M$ and for a Borel subset $A\subset M,$   $\|R\|_A$  denotes the mass of $R$  on $A$ (see Subsection \ref{SS:currents}).
 
 \item If $T$ is moreover  closed (resp.  pluriharmonic), then $(T_n)_{n=1}^\infty$  can be  chosen to be positive  closed (resp. positive pluriharmonic), and  hence
 $$
 \nu_\AB(T,B)=\| \widetilde T \|(\Pi^{-1}_p(B)).
 $$
 \end{enumerate}
\end{theorem}
Theorem \ref{T:AB-2} (4)  says   that  when $T$ is positive pluriharmonic (resp. positive closed), the     Alessandrini--Bassanelli's  Lelong  number $\nu_\AB(T,B)$ is  equal to the mass of the limiting pluriharmonic (resp. closed) current $\widetilde T$ on the exceptional fiber 
 on $B.$ In particular, this mass is independent of  limiting currents. 
 
Using this  geometric  interpretation (Theorem \ref{T:AB-2} (3)) for the case when $V$ is  a  single point, the  following  result  is established.
\begin{theorem}\label{T:AB-3} {\rm  (Siu \cite[Section 11]{Siu} for positive closed currents,   Alessandrini--Bassanelli \cite[Theorem II]{AlessandriniBassanelli96} for positive plurisubharmonic  currents)}
 Let $F:\ \Omega\to\Omega'$ be  a  biholomorphic map between open subsets of $\C^k.$
 If $T$ is a  positive  plurisubharmonic $(p,p)$-current on $\Omega$ and $x\in\Omega,$ then
 $$
 \nu(T,x)=\nu(F_*T,F(x)).
 $$
\end{theorem}
Hence, the limit \eqref{e:Lelong-number-point} does not depend on the choice of coordinates even  for   positive  plurisubharmonic currents. 
So,  the Lelong number  of a  positive plurisubhamonic   current at a single point is an intrinsic notion.
 
 Although  the  assumption   on  the  pair of manifolds $(X,V)$ in Theorem  \ref{T:AB-1} is  quite  restrictive and  this theorem   provides only one  Lelong number, Alessandrini--Bassanelli theory may be regarded as the first  effort to elaborate  the notion of  numerical Lelong numbers  when the dimension of $V$ is positive.
\subsection{Main purpose of the  article}

The  main purpose of  this  work is to create   a unified framework  where we can develop and  generalize   both the  above mentioned theories.  
There  are  two concrere tasks.
The  first one
 is  to generalize the notion of     Dinh--Sibony \cite{DinhSibony18} on tangent  and  density  currents  
\begin{itemize} 
\item[$\bullet$] for    a  very  general and natural  class of currents:
the  positive plurisubharmonic  currents;
\item   for a  general and natural context of a piecewise smooth  open set $B\subset V:$ 
studying the  tangent  currents to $T$ along 
$B.$
\end{itemize}
The  second  task   is  to 
generalize   the   notion of Alessandrini--Bassanelli  \cite{AlessandriniBassanelli96} on   Lelong numbers, 
and the  results of   Siu  \cite{Siu} and of Alessandrini--Bassanelli  \cite{AlessandriniBassanelli96} on geometric characterizations of Lelong numbers
to the above  contexts.
 Following  the  tradition of  Lelong  \cite{Lelong} and  Skoda \cite{Skoda},   we   formulate  some  natural  and   numerical  Lelong numbers.  
So  our  viewpoint  which is  close  to that  of  Alessandrini--Bassanelli seems to be   quite  different from    Dinh-Sibony's  viewpoint  of  defining  some cohomology  classes  as  Lelong numbers.

A  novelty of our  work  is  that   our  approach is  technically different  from 
those of   Dinh--Sibony. Indeed,
these authors rely on a  cohomological calculus for positive closed currents.
When   neither  the currents  in  questions are  closed  nor their supports are compactly intersected with $V,$
this method does not  seems  to be  applicable.
To overcome   this  basic difficulty,   our new  key tool to study  the  tangent  currents to $T$ along 
an open subset $B\subset V$  is     some Lelong--Jensen type  formulas, which are applied    at the limit  on $B\subset V$, that is,   on an  infinitesimally  small 
tubular neighborhood
of $B$ in the  normal  bundle  to $V$ in $X.$   So  our approach is  close  to  that of  Alessandrini--Bassanelli.
But our Lelong--Jensen type formulas are more general than theirs. Indeed, even  in   their context where the normal bundle $\E$ is  trivial,  our formulas are more general and they  can not be obtained  from slicing method.    In order to make  our machinary work,  we introduce  new  classes of currents
 which  satisfy   mild reasonable approximation conditions.

\subsection{Approximations  and new classes of currents} \label{SS:New-classes-currents}
Now we are in the position to define  the needed  notions of approximations.

 \begin{definition}\label{D:approximable}
  \rm 
    Let $m,m'\in\N$ with $m\geq m'.$  Let  $W\subset U\subset X$ be two open subsets.
  Let $T$ be  a  positive  $(p,p)$-current  defined  on an open  set containing  $U.$ 

 \smallskip
 
 \noindent  (1) We say that  $T$ is  {\it approximable  on $U$ by  $\Cc^m$-smooth positive plurisubharmonic forms} 
  and write $T\in\SH^{m}_p(U)$
   if there  is  a sequence of $\Cc^m$-smooth
  positive plurisubharmonic 
  $(p,p)$-forms $(T_n)_{n=1}^\infty$  defined on $U$  
  such that
  \begin{enumerate}
\item[(i)] the masses $\|T_n\|$ on $U$   are uniformly bounded; 
\item[(ii)]   $T_n$ converge  weakly to  $T$    on $U$  as $n$ tends to infinity.
   \end{enumerate}
   If moreover,  the following condition is  fulfilled:
   \begin{itemize}
    \item [(iii-a)]  the restrictions  of the  forms $T_n$ on $W$   are of uniformly   bounded $\Cc^{m'}$-norm;
   \end{itemize}
then  we  say that  $T$ is  {\it approximable  on $U$ by  $\Cc^m$-smooth positive plurisubharmonic forms with  $\Cc^{m'}$-control on $W,$}
and write $T\in\SH^{m,m'}_p(U,W).$
   
    If moreover,  the following condition is  fulfilled:
   \begin{itemize}
    \item [(iii-b)]   $\supp (T_n)\cap W=\varnothing$   for $n\geq 1;$
   \end{itemize}
then  we  say that  $T$ is  {\it approximable  on $U$ by  $\Cc^m$-smooth positive plurisubharmonic forms with support outside $W,$}
and write $T\in \SH^{m}_p(U,W,\comp).$

   We     say  that $(T_n)_{n=1}^\infty$ is a {\it sequence of approximating  forms} for  $T$
   as an element of   $\SH^{m}_p(U)$ in the  first case (resp. as an element of  $\SH^{m,m'}_p(U,W)$ in the second case, resp.
   as an element of  $\SH^{m}_p(U,W,\comp)$ in the third case).

 \noindent  (2) Similarly,  we say that  $T$ is  {\it approximable  on $U$ by  $\Cc^m$-smooth positive pluriharmonic (resp. positive closed) forms} and   we write $T\in\PH^{m}_p(U)$   (resp.  $T\in  \CL^{m}_p(U)$) 
   if,  the approximating forms  $T_n$ satisfying  conditions (i)--(ii) in Definition \ref{D:approximable} (1) are   positive pluriharmonic (resp. positive closed).

 \smallskip

 \noindent  (3)  We say that  $T$ is  {\it approximable  on $U$ by  $\Cc^m$-smooth positive pluriharmonic (resp. positive closed) forms with  $\Cc^{m'}$-control on $W$} and   we write $T\in\PH^{m,m'}_p(U,W)$   (resp.  $T\in  \CL^{m,m'}_p(U,W)$) 
   if,   the approximating forms  $T_n$  satisfying  conditions (i)--(ii)--(iii-a)  in Definition \ref{D:approximable} (1) are   positive pluriharmonic (resp. positive closed).

 \smallskip

 \noindent  (4)  We say that  $T$ is  {\it approximable  on $U$ by  $\Cc^m$-smooth positive pluriharmonic (resp. positive closed) forms with  support outside $W$} and   we write $T\in\PH^{m}_p(U,W,\comp)$   (resp.  $T\in  \CL^{m}_p(U,W,\comp)$) 
   if,   the approximating forms  $T_n$  satisfying  conditions (i)--(ii)--(iii-b)   in Definition \ref{D:approximable} (1) are   positive pluriharmonic (resp. positive closed).

  \end{definition}

  Next, we introduce  some suitable classes  of currents. Recall  that  $X$ is a complex manifold of dimension $k$
  and $V\subset X$ is a submanifold of dimension $1\leq l<k.$
 
 \begin{definition}\rm \label{D:classes}
   Let $B$ be a relatively compact open subset of $V.$  Let $m,m'\in\N$ with $m\geq m'.$

 \smallskip

 \noindent  (1) We say that  $T$ is  {\it approximable  along $B$ by  $\Cc^m$-smooth positive plurisubharmonic forms} 
   if, there are  an open  neighborhood $U$ of $\overline B$  in $X$ 
  such that  $T\in\SH^{m}_p(U).$
 We denote  by  $\SH^{m}_p( B)$  the class of all positive $(p,p)$-currents which are approximable  along $B$ by  $\Cc^m$-smooth positive plurisubharmonic forms.

   Similarly,  we say that  $T$ is  {\it approximable  along $B$ by  $\Cc^m$-smooth positive pluriharmonic (resp. positive closed) forms}   
   if, in the  above definition  the approximating forms  $T_n$ are   positive pluriharmonic (resp. positive closed).
    We denote  by  $\PH_p^{m}( B)$  the class of all positive $(p,p)$-currents which are approximable  along $B$ by  $\Cc^m$-smooth positive pluriharmonic forms.
 Analogously, we denote  by  $\CL_p^{m}( B)$  the class of all positive $(p,p)$-currents which are approximable  along $B$ by  $\Cc^m$-smooth positive closed forms.

 \smallskip
 
 \noindent  (2) We say that  $T$ is  {\it approximable  along $B$ by  $\Cc^m$-smooth positive plurisubharmonic forms with  $\Cc^{m'}$-control on boundary} 
   if, there are  an open  neighborhood $U$ of $\overline B$  in $X$ and  an open neighborhood  $W $ of $\partial B$ in $X$  with $W\subset U$
  such that  $T\in\SH^{m,m'}_p(U,W).$
 We denote  by  $\SH^{m,m'}_p( B)$  the class of all positive $(p,p)$-currents which are approximable  along $B$ by  $\Cc^m$-smooth positive plurisubharmonic forms with  $\Cc^{m'}$-control on boundary.

   Similarly,  we say that  $T$ is  {\it approximable  along $B$ by  $\Cc^m$-smooth positive pluriharmonic (resp. positive closed) forms with  $\Cc^{m'}$-control on boundary}   
   if, in the  above definition  the approximating forms  $T_n$  are   positive pluriharmonic (resp. positive closed).

 We denote  by  $\PH_p^{m,m'}( B)$  the class of all positive $(p,p)$-currents which are approximable  along $B$ by  $\Cc^m$-smooth positive pluriharmonic forms with  $\Cc^{m'}$-control on boundary.
 Analogously, we denote  by  $\CL_p^{m,m'}( B)$  the class of all positive $(p,p)$-currents which are approximable  along $B$ by  $\Cc^m$-smooth positive closed forms with  $\Cc^{m'}$-control on boundary.

 \smallskip

 \noindent  (3) We say that  $T$ is  {\it approximable  along $B$ by  $\Cc^m$-smooth positive plurisubharmonic forms with  compact support along $B$} 
   if, there are  an open  neighborhood $U$ of $\overline B$  in $X$ and  an open neighborhood  $W $ of $\partial B$ in $X$  with $W\subset U$
  such that  $T\in\SH^{m}_p(U,W,\comp).$
 We denote  by  $\SH^{m}_p( B,\comp)$  the class of all positive $(p,p)$-currents which are approximable  along $B$ by  $\Cc^m$-smooth positive plurisubharmonic forms with  compact support along $B$.

   Similarly,  we say that  $T$ is  {\it approximable  along $B$ by  $\Cc^m$-smooth positive pluriharmonic (resp. positive closed) forms with  compact support along $B$}   
   if, in the  above definition  the approximating forms  $T_n$  are   positive pluriharmonic (resp. positive closed).

 We denote  by  $\PH_p^{m}( B,\comp)$  the class of all positive $(p,p)$-currents which are approximable  along $B$ by  $\Cc^m$-smooth positive pluriharmonic forms with  compact support along $B.$
 Analogously, we denote  by  $\CL_p^{m}( B,\comp)$  the class of all positive $(p,p)$-currents which are approximable  along $B$ by  $\Cc^m$-smooth positive closed forms with  compact support along $B$.

   \end{definition}
   
    \begin{definition} \label{D:approximable-bis}\rm
  Let $K$ be a  relatively compact  subset of $V.$  Let $m,m'\in\N$ with $m\geq m'.$

 \smallskip

 \noindent  (1) We say that  a current $T$ is  {\it approximable along $K$ by   $\Cc^m$-smooth positive plurisubharmonic forms (resp.   $\Cc^m$-smooth positive pluriharmonic forms, resp.    $\Cc^m$-smooth positive closed forms)} if, there is a relatively compact  open neighborhood $B$ of $K$ in $V$
  such  that  $T\in \SH^{m}_p(B)$   (resp.  $T\in \PH^{m}_p(B),$ resp. $T\in\CL^{m}_p(B)$).  In other words,
  $$
   \SH_p^{m}(K):=\bigcup_{B\ \text{open in }\ V:\  K\subset B}  \SH_p^{m}(B),
   $$
   and similarly  for  $\PH_p^{m}(K)$ and $\CL_p^{m}(K).$

 \smallskip

 \noindent  (2)
  We say that  a current $T$ is  {\it approximable along $K$ by   $\Cc^m$-smooth positive plurisubharmonic forms (resp.   $\Cc^m$-smooth positive pluriharmonic forms, resp.    $\Cc^m$-smooth positive closed forms)  with  $\Cc^{m'}$-control on boundary} if, there is a relatively compact  open neighborhood $B$ of $K$ in $V$
  such  that  $T\in \SH^{m,m'}_p(B)$   (resp.  $T\in \PH^{m,m'}_p(B),$ resp. $T\in\CL^{m,m'}_p(B)$).  In other words,
  $$
   \SH_p^{m,m'}(K):=\bigcup_{B\ \text{open in }\ V:\  K\subset B}  \SH_p^{m,m'}(B),
   $$
   and similarly  for  $\PH_p^{m,m'}(K)$ and $\CL_p^{m,m'}(K).$

 \smallskip

 \noindent  (3)
  We say that  a current $T$ is  {\it approximable along $K$ by   $\Cc^m$-smooth positive plurisubharmonic forms (resp.   $\Cc^m$-smooth positive pluriharmonic forms, resp.    $\Cc^m$-smooth positive closed forms)   with  compact support along $B$} if, there is a relatively compact  open neighborhood $B$ of $K$ in $V$
  such  that  $T\in \SH^{m}_p(B,\comp)$   (resp.  $T\in \PH^{m}_p(B,\comp),$ resp. $T\in\CL^{m}_p(B,\comp)$).  In other words,
  $$
   \SH_p^{m}(K,\comp):=\bigcup_{B\ \text{open in }\ V:\  K\subset B}  \SH_p^{m}(B,\comp),
   $$
   and similarly  for  $\PH_p^{m}(K,\comp)$ and $\CL_p^{m}(K,\comp).$
   
 \end{definition}
    
\subsection{Statement of the main general results}\label{SS:Main-General-Results}
 Recall  that  $X$ is a complex manifold of dimension $k$
  and $V\subset X$ is a submanifold of dimension $1\leq l<k.$
  Fix $0\leq p\leq k$ and  define $\lowm$ and $\upm$  by  \eqref{e:m}.   
 
The  vector bundle $\E$ (that is, the normal  bundle  to $V$ in $X$) is  endowed with a Hermitian  metric $h.$
Several  notions appearing   in the  following   theorems such as  (strongly) admissible maps, generalized Lelong numbers, tangent currents,  etc.  will be  defined in later  sections. In particular, the mass indicator $\nu_j(T,B,r,\tau ,\omega,h)$   (resp.  the mass indicator $\kappa_j(T,B,r,\tau ,\omega,h)$)  appearing in Theorem \ref{T:main_1} as well as    Theorems  \ref{T:main_2}, \ref{T:main_1'} and \ref{T:main_2'} below
 are naturally  defined using a  canonical average of $T$  along $B$   in the spirit of  the model formula  
 \eqref{e:Lelong-number-point} for  average mean  (resp.     the model formula  \eqref{e:Lelong-number-point-bis}  for  logarithmic mean).

 The first main result of the  article is  the following
 
 \begin{theorem}\label{T:main_1} {\rm (Tangent Theorem I) }
  Let $X,$ $V$ be as  above  and  suppose that $(V,\omega)$ is  K\"ahler,  
  and  that  $B$ is   a piecewise $\Cc^2$-smooth open subset  of $V$ and that  there exists a strongly   admissible map for $B.$
  Let $T$ be  a   positive plurisubharmonic  $(p,p)$-current  on a neighborhood of $\overline B$ in $X$ such that  $T=T^+-T^-$  for some $T^\pm\in\SH_p^{3,3}( B).$
  Then  the following  assertions  hold:
  \begin{enumerate}
  \item For every $\lowm\leq j\leq \upm,$  the following limit  exists and is  finite
  $$    \nu_j(T,B,\omega,h):=\lim\limits_{r\to 0+}\nu_j(T,B,r,\tau ,\omega,h)                              $$
   for  all strongly  admissible maps  $\tau$     for $B$ and for all  Hermitian  metrics  $h$ on $\E.$   
   
\item The  real numbers  $\nu_j(T,B,\omega,h)$ are intrinsic, that is, they are
  independent of  the choice  of   $\tau.$ 
  
  \item  The following  equality holds
  $$\lim\limits_{r\to 0+}\kappa_j(T,B,r,\tau,\omega,h)=\nu_j(T,B,\omega,h)$$ 
  for all $\lowm\leq j\leq \upm$ with $j>l-p,$  and  for all strongly  admissible maps $\tau$  for $B$ and for all Hermitian  metrics $h$ on $\E.$
  \item   $\nu_{\upm}(T,B,\omega,h)$ is  nonnegative. Moreover,  it   is  totally intrinsic, i.e.  it is independent of  the choice  of  both  $\tau$ and $h.$ So  we will denote it simply by   $\nu_{\upm}(T,B,\omega).$
  Moreover, it has  a  geometric meaning in the sense of Siu and Alessandrini--Bassanelli (see Theorem \ref{T:AB-2}).

  \item  If  $\tau$ is a holomorphic admissible map  and if  $T^\pm$ belong  only  to  the class $\SH_p^{2}(\overline B),$
  then the above four assertions (1)--(4) still hold for $j=\upm.$
   
   \item 
There  exists   tangent currents to $T$  along $B,$   and all tangent  currents $T_\infty$  are  positive plurisubharmonic   on $\pi^{-1}(B)\subset \E.$
Moreover,  $T_\infty$ are partially $V$-conic pluriharmonic on $\pi^{-1}(B)\subset \E$ in the  sense   that
the current   $T_\infty\wedge\pi^*(\omega^\lowm)$ is $V$-conic pluriharmonic on $\pi^{-1}(B)\subset \E.$

\item If  instead of the above  assumption on $T,$  we assume  that  $T$ is a   positive pluriharmonic  $(p,p)$-current  on a neighborhood of $\overline B$ in $X$ such that  $T=T^+-T^-$  for some $T^\pm\in\PH_p^{2,2}( B),$  then  all the above  assertions still  hold and  moreover every  tangent  current $T_\infty$ is also   $V$-conic pluriharmonic  on $\pi^{-1}(B)\subset \E.$
  \end{enumerate}
 \end{theorem}

 \begin{remark}\rm
 There are at least two interpretations of the limit
  $$\nu_j(T,B,\omega,h)=\lim\limits_{r\to 0+}\kappa_j(T,B,r,\tau,\omega,h)$$ 
 which correspond  to the model  interpretations  \eqref{e:Lelong-number-point-bisbis(1)}  and  \eqref{e:Lelong-number-point-bisbis(2)}.
 \end{remark}
 \begin{remark}\rm
At the first glance    the condition  $T=T^+-T^-$  looks artificial. However, it is unavoidable in practice, see Theorem \ref{T:approx-admiss} below. 
 \end{remark}
 
 Our second main result deals with positive closed  currents.  
 \begin{theorem}\label{T:main_2} {\rm (Tangent Theorem II) }
  Let $X,$ $V$ be as  above. Assume that there is a  Hermitian metric $\omega$ on  $V$      for which  $\ddc\omega^j=0$  for $1\leq  j\leq \upm-1.$ Assume
  also that  $B$ is   a piecewise $\Cc^2$-smooth open subset  of $V$ and that  there exists a strongly   admissible map for $B.$
  Let $T$ be  a   positive closed  $(p,p)$-current  on a neighborhood of $\overline B$ in $X$ such that  $T=T^+-T^-$  for some $T^\pm\in\CL_p^{2,2}( B).$
  Then  the following  assertions  hold:
   \begin{enumerate}
   \item For every $\lowm\leq j\leq \upm,$  the following limit  exists and is  finite
  $$    \nu_j(T,B,\omega,h):=\lim\limits_{r\to 0+}\nu_j(T,B,r,\tau ,\omega,h)                              $$
   for  all strongly  admissible maps  $\tau$     for $B$ and for all  Hermitian  metrics  $h$ on $\E.$   
   
\item The  real numbers  $\nu_j(T,B,\omega,h)$ are intrinsic, that is, they are
  independent of  the choice  of   $\tau.$ 
  
  \item  The following  equality holds
  $$\lim\limits_{r\to 0+}\kappa_j(T,B,r,\tau,\omega,h)=\nu_j(T,B,\omega,h)$$ 
 for all $\lowm\leq j\leq \upm$ with $j>l-p,$  and for all strongly  admissible maps $\tau$  for $B$ and for all Hermitian  metrics $h$ on $\E.$
  \item   $\nu_{\upm}(T,B,\omega,h)$ is  nonnegative. Moreover,  it is  totally intrinsic, i.e.  it is independent of  the choice  of both   $\tau$ and $h.$ So  we will denote it simply by   $\nu_{\upm}(T,B,\omega).$
  Moreover, it has  a  geometric meaning in the sense of Siu and Alessandrini--Bassanelli (see Theorem \ref{T:AB-2}).

  \item  If  $\tau$ is a holomorphic admissible map  and if  $T^\pm$ belong  only  to  the class $\CL_p^{1,1}(\overline B),$
  then the above four assertions (1)--(4) still hold.

   \item 
  There  exist   tangent currents  to  $T$  along $B$   and all  tangent  currents  $T_\infty$  are $V$-conic  positive  closed   on $\pi^{-1}(B)\subset \E.$
 
 \item  If  instead of the above  assumption on $\omega$ and  $T,$  we assume  that  the form $\omega$  is  K\"ahler and $T$ is a   positive closed  $(p,p)$-current  on a neighborhood of $\overline B$ in $X$ such that  $T=T^+-T^-$  for some $T^\pm\in\CL_p^{1,1}( B),$  then all the above  assertions still hold.  
 If  moreover $\tau$ is  holomorphic and  $T=T^+-T^-$  for some $T^\pm\in\CL_p^{1}( B),$ then the above four assertions (1)--(4) still hold for $j=\upm.$
  \end{enumerate}
 \end{theorem}

 \begin{definition}\label{D:Lelong-numbers}\rm 
 The value $\nu_j(T,B,\omega,h)$  is called {\it the $j$-th (generalized)  Lelong number of $T$ along $B.$}
 The   set of real numbers  $\{\nu_j(T,B,\omega,h):$\ $\lowm\leq j\leq \upm\}$  are called {\it  the Lelong numbers of $T$ along $B.$}  
The nonnegative  number $\nu_\upm(T,B,\omega):=\nu_\upm(T,B,\omega,h)$ is  called  {\it  the top Lelong numbers of $T$ along $V,$}
it is  also denoted by  $\nu_\top(T,B,\omega).$
\end{definition}

\begin{remark}\rm \label{R:Top-Lelong-vs-AB}
 In the context of    Theorem \ref{T:AB-1} where the admissible map is  identity, $\omega$ is the canonical K\"ahler metric on the linear subspace $V$
 and $h$  is the canonical Euclidean metric  on $\C^{k-l},$  the Alessandrini-Bassanelli's  Lelong number of $T$ coincides with our top Lelong  number of $T,$  that is, $\nu_\AB(T,B)=\nu_\top(T,B,\omega).$ 
 
 Our Lelong numbers $\nu_j(T,B,\omega,h)$ are indexed  by the same set  $\{j:\ \lowm\leq j\leq \upm\}$ as  Dinh--Sibony's cohomology classes  $\kappa_j(T)$  which  were previously introduced in Theorem \ref{T:Dinh-Sibony-first} (3).  In  a forthcoming work we will study the  relation between our Lelong numbers and  Dinh--Sibony's cohomology classes.
\end{remark}

\smallskip

 We  are  particularly interested  in   the  special but very important   case  where $\supp( T)\cap V$ is compact in $V.$
 In this case
 we can  choose  any piecewise smooth open  neighborhood $B$  of  $\supp( T)\cap V$ in $V$
 and define simply
 \begin{equation}\label{e:def_nu_j_T,V}
 \nu_j(T,V,\omega,h):=\nu_j(T,B,\omega,h).
 \end{equation}
 We will see later that   this definition is independent of the choice of  such a $B.$
 The  above main results yield following two important    applications. 
 The  first consequence  is for  positive  plurisubharmonic  currents.
 
 \begin{theorem}\label{T:main_1'} {\rm (Tangent Theorem I') }
  Let $X,$ $V$ be as  above and  suppose that $(V,\omega)$ is  K\"ahler.  Assume   that there   exists a strongly admissible map for $V.$
  Let $T$ be  a  positive plurisubharmonic  $(p,p)$-current  on $X$  
  such that
  $\supp (T )\cap V$ is compact. Assume  in addition that  on an open  neighborhood of $\supp(T)\cap V$ in $X,$ we have    $T=T^+-T^-$  for some $T^\pm\in\SH_p^{3}(\supp(T)\cap V,\comp).$
  Then  the following  assertions  hold:
  \begin{enumerate}
  \item For every $\lowm\leq j\leq \upm,$  the following limit  exists and is  finite
  $$    \nu_j(T,V,\omega,h):=\lim\limits_{r\to 0+}\nu_j(T,B,r,\tau ,\omega,h).                                $$
   Here, $B$ is a  piecewise smooth open  neighborhood   of  $\supp( T)\cap V$ in $V,$
   $\tau$  is a strongly   admissible map for $B$  and $h$ is as usual a  Hermitian  metric on $\E.$   
   Moreover, for all $\lowm\leq j\leq \upm$ with $j>l-p,$  we  also have
   $$
     \nu_j(T,V,\omega,h)=\lim\limits_{r\to 0+}\kappa_j(T,B,r,\tau,\omega,h).
   $$
   
\item The real numbers   $\nu_j(T,V,\omega,h)$   are intrinsic, that is, they are
  independent of  the choice  of $B$ and   $\tau.$ 
  
  \item   $\nu_{\upm}(T,V,\omega,h)$ is  nonnegative.  Moreover,  it  is totally intrinsic, i.e. it is independent of the
choice of  $B,$ $\tau$ and $h.$  So  we  denote it simply by $\nu_{\upm}(T,V,\omega).$  Moreover, it has a 
geometric meaning in the sense of Siu and
Alessandrini–Bassanelli (see Theorem \ref{T:AB-1}).

   \item 
There  exists   tangent currents  to $T$  along $V$   and  all tangent  currents   $T_\infty$ are  positive plurisubharmonic   on $\E.$
Moreover,  $T_\infty$ are partially $V$-conic pluriharmonic on $ \E$ in the  sense   that
the current   $T_\infty\wedge\pi^*(\omega^\lowm)$ is $V$-conic pluriharmonic on $ \E.$

\item If  instead of the above  assumption on $T,$  we assume  that  $T$ is a   positive pluriharmonic  $(p,p)$-current   in $X$ such that   $\supp (T )\cap V$ is compact and  that  on an open  neighborhood of $\supp(T)\cap V$ in $X,$ we have    $T=T^+-T^-$   for some $T^\pm\in\PH_p^{2}( \supp(T)\cap V,\comp),$  then  all the above  assertions still hold and  morever every  tangent  current $T_\infty$ is also   $V$-conic positive  pluriharmonic  on $\E.$
   
  \end{enumerate}
 \end{theorem}

 The  second consequence  is for  positive  closed  currents.
 
 \begin{theorem}\label{T:main_2'} {\rm (Tangent Theorem II') }
  Let $X,$ $V$ be as  above. Let $T$ be  a  positive closed  $(p,p)$-current  on $X$ 
  such that
  $\supp( T) \cap V$ is compact.   Assume   that  on an open  neighborhood of $\supp(T)\cap V$ in $X,$ we have    $T=T^+-T^-$  for some $T^\pm\in\CL_p^{2}(\supp(T)\cap V,\comp).$  
  Assume in addition that there is a  Hermitian metric $\omega$ on  $V$      for which  $\ddc\omega^j=0$  for $1\leq  j\leq \upm-1.$   
  Then the following assertions hold:
   \begin{enumerate}
  \item For every $\lowm\leq j\leq \upm,$  the following limit  exists and is  finite
  $$    \nu_j(T,V,\omega,h):=\lim\limits_{r\to 0+}\nu_j(T,B,r,\tau,\omega,h ).                                $$
   Here, $B$ is a  piecewise smooth open  neighborhood   of  $\supp( T)\cap V$ in $V,$ $\tau$  is a strongly admissible   map and  $h$ is a  Hermitian  metric on $\E.$ Moreover, for all $\lowm\leq j\leq \upm$ with $j>l-p,$  we  also have
   $$
     \nu_j(T,V,\omega,h)=\lim\limits_{r\to 0+}\kappa_j(T,B,r,\tau,\omega,h).
   $$

   \item The real numbers 
   $\nu_j(T,V,\omega,h)$  are intrinsic, that is, they are
  independent of  the choice  of $B$ and   $\tau.$  
  \item   $\nu_{\upm}(T,V,\omega,h)$ is  nonnegative. Moreover, it is totally intrinsic, i.e. it is independent of the
choice of $B,$ $\tau$ and $h.$   So  we  denote it simply by $\nu_{\upm}(T,V,\omega).$  Moreover, it has a 
geometric meaning in the sense of Siu and
Alessandrini–Bassanelli (see Theorem \ref{T:AB-1}).
  
   \item 
  There  exist   tangent currents  to  $T$  along $V$   and  all tangent   currents   $T_\infty$ are $V$-conic  positive  closed   on $\E.$

\item If  instead of the above  assumption on $\omega$ and  $T,$  we assume  that the form $\omega$ is  K\"ahler and  $T$ is a   positive  closed  $(p,p)$-current   in $X$ such that   $\supp (T )\cap V$ is compact and  that  on an open  neighborhood of $\supp(T)\cap V$ in $X,$ we have    $T=T^+-T^-$   for some $T^\pm\in\CL_p^{1}( \supp(T)\cap V,\comp),$  then  all the above  assertions still hold.  
   
  \end{enumerate}
 \end{theorem}

 \begin{remark}\rm
 Observe that the  condition  on $V$  in Theorem  \ref{T:main_2} is   weaker  than  that in Theorem \ref{T:main_1}.
 Namely, in  Theorem  \ref{T:main_2} we only require  $\ddc\omega^j=0$  for $1\leq  j\leq \upm-1,$  whereas  in Theorem \ref{T:main_1}  we assume that $\omega$ is  K\"ahler.  
 \end{remark}
 \begin{definition}\label{D:Lelong-numbers-compact-supp}\rm 
 The value $\nu_j(T,V,\omega,h)$  is called {\it the $j$-th Lelong number of $T$ along $V.$}
 The   set of real numbers  $\{\nu_j(T,V,\omega,h):$\ $\lowm\leq j\leq \upm\}$  are called {\it  the Lelong numbers of $T$ along $V.$}  
The nonnegative  number $\nu_\upm(T,V,\omega):=\nu_\upm(T,V,\omega,h)$ is  called  {\it  the top Lelong numbers of $T$ along $V,$}
it is  also denoted by  $\nu_\top(T,V,\omega).$
\end{definition}

There  are  two  assumptions appearing in the above main theorems, namely, the (strongly) admissible maps  and the  approximation
of  positive plurisubharmonic (resp. positive pluriharmonic, resp.  positive closed) $(p,p)$-currents.
The last main result shows  that when $X$ is  K\"ahler, these  conditions are fulfilled. 
This is a  consequence of   Appendix \ref{A:admissible_maps}  and Appendix \ref{A:approximation}. 
 \begin{theorem}\label{T:approx-admiss}
  Let $X,$ $V$ be as  above.   Assume that  $X$ is  K\"ahler.  Then,  for every  relatively compact open set $B\subset V,$
  the  following assertions hold:
  \begin{enumerate}
  \item There is a strongly  admissible map  for $B.$
   \item Let $m,m'\in\N$ with $m\geq m'.$  
   Let $T$ be a   positive plurisubharmonic (resp.  positive  pluriharmonic, resp.  positive closed) $(p,p)$-current   on $X$ which satisfies the following conditions (i)--(ii):
   \begin{itemize}
   \item[(i)] $T$ is  of class $\Cc^{m'}$ near  $\partial B;$
   \item[(ii)] There is   a relatively compact open subset $\Omega$ of $X$ with  $B\Subset \Omega$
   and $dT$ is of class $\Cc^0$ near $\partial\Omega.$
   \end{itemize}
  Then  $T$ can be written  in an open neighborhood of $\overline B$ in $X$ as $T=T^+-T^-$  for some $T^\pm  \in  \SH^{m,m'}_p( B)$ (resp. $T^\pm  \in  \PH^{m,m'}_p(\overline B),$ $T^\pm  \in  \CL^{m,m'}_p( B)$).
\end{enumerate}  \end{theorem}

 \begin{remark}\rm
  \label{R:Dinh-Sibony}
  In view of Theorem  \ref{T:approx-admiss} below, the assumption of  Theorem  \ref{T:main_2'} is close to that of  Theorem \ref{T:Dinh-Sibony-first}.  So  Theorem  \ref{T:main_2'} may be regarded  as a numerical  complement to the original result of Dinh--Sibony  when  the ambient manifold $X$ is compact K\"ahler.
 \end{remark}

\begin{remark}\label{R:Vu}\rm
 It is  worth  noting that Theorem \ref{T:main_2'} (4) improves, in some  sense,  Vu's recent result \cite[Theorem 1.1]{Vu19a} (see Theorem \ref{T:Vu} below).  Indeed, instead of  the  Hermitian form $\omega$ on $V,$ Vu  assumes  
the following slightly stronger condition:
 there is a Hermitian  metric $\hat\omega$ on $X$ 
 for which  $\ddc\hat\omega^j=0$ on $V$ for $1\leq  j\leq k-p-1.$  Setting $\omega:=\hat\omega|_V,$ we  get the Hermitian  metric $\omega$ needed for Theorem \ref{T:main_2}.
 
 However, Vu does not need   that  $T$ is approximable by the difference of positive closed smooth forms  along $\supp(T)\cap V.$
 Moreover, he only needs that $\tau$ is  an admissible map,  whereas $\tau $ is  strongly  admissible in 
 Theorem \ref{T:main_2'}.  
 \end{remark}

 To  end  this  subsection, we  record the  following two corollaries which capture the essential points of the above  main results in 
 the special  but important context where the ambient manifold $X$ is  K\"ahler. We think that these  explicit statements will be  useful in practice.
 
 Our first main corollary concerns positive plurisubharmonic and  positive pluriharmonic currents.  
 \begin{corollary}\label{C:main_1} {\rm (Tangent Corollary I) }
  Let $X,$ $V$ be as  above  and  suppose that $X$ is  K\"ahler and that  $(V,\omega)$ is  K\"ahler,  
  and  that  $B$ is   a piecewise $\Cc^2$-smooth open subset  of $V.$ 
  Let $T$ and $T^\pm$  be  three  positive plurisubharmonic  $(p,p)$-currents  on a neighborhood of $\overline B$ in $X$ such that $T=T^+-T^-$ and that  
  \begin{itemize}
  \item [(i)] $T^\pm $  is of class $\Cc^3$
  in a neighborhood of $\partial  B$ in $X;$
   \item[(ii)] There is   a relatively compact open subset $\Omega$ of $X$ with  $B\Subset \Omega$
   and $dT^\pm$ is of class $\Cc^0$ near $\partial\Omega.$
  \end{itemize}
  Then  the following  assertions  hold:
  \begin{enumerate}
  \item For every $\lowm\leq j\leq \upm,$  the following limit  exists and is  finite
  $$    \nu_j(T,B,\omega,h):=\lim\limits_{r\to 0+}\nu_j(T,B,r,\tau ,\omega,h)                              $$
   for  all strongly  admissible maps  $\tau$     for $B$ and for all  Hermitian  metrics  $h$ on $\E.$   
   
\item The  real numbers  $\nu_j(T,B,\omega,h)$ are intrinsic, that is, they are
  independent of  the choice  of   $\tau.$ 
  
  \item  The following  equality holds
  $$\lim\limits_{r\to 0+}\kappa_j(T,B,r,\tau,\omega,h)=\nu_j(T,B,\omega,h)$$ 
   for all $\lowm\leq j\leq \upm$ with $j>l-p,$ and for all strongly  admissible maps $\tau$  for $B$ and for all Hermitian  metrics $h$ on $\E.$
  \item   $\nu_{\upm}(T,B,\omega,h)$ is  nonnegative. Moreover, is  totally intrinsic, i.e.  it is independent of  the choice  of  both $\tau$ and $h.$ So we will denote it simply by $\nu_{\upm}(T,B,\omega).$ 
  Moreover, it has  a  geometric meaning in the sense of Siu and Alessandrini--Bassanelli (see Theorem \ref{T:AB-2}).

  \item  If  $\tau$ is a holomorphic admissible map  
  then the above four assertions (1)--(4) still hold for $j=\upm.$
   
   \item 
There  exists   tangent currents to $T$  along $B,$   and all tangent  currents $T_\infty$  are  positive plurisubharmonic   on $\pi^{-1}(B)\subset \E.$
Moreover,  $T_\infty$ are partially $V$-conic pluriharmonic on $\pi^{-1}(B)\subset \E$ in the  sense   that
the current   $T_\infty\wedge\pi^*(\omega^\lowm)$ is $V$-conic pluriharmonic on $\pi^{-1}(B)\subset \E.$

\item If  instead of the above  assumption on $T,$  we assume  that  $T$ and $T^\pm$  are three   positive pluriharmonic  $(p,p)$-currents  on a neighborhood of $\overline B$ in $X$ such that $T=T^+-T^-$ and  that  \begin{itemize} \item  $T^\pm $  is of class $\Cc^2$
  in a neighborhood of $\partial  B$ in $X;$ 
  \item there is   a relatively compact open subset $\Omega$ of $X$ with  $B\Subset \Omega$
   and $dT^\pm$ is of class $\Cc^0$ near $\partial\Omega,$
   \end{itemize}
  then  all the above  assertions still  hold and  moreover every  tangent  current $T_\infty$ is also   $V$-conic pluriharmonic  on $\pi^{-1}(B)\subset \E.$
  \end{enumerate}
 \end{corollary}
 
 Our second main corollary deals with positive closed  currents.  
 \begin{corollary}\label{C:main_2} {\rm (Tangent Corollary II) }
  Let $X,$ $V$ be as  above. Assume that $X$ is   K\"ahler and  that there is a  Hermitian metric $\omega$ on  $V$      for which  $\ddc\omega^j=0$  for $1\leq  j\leq \upm-1.$ Assume
  also that  $B$ is   a piecewise $\Cc^2$-smooth open subset  of $V.$  
  Let $T$ and $T^\pm$  be  three   positive closed  $(p,p)$-currents  on a neighborhood of $\overline B$ in $X$ such that $T=T^+-T^-$ and that    $T^\pm $  is of class $\Cc^2$
  in a neighborhood of $\partial  B$ in $X.$ 
  Then  the following  assertions  hold:
   \begin{enumerate}
   \item For every $\lowm\leq j\leq \upm,$  the following limit  exists and is  finite
  $$    \nu_j(T,B,\omega,h):=\lim\limits_{r\to 0+}\nu_j(T,B,r,\tau ,\omega,h)                              $$
   for  all strongly  admissible maps  $\tau$     for $B$ and for all  Hermitian  metrics  $h$ on $\E.$   
   
\item The  real numbers  $\nu_j(T,B,\omega,h)$ are intrinsic, that is, they are
  independent of  the choice  of   $\tau.$ 
  
  \item  The following  equality holds
  $$\lim\limits_{r\to 0+}\kappa_j(T,B,r,\tau,\omega,h)=\nu_j(T,B,\omega,h)$$ 
  for all $\lowm\leq j\leq \upm$ with $j>l-p,$ and  for all strongly  admissible maps $\tau$  for $B$ and for all Hermitian  metrics $h$ on $\E.$
  \item   $\nu_{\upm}(T,B,\omega,h)$ is  nonnegative. Moreover,  it is  totally intrinsic, i.e.  it is independent of  the choice  of both    $\tau$ and $h.$  So we will denote it simply by $\nu_{\upm}(T,B,\omega).$
  Moreover, it has  a  geometric meaning in the sense of Siu and Alessandrini--Bassanelli (see Theorem \ref{T:AB-2}).

  \item  If  $\tau$ is a holomorphic admissible map,
  then the above four assertions (1)--(4) still hold for $j=\upm.$

   \item 
  There  exist   tangent currents  to  $T$  along $B$   and all  tangent  currents  $T_\infty$  are $V$-conic  positive  closed   on $\pi^{-1}(B)\subset \E.$
 
 \item  If  instead of the above  assumption on $\omega$ and  $T,$  we assume  that  the form $\omega$  is  K\"ahler and $T$ is a   positive closed  $(p,p)$-current  on a neighborhood of $\overline B$ in $X$ such that     $T $  is of class $\Cc^1$
  in a neighborhood of $\partial  B$ in $X,$  then all the above  assertions still hold.  
 If  moreover $\tau$ is  holomorphic, then the above four assertions (1)--(4) still hold for $j=\upm.$
  \end{enumerate}
 \end{corollary}

 \subsection{Organization of the  article}
 
 The  article is divided into two parts.
 The first part, which covers the first 16  sections,  is devoted to the generalized Lelong numbers.
 The  second part, which consists of the last 8 sections,  studies the geometric  charaterizations of these  characteristic numbers.
 More concretely,   the  article is organized  as follows.
 
 \smallskip
 
 In  Section \ref{S:preparatory_results} below   we set up the  background and  introduce some  main objects, important  definitions which will be  used throughout the  article.  More specifically,  we  first  recall some basic  definitions of  Dinh-Sibony  \cite{DinhSibony18} such as normal vector bundles,  $V$-conic  currents, admissible maps, tangent currents  and   review  quickly their results  as well as a recent result of  Vu \cite{Vu19a}.  Next we  introduce  our definition of strongly admissible maps. The section is concluded with the introduction of some fundamental forms and  our notion of tubes  which generalizes
 that of Alessandrini--Bassanelli given in  \eqref{e:Tube-AB}.

 In  Section \ref{S:Lelong}
 we  state the first   collection  of main  results. These  results  are, in some sense,  more specialized and more detailed  than  the main general 
 results
 stated  in Subsection \ref{SS:Main-General-Results}. This  collection can be  divided into  two groups.
 The  first  group consists of four theorems which  consider positive closed  currents.
 The  second group consists of three theorems which  handle positive plurisubhamonic currents.
 The section is ended  with a short interpretation of our result in the context of Alessandrini--Bassanelli \cite{AlessandriniBassanelli96}.
 Even in this context, we obtain   new   results.
 
 In Section \ref{S:Lelong-Jensen}  we present the main tool  of our   method:  {\it   Lelong--Jensen formulas for   tubes  in a  vector bundle.}
 These formulas  arise in connection with  the  generalization of the classical Lelong-Jensen formula  for a ball in $\C^k.$
 Both Lelong's and  Skoda's formulations  of the Lelong number at a single point  rely on   the latter formula.
 In  comparison  with  Euclidean  balls,   our  tubes have  not only  horizontal  boundary, but also  vertical  boundary.
 A typical  feature of  our new  Lelong-Jensen formulas  for  tubes is  the presence of  vertical  boundary terms which are
linked  to  the  vertical  boundary.    This  section is  devoted to    Lelong-Jensen formulas for tubes in     abstract  context as well as  
in  concrete applications.  The formulas  are, in fact,  applied to  various  objects:
closed  currents,   currents with compact support, currents which are full in vertical directions   etc.
We  also give estimates  for the vertical boundary terms.

In Section  \ref{S:Basic-forms-and-CV-test}
we introduce  some basic  forms  for  the    bundle $\E.$ 
 We also prove  a  convergence test. These forms and  this  test will be  used  throughout   this  work.

 Section \ref{S:closed-holomorphic}  gives the proof of some of  the main theorems  in a   special  situation. Namely, 
 we deal  with positive  closed  currents  and  we assume  that there exists a  holomorphic  admissible map.
 Here,  some basic  ideas are carefully  explained  in  such  a  particular case. This  case is  simpler and it suggests also  how to manage the general case  of positive  closed  currents  with  non-holomorphic  admissible maps. 

 Section \ref{S:Regularization} develop technical  tools  which will be  used throughout the article.
 We first  introduce the Extended Standing Hypothesis which is  a technically complete version of the Standing Hypothesis.
 Next, we introduce  the representative current $T^\hash$ living on $\E$ of a  positive  current $T$ living on $X.$
 The section is then devoted to the  study of admissible estimates, that is,  estimates  which are related  to  admissible maps.
 Basic   individual admissible  estimates are obtained before more sophisticated and abstract estimates for  wedge-product are  established.

 In Section \ref{S:Positive-closed-currents}  we prove three of  the four theorems  in  the first  collection of the main results which  concern positive  closed currents.
 The two first  subsections introduce and study some new mass indicators. The major result of the  section is Theorem \ref{T:Lc-finite} on the  finiteness of these
 mass indicators. Using this  result, we  establish  in the  third subsection the existence of the generalized Lelong  numbers. The fourth subsection is then  devoted to  the proof that these  characteristic numbers  are in fact   independent of the choice of a (strongly) admissible map.
 As a by-product,  we introduce  some  variants $\hat\nu$ of the Lelong numbers $\nu$  which are always  non-negative.

 It is  classical (see \cite{Lelong,Lelong68}) that the  Lelong mean $\nu(T,x,r)$  of  a positive closed current $T$ at a point $x$
 (see formula \eqref{e:Lelong-number-point}) is a non-negative valued   increasing function in the radius $r.$ 
 Section \ref{S:quasi-monotone}  establishes   analogous  properties  for 
 the  generalized  Lelong numbers  of a  positive  closed  current.  In this new general context, we only achieve 
 a quasi-positivity and a quasi-monotonicity  of a positive linear combination  of  the  generalized  Lelong numbers. However,  this  seems the best property that  we may hope for.
 
 Section  \ref{S:psh-holomorphic} studies positive plurisubharmonic currents in   a    special  setting. Namely, 
  we suppose   that there exists a  holomorphic  admissible map.
 In comparison with the case of  positive  closed  currents  with a holomorphic  admissible map  treated in  Section  \ref{S:closed-holomorphic},
 the new difficulty  here is how to deal  with  the $\ddc$-part (that is,  the current $\ddc T$)  of a given positive plurisubharmonic  current $T.$  
 Basic  ideas and techniques  are well presented   in order to tackle  with   this obstacle. 
 This  study suggests us   how to treat the general case  of positive  plurisubharmonic  currents  with  non-holomorphic  admissible maps.

 Section  \ref{S:Admissible-psh}
sets up the necessary machinary for   admissible  estimates   and for positive plurisubharmonic  currents $T.$ 
These  estimates are much more  difficult than those  for positive  closed  currents since we have to  deal with  the influence of (strongly) admissible maps on  the  curvature current $\ddc T,$ whereas  this term  vanishes  automatically  when $T$ is closed. 
 In the  two first subsections, we introduce pointwise   admissible  estimates, negligible  test forms. Using  this, we next develop  basic   volume estimates which relate
 the values  of $T$ on  test forms  to  the  generalized  Lelong numbers.
 Next, we   establish a  basic  boundary  formula  using Stokes' Theorem (see Proposition \ref{P:Stokes}). Let us explain briefly this subsection.
  Let $\tau$ be a (strongly) admissible map and let $T$ be  a $(p,p)$-current of order $0$  such that $\ddc T$ is  also  a current of order $0.$ 
 This formula expresses the  difference $\ddc  (\tau_*T) -\tau_*(\ddc T)$ on a  tube in terms of some boundary  integrals.
 Note that understanding   the  above  difference is a  key    problem,   since   in general there is no holomorphic admissible map
 and hence the above  difference is in general  non trivial. The last three subsections are then devoted to  estimate  these  boundary terms  using
 the so-called  boundary  differential operators. Combining
 all the tools   developed so far,  Proposition \ref{P:basic-bdr-estimates} is the major technical result of this  section.

  Section \ref{S:Positive-psh-currents} has two  purposes. The first  one is  to establish
  some  abstract estimates on  the  difference $\ddc  (\tau_*T) -\tau_*(\ddc T)$ on a  tube in terms of some mass indicators modelling  the Lelong average means of $T$ and of $\ddc T.$
  For this  purpose we rely on   the result of Section  \ref{S:Admissible-psh}.
 Roughly speaking,  these estimates say that   this  difference is small when  the radius of the tube is  small.
 Using  these  inequalities,  the second  purpose of the section is  to achieve  a quasi-positivity and a quasi-monotonicity of  the generalized  Lelong numbers  of a  positive  plurisubharmonic  current. 
 So this  is  a generalization of  Section  \ref{S:quasi-monotone}.

 In Section \ref{S:finitness-Kc-Lc}  we prove  the last theorem  (which is also the  most  important one) in  the second collection of main results: Theorem \ref{T:Lelong-psh}.
 This  result treats  the class of positive plurisubharmonic currents with non-holomorphic admissible maps.
 In the  first  subsection we  introduce  some new global mass indicators which capture not only the mass of $T$ but also  the mass of its curvature $\ddc T .$   Using   Lelong-Jensen  formulas, we  study  these  mass  indicators in the two next subsections.
 The main result of the  section is Theorem \ref{T:Lc-finite-psh} on the  finiteness of these
 mass indicators. 
 As a by-product,  we show  in Theorem \ref{T:vanishing-Lelong} that if $T$ is a positive  plurisubharmonic  current with some suitable  additional assumptions, then all  the generalized Lelong  number of the positive  closed current $\ddc T$ vanish. 
 Based on this  development, the last subsection  is devoted to the existence of the generalized Lelong  numbers.

 Section  \ref{S:non-Kaehler-metrics} is devoted to the proof of the last theorem in  the first  collection of main results: Theorem \ref{T:Lelong-closed}. This  theorem deals with the class of positive closed currents, but  the admissible map in question is not  holomorphic, and the metric on the submanifold $V$ in question
 is not K\"ahler. The non-K\"ahlerity  of  the metric forces  us  to adapt the  method developed in Sections 
  \ref{S:Positive-psh-currents}
 and \ref{S:finitness-Kc-Lc} in a rather delicate  situation.  More  concretely, we achieve 
  some  abstract estimates on  the  difference $\dbar  (\tau_*T) -\tau_*(\dbar T)$ on a  tube in terms of some mass indicators modelling  the Lelong average means of $T.$ 
 These estimates assert that   this  difference is small when  the radius of the tube is  small.
 
Section \ref{S:Existence-tangent-currents} establishes the existence of tangent  currents 
in the following  three classes of currents: 
positive closed currents, positive pluriharmonic  currents and plurisubharmonic   currents. Here, the idea is to combine a  local analysis  and  the  finiteness  of the mass indicators modeling  the Lelong numbers
 which was previously obtained in Theorem \ref{T:Lc-finite}  and Theorem \ref{T:Lc-finite-psh}.

 Section \ref{S:V-conic-and-plurihar}, which is the last section of  Part 1,  describes  basic properties of  the  tangent currents such as the $V$-conicity  and the (partial)-pluriharmonicity.
This is  a  consequence of  our  Lelong-Jensen formulas for vector bundles which are  applied to  the  tangent currents.

\smallskip

Part 2 of the  article investigates the  geometric  characterizations of the generalized  Lelong numbers in the  spirit of  Siu \cite{Siu} and Alessandrini--Bassanelli  \cite{AlessandriniBassanelli96}.
Section  
\ref{S:Grassmannian} introduces Grassmannian  bundles associated to the normal  vector bundle $\E$ as well as   some canonical  projection
$\Pi_j,$ $\Pr_j$ and  some canonical  vertical  forms  $\alpha_\ver,$ $\beta_\ver,$ $\Upsilon_j.$ This notion  is  a generalized version of
the  blow-ups. We  reformulate  some  important  identities relating
these  objects.
These  identities  are
due  to Siu \cite{Siu}  in  the context of a single point and  to  Alessandrini--Bassanelli  \cite{AlessandriniBassanelli96} in the context of a linear complex  subspace.

Section  
 \ref{S:Flatness} begins  with  a review of  basic  notions and results in Complex  Geometry regarding $\C$-flat currents, $\C$-normal currents
 and their extension properties through analytic subsets such as   Federer type  theorems  (see \cite{Bassanelli, AlessandriniBassanelli93,AlessandriniBassanelli96, Sibony}).   Next,   we deal with the extension property of   some  currents in the presence of  a holomorphic admissible  map.  
 
 Based on the previous section, Section  \ref{S:Charac-pos-closed-and-plurihar-I} characterizes the generalized  Lelong numbers geometrically  when the admissible map in question is  holomorphic and  the currents in question are  either positive closed  or positive  pluriharmonic.

 In order to treat the general case of non-holomorphic admissible maps, 
  Section  \ref{S:Charac-pos-closed-and-plurihar-II} uses   the  finiteness  of the mass indicators modeling  the Lelong numbers. 
  Consequently, we  can show that our currents have similar extension properties as in the case of holomorphic admissible maps.
 Based on this  remarkable fact,  we give geometric characterizations of the generalized  Lelong numbers for two classes of currents: positive  closed currents and  positive  pluriharmonic  currents. Roughly  speaking, each (generalized)  Lelong number of a  current $T$ in each one of these two classes is expressed as  the mass of a  suitable cut-off current on the exceptional fiber in a suitable  Grassmannian   bundle. 


 Section \ref{S:Charac-psh-hol}
 extends the result of  Section  \ref{S:Charac-pos-closed-and-plurihar-II}  
 to the top Lelong  number  for the  class of  positive  plurisubharmonic  currents. In this  general  context,
 there is a new  phenomenon: to each  positive  plurisubharmonic  current $T$  we associate two cut-off currents $T^{(1)}$ and $S^{(0)}.$
 The top Lelong number of $T$ is expressed as  the sum of  the masses of these two cut-off currents on suitable exceptional fibers in  corresponding  two Grassmannian   bundles.  The analysis of the currents considered in this  section is much harder than that of the previous section.

Section  
 \ref{S:Charac-psh}  completes  Section \ref{S:Charac-psh-hol} by treating  the general case of non-holomorphic admissible maps.
 
Section \ref{S:Top-Lelong-strong-intrinsic} states and proves the second collection of main results: if the current  $T$ is positive closed
(resp.  positive  pluriharmonic, resp. positive  plurisubharmonic) and it satisfies some suitable approximation property, then 
its top  Lelong number is totally intrinsic.
 The proof is based on the geometric  description of the top Lelong number obtained in the previous  two sections.

 Section \ref{S:Proofs} gives the proof of the  main general theorems stated  in Subsection \ref{SS:Main-General-Results} by combining  the two collections  of main results.
We discuss the particular case where $\dim V=0,$ that is, $V$ is a single point.  We also study the dependence of the generalized Lelong numbers on
the   Hermitian metric $\omega$ on $V$ and on the   Hermitian  metric $h$ on the normal  bundle $\E.$
The section concludes  with  some open questions and further remarks.

 The articles ends with two appendices.
 
 Appendix A construct  strongly  admissible maps when the ambient manifold is  K\"ahler following  the idea of  Dinh--Sibony in \cite{DinhSibony18b}.
 Appendix B discusses  various  approximation results for  three classes of currents: positive closed currents,  positive pluriharmonic  currents and  positive plurisubharmonic  currents. The latter  appendix is inspired by  another work of  Dinh--Sibony in \cite{DinhSibony04}. 
 
 \smallskip
 
\noindent
{\bf Acknowledgments. }   The    author   acknowledges support from the Labex CEMPI (ANR-11-LABX-0007-01)
and from  the project QuaSiDy (ANR-21-CE40-0016).
The paper was partially prepared 
during the visit of the  author at the Vietnam  Institute for Advanced Study in Mathematics (VIASM). He would like to express his gratitude to this organization 
for hospitality and  for  financial support.

 \part{The generalized Lelong numbers and the tangent theorems}

\section{Preparatory results}  \label{S:preparatory_results}

\subsection{Currents  and positive  currents}\label{SS:currents}

 Let $M$ be a  complex manifold of dimension $k.$
A $(p,p)$-form on  $M$  is {\it
  positive} if it can be written at every point as a combination with
positive coefficients of forms of type
$$i\alpha_1\wedge\overline\alpha_1\wedge\ldots\wedge
i\alpha_p\wedge\overline\alpha_p$$
where the $\alpha_j$ are $(1,0)$-forms. A $(p,p)$-current or a $(p,p)$-form $T$ on $M$ is
{\it weakly positive} if $T\wedge\varphi$ is a positive measure for
any smooth positive $(k-p,k-p)$-form $\varphi$. A $(p,p)$-current $T$
is {\it positive} if $T\wedge\varphi$ is a positive measure for
any smooth weakly positive $(k-p,k-p)$-form $\varphi$. 
If $M$ is given with a Hermitian metric $\beta$ and $T$ is  a positive  $(p,p)$-current on $M,$
$T\wedge \beta^{k-p}$ is a positive measure on
$M$. The mass of  $T\wedge \beta^{k-p}$
on a measurable set $A$ is denoted by $\|T\|_A$ and is called {\it the mass of $T$ on $A$}.
{\it The mass} $\|T\|$ of $T$ is the total mass of  $T\wedge \beta^{k-p}$ on $M.$
A $(p, p)$-current $T$ on $M$  is {\it strictly positive} if  we have locally $T \geq \epsilon \beta^p ,$ i.e., $T -\epsilon \beta^p$
is positive, for some constant $\epsilon > 0.$ The definition does not depend on the
choice of $\beta.$
 
 Let $T$ be  a current of bidegree $(p,p)$ on an open set $U\subset\C^k.$  Write 
\begin{equation}\label{e:expression_currents}T=i^{p^2}\sum T_{I,J} dx_I\wedge d\overline x_J\quad\text{with}\quad T_{I,J}\quad \text{a distribution on}\quad U,
\end{equation}
the sum being taken all  over all  multi-indices $I,J$ with $|I|=|J|=p.$
Here, for a multi-index $I=(i_1,\ldots,i_p)$  with $1\leq i_1\leq \ldots\leq i_p\leq k,$  $|I|$ denotes the length $p$ of $I,$   $dx_I$ denotes $dx_{i_1}\wedge \ldots  dx_{i_p}$  and
 $d\bar x_I$ denotes $d\bar x_{i_1}\wedge \ldots  d\bar x_{i_p}.$

\begin{proposition}\label{P:Demailly}{\rm (see  e.g.  \cite[Proposition 1.14]{Demailly}) }
 Let  $T=i^{p^2}\sum T_{I,J} dx_I\wedge d\overline x_J$ of bidegree $(p,p)$ be  a  positive  current on an open  set in $\C^k.$  Then its coefficients $T_{I,J}$ are complex measures
 and satisfy  $\overline T_{I,J}=T_{J,I}$ for all  multi-indices $|I|=|J|=p.$  Moreover,  $T_{I,I}\geq 0,$ and the absolute values $|T_{I,J}|$
 of the measure $T_{I,J}$ satisfy the  inequality
 $$
 \lambda_I\lambda_J|T_{I,J}|\leq 2^{k-p}\sum_M  \lambda_M^2 T_{M,M},  \qquad  I\cap J\subset M\subset I\cup J,
 $$
 where $\lambda_j\geq 0$   are  arbitrary   coefficients and $\lambda_I=\prod_{j\in I}   \lambda_j.$
\end{proposition}

 The following  elementary lemma whose  proof is left to the interested reader says that  any set of  positive currents with uniformly bounded mass is weakly relatively compact in the weak-$\star$  topology.
 \begin{lemma}
 \label{L:subsequence}
 Let $(R_n)_{n=0}^\infty$ be  a  sequence of positive $(p,p)$-currents on an  open set $\Omega\subset \C^k$ such that
 $$
 \sup_{n\in\N} \int_\Omega  R_n\wedge (\ddc \|x\|^2)^{k-p}  <\infty.
 $$
 Then there exists a subsequence $(R_{N_n})_{n=0}^\infty$ and a positive current $R$ on $\Omega$  such that $\lim_{n\to\infty}R_{N_n}=R$ weakly in $\Omega.$
\end{lemma}

 Let $R$ be  a current with measure coefficients (or equivalently,  of order $0$) on an open set $\Omega$  in a complex  manifold $X$ of dimension $k.$
 Let $W$ be a relatively compact open subset of $\Omega$ and $\Phi$ a smooth test form  on $\Omega,$ we will  write
 \begin{equation}\label{e;cut-off}
 \int_W  R\wedge \Phi:=  \langle  R,\ind_W\Phi\rangle,
 \end{equation}
 where $\ind_W$ is  the  characteristic  function of $W.$ Let $(R_n)_{n=1}^\infty$ be a  sequence of positive  currents on $\Omega$ such that
 $\lim_{n\to\infty}  R_n=R$  weakly on  $\Omega,$ then  we  see that
 \begin{equation}\label{e:continuity-cut-off}
 \lim_{n\to\infty}\int_W  R_n\wedge \Phi=\int_W  R\wedge \Phi
 \end{equation}
 for every smooth test form $\Phi$  on $\Omega$ and  every relatively compact open subset $W\subset \Omega$ with $\|R\|(\partial W)=0.$
 Here, $\partial W$ is  the topological boundary of $W$ and $\|R\|$ is  the mass-measure of $R.$
 Consequently, if $K$  is a compact subset of $\Omega$ and $(W_i)_{i\in I}$ is a family of open subsets of $\Omega$ such that
 $K\subset  W_i$ for all $i\in I$ and $\partial W_i\cap \partial W_{j}=\varnothing$ for $i\not=j,$ then  we have
 \begin{equation} \label{e:except-countable}
   \lim_{n\to\infty}\int_{W_i}  R_n\wedge \Phi=\int_{W_i}  R\wedge \Phi
 \end{equation}
for every smooth test form $\Phi$  on $\Omega$ and  every $i\in I$  except  for a countable subset of $I.$

 In this  article we are concerned with the following  notion of weak convergence of quasi-positive currents.
 \begin{definition}\label{D:quasi-positivity}\rm
  We say that   a    current $R$ defined on $\Omega$  is {\it  quasi-positive} if, for every $x\in \Omega,$ there are an open neighborhood $\Omega_x$ of $x$  in $\Omega$ and  a $\Cc^1$-diffeomorphism $\tau_x$
 of $\Omega_x$ such that $\tau_x^*R$ is a  positive current.
 
 We say that   a  sequence of  currents $(R_n)_{n=1}^\infty$ {\it converge  in the sense of quasi-positive currents  on $\Omega$ to a current $R$}
  if   for every $x\in \Omega,$ there are an open neighborhood $\Omega_x$ of $x$  in $\Omega$ and  a $\Cc^1$-diffeomorphism $\tau_x$
 of $\Omega_x$  and two sequences of  positive  currents $(T^\pm_n)_{n=1}^\infty$  on $\Omega_x$   such that all currents $\tau_x^*(R_n-R)= T^+_n-T^-_n$   and that both sequences $T_n^\pm$ converge  weakly to  a common positive current  $T$ on $\Omega_x.$  
 \end{definition}
 
 The relevance of this notion is justified by the following simple result.
 \begin{lemma}\label{L:quasi-positivity}
  If  a  sequence of  currents $(R_n)_{n=1}^\infty$  converge  in the sense of quasi-positive currents  on $\Omega$ to a current $R,$ then 
  both  \eqref{e:continuity-cut-off} and \eqref{e:except-countable}  hold.
 \end{lemma}
 \proof
 Since  problem  is  local, we  are reduced to the 
 situation   where there is a   $\Cc^1$-diffeomorphism $\tau$
 on $\Omega$ such that all currents $\tau^*R_n$ are positive. Applying  \eqref{e:continuity-cut-off} and  \eqref{e:except-countable} to $\tau^*R_n,$ 
 the result follows.
 \endproof

\subsection{Normal bundle and admissible maps}\label{SS:admissible-maps}

Let $X$  be a complex manifold of dimension $k.$   Let  $V$ be a  smooth complex submanifold of
$X$ of dimension $l.$  Let $\E$ be the normal vector bundle to $V$ in $X.$ 

Consider a point $x\in V.$ If $\Tan_x(X)$  and $\Tan_x(V)$ denote, respectively, the tangent spaces of $X$ and of $V$ at $x,$
the  fiber $\E_x$ of $\E$ over $x$ is canonically  identified  with the quotient space $\Tan_x(X)/\Tan_x(V).$

For $\lambda \in \C^\ast ,$ let $A_\lambda :\ \E \to  \E$ be the multiplication by $\lambda$ in fibers of $\E,$
that is, 
\begin{equation}\label{e:A_lambda} A_\lambda(y):=\lambda y\qquad\text{for}\qquad y\in \E.
\end{equation}
A  current $T$ on $\E$ is  said to be {\it $V$-conic} if $T$ is  invariant under  the action of  $A_\lambda,$ that is,
$(A_\lambda)_*T=T$ for all $\lambda\in\C^*.$


The  following notion,  introduced  by  Dinh--Sibony \cite{DinhSibony18}, plays a vital role in their tangent theory for positive closed currents.  
\begin{definition}\label{D:admissible-maps} {\rm  (See \cite[Definitions 2.15 and  2.18]{DinhSibony18}) } \rm   
 Let $B$ be a relatively compact nonempty open subset of $V.$
An admissible map   along  $B$  is   a $\Cc^1$-smooth diffeomorphism $\tau$  from an open neighborhood  $U$ of    $\overline{B}$ in $X$
onto an open neighborhood of $B\subset V$ in $\E$ (where $V$ is  identified with the zero section $0_\E$) which is identity on an open neighborhood of  $\overline B\subset V$ such that the  endomorphism   on $\E$ induced  by the restriction
of the differential $d\tau$ to ${\overline B}$ is identity.  

In local coordinates, we can describe an admissible
map $\tau$ as follows: for every  point $x\in V\cap U,$  for every local
chart $y=(z,w)$  on a neighborhood $W$ of   $x$ in $U$  with  $V\cap W=\{z=0\}$, 
 we have
\begin{equation}\label{e:admissible-maps}
\tau (y) = \big( z + O(\|z\|^2), w+ O(\|z\|\big) ,
\end{equation}
and
\begin{equation}
d\tau  (y) = \big( dz + \widetilde O(\|z\|^2), dw + \widetilde O(\|z\| ) \big),
\end{equation}
as $z \to  0$ where for every positive integer $m,$  $\widetilde O(\|z \|^m )$ denotes the sum of $1$-forms with
$O(\|z\|^m )$-coefficients and a linear combination of $dz,$ $d\bar z$  with $O(\|z\|^{m-1} )$-coefficients.
\end{definition}
It  is worthy noting that in  \cite{DinhSibony18}   Dinh--Sibony  use the terminology {\it almost-admissible}  for those
maps satisfying  Definition \ref{D:admissible-maps}.
In general, $\tau$ is not holomorphic. When $U$ is a small enough local chart,
we can choose a holomorphic  admissible map by using suitable holomorphic coordinates
on $U .$   For the  global situation, the following  result gives a positive answer.
\begin{theorem}{\rm (\cite[Lemma 4.2]{DinhSibony18})}
For every  compact  subset $V_0\subset V,$ there always exists  an admissible map $\tau$ defined on a small enough tubular neighborhood  $U$ of $V_0$ in $X.$
\end{theorem}
  
In order  to  develop   a quantitative theory of  tangent  and density currents for  positive plurisubharmonic currents, the following notion, which is greatly inspired
by Dinh--Sibony \cite[Proposition 3.8]{DinhSibony18b}, is needed.  
\begin{definition}
 \label{D:Strongly-admissible-maps}\rm  Let $B$ be a relatively compact nonempty open subset of $V.$
   A {\it strongly admissible} map  along  $B$ is   a $\Cc^2$-smooth diffeomorphism $\tau$  from an open neighborhood $U$ of $\overline B$  in $X$
onto an open neighborhood of $V \cap U$ in $\E$  such that for every  point $x\in V\cap U,$  for every local
chart $y=(z,w)$  on a neighborhood $W$ of   $x$ in $U$  with  $V\cap W=\{z=0\}$,
 we have   
 \begin{eqnarray*}
 \tau_{j}(z,w)&= &  z_j+\sum_{p,q=1}^{k-l}  a_{pq}(w) z_pz_{q}  +O(\|z\|^3)\quad\text{for}\quad 1\leq j\leq  k-l,\\
 \tau_{j}(z,w)&= &  w_{j-(k-l)} +\sum_{p=1}^{k-l} b_p(w)z_p 
 +O(\|z\|^2)\quad\text{for}\quad k-l< j\leq  k.
 \end{eqnarray*}
 Here,  we  write $\tau(y)=(\tau_{1}(y),\ldots,\tau_{k-l}(y),\tau_{k-l+1}(y),  \ldots \tau_{k}(y))\in\C^k,$  and $a_{pq},$  $b_p,$ 
 are $\Cc^2$-smooth functions depending only on $w.$
 In other words, if we write $\tau(z,w)=(z',w')\in\C^{k-l}\times \C^l,$ then
 \begin{eqnarray*}
 z'&= &  z+  z A z^T  +O(\|z\|^3),\\
 w'&= &  w + Bz 
 +O(\|z\|^2),
 \end{eqnarray*}
 where  $A$ is a $(k-l)\times(k-l)$-matrix and $B$ is  a $l\times(k-l)$-matrix whose entries are $\Cc^2$-smooth functions in $w,$ $z^T$ is the transpose of $z,$
\end{definition}

Observe  that  a  strongly  admissible map is necessarily  admissible in the sense of Definition \ref{D:admissible-maps}.
 On the  other hand,  holomorphic  admissible maps  are always strongly   admissible.


\subsection{Tangent  currents and  known results} \label{SS:tangent-currents}
For  every  current $T$  of order $0$  on  the open set $U$  given by Definitions \ref{D:admissible-maps} or \ref{D:Strongly-admissible-maps}, let $\U:=\tau(U)$ and consider the
family of  currents  of order $0$ parameterized by $\lambda\in \C^* :$
\begin{equation}\label{e:T_lambda}  T_\lambda:= (A_\lambda)_*(\tau_*  T)\qquad\text{ on}\qquad \E|_{V \cap \U}.
 \end{equation}
\begin{definition}\label{D:tangent-currents} {\rm (\cite{DinhSibony18})}   \rm 
Let $B$ be an  open subset of $V.$
A tangent current $T_\infty$ of $T$ along $B$ is a  current on
$\pi^{-1}(B)\subset \E$ such that there are a sequence $(\lambda_n ) \subset \C^*$ converging to $\infty$  and a collection of 
admissible maps $\tau_\ell :\ U_\ell \to \U_\ell:=\tau_\ell(U_\ell)\subset \E $ for $\ell \in L,$ where $L$ is an index set, which  satisfy the following two properties:
\begin{enumerate}
\item[(i)]  $(U_\ell)_{\ell\in L}$ covers $B,$ that is,
$B \subset 
\bigcup\limits_{\ell\in L} U_\ell ;$
\item  [(ii)]  the masses  of the  currents  $T_{\lambda_n,\ell} $ 
are uniformly bounded  on compact subsets of  $\pi^{-1}(\U_\ell \cap B );$
\item[(iii)] the  following limit exists 
$$T_\infty := \lim\limits_{n\to\infty}  T_{\lambda_n,\ell}
\quad\text{on}\quad \pi^{-1}(\U_\ell \cap B )\quad\text{for all}\quad \ell \in  L.$$
Here, $ T_{\lambda_n,\ell}$ is  given  by \eqref{e:T_lambda} associated  to  the admissible map $\tau_\ell$ and  to $\lambda:=\lambda_n.$
\end{enumerate}
\end{definition}
We  record here  basic properties of tangent  currents.
\begin{theorem}\label{T:Dinh-Sibony} {\rm  (Dinh-Sibony  \cite{DinhSibony18})}
Assume  that $X$ is K\"ahler and $\supp(T) \cap V$ is compact. Then for every positive  closed  $(p,p)$-current $T$  on $X,$
the following  assertions hold:
\begin{enumerate}
 \item The masses  of the  currents  $T_{\lambda,\tau_\ell} $ with $\lambda\in\C^*$  are uniformly bounded  on compact subsets of  $\pi^{-1}(\U_\ell \cap 0_\E ).$
 In particular, the cluster limits  of   $T_{\lambda,\tau_\ell} $ as $\lambda\to\infty$ always exist.
 
 \item   $T_{\lambda,\tau_\ell} -T_{\lambda,\tau_{\ell'}} $ tends weakly to $0$  as $\lambda\to\infty$ on the  overlap   $\pi^{-1}(\U_\ell \cap \U_{\ell'}\cap 0_\E ).$
\item   If  the limit (iii) of Definition \ref{D:tangent-currents} holds for a sequence $(\lambda_n),$ then it still holds for this  sequence  when we replace $(\tau_\ell)_{\ell\in L}$  by another  collection
of admissible maps. In other words, the tangent limits are independent of the  choice of admissible maps.
\item Every tangent current $T_\infty$ is positive closed $V$-conic of bidegree $(p,p).$  
\end{enumerate}
\end{theorem}
It  is  interesting  to mention the following  improvement  where $X$ need not to be K\"ahler. 
\begin{theorem}\label{T:Vu} {\rm  (Vu \cite{Vu19a})} Assume that 
 there is a Hermitian  metric $\hat\omega$ on $X$ 
 for which  $\ddc\hat\omega^j=0$ on $V$ for $1\leq  j\leq k-p-1.$ Then the conclusion of   Theorem \ref{T:Dinh-Sibony} still holds for every positive  closed  $(p,p)$-current $T$  on $X$
 such that  $\supp(T) \cap V$ is compact.
\end{theorem}

\subsection{Function $\varphi$ and forms $\alpha$   and $\beta$ and  tubes}
 In this  subsection we  introduce  three  important objects  which will be used throughout the  article. Let $B\Subset V_0\Subset V$ be two open subsets of $V.$ 
  Denote by $\pi:\ \E\to V$ the canonical projection.
Consider a Hermitian metric  $h=\|\cdot\|$  on the  vector bundle  $\E_{\pi^{-1}(V_0)}$  and    let   $\varphi:\ \E_{\pi^{-1}(V_0)}\to \R^+$ be the function defined by  
\begin{equation}\label{e:varphi-spec}
 \varphi(y):=\|y\|^2\qquad \text{for}\qquad  y\in \pi^{-1}(V_0)\subset \E.
\end{equation}
Consider also the following  closed  $(1,1)$-forms on  $ \pi^{-1}(V_0)\subset \E $
\begin{equation}\label{e:alpha-beta-spec}
 \alpha:=\ddc\log\varphi\quad\text{and}\quad \beta:= \ddc\varphi.
\end{equation} 
So, for every $x\in V_0\subset  X$  the  metric $\| \cdot\|$  on the fiber $\E_x\simeq \C^{k-l}$ is  an Euclidean metric (in a suitable basis). In particular, we have 
\begin{equation}\label{e:varphi_bis-spec}  \varphi(\lambda y)=|\lambda|^2\varphi(y)\qquad\text{for}\qquad  y\in \pi^{-1}(V_0)\subset \E,\qquad\lambda\in\C. 
\end{equation}
For $r>0$ 
consider the  following {\it tube with base $B$ and radius $r$}
\begin{equation}
\label{e:tubular-nbh-0}
\Tube(B,r):=\left\lbrace y\in \E:\    \pi(y)\in B\quad\text{and}\quad  \|y\|<r  \right\rbrace.
\end{equation}
So   this is  a  natural generalization of Euclidean  tubes   considered by  Alessandrini--Bassanelli  in \eqref{e:Tube-AB}.
For for  all $0\leq s<r<\infty,$  define also the  corona tube
\begin{equation}\label{e:tubular-corona-0}\Tube(B,s,r):=\left\lbrace y\in \E:\   \pi(y)\in B\quad\text{and}\quad   s<\|y\|<r \right\rbrace.
\end{equation}
Since  $V_0\Subset V,$ there  is  a constant $c>0$ large enough such that 
$c\pi^*\omega+\beta$
is  positive on $\pi^{-1}(V_0).$  Moreover,   the latter  form defines
a K\"ahler metric there  if  $\omega$ is  K\"ahler on $V_0.$

\section{Lelong numbers  and first  collection of  main results}\label{S:Lelong}

\subsection{Standing  Hypothesis, global setting and generalized (main) Lelong  numbers}\label{SS:Global-setting}
We  keep the notation  introduced in   Sections  \ref{S:Intro} and \ref{S:preparatory_results}.
More specifically, we  assume the following 

\smallskip
\noindent {\bf Standing Hypothesis.} {\it 
Let $X$ be a complex  manifold of dimension $k.$
Let $V\subset  X$ be   a submanifold  of dimension $l$ and $B\subset V$ a  relatively compact piecewise  $\Cc^2$-smooth open subset.
Let  $V_0$ be  a relatively compact  open subset of $V$ such that  $B\Subset V_0.$ Let $\omega$ be a Hermitian form on $V.$
Let $\tau:\ U\to\tau(U)$ be  an admissible  map along $B$
 from an open neighborhood $U$ of $\overline B$ in $X.$ 
 Let $\bfr$ be small enough such that $\Tube(B,\bfr)\subset  \tau(U),$  see \eqref{e:tubular-nbh-0}.
 Fix $0\leq p\leq k.$
 Let $T$ be a real  current of degree $2p$ and  of  order $0$    on    $U.$
 }

\smallskip
\noindent {\bf Convention.} {\it Throughout the  first part of the  article, for the sake of simplicity   we will
omit the dependence of  the  mass indicators $\nu_j$ and $\kappa_j$ below  on  the  Hermitian form $\omega$ on $V$ and the  Hermitian metric $h$ on  $\E_{\pi^{-1}(V_0)}\subset \E.$ For example,  we  will write $\nu_j(T,B,r,\tau)$  (resp. $  \kappa_j(T,B,r,\tau)$)
instead  of  $\nu_j(T,B,r,\tau,\omega,h)$  (resp. $  \kappa_j(T,B,r,\tau,\omega,h)$).} 
 
 \smallskip

 Recall  from \eqref{e:m} that $\upm:= \min(l,k-p)$ and $
  \lowm:=\max(0,l-p).$  Let $\alpha$ and $\beta$  be the $(1,1)$-forms on $\pi^{-1}(V_0)\subset \E$   given by \eqref{e:alpha-beta-spec}.
 For $0\leq j\leq \upm$ and  $0<r\leq \bfr,$   consider 
\begin{equation}\label{e:Lelong-numbers}
 \nu_j(T,B,r,\tau):=  {1\over r^{2(k-p-j)}}\int_{\Tube(B,r)} (\tau_*T)\wedge \pi^*(\omega^j) \wedge \beta^{k-p-j}.  
\end{equation} 
When $j=\upm$  we  also denote  $\nu_\upm(T,B,r,\tau)$    by $\nu_\top(T,B,r,\tau).$

Let  $0\leq j\leq \upm.$ For  $0<s<r\leq \bfr,$   consider 
\begin{equation}\label{e:Lelong-corona-numbers}
 \kappa_j(T,B,s,r,\tau):=   \int_{\Tube(B,s,r)} (\tau_*T)\wedge \pi^*(\omega^j) \wedge \alpha^{k-p-j}. 
\end{equation} 
Let $0<r\leq\bfr.$   Consider
\begin{equation}\label{e:Lelong-log-bullet-numbers}
 \kappa^\bullet_j(T,B,r,\tau):= \limsup\limits_{s\to0+}   \kappa_j(T,B,s,r,\tau). 
\end{equation} 
We  also consider
\begin{equation}\label{e:Lelong-log-numbers}
 \kappa_j(T,B,r,\tau):= \int_{\Tube(B,r)} (\tau_*T)\wedge \pi^*(\omega^j) \wedge \alpha^{k-p-j},  
\end{equation}
provided that the RHS side makes sense according  to the  following definitions.

 

\begin{definition}\label{D:Lelong-log-numbers(1)}\rm
We say that  \eqref{e:Lelong-log-numbers} holds in the spirit of \eqref{e:Lelong-number-point-bisbis(1)}
if $T=T^+-T^-$  in an open neighborhood of $\overline B$ in $X$ and $T^\pm\in\SH^{m,m'}_p(B)$ (resp.  $T^\pm\in\PH^{m,m'}_p(B),$  resp.
$T^\pm\in\CL^{m,m'}_p(B)$ for some  suitable integers $0\leq m'\leq m$) with the corresponding sequences of approximating forms   $(T^\pm_n)_{n=1}^\infty,$
and  for any  such  forms  $(T^\pm_n),$  the two  limits  on the  following  RHS exist and  are finite   
\begin{equation}\label{e:Lelong-log-numbers(2)}
\kappa_j(T,B,r,\tau):= \lim\limits_{n\to\infty} \kappa_j(T^+_n,B,r,\tau)-\lim\limits_{n\to\infty}\kappa_j(T^-_n,B,r,\tau),
\end{equation}
and the  value on the RHS is  independent of the choice of $(T^\pm_n)_{n=1}^\infty.$
\end{definition}

\begin{definition}\label{D:Lelong-log-numbers(2)}\rm
We say that  \eqref{e:Lelong-log-numbers} holds in the spirit of \eqref{e:Lelong-number-point-bisbis(2)}
if  the    limit  on the  following  RHS exists and  is finite   
\begin{equation}\label{e:Lelong-log-numbers(1)}
\kappa_j(T,B,r,\tau):= \lim\limits_{\epsilon\to 0+}  \int_{\Tube(B,r)} (\tau_*T)\wedge \pi^*(\omega^j) \wedge \alpha^{k-p-j}_\epsilon.
\end{equation}
Here, the  smooth  form $\alpha_\epsilon$ is given  by \eqref{e:alpha-beta-eps} below.
\end{definition}

When $j=\upm$  we  also denote  $\kappa_\upm(T,B,s,r,\tau)$ (resp.    $\kappa_\upm(T,B,r,\tau)$)   by $\nu_\top(T,B,s,r,\tau)$
 (resp.    $\kappa_\top(T,B,r,\tau)$).

 \subsection{Intermediate  average means}\label{SS:Intermediate-Lelong-means}

 For $0\leq j\leq \upm,$ $0\leq q\leq k-l$  and  $0<r\leq \bfr,$   consider 
\begin{equation}\label{e:inter-Lelong-numbers}
 \nu_{j,q}(T,B,r,\tau):=  {1\over r^{2q}}\int_{\Tube(B,r)} (\tau_*T)\wedge \pi^*(\omega^j) \wedge \beta^{k-p-j}.  
\end{equation}

Let  $0\leq j\leq \upm$ and $0\leq q\leq \min(k-l,k-p-j).$ For  $0<s<r\leq \bfr,$   consider 
\begin{equation}\label{e:Lelong-inter-corona-numbers}
 \kappa_{j,q}(T,B,s,r,\tau):=   \int_{\Tube(B,s,r)} (\tau_*T)\wedge \pi^*(\omega^j)\wedge \beta^{k-p-j-q} \wedge \alpha^{q}. 
\end{equation} 
Let $0<r\leq\bfr.$   Consider
\begin{equation}\label{e:Lelong-inter-log-bullet-numbers}
 \kappa^\bullet_{j,q}(T,B,r,\tau):= \limsup\limits_{s\to0+}   \kappa_{j,q}(T,B,s,r,\tau). 
\end{equation} 
We  also consider
\begin{equation}\label{e:Lelong-inter-log-numbers}
 \kappa_{j,q}(T,B,r,\tau):= \int_{\Tube(B,r)} (\tau_*T)\wedge \pi^*(\omega^j)\wedge \beta^{k-p-j-q} \wedge \alpha^{q},  
\end{equation}
provided that the RHS side makes sense according  to  Definitions \ref{D:Lelong-log-numbers(1)} and \ref{D:Lelong-log-numbers(2)}.

 \begin{remark}\label{R:Lelong-inter}\rm
  For $\lowm\leq j\leq\upm,$  we have  
  \begin{eqnarray*}\nu_j(T,B,r,\tau)&=&\nu_{j,k-p-j}(T,B,r,\tau)
\quad\text{and}\quad  \nu_{j,q}(T,B,r,\tau)=r^{2(k-p-j-q)}\nu_{j,k-p-j}(T,B,r,\tau),\\
 \kappa_j(T,B,r,\tau)&=&\kappa_{j,k-p-j}(T,B,r,\tau).
\end{eqnarray*}
For $0\leq j\leq\lowm,$   we have
\begin{equation*} \nu_j(T,B,r,\tau)=\nu_{j,k-l}(T,B,r,\tau)\qquad\text{ and}\qquad 
  \kappa_j(T,B,r,\tau)=\kappa_{j,k-l}(T,B,r,\tau).\end{equation*}

 \end{remark}

\subsection{First collection of main  results}\label{SS:First-collection}
 
 The  main purpose  of this section is to state    the following seven theorems. The first four deal with positive closed  currents,  whereas
 the last three are  devoted to positive pluriharmonic  currents and  positive  plurisubharmonic  currents.

 The  first theorem only deals with the top degree $j=\upm$ and with  a holomorphic admissible map $\tau,$
 but  it does not require  any  condition on the support of $T$ nor on the Hermitian metric $\omega.$

\begin{theorem}\label{T:top-Lelong-closed}
  We  keep the     Standing Hypothesis and  assume that $p>0.$ 
Suppose in addition that     $T=T^+-T^-$ on an open neighborhood  of $\overline B$ in $X$
with $T^\pm$ in the class $\CL_p^{1}( B).$
  Then, for every holomorphic admissible map $\tau,$ the following  assertions  hold: 
  \begin{enumerate}
  \item  For $0<r_1<r_2\leq \bfr,$
  $$ \nu_\top(T,B,r_2,\tau)-\nu_\top(T,B,r_1,\tau)=\kappa_\top(T,B,r_1,r_2,\tau),
  $$
  and all three numbers are nonnegative real numbers. In particular, the function $r\mapsto \nu_\top(T,B,r,\tau),$ defined for  $r\in (0,\bfr)$  with non-negative values,  is  increasing.

  \item  The following limit  exists
  $$    \nu_\top(T,B,\tau):=\lim\limits_{r\to 0+}\nu_\top(T,B,r,\tau )        ,                   $$
  and
$\nu_\top(T,B,\tau)$  is  a  nonnegative real number.

 \item       
  $
  \lim\limits_{r\to 0+}\kappa^\bullet_\top(T,B,r,\tau )=0.                                
  $
 \item    The following limit holds in the sense of Definitions \ref {D:Lelong-log-numbers(1)}  and  \ref {D:Lelong-log-numbers(2)}:
  $$    \lim\limits_{r\to 0+}\kappa_\top(T,B,r,\tau )=\nu_\top(T,B,\tau).                                $$
  
\end{enumerate}
 \end{theorem}

 The  second theorem  deals with all degrees $\lowm\leq j\leq \upm$  (eventually  with  all degrees $0\leq j\leq \upm$),   but with  a holomorphic admissible map $\tau,$  and it requires  a control  of approximation of  $T$ on the boundary and a condition on the Hermitian form $\omega.$
 
 \begin{theorem}\label{T:Lelong-closed-all-degrees}
   We  keep the    Standing Hypothesis. Suppose that $ \ddc\omega^j=0$ on $V_0$ for all  $1\leq j\leq   \upm-1.$
Suppose in addition that    the current   $T$  is  positive closed  and $T=T^+-T^-$ on an open neighborhood  of $\overline B$ in $X$
with $T^\pm$ in the class $\CL_p^{1,1}( B).$
  Then, for every holomorphic admissible map $\tau,$ the following  assertions hold for  $\lowm\leq j\leq \upm:$ 
  \begin{enumerate}
  
   \item  For $0<r_1<r_2\leq\bfr,$
  $$ \nu_j(T,B,r_2,\tau)-\nu_j(T,B,r_1,\tau)=\kappa_j(T,B,r_1,r_2,\tau)+O(r_2),
  $$
  where $|O(r_2)|\leq cr_2$  for a  constant $c>0$ which depends only on $T,X,V, B,\omega$ but which does not depend on $r_2.$
   \item  The limit 
  $    \nu_j(T,B,\tau):=\lim\limits_{r\to 0+}\nu_j(T,B,r,\tau )$ exists and $ \nu_j(T,B)\in\R.$
  \item   
  $
  \lim\limits_{r\to 0+}\kappa^\bullet_j(T,B,r,\tau )=0      .
  $
  \item For all $\lowm\leq j\leq \upm$ with $j>l-p,$ the following limit holds in the sense of Definition \ref {D:Lelong-log-numbers(1)} and    Definition \ref {D:Lelong-log-numbers(2)}:
  $    \lim\limits_{r\to 0+}\kappa_j(T,B,r,\tau )=\nu_j(T,B,\tau).   $
  
  \item  Suppose in addition that  $\supp(T^\pm_n)\cap V\subset  B$ for $n\geq 1.$ Then the above assertions (1)--(5) also hold
  for all $0\leq j\leq \upm.$  Moreover, the following  stronger  version of assertion (1) also holds:
  For  $0\leq j\leq \upm$  and   $0<r_1<r_2\leq\bfr,$ 
  $$ \nu_j(T,B,r_2,\tau)-\nu_j(T,B,r_1,\tau)=\kappa_j(T,B,r_1,r_2,\tau).
  $$
  \item   If  moreover $\omega$ is  K\"ahler, then  all the above assertions (1)--(5) still hold if 
    $T=T^+-T^-,$ where  $T^\pm$ only belong to  the class $\CL_p^{1,0}(B).$
\end{enumerate}
 \end{theorem}

 The  third   theorem  deals with all  degree $\lowm\leq j\leq \upm,$   with  a non-holomorphic admissible map $\tau,$  but it requires  a control  of approximation of  $T$ on the boundary.
  
 \begin{theorem}\label{T:Lelong-closed}
   We  keep the     Standing Hypothesis. Suppose that $ \ddc\omega^j=0$ on $B$ for all  $1\leq j\leq   \upm-1.$
Suppose in addition that    the current   $T$  is  positive closed  and $T=T^+-T^-$ on an open neighborhood  of $\overline B$ in $X$
with $T^\pm$ in the class $\CL_p^{2,2}(B).$
  Then, for every  strongly admissible map $\tau,$ the following  assertions hold for  $\lowm\leq j\leq \upm:$
  \begin{enumerate}
  
   \item  For $0<r_1<r_2\leq\bfr,$
  $$ \nu_j(T,B,r_2,\tau)-\nu_j(T,B,r_1,\tau)=\kappa_j(T,B,r_1,r_2,\tau)+O(r_2).
  $$
  \item   The  limit  $    \nu_j(T,B,\tau):=\lim\limits_{r\to 0+}\nu_j(T,B,r,\tau )$ exists and $  \nu_j(T,B,\tau)\in\R .$
   \item  
  $
  \lim\limits_{r\to 0+}\kappa^\bullet_j(T,B,r,\tau )=0 .
  $
  \item For all $\lowm\leq j\leq \upm$ with $j>l-p,$ the following limit holds in the sense of   Definition \ref {D:Lelong-log-numbers(1)} and Definition \ref {D:Lelong-log-numbers(2)}:
$    \lim\limits_{r\to 0+}\kappa_j(T,B,r,\tau )=\nu_j(T,B,\tau).                                $
  
  \item  $\nu_j(T,B,\tau)$  is
  independent of  the choice  of   $\tau.$

  \item $\nu_\top(T,B,\tau)$ is a nonnegative real number.
\end{enumerate}
 \end{theorem}

 The  fourth   theorem  deals with all  degree $\lowm\leq j\leq \upm,$   with  a non-holomorphic admissible map $\tau,$  but it requires the K\"ahlerity of the  metric $\omega$ and  a control  of approximation of  $T$ on the boundary.
  
 \begin{theorem}\label{T:Lelong-closed-Kaehler}
   We  keep the     Standing Hypothesis. Suppose that $ \omega$ is  K\"ahler.
Suppose in addition that    the current   $T$  is  positive closed  and $T=T^+-T^-$ on an open neighborhood  of $\overline B$ in $X$
with $T^\pm$ in the class $\CL_p^{1,1}(B).$
  Then, for every  strongly admissible map $\tau,$ the following  assertions hold for  $\lowm\leq j\leq \upm:$
  \begin{enumerate}
  
    \item  For $0<r_1<r_2\leq\bfr,$
  $$ \nu_j(T,B,r_2,\tau)-\nu_j(T,B,r_1,\tau)=\kappa_j(T,B,r_1,r_2,\tau)+O(r_2).
  $$
  \item   The  limit  $    \nu_j(T,B,\tau):=\lim\limits_{r\to 0+}\nu_j(T,B,r,\tau )$ exists and $  \nu_j(T,B,\tau)\in\R .$
   \item  
  $
  \lim\limits_{r\to 0+}\kappa^\bullet_j(T,B,r,\tau )=0 .
  $
  \item For all $\lowm\leq j\leq \upm$ with $j>l-p,$ the following limit holds in the sense of   Definition \ref {D:Lelong-log-numbers(1)}  and   Definition \ref {D:Lelong-log-numbers(2)}:
$    \lim\limits_{r\to 0+}\kappa_j(T,B,r,\tau )=\nu_j(T,B,\tau).                                $
  
  \item  $\nu_j(T,B,\tau)$  is
  independent of  the choice  of   $\tau.$

  \item $\nu_\top(T,B,\tau)$ is a nonnegative real number.
\end{enumerate}
 \end{theorem}

 Now  we come  to the  three theorems on positive pluriharmonic  currents and positive plurisubharmonic currents. 
 The  first theorem only deals with the top degree $j=\upm$ and with  a holomorphic admissible map $\tau,$
 but  it does not require  any  condition on the support of $T.$

\begin{theorem}\label{T:top-Lelong-psh}
  We  keep the     Standing Hypothesis. Suppose that $\omega$ is  K\"ahler.
  Suppose in addition that $T=T^+-T^-$ on an open neighborhood  of $\overline B$ in $X$
with $T^\pm$ in the class $\SH_p^{2}(\overline B).$
  Then, for every holomorphic admissible map $\tau,$ the following  assertions hold.  
  \begin{enumerate}
   \item  For $0<r_1<r_2\leq \bfr,$
  $$ \nu_\top(T,B,r_2,\tau)-\nu_\top(T,B,r_1,\tau)\geq\kappa_\top(T,B,r_1,r_2,\tau),
  $$
  and all three numbers are nonnegative real numbers. In particular, the function $r\mapsto \nu_\top(T,B,r,\tau),$ defined for  $r\in (0,\bfr)$  with non-negative values,  is  increasing.
  
  \item  The following limit  exists
  $$    \nu_\top(T,B,\tau):=\lim\limits_{r\to 0+}\nu_\top(T,B,r,\tau )        ,                   $$
  and
$\nu_\top(T,B,\tau)$  is  a  nonnegative real number.

  \item       
  $
  \lim\limits_{r\to 0+}\kappa^\bullet_\top(T,B,r,\tau )=0.                                
  $
  
  \item   As  a positive closed current, $\ddc T$ satisfies 
  $\nu_\top(\ddc T,B,\tau)=0.$
 \item    The following limit holds in the sense  of Definition  \ref {D:Lelong-log-numbers(1)}  and
 \ref {D:Lelong-log-numbers(2)}:
  $$    \lim\limits_{r\to 0+}\kappa_\top(T,B,r,\tau )=\nu_\top(T,B,\tau).                                $$
  \item  $\nu_\top(T,B,\tau)$  is
  independent of  the choice  of a holomorphic admissible map  $\tau.$ 
  
\end{enumerate}
 \end{theorem}

 The  second theorem  deals with all degrees $\lowm\leq j\leq \upm,$   but with  a 
 holomorphic admissible map $\tau,$  and it requires a positivity of the basic  forms $\alpha,$ $\beta$ as well
 as   a control  of approximation of  $T$ on the boundary.
 
 \begin{theorem}\label{T:Lelong-psh-all-degrees}
We  keep the     Standing Hypothesis. Suppose that $\omega$ is  K\"ahler and the forms $\alpha,$ $\beta$   are positive.
   Suppose in addition that  the current   $T$  is  positive plurisubharmonic  
  and $T=T^+-T^-$ on an open neighborhood  of $\overline B$ in $X$
with $T^\pm$ in the class $\SH_p^{2,2}(B).$
  Then, for every holomorphic admissible map $\tau,$ the following  assertions hold for  $\lowm\leq j\leq \upm:$ 
  \begin{enumerate}
  \item
  For $0<r_1<r_2\leq \bfr,$
  $$ \nu_j(T,B,r_2,\tau)-\nu_j(T,B,r_1,\tau)\geq\kappa_j(T,B,r_1,r_2,\tau)+O(r_2),
  $$
  and all three numbers $ \nu_j(T,B,r_2,\tau),$ $\nu_j(T,B,r_1,\tau),$ $\kappa_j(T,B,r_1,r_2,\tau)$
  are nonnegative real numbers.

  \item  The limit 
  $    \nu_j(T,B,\tau):=\lim\limits_{r\to 0+}\nu_j(T,B,r,\tau )$ exists and $ \nu_j(T,B,\tau)\in\R.$
  
  \item   
  $
  \lim\limits_{r\to 0+}\kappa^\bullet_j(T,B,r,\tau )=0      .
  $

  \item   As  a positive closed current, $\ddc T$ satisfies 
  $\nu_\top(\ddc T,B,\tau)=0.$
  \item  For all $\lowm\leq j\leq \upm$ with $j>l-p,$ the following limit holds in the sense   of   Definition \ref {D:Lelong-log-numbers(1)}  and  Definition \ref {D:Lelong-log-numbers(2)}:
  $    \lim\limits_{r\to 0+}\kappa_j(T,B,r,\tau )=\nu_j(T,B,\tau).   $
  \item  $\nu_j(T,B,\tau)$  is
  independent of  the choice  of   $\tau.$
  
  \item  If $\supp(T^\pm_n)\cap V\subset  B$ for $n\geq 1,$ then the above assertions (1)--(4) also hold
  for all $0\leq j\leq \upm.$
\end{enumerate}
 \end{theorem}

 The  third   theorem  deals with all  degree $\lowm\leq j\leq \upm,$   with  a non-holomorphic admissible map $\tau,$  but it requires  a control  of approximation of  $T$ on the boundary.

 \begin{theorem}\label{T:Lelong-psh}
   We  keep the     Standing Hypothesis. Suppose that $\omega$ is  K\"ahler.
     Suppose in addition that  the current   $T$  is  positive plurisubharmonic  
  and $T=T^+-T^-$ on an open neighborhood  of $\overline B$ in $X$
with $T^\pm$ in the class $\SH_p^{3,3}(B).$
 Then, for every  strongly admissible map $\tau,$ the following  assertions hold for  $\lowm\leq j\leq \upm:$
  \begin{enumerate}
  
  \item   The  limit  $    \nu_j(T,B,\tau):=\lim\limits_{r\to 0+}\nu_j(T,B,r,\tau )$ exists and $  \nu_j(T,B,\tau)\in\R .$
  
   \item  
  $
  \lim\limits_{r\to 0+}\kappa^\bullet_j(T,B,r,\tau )=0 .
  $
  \item For all $\lowm\leq j\leq \upm$ with $j>l-p,$ the following limit holds 
in the sense   of   Definition \ref {D:Lelong-log-numbers(1)}   and    Definition 
\ref {D:Lelong-log-numbers(2)}: $    \lim\limits_{r\to 0+}\kappa_j(T,B,r,\tau )=\nu_j(T,B,\tau).                                $
  
   \item $\nu_\top(T,B,\tau)$ is a nonnegative real number.
  \item  $\nu_j(T,B,\tau)$  is
  independent of  the choice  of   $\tau.$ 
  \item  If instead of the above assumption on $ T,$ we assume that   T is a positive pluriharmonic $(p,p)$-current on a neighborhood of $\overline B$ in X such that
$T=T^+-  T^-$ for some $T^\pm\in \PH^{2,2}_p(B),$
 then all the above assertions still hold.
\end{enumerate}
 \end{theorem}

  \begin{definition}\label{D:Lelong-numbers-bis}
  \rm  For  $0\leq j\leq  \upm,$ the real number 
 $\nu_j(T,B,\tau)$ (if it is  well-defined) is  called  the $j$-th Lelong number of $T$ along  $B.$
 Since  by the above  theorems, $\nu_j(T,B,\tau)$ is  independent of the choice of a strongly admissible map $\tau,$
 we  will  denote  it simply by  $\nu_j(T,B).$
 
The   set  of all  well-defined  Lelong numbers $\{\nu_j(T,B):$\ $0\leq j\leq \upm\}$    are called {\it the Lelong numbers of $T$ along
$B.$} 
The nonnegative  number $\nu_\upm(T,B)$ is  called  {\it  the top Lelong numbers of $T$ along $B,$}
it is  also denoted by  $\nu_\top(T,B).$
  \end{definition}

The functions $\nu_j$  and $\kappa_j$ introduced  in \eqref{e:Lelong-numbers} and \eqref{e:Lelong-corona-numbers} enjoy the following   simple  scaling property.

\begin{proposition}\label{P:scale}
 For  every $0<s<r\leq\bfr$ and $\lambda\geq 1$ and  $0\leq j\leq \upm,$  we have that 
\begin{eqnarray*}
\nu_j(T,B,{r\over  \lambda},\tau)=\nu_j(   (A_\lambda)_*(\tau_*T),B, r,\id),\\
\kappa_j(T,B,{s\over \lambda},{r\over \lambda},\tau)=\kappa_j( (A_\lambda)_*(\tau_*T), B,r,\id).
\end{eqnarray*}
\end{proposition}
\proof
 By \eqref{e:varphi-spec}
 and \eqref{e:alpha-beta-spec},  we have that
 $$(A_\lambda)^*\beta^{k-p-j}=|\lambda|^{2(k-p-j)} \beta\qquad\text{and}\qquad  (A_\lambda)^* \big(\pi^*(\omega^j)\big)= \pi^*(\omega^j).$$
 Using this, we get that 
 \begin{multline*}
   {\lambda^{2(k-p-j)}\over r^{2(k-p-j)}}\int_{\Tube(B,{r\over \lambda})} (\tau_*T)\wedge \pi^*(\omega^j)\wedge \beta^{k-p-j}\\=
   {1\over  r^{2(k-p-j)}}\int_{\Tube(B,r)} (A_\lambda)_*(\tau_*T)\wedge
   \pi^*(\omega^j)\wedge \beta^{k-p-j}.
 \end{multline*}
 Hence, 
the first identity follows. The  second one  can be proved  in the  same way. 
\endproof

When  $X$  is  a   holomorphic vector bundle and $V$ is  the  base space, we see that $\E=X$ and 
we can choose $\tau=\id,$ and the above  proposition gives the following 
useful identity

\begin{corollary} \label{C:scale} Let $T$   be as  above. 
Then  for  every $0<s<r\leq\bfr$ and $\lambda\geq 1$ and  $0\leq j\leq \upm,$  we have that 
\begin{eqnarray*}
\nu_j(T,B,{r\over  \lambda},\id)=\nu_j(   (A_\lambda)_*T,B, r,\id),\\
\kappa_j(T,B,{s\over \lambda},{r\over \lambda},\id)=\kappa_j( (A_\lambda)_*T, B,r,\id).
\end{eqnarray*}
\end{corollary}

\subsection{Local setting}\label{SS:Local-setting}

We  explain some of the above  main results   in  the local setting   of Alessandrini--Bassanelli \cite{AlessandriniBassanelli96}.

Let $T$ be a $(p,p)$-current of order $0$  defined on  an open neighborhood $U$ of $0$ in $\C^k.$
We use the coordinates $(z,w)\in\C^{k-l}\times \C^l.$ 
We may assume that  $U$ has the form $U=U'\times U'',$ where $U'$ (resp. $U'')$ are open neighborhood of $0'$ in $\C^{k-l}$ of  ($0''$ in $\C^l$).
Let $V=\{z=0\}=U''$ and let $B=B_w\Subset U''$ be a domain with  piecewise  $\Cc^2$-smooth boundary    and $\bfr>0$  such that  $\{\|z\|<\bfr\}\times B\Subset U.$ 
Consider  the trivial  vector bundle $\pi:\ \E \to  U''$ with  $\E\simeq  \C^{k-l}\times U''.$ For $\lambda\in\C^*,$  let $a_\lambda:\ \E\to \E$ be the multiplication by  $\lambda$
on fibers, that is, 
$a_\lambda(z,w):=(\lambda z,w)$ for $(z,w)\in \E.$
The admissible map $\tau$ in this setting is simply the identity $\id.$

Consider the positive closed $(1,1)$-forms
\begin{equation}\label{e:alpha-beta-Upsilon-local}\beta=\omega_z:=\ddc \|z\|^2\quad\text{and}\quad  \omega=\omega_w:=\ddc\|w\|^2\quad\text{and}\quad \alpha=\Upsilon_z:=\ddc \log\|z\|^2.
\end{equation}
Let $\lowm\leq j\leq \upm.$ For  $0<r<\bfr,$ consider the quantity 
\begin{equation}\label{e:Lelong-numbers-local}
 \nu_j(T,B,r):=  {1\over r^{2(k-p-j)}}\int_{\|z\|<r,\ w\in B} T\wedge \omega_w^j \wedge \omega_z^{k-p-j}.
\end{equation}
 For  $0<s<r\leq \bfr,$   consider 
\begin{equation}\label{e:Lelong-corona-numbers-local}
 \kappa_j(T,B,s,r):=   \int_{s<\|z\|<r,\ w\in B}  T\wedge \omega_w^j \wedge \Upsilon_z^{k-p-j}. 
\end{equation} 
Let $0<r\leq \bfr.$   Consider
\begin{equation}\label{e:Lelong-log-bullet-numbers-local}
 \kappa^\bullet_j(T,B,r):= \limsup\limits_{s\to0+}   \kappa_j(T,B,s,r). 
\end{equation} 
We  also consider
\begin{equation}\label{e:Lelong-log-numbers-local}
 \kappa_j(T,B,r):=  \int_{\|z\|<r,\ w\in B}  T\wedge \omega_w^j \wedge \Upsilon_z^{k-p-j},  
\end{equation}
provided that the right hand  side makes sense in the sense of
of \eqref{e:Lelong-number-point-bisbis(1)} or
\eqref{e:Lelong-number-point-bisbis(2)}. 
 
 As  an immediate consequence of  Theorems \ref{T:top-Lelong-psh} and \ref{T:Lelong-psh-all-degrees}, we get the following result which is  in part more general than  Theorem \ref{T:AB-1}. The latter theorem  only  gives the top Lelong number.

\begin{corollary}
 We  keep the     Standing Hypothesis and   the  above  assumption in the local setting.
 Suppose in addition that  the current   $T$  is  positive plurisubharmonic  
  and $T=T^+-T^-$ on an open neighborhood $\Omega$ of $\overline B$ in $U$
with either  of  the following      conditions: 
\begin{itemize}
\item[(a)] If $j=\upm$  then  we require  that the currents $T^\pm$
  are approximable on   $U$ by  some $\Cc^2$-smooth positive  plurisubharmonic  forms $(T^\pm_n)_{n=1}^\infty$
in the  following sense: 
\begin{itemize}
\item [(a-i)] $T^\pm_n$ converge weakly to the current $T^\pm$ on $U ;$
\item [(a-ii)] the masses $\|T^\pm_n\|$ on $U$ are uniformly bounded.
\end{itemize}
\item[(b)]  If $0 \leq j<\upm$  then we require that 
$T^\pm$ belong to  the class $\SH_p^{2,1}(B).$
\end{itemize}
  Then, for every holomorphic admissible map $\tau,$ the following  assertions hold for  $\lowm\leq j\leq \upm:$ 
  \begin{enumerate}
  \item   
  $
  \lim\limits_{r\to 0+}\kappa^\bullet_j(T,B,r)=0      .
  $
  \item  The limit 
  $    \nu_j(T,B):=\lim\limits_{r\to 0+}\nu_j(T,B,r )$ exists and $ \nu_j(T,B)\in\R.$
  \item  If $j>l-p,$ then
  $    \lim\limits_{r\to 0+}\kappa_j(T,B,r )=\nu_j(T,B).   $
  
 \item  $\nu_\top(T,B)$  is  a  nonnegative real number.
  \item  If $\supp(T^\pm_n)\cap V\subset  B$ for $n\geq 1,$ then the above assertions (1)--(3) also hold
  for all $0\leq j\leq \upm.$
  \end{enumerate}
  
\end{corollary}



\section{Lelong-Jensen  formulas   for  vector bundles} \label{S:Lelong-Jensen}

In this  section we  introduce the main tool of this  article: {\it  Lelong-Jensen  formulas   for  vector bundles.}   These formulas  play
a key role  throughout   this  work.

\subsection{Tubes, horizontal and  vertical boundaries}
Let $V$ be  a complex manifold of dimension $l.$ Let $\E$ be a holomorphic  bundle of rank $k-l$  over $V.$
So $\E$ is a complex manifold of dimension $k.$  Denote by $\pi:\ \E\to V$ the canonical projection.
Let  $B$ be a relatively compact  open set of $V$ with piecewice $\Cc^2$-smooth boundary.
Let $\U$ be an open neighborhood of $\overline B$ in $\E.$
Let $\bfr\in \R^+_*\cup \{\infty\}$ and  $0\leq r_0<\bfr.$  Let $\varphi:\ \U\to  [0,\infty)$ be a  $\Cc^2$-smooth  function   such that
\begin{itemize}
 \item $\varphi(y)=r^2_0$ for  $y\in\U\cap V$ and $\varphi(y)>r^2_0$ for $y\in \U\setminus V;$
 \item for every $r\in(r_0, \bfr],$  the set $\{y\in\U:\ \varphi(y)=r^2\}$ is a connected nonsingular real hypersurface of $\U$
 which intersects  the real hypersurface $\pi^{-1}(\partial B)\subset \E$ transversally.
\end{itemize}
Consider also the following  closed  $(1,1)$-forms on  $\U$
\begin{equation}\label{e:alpha-beta}
 \alpha:=\ddc\log\varphi\quad\text{and}\quad \beta:= \ddc\varphi.
\end{equation} 
Let $r>0$ and   $B\Subset V$  an open set.
Consider the  following {\it tube with base $B$ and radius $r$}
\begin{equation}
\label{e:tubular-nbh}
\Tube(B,r):=\left\lbrace y\in \E:\   \varphi(y)<r^2  \right\rbrace.
\end{equation}
For   all $r_0\leq r<s\leq \bfr,$  define
\begin{equation}\label{e:tubular-corona}\Tube(B,r,s):=\left\lbrace y\in \E:\   \pi(y)\in B\quad\text{and}\quad   r^2<\varphi(y)<s^2  \right\rbrace.
\end{equation}
Note that the boundary   $\partial\Tube(B,r)$ can be decomposed as the disjoint union  of the {\it  vertical boundary}   $\partial_\ver\Tube(B,r)$ and
the  {\it horizontal  boundary}    $\partial_\hor\Tube(B,r)$, where
\begin{eqnarray*}
 \partial_\ver\Tube(B,r)&:=&  \left\lbrace y\in \E:\   \pi(y)\in \partial B\quad\text{and}\quad  \varphi(y)\leq r^2  \right\rbrace ,\\
 \partial_\hor\Tube(B,r)&:=&  \left\lbrace y\in \E:\   \pi(y)\in B\quad\text{and}\quad  \varphi(y)=r^2  \right\rbrace .
\end{eqnarray*}
Under the above assumption on $\varphi,$ we see that $\Tube(B,r)$ is a  manifold  with piecewise $\Cc^2$-smooth  boundary
for  every $r\in[r_0,\bfr].$  When  $\partial B=\varnothing,$  we have $\partial_\ver\Tube(B,r)=\varnothing.$

\subsection{Abstract  formulas}
\begin{notation}
 \label{N:principal}\rm 
  Let $S$ be  a current of bidegree $2p$  defined on $\Tube(B,\bfr)\subset \E.$
  We  denote by $S^\sharp$  or equivalently $(S)^\sharp$ its component of bidegree $(p,p).$
\end{notation}
We are in the position to state  and prove the first  Lelong-Jensen  formulas   for  vector bundles.
\begin{theorem}\label{T:Lelong-Jensen}
Let $r\in(r_0,\bfr]$ and   $B\Subset V$  a relatively compact  open set with piecewice $\Cc^2$-smooth boundary.
Let $S$ be a real current of dimension $2q$ on  a neighborhood of $\overline\Tube(B,r)$ such that $S$ and $\ddc S$ are of order  $0$
and  that  $S$ is of class $\Cc^1$ near $\partial_\ver \Tube( B, r).$  
Suppose that there is a sequence of $\Cc^2$-smooth forms of dimension $2q$ $(S_n)_{n=1}^\infty$ defined on  a neighborhood of $\overline\Tube(B,r)$ 
such  that
\begin{itemize}
 \item[(i)] 
$S_n$ converge to $S$  in the sense of quasi-positive currents   on a neighborhood of  $\overline\Tube(B,r)$ as $n$ tends to infinity (see Definition \ref{D:quasi-positivity});
\item [(ii)]  $\ddc S_n$ converge to $\ddc S$  in the sense of quasi-positive currents   on a neighborhood of  $\overline\Tube(B,r)$ as $n$ tends to infinity;
\item[(iii)] there is an  open neighborhood of $\partial_\ver\Tube(B,r)$  on which  $S_n$  converge to $S$  in  $\Cc^1$-norm.
\end{itemize}
Then the  following two assertions hold: 
\begin{enumerate} \item 
The following  four sub-assertions hold:
\begin{itemize} \item[(1-i)] For  all $r_1,r_2\in  (r_0,r]$ with  $r_1<r_2$  except for  a  countable set of   values,  we have that  
\begin{equation}\label{e:Lelong-Jensen}\begin{split}
 {1\over  r_2^{2q}} \int_{\Tube(B,r_2)} S\wedge \beta^q- {1\over  r_1^{2q}} \int_{\Tube(B,r_1)} S\wedge \beta^q
 = \lim\limits_{n\to\infty}\Vc(S_n,r_1,r_2)+   \int_{\Tube(B,r_1,r_2)} S\wedge \alpha^q\\
 +  \int_{r_1}^{r_2} \big( {1\over t^{2q}}-{1\over r_2^{2q}}  \big)2tdt\int_{\Tube(B,t)} \ddc S\wedge \beta^{q-1} 
 +  \big( {1\over r_1^{2q}}-{1\over r_2^{2q}}  \big) \int_{r_0}^{r_1}2tdt\int_{\Tube(B,t)} \ddc S\wedge \beta^{q-1}.
 \end{split}
\end{equation} 
Here the vertical boundary term  $\Vc(S,r_1,r_2)$  for a $\Cc^1$-smooth form $S$  is  given by the following formula,  where $S^\sharp$ denotes, 
according to Notation \ref{N:principal}, the component of bidimension $(q,q)$ of the current $S:$
\begin{equation}\label{e:vertical-boundary-term}
\begin{split}
 \Vc(S,r_1,r_2)&:=-\int_{r_1}^{r_2} \big( {1\over t^{2q}}-{1\over r_2^{2q}}  \big)2tdt\int_{\partial_\ver\Tube(B,t)} \dc S^\sharp\wedge \beta^{q-1} 
 \\ &-  \big( {1\over r_1^{2q}}-{1\over r_2^{2q}}  \big) \int_{r_0}^{r_1}2tdt\int_{\partial_\ver\Tube(B,t)} \dc S^\sharp\wedge \beta^{q-1}
 +{1\over r_2^{2q}} \int_{\partial_\ver \Tube(B,r_2)}\dc\varphi\wedge S^\sharp\wedge \beta^{q-1} \\&-{1\over r_1^{2q}} \int_{\partial_\ver \Tube(B,r_1)}\dc\varphi\wedge S^\sharp\wedge \beta^{q-1}
  -\int_{\partial_\ver \Tube(B,r_1,r_2)}\dc\log\varphi\wedge S^\sharp\wedge \alpha^{q-1}.
  \end{split}
\end{equation}
 \item[(1-ii)] If $S$ is a  $\Cc^2$-smooth form, then identity \eqref{e:Lelong-Jensen}  (with $S_n:=S$ for $n\geq 1$) holds   for  all $r_1,r_2\in  (r_0,r]$ with  $r_1<r_2.$
 
 \item[(1-iii)] If for  all $n,$ $S_n$ is  a $\Cc^1$-smooth form of bidimension $(q,q)$  such that $\ddc S_n=0,$  then identity \eqref{e:Lelong-Jensen} holds  for  all $r_1,r_2\in  (r_0,r]$ with  $r_1<r_2$  except for  a  countable set of   values.

 \item[(1-iv)] If $S$ is  a $\Cc^1$-smooth form of bidimension $(q,q)$  such that $\ddc S=0,$  then identity \eqref{e:Lelong-Jensen}  (with 
 $S_n:=S$  for $n\geq 1$) holds  for  all $r_1,r_2\in  (r_0,r]$ with  $r_1<r_2.$  
 
\end{itemize}
\item  Assume  that $r_0>0.$ Then  the following four sub-assertions hold:
\begin{itemize}
\item  [(2-i)] For  all $r_2\in  (r_0,r]$  except for  a  countable set of   values,  we have that  
\begin{equation}\label{e:Lelong-Jensen-bis}\begin{split}
 {1\over  r_2^{2q}} \int_{\Tube(B,r_2)} S\wedge \beta^q
 &= \lim\limits_{n\to\infty}\Vc(S_n,r_2)+   \int_{\Tube(B,r_2)} S\wedge \alpha^q\\
 &+  \int_{r_0}^{r_2} \big( {1\over t^{2q}}-{1\over r_2^{2q}}  \big)2tdt\int_{\Tube(B,t)} \ddc S\wedge \beta^{q-1} 
 .
 \end{split}
\end{equation} 
Here the vertical boundary term  $\Vc(S,r_2)$  for a $\Cc^1$-smooth form $S$ is  given by the following formula:
\begin{equation}\label{e:vertical-boundary-term-bisbis}
\begin{split}
 \Vc(S,r_2)&:=-\int_{r_0}^{r_2} \big( {1\over t^{2q}}-{1\over r_2^{2q}}  \big)2tdt\int_{\partial_\ver\Tube(B,t)} \dc S^\sharp\wedge \beta^{q-1} 
 \\ &
 +{1\over r_2^{2q}} \int_{\partial_\ver \Tube(B,r_2)}\dc\varphi\wedge S^\sharp\wedge \beta^{q-1}
  -\int_{\partial_\ver \Tube(B,r_2)}\dc\log\varphi\wedge S^\sharp\wedge \alpha^{q-1}.
  \end{split}
\end{equation}
\item[(2-ii)] If $S$ is a  $\Cc^2$-smooth form, then identity \eqref{e:Lelong-Jensen-bis} holds  for  all $r_2\in  (r_0,r].$

 \item[(1-iii)] If for  all $n,$ $S_n$ is  a $\Cc^1$-smooth form of bidimension $(q,q)$  such that $\ddc S_n=0,$  then identity \eqref{e:Lelong-Jensen-bis} holds  for  all $r_2\in  (r_0,r]$   except for  a  countable set of   values.

 \item[(1-iv)] If $S$ is  a $\Cc^1$-smooth form of bidimension $(q,q)$  such that $\ddc S=0,$  then identity \eqref{e:Lelong-Jensen-bis}  (with 
 $S_n:=S$  for $n\geq 1$) holds  for  all $r_2\in  (r_0,r].$

\end{itemize}
\end{enumerate}
\end{theorem}

For the  proof of this  theorem the following two lemmas are needed.

\begin{lemma}\label{L:Lelong-Jensen_1}
For  every $t\in(r_0,\bfr],$ let  $j_t:\  \partial_\hor\Tube(B,t)\to \E$ be the  canonical  injection of the  real submanifold  $\partial_\hor\Tube(B,t)$ into $\E.$
Then  we have
$$
j^*_t(\alpha)={1\over t^2} j^*_t(\beta).
$$
\end{lemma}
\proof
Since we have
$$
j^*_t(\partial\varphi)+j^*_t(\dbar \varphi)=j^*_t(d\varphi)=d(\varphi\circ j_t)=0,
$$
it follows that 
$$
j^*_t(\partial\varphi)\wedge  j^*_t(\dbar\varphi)=0.
$$
On the other hand, a straightforward computation shows that
$$
\ddbar \log\varphi ={1\over\varphi} \ddbar \varphi- {1\over\varphi^2} \partial\varphi\wedge\dbar\varphi.
$$
Hence,
$$
j^*_t(\ddbar \log\varphi)={1\over t^2} j^*_t(\ddbar \varphi)\qquad\text{and}\qquad j^*_t(\ddc \log\varphi)={1\over t^2} j^*_t\beta.
$$
\endproof

\begin{lemma}\label{L:Lelong-Jensen_2}
Let $u$  be  a function  and  $\gamma$   a smooth  form of bidegree $(q-1,q-1).$
Let $S$ be  a smooth form of  bidimension $(q,q).$ Then
 we have
$$
\dc u\wedge dS \wedge \gamma=-du\wedge \dc S\wedge\gamma.
$$
\end{lemma}
\proof
By bidegree consideration, we have  that 
\begin{eqnarray*}
 \dc u\wedge dS \wedge \gamma&=&(i\dbar u\wedge \partial S-i\partial u\wedge \dbar S)\wedge\gamma,\\
 d u\wedge \dc S \wedge \gamma&=&(i\partial u\wedge \dbar S-i\dbar u\wedge \partial  S)\wedge\gamma.
 \end{eqnarray*}
The result follows.
\endproof

\proof[Proof of Theorem \ref{T:Lelong-Jensen}]
First  we  assume  that $S$ is  a $\Cc^2$-smooth form of bidimension $(q,q).$

We  will prove sub-assertion (1-ii). Write
\begin{equation}\label{e:Lelong-Jensen-1}
\int_{\Tube(B,r_1,r_2)} S\wedge \alpha^q=\int_{\Tube(B,r_1,r_2)}d\left\lbrack (\dc\log\varphi)\wedge S\wedge \alpha^{q-1}\right\rbrack+ 
\int_{\Tube(B,r_1,r_2)} (\dc\log\varphi)\wedge dS\wedge \alpha^{q-1}.
\end{equation}
Consider the  quantity
\begin{equation}\label{e:Lelong-Jensen-2}
 J:= \int_{\Tube(B,r_1,r_2)}(\dc\log\varphi)\wedge dS\wedge \alpha^{q-1},
\end{equation}
and the following one  for $r\in[r_0,\bfr]$:
\begin{equation}\label{e:Lelong-Jensen-I(r)}
 I(r):= \int_{\partial \Tube(B,r)}\dc\log\varphi\wedge S\wedge \alpha^{q-1}.
\end{equation}
Consequently,
Stokes' formula applied  in  \eqref{e:Lelong-Jensen-1} to the manifold  with boundary $\overline\Tube(B,r_1,r_2)$ using  the algebraic  identity
$\partial  \Tube (B,r_1,r_2)=\partial  \Tube (B,r_2)- \partial  \Tube (B,r_1)$  gives that 
\begin{equation}\label{e:Lelong-Jensen-3}
 \int_{\Tube(B,r_1,r_2)} S\wedge \alpha^q=I(r_2)-I(r_1)+J.
\end{equation}
  Using   \eqref{e:Lelong-Jensen-I(r)} and  the identity   $\partial \Tube(B,r)=\partial_\hor \Tube(B,r)\cup \partial_\ver \Tube(B,r)$,
we see that  $I(r)$ is  equal to
\begin{eqnarray*} 
 &&\int_{\partial_\hor \Tube(B,r)}\dc\log\varphi\wedge S\wedge \alpha^{q-1}+\int_{\partial_\ver \Tube(B,r)}\dc\log\varphi\wedge S\wedge \alpha^{q-1}\\
 &=&{1\over r^{2q}} \int_{\partial_\hor \Tube(B,r)}\dc\varphi\wedge S\wedge \beta^{q-1}+\int_{\partial_\ver \Tube(B,r)}\dc\log\varphi\wedge S\wedge \alpha^{q-1}\\
 &=&  {1\over r^{2q}} \int_{\partial \Tube(B,r)}\dc\varphi\wedge S\wedge \beta^{q-1} -{1\over r^{2q}} \int_{\partial_\ver \Tube(B,r)}\dc\varphi\wedge S\wedge \beta^{q-1}\\
 &+&\int_{\partial_\ver \Tube(B,r)}\dc\log\varphi\wedge S\wedge \alpha^{q-1},
\end{eqnarray*}
where  for the  first  integral in the second line  we have applied  Lemma \ref{L:Lelong-Jensen_1}, and for the  third line  we have used the algebraic identity  $\partial_\hor \Tube(B,r)=\partial \Tube(B,r)- \partial_\ver \Tube(B,r).$

Stokes' formula applied  to the first integral  of the last line  
gives that
\begin{eqnarray*}
 I(r)&=&{1\over r^{2q}} \int_{ \Tube(B,r)} S\wedge \beta^q-  
  {1\over r^{2q}} \int_{\Tube(B,r)}\dc\varphi\wedge dS\wedge \beta^{q-1}\\
  &-&{1\over r^{2q}} \int_{\partial_\ver \Tube(B,r)}\dc\varphi\wedge S\wedge \beta^{q-1}+\int_{\partial_\ver \Tube(B,r)}\dc\log\varphi\wedge S\wedge \alpha^{q-1}.
\end{eqnarray*}
Next,  applying  Lemma \ref{L:Lelong-Jensen_2}  to $u:=\varphi,$ $\gamma:=\beta^{q-1}$  in the second integral on the RHS yields that
\begin{eqnarray*}
 I(r)&=&{1\over r^{2q}} \int_{ \Tube(B,r)} S\wedge \beta^q+
  {1\over r^{2q}} \int_{\Tube(B,r)}d\varphi\wedge \dc S\wedge \beta^{q-1}\\
   &-&{1\over r^{2q}} \int_{\partial_\ver \Tube(B,r)}\dc\varphi\wedge S\wedge \beta^{q-1}+\int_{\partial_\ver \Tube(B,r)}\dc\log\varphi\wedge S\wedge \alpha^{q-1}.
\end{eqnarray*}
Applying Fubini's theorem (see  \cite[4.3.2., (1)]{Federer}) or a variant  (see \cite[7.2.]{Siu}) to the second integral on the RHS and 
using  that $d\varphi=2tdt$ for $\varphi=t^2,$   we get that
\begin{eqnarray*}
 I(r)&=&{1\over r^{2q}} \int_{ \Tube(B,r)} S\wedge \beta^q+
  {1\over r^{2q}}\int_{r_0}^r 2tdt \int_{\partial_\hor \Tube(B,t)} \dc S\wedge \beta^{q-1}\\
  &-&{1\over r^{2q}} \int_{\partial_\ver \Tube(B,r)}\dc\varphi\wedge S\wedge \beta^{q-1}+\int_{\partial_\ver \Tube(B,r)}\dc\log\varphi\wedge S\wedge \alpha^{q-1} .
\end{eqnarray*}
Since  $\partial_\hor \Tube(B,t)=\partial \Tube(B,t)- \partial_\ver \Tube(B,t),$ it follows that
\begin{eqnarray*}
 I(r)&=&{1\over r^{2q}} \int_{ \Tube(B,r)} S\wedge \beta^q+
  {1\over r^{2q}}\int_{r_0}^r 2tdt \int_{\partial \Tube(B,t)} \dc S\wedge \beta^{q-1}\\
  &-&{1\over r^{2q}}\int_{r_0}^r 2tdt \int_{\partial_\ver \Tube(B,t)} \dc S\wedge \beta^{q-1}\\
  &-&{1\over r^{2q}} \int_{\partial_\ver \Tube(B,r)}\dc\varphi\wedge S\wedge \beta^{q-1}+\int_{\partial_\ver \Tube(B,r)}\dc\log\varphi\wedge S\wedge \alpha^{q-1} .
\end{eqnarray*}
Applying  Stokes' formula to the inner  integral of the  first double  integral on the RHS, the last line is  equal to
\begin{equation}\label{e:Lelong-Jensen-4}
\begin{split}
 I(r)&={1\over r^{2q}} \int_{ \Tube(B,r)} S\wedge \beta^q+
  {1\over r^{2q}}\int_{r_0}^r 2tdt \int_{\Tube(B,t)} \ddc S\wedge \beta^{q-1}\\
   &-{1\over r^{2q}}\int_{r_0}^r 2tdt \int_{\partial_\ver \Tube(B,t)} \dc S\wedge \beta^{q-1}\\
  &-{1\over r^{2q}} \int_{\partial_\ver \Tube(B,r)}\dc\varphi\wedge S\wedge \beta^{q-1}+\int_{\partial_\ver \Tube(B,r)}\dc\log\varphi\wedge S\wedge \alpha^{q-1} .
  \end{split}
\end{equation}
Rewrite  \eqref{e:Lelong-Jensen-2} using Lemma \ref{L:Lelong-Jensen_2} with $u:=\log\varphi,$ $\gamma:=\alpha^{q-1},$
\begin{equation*}
 J= -\int_{\Tube(B,r_1,r_2)}(d\log\varphi)\wedge \dc S\wedge \alpha^{q-1}.
\end{equation*}
By  Fubini's theorem, we get that
\begin{equation*}
 J= -\int_{r_1}^{r_2}{2dt\over t}  \int_{\partial_\hor \Tube(B,t)} \dc S\wedge \alpha^{q-1}.
\end{equation*}
By Lemma \ref{L:Lelong-Jensen_1}   applied to $\partial_\hor \Tube(B,r)$ and the  equality $\varphi(y)=t^2$  for $y\in  \partial_\hor \Tube(B,t),$  we obtain that 
\begin{eqnarray*}
 J&=& -\int_{r_1}^{r_2}{2tdt\over t^{2q}}  \int_{\partial_\hor \Tube(B,t)} \dc S\wedge \beta^{q-1}=-\int_{r_1}^{r_2}{2tdt\over t^{2q}}  \int_{\partial \Tube(B,t)} \dc S\wedge \beta^{q-1}\\
 &+&\int_{r_1}^{r_2}{2tdt\over t^{2q}}  \int_{\partial_\ver \Tube(B,t)} \dc S\wedge \beta^{q-1},
\end{eqnarray*}
where the second equality holds since $\partial_\hor \Tube(B,t)=\partial\Tube(B,t) -\partial_\ver \Tube(B,t).$
Stokes' formula applied  to the first integral on the RHS   gives that
\begin{equation}\label{e:Lelong-Jensen-5}
 J= -\int_{r_1}^{r_2}{2tdt\over t^{2q}}  \int_{ \Tube(B,t)} \ddc S\wedge \beta^{q-1}+\int_{r_1}^{r_2}{2tdt\over t^{2q}}  \int_{\partial_\ver \Tube(B,t)} \dc S\wedge \beta^{q-1}.
\end{equation}
This, combined with \eqref{e:Lelong-Jensen-3}--\eqref{e:Lelong-Jensen-4}, implies that
\begin{eqnarray*}
 \int_{\Tube(B,r_1,r_2)} S\wedge \alpha^q&=&{1\over  r_2^{2q}} \int_{\Tube(B,r_2)} S\wedge \beta^q- {1\over  r_1^{2q}} \int_{\Tube(B,r_1)} S\wedge \beta^q-\int_{r_1}^{r_2}{2tdt\over t^{2q}}  \int_{ \Tube(B,t)} \ddc S\wedge \beta^{q-1}
 \\
 &+&{1\over r_2^{2q}}\int_{r_0}^{r_2}2tdt  \int_{ \Tube(B,t)} \ddc S\wedge \beta^{q-1}
 -{1\over r_1^{2q}}\int_{r_0}^{r_1}2tdt  \int_{ \Tube(B,t)} \ddc S\wedge \beta^{q-1}\\
  &-&{1\over r_2^{2q}}\int_{r_0}^{r_2} 2tdt \int_{\partial_\ver \Tube(B,t)} \dc S\wedge \beta^{q-1}+
  {1\over r_1^{2q}}\int_{r_0}^{r_1} 2tdt \int_{\partial_\ver \Tube(B,t)} \dc S\wedge \beta^{q-1}
  \\
  &-&{1\over r_2^{2q}} \int_{\partial_\ver \Tube(B,r_2)}\dc\varphi\wedge S\wedge \beta^{q-1}+{1\over r_1^{2q}} \int_{\partial_\ver \Tube(B,r_1)}\dc\varphi\wedge S\wedge \beta^{q-1}\\
  &+&\int_{\partial_\ver \Tube(B,r_2)}\dc\log\varphi\wedge S\wedge \alpha^{q-1}-\int_{\partial_\ver \Tube(B,r_1)}\dc\log\varphi\wedge S\wedge \alpha^{q-1}\\
& +&\int_{r_1}^{r_2}{2tdt\over t^{2q}}  \int_{\partial_\ver \Tube(B,t)} \dc S\wedge \beta^{q-1}
 .
\end{eqnarray*}
So formula \eqref{e:Lelong-Jensen} holds for  all $r_1,r_2\in  [r_0,r]$ with  $r_1<r_2.$
This  completes the proof of assertion (1-ii)  for the case  when   $S$ is  a $\Cc^2$-smooth form  of bidimension $(q,q).$

We  turn to the proof of assertion (2).  Roughly speaking,  assertion (1) for $r_1:=r_0$ becomes  assertion (2).
More precisely, consider the  quantity
\begin{equation}\label{e:Lelong-Jensen-2bis}
 \widetilde J:= \int_{\Tube(B,r_2)}(\dc\log\varphi)\wedge dS\wedge \alpha^{q-1}.
\end{equation}
Consequently,
Stokes' formula applied  in  \eqref{e:Lelong-Jensen-1} to the manifold  with boundary $\overline\Tube(B,r_2)$ using  the algebraic  identity
$\partial  \Tube (B,r_1,r_2)=\partial  \Tube (B,r_2)- \partial  \Tube (B,r_1)$  gives that 
\begin{equation}\label{e:Lelong-Jensen-3bis}
 \int_{\Tube(B,r_2)} S\wedge \alpha^q=I(r_2)+\widetilde J.
\end{equation}
Arguing as  in the  proof of \eqref{e:Lelong-Jensen-5}, we obtain that
\begin{equation}
 \widetilde J= -\int_{r_0}^{r_2}{2tdt\over t^{2q}}  \int_{ \Tube(B,t)} \ddc S\wedge \beta^{q-1}+\int_{r_0}^{r_2}{2tdt\over t^{2q}}  \int_{\partial_\ver \Tube(B,t)} \dc S\wedge \beta^{q-1}.
\end{equation}
This, combined with \eqref{e:Lelong-Jensen-2bis}, \eqref{e:Lelong-Jensen-3bis} and \eqref{e:Lelong-Jensen-4}, implies  assertion (2).
Hence, we have  proved the theorem  for the case  when   $S$ is  a $\Cc^2$-smooth form  of bidimension $(q,q).$ 

Next, we treat the case when   $S$ is  a $\Cc^2$-smooth form  of dimension $2q.$ 
We only give the proof of assertion (1) since  assertion (2) can be proved similarly.
In this  case  we only need  to apply   the  previous case to $S^\sharp$
and  observe  that by a  consideration of bidegree in formula \eqref{e:Lelong-Jensen},  $S^\sharp$ can be replace  by $S$ except for the vertical boundary
term $\Vc(S,r_1,r_2).$  Hence, this case is done.

Now  we pass to the more  general case where  $S$ is  a current  of dimension $2q$  with  an approximating $\Cc^2$-smooth forms $S_n$ of dimension $2q$
as  in the  hypothesis.   We apply the  previous  case to each   form $S_n$ and then we take  the limit in each term  of formula \eqref{e:Lelong-Jensen}.
Arguing as in the proof of \eqref{e:except-countable} in  Lemma \ref{L:quasi-positivity}, we conclude that formula  \eqref{e:Lelong-Jensen} holds for  all $r_1,r_2\in  [r_0,r]$ with  $r_1<r_2$  except for a countable set of values.
\endproof

The next theorem deals with the special case where  the current  is  approximable by smooth {\bf closed}  forms with  control on the boundary.
Here, we gain the  smoothness.

\begin{theorem}\label{T:Lelong-Jensen-closed}
Let $r\in[r_0,\bfr]$ and  let $S$ be a real closed current of dimension $2q$ on  a neighborhood of $\overline\Tube(B,r).$ 
%
%
Suppose that there is a sequence of $\Cc^1$-smooth closed forms of dimension $2q:$  $(S_n)_{n=1}^\infty$ defined on  a neighborhood of $\overline\Tube(B,r)$ 
such  that
$S_n$ converge to $S$  in the sense of quasi-positive currents   on a neighborhood of  $\overline\Tube(B,r)$ as $n$ tends to infinity (see Definition \ref{D:quasi-positivity}).

Then the following  two assertions  hold:
\begin{enumerate} \item  The  following two  sub-assertions hold:
\begin{itemize}
 \item [(1-i)]
For  all $r_1,r_2\in  [r_0,r]$ with  $r_1<r_2$  except for  a  countable set of   values, we have that  
\begin{equation}\label{e:Lelong-Jensen-closed}
 {1\over  r_2^{2q}} \int_{\Tube(B,r_2)} S\wedge \beta^q- {1\over  r_1^{2q}} \int_{\Tube(B,r_1)} S\wedge \beta^q
 =\lim\limits_{n\to\infty} \Vc(S_n,r_1,r_2)+   \int_{\Tube(B,r_1,r_2)} S\wedge \alpha^q.
\end{equation} 
Here the vertical boundary term  $\Vc(S,r_1,r_2)$ for  a  continuous form $S$  is  given by
\begin{equation}\label{e:vertical-boundary-term-closed}
\begin{split}
 \Vc(S,r_1,r_2)&:=
 {1\over r_2^{2q}} \int_{\partial_\ver \Tube(B,r_2)}\dc\varphi\wedge S\wedge \beta^{q-1}-{1\over r_1^{2q}} \int_{\partial_\ver \Tube(B,r_1)}\dc\varphi\wedge S\wedge \beta^{q-1}\\
  &-\int_{\partial_\ver \Tube(B,r_1,r_2)}\dc\log\varphi\wedge S\wedge \alpha^{q-1}.
  \end{split}
\end{equation}
\item [(1-ii)] If $S$ is a closed $\Cc^1$-smooth form, then identity \eqref{e:Lelong-Jensen-closed} (with $S_n:= S$ for $n \geq 1$) holds
for  all $r_1,r_2\in  [r_0,r]$ with  $r_1<r_2.$
\end{itemize}
\item If $r_0>0,$  then the following  two sub-assertions hold:
\begin{itemize}
\item[(2-i)]  For  all $r_2\in  [r_0,r]$   except for  a  countable set of   values, we have that  
\begin{equation}\label{e:Lelong-Jensen-closed-bis}
 {1\over  r_2^{2q}} \int_{\Tube(B,r_2)} S\wedge \beta^q
 = \lim\limits_{n\to\infty} \Vc(S_n,r_2)+   \int_{\Tube(B,r_2)} S\wedge \alpha^q.
\end{equation} 
Here the vertical boundary term  $\Vc(S,r_2)$ for  a  continuous form $S$ is  given by
\begin{equation}\label{e:vertical-boundary-term-closed-bisbis}
 \Vc(S,r_2):=
 {1\over r_2^{2q}} \int_{\partial_\ver \Tube(B,r_2)}\dc\varphi\wedge S\wedge \beta^{q-1}-\int_{\partial_\ver \Tube(B,r_2)}\dc\log\varphi\wedge S\wedge \alpha^{q-1}.
\end{equation}
\item [(2-ii)] If $S$ is a closed $\Cc^1$-smooth form, then identity \eqref{e:Lelong-Jensen-closed-bis} (with $S_n:= S$ for $n \geq 1$) holds
for  all $r_2\in  [r_0,r].$ 
\end{itemize}
\end{enumerate}
\end{theorem}

\proof
We only treat the case where $S$ is  a $\Cc^1$-smooth form  on  a neighborhood of $\overline\Tube(B,r)$
and we only give the proof of assertion (1).
We  follow  the proof of Theorem \ref{T:Lelong-Jensen} by making the following   observation. All terms  containing $dS,$  $\ddc S$
vanishes, for example $J$ in \eqref{e:Lelong-Jensen-2}. Moreover, we do not use Lemma \ref{L:Lelong-Jensen_2}.
Consequently, instead of \eqref{e:Lelong-Jensen-4}, we get that
\begin{equation*}
 I(r)={1\over r^{2q}} \int_{ \Tube(B,r)} S\wedge \beta^q
  -{1\over r^{2q}} \int_{\partial_\ver \Tube(B,r)}\dc\varphi\wedge S\wedge \beta^{q-1}
  +\int_{\partial_\ver \Tube(B,r)}\dc\log\varphi\wedge S\wedge \alpha^{q-1}  .
\end{equation*}
This, coupled with \eqref{e:Lelong-Jensen-3} and $J=0$ and  assumption (iii) implies the desired conclusion.
\endproof

Now  we consider  a special case where the vertical boundary term $\Vc(S,r_1,r_2)$ defined  in \eqref{e:vertical-boundary-term} vanishes.   

\begin{theorem}\label{T:Lelong-Jensen-compact-support}
 Let $r\in\R^+_*$ and   let $S$ be a real current of dimension $2q$ on  a neighborhood of $\overline\Tube(B,r)$ such that $S$ and $\ddc S$ are of order  $0.$
Suppose that there is a sequence of $\Cc^2$-smooth $2q$-forms $(S_n)_{n=1}^\infty$ defined on  a neighborhood of $\overline\Tube(B,r)$ 
such  that 
\begin{itemize}
 \item[(i)] 
$S_n$ converge to $S$  in the sense of quasi-positive currents   on a neighborhood of  $\overline\Tube(B,r)$ as $n$ tends to infinity;
\item [(ii)]  $\ddc S_n$ converge to $\ddc S$  in the sense of quasi-positive currents   on a neighborhood of  $\overline\Tube(B,r)$ as $n$ tends to infinity;
\item[(iii)]  the following  equalities hold:
\begin{eqnarray*} \lim_{n\to\infty}\|S^\sharp_n\| (\partial_\ver \Tube( B, r))&=&0\quad\text{and}\quad  \lim_{n\to\infty}\|\partial S^\sharp_n\| (\partial_\ver \Tube( B, r))=0\\
 \qquad \text{and}\quad  \lim_{n\to\infty}\|\dbar S^\sharp_n\| (\partial_\ver \Tube( B, r))&=&0.
 \end{eqnarray*}
\end{itemize}  
Then,  for  all $r_1,r_2\in [r_0,r]$  with  $r_1<r_2,$ except for a   countable set of values,
the vertical boundary term $\Vc(S,r_1,r_2)$ vanishes
and 
\begin{multline*}
 {1\over  r_2^{2q}} \int_{\Tube(B,r_2)} S\wedge \beta^q- {1\over  r_1^{2q}} \int_{\Tube(B,r_1)} S\wedge \beta^q
 =\int_{\Tube(B,r_1,r_2)} S\wedge \alpha^q\\
 +  \int_{r_1}^{r_2} \big( {1\over t^{2q}}-{1\over r_2^{2q}}  \big)2tdt\int_{\Tube(B,t)} \ddc S\wedge \beta^{q-1} 
 +  \big( {1\over r_1^{2q}}-{1\over r_2^{2q}}  \big) \int_{0}^{r_1}2tdt\int_{\Tube(B,t)} \ddc S\wedge \beta^{q-1}.
\end{multline*}
 In particular,  when  
$\supp (S^\sharp_n)\cap \partial_\ver \Tube( B, r)=\varnothing$ for all $n\geq 1,$ 
then condition  (iii) above is  automatically satisfied and  the above  formula  holds  whenever conditions (i)--(ii)  are fulfilled.
\end{theorem}
\proof
The  second assertion of  the  theorem   follows immediately from  the first one.

We now prove  the  first assertion. Applying   assumption (iii)  to  formula \eqref{e:vertical-boundary-term} and  \eqref{e:continuity-cut-off} in Lemma \ref{L:quasi-positivity}, we conclude that  $\lim_{n\to\infty}\Vc(S_n,r_1,r_2)=0.$
Hence, the first  assertion follows from Theorem \ref{T:Lelong-Jensen}.
\endproof

We also need a version of Theorem \ref{T:Lelong-Jensen-compact-support}  when the current is approximable by smooth closed forms.
\begin{theorem}\label{T:Lelong-Jensen-closed-compact-support}
 Let $r\in\R^+_*$ and   let $S$ be a real closed  current of dimension $2q$ on  a neighborhood of $\overline\Tube(B,r).$ 
Suppose that there is a sequence of $\Cc^1$-smooth closed $2q$-forms $(S_n)_{n=1}^\infty$ defined on  a neighborhood of $\overline\Tube(B,r)$ 
such  that 
\begin{itemize}
 \item[(i)] 
$S_n$ converge to $S$  in the sense of quasi-positive currents   on a neighborhood of  $\overline\Tube(B,r)$ as $n$ tends to infinity;
\item[(ii)]  the following  equality holds
$ \lim_{n\to\infty}\|S_n\| (\partial_\ver \Tube( B, r))=0.$
\end{itemize}  
Then,  for  all $r_1,r_2\in [r_0,r]$  with  $r_1<r_2,$ except for a   countable set of values,
the vertical boundary term $\Vc(S,r_1,r_2)$ vanishes
and 
\begin{equation*}
 {1\over  r_2^{2q}} \int_{\Tube(B,r_2)} S\wedge \beta^q- {1\over  r_1^{2q}} \int_{\Tube(B,r_1)} S\wedge \beta^q
 =\int_{\Tube(B,r_1,r_2)} S\wedge \alpha^q. 
\end{equation*}
  In  particular,   when $\supp (S_n)\cap \partial_\ver \Tube( B, r)=\varnothing$ for all $n\geq 1,$
  then condition  (ii) above is  automatically satisfied and  the above  formula  holds  whenever condition (i) is fulfilled.
\end{theorem}
\proof
The  second assertion of  the  theorem   follows immediately from  the first one.

We now  prove  the first one.
 Applying   assumption (ii) to the formula  of  $\Vc(S_n^,r_1,r_2)$
given by \eqref{e:vertical-boundary-term-closed} and \eqref{e:continuity-cut-off} in Lemma \ref{L:quasi-positivity}, we conclude that  $\lim_{n\to\infty}\Vc(S_n,r_1,r_2)=0.$
Hence, the result follows from   Theorem \ref{T:Lelong-Jensen-closed}.
\endproof
For the remainder of the  section, we fix $0\leq p\leq k$ and  recall from \eqref{e:m} that
\begin{equation*}\upm:= \min(l,k-p)\qquad\text{and}\qquad
  \lowm:=\max(0,l-p).
  \end{equation*}
  
 As  an immediate consequence,  we obtain in a particular  situation  the following  Lelong-Jensen formula without  boundary support condition.
\begin{corollary}\label{C:Lelong-Jensen}
  Let $r\in\R^+_*$ and
    let $\omega$ be a  smooth $(1,1)$-form  defined  on a neighborhood of $\overline B$   in $V.$
   Let $T$ be a  real current of degree $2p$  and $(T_n)_{n=1}^\infty$    a sequence of $\Cc^2$-smooth $2p$-forms
defined on  a neighborhood of $\overline\Tube(B,r)$ satisfying  the following properties: 
   \begin{itemize}
\item[(i)] 
$T_n$ converge to $T$  in the sense of quasi-positive currents   on a neighborhood of  $\overline\Tube(B,r)$ as $n$ tends to infinity (see Definition \ref{D:quasi-positivity});
\item [(ii)]  $\ddc T_n$ converge to $\ddc T$  in the sense of quasi-positive currents   on a neighborhood of  $\overline\Tube(B,r)$ as $n$ tends to infinity.
  \end{itemize} 
   Then    the conclusion of  Theorem  \ref{T:Lelong-Jensen-compact-support} holds
  with the  $(2p+2\upm)$-current $S:=T\wedge \pi^*\omega^\upm$  and $q:=k-(p+\upm).$  
\end{corollary}

\proof  Consider the $\Cc^2$-smooth  $(2p+2\upm)$-forms
 $S_n=T_n\wedge \pi^*\omega^\upm.$
Consider  a small neighborhood $V(x_0)$ of  an arbitrary  point $x_0\in \partial_\ver \Tube(B, r),$  where in a local chart $V(x_0)\simeq \D^l$ and  $\E|_{V(x_0)}\simeq \C^{k-l}\times \D^l.$
For $y\in \E|_{V(x_0)},$ write $y=(z,w).$  We  will 
prove the  following 

\noindent {\bf Fact.} {\it  $S_n$ is  of bidegree $(l,l)$  in $dw,$ $d\bar w.$}

Indeed, there are two cases  to consider.

If  $\upm=l,$  then  clearly $ \pi^*\omega^\upm$ is  of bidegree $(l,l)$  in $dw,$ $d\bar w,$
and the  above fact  follows  because of the above  formula of  $S_n.$

Otherwise, we have  $k-p<l$ and   $\upm=k-p.$  In this  case  $ \pi^*\omega^\upm$ is  of bidegree $(k-p,k-p)$ and  every component  of $T_n$ should contain $dw_I\wedge d\bar w_J$   with $|I|+|J|\geq 2(p-k+l).$ 
Since $S_n=T_n\wedge \pi^*\omega^\upm$, we see that $S_n$ can be factorized  by   $ \pi^*\omega^l,$ and hence the  above fact also follows in this  last case.
 
Since  $\dim_\R (\partial B)=2l-1,$  it follows from the  above fact that 
$$ \|S_n\| (\partial_\ver \Tube( B, r))=0\quad\text{and}\quad \|\partial S_n\| (\partial_\ver \Tube( B, r))=0 \quad\text{and}\quad \|\dbar S_n\| (\partial_\ver \Tube( B, r))=0.$$
Hence,   Theorem \ref{T:Lelong-Jensen-compact-support} gives the desired conclusion.
\endproof

\begin{corollary}\label{C:Lelong-Jensen-closed}
  Let $r\in\R^+_*$ and
    let $\omega$ be a  smooth $(1,1)$-form  defined  on a neighborhood of $\overline B$   in $V.$
   Let $T$ be a  real current of degree $2p$  and $(T_n)_{n=1}^\infty$    a sequence of $\Cc^1$-smooth $2p$-forms
defined on  a neighborhood of $\overline\Tube(B,r)$  such that
$T_n$ converge to $T$  in the sense of quasi-positive currents   on a neighborhood of  $\overline\Tube(B,r)$ as $n$ tends to infinity (see Definition \ref{D:quasi-positivity}).
   Then    the conclusion of  Theorem  \ref{T:Lelong-Jensen-closed-compact-support} holds
  with the  $(2p+2\upm)$-current $S:=T\wedge \pi^*\omega^\upm$  and $q:=k-(p+\upm).$  
\end{corollary}

\proof
We argue as in the proof of Corollary  \ref{C:Lelong-Jensen} replacing  Theorem   \ref{T:Lelong-Jensen-compact-support}  by Theorem   \ref{T:Lelong-Jensen-closed-compact-support}.
\endproof

\subsection{Applications}
Consider a Hermitian metric  $\|\cdot\|$  on the  vector bundle  $\E$  and    let   $\varphi:\ \E\to \R^+$ be the function defined by  
\begin{equation}\label{e:varphi}
 \varphi(y):=\|y\|^2\qquad \text{for}\qquad  y\in \E.
\end{equation}
Since for every $x\in X$  the  metric $\| \cdot\|$  on the fiber $\E_x\simeq \C^{k-l}$ is  an Euclidean metric (in a suitable basis), we have 
\begin{equation}\label{e:varphi_bis}  \varphi(\lambda y)=|\lambda|^2\varphi(y)\qquad\text{for}\qquad  y\in \E,\qquad\lambda\in\C. 
\end{equation}
In this case  where we have $r_0=0$ and $\bfr=\infty,$ and the cooresponding tubes as well as  the corresponding  corona tubes are already  defined  
in \eqref{e:tubular-nbh-0} and \eqref{e:tubular-corona-0}.
Unless otherwise specified  we   consider  mainly  these  tubes and  corona tubes in this  work.

In some places we also consider the  following variant  of $\varphi$  in the spirit of \eqref{e:Lelong-number-point-bisbis(2)}:
for  every $\epsilon>0,$  set
\begin{equation}\label{e:varphi_eps}
 \varphi_\epsilon(y):=\|y\|^2+\epsilon^2\qquad \text{for}\qquad  y\in \E.
\end{equation}
In this case  where we have $r_0=\epsilon$ and $\bfr=\infty.$
Following the model \eqref{e:alpha-beta}, consider also the following  closed  $(1,1)$-form  for each $\epsilon >0$ on  $\U:$
\begin{equation}\label{e:alpha-beta-eps}
 \alpha_\epsilon:=\ddc\log\varphi_\epsilon\quad\text{and note that }\quad \beta= \ddc\varphi_\epsilon.
\end{equation}  
The  following   result which  will play a key role  for proving logarithmic   interpretation  version in the spirit of  \eqref{e:Lelong-number-point-bisbis(2)}.
\begin{theorem}\label{T:Lelong-Jensen-eps}
Let $\bfr\in\R^+_*$  and   $B\Subset V$  a relatively compact  open set with piecewice $\Cc^2$-smooth boundary.
Let $S$ be a real current of dimension $2q$ on  a neighborhood of $\overline\Tube(B,\bfr).$ 
%
%
Suppose that there is a sequence of $\Cc^2$-smooth forms of dimension $2q$ $(S_n)_{n=1}^\infty$ defined on  a neighborhood of $\overline\Tube(B,\bfr)$ 
such  that
\begin{itemize}
 \item[(i)] 
$S_n$ converge to $S$  in the sense of quasi-positive currents   on a neighborhood of  $\overline\Tube(B,\bfr)$ as $n$ tends to infinity (see Definition \ref{D:quasi-positivity});
\item [(ii)]  $\ddc S_n$ converge to $\ddc S$  in the sense of quasi-positive currents   on a neighborhood of  $\overline\Tube(B,\bfr)$ as $n$ tends to infinity.
\end{itemize}
Then  the following  two assertions hold:
\begin{enumerate}
 \item 
For  all  $r\in (0,\bfr)$ and   $\epsilon\in  (0,r)$      except for a  countable  set of  values, we have that
\begin{multline*}
 {1\over  (r^2+\epsilon^2)^q} \int_{\Tube(B,r)} S\wedge \beta^q
 =\lim\limits_{n\to\infty} \Vc_\epsilon(S_n,r)+   \int_{\Tube(B,r)} S\wedge \alpha^q_\epsilon\\
 +  \int_{0}^{r} \big( {1\over (t^2+\epsilon^2)^q}-{1\over (r^2+\epsilon^2)^q}  \big)2tdt\int_{\Tube(B,t)} \ddc S\wedge \beta^{q-1} .
\end{multline*} 
Here the vertical boundary term  $\Vc_\epsilon(S,r)$ for  a  $\Cc^1$-smooth form $S$ is  given by
\begin{equation}\label{e:vertical-boundary-term-eps}
\begin{split}
 \Vc_\epsilon(S,r)&:=-\int_{0}^{r} \big( {1\over (t^2+\epsilon^2)^q}-{1\over (r^2+\epsilon^2)^q}  \big)2tdt\int_{\partial_\ver\Tube(B,t)} \dc S^\sharp\wedge \beta^{q-1} 
  \\
 &+{1\over  (r^2+\epsilon^2)^q} \int_{\partial_\ver \Tube(B,r)}\dc\varphi\wedge S^\sharp\wedge \beta^{q-1}
  -\int_{\partial_\ver \Tube(B,r)}\dc\log\varphi_\epsilon\wedge S^\sharp\wedge \alpha_\epsilon^{q-1}.
  \end{split}
\end{equation}

\item  If $S$  is a $\Cc^2$-smooth form, then the above identity (with $S_n:=S$ for $n\geq 1$) holds  for  all  $r\in (0,\bfr)$ and 
$\epsilon\in  (0,r).$ 

\end{enumerate}
\end{theorem}
\proof
Note that  $\Tube (B,r)=\{y\in \E:\  \varphi_\epsilon<r^2+\epsilon^2\}.$  Note also  by    \eqref{e:varphi_eps}
that  $\dc \varphi_\epsilon=\dc\varphi.$ Consequently, 
the  result  follows  from  Theorem  \ref{T:Lelong-Jensen} (2) applied  to $\varphi_\epsilon,$ $\alpha_\epsilon$ and $\beta$
given by    \eqref{e:varphi_eps}--\eqref{e:alpha-beta-eps}, and  to $r_1:=\epsilon$ and 
 $r_2:=\sqrt{r^2+\epsilon^2}.$ 
\endproof

We record a version of Theorem \ref{T:Lelong-Jensen-eps} for {\bf closed} currents.

\begin{theorem}\label{T:Lelong-Jensen-closed-eps}
Let $\bfr\in\R^+_*$ and  let $S$ be a real closed current of dimension $2q$ on  a neighborhood of $\overline\Tube(B,\bfr).$ 
%
%
Suppose that there is a sequence of $\Cc^1$-smooth closed forms of dimension $2q:$  $(S_n)_{n=1}^\infty$ defined on  a neighborhood of $\overline\Tube(B,\bfr)$ 
such  that
%
$S_n$ converge to $S$  in the sense of quasi-positive currents   on a neighborhood of  $\overline\Tube(B,\bfr)$ as $n$ tends to infinity (see Definition \ref{D:quasi-positivity}).
%
%
Then the  following two  assertions hold:
\begin{enumerate}
 \item 
For  all $r\in (0,\bfr]$ and   $\epsilon\in  (0,r)$  except for a countable set of values,   we have that 
\begin{equation}\label{e:Lelong-Jensen-closed-eps}
 {1\over  (r^2+\epsilon^2)^{q}} \int_{\Tube(B,r)} S\wedge \beta^q
 = \lim\limits_{n\to\infty} \Vc_\epsilon(S_n,r)+   \int_{\Tube(B,r)} S\wedge \alpha^q_\epsilon.
\end{equation} 
Here the vertical boundary term  $\Vc_\epsilon(S,r)$ for a  continuous   form $S$  is  given by
\begin{equation}\label{e:vertical-boundary-term-closed-eps}
 \Vc_\epsilon(S,r):=
 {1\over (r^2+\epsilon^2)^{q}} \int_{\partial_\ver \Tube(B,r)}\dc\varphi\wedge S\wedge \beta^{q-1}-
  \int_{\partial_\ver \Tube(B,r)}\dc\log\varphi_\epsilon\wedge S\wedge \alpha^{q-1}_\epsilon.
\end{equation}

\item  If $S$  is a closed  $\Cc^1$-smooth form, then the above identity (with $S_n:=S$ for $n\geq 1$) holds  for  all  $r\in (0,\bfr)$ and 
$\epsilon\in  (0,r).$
\end{enumerate}
\end{theorem}
\proof
Note that  $\Tube (B,r)=\{y\in \E:\  \varphi_\epsilon<r^2+\epsilon^2\}.$ Consequently, 
the  result  follows  from  Theorem  \ref{T:Lelong-Jensen-closed} (2) applied  to $\varphi_\epsilon$  and $r_2:=\sqrt{r^2+\epsilon^2}.$ 
\endproof

Now  we consider  a special case where the vertical boundary term $\Vc_\epsilon(S,r)$ defined  in \eqref{e:vertical-boundary-term-eps} vanishes.   

\begin{theorem}\label{T:Lelong-Jensen-compact-support-eps}
 Let $\bfr\in\R^+_*$ and   let $S$ be a real current of dimension $2q$ on  a neighborhood of $\overline\Tube(B,\bfr)$ such that $S$ and $\ddc S$ are of order  $0.$
Suppose that there is a sequence of $\Cc^2$-smooth $2q$-forms $(S_n)_{n=1}^\infty$ defined on  a neighborhood of $\overline\Tube(B,\bfr)$ 
such  that 
\begin{itemize}
 \item[(i)] 
$S_n$ converge to $S$  in the sense of quasi-positive currents   on a neighborhood of  $\overline\Tube(B,\bfr)$ as $n$ tends to infinity;
\item [(ii)]  $\ddc S_n$ converge to $\ddc S$  in the sense of quasi-positive currents   on a neighborhood of  $\overline\Tube(B,\bfr)$ as $n$ tends to infinity;
\item[(iii)]  the following  equalities hold:
\begin{eqnarray*} \lim_{n\to\infty}\|S^\sharp_n\| (\partial_\ver \Tube( B, \bfr))&=&0\quad\text{and}\quad  \lim_{n\to\infty}\|\partial S^\sharp_n\| (\partial_\ver \Tube( B, \bfr))=0\\
 \qquad \text{and}\quad  \lim_{n\to\infty}\|\dbar S^\sharp_n\| (\partial_\ver \Tube( B, \bfr))&=&0.
 \end{eqnarray*}
\end{itemize}  
 Then,  for  all $r\in (0,\bfr]$ and   $\epsilon\in  (0,r),$   
the vertical boundary term $\Vc_\epsilon(S,r)$ vanishes
and 
\begin{multline*}
  {1\over  (r^2+\epsilon^2)^q} \int_{\Tube(B,r)} S\wedge \beta^q
 =   \int_{\Tube(B,r)} S\wedge \alpha^q_\epsilon\\
 +  \int_{0}^{r} \big( {1\over (t^2+\epsilon^2)^q}-{1\over (r^2+\epsilon^2)^q}  \big)2tdt\int_{\Tube(B,t)} \ddc S\wedge \beta^{q-1} .
\end{multline*}
 In particular,  when  
$\supp (S_n)\cap \partial_\ver \Tube( B, \bfr)=\varnothing$ for all $n\geq 1,$ 
then  the above  formula  holds.
\end{theorem}
\proof
We combine  the proofs of Theorem \ref{T:Lelong-Jensen-compact-support} and Theorem \ref{T:Lelong-Jensen-eps}.
\endproof

Now  we consider  a special case where the vertical boundary term $\Vc_\epsilon(S,r)$ defined  in \eqref{e:vertical-boundary-term-closed-eps} vanishes.   

\begin{theorem}\label{T:Lelong-Jensen-compact-support-closed-eps}
 Let $\bfr\in\R^+_*$ and   let $S$ be a real current of dimension $2q$ on  a neighborhood of $\overline\Tube(B,\bfr).$
Suppose that there is a sequence of closed $\Cc^1$-smooth $2q$-forms $(S_n)_{n=1}^\infty$ defined on  a neighborhood of $\overline\Tube(B,\bfr)$ 
such  that 
\begin{itemize}
 \item[(i)] 
$S_n$ converge to $S$  in the sense of quasi-positive currents   on a neighborhood of  $\overline\Tube(B,\bfr)$ as $n$ tends to infinity;
\item [(ii)]   the following  equality hold:
$ \lim_{n\to\infty}\|S_n\| (\partial_\ver \Tube( B, \bfr))=0.$
\end{itemize}  
 Then,  for  all $r\in (0,\bfr)$ and   $\epsilon\in  (0,r),$   
the vertical boundary term $\Vc_\epsilon(S,r)$ vanishes
and 
\begin{equation*}
  {1\over  (r^2+\epsilon^2)^q} \int_{\Tube(B,r)} S\wedge \beta^q
 =   \int_{\Tube(B,r)} S\wedge \alpha^q_\epsilon .
\end{equation*}
 In particular,  when  
$\supp (S_n)\cap \partial_\ver \Tube( B, \bfr)=\varnothing$ for all $n\geq 1,$ 
then  the above  formula  holds.
\end{theorem}
\proof
We combine  the proofs of Theorem \ref{T:Lelong-Jensen-closed-compact-support} and Theorem \ref{T:Lelong-Jensen-closed-eps}.
\endproof

The remaining    of this  subsection is  devoted to some   estimates of  the  terms  in Lelong-Jensen formulas  when
the  current in question  is a $\Cc^m$-smooth form.

\begin{lemma}\label{L:vertical-boundary-terms}
Let $\bfr\in\R^+_*$ and   
let $S$ be a real current of dimension $2q$ on  a neighborhood of $\overline\Tube(B,\bfr)$ such that $S$ and $\ddc S$ are of order  $0.$ Suppose that $q\leq k-l.$
\begin{enumerate}
\item    Assume that $S$ is  continuous  near $\overline B$ in $\E.$
\begin{itemize}
\item[(1a)]  Then the following  limit $$\lim\limits_{r\to 0+}  {1\over  r^{2(k-l)}}\int_{\Tube(B,r)}  S\wedge \beta^{q}$$ exists
and is  finite. If  moreover $S(y)$ is  a  positive form for all $y\in \overline B\subset \E,$   then  this limit is nonnegative. 

\item[(1b)] If moreover $S$ is  of class $\Cc^1$   near $\overline B$ in $\E,$
then we have the following asymptotic estimate
$$ {1\over  r^{2(k-l)}}\int_{\Tube(B,r)}  S\wedge \beta^{q}=O(r)+\lim\limits_{s\to 0+}  {1\over  s^{2(k-l)}}\int_{\Tube(B,s)}  S\wedge \beta^{q}.$$

\item[(1c)]
If $S$ is  of class $\Cc^2$   near $\overline B$ in $\E,$
then the following  limit $$\lim\limits_{r\to 0+}  {1\over  r^{2(k-l)}}\int_{\Tube(B,r)}  \ddc S\wedge \beta^{q-1}$$ exists
and is  finite.
\end{itemize}

\item  If $S$ is  continuous  near $\overline B$ in $\E,$
then  there is a constant $c>0$ depending only on $S$  such that  for $ 0<r\ll \bfr,$ 
\begin{equation*}
 \big|\int_{\Tube(B,r)} S\wedge \alpha^q\big|\leq  cr^{\max(2(k-l-q),1)}.
 \end{equation*}
 If  If $S'$ is  continuous $(2q-1)$-form near $\overline B$ in $\E,$
then  there is a constant $c>0$ depending only on $S'$  such that  for $ 0<r\ll \bfr,$ 
\begin{equation*}
 \big|\int_{\partial\Tube(B,r)} S'\wedge \alpha^q\big|\leq  cr^{\max(2(k-l-q),1)}.
 \end{equation*}
 \item  If  $S$ is continuous near $\partial_\ver \Tube( B, \bfr),$  then the following  limit
 $$\lim\limits_{r\to 0+} {1\over  r^{2(k-l)}}\int_{\partial_\ver \Tube(B,r)}\dc\varphi\wedge S\wedge \beta^{q-1}$$
 exists
and is  finite.
 If  $S$ is of class $\Cc^1$ near $\partial_\ver \Tube( B, \bfr),$  then  the following  limit
 \begin{equation*}
\lim\limits_{r\to 0+}{1\over  r^{2(k-l)}}\int_{\partial_\ver\Tube(B,r)} \dc S\wedge \beta^{q-1}
\end{equation*}
exists
and is  finite.
\item   If  $S$ is continuous near $\partial_\ver \Tube( B, \bfr),$  then 
\begin{equation*}\big|\int_{\partial_\ver \Tube(B,r)}\dc\log\varphi\wedge S\wedge \alpha^{q-1}\big|\leq  cr^{2(k-l-q)+1}.
 \end{equation*}

\end{enumerate}
\end{lemma}
\proof
\noindent {\bf Assertion (1).} Using  the  partition  of unity $(\theta_\ell)$ introduced in Section  \ref{S:Regularization}, we may   suppose without loss  of generality that
$\pi(\supp(S))$ is  compactly supported in a small open neighborhood 
$V(x_0)$ in $V$ of  a given point $x_0\in V,$  where in a local chart $V(x_0)\simeq \D^l$ and  $\E|_{V(x_0)}\simeq \C^{k-l}\times \D^l.$
For $y\in \E|_{V(x_0)},$ write $y=(z,w).$ 
Consider the function  $R$ given  by
$$S\wedge \beta^q=R(z,w)(i^{l}dw\wedge  d\bar w)\wedge (i^{k-l}dz\wedge d\bar z).$$
Let $\pi:\ \C^{k-l}\setminus \{0\}\to  \P^{k-l-1},$  $z\mapsto  \pi(z):=[z]$   be the canonical projection.
Let $\omega_\FS$ be the  Fubini-Study  form  on  $\P^{k-l-1}.$
There is a smooth function  $h:\  \D^l\times \P^{k-l-1}\to (0,\infty)$ such that
$$
\varphi(z,w)=h([z],w)^2\|z\|^2\qquad\text{for}\qquad  z\in\C^{k-l}\setminus  \{0\},\ w\in \D^l.
$$
We have 
\begin{eqnarray*}
& &\lim\limits_{r\to 0+}  {1\over  r^{2(k-l)}}\int_{\Tube(B,r)}  S\wedge \beta^{q}=\lim\limits_{r\to 0+}  {1\over  r^{2(k-l)}}\int_{(z,w):\ h(w,[z])\|z\|<r}  R(z,w)(i^{l}dw\wedge  d\bar w)\wedge (i^{k-l}dz\wedge d\bar z)\\
&=&\lim\limits_{r\to 0+}  {1\over  r^{2(k-l)}}\int_{w\in\D^l} \big(r^{2(k-l-1)}\int_{ \P^{k-l-1}}  {\pi r^2\over h(w,[z])^2}  \omega^{k-l-1}_\FS([z])\big)  R(0,w)i^{l}dw\wedge  d\bar w\\
&=&  \int_{w\in\D^l} \big(\int_{\P^{k-l-1}}  {\pi\over h(w,[z])^2}  \omega^{k-l-1}_\FS([z])\big)  R(0,w)i^{l}dw\wedge  d\bar w.
\end{eqnarray*}
Hence, the desired limit  exists and is finite.

 Consider  the case  where    $S(y)$  is a positive form  for all $y\in\overline B.$ 
 By 
Lemma \ref{L:hat-alpha-beta} (2) below
  there is  a  constant $c_1>0$ large enough such that 
$
  \hat\beta:=c_1\varphi\cdot  \pi^*\omega+\beta
 $
is  positive  on $\pi^{-1}(\overline B)$ and is strictly positive on $\pi^{-1}(\overline B)\setminus \overline B.$
Consider the function $\widehat R$ given by
 $$ S\wedge \hat\beta^q=\widehat R(z,w)(i^{l}dw\wedge  d\bar w)\wedge (i^{k-l}dz\wedge d\bar z).$$
 So $\widehat R(0,w)$ is  non-negative for $w\in\overline B.$  Arguing as in the previous paragraph, we see that
 \begin{equation*}
\lim\limits_{r\to 0+}  {1\over  r^{2(k-l)}}\int_{\Tube(B,r)}  S\wedge \hat\beta^{q}= 
 \int_{w\in\D^l} \big(\int_{\P^{k-l-1}}  {\pi\over h(w,[z])^2}  \omega^{k-l-1}_\FS([z])\big)  \widehat R(0,w)i^{l}dw\wedge  d\bar w.
\end{equation*}
 Since the last  double integral is nonnegative, we  infer that the  expression on the LHS in the last
 line is also non-negative. Expanding this  expression, we get that
  $$\lim\limits_{r\to 0+}  {1\over  r^{2(k-l)}}\int_{\Tube(B,r)}  S\wedge \beta^{q}+\sum\limits_{j=1}^q {q\choose j} c_1^j\lim\limits_{r\to 0+}  {1\over  r^{2(k-l)}}\int_{\Tube(B,r)} \varphi^j S\wedge\pi^*\omega^j\wedge  \beta^{q-j} \geq 0.$$
 On the other hand,  since $\varphi\lesssim r^2$ on $\Tube(B,r),$   we deduce  from  the first part of assertion (1a)  that
 \begin{equation*}
    \big|{1\over  r^{2(k-l)}}\int_{\Tube(B,r)} \varphi^j S\wedge\pi^*\omega^j\wedge  \beta^{q-j}\big|\lesssim r^{2j}\qquad\text{for}\qquad j\geq 1.
  \end{equation*}
So 
 all terms in the above sum are zero. This implies that the  limit in front of  the  above  sum    is  nonnegative.
The proof of  assertion (1a) is thereby completed.

When  $S$ is  of class $\Cc^1$ near $\overline B$ in $\E,$  we see that $R(z,w)=R(0,w)+O(\|z\|).$
Arguing as in the proof of  assertion (1a) and using this  asymptotic  expression instead of the limit process, assertion (1b) follows.

 Assertion (1c) can be proved  in the  same way as that of assertion (1a) replacing
$S$  by $\ddc S.$

\noindent {\bf Assertion (2).} 
We will only give the proof of the first  inequality since the second one can be shown  similarly using the proof of assertion (3) below. 
We  need  some  estimates which will be  established in Section \ref{S:Basic-forms-and-CV-test}. By
\eqref{e:varphi-new-exp} 
there is a smooth function  $A:\  \D^l\to \GL(\C,k-l)$ such that
\begin{equation*}
\varphi(z,w)= \| A(w)z\|^2\qquad\text{for}\qquad  z\in\C^{k-l},\ w\in \D^l.
\end{equation*}
We see easily that there is a constant $c>1$  such that
\begin{equation*}
c^{-1}\leq \| A(w)\|\leq c,\ w\in \D^l.
\end{equation*}
This,  combined with 
 the  second  equality of \eqref{e:alpha-beta-local},
implies  that
\begin{equation*}
I:=\big|\int_{\Tube(B,r)} S\wedge \alpha^q\big|\lesssim \sum_{j=0}^q  \int_{(z,w)\in\C^{k-l}\times \D^l:\ \|z\|<r}  \|z\|^{-(q-j)}(\ddc\|w\|^2)^l\wedge  (\ddc\|z\|^2)^{k-l-j}\wedge \omega^{j}_\FS([z]).
\end{equation*}
Recall  from   the hypothesis  that  $q\geq k-l.$
Since  $\omega_\FS^{k-l}([z])=0$   and $\omega_\FS([z])\lesssim \|z\|^{-2}  (\ddc\|z\|^2),$  we  see that
\begin{eqnarray*}
I&\lesssim &\sum_{j=0}^{\min (q,k-l-1)} \int_{z\in\C^{k-l}:\ \|z\|<r}\|z\|^{-(q-j)}  (\ddc\|z\|^2)^{k-l-j}\wedge \omega^{j}_\FS([z])\\
&\lesssim& \sum_{j=0}^{\min (q,k-l-1)} \int_{z\in\C^{k-l}:\ \|z\|<r}\|z\|^{-q-j}  (\ddc\|z\|^2)^{k-l} \lesssim  \int_{z\in\C^{k-l}:\ \|z\|<r}\|z\|^{-\min(2q,2k-2l-1)}  (\ddc\|z\|^2)^{k-l}.
\end{eqnarray*}
This  proves assertion (2).

\noindent {\bf Assertion (3).}  Observe that near $x_0,$ 
$$\partial_\ver\Tube(B,r)=\left\lbrace (z,w)\in\C^{k-l}\times \D^l:\ \|z\|<r\quad\text{and}\quad w\in\partial B\right\rbrace.$$ For $y\in \partial_\ver\Tube(B,r)\cap \E|_{V(x_0)},$ write $y=(z,w).$ 
Write  $$\dc\phi\wedge S\wedge \beta^{q-1}=R(z,w)  d\sigma(w)\wedge (i^{k-l}dz\wedge d\bar z),$$
where $d\sigma(w)$ is  the  volume form on $\partial B\cap V(x_0).$
Using this, we argue as in the  proof of assertion (1).  Hence,  the first  limit of assertion (3)  follows.
The  second  one can be proved  similarly.

\noindent {\bf Assertion (4).}  It follows from  
\eqref{e:varphi-new-exp} that there is a constant $c>1$  such that
\begin{equation*}
\dc\log\varphi(z,w)=\sum  O(\|z\|^{-1})  dz_p + O(\|z\|^{-1})  d\bar z_{p'}+  O(1)dw_j+O(1)d\bar w_{j'},
\end{equation*}
the sum being taken over all $1\leq p,p'\leq k-l$ and $1\leq j,j'\leq l.$
Using this inequality   and arguing as in  the proof of assertion (2),  we  see that
\begin{equation*}
I:=\big|\int_{\partial_\ver \Tube(B,r)} \dc\log\varphi\wedge S\wedge \alpha^{q-1}\big|\lesssim \sum_{j=0}^{q-1}  \int_{(z,w)\in\C^{k-l}\times \partial B:\ \|z\|<r}  \|z\|^{-(q-j)} d\sigma(w)\wedge  (\ddc\|z\|^2)^{k-l-j}\wedge \omega^{j}_\FS([z]).
\end{equation*}
Since   $q\leq k-l$  and $\omega_\FS\leq \|z\|^{-2}  (\ddc\|z\|^2),$  we  see that
\begin{eqnarray*}
I&\lesssim &\sum_{j=0}^{q-1} \int_{z\in\C^{k-l}:\ \|z\|<r}\|z\|^{-(q-j)}  (\ddc\|z\|^2)^{k-l-j}\wedge \omega^{j}_\FS([z])\\
&\lesssim& \sum_{j=0}^{q-1} \int_{z\in\C^{k-l}:\ \|z\|<r}\|z\|^{-q-j}  (\ddc\|z\|^2)^{k-l} \lesssim  \int_{z\in\C^{k-l}:\ \|z\|<r}\|z\|^{-2q+1}  (\ddc\|z\|^2)^{k-l}.
\end{eqnarray*}
This  proves assertion (4). 

\endproof

Here is a version of   Theorem  \ref{T:Lelong-Jensen} for smooth  forms when the  minor radius $r_1$  becomes  infinitesimally small.
\begin{theorem}\label{T:Lelong-Jensen-smooth}
Let $\bfr\in\R^+_*$  and let $S$ be a $\Cc^2$-smooth  form of dimension $2q$ on  a neighborhood of $\overline\Tube(B,\bfr).$
Suppose that $q\leq k-l.$
\begin{enumerate}
\item
Then,  for  all $0<r\leq\bfr,$
\begin{equation}\label{e:Lelong-Jensen-smooth}\begin{split}
& {1\over  r^{2q}} \int_{\Tube(B,r)} S\wedge \beta^q -\lim_{s\to 0+} {1\over  s^{2q}} \int_{\Tube(B,s)} S\wedge \beta^q
 =   \int_{\Tube(B,r)} S\wedge \alpha^q\\
 &+  \int_{0}^{r} \big( {1\over t^{2q}}-{1\over r^{2q}}  \big)2tdt\int_{\Tube(B,t)} \ddc S\wedge \beta^{q-1} 
 + \Vc(S,r).
 \end{split}
\end{equation} 
Here, the vertical boundary term  $\Vc(S,r)$ is  given by 
\begin{equation}\label{e:vertical-boundary-term-bis}
\begin{split}
 \Vc(S,r)&:=-\int_{0}^{r} \big( {1\over t^{2q}}-{1\over r^{2q}}  \big)2tdt\int_{\partial_\ver\Tube(B,t)} \dc S^\sharp\wedge \beta^{q-1} 
 -\int_{\partial_\ver \Tube(B,r)}\dc\log\varphi\wedge S^\sharp\wedge \alpha^{q-1}\\
 &+\big( {1\over r^{2q}} \int_{\partial_\ver \Tube(B,r)}\dc\varphi\wedge S^\sharp\wedge \beta^{q-1}-\lim_{s\to 0+} {1\over s^{2q}} \int_{\partial_\ver \Tube(B,s)}\dc\varphi\wedge S^\sharp\wedge \beta^{q-1}\big)
  .
  \end{split}
\end{equation}
\begin{itemize}
\item If   $q<k-l,$ then $\lim_{s\to 0+} {1\over  s^{2q}} \int_{\Tube(B,s)} S\wedge \beta^q=0.$
 \item If   $q=k-l$ and  $S(y)$ is a positive form  for all $y\in \overline B,$   then  the last  limit is nonnegative.
\end{itemize}

\item   Suppose  in addition that  $\supp(S) \cap \partial_\ver \Tube( B, \bfr)=\varnothing.$   
Then,  for  all $0<r<\bfr,$
\begin{eqnarray*}
 &&{1\over  r^{2q}} \int_{\Tube(B,r)} S\wedge \beta^q-\lim_{s\to 0+} {1\over  s^{2q}} \int_{\Tube(B,s)} S\wedge \beta^q\\
& =&\int_{\Tube(B,r)} S\wedge \alpha^q
 + \int_{0}^{r} \big( {1\over t^{2q}}-{1\over r^{2q}}  \big)2tdt\int_{\Tube(B,t)} \ddc S\wedge \beta^{q-1} 
 .
\end{eqnarray*}
\end{enumerate}
\end{theorem}
\proof
Assertion (2) is an immediate consequence of assertion (1).

Assertion (1)  follows from combining  Theorem \ref{T:Lelong-Jensen} for a $\Cc^2$-smooth form $S$ and for  $0<r_1<r_2:=r$ and Lemma \ref{L:vertical-boundary-terms}  for $r:=r_1.$ Indeed, in formulas \eqref{e:Lelong-Jensen} and 
\eqref{e:vertical-boundary-term} we apply  Lemma \ref{L:vertical-boundary-terms} when  $r_1$ tends to $0.$ 
\endproof
 
Here is a version of   Theorem  \ref{T:Lelong-Jensen-closed} for smooth closed forms when the  minor radius $r_1$  becomes  infinitesimally small.
\begin{theorem}\label{T:Lelong-Jensen-smooth-closed}
Let $\bfr\in\R^+_*$  and let $S$ be a $\Cc^1$-smooth closed form of dimension $2q$ on  a neighborhood of $\overline\Tube(B,\bfr).$
Suppose that $q\leq k-l.$
\begin{enumerate}
\item
Then,  for  all $0<r\leq\bfr,$
\begin{equation}\label{e:Lelong-Jensen-smooth-closed}
 {1\over  r^{2q}} \int_{\Tube(B,r)} S\wedge \beta^q- \lim_{s\to 0+}{1\over  s^{2q}} \int_{\Tube(B,s)} S\wedge \beta^q
 = \Vc(S,r)+   \int_{\Tube(B,r)} S\wedge \alpha^q.
\end{equation} 
Here the vertical boundary term  $\Vc(S,r)$ is  given by
\begin{equation}\label{e:vertical-boundary-term-closed-bis}
\begin{split}
 \Vc(S,r)&:=\big(
 {1\over r^{2q}} \int_{\partial_\ver \Tube(B,r)}\dc\varphi\wedge S\wedge \beta^{q-1}-\lim_{s\to 0+}{1\over s^{2q}} \int_{\partial_\ver \Tube(B,s)}\dc\varphi\wedge S\wedge \beta^{q-1}\big)\\
  &-\int_{\partial_\ver \Tube(B,r)}\dc\log\varphi\wedge S\wedge \alpha^{q-1}.
  \end{split}
\end{equation}
\begin{itemize}
\item If  $q<k-l$ then $\lim_{s\to 0+} {1\over  s^{2q}} \int_{\Tube(B,s)} S\wedge \beta^q=0.$
\item If  $q=k-l$ and  $S(y)$  is a  positive  form for all $y\in \overline B,$   then  the last  limit is nonnegative.
\end{itemize}
\item  Suppose  in addition that  $\supp(S) \cap \partial_\ver \Tube( B, \bfr)=\varnothing.$   
Then,  for  all $0<r<\bfr,$
\begin{equation*}
 {1\over  r^{2q}} \int_{\Tube(B,r)} S\wedge \beta^q- \lim_{s\to 0+}{1\over  s^{2q}} \int_{\Tube(B,s)} S\wedge \beta^q
 =   \int_{\Tube(B,r)} S\wedge \alpha^q.
\end{equation*} 
\end{enumerate}
\end{theorem}
\proof
Assertion (2) is an immediate consequence of assertion (1).

Assertion (1)  follows from combining  Theorem \ref{T:Lelong-Jensen-closed} for a $\Cc^1$-smooth form $S$ and for  $0<r_1<r_2:=r$ and Lemma \ref{L:vertical-boundary-terms}  for $r:=r_1.$ Indeed, in formulas \eqref{e:Lelong-Jensen-closed} and \eqref{e:vertical-boundary-term-closed} we apply  Lemma \ref{L:vertical-boundary-terms} when  $r_1$ tends to $0.$ 
\endproof

Finally,  we conclude the section with two  asymptotic   Lelong--Jensen formulas.

\begin{theorem}\label{T:vertical-boundary-terms}
Let $\bfr\in\R^+_*$ and   $0\leq q\leq k-l.$ Let $S$ be a real current of dimension $2q$ on  a neighborhood of $\overline\Tube(B,\bfr)$ such that $S$ and $\ddc S$ are of order  $0.$
Suppose that there is a sequence of $\Cc^2$-smooth forms of dimension $2q$ $(S_n)_{n=1}^\infty$ defined on  a neighborhood of $\overline\Tube(B,\bfr)$ 
such  that
\begin{itemize}
 \item[(i)] 
$S_n$ converge to $S$  in the sense of quasi-positive currents   on a neighborhood of  $\overline\Tube(B,\bfr)$ as $n$ tends to infinity (see Definition \ref{D:quasi-positivity});
\item [(ii)]  $\ddc S_n$ converge to $\ddc S$  in the sense of quasi-positive currents   on a neighborhood of  $\overline\Tube(B,\bfr)$ as $n$ tends to infinity;
\item[(iii)] there is an  open neighborhood of $\partial_\ver\Tube(B,\bfr)$  on which the $\Cc^1$-norms of  $S_n$  are uniformly bounded.
\end{itemize}
Then,  for  all $s,r\in  (0,\bfr]$ with  $s<r$  except for  a  countable set of   values,   formula \eqref{e:Lelong-Jensen} for $r_1:=s,$ $r_2:=r$
(resp.  formula  \eqref{e:Lelong-Jensen-smooth})
holds with 
\begin{equation*}
 |\Vc(S,s,r)|\leq cr\qquad\big(\text{resp.}\qquad  |\Vc(S,r)|\leq  cr\quad\big),
\end{equation*}
where $c$ is a constant  independent of $s,r.$
\end{theorem}
\proof
Combining  Theorem \ref{T:Lelong-Jensen} (resp.  Theorem \ref{T:Lelong-Jensen-smooth})
and Lemma \ref{L:vertical-boundary-terms}, the result follows.
\endproof

\begin{theorem}\label{T:vertical-boundary-closed}
Let $\bfr\in\R^+_*$ and   $0\leq q\leq k-l.$  Let $S$ be a real closed current of dimension $2q$ on  a neighborhood of $\overline\Tube(B,\bfr).$ 
Suppose that there is a sequence of $\Cc^1$-smooth closed forms of dimension $2q:$  $(S_n)_{n=1}^\infty$ defined on  a neighborhood of $\overline\Tube(B,\bfr)$ 
such  that
\begin{itemize}
 \item[(i)] 
$S_n$ converge to $S$  in the sense of quasi-positive currents   on a neighborhood of  $\overline\Tube(B,\bfr)$ as $n$ tends to infinity (see Definition \ref{D:quasi-positivity});
\item[(ii)] there is an  open neighborhood of $\partial_\ver\Tube(B,\bfr)$  on which the $\Cc^m$-norms of  $S_n$  are uniformly bounded, where
$m=0$ if  $q<k-l$ and $m=1$ if $q=k-l.$ 
\end{itemize}
Then,  for  all $s,r\in  (0,\bfr]$ with  $s<r$  except for  a  countable set of   values,  
formula \eqref{e:Lelong-Jensen-closed}  for $r_1:=s,$ $r_2:=r$  (resp. formula \eqref{e:Lelong-Jensen-smooth-closed}) holds
  with 
\begin{equation*}
 |\Vc(S,s,r)|\leq cr\qquad\big(\text{resp.}\qquad  |\Vc(S,r)|\leq  cr\quad\big),
\end{equation*}
where $c$ is a constant independent of $s,r.$
\end{theorem}
\proof
Combining  Theorem \ref{T:Lelong-Jensen-closed} (resp.  Theorem \ref{T:Lelong-Jensen-smooth-closed})
and Lemma \ref{L:vertical-boundary-terms}, the result follows.
\endproof

\begin{theorem}\label{T:vertical-boundary-eps}
 We keep the hypothesis and the notation of Theorem \ref{T:Lelong-Jensen-eps} (resp.  Theorem \ref{T:Lelong-Jensen-closed-eps}).
 Then there  is  a  constant $c$ depending only on $S$ such that 
 for all  $r\in (0,\bfr]$  and $\epsilon\in(0,r),$ the  following assertions  hold:
 \begin{enumerate}
  \item If $q<k-l,$ then 
 $
 |\Vc_\epsilon(S,r)|\leq  cr.$
 
 \item If $  q=k-l$ and  we  are in the assumption of Theorem \ref{T:Lelong-Jensen-eps}, then
 
 $$ \big|\Vc_\epsilon(S,r)-{
1\over r^{2q}} \int_{\partial_\ver \Tube(B,r)}\dc\varphi\wedge S^\sharp\wedge \beta^{q-1}\big|
\leq cr 
.$$
\item If $  q=k-l$ and  we  are in the assumption of Theorem \ref{T:Lelong-Jensen-closed-eps}, then
 
 $$ \big|\Vc_\epsilon(S,r)-{
1\over r^{2q}} \int_{\partial_\ver \Tube(B,r)}\dc\varphi\wedge S\wedge \beta^{q-1}\big|
\leq cr 
.$$

 \end{enumerate}

\end{theorem}
\proof
Combining  Theorem \ref{T:Lelong-Jensen-eps} (resp.  Theorem \ref{T:Lelong-Jensen-closed-eps})
and Lemma \ref{L:vertical-boundary-terms}, the result follows.
\endproof
\section{Forms   $\hat\alpha,$ $\alpha_\ver$ and   $\hat\beta,$    $\beta_\ver$  and a  convergence test}
\label{S:Basic-forms-and-CV-test}
 In this  section we introduce  some basic  forms  for  the    bundle $\E$ which is  the normal bundle to  $V$ in $X.$
 We also prove  a  convergence test. They will be  used  throughout   this  work.
 We keep the Standing Hypothesis introduced in Subsection \ref{SS:Global-setting}.

\subsection{Forms $\alpha_\ver$ and $\beta_\ver$}
 
 Since  the  transition  functions  of the holomorphic vector bundle $\E$ are holomorphic, 
 the vertical operators $\partial_\ver,$ $\dbar_v$ 
which are the  restrictions of the usual operators $\partial$ and $\dbar$ on fibers of $\E$  are well-defined. More precisely, for a smooth
form $\Phi$ on an open set $\Omega$ in $\E,$ we can define
\begin{equation}\label{e:d-and-dbar-ver}
 \partial_\ver\Phi(y):=\partial |_{\E_{\pi(y)}}\Phi (y) \qquad\text{and}\qquad \dbar_\ver\Phi(y):=\dbar|_{\E_{\pi(y)}}\Phi(y)\qquad\text{for}\qquad y\in \Omega.
\end{equation}
So  the vertical operators  $d_\ver$ and  $\ddcv$ are  also well-defined  by the formulas \begin{equation}
\label{e:ddc-ver}
d_\ver \Phi:=\partial_\ver \Phi+\dbar_\ver \Phi\quad\text{and}\quad \ddcv\Phi:={i\over \pi}  \partial_\ver\dbar_\ver\Phi.                           
                          \end{equation}  
Consider for $ y\in \E,$
\begin{equation}\label{e:alpha-beta-ver}
 \alpha_\ver(y):=\ddcv \log\varphi (y)=\ddc|_{\E_{\pi(y)}}\log\varphi (y) \qquad\text{and}\qquad \beta_\ver(y):=\ddcv \varphi(y)=\ddc|_{\E_{\pi(y)}}\varphi(y),
\end{equation}
where  $\ddc|_{\E_{\pi(y)}}$  is restriction  of  the operator $\ddc$ on the fiber $\E_{\pi(y)}.$  
Observe  that  both
 $\alpha_\ver$  and $\beta_\ver$  are positive  $(1,1)$-forms on $\E.$
However, they are not necessarily closed.

\subsection{Analysis  in local coordinates}
Since $V_0\Subset  V,$ we only need to  prove  a local result near a given point $y_0\in V_0.$
We use the coordinates $(z,w)\in\C^{k-l}\times \C^l$  around a neighborhood $U$ of $y_0$  such that $y_0=0$  in  these coordinates. 
We may assume that  $U$ has the form $U=U'\times U'',$ where $U'$ (resp. $U'')$ are open neighborhood of $0'$ in $\C^{k-l}$ of  ($0''$ in $\C^l$)
and  $V=\{z=0\}\simeq U''.$ Moreover,  we  may assume  that $U''=(2\D)^l.$  
Consider  the trivial  vector bundle $\pi:\ \E \to  U''$ with  $\E\simeq  \C^{k-l}\times U''.$ 
Let $\pi_\FS:\ \C^{k-l}\setminus \{0\}\to  \P^{k-l-1},$  $z\mapsto  \pi_\FS(z):=[z]$   be the canonical projection.
Let $\omega_\FS$ be the  Fubini-Study  form  on  $\P^{k-l-1}.$  So
\begin{equation}\label{e:FS}
\pi^*_\FS  (\omega_\FS([z]))=\ddc (\log{\|z\|^2})\qquad\text{for}\qquad  z\in\C^{k-l}\setminus \{0\}.
\end{equation}
There is a smooth function  $A:\  \D^l\to \GL(\C,k-l)$ such that
\begin{equation}\label{e:varphi-new-exp}
\varphi(z,w)= \| A(w)z\|^2\qquad\text{for}\qquad  z\in\C^{k-l},\ w\in \D^l.
\end{equation}
It follows  from \eqref{e:alpha-beta-ver} and  \eqref{e:varphi-new-exp} that  
\begin{equation}\label{e:tilde-alpha-beta-local-exp}
\alpha_\ver(z,w)=  A(w)^* [\ddc \log{\|z\|^2}]\quad\text{and}\quad    \beta_\ver(z,w)=  A(w)^* [\ddc \|z\|^2]\quad\text{for}\quad  z\in\C^{k-l},\ w\in \D^l.
\end{equation}
We place ourselves on  an open set of $\C^{k-l}$ defined by $z_{k-l}\not=0.$
We   may assume without loss of generality that
\begin{equation}\label{e:max-coordinate} 2|z_{k-l}| > \max\limits_{1\leq j\leq k-l}|z_j|.
\end{equation}
and use the projective coordinates
\begin{equation}\label{e:homogeneous-coordinates}
\zeta_1:={z_1\over z_{k-l}},\ldots, \zeta_{k-l-1}:={z_{k-l-1}\over z_{k-l}},\quad \zeta_{k-l}=z_{k-l}.
\end{equation}
In the coordinates  $\zeta=(\zeta_1,\ldots,\zeta_{k-l})=(\zeta',\zeta_{k-l}),$ the form $\omega_\FS([z])$  can be  rewritten as  
\begin{equation}\label{e:FS-zeta} \omega_\FS([z])= \ddc \log{ (1+|\zeta_1|^2+\cdots+|\zeta_{k-l-1}|^2)},
\end{equation}
and a direct computation shows that 
\begin{equation}\label{e:omega'-zeta'}
\omega_\FS([z])\approx  (1+\|\zeta'\|^2)^{-2}\omega'(\zeta'),\quad\text{where}\quad \omega'(\zeta'):=\ddc (|\zeta_1|^2+\cdots+|\zeta_{k-l-1}|^2).
\end{equation}
Since  $|\zeta_j|<2$ for $1\leq j\leq k-l-1$ by  \eqref{e:max-coordinate}, it follows from  \eqref{e:omega'-zeta'} and  the first equality in \eqref{e:tilde-alpha-beta-local-exp} that
\begin{equation}\label{e:omega_FS-vs-omega'}
\omega_\FS([z])\approx \omega'(\zeta')\approx \alpha_\ver.
\end{equation}

\subsection{Forms $\hat\alpha$ and $\hat\beta$}
The next result  shows  that $\alpha$ and $\beta$ are,  in some sense,  nearly  positive forms
on $\pi^{-1}(V_0)\subset \E.$
Namely, their following  variants  $\hat\alpha,$ $\hat\alpha'$  and $\hat\beta$  are positive.
This positivity plays  a  crucial  role    in the sequel.

\begin{lemma}\label{L:hat-alpha-beta}
\begin{enumerate}
\item  We have the following  expressions
\begin{equation}\label{e:alpha-beta-local}
 \begin{split}
 \beta(z,w)&=   A^*(w)(\sum_{p=1}^{k-l} idz_p\wedge d\bar z_{p})+ \sum O(\|z\|) dz_p\wedge d\bar w_{q'} +   O(\|z\|) d\bar z_{p'}\wedge d w_{q} + O(\|z\|^2) dw_q\wedge d\bar w_{q'},\\
 \alpha(z,w)&=    A(w)^* [\ddc \log{\|z\|^2}]+ \sum O(\|z\|^{-1}) dz_p\wedge d\bar w_{q'} +   O(\|z\|^{-1}) d\bar z_{p'}\wedge d w_{q} + O(1) dw_q\wedge d\bar w_{q'}. 
 \end{split}
 \end{equation}
 Here, in the  first sum $A(w)$ is  regarded as a $\C$-linear endomorphism of $\C^{k-l},$ and   the other sums are taken over   $ 1\leq p,p'\leq k-l$ and $1\leq q,q'\leq l.$
\item  There is  a  constant $c_1>0$ large enough such that 
\begin{equation}\label{e:hat-beta}
  \hat\beta:=c_1\varphi\cdot  \pi^*\omega+\beta
 \end{equation}
is  positive  on $\pi^{-1}(V_0)$ and is strictly positive on $\pi^{-1}(V_0)\setminus V_0,$
and 
\begin{equation}\label{e:hat-alpha'}
\hat\alpha':= c_1 \pi^*\omega+ \alpha
\end{equation}
satisfies 
\begin{equation}\label{e:hat-alpha'-vs-alpha_ver}
\hat\alpha'\geq  c^{-1}_1\alpha_\ver.
\end{equation}
In particular,   $\hat\alpha'$ is positive  on $\pi^{-1}(V_0).$
\item 
  For  every $r>0,$ there are constants $c_2,c_3>0$ such that on $\Tube(V_0,r)\setminus V_0,$  
\begin{equation}\label{e:hat-alpha}
\hat\alpha:=\hat\alpha'+c_2\beta= c_1 \pi^*\omega+ \alpha+c_2\beta
\end{equation}
is   strictly positive, and 
\begin{equation}\label{e:hat-alpha-vs-alpha_ver}
\hat\alpha\geq  c^{-1}_1\alpha_\ver,
\end{equation}
and
\begin{equation}\label{e:hat-alpha-vs-hat-beta}
\varphi\hat\alpha\leq  c_3\hat\beta.
\end{equation}
\item 
  For  every $r>0,$ there are constants $c_3>0$ such that  on $\Tube(V_0,r),$
  \begin{equation}\label{e:hat-beta-vs-beta_ver}
\hat\beta\geq  c^{-1}_1\beta_\ver,
\end{equation}
  and on $\Tube(V_0,r)\setminus V_0,$  
\begin{equation}\label{e:tilde-alpha-vs-hat-beta}
\varphi\alpha_\ver\leq c_3\hat\beta.
\end{equation}
\end{enumerate}

\end{lemma}
\proof

\noindent{\bf  Proof of \eqref{e:alpha-beta-local}.} Its  proof  follows  from formulas  \eqref{e:alpha-beta-spec} and  expression \eqref{e:varphi-new-exp}.

\noindent{\bf  Proof of \eqref{e:hat-beta}.}  By  expression \eqref{e:varphi-new-exp}, there is a constant $c>1$ such that 
$$
\varphi(z,w)\geq c\|z\|^2 \quad\text{and}\quad A^*(w)(\sum_{p=1}^{k-l} idz_p\wedge d\bar z_{p})\geq   c(\sum_{p=1}^{k-l} idz_p\wedge d\bar z_{p})\quad\text{and}\quad  \omega(w)\geq c (\sum_{q=1}^{l} idw_q\wedge d\bar w_q).
$$
Using this and  the first  inequality of assertion  (1),  we get that for $c_1>0,$
\begin{eqnarray*}
  \hat\beta&=&c_1\varphi\cdot  \pi^*\omega+\beta \geq  c^2c_1\|z\|^2 (\sum_{q=1}^{l} idw_q\wedge d\bar w_q)
  +c(\sum_{p=1}^{k-l} idz_p\wedge d\bar z_{p})\\
  &+&\sum_{p,q'} \varphi_{p,q'}(z,w) dz_p\wedge d\bar w_{q'} +\sum_{p',q}   \varphi_{p',q}(z,w) d\bar z_{p'}\wedge d w_{q} +f(z,w) \sum_{q=1}^{l} idw_q\wedge d\bar w_q,
 \end{eqnarray*}
where $\varphi_{p,q'},$ $\varphi_{p',q}$  are  complex-valued functions and $f(z,w)$ is a real-valued function such that    $|\varphi_{p,q'}(z,w)|\leq    c'\|z\|$ and $|\varphi_{p',q}(z,w)|\leq    c'\|z\|$  and  $f(z,w)\geq  -c'$  for some  constant $c'>0.$
By  Cauchy--Schwarz inequality,  we see that for $c_1>0$ large  enough, more  precisely,
when  $(c_1c^2-c' )c>4l(k-l)c'^2,$   $\hat\beta(z,w)$ is positive for all $(z,w)\in \Tube(V_0,r)$ and is  strictly positive  outside $z=0.$

\noindent{\bf  Proof of \eqref{e:hat-alpha'-vs-alpha_ver}.} 
We use the  homogeneous coordinates introduced  in \eqref{e:homogeneous-coordinates}.
We infer from \eqref{e:alpha-beta-spec} and  \eqref{e:varphi-new-exp} that
\begin{eqnarray*}
 \alpha&=&\ddcwz \log{\| A(w)z)\|^2}=\ddcwzetap\log{\| A(w)(\zeta',1)\|^2}\\
 &=&\ddczetap\log{\| A(w)(\zeta',1)\|^2}+\ddcw\log{\| A(w)(\zeta',1)\|^2}+\partial_w\dbar_{\zeta'}\log{\| A(w)(\zeta',1)\|^2}\\
 &+&\dbar_w\partial_{\zeta'}\log{\| A(w)(\zeta',1)\|^2}.
\end{eqnarray*}
Since there are constants $c,c'>0$ such that  $c'\leq\|A(w)\|\leq c,$ we see that
the first term in the last line is  equivalent to  $\ddc_{\zeta'}\log{\| (\zeta',1)\|^2},$ which is
in turn  equivalent to $\omega'(\zeta')$  by \eqref{e:FS-zeta} 
and \eqref{e:omega'-zeta'}.

As $|\zeta_j|<2$ for $1\leq j\leq k-l-1,$  a straightforward computation shows  that the  sum of the other 3 terms
is a smooth differential form $\Phi(\zeta',w)$  with bounded coefficients.
Consequently, by  Cauchy-Schwarz inequality, when $c_1>0$ is  large  enough, we obtain 
$$
\hat\alpha'= c_1 \pi^*\omega+ \alpha\geq  c c_1(\sum_{q=1}^{l} idw_q\wedge d\bar w_q)+c'\omega'(\zeta')+\Phi(\zeta',w)\geq {c'\over2} \omega'(\zeta')\approx 
\alpha_\ver,
$$
where   the last inequality  follows from \eqref{e:omega_FS-vs-omega'}.
This proves \eqref{e:hat-alpha'-vs-alpha_ver}.

\noindent{\bf  Proof of \eqref{e:hat-alpha}.}
Let $c'_1$ be  a constant   which satisfies both \eqref{e:hat-beta} and \eqref{e:hat-alpha'}-\eqref{e:hat-alpha'-vs-alpha_ver}  when  $c_1$ therein
is  replaced by $c'_1.$ 
Let  $c_1:=2c'_1$ and $c_2:= {c'_1\over  r^2},$  and    set $\hat\alpha':= c_1\pi^*\omega+\alpha.$ We have
for $(z,w)\in  \Tube(V_0,r)$  that
\begin{equation*}
\hat\alpha:=\hat\alpha'+c_2\beta= 2c'_1 \pi^*\omega+ \alpha+c_2\beta=(c'_1 \pi^*\omega+ \alpha)
+{c'_1\over r^2}( r^2\pi^*\omega+\beta)\geq c_2 ( c'_1\varphi\pi^*\omega+\beta).
\end{equation*}
Since we know  by \eqref{e:hat-beta} that  the last form is  strictly positive on $\pi^{-1}(V_0)\setminus V_0,$
the proof of  \eqref{e:hat-alpha} is thereby completed.

\noindent{\bf  Proof of \eqref{e:hat-alpha-vs-alpha_ver}.} It is  similar to the proof of \eqref{e:hat-alpha'-vs-alpha_ver}.

\noindent{\bf  Proof of \eqref{e:hat-alpha-vs-hat-beta}.}
Using   \eqref{e:alpha-beta-local} and  applying  Cauchy-Schwarz inequality,  there are constants $c',c''>0$ such that
$$\varphi\alpha\leq   c'\varphi\pi^*\omega +c' \sum_{p=1}^{k-l} idz_p\wedge d\bar z_p \leq c''\varphi\pi^*\omega+c''\beta.$$
Therefore, for $c_3>0$ large  enough, we  obtain that
\begin{equation*}
 \varphi\hat\alpha=c_1\varphi \pi^*\omega+c_2\varphi\beta+\varphi\alpha\leq (c_1+c'') \varphi \pi^*\omega+(c_2\varphi+c'')\beta\leq  c_3\hat\beta.
\end{equation*}

\noindent{\bf  Proof of \eqref{e:hat-beta-vs-beta_ver}.} 
Using \eqref{e:alpha-beta-local}  and \eqref{e:tilde-alpha-beta-local-exp}  and  applying  Cauchy-Schwarz inequality,
we see that for $c_1>1$ large enough, there is  $c'>0$ such that
$$
\hat \beta =c_1\varphi\cdot \pi^*\omega+\beta\geq  c'\sum_{p=1}^{k-l}idz_p\wedge d\bar z_p\geq  c^{-1}_1\beta_\ver.
$$
\noindent{\bf  Proof of \eqref{e:tilde-alpha-vs-hat-beta}.}   By \eqref{e:tilde-alpha-beta-local-exp} we have for a  large constant $c_3>1$ that
$$
\varphi\alpha_\ver\lesssim  \|z\|^2 \ddc \log{\|z\|^2}\lesssim \sum_{p=1}^{k-l}  idz_p\wedge d\bar z_p\lesssim   c_3\hat\beta.
$$
This completes the proof. \endproof

\subsection{A convergence test}
The  following  elementary  result will be  repeatedly used in this  work.
\begin{lemma}\label{L:elementary}
Let $0<r_1<r_2\leq\bfr.$   Consider  two functions
 $f:\ (0,\bfr]\to \R$  and  $\epsilon:\   [\bfr^{-1},\infty)\to (0,\infty),$ $ \lambda\mapsto\epsilon_\lambda$
  such that
 \begin{itemize}\item[(i)]   there are  two constants $c>0$ and $N\in\N$
 such that  if $2^n\leq  \lambda <2^{n+1}$ and $2^{n-N}>\bfr^{-1},$   then 
 $\epsilon_\lambda\leq c\sum_{j=- N}^N \epsilon_{2^{n+j}};$ 
 \item[(ii)]  $\sum_{n\in\N: 2^n\geq \bfr^{-1}}\epsilon_{2^n}<\infty;$
 \item[(iii)] For $r\in(r_1,r_2), $  we have $
 f({r\over \lambda})-f({r_1\over\lambda})\geq  -\epsilon_\lambda.$  
 \end{itemize}
 \begin{enumerate}
  \item Then we have 
$\lim_{r\to 0} f(r)=\liminf_{r\to 0} f(r)\in\R\cup\{-\infty\}.$

\item If instead of condition (iii) we  have the following stronger  condition (iii'): 
$$ |f({r_2\over \lambda})-f({r_1\over\lambda})|\leq  \epsilon_\lambda,$$
 then $\lim_{r\to 0} f(r)=\liminf_{r\to 0} f(r)\in\R,$ that is, the last limit is  finite.

\end{enumerate}
\end{lemma}
\proof
Set
 $
 \nu:=\liminf\limits_{r\to 0+} f(r)\in \R\cup\{-\infty\}.
 $
 So  there is a decreasing  sequence  $s_n$  such that  $s_n\to 0$ and  $\lim\limits_{n\to\infty} f(s_n)=\nu.$
 Using the  hypothesis (i)-(ii)-(iii) one can show that for $0<r<s_n,$
 \begin{eqnarray*}
 f(r)-f(s_n)&=&f(r)-f(2r)+f(2r)-f(2^2r) +\ldots+  f(2^{M-1}r) -f(2^Mr)+ f(2^Mr)-f(s_n)\\
 &\leq&   
 c \sum_{k=0 }^M\epsilon_{2^{-k} r^{-1}r_2^{-1}},
 \end{eqnarray*}
 where $M$ is  the largest  non-negative integer such that $2^Mr\leq s_n.$ 
 Using the hypothesis (i)-(ii)-(iii) again we see that the sum on the RHS tends to $0$  as $n$ tends to infinity. This proves  assertion (1).
 
 The hypothesis (i)-(ii)-(iii') also  shows that $ |f({r_2\over \lambda})-f({r_2\over 2\lambda})|\leq  c\epsilon_\lambda.$
 So 
 $$
 \sum_{n=0}^\infty |f({r_2\over 2^n})-f({r_2\over 2^{n+1}})|\leq  c\sum_{n=0}^\infty\epsilon_{2^n}<\infty.
 $$
 Hence,  $\liminf_{n\to \infty} f({r_2\over  2^n})\in\R.$  This, combined  with assertion (1), implies  assertion (2).
\endproof

\section{Positive closed currents  and holomorphic admissible maps}\label{S:closed-holomorphic}

In  this  section we deal with positive  closed currents together with holomorphic admissible maps, and   we  prove  Theorem  \ref{T:top-Lelong-closed}  and  then  Theorem \ref{T:Lelong-closed-all-degrees}.
 This section  may be regarded   as  a preparation  for  the  proof of Theorems \ref{T:Lelong-closed}, where    the  general situation  with non-holomorphic  admissible maps  will be  investigated.

We  keep  the  global  setting  of Subsection  \ref{SS:Global-setting}, in particular, the Standing Hypothesis.
We  also suppose   in addition that $T$ is a positive closed on $X,$ $\tau$ is a holomorphic  admissible map,
and   $\omega$ is  a Hermitian  form $\omega$ on $V.$ 

\subsection{Top Lelong  number}
   This   subsection is  devoted  to the proof of Theorem \ref{T:top-Lelong-closed}.

\proof[Proof of  assertion (1) of Theorem \ref{T:top-Lelong-closed}]

Consider  a small neighborhood $V(y_0)$ of  an arbitrary  point $y_0\in  \Tube(B, r_0),$  where in a local chart $V(y_0)\simeq \D^l$ and  $\E|_{V(y_0)}\simeq \C^{k-l}\times \D^l.$
For $y\in \E|_{V(y_0)},$ write $y=(z,w).$ Since   $\upm=\min(l,k-p)$ and $T$ is  of bidegree $(p,p)$ and $\tau$ is  holomorphic,  we  argue as in
the proof of the  Fact  in  Corollary \ref{C:Lelong-Jensen} that  $\tau_*T\wedge \pi^*\omega^\upm$ is  of full bidegree $(l,l)$  in $dw,$ $d\bar w.$
Consequently, since  $\tau$ is  holomorphic and   $T$ is positive closed,    
it follows that
$$
 d ( \tau_*T\wedge \pi^*\omega^\upm)=d(\tau_* T)\wedge \pi^*\omega^\upm= \tau_* (dT)\wedge \pi^*\omega^\upm =  0.
$$
So  $\tau_*T\wedge \pi^*\omega^\upm$  is a  positive  closed  current.
For $0<r_1<r_2<\bfr,$  Theorem    \ref{T:Lelong-Jensen-closed-compact-support} and 
Corollary \ref{C:Lelong-Jensen-closed}
applied  to    this  current 
gives
\begin{equation*}
  \nu_\top(T,B,r_2,\tau)- \nu_\top(T,B,r_1,\tau)=\int_{\Tube(B,r_1,r_2)}\tau_*T\wedge \pi^*(\omega^\upm)\wedge\alpha^{k-p-\upm} =\kappa_\top(T,B,r_1,r_2,\tau).
\end{equation*}
 Hence,  the  identity of assertion (1) follows.

It remains  to show  that $\nu_\top(T,B,r_1,\tau),$  $\nu_\top(T,B,r_2,\tau)$ and $\kappa_\top(T,B,r_1,r_2,\tau)$ are non-negative.
 As previously observed,   $T\wedge \pi^*\omega^\upm$ is  of full bidegree $(l,l)$  in $dw,$ $d\bar w.$
Consequently, we infer from \eqref{e:hat-beta} and  \eqref{e:hat-alpha'} that
\begin{eqnarray*}
\tau_*T\wedge \pi^*(\omega^m)\wedge\alpha^{k-p-\upm}
&=&\tau_*T\wedge \pi^*(\omega^\upm)\wedge(\hat\alpha')^{k-p-\upm},\\
\tau_* T\wedge \pi^*\omega^\upm\wedge\beta^{k-p-\upm}&=&\tau_* T\wedge \pi^*\omega^\upm\wedge\hat\beta^{k-p-\upm}.
\end{eqnarray*}
Therefore, we deduce from  \eqref{e:Lelong-corona-numbers} that
\begin{equation*}
  \kappa_\top(T,B,r_1,r_2,\tau)= \int_{\Tube(B,r_1,r_2)}\tau_*T\wedge \pi^*(\omega^\upm)\wedge(\hat\alpha')^{k-p-\upm}. 
\end{equation*}
Moreover, by \eqref{e:Lelong-numbers} we also get that
\begin{equation*}
 \nu_\top(T,B,r,\tau)=   
 {1\over r^{2(k-p-\upm)}}\int_{\Tube(B,r)} (\tau_*T)\wedge \pi^*(\omega^\upm)
 \wedge \hat\beta^{k-p-\upm}.  
\end{equation*} 
Since $T$ is a  positive  current, and by Lemma \ref{L:hat-alpha-beta} $\omega,$ $\hat\alpha',$  $\hat\beta$ are positive forms,
and by the  hypthesis the map $\tau$   is holomorphic, the RHS of the last two  equations   are  $\geq 0.$
Hence, $\nu_\top(T,B,r,\tau)$  and 
$\kappa_\top(T,B,r_1,r_2,\tau)$ are non-negative.   This, combined  with  the identity of assertion (1), 
 show that $r\mapsto  \nu_\top(T,B,r,\tau)$  is increasing for $r\in(0,\bfr].$ This completes the proof of  assertion (1).
\endproof


\proof[Proof of  assertion (2) of Theorem \ref{T:top-Lelong-closed}]
Since  we know  by assertion (1) that the non-negative function $r\mapsto  \nu_\top(T,B,r,\tau)\geq 0$ is increasing for $r\in(0,\bfr),$
assertion (2) follows.
\endproof


\proof[Proof of  assertion (3) of Theorem \ref{T:top-Lelong-closed}]
By \eqref{e:Lelong-log-bullet-numbers}  and  the identity of assertion (1), we have 
\begin{eqnarray*}\
 0\leq \kappa^\bullet_\top(T,B,r,\tau)= \limsup\limits_{s\to0+}   \kappa_\top(T,B,s,r,\tau)&=& \nu_\top(T,B,r,\tau)-\liminf_{s\to0+} \nu_\top(T,B,s,\tau)\\
 &=&  \nu_\top(T,B,r,\tau)- \nu_\top(T,B,\tau),
\end{eqnarray*} 
where the last equality holds by assertion (2).   Consequently, we  infer  from  assertion  (2) again that
\begin{equation*}
 \lim\limits_{r\to0+}\kappa^\bullet_\top(T,B,r,\tau)= \lim\limits_{r\to0+}\nu_\top(T,B,r,\tau)- \nu_\top(T,B,\tau)=0.
\end{equation*}
\endproof


\proof[Proof of  assertion (4) of Theorem \ref{T:top-Lelong-closed}]
First, 
we  will prove  the interpretation of assertion (4)  in the spirit  of \eqref{e:Lelong-number-point-bisbis(1)}. 
Since $p>0$ and $l<k,$  it follows  from \eqref{e:m}  that $k-p-\upm<k-l.$
Therefore, we are in the position to apply  Theorem  \ref{T:Lelong-Jensen-smooth-closed}  to the case where $q=k-p-\upm<k-l.$ 
Hence, we get that 
$$
\nu_\top( T^\pm_n,B,r,\tau)=\kappa_\top(T^\pm_n,B,r,\tau)+\Vc(\tau_*T^\pm_n\wedge \pi^*(\omega^\upm),r).
$$
On the other hand,  we deduce  from \eqref{e:vertical-boundary-term-closed-bis} and  the fact that   $\tau_*T^\pm_n\wedge \pi^*\omega^\upm$ is  of full bidegree $(l,l)$  in $dw,$ $d\bar w$  that $\Vc( \tau_*T^\pm_n\wedge \pi^*(\omega^\upm),r)=0$  since all  the integrals involved in this  term  are performed over $\partial_\ver \Tube(B,r)$ which is a manifold of real dimension $2l-1$ in $w.$  
Consequently, by Lemma \ref{L:quasi-positivity}, we have,  for all $0<r<\bfr$ except for a  countable set of values,
\begin{eqnarray*} 
\kappa_\top(T,B,r,\tau)&:=& \lim\limits_{n\to\infty} \kappa_\top(T^+_n-T^-_n,B,r,\tau)=\lim\limits_{n\to\infty}\nu_\top(T^+_n,B,r,\tau)-\lim\limits_{n\to\infty}\nu_\top(T^-_n,B,r,\tau)\\
&=&\nu_\top(T^+,B,r,\tau)-\nu_\top(T^-,B,r,\tau)=\nu_\top(T,B,r,\tau).
\end{eqnarray*}
This, combined  with assertion (1),  implies  the  desired interpretation  according to  Definition \ref{D:Lelong-log-numbers(1)}.

Second,  we  will prove  the interpretation of assertion (4)  in the spirit  of \eqref{e:Lelong-number-point-bisbis(2)}.
To start  with,  we  fix $0<r\leq\bfr$ and let  $0<\epsilon<r.$  Theorem \ref{T:Lelong-Jensen-closed-eps} applied  to  $\tau_*T\wedge \pi^*(\omega^\upm)$  gives
 $$
 {1\over  (r^2+\epsilon^2)^{k-p-\upm}} \int_{\Tube(B,r)} \tau_*T\wedge \pi^*(\omega^\upm)\wedge \beta^{k-p-\upm}
 = \Vc_\epsilon( \tau_*T\wedge \pi^*(\omega^\upm),r)+   \int_{\Tube(B,r)} \tau_*T\wedge \pi^*(\omega^\upm)\wedge \alpha^{k-p-\upm}_\epsilon.
$$
Now  we let $\epsilon$ tend to $0.$  Then the LHS  tends to $\nu_\top(T,B,r,\tau).$  On the other hand,  we deduce  from \eqref{e:vertical-boundary-term-closed-eps} and  the fact that   $T\wedge \pi^*\omega^\upm$ is  of full bidegree $(l,l)$  in $dw,$ $d\bar w$  that $\Vc_\epsilon( \tau_*T\wedge \pi^*(\omega^\upm),r)=0.$  Consequently,    the second  term on the RHS  tends to $\nu_\top(T,B,r,\tau)$  as $\epsilon$ tends to $0+.$
Hence,  by assertion (2),  the  desired  interpretation  according to  Definition \ref{D:Lelong-log-numbers(2)} follows.
\endproof




\subsection{Other Lelong numbers}

Introduce  the following  mass indicators, for a  positive  current $T$  of bidegree $(p,p)$ defined on $X$ and  for $0\leq j\leq \upm$ and  for $0<s<r<\bfr,$ 
 \begin{equation}\label{e:mass-indicators}
 \begin{split}
 \hat\nu_j(T,B,r,\tau)&:={1\over  r^{2(k-p-j)}}\int\limits_{(\Tube(B,r))} \tau_*T\wedge\pi^*\omega^j\wedge  (\beta+c_1r^2\pi^*\omega)^{k-p-j},\\
 \hat\kappa^{'\bullet}_j(T,B,r,\tau)&:=\int\limits_{\Tube(B,r) \setminus V  } \tau_*(T)
 \wedge\pi^*\omega^j\wedge  (\hat\alpha')^{k-p-j},\\
 \hat\kappa'_j(T,B, s,r,\tau)&:=\int\limits_{\Tube(B,s,r)} \tau_*(T)
 \wedge\pi^*\omega^j\wedge  (\hat\alpha')^{k-p-j}.
 \end{split}
 \end{equation}
We also  write $\hat\nu_\top(T,r,\tau),$ $\hat\kappa^{'\bullet}_\top(T,r,\tau),$  $\hat\kappa'_\top(T,r,s,\tau)$ instead of $\hat\nu_\upm(T,r,\tau),$ $\hat\kappa^{'\bullet}_\upm(T,r,\tau),$  $\hat\kappa'_\upm(T,r,s,\tau)$ respectively. 

\begin{lemma}\label{L:hat-nu}
For  $0<r<\bfr,$  we have that $\hat\nu_j(T,B,r)\geq 0$ 
and 
 \begin{equation} \label{e:hat-nu-vs-nu}
   \hat\nu_j(T,B,r)=\sum\limits_{q=0}^{k-p-j}  {k-p-j\choose q}c_1^q \nu_{j+q}(T,B,r,\tau).  
 \end{equation}
 \end{lemma}
 \proof
  By   Lemma \ref{L:hat-alpha-beta} (1), $\hat\beta$ and  $\beta+c_1r^2\pi^*\omega$ are smooth positive forms. This, combined with the positivity of the current $T$     and  the  explicit formula of  $\hat\nu_j(T,B,r)$ 
  in  \eqref{e:mass-indicators}, implies  that 
  this real number is  non-negative.
  
  Using a binomial  expansion
  $$
   (\beta+c_1r^2\pi^*\omega)^{k-p-j}=\sum_{q=0}^{k-p-j}  {k-p-j\choose q}c_1^q r^{2q} \pi^*(\omega^q)\wedge \beta^{k-p-j-q},
  $$
the  equality of the lemma follows  from   \eqref{e:Lelong-numbers} and the  explicit formula of $\hat\nu_j(T,B,r)$  in \eqref{e:mass-indicators}.
 \endproof
\begin{definition}\label{D:sup}\rm
Fix an open neighborhood $\bfU$ of $\overline B$ and an open neighborhood $\bfW$ of $\partial B$ in $X$ with $\bfW\subset \bfU.$
Let $\widetilde\CL^{1,1}_p(\bfU,\bfW)$ be the  set of all $T\in \CL^{1,1}_p(\bfU,\bfW)$  whose  a sequence of approximating  forms $(T_n)_{n=1}^\infty$
satisfies the following   condition:
 \begin{equation}\label{e:unit-CL-1,1} \|T_n\|_{\bfU}\leq  1\qquad\text{and}\qquad  \| T_n\|_{\Cc^1(\bfW)}\leq 1.\end{equation}

Let $\Mc(T)$ be  a mass indicator of a current $T$ and  $\Fc$  a  class of currents. 
 We denote  by $\sup_{T\in\Fc}\Mc(T)$  the supremum of $\Mc(T)$  when  $T$ is taken over $\Fc.$
\end{definition}

 As  an immediate consequence of  Theorem \ref{T:top-Lelong-closed}, we get  the following   finiteness  for the above mass indicators.
\begin{corollary}\label{C:finiteness-top-mass-ind}
 \begin{enumerate}
  \item the function $(0,\bfr)\ni r\mapsto \hat\nu_\top(T,r,\tau)\in\R^+$ is increasing and $\lim_{r\to 0+} \hat\nu_\top(T,r,\tau)\in\R^+$
  and    $\sup_{T\in  \widetilde\CL^{1,1}_p(\bfU,\bfW)} \hat\nu_\top(T,\bfr,\tau)<\infty.$   
  \item $\sup_{T\in  \widetilde\CL^{1,1}_p(\bfU,\bfW),\ r \in(0,\bfr]}\hat\kappa^{'\bullet}_\top(T,r,\tau)<\infty.$   
 \end{enumerate}
\end{corollary}

  To prove  Theorem \ref{T:Lelong-closed-all-degrees}  we reformulate it in a more  technical way, which allows
  us to make an induction argument.

\begin{theorem}\label{T:Lelong-closed-all-degrees-bis}
 We  keep the   the    assumption of   Theorem  \ref{T:Lelong-closed-all-degrees}.
 Then
 the following  assertions hold.  
  \begin{enumerate}
  \item[(1)--(5)] The corresponding  assertions  (1)--(5) of Theorem \ref{T:Lelong-closed-all-degrees} hold.
  
  For  the  remaining two assertions, we assume that  $\lowm\leq j\leq \upm.$
  If moreover,  we are in 
  the hypothesis of  assertion (5),  then  we assume that   $0\leq j\leq \upm.$
  
  \item [(6)]  
   $\hat\kappa^{'\bullet}_j(T,\bfr,\tau)<\infty.$  
   
   \item  [(7)]  $\sup_{T\in  \widetilde\CL^{1,1}_p(\bfU,\bfW),\ r\in (0,\bfr]} \hat\nu_j(T,r,\tau)<\infty.$ 
\end{enumerate}
 \end{theorem}
\proof[End of the proof of Theorem \ref{T:Lelong-closed-all-degrees}]
It  follows from assertions (1)--(6) of  Theorem \ref{T:Lelong-closed-all-degrees-bis}.
\endproof
\proof[Proof of  assertion (1) of Theorem \ref{T:Lelong-closed-all-degrees-bis}]

First assume  that the current $T$ is a closed $\Cc^1$-smooth form. Recall  from the  hypothesis that  $\tau$ is  holomorphic and
 the   identity   $\ddc\omega^j=0$ holds on $B$ for   $1\leq j\leq   \upm-1.$  Therefore,   we have for $1\leq j\leq   \upm-1$ that
\begin{eqnarray*}
\ddc [(\tau_*T)\wedge \pi^*\omega^{j}]&=&{i\over \pi}(\partial\tau_*T)\wedge (\dbar\pi^*\omega^{j})-{i\over \pi} (\dbar\tau_*T)\wedge (\partial\pi^*\omega^{j})\\
&=&{i\over\pi}(\tau_*\partial T)\wedge (\dbar \pi^*\omega^{j})-{i\over \pi} (\tau_*\dbar T)\wedge (\partial\pi^*\omega^{j})=0 .
\end{eqnarray*}
Recall   from the proof of  assertion (1) of Theorem \ref{T:top-Lelong-closed} that the equality
\begin{equation}\label{e:ddc-equal-zero} \ddc [(\tau_*T)\wedge \pi^*\omega^{j}]=0
\end{equation} also  holds for $j=\upm,$ and hence for all $0\leq j\leq \upm.$ 
Applying Theorem  \ref{T:Lelong-Jensen}  to  $\tau_*T\wedge \pi^*(\omega^{j})$ with $r_0=0$ and using   the  above  equality,  we get, for $0<r_1<r_2<\bfr$ except  for a countable set of values,  that
\begin{equation}\label{e:Lelong-closed-all-degrees-bis-asser(1)}
  \nu_{j}(T,B,r_2,\tau)- \nu_{j}(T,B,r_1,\tau)= \int_{\Tube(B,r_1,r_2)}\tau_*T\wedge \pi^*(\omega^{j})\wedge\alpha^{k-p-j}+\lim\limits_{n\to\infty}\Vc(\tau_*T_n\wedge \pi^*(\omega^{j}),r_1,r_2).
\end{equation}
On  the other hand,  since $j\geq \lowm$ we  get that  $k-p-j\leq  k-l.$  Therefore, we can apply
Theorem \ref{T:vertical-boundary-closed} to the current $\tau_*T\wedge \pi^*(\omega^{j}),$ which gives  that  $\Vc(\tau_*T\wedge \pi^*(\omega^{j}),r_1,r_2)=O(r_2).$ This proves  assertion (1) in  the special case where $T$ is $\Cc^1$-smooth.

Now  we consider the general case where $T$ is a general positive closed  $(p,p)$-current such that $T=T^+-T^-,$ where $T^\pm$ are approximable
along $B\subset V$ by positive closed $\Cc^1$-smooth  $(p,p)$-forms $(T^\pm_n)$  with $\Cc^1$-control on boundary. 
So $T^+_n\to T^+$ and $T^-_n\to T^-$ as $n$ tends to infinity.  
By the  previous case applied  to $T^\pm_n,$ we get that
\begin{equation*}
\nu_j( T^\pm_n,B,r_2,\tau)- \nu_j(T^\pm_n,B,r_1,\tau)=\kappa_j(T^\pm_n,B,r_1,r_2,\tau) +O(r_2).
\end{equation*}
Letting  $n$ tend to infinity, we infer that
\begin{equation*}
\nu_j( T^\pm,B,r_2,\tau)- \nu_j(T^\pm,B,r_1,\tau)=\kappa_j(T^\pm,B,r_1,r_2,\tau)+O(r_2).
\end{equation*}
This implies assertion (1)  since $T=T^+-T^-.$
\endproof

The  remaining  assertions of Theorem \ref{T:Lelong-closed-all-degrees-bis} will be proved by decreasing  induction on $j\in[0,\upm].$

The theorem   for $j=\upm$ is a  consequence of Theorem  \ref{T:top-Lelong-closed}.
Suppose that the theorem  is true  for  all $j$ such that  $j_0<j\leq \upm,$ where $j_0$ is a given integer with $0\leq j_0  <\upm.$  We need to show that the theorem is also true for $j=j_0.$

The plan of the proof is  as follows. We first establish some preliminary results, next we prove  assertion (6) for $j_0,$
next we prove  assertions (2)--(5)  for $j_0,$ and  finally we prove assertion (7) for $j_0.$


\begin{lemma}\label{L:inductive-holo-Lelong-numbers}
For  every $0\leq j\leq  \upm$ and for all $r_1,\ r_2\in  (0,\bfr]$ with $r_1<r_2$  except for a  countable of values,  we have
 \begin{multline*}
  \int_{\Tube(B,r_1,r_2)}(\tau_*T)\wedge \pi^*(\omega^{j})\wedge(\hat\alpha')^{k-p-j}
  = \nu_{j}(T,B,r_2,\tau)-\nu_{j}(T,B,r_1,\tau) \\
  +\sum_{q=1}^{\upm-j}
  {k-p-j\choose q}c_1^q\big(\nu_{j+q}(T,B,r_2,\tau)-  \nu_{j+q}(T,B,r_1,\tau)\big)+O(r_2).
 \end{multline*}
\end{lemma}
\proof
It follows from \eqref{e:hat-alpha'} that
\begin{equation*} (\hat \alpha')^{k-p-j}-\alpha^{k-p-j}=\sum_{q=1}^{k-p-j} {k-p-j 
\choose q} c^q_1\pi^*(\omega^q)\wedge \alpha^{k-p-j-q}.
\end{equation*}
So we get that 
\begin{multline*}
 \int_{\Tube(B,r_1,r_2)}(\tau_*T)\wedge \pi^*(\omega^{j})\wedge(\hat\alpha')^{k-p-j}  = \int_{\Tube(B,r_1,r_2)}(\tau_*T)\wedge \pi^*(\omega^{j})\wedge\alpha^{k-p-j}\\
  +  \sum_{q=1}^{k-p-j} {k-p-j\choose q} c^{q}_1\int_{\Tube(B,r_1,r_2)}(\tau_*T)\wedge  \pi^*(\omega^{j+q})\wedge 
  \alpha^{k-p-j-q}.
\end{multline*} 
On the  other hand, by assertion (1) of Theorem \ref{T:Lelong-closed-all-degrees-bis}, we have for $0\leq q\leq k-p-j$ that
\begin{equation*}
  \nu_{j+q}(T,B,r_2,\tau)- \nu_{j+q}(T,B,r_1,\tau)=  \int_{\Tube(B,r_1,r_2)}(\tau_*T)\wedge \pi^*(\omega^{j+q})\wedge\alpha^{k-p-j-q}+O(r_2).
\end{equation*}
This completes the proof.
\endproof

\proof[Proof of  assertion (6) of Theorem \ref{T:Lelong-closed-all-degrees-bis}]
Let $T$ be  a  $\Cc^1$-smooth positive  current and  let $0<r<\bfr.$
Applying  Lemma  \ref{L:inductive-holo-Lelong-numbers} to $T$ and $0<r_1<r_2\leq\bfr,$ we  get that
\begin{multline*}
  \int_{\Tube(B,r_1,r_2)}(\tau_*T)\wedge \pi^*(\omega^{j_0})\wedge(\hat\alpha')^{k-p-j_0}
  = \nu_{j_0}(T,B,r_2,\tau)-\nu_{j_0}(T,B,r_1,\tau) \\
  +\sum_{q=1}^{\upm-j_0}
  {k-p-j_0\choose q}c^q_1\big(\nu_{j_0+q}(T,B,r_2,\tau)-  \nu_{j_0+q}(T,B,r_1,\tau)\big)+O(r_2).
 \end{multline*}
 This, combined  with \eqref{e:hat-nu-vs-nu}, implies  that
 \begin{equation*}
  \int_{\Tube(B,r_1,r_2)}(\tau_*T)\wedge \pi^*(\omega^{j_0})\wedge(\hat\alpha')^{k-p-j_0}
  = \hat\nu_{j_0}(T,B,r_2,\tau)-\hat\nu_{j_0}(T,B,r_1,\tau)  +O(r_2).
 \end{equation*}
By Lemma \ref{L:hat-alpha-beta} (see \eqref{e:hat-alpha'}), the  form $\hat\alpha'$  is positive  smooth outside $V.$  Moreover, $\tau$ is holomorphic  and the current  $T$ is  positive.
Hence, the LHS is  $\geq 0.$ On the  other hand, by  Lemma  \ref{L:hat-nu} $\hat\nu_{j_0}(T,B,r_1,\tau)\geq 0.$
We infer that 
\begin{equation*}
0\leq \int_{\Tube(B,r_1,r_2)}(\tau_*T)\wedge \pi^*(\omega^{j_0})\wedge(\hat\alpha')^{k-p-j_0}\leq\hat\nu_{j_0}(T,B,r_2,\tau) +O(r_2).
 \end{equation*}
  Letting $r_2\to \bfr-$ and $r_1\to  0+$ and using the hypothesis of induction, we see that the RHS is   finite.
  Hence,  so is  the LHS.  This completes the proof of  assertion (6) for $j=j_0$ but only  for  every   $\Cc^1$-smooth positive form $T$ on
  $\bfU$  satisfying \eqref{e:unit-CL-1,1} (with $T$ in place of $T_n$ therein).
  
Now  let $T\in \widetilde\CL^{1,1}_p(\bfU,\bfW)$ be a general current with a sequence of approximating  forms $(T_n)_{n=1}^\infty$  satisfying 
\eqref{e:unit-CL-1,1}.  We have  demonstrated  that for every $n\geq 1,$
$$
 \hat\kappa^{'\bullet}_{j_0}(T_n,\bfr,\tau) \leq M<\infty.
$$
Since the forms in the integral formula  of $ \hat\kappa^{'\bullet}_{j_0}(T,\bfr,\tau)$ are positive and $T_n$ converge weakly to $T,$
we infer that $\hat\kappa^{'\bullet}_{j_0}(T,\bfr,\tau) \leq M.$  This completes the proof of  assertion (6) for $j=j_0.$ 
\endproof
\proof[Proof of  assertion (2) of Theorem \ref{T:Lelong-closed-all-degrees-bis}]  Fix $r_1,\ r_2\in  (0,\bfr]$ with $r_1<r_2.$
We prove the  following fact  by  decreasing induction on $j.$
\\
{\noindent \bf Fact.} {\it Assertion (2)  as well as  inequality
\begin{equation}\label{e:holo-finite-series}
 \sum_{n=0}^\infty\big|\nu_{j}(T,B,{r_2\over 2^n },\tau)- \nu_{j}(T,B,{r_1\over 2^n},\tau)\big|<\infty
\end{equation}
hold for $j.$
}

 Assertion (2) for $j=\upm$ is a  consequence of Theorem  \ref{T:top-Lelong-closed} (2). To  prove the above fact for $j=\upm,$
 it remains to establish \eqref{e:holo-finite-series} for  $j=\upm.$  By Theorem  \ref{T:top-Lelong-closed} (1), we  see that
 $$
 \sum_{n=0}^\infty\big|\nu_{\upm}(T,B,{r_2\over 2^n },\tau)- \nu_{\upm}(T,B,{r_1\over 2^n},\tau)\big|=\sum_{n=0}^\infty\kappa_\upm(T,B,{r_1\over 2^n },{r_2\over 2^n},\tau)\leq  c\kappa^\bullet_\upm(T,B,r_2).
 $$
By  Theorem  \ref{T:top-Lelong-closed} (3) the last term is finite. Hence, \eqref{e:holo-finite-series} for  $j=\upm$ follows.

Suppose  that the  fact is true  for all $j$  such that  $j_0<j\leq \upm,$ where $j_0$ is  a given non-negative integer  with
$0 \leq j_0< \upm.$  We need to show that  the fact is  also true  for $j=j_0.$
Let $\lambda\geq 1.$ 
 By Lemma \ref{L:inductive-holo-Lelong-numbers}, 
  we have
 \begin{equation*}
 | \nu_{j_0}(T,B,r_2/\lambda ,\tau)-\nu_{j_0}(T,B,r_2/\lambda,\tau)| \leq  \epsilon_{\lambda}:=\epsilon'_\lambda+\epsilon''_\lambda,
  \end{equation*}
 where    $\epsilon'_\lambda:=\int_{\Tube(B,r_1/\lambda,r_2/\lambda)}(\tau_*T)\wedge \pi^*(\omega^{j_0})\wedge(\hat\alpha')^{k-p-j_0}.$  
 \begin{equation*}
 \epsilon''_\lambda :=\sum_{q=1}^{\upm-j_0}
  {k-p-j\choose q}c^q_1\big|\nu_{j_0+q}(T,B,r_2/\lambda,\tau)-  \nu_{j_0+q}(T,B,r_1/\lambda,\tau)\big| +O({r_2\over \lambda}).
 \end{equation*}  
Observe that   there is a constant $c>0$ depending on $\lambda,$ $r_1$ and $r_2$  such that
$$
\sum_{n=0}^\infty \epsilon'_{2^n \lambda}\leq c  \int_{\Tube(B,\bfr)}(\tau_*T)\wedge \pi^*(\omega^{j_0})\wedge(\hat\alpha')^{k-p-j_0}<\infty,
$$
where the finiteness of the last integral  holds by   assertion (6).  On  the other hand,
 by  the  inductive hypothesis  of inequality \eqref{e:holo-finite-series}, we  see easily that 
$
\sum_{n=0}^\infty \epsilon''_{2^n \lambda} <\infty.
$
Therefore, we obtain that  $
\sum_{n=0}^\infty \epsilon_{2^n \lambda} <\infty.
$
Consequently, by Lemma  \ref{L:elementary},  $\lim\limits_{\lambda\to\infty} \nu_{j_0}(T,B,r/\lambda ,\tau)$
exists and  is  finite.
This  proves  assertion (2) for $j=j_0.$
\endproof

\proof[Proof of  assertion (3) of Theorem \ref{T:Lelong-closed-all-degrees-bis}] 
By   \eqref{e:hat-alpha'},  we have that
$
\alpha:=\hat\alpha'-c_1\pi^*\omega.
$
Inserting   this into   \eqref{e:Lelong-log-numbers} and using  the third formula in \eqref{e:mass-indicators}, we get that  for $0<s<r\leq \bfr,$
 \begin{equation*}\begin{split} 
  \kappa_j(T,B, s,r,\tau)&=\int\limits_{\Tube(B,s,r)} \tau_*(T)
 \wedge\pi^*\omega^j\wedge  (\hat\alpha'-c_1\pi^*\omega)^{k-p-j}\\
&=\sum_{q=0}^{k-p-j}{k-p-j \choose  q}\int\limits_{\Tube(B,s,r)} \tau_*(T)
 \wedge\pi^*\omega^j\wedge  (\hat\alpha')^{k-p-j-q} \wedge   (-1)^qc^q_1\pi^*(\omega^q)\\
 &= \sum_{q=0}^{k-p-j}  (-1)^qc^q_1{k-p-j \choose  q} 
 \hat\kappa'_{j+q}(T,B, s,r,\tau).
 \end{split}\end{equation*}
Using this and  the fact that $\alpha'\geq 0$ outside $V,$   we infer that
\begin{equation}\label{e:kappa-vs-kappa'}|\kappa_j(T,B,s,r,\tau)|\leq c\sum_{q=0}^{\upm-j} \hat\kappa'_{j+q}(T,B,s,r,\tau).\end{equation}
By  assertion (6),  the  RHS tends to $0$  as $r$ tends to $0.$  So
$$
\lim_{r\to 0,\ s<r}|\kappa_j(T,B,s,r,\tau)|=0.
$$
This, coupled with  \eqref{e:Lelong-log-bullet-numbers}, implies that 
\begin{equation*}
  \lim _{r\to 0}\kappa^\bullet_j(T,B,r,\tau)= \lim_{r\to 0}\big(\limsup\limits_{s\to0+}   \kappa_j(T,B,s,r,\tau)\big)  =0,
  \end{equation*}
  as  desired.
\endproof


\proof[Proof of  assertion (4) of Theorem \ref{T:Lelong-closed-all-degrees-bis}]
First, 
we  will prove  the interpretation of assertion (4)  in the spirit  of \eqref{e:Lelong-number-point-bisbis(1)}.  Since  $q:=k-p-j<k-l,$
we infer from   Theorems  \ref{T:Lelong-Jensen-smooth}  and \ref{T:vertical-boundary-terms}  that
$$
\kappa_j(T^\pm_n,B,r,\tau)=\nu_j( T^\pm_n,B,r,\tau)+O(r).
$$
Consequently,
\begin{eqnarray*}
\kappa_j(T,B,r,\tau)&:=& \lim\limits_{n\to\infty} \kappa_j(T^+_n-T^-_n,B,r,\tau)=\lim\limits_{n\to\infty}\nu_j(T^+_n,B,r,\tau)-\lim\limits_{n\to\infty}\nu_j(T^-_n,B,r,\tau)+O(r)\\
&=&\nu_j(T^+,B,r,\tau)-\nu_j(T^-,B,r,\tau)+O(r) =\nu_j(T,B,r,\tau)+O(r).
\end{eqnarray*}
This  implies  the  desired interpretation  according to  Definition \ref{D:Lelong-log-numbers(1)}.

Second,  we  will prove  the interpretation of assertion (4)  in the spirit  of \eqref{e:Lelong-number-point-bisbis(2)}.
To start  with,  we  fix $0<r<\bfr$ and let  $0<\epsilon<r.$  Theorem \ref{T:Lelong-Jensen-eps} applied  to  $\tau_*T\wedge \pi^*(\omega^j)$  and using identity \eqref{e:ddc-equal-zero} gives
 $$
 {1\over  (r^2+\epsilon^2)^{k-p-j}} \int_{\Tube(B,r)} \tau_*T\wedge \pi^*(\omega^j)\wedge \beta_\epsilon^{k-p-j}
 = \lim\limits_{n\to\infty} \Vc_\epsilon( \tau_*T_n\wedge \pi^*(\omega^j),r)+   \int_{\Tube(B,r)} \tau_*T\wedge \pi^*(\omega^j)\wedge \alpha^{k-p-j}_\epsilon.
$$
Now  we let $\epsilon$ tend to $0.$  Then the LHS  tends to $\nu_j(T,B,r,\tau).$  On the other hand,  we deduce  from \eqref{e:vertical-boundary-term-eps} and  Theorem \ref{T:vertical-boundary-eps}  that $\Vc_\epsilon( \tau_*T_n\wedge \pi^*(\omega^j),r)=O(r).$  Consequently,    the second  term on the RHS  tends to $\nu_j(T,B,r,\tau)+O(r).$
This  proves the  desired  interpretation  according to  Definition \ref{D:Lelong-log-numbers(2)}.
\endproof



\proof[Proof of  assertion (5) of Theorem \ref{T:Lelong-closed-all-degrees-bis}]
Recall from  the hypothesis of this  assertion  that  $\supp(T^\pm_n)\cap V\subset  B$ for $n\geq 1.$
 First  we   explain how  to prove the following  stronger  version of assertion (1) also holds:
  For  
  \begin{equation}\label{e:stronger-version-asser(1)-closed} \nu_j(T,B,r_2,\tau)-\nu_j(T,B,r_1,\tau)=\kappa_j(T,B,r_1,r_2,\tau)\quad\text{for}\quad 0\leq j\leq \upm\quad\text{and}\quad 0<r_1<r_2<\bfr.
  \end{equation}
  We argue  as  in the proof  of assertion (1) of Theorem \ref{T:Lelong-closed-all-degrees-bis}.
  However in \eqref{e:Lelong-closed-all-degrees-bis-asser(1)} we apply  Theorem \ref{T:Lelong-Jensen-closed-compact-support}
  instead of Theorem \ref{T:Lelong-Jensen}. Consequently,  in the  present context we get  $\Vc(\tau_*T_n\wedge \pi^*(\omega^l),r_1,r_2)=0$
  and the above identity follows. Note that in  the present context we do not need Theorem \ref{T:vertical-boundary-closed}. That is why our result is valid
  for  $0\leq j\leq \upm.$
  
  Using the above stronger version of  assertion (1) of Theorem \ref{T:Lelong-closed-all-degrees-bis},  we argue as in  the proof of assertions
  (2)--(4) of this theorem  in order to extend the validity of these assertions to   $0\leq j\leq \upm.$
  \endproof
  

\proof[Proof of  assertion (7) of Theorem \ref{T:Lelong-closed-all-degrees-bis}]
We prove  assertion (7)  for $j_0$   using  assertion (6) for all $j$ with $j_0\leq j\leq \upm.$ Let $0<s<\bfr.$ Applying 
identity \eqref{e:stronger-version-asser(1)-closed} to $j_0$ and $r_1:=s,$ $r_2:=\bfr,$  we  get 
$$\nu_{j_0}(T,B,\bfr,\tau)-\nu_{j_0}(T,B,s,\tau)=\kappa_{j_0}(T,B,s,\bfr,\tau).$$
By \eqref{e:kappa-vs-kappa'} we know that
\begin{equation*}|\kappa_{j_0}(T,B,s,\bfr,\tau)|\leq c\sum_{q=0}^{\upm-j_0} \hat\kappa'_{j_0+q}(T,B,s,\bfr,\tau).
\end{equation*}
Using  assertion (6) for all $j$ with $j_0\leq j\leq \upm,$ we see that the expression on the RHS is uniformly bounded independently of $s.$
Consequently, $\nu_{j_0}(T,B,s,\tau)$ is uniformly bounded independently of $s.$
On the other hand,  applying  Lemma \ref{L:hat-nu}, we obtain that
$$
\limsup_{s\to0+} \hat\nu_{j_0}(T,B,s,\tau)=\limsup_{s\to0+} \nu_{j_0}(T,B,s,\tau)+\sum\limits_{q=1}^{k-p-j_0}  {k-p-j_0\choose q}c_1^q \nu_{j_0+q}(T,B,\tau).
$$
Since the expression on the RHS is  bounded,  so is the LHS.  This  proves assertion (7) for $j_0.$
  \endproof
  
\section{Regularization of currents and admissible estimates}\label{S:Regularization}


\subsection{Extended Standing Hypothesis}\label{SS:Ex-Stand-Hyp}

In this subsection we  introduce a  standard   setting for further technical  developments.
We keep the Standing Hypothesis  formulated in Subsection \ref{SS:Global-setting}.

Let $B$ be a  relatively compact piecewise $\Cc^2$-smooth open
subset. Let $V_0$ be a relatively compact open subset of $V$ such that $B \Subset V_0 .$
Consider  a  strongly admissible map $\tau:\ \bfU\to\tau(\bfU)$  along $B,$  with $\bfU$ a neighborhood of $\overline B$ in $X.$
By shrinking $\bfU$ if necessary, we may   fix a finite collection $\Uc=(\bfU_\ell,\tau_\ell)_{1\leq \ell\leq \ell_0} ,$
 of holomorphic admissible maps for $\bfU.$  More precisely,  
there is a  finite cover of $\overline \bfU$ by open subsets $\bfU_\ell,$ $1\leq \ell\leq  \ell_0,$ of $X$
such that   there is   a holomorphic coordinate system on $\overline \bfU_\ell$ in $X$ and $\bfU_\ell$  is  biholomorphic  to  $\U_\ell:=\tau_\ell(\bfU_\ell)\subset \E$
by a  holomorphic admissible map $\tau_\ell.$ By  choosing $\bfr>0$ small enough, we may assume  without loss of generality that $\overline\Tube(B,\bfr)\Subset \U:=\bigcup_{\ell=1}^{\ell_0} \U_\ell.$ 
Fix  a partition of unity  $(\theta_\ell)_{1\leq \ell\leq \ell_0}$ subordinate   to the open cover  $(\bfU_\ell\cap V)_{1\leq \ell\leq \ell_0}$   of $\overline{\bfU\cap V}$  in $V$  such that $\sum_{1\leq \ell\leq \ell_0}  \theta_\ell=1$ on an open neighborhood of $\overline {\bfU\cap V} \subset V.$ We  may  assume  without loss of generality  that there are open  subsets
$\widetilde V_\ell\subset V$  for $1\leq \ell\leq \ell_0$ such that
\begin{equation}\label{e:supp} \supp(\theta_l)\subset  \widetilde V_\ell\Subset \bfU_\ell\cap V \quad\text{and}\quad 
\tau( \widetilde V_\ell)\Subset  \U_\ell\quad\text{and}\quad \pi^{-1}(\supp(\theta_\ell))\cap \U\subset \U_\ell.
\end{equation}
 For $1\leq \ell\leq \ell_0$ set 
 \begin{equation}\label{e:tilde-tau_ell}\tilde \tau_\ell:=\tau\circ\tau_\ell^{-1}.
 \end{equation}
  So  $\tilde \tau_\ell$ defines a map  from  $\U_\ell\subset \E$  onto  $\tau(\bfU_\ell)\subset \E.$

  We also assume  that  for every $1\leq\ell\leq\ell_0,$  there is a local  coordinate system $y=(z,w)$  on $\U_\ell$  with  $V\cap \U_\ell=\{z=0\}.$
  
  $\Uc=(\bfU_\ell,\tau_\ell)_{1\leq \ell\leq \ell_0}$ is  said to be  a  {\it covering  family of  holomorphic admissible maps for $B.$}
  
  Now  we  formulate the
  
  \noindent  {\bf  Extended  Standing  Hypothesis.}  {\it   This means that  we assume   the usual  Standing Hypothesis (introduced  in Subsection \ref{SS:Global-setting}) and      a  covering  family $\Uc=(\bfU_\ell,\tau_\ell)_{1\leq \ell\leq \ell_0}$ of  holomorphic admissible maps for $B$  as described above.} 
\subsection{Representative current, regularization of currents and an elementary lemma}\label{SS:Regularizations}

Recall that $\pi:\  \E\to V$ is  the canonical projection.
For  every current $S$ defined  on $\bfU_\ell\subset X,$ we    denote  by $S_\epsilon,$ or equivalently $(S)_\epsilon,$  with $\epsilon>0,$  a  family of  forms which  regularize $S$ by convolution.
\begin{definition}\label{D:T-hash}\rm

Let  $ T$ be a current  defined on $\bfU.$ Consider the current  $T^{\#}$ defined  on $\U$ by the   following formula:
\begin{equation}\label{e:T-hash}
T^{\#}:=\sum_{\ell=1}^{\ell_0}  (\pi^*\theta_\ell)\cdot (\tau_\ell)_*( T|_{\bfU_\ell}).
\end{equation}
By  \eqref{e:supp},  $T^{\#}$ is  well-defined.
We also  consider the smooth regularizing   forms  $(T_\epsilon)_{\epsilon>0}$ on $\U$ defined by
\begin{equation}\label{e:T_eps}
T_\epsilon:= \sum_{\ell=1}^{\ell_0}(\pi^*\theta_\ell)\cdot(\tau_\ell)_*\big(( T|_{\bfU_\ell})_\epsilon\big).
\end{equation}
\end{definition}
The following observation is  an immediate consequence of Definition  \ref{D:T-hash}.
\begin{lemma}\label{L:regularization}
Let  $ T$ be a current of bidegree $(p,p)$ defined on $\bfU.$
\begin{enumerate}
 \item the forms $T_\epsilon$ are smooth   of bidegree $(p,p),$ and  $T_\epsilon$ converge to $T$ weakly    on $X$
 and  $T_\epsilon^{\#}$ converge to $T^{\#}$ weakly    on $\E$  as $\epsilon$ tends
 to $0.$
 \item  If moreover $T$ is  positive,  then so are the forms $T_\epsilon$ and the current  $T^{\#}.$
\end{enumerate}

\end{lemma}

\begin{lemma}\label{L:basic-T-hash}
For  every current $R$ on $\U,$ the following identity holds
 $$\langle \tau_*T,R\rangle -\langle T^\hash, R\rangle =\sum_{\ell=1}^{\ell_0}\langle(\tau_\ell)_* T,   (\tilde\tau_\ell)^*[(\pi^*\theta_\ell )R]-[(\pi^*\theta_\ell )R]\rangle .$$
\end{lemma}
\proof
Since  $\sum_{\ell=1}^{\ell_0}  \pi^*\theta_\ell  =1$  on an open  neighborhood  of $\pi^{-1}(\overline{\bfU\cap V})\subset \pi^{-1}(V),$  we have
$$\tau_*T=\sum_{\ell=1}^{\ell_0}  \pi^*\theta_\ell\cdot\tau_* T\qquad\text{on}\qquad \U.$$
So
$$
\langle \tau_*T,R\rangle -\langle T^\hash, R\rangle =\sum_{\ell=1}^{\ell_0}   ( \langle \tau_*T,(\pi^*\theta_\ell)  R\rangle
- \langle (\tau_\ell)_*( T|_{\bfU_\ell}),(\pi^*\theta_\ell )R\rangle ).
$$
Writing   $\tau_*T=(\tau\circ \tau_\ell^{-1})_*(\tau_\ell)_* T=(\tilde\tau_\ell)_*(\tau_\ell)_* T$ on $\U_\ell,$
we get that
\begin{eqnarray*}
\langle \tau_*T,R\rangle -\langle T^\hash, R\rangle &=&\sum_{\ell=1}^{\ell_0}  \big( \langle (\tilde\tau_\ell)_*(\tau_\ell)_* T,(\pi^*\theta_\ell )R\rangle - \langle (\tau_\ell)_*( T|_{\bfU_\ell}),(\pi^*\theta_\ell )R\rangle \big)\\
&=&\sum_{\ell=1}^{\ell_0}  \big( \langle(\tau_\ell)_* T, (\tilde\tau_\ell)^*[(\pi^*\theta_\ell )R]\rangle - \langle 
(\tau_\ell)_*( T),(\pi^*\theta_\ell )R\rangle \big)\\
&=& \sum_{\ell=1}^{\ell_0}  \langle(\tau_\ell)_* T, (\tilde\tau_\ell)^*[(\pi^*\theta_\ell )R]-[(\pi^*\theta_\ell )R]\rangle ,
\end{eqnarray*}
which implies the desired identity.
\endproof

The  following  elementary lemma will be repeatedly  used in the sequel.

\begin{lemma}\label{L:difference}
Let  $U$ be an open neighborhood of $0$ in $\C^n.$
Let $\Ic$ be a nonempty  finite  index set.
For every $I\in \Ic,$  there are  a  number $p_I\in \N$ and  $2p_I$  continuous forms
 $f_{I1},\ldots, f_{Ip_I}$ and $\tilde f_{I1},\ldots, \tilde f_{Ip_I}$   on $U.$ 
Set $$S:= \sum_{I\in\Ic}  f_{I1}\wedge  \ldots\wedge f_{Ip_I}\quad\text{and}\quad \widetilde 
S:= \tilde f_{I1}\wedge  \ldots\wedge \tilde f_{Ip_I}.$$
\begin{enumerate}
\item Then  we have
 $$
 \widetilde S-S=\sum_{I\in\Ic}\sum_{ J\subset \{1,\ldots ,p_I\}:\ J\not=\varnothing}  (f_{I1})_J\wedge  \ldots\wedge (f_{Ip_I})_J,
 $$
 where  for $1\leq j\leq p_I,$   
 $$
 (f_{Ij})_J:=\begin{cases}
           f_{Ij},&\text{if}\quad j\not\in J;\\
           \tilde f_{Ij}-f_{Ij},&\text{otherwise}.
          \end{cases}
 $$  
\item (Application) Let $\tau:\ U\to U$  be  a  $\Cc^1$-smooth functions and suppose that  $\tilde f_{Ij}=\tau^*f_{Ij}$ for $I\in\Ic$ and  $1\leq j\leq p_I.$ Then  the above conclusion  holds and  $\widetilde S=\tau^*S.$
\end{enumerate}
 \end{lemma}
\proof
  For every $I\in\Ic$ and  $1\leq j\leq p_I,$ write  $g_{Ij}:=\tilde f_{Ij}-f_{Ij}.$  Observe  that
\begin{eqnarray*}
\widetilde S-S&=&\sum_{I\in\Ic}\tilde f_{I1}\wedge  \ldots\wedge \tilde f_{Ip_I}-f_{I1}\wedge  \ldots\wedge f_{Ip}
\\
&=&\sum_{I\in\Ic}(g_{I1}+f_{I1})\wedge \ldots (g_{Ip}+f_{Ip})- f_{I1}\wedge  \ldots\wedge f_{Ip}\\
&=&\sum_{I\in\Ic}\sum_{ J\subset \{1,\ldots ,p_I\}:\ J\not=\varnothing}  (f_{I1})_J\wedge  \ldots\wedge (f_{Ip_I})_J.
\end{eqnarray*}
This proves assertion (1).

Assertion (2) is an immediate consequence of assertion (1).
\endproof
\subsection{Admissible estimates}\label{SS:admissible-estimates}

Admissible estimates   are those  estimates which are related to   admissible maps.
This  subsection provides  necessary admissible estimates.  

Let $\U$ be an open neighborhood of $0$ in $\C^k.$
We use the local coordinates $y=(z,w)\in\C^{k-l}\times \C^l$ on $\U.$

The following  notion  will be  needed  in order  to obtain   admissible estimates.
\begin{definition}
 \label{D:precsim}
 \rm  Let $\Gamma$ be  a form of degree $2$ and $S$ a positive $(1,1)$-form  defined on $\U.$ 
 For  $(p,q)\in \{(0,2),(1,1),(2,0)\},$ $\Gamma^{p,q}$ denotes the component of bidegree $(p,q)$ of $\Gamma.$
 So $\Gamma^{1,1}=\Gamma^\sharp$  according to  Notation \ref{N:principal}.
 
 We  write $\Gamma\trianglelefteq  S$  
 if there is a constant $c>0$ such that   the  following   two inequalities  hold for $y\in\U:$
 $$
  \Gamma^{0,2}(y)\wedge \overline{\Gamma^{0,2}}(y)\leq c S^2(y)\quad\text{and}\quad  \Gamma^{2,0}(y)\wedge \overline{\Gamma^{2,0}}(y)\leq cS^2(y).
 $$
 \end{definition}
\begin{notation}\label{N:pm}
 \rm  Let $\Gamma$ and  $S$ be  two  real $(1,1)$-forms   defined on $\U.$
 
 We  write  $\Gamma\lesssim S$ if there is a constant $c>0$ such that  $\Gamma\leq c S.$
 We write  $\pm \Gamma\lesssim   S$  if  we  have both $ \Gamma\lesssim  S$ and  $- \Gamma\lesssim  S.$
 
 We  write $\Gamma\approx  S$ if we  have both $ \Gamma\lesssim  S$ and  $S\lesssim  \Gamma.$
\end{notation}

\begin{definition}\label{D:Hc}
 \rm Let $\Hc=\Hc(\U)$ be the class of all   real $(1,1)$-forms $H$ on $\U$  which  can be written as
 $$
 H=\sum f_{pq'} dz_p\wedge d\bar w_{q'}+ \sum g_{p'q} d\bar z_{p'}\wedge d w_q,
 $$
 where  $f_{pq'}$ and  $g_{p'q}$ are  bounded functions. 
\end{definition}

Now  we place ourselves under the Extended Standing  Hypothesis at the beginning of this  section.
Since $\tau$ is  strongly  admissible, we infer  from  Definition \ref{D:Strongly-admissible-maps} that the  following    estimates of $1$-forms  for 
the components of $\tau=(s_1,\ldots, s_k)$ in the  local coordinates $y=(z,w).$
Note  that $s_j= \tau^*z_j$  for $1\leq j\leq k-l$  and  $s_j= \tau^*w_{j-k+l}$  for  $k-l<j\leq k.$
 \begin{equation}\label{e:diff-dz}
 d(\tau^*z_j)-dz_j=\sum_{p=1}^{k-l}  O(\|z\|) dz_p +  O(\|z\|^2)\quad\text{and}\quad 
 d(\tau^*\bar z_j)-d\bar z_j=\sum_{p=1}^{k-l}  O(\|z\|) d\bar z_p +  O(\|z\|^2).
 \end{equation}
 \begin{equation}\label{e:diff-dw}
 d(\tau^*w_m)-dw_m=\sum_{p=1}^{k-l}  O(1) dz_p + O(\|z\|)\quad\text{and}\quad
 d(\tau^*\bar w_m)-d\bar w_m=\sum_{p=1}^{k-l}  O(1) d\bar z_p +  O(\|z\|).
 \end{equation}
Using  this we  infer   the  following  estimates for the    change  under $\tau$ of  a $\Cc^1$-smooth function  and of the basic $(1,1)$-forms  $\pi^*\omega,$ $\beta,$ $\hat\beta.$ 
\begin{proposition}\label{P:basic-admissible-estimates-I}
There  are constants $c_3,c_4>0$ such that   $ c_3 \pi^*\omega +c_4\beta\geq 0$  on  $\pi^{-1}(V_0)\subset \E$  and that
for every $1\leq \ell \leq \ell_0,$
 the following inequalities hold  on $\U_\ell \cap \Tube(B,\bfr):$
\begin{enumerate}

\item   $|\tilde \tau_\ell^*(\varphi)-\varphi|\leq c_3 \varphi^{3\over 2},$ and   $|\tilde \tau_\ell^*(f)-f|\leq c_3 \varphi^{1\over 2}$  for every  $\Cc^1$-smooth  function $f$ on  $\Tube(B,\bfr);$
 \item $\pm\big(   \tilde \tau_\ell^*(\pi^*\omega)  -\pi^*\omega-H\big)^\sharp \lesssim  c_3 \varphi^{1\over 2}  \pi^*\omega +c_4\varphi^{1\over 2} \beta,   $
  and 
 $\tilde\tau_\ell^*(\pi^*\omega)  -\pi^*\omega\trianglelefteq   c_3 \varphi^{1\over 2}  \pi^*\omega +c_4\varphi^{1\over 2} \beta;$
 \item $\pm\big(   \tilde \tau_\ell^*(\beta) -\beta  \big)^\sharp \lesssim c_3 \phi^{3\over 2}\cdot  \pi^*\omega +c_4\phi^{1\over 2}\cdot \beta,$ and 
 $\pm\big(   \tilde \tau_\ell^*(\beta) -\beta  \big) \trianglelefteq c_3 \phi^{3\over 2}\cdot  \pi^*\omega +c_4\phi^{1\over 2}\cdot \beta;$
   \item $\pm\big(   \tilde \tau_\ell^*(\hat\beta) -\hat\beta  \big)^\sharp \lesssim c_3 \phi^{3\over 2}\cdot  \pi^*\omega +c_4\phi^{1\over 2}\cdot \hat\beta,$ and 
 $\pm\big(   \tilde \tau_\ell^*(\hat \beta) -\hat\beta  \big) \trianglelefteq c_3 \phi^{3\over 2} \cdot  \pi^*\omega +c_4\phi^{1\over 2}\cdot \hat\beta.$
\end{enumerate}
Here, in the first  inequalities of  (2)-(3)-(4), $H$ is  some  form in the class $\Hc$ given in Definition  \ref{D:Hc}.
\end{proposition}
\proof 
\noindent {\bf Proof of   assertion (1).} It follows from \eqref{e:admissible-maps}.

Let $S$  be the  positive $(1,1)$-form on the  RHS of each  inequality of the above three assertions (2),(3) and (4). Let $\Gamma$  be the $2$-form  on  the corresponding   LHS.
We may assume  that the sign  on the LHS is plus $+.$  The remaining  case when the sign is minus $-$ can be treated similarly.

\noindent {\bf Proof of   assertion (2).} 
Using  \eqref{e:varphi-new-exp} and \eqref{e:alpha-beta-local} we may assume  without loss of generality that 
$$\Gamma=\tau^*(idw_q\wedge dw_{q'})- idw_q\wedge dw_{q'}\quad\text{and}\quad
S=\|z\|\big( \sum_{p=1}^{k-l} idz_p\wedge d\bar z_p+
  \sum_{m=1}^{l} idw_m\wedge d\bar w_m\big).$$
By \eqref{e:diff-dw} we see that
\begin{multline*}\Gamma=\big( \sum_{p=1}^{k-l} O(1)  dz_p +O(\|z\|) \big)\wedge d\bar w_{q'}+ \big( \sum_{p=1}^{k-l} O(1) d\bar z_p +O(\|z\|) \big) \wedge d w_{q} + \big( \sum_{p,p'=1}^{k-l} O(1) dz_p\wedge d\bar z_{p'} \big)\\
+   O(\|z\|)\{dz,d\bar z\}+O(\|z\|^2).
\end{multline*}
So there is  a $(1,1)$-form  $H\in\Hc(\U)$
such that $
\Gamma^{1,1}=H+O(\|z\|).$
 This  implies that $\pm( \Gamma^{1,1}- H)\lesssim  S.$  
 The first  inequality of assertion (2) follows.
 
 On the other hand, we also see that  
 \begin{equation} \label{e:Gamma-0,2}\Gamma^{0,2}=O(\|z\|).
 \end{equation}
  Therefore, there are constants $c'', c>0$ such that
 $$
  \Gamma^{0,2}(y)\wedge \overline{\Gamma^{0,2}}(y)\leq  c''\|z\|^2 \big( \sum_{p=1}^{k-l} idz_p\wedge d\bar z_p+
  \sum_{m=1}^{l} idw_m\wedge d\bar w_m\big)^2 \leq c S^2(y).
 $$
 This proves $
  \Gamma^{0,2}(y)\wedge \overline{\Gamma^{0,2}}(y)\leq c S^2(y).$  The inequality  $
  \Gamma^{2,0}(y)\wedge \overline{\Gamma^{2,0}}(y)\leq c S^2(y)$ can be proved similarly.
  Hence, the second  inequality of assertion (2) follows.
  
\noindent {\bf Proof of assertion (3).}  
First  we prove the first  inequality  of assertion (3). Using  \eqref{e:varphi-new-exp} and \eqref{e:alpha-beta-local} and applying  the Cauchy--Schwarz inequality,  we may assume  without loss of generality that 
$$\Gamma= \Gamma_1+\Gamma_2+\Gamma_3\quad\text{and}\quad
S=\|z\|\big( \sum_{p=1}^{k-l} idz_p\wedge d\bar z_p\big)+
  \|z\|^3\big(\sum_{m=1}^{l} idw_m\wedge d\bar w_m\big),$$
  where, for some $\Cc^1$-smooth functions $f(z,w)=O(\|z\|)$ and   $g(z,w)=O(\|z\|^2),$
  \begin{multline*}
  \Gamma_1:=\tilde\tau_\ell^*\big( A(w)idz_p\wedge d\bar z_p\big)-\big( A(w)idz_p\wedge d\bar z_p\big),\quad  \Gamma_2:= \tilde\tau_\ell^*\big( fdz_p\wedge d\bar w_{q'}\big)-\big( f dz_p\wedge d\bar w_{q'}\big),\\
   \Gamma_3:= \tilde\tau_\ell^*\big( g dw_q\wedge d\bar w_{q'}\big)-\big( g dw_q\wedge d\bar w_{q'}\big).
  \end{multline*}
By \eqref{e:diff-dz} we see that
\begin{equation}
 \label{e:Gamma_1}
 \Gamma_1=R_{11}+R_{12}+R_{13},
\end{equation}
 where 
 \begin{itemize} \item $R_{11}$ is  a $(1,1)$-form in $\{dz,d\bar z\}$ with coefficients of order $O(\|z\|);$
 \item $R_{12}$  is a $2$-form  in $\{dz,dw, d\bar z,d\bar w \}$  with coefficients  of order $O(\|z\|^2)$ such that the total  degree of  $\{dw,d\bar w\}$
 in  each term in the representation 
 of $R_{12}$  is $1;$
\item  $R_{13}$ is  a $2$-form  in  $\{dz,dw ,d\bar z,d\bar w \}$  with coefficients  of order $O(\|z\|^3).$
\end{itemize}
We  handle  $\Gamma_3$  in the same way.   Since  $g(z,w)= O(\|z\|^2),$
it follows  that  $\tilde\tau_\ell^*g-g= O(\|z\|^3)$  by assertion (1). Hence,
\begin{equation}
 \label{e:Gamma_3}
 \Gamma_3=R_{33},
\end{equation}
where $R_{33}$ has the same description as $R_{13}$ in \eqref{e:Gamma_1}.

Similarly,  Since  $f(z,w)= O(\|z\|),$ we deduce from  \eqref{e:diff-dz}-\eqref{e:diff-dw} that
\begin{equation}
 \label{e:Gamma_2}
 \Gamma_2=R_{21}+R_{22}+R_{23},
\end{equation}
where $R_{2j}$ has the same description as $R_{1j}$ in \eqref{e:Gamma_1}.
 
 By  \eqref{e:Gamma_1}--\eqref{e:Gamma_3}--\eqref{e:Gamma_2} and  using the Cauchy--Schwarz inequality,
 we  get that $\Gamma^{1,1}_j\lesssim  S$ for $1\leq j\leq 3.$
This completes the proof of the  first inequality of assertion (3).

It  also follows  from  \eqref{e:Gamma_1}--\eqref{e:Gamma_3}--\eqref{e:Gamma_2} that for $1\leq j\leq 3,$
\begin{equation}\label{e:Gamma_2,0}
 \Gamma^{2,0}_j=R'_j+R''_j\quad \text{and}\quad  \Gamma^{0,2}_j=\overline{\Gamma}^{2,0}_j=\overline R'_j+\overline R''_j,
\end{equation}
where  $R'_j$ (resp.  $R''_j)$  is  a $(2,0)$-form  which   has the same description as $R_{12}$  (resp. $R_{13}$)  in \eqref{e:Gamma_1}. 
Using \eqref{e:Gamma_2,0} and the above  expression of $S,$  and  applying the Cauchy-Schwarz  inequality,
we see that $\Gamma \trianglelefteq  S.$  Hence, the proof of the  second inequality of assertion (3)
is  complete.

\noindent {\bf Proof of assertion (4).}
Using  formula  \eqref{e:hat-beta} and applying  assertions (1), (2) and (3), we get the  desired conclusion.
 \endproof

Using \eqref{e:diff-dz}--\eqref{e:diff-dw}, we also infer   the  following  estimates for the    change  under $\tau$ of   the basic $(1,1)$-forms  $\alpha$ and  $\hat\alpha.$

 In the sequel, when we use the coordinate system $(w,\zeta',t)$  with $\zeta=(\zeta',t)$ given in \eqref{e:homogeneous-coordinates} we adopt  the  following notation for $n\in\N:$
 \begin{equation}\label{e:w,zeta',t}\begin{split} 
 O(t^n)dw&= \sum_{q=1}^l O(t^n)dw_q,\qquad O(t^n) d\zeta'= \sum_{p=1}^{k-l-1}O(t^n)d\zeta_p,\\
 O(t^n)dw\wedge dw&:=\sum_{q,q'=1}^l O(t^n)dw_q\wedge dw_{q'},\qquad O(t^n)dw\wedge d\zeta':=\sum_{p=1}^{k-l-1}
 \sum_{q=1}^{l} O(t^n)dw_q\wedge d\zeta_{p},\\
 O(t^n)dw\wedge dt&:=\sum_{q=1}^l O(t^n)dw_q\wedge dt,\\
  O(t^n)d\zeta'\wedge d\zeta'&:=\sum_{p,p'=1}^{k-l-1} O(t^n)d\zeta_p\wedge d\zeta_{p'},\quad  O(t^n)d\zeta'\wedge dt:=\sum_{p=1}^{k-l-1} O(t^n)d\zeta_p\wedge dt.
 \end{split}
 \end{equation}
\begin{proposition}\label{P:basic-admissible-estimates-II}
There  are constants $c_3,c_4>0$  such that   the conclusion of Proposition \ref{P:basic-admissible-estimates-I}  holds and that
for every $1\leq \ell \leq \ell_0,$
 the following inequalities hold  on $\U_\ell \cap \Tube(B,\bfr):$
\begin{enumerate}

\item  
$\pm\big(   \tilde \tau_\ell^*(\alpha) -\alpha  \big)^\sharp\lesssim   c_3 \pi^*\omega +c_4\beta+c_3\varphi^{1/2}\alpha$ and 
$  \big(   \tilde \tau_\ell^*(\alpha) -\alpha  \big)\trianglelefteq   c_3 \pi^*\omega +c_4\beta+c_3\varphi^{1/4}\alpha ;$
\item 
$\pm\big(  \tilde \tau_\ell^*(\hat\alpha) -\hat\alpha -H\big)^\sharp \lesssim  c_3 \pi^*\omega +c_4\hat\beta +c_3\varphi^{1/2}\hat\alpha$ and 
$\big(  \tilde \tau_\ell^*(\hat\alpha) -\hat\alpha  \big) \trianglelefteq   c_3 \pi^*\omega +c_4\hat\beta +c_3\varphi^{1/4}\hat\alpha.$   Here,  $H$ is  some  form in the class $\Hc$ given in Definition  \ref{D:Hc}.           
\end{enumerate} 
\end{proposition}

\proof
We use the  homogeneous  coordinates $\zeta=(\zeta',\zeta_{k-l})$  introduced in  \eqref{e:max-coordinate}--\eqref{e:homogeneous-coordinates}. 
For  $(z,w)\in  \U_l\cap \Tube(B,\bfr),$ write $(\tilde z,\tilde w)=\tilde\tau_\ell(z,w).$
Let $\tilde\zeta$  be the homogeneous  coordinate   of $\tilde z.$  For the  sake of simplicity
we  will  also  write  $t$ (resp.  $\tilde t$) instead of $\zeta_{k-l}=z_{k-l}$  (resp.   $\tilde\zeta_{k-l}=\tilde z_{k-l}$).

Write
\begin{equation}\label{e:admissible-alpha(1)}
       \begin{split}
 \alpha&=\ddcwzeta\log{\| A(w)(z)\|^2}=\ddcw\log{\| A(w)(z)\|^2}\\
 &+\partial_w\dbar_{\zeta}\log{\| A(w)(z)\|^2}+\dbar_w\partial_{\zeta}\log{\| A(w)(z)\|^2} +\ddczeta\log{\| A(w)(z)\|^2}\\
 &= I_1+I_2+I_3+I_4.
\end{split}          
\end{equation}
We also write
\begin{equation}\label{e:admissible-alpha(2)}
       \begin{split}
 (\tilde\tau_\ell)^*\alpha&=(\tilde\tau_\ell)^*[\ddcwzeta\log{\| A(w)(z)\|^2}]=(\tilde\tau_\ell)^*[\ddcw\log{\| A(w)(z)\|^2}]\\
 &+(\tilde\tau_\ell)^*[\partial_w\dbar_{\zeta}\log{\| A(w)(z)\|^2}]+(\tilde\tau_\ell)^*[\dbar_w\partial_{\zeta}\log{\| A(w)(z)\|^2}] +(\tilde\tau_\ell)^*[\ddczeta\log{\| A(w)(z)\|^2}]\\
 &= \widetilde I_1+\widetilde I_2+\widetilde I_3+\widetilde I_4.
\end{split}          
\end{equation}
We will show that  for $1\leq j\leq 4,$
\begin{equation}\label{e:widetildeI_j-I_j}
\widetilde I_j-I_j=f_j+g_j +h_{j1} d\zeta'\wedge d\zeta'+  h_{j2} d\bar\zeta'\wedge d\bar\zeta',
\end{equation}
where \begin{itemize}\item $f_j$ is  a  bounded $2$-form  in $\{dw,dt\}$  and their complex  conjugates;\item  $g_j$ is a   $2$-form  in $\{dw,d\zeta', dt\}$ and their complex conjugates with coefficients  of order $O(t),$ and there is  no term of the forms  $d\zeta'\wedge d\zeta'$ or  $d\bar\zeta'\wedge d\bar\zeta';$
       \item $h_{j1}$ and $h_{j2}$  are functions  of ordet  $O(t^2).$
      \end{itemize}

Taking  \eqref{e:widetildeI_j-I_j} for  granted  we  are in the  position to prove assertion (1).

\noindent {\bf Proof of the first  inequality of assertion (1).}
 Using
\eqref{e:admissible-alpha(1)}--\eqref{e:admissible-alpha(2)} and \eqref{e:widetildeI_j-I_j}, we  get that
\begin{equation}  \label{e:widetildeI-I}(\tilde\tau_\ell)^*\alpha -\alpha=f+g,
\end{equation}  
where  $f$ is  a  bounded $2$-form  in $\{dw,dt\}$  and their complex  conjugates,  $g$ is a   $2$-form  in $\{dw,d\zeta',dt\}$ and their complex conjugates with coefficients  of order $O(t).$   
Applying the Cauchy--Schwarz inequality, we  infer from  the last equality that
\begin{equation*}
  -(\omega(w)+ |t|\omega'(\zeta')+  idt\wedge  d\bar t)\lesssim [(\tilde\tau_\ell)^*\alpha -\alpha]^\sharp \lesssim  \omega(w)+   |t|\omega'(\zeta')+  idt\wedge  d\bar t,
\end{equation*}
where we recall that  $\omega(z):=\sum_{p=1}^{k-l} idz_p\wedge d\bar z_p$ and $\omega(w):=\sum_{q=1}^{l} idw_q\wedge d\bar w_q.$
On the one hand, we have  
$$idt\wedge d\bar t\leq \omega(z)\lesssim \hat \beta.$$  
On the other hand, since   $|t|^2\approx \|z\|^2\approx\varphi(z,w),$ we deduce from  \eqref{e:omega'-zeta'}--\eqref{e:omega_FS-vs-omega'} and \eqref{e:tilde-alpha-beta-local-exp} and \eqref{e:hat-alpha'}--\eqref{e:hat-alpha'-vs-alpha_ver}  and \eqref{e:tilde-alpha-vs-hat-beta} that
$$
|t|\omega'(\zeta')\approx \varphi^{1\over 2}\omega_\FS([z])\approx \varphi^{1\over 2}\alpha_\ver\lesssim  \varphi^{1\over 2}\alpha+ c_1\varphi^{1\over 2}\pi^*\omega.
$$
Putting the last three estimates together,  the first inequality of assertion (1) follows.

\noindent {\bf Proof of the second  inequality of assertion (1).}
Let $\Gamma:=(\tilde\tau_\ell)^*\alpha -\alpha.$
We infer  from    \eqref{e:widetildeI-I} that
$$
\Gamma^{2,0}\wedge  \overline{\Gamma^{2,0}}=F+G,
$$
where $F$ is  a  bounded $(2,2)$-form  in $\{dw,dt\}$  and their complex  conjugates,  $G$ is a   $(2,2)$-form  in
$\{dw,d\zeta',dt\}$ and their complex conjugates with coefficients  of order $O(t).$
Applying  Cauchy--Schwarz inequality, we  infer from  the last equality that
\begin{equation*}
   \Gamma^{2,0}\wedge  \overline{\Gamma^{2,0}} \lesssim  \big(\omega(w)+   |t|^{1\over 2}\omega'(\zeta')+  idt\wedge  d\bar t\big)^2.
\end{equation*}
We obtain  in the  same way    as  in the  proof  of the  first inequality  of assertion (1)  that
$$ \Gamma^{0,2}\wedge  \overline{\Gamma^{0,2}}\lesssim (2c_1 \pi^*\omega +\varphi^{1\over 4} \alpha+\beta)^2.$$
 This  implies the second inequality of assertion (1).

\noindent {\bf Proof of the first  inequality of assertion (2).}
We  deduce from  \eqref{e:hat-alpha}  that \begin{equation}\label{e:tau-hat-alpha-minus-alpha}
\tilde\tau^*_\ell (\hat\alpha)-\hat\alpha= c_1[\tilde\tau^*_\ell ( \pi^*\omega) -\pi^*\omega ]+ [\tilde\tau^*_\ell (\alpha)-\alpha]+c_2
[\tilde\tau^*_\ell(\beta)-\beta].
\end{equation}
This, combined  with  the first inequalities of  Proposition \ref{P:basic-admissible-estimates-I} (2)--(3) and Proposition  \ref{P:basic-admissible-estimates-II} (1), yields  the desired result.

\noindent {\bf Proof of  the  second  inequality of assertion (2).}
Let $\Gamma:=(\tilde\tau_\ell)^*(\hat \alpha) -\hat\alpha.$
We infer  from    \eqref{e:tau-hat-alpha-minus-alpha},   \eqref{e:widetildeI-I}  and  the proof of  Proposition \ref{P:basic-admissible-estimates-I} (2)--(3)  that
$$
\Gamma^{2,0}\wedge  \overline{\Gamma^{2,0}}=F+G,
$$
where $F$ is  a  bounded $(2,2)$-form  in $\{dw,dt\}$  and their complex  conjugates,  $G$ is a   $(2,2)$-form  in
$\{dw,d\zeta',dt\}$ and their complex conjugates with coefficients  of order $O(t).$ Finally,
 we  proceed  as  in  the proof of  the second inequality of assertion (1).
\endproof

\noindent {\bf End of the proof of Proposition \ref{P:basic-admissible-estimates-II}.}
It remains to prove \eqref{e:widetildeI_j-I_j} for all $1\leq j\leq 4.$
To  prove \eqref{e:widetildeI_j-I_j} for $j=1,$   
observe  that 
 $\ddcw\log{\| A(w)(\zeta',1)\|^2}]$ is  a smooth form of bidegree $(1,1)$ in $(dw,d\bar w)$
and of bidegree $(0,0)$ in $(d\zeta,d\bar \zeta).$ 
Next,  since $\tilde\tau_\ell$ is strongly  admissible  and $|t|\approx \|z\|,$ we see  that 
\begin{equation}\label{e:tau-(z,w)}
 (\tilde \zeta ,\tilde w)-(\zeta,w)=O(z)=O(t),
\end{equation}
where we recall that $(\tilde z,\tilde w)=\tau(z,w)$  and  $\zeta$ (resp. $\tilde \zeta$) is the homogeneous
coordinate of $z$ (resp. $\tilde z$) according   to  \eqref{e:homogeneous-coordinates}. Moreover,
$$\tilde\tau_\ell^*(dw_q)- dw_q =\widetilde O(\|z\|)= \sum_{q'=1}^l(O(t) dw_{q'} + O(t) d\bar w_{q'})+  \sum_{p=1}^l( O(1) dz_{p} + O(1) d\bar z_p) .$$
We have, for $1\leq p\leq k-l-1,$ that 
\begin{equation}\label{e:dz_p-vs-zeta'-and-t}  dz_p=d(\zeta_p t)=\zeta_pdt +td\zeta_p
 =O(t)
\end{equation}
 because $|\zeta_p|\leq 2|t|$  by  \eqref{e:max-coordinate}.
Consequently,
\begin{equation}\label{e:tau-dw_q}
\tilde\tau_\ell^*(dw_q)- dw_q = O(t) dw +O(t)d\bar w+ O(t) d\zeta'+O(t) dt+O(t) d\bar t+O(t^2) d\bar\zeta' .
\end{equation}
We have  the same  expression for $
\tilde\tau_\ell^*(d\bar w_q)- d\bar w_q.$  

Combining  estimates \eqref{e:tau-(z,w)} and \eqref{e:tau-dw_q} and applying Lemma \ref{L:difference}, we  infer from the equality  
$$
\widetilde I_1-I_1=(\tilde\tau_\ell)^*[\ddcw\log{\| A(w)(\zeta',1)\|^2}] -[\ddcw\log{\| A(w)(\zeta',1)\|^2}]
$$
  that 
\eqref{e:widetildeI_j-I_j} holds for $j=1.$

The following two lemmas are needed.
\begin{lemma}\label{L:tilde-tau-dzeta_j}  For $1\leq j\leq   k-l-1,$ we have
 $$d\tilde\zeta_j-d\zeta_j=(\tilde\tau^*_\ell)(d\zeta_j) -d\zeta_j=O(t) dw +O(t)d\bar w+ O(t) d\zeta'+O(t) dt+O(t) d\bar t+O(t^2) d\bar\zeta'.$$
 Moreover,  $d\tilde t-dt =O(t) dw +O(t)d\bar w+ O(t) d\zeta'+O(t) dt+O(t) d\bar t+O(t^2) d\bar\zeta'.$
\end{lemma}
\proof
Since by  \eqref{e:homogeneous-coordinates} $\zeta_j={z_j\over  t}$ and $\tilde\tau_\ell$ is  strongly  admissible,  it 
follows from  \eqref{e:diff-dz} that
\begin{equation*}
 d\tilde\zeta_j-d\zeta_j= {(t+O(t^2))(dz_j+\widetilde O(t^2))-(z_j+O(t^2))(dt+\widetilde O(t^2))\over (t+O(t^2))^2}- {tdz_j-z_jdt\over t^2},
\end{equation*}
where the notation $\widetilde O$ is introduced in Definition \ref{D:admissible-maps}.
This, combined  with   \eqref{e:dz_p-vs-zeta'-and-t}, implies  the first  estimate of the lemma.

The  second  estimate  follows  from the  second identity  in  \eqref{e:z-t} below.
\endproof

Combining  estimates \eqref{e:tau-(z,w)}, \eqref{e:tau-dw_q} and applying Lemmas
\ref{L:difference}  and \ref{L:tilde-tau-dzeta_j}, we  infer from the equality  
$$
\widetilde I_2-I_2=(\tilde\tau_\ell)^*[\partial_{w}\dbar_\zeta\log{\| A(w)(\zeta',1)\|^2}] -[\partial_{w}\dbar_\zeta\log{\| A(w)(\zeta',1)\|^2}]
$$
  that  \eqref{e:widetildeI_j-I_j} holds for $j=2.$
  
Similarly, we can show that  \eqref{e:widetildeI_j-I_j}  also holds for $j=3.$

It remains to us   to show  that  \eqref{e:widetildeI_j-I_j}   holds for $j=4.$ Write
\begin{equation}\label{e:widetildeI_4-I_4}
\widetilde {I}_4-I_4=  I'_4+ I''_4+I'''_4, 
\end{equation}
where,   recalling that $(\tilde z,\tilde w)=\tilde\tau_\ell(z,w),$ we have 
\begin{eqnarray*}
 I'_4&:=& (\tilde\tau_\ell)^*[\ddczeta \log{\|A(w)(z)\|}] -\ddczeta [ (\tilde\tau_\ell)^*(  \log{\|A(w)( z)\|^2})],\\
 I''_4&:=&  \ddczeta [ (\tilde\tau_\ell)^*(  \log{\|A(w)( z)\|^2})] -\ddczeta [\log{\|A(w)(\tilde z)\|^2}],\\
 I'''_4&:=&\ddczeta[\log {\|A(w) (\tilde z)\|^2}]  -  \ddczeta  (\log {\|A(w)( z)\|^2}).
\end{eqnarray*}
By Lemma \ref{L:tilde-tau-dzeta_j}, we get \begin{equation}\label{e:I'_4} I'_4=O(t)\quad\text{and}\quad
 (I'_4)^{2,0}= O(t) dw\wedge dw + O(t) dw\wedge d\zeta'+O(t) dw \wedge dt+ O(t) d\zeta'\wedge dt +  O(t^2)d\zeta'\wedge d\zeta'
 .
                            \end{equation}
By Lemma \ref{L:difference}, $I''_4$  contains  the terms  which are of order $\tau^*A(w)-A(w).$
 By \eqref{e:tau-(z,w)},  we  conclude that $I''_4=O(t).$
Consequently, in order to  prove  \eqref{e:widetildeI_j-I_j}   for $j=4,$ we only need to show that
\begin{equation}\label{e:admissible-alpha(4)}
 -(\omega(w)+ |t|\omega'(\zeta')+  idt\wedge  d\bar t)\lesssim  (I'''_4)^\sharp
 \lesssim  \omega(w)+   |t|\omega'(\zeta')+  idt\wedge  d\bar t.
\end{equation}
Write
\begin{equation}\label{e:I'''_4}
 \begin{split}
 I'''_4&=\big(\ddczeta [\log {\|A(w) (\tilde \zeta',1)\|^2}]-\ddczeta  [\log {\|A(w)( \zeta',1)\|^2}]\big)+
 \big(\ddczeta [(\log {|\tilde t |^2})]-\ddczeta  [\log {|t|^2}]\big)\\
 &:= I_{41}+I_{42}.
 \end{split}
\end{equation}
Therefore, inequality \eqref{e:admissible-alpha(4)}, and hence  inequality \eqref{e:widetildeI_j-I_j}  for $j=4,$ is a consequence of the following two lemmas.

\begin{lemma}\label{L:I_41}  The following estimate holds: $
I_{41}=O(t)+ O(1)dt\wedge d\bar t.
$ 
\end{lemma}
\proof
Recall  that $|\zeta_j|,|\zeta'_j|<3$ for  $1\leq j\leq k-l-1.$
 Therefore,  $\ddczetat  [\log {\|A( w)(\tilde \zeta',1)\|^2}]=\ddczetap  [\log {\|A( w)(\tilde \zeta',1)\|^2}] $ is a smooth  function.
Moreover, the matrix-valued  function $A(w):\ \overline{\D}^l\to \GL(\C,k-l)$ is  smooth.

 We will  prove the  following  two facts.  The first fact  says that for every $f$   among $\{t,\bar t,\zeta_1,\bar \zeta_1,\ldots,\zeta_{k-l-1}, \bar \zeta_{k-l-1}\}$
and  for  every  $D$ among $\{ \id, \partial_{\zeta'},\partial_t, \dbar_{\zeta'},\dbar_t, \ddczetap,\partial_{\zeta'}\dbar_t, \partial_t\dbar_{\zeta'} \},$
the  following  inequality holds:
\begin{equation}\label{e:admissible-zeta',t}
 D\tilde  f  -Df  =O(t)+O(1) dt +O(1)d\bar t.
\end{equation}
The  second  fact  says   that for  every $f$   among $\{t,\bar t,\zeta_1,\bar \zeta_1,\ldots,\zeta_{k-l-1}, \bar \zeta_{k-l-1}\},$
the  following  inequality holds:
\begin{equation}\label{e:admissible-ddc_t}
 \ddct\tilde  f  -\ddct f  =O(1).
\end{equation}
Assuming  \eqref{e:admissible-zeta',t}--\eqref{e:admissible-ddc_t}  we  resume the proof of the lemma. First  we  apply  the   equality
$$
\ddc\log\phi={1\over  \phi} \ddc\phi-{i\over \pi\phi^2}  \partial \phi\wedge \dbar\phi
$$
to  $\phi:=\| A(w) (\zeta',1)\|^2$  and then to $\phi:=\|A(w)(\tilde\zeta',1)\|^2.$
Next, using  \eqref{e:admissible-zeta',t}--\eqref{e:admissible-ddc_t} and the identities
$$
\ddc={i\over \pi}\ddbar\quad\text{and}\quad  \partial_{\zeta',t}=\partial_{\zeta'}+\partial_t\quad\text{and}\quad
\dbar_{\zeta',t}=\dbar_{\zeta'}+\dbar_t,
$$
we    apply
 Lemma \ref{L:difference}
 to  $I_{41}.$
Consequently, we  get the  desired conclusion of the lemma.

It remains  to prove  \eqref{e:admissible-zeta',t}--\eqref{e:admissible-ddc_t}.  We use the  homogeneous  coordinates $\zeta,$ $\tilde \zeta,$
and recall that $t=\zeta_{k-l}=z_{k_l}$  and 
$\tilde t=\tilde\zeta_{k-l}=\tilde z_{k_l}.$  Since  $\tilde\tau_\ell$ is   strongly    admissible,
it follows from  Definition  \ref{D:Strongly-admissible-maps} that we can  write
\begin{equation}\label{e:z-t}\begin{split}
\tilde z_j&=t\zeta_j+   a_jt^2+ \sum_{p=1}^{k-l-1} a_{jp} t^2\zeta_p  +\sum_{p,q=1}^{k-l-1}a_{jpq}t^2\zeta_p \zeta_{q} +O(t^3),\\
\tilde t&= t+ a'_0t^2+  \sum_{p=1}^{k-l-1} a'_{p} t^2  \zeta_p+\sum_{p,q=1}^{k-l-1}a'_{pq}t^2\zeta_p \zeta_{q}  +O(t^3).
\end{split} \end{equation}
Here,  $ 1\leq j\leq k-l-1$ and  $a_j,a_{jp},a_{jpq},a'_0,a'_p,a'_{pq}\in\C$  are some constants.
 Therefore,  
 \begin{equation*}
  \tilde\zeta_j  -\zeta_j={\tilde z_j\over  \tilde t}-\zeta_j= b_jt+ \sum_{p=1}^{k-l-1} b_{jp} t  \zeta_p+ \sum_{p,q=1}^{k-l-1}b_{jpq}t\zeta_p \zeta_{q}  +{O(t^3)\over  t},
 \end{equation*}
where  $b_j,b_{jp},b_{jpq} \in\C$  are some constants.
  Note that  $O(t^3)$ is  a $\Cc^2$-function in $\zeta',t$ and hence the function
${O(t^3)\over t}$ is of class $\Cc^{1,1}.$
Using this explicit   expression and Lemma \ref{L:difference}  we  can check \eqref{e:admissible-zeta',t}--\eqref{e:admissible-ddc_t}.
The proof of Lemma \ref{L:I_41} is thereby completed.
 \endproof

\begin{lemma}\label{L:I_42}  $I_{42}$ is  a bounded $(1,1)$-form in $dt,$ $d\bar t.$
\end{lemma}
\proof By hypothesis,   $\tau$ is  strongly   admissible. Therefore,  $\tilde\tau_\ell$ is  also  strongly    admissible. Hence,  we can  write, by the second equality of \eqref{e:z-t},
\begin{equation*}
\tilde t= t+ t\big(a_0t+\sum_{p=1}^{k-l-1}a_p\zeta_p \big)+O(t^3)
\end{equation*}
for some  constants $a_1,\ldots,a_{k-l-1}\in\C.$   So 
\begin{equation*}
 I_{42}=\ddct \log{\big|{\tilde t\over t}\big|^2}=2\ddct\log{\big| 1+ \big(a_0t+\sum\limits_{p=1}^{k-l-1}a_j\zeta_j \big)+{O(t^3)\over t}\big|}.
\end{equation*}
We have the classical  Taylor expansion $\log{|1+\xi|}=\Re\big( \sum_{n=1}^\infty (-1)^{n-1}{\xi^n\over n}  \big)$  for $\xi$ close to   the point $0\in\C$ and the function
${O(t^3)\over t}$ is of class $\Cc^{1,1}.$  Therefore,  we infer  that the form $I_{42}$ is  a bounded $(1,1)$-form in $dt,$ $d\bar t.$ This completes the proof.   
\endproof

Combining  equality \eqref{e:widetildeI_4-I_4}, estimate \eqref{e:I'_4}  and the fact that  $I''_4=O(t),$
equality \eqref{e:I'''_4}, Lemma \ref{L:I_41}  and \ref{L:I_42} we see that
\eqref{e:widetildeI_j-I_j} holds for $j=4.$
This completes  the proof of   Proposition \ref{P:basic-admissible-estimates-II}.
\endproof

\begin{remark}
 \rm  In both  Lemmas \ref{L:I_41} and \ref{L:I_42} we have made full use of the  assumption that
 $\tau$ is   strongly  admissible.
\end{remark}

The following  notion generalizes Definition
 \ref{D:precsim}  to   a collection of finite  $2$-forms. It will be  needed  in order  to obtain   admissible estimates.
\begin{definition}
 \label{D:precsim-bis}
 \rm  Let $\Gamma:=(\Gamma_1,\ldots, \Gamma_n)$ be a collection of  $n$ forms of degree $2$ and $S:=(S_1,\ldots,S_n)$ be  a collection of $n$ positive $(1,1)$-forms  defined on $\U.$ 
 We  write $\Gamma\trianglelefteq  S$  
 if there is a constant $c>0$ such that   the  following   two inequalities  hold for all $1\leq p,q\leq n$ and  $y\in\U:$
 $$
  \pm\Re\big[\Gamma^{2,0}_p(y)\wedge \overline{\Gamma^{2,0}_q}(y)\big]\leq c (S_p\wedge  S_q)(y)\ \quad\text{and}\quad \pm\Im\big[ \Gamma^{2,0}_p(y)\wedge \overline{\Gamma^{2,0}_q}(y)\big]\leq c (S_p\wedge  S_q)(y).
 $$
 Here,
 $\Gamma^{2,0}_p$  (resp.    $\Gamma^{0,2}_p$)  denotes the component of bidegree $(2,0)$ (resp.  $(0,2)$) of $\Gamma_p.$
\end{definition}

\begin{theorem}
  \label{T:basic-admissible-estimates}
There  are constants $c_3,c_4>0$ such that   $ c_3 \pi^*\omega +c_4\beta\geq 0$  on  $\pi^{-1}(V_0)\subset \E$  and that
for every $1\leq \ell \leq \ell_0,$
 $\Gamma \trianglelefteq  S$ on $\U_\ell \cap \Tube(B,\bfr),$
where  
\begin{eqnarray*}
\Gamma&:=& (\Gamma_1,\Gamma_2,\Gamma_3)\quad\text{and}\quad S:=(S_1,S_2,S_3),\\
 \Gamma_1&:=& \tilde\tau_\ell^*(\pi^*\omega)  -\pi^*\omega  \quad\text{and}\quad S_1:=   c_3 \varphi^{1\over 2}  \pi^*\omega +c_4\varphi^{1\over 2} \beta,\\
 \Gamma_2&:=&    \tilde \tau_\ell^*(\hat\beta) -\hat\beta   \quad\text{and}\quad S_2:= c_3 \phi^{3\over 2}\cdot  \pi^*\omega +c_4\phi^{1\over 2}\cdot \hat\beta,\\
  \Gamma_3&:=&\tilde \tau_\ell^*(\hat\alpha) -\hat\alpha  \quad\text{and}\quad S_3:=  c_3 \pi^*\omega +c_4\hat\beta +c_3\varphi^{1/4}\hat\alpha.
\end{eqnarray*}
 
\end{theorem}
\proof
We will  express  the  forms in terms of $(w,\zeta,t)$ and $\{dw,d\bar w,d\zeta', d\bar \zeta', dt, d\bar t\}.$
By estimates \eqref{e:Gamma-0,2} and \eqref{e:Gamma_2,0},  we get  that 
\begin{eqnarray*}
\Gamma_1^{2,0}&=&O(t)dw\wedge dw+O(t)dw\wedge dt +O(t^2)dw\wedge d\zeta'+O(t^2)dt\wedge d\zeta'+O(t^3)d\zeta'\wedge d\zeta',\\
 \Gamma_2^{2,0}&=&O(t^3) dw\wedge dw+O(t^2) dw\wedge dt+  O(t^3) dw\wedge d\zeta'+ O(t^3) dt\wedge d\zeta'+O(t^4)d\zeta'\wedge d\zeta'.
 \end{eqnarray*}
 Moreover, by  \eqref{e:widetildeI-I} and \eqref{e:tau-hat-alpha-minus-alpha},
 we obtain  that
 \begin{equation*}
 \Gamma_3^{2,0}=O(1) d w\wedge d w+O(1) d w\wedge d t+O(t)dw\wedge d\zeta'+O(t)dt\wedge d\zeta'+O(t^2)d\zeta'\wedge d\zeta'.
\end{equation*} 
On the other hand, combining 
\eqref{e:omega'-zeta'} and \eqref{e:alpha-beta-local}
yields
\begin{eqnarray*}
S_1&\gtrsim& |t|\omega(w)+|t|^3\omega'(\zeta')+i|t| dt\wedge d\bar t,\\
 S_2&\gtrsim& | t|^3\omega(w)+|t|^3\omega'(\zeta')+|t|i dt\wedge d\bar t.
 \end{eqnarray*}
 Moreover, combining estimates \eqref{e:hat-alpha-vs-alpha_ver}, \eqref{e:hat-alpha'-vs-alpha_ver}
 and \eqref{e:omega_FS-vs-omega'} 
 yields
 \begin{equation*}
 S_3\gtrsim\omega(w)+|t|^{1\over 2}\omega'(\zeta')+ idt\wedge d\bar t.
\end{equation*}
We will prove that  for $1\leq p\leq q\leq 3,$ the following  inequalities hold:
\begin{equation}\label{e:admissible-estimates-p,q}
  \pm\Re\big[\Gamma^{2,0}_p(y)\wedge \overline{\Gamma^{2,0}_q}(y)\big]\leq c (S_p\wedge  S_q)(y)\ \quad\text{and}\quad \pm\Im\big[ \Gamma^{2,0}_p(y)\wedge \overline{\Gamma^{2,0}_q}(y)\big]\leq c (S_p\wedge  S_q)(y).
 \end{equation}
 By Proposition \ref{P:basic-admissible-estimates-II},  inequalities 
\eqref{e:admissible-estimates-p,q} hold for $p=q.$ Therefore, we only need to  prove \eqref{e:admissible-estimates-p,q}  for $(p,q)\in\{(1,2),(1,3),(2,3)\}.$

 \noindent {\bf Proof of  \eqref{e:admissible-estimates-p,q} for $(p,q)=(1,2).$} Using the above  estimates, we  see that
\begin{eqnarray*}
 \Gamma_1^{2,0}\wedge \Gamma_2^{0,2}
 &=& O(t^4) dw\wedge dw\wedge d\bar w\wedge d\bar w+O(t^3)  dw\wedge dw\wedge d\bar w\wedge d\bar t\\
 &+&O(t^4)  dw\wedge dw\wedge d\bar w\wedge d\bar \zeta'
 +O(t^4)  dw\wedge dw\wedge d\bar t\wedge d\bar \zeta'+ O(t^5)  dw\wedge dw\wedge d\bar \zeta'\wedge d\bar \zeta'
 \\
 &+&     O(t^4) dw \wedge dt\wedge d\bar w\wedge d\bar w+      O(t^3) dw \wedge dt\wedge d\bar w\wedge d\bar t\\
 &+& O(t^4)  dw \wedge dt\wedge  d\bar w\wedge d\bar \zeta' +O(t^4) dw\wedge dt \wedge d\bar t\wedge d\bar \zeta'+  O(t^5) 
 dw\wedge dt \wedge d\bar \zeta'\wedge d\bar \zeta'\\
 &+&O(t^5)dw\wedge d\zeta'\wedge d\bar w\wedge d\bar w+O(t^4)dw\wedge d\zeta'\wedge d\bar w\wedge d\bar t\\
 &+& O(t^5)dw\wedge d\zeta'\wedge d\bar w\wedge d\bar \zeta'+O(t^5)dw\wedge d\zeta'\wedge d\bar t\wedge d\bar \zeta'+O(t^6)dw\wedge d\zeta'\wedge d\bar \zeta'\wedge d\bar \zeta'\\
 &+& O(t^5)dt\wedge d\zeta'\wedge d\bar w\wedge d\bar w+O(t^4)dt\wedge d\zeta'\wedge d\bar w\wedge d\bar t\\
 &+& O(t^5)dt\wedge d\zeta'\wedge d\bar w\wedge d\bar \zeta'+O(t^5)dt\wedge d\zeta'\wedge d\bar t\wedge d\bar \zeta'+O(t^5)dt\wedge d\zeta'\wedge 
 d\bar \zeta'\wedge d\bar \zeta'\\
 &+& O(t^6)d\zeta'\wedge d\zeta'\wedge d\bar w\wedge d\bar w+O(t^5)d\zeta'\wedge d\zeta'\wedge d\bar w\wedge d\bar t\\
 &+& O(t^6)d\zeta'\wedge d\zeta'\wedge d\bar w\wedge d\bar \zeta'+O(t^6)d\zeta'\wedge d\zeta'\wedge d\bar t\wedge d\bar \zeta'+O(t^7)d\zeta'\wedge d\zeta'\wedge d\bar \zeta'\wedge d\bar \zeta'.
\end{eqnarray*}
 Moreover, we also have that
\begin{eqnarray*}
  S_1\wedge S_2&\gtrsim & |t|^4\omega(w)^2+ |t|^2i dt\wedge d\bar t\wedge  \omega(w)+|t|^4\omega(w)\wedge \omega'(\zeta') \\
 &+& |t|^4i dt\wedge d\bar t\wedge  \omega'(\zeta')+|t|^6 \omega'(\zeta')^2.
\end{eqnarray*}
Now we  treat the two terms  with exact order $O(t^3)$ on the RHS of the  expression for $\Gamma_1^{2,0}\wedge \Gamma_2^{0,2}.$
Applying the Cauchy--Schwarz inequality yields that
\begin{eqnarray*}
 O(t^3) dw \wedge dw\wedge d\bar w\wedge d\bar t&\lesssim&   |t|^3i (|t|\omega(w)+|t|^{-1}dt\wedge d\bar t) \wedge  \omega(w)\lesssim S_1\wedge S_2,\\
 O(t^3) dw \wedge dt\wedge d\bar w\wedge d\bar t&\lesssim&   |t|^3i dt\wedge d\bar t\wedge  \omega(w)\lesssim S_1\wedge S_2.
\end{eqnarray*}
Next, we  treat all  terms $I$  with order at least $O(t^4)$ on the RHS of the  expression for $\Gamma_1^{2,0}\wedge \Gamma_2^{0,2}$ such that $I$  contains
neither the factor $d\zeta'\wedge d\zeta'$ nor the factor  $d\bar\zeta'\wedge d\bar\zeta'.$ There are 14  such terms.
Applying the Cauchy--Schwarz inequality yields that
\begin{equation*}
 I\lesssim |t|^4\omega(w)^2+ |t|^2i dt\wedge d\bar t\wedge  \omega(w)+|t|^4\omega(w)\wedge \omega'(\zeta') 
 + |t|^4i dt\wedge d\bar t\wedge  \omega'(\zeta')\lesssim S_1\wedge S_2.
\end{equation*}
Next, we  treat  all terms  $I$  among the 9 remaining  terms on the RHS of the  expression for $\Gamma_1^{2,0}\wedge \Gamma_2^{0,2}$
such that $ I$ is  of order at least $O(t^6).$ There are 5 such terms.
Applying the Cauchy--Schwarz inequality yields that
\begin{equation*}
 I\lesssim |t|^6\omega(w)^2+ |t|^6i dt\wedge d\bar t\wedge  \omega(w)+|t|^6\omega(w)\wedge \omega'(\zeta') 
 + |t|^6i dt\wedge d\bar t\wedge  \omega'(\zeta')+ |t|^6 \omega'(\zeta')^2 \lesssim S_1\wedge S_2.
\end{equation*}
Finally, we treat the last 4 terms. They are all of order $O(t^3).$ Applying the Cauchy--Schwarz inequality yields that
\begin{eqnarray*}
O(t^5) dw\wedge  dw \wedge  d\bar\zeta' \wedge d\bar\zeta'&\lesssim& |t|^5( |t|\omega'(\zeta')+ |t|^{-1} \omega(w) ) \wedge ( \omega'(\zeta')+  \omega(w) ) \lesssim S_1\wedge S_2,\\
O(t^5) dw\wedge  dt \wedge  d\bar\zeta' \wedge d\bar\zeta'&\lesssim& |t|^5( |t|\omega'(\zeta')+ |t|^{-1} \omega(w) ) \wedge ( \omega'(\zeta')+  idt\wedge d\bar t ) \lesssim S_1\wedge S_2,\\
O(t^5) dt\wedge  d\zeta' \wedge  d\bar\zeta' \wedge d\bar\zeta'&\lesssim& |t|^5 \omega'(\zeta') \wedge ( |t|\omega'(\zeta')+ |t|^{-1} idt\wedge d\bar t ) \lesssim S_1\wedge S_2,\\
 O(t^5)  d\zeta' \wedge d\zeta'\wedge d\bar w\wedge d\bar t&\lesssim& |t|^5( |t|\omega'(\zeta')+ |t|^{-1} \omega(w) ) \wedge (i dt\wedge d\bar t+  \omega'(\zeta'))\lesssim S_1\wedge S_2.
\end{eqnarray*}
In summary, we have shown
inequality  \eqref{e:admissible-estimates-p,q} for $(p,q)=(1,2).$

 \noindent {\bf Proof of  \eqref{e:admissible-estimates-p,q} for $(p,q)=(1,3).$} Using the above  estimates, we  see that
\begin{eqnarray*}
 \Gamma_1^{2,0}\wedge \Gamma_3^{0,2}
 &=& O(t) dw\wedge dw\wedge d\bar w\wedge d\bar w+O(t)  dw\wedge dw\wedge d\bar w\wedge d\bar t\\
 &+&O(t^2)  dw\wedge dw\wedge d\bar w\wedge d\bar \zeta'
 +O(t^2)  dw\wedge dw\wedge d\bar t\wedge d\bar \zeta'+ O(t^3)  dw\wedge dw\wedge d\bar \zeta'\wedge d\bar \zeta'
 \\
 &+&     O(t) dw \wedge dt\wedge d\bar w\wedge d\bar w+      O(t) dw \wedge dt\wedge d\bar w\wedge d\bar t\\
 &+& O(t^2)  dw \wedge dt\wedge  d\bar w\wedge d\bar \zeta' +O(t^2) dw\wedge dt \wedge d\bar t\wedge d\bar \zeta'+  O(t^3) 
 dw\wedge dt \wedge d\bar \zeta'\wedge d\bar \zeta'\\
 &+&O(t^2)dw\wedge d\zeta'\wedge d\bar w\wedge d\bar w+O(t^2)dw\wedge d\zeta'\wedge d\bar w\wedge d\bar t\\
 &+& O(t^3)dw\wedge d\zeta'\wedge d\bar w\wedge d\bar \zeta'+O(t^3)dw\wedge d\zeta'\wedge d\bar t\wedge d\bar \zeta'+O(t^4)dw\wedge d\zeta'\wedge d\bar \zeta'\wedge d\bar \zeta'\\
 &+& O(t^2)dt\wedge d\zeta'\wedge d\bar w\wedge d\bar w+O(t^2)dt\wedge d\zeta'\wedge d\bar w\wedge d\bar t\\
 &+& O(t^3)dt\wedge d\zeta'\wedge d\bar w\wedge d\bar \zeta'+O(t^3)dt\wedge d\zeta'\wedge d\bar t\wedge d\bar \zeta'+O(t^4)dt\wedge d\zeta'\wedge 
 d\bar \zeta'\wedge d\bar \zeta'\\
 &+& O(t^3)d\zeta'\wedge d\zeta'\wedge d\bar w\wedge d\bar w+O(t^3)d\zeta'\wedge d\zeta'\wedge d\bar w\wedge d\bar t\\
 &+& O(t^4)d\zeta'\wedge d\zeta'\wedge d\bar w\wedge d\bar \zeta'+O(t^4)d\zeta'\wedge d\zeta'\wedge d\bar t\wedge d\bar \zeta'+O(t^5)d\zeta'\wedge d\zeta'\wedge d\bar \zeta'\wedge d\bar \zeta'.
\end{eqnarray*}
 Moreover, we also have that
\begin{eqnarray*}
  S_1\wedge S_3&\gtrsim & |t|\omega(w)^2+ |t|i dt\wedge d\bar t\wedge  \omega(w)+|t|^{3\over 2}\omega(w)\wedge \omega'(\zeta') \\
 &+& |t|^{3\over 2}i dt\wedge d\bar t\wedge  \omega'(\zeta')+|t|^{7\over 2} \omega'(\zeta')^2.
\end{eqnarray*}
Now we  treat the four terms  with exact order $O(t)$ on the RHS of the  expression for $\Gamma_1^{2,0}\wedge \Gamma_3^{0,2}.$
Applying the Cauchy--Schwarz inequality yields that
\begin{eqnarray*}
 O(t) dw \wedge dw\wedge d\bar w\wedge d\bar w&\lesssim&   |t| \omega(w)^2 \lesssim S_1\wedge S_3,\\
 O(t) dw \wedge dw\wedge d\bar w\wedge d\bar t&\lesssim&   |t|   (\omega(w)+ idt\wedge d\bar t)\wedge  \omega(w)\lesssim S_1\wedge S_3,\\
 O(t) dw \wedge dt\wedge d\bar w\wedge d\bar w&\lesssim&   |t|   (\omega(w)+ idt\wedge d\bar t)\wedge  \omega(w)\lesssim S_1\wedge S_3,\\
 O(t) dw \wedge dt\wedge d\bar w\wedge d\bar t&\lesssim&   |t|   ( idt\wedge d\bar t)\wedge  \omega(w)\lesssim S_1\wedge S_3.
\end{eqnarray*}
Next, we  treat all  terms $I$  with order at least $O(t^{3\over 2})$ on the RHS of the  expression for $\Gamma_1^{2,0}\wedge \Gamma_3^{0,2}$ such that $I$  contains
neither the factor $d\zeta'\wedge d\zeta'$ nor the factor  $d\bar\zeta'\wedge d\bar\zeta'.$ There are 12  such terms.
Applying the Cauchy--Schwarz inequality yields that
\begin{equation*}
 I\lesssim |t|^{3\over 2}\omega(w)^2+ |t|^{3\over 2}i dt\wedge d\bar t\wedge  \omega(w)+|t|^{3\over 2}\omega(w)\wedge \omega'(\zeta') 
 + |t|^{3\over 2}i dt\wedge d\bar t\wedge  \omega'(\zeta')\lesssim S_1\wedge S_3.
\end{equation*}
Next, we  treat  all terms  $I$  among the 9 remaining  terms on the RHS of the  expression for $\Gamma_1^{2,0}\wedge \Gamma_3^{0,2}$
such that $ I$ is  of order at least $O(t^4).$ There are 5 such terms.
Applying the Cauchy--Schwarz inequality yields that
\begin{equation*}
 I\lesssim |t|^4\omega(w)^2+ |t|^4i dt\wedge d\bar t\wedge  \omega(w)+|t|^4\omega(w)\wedge \omega'(\zeta') 
 + |t|^4i dt\wedge d\bar t\wedge  \omega'(\zeta')+ |t|^4 \omega'(\zeta')^2 \lesssim S_1\wedge S_3.
\end{equation*}
Finally, we treat the last 4 terms. They are all of order $O(t^3).$ Applying the Cauchy--Schwarz inequality yields that
\begin{eqnarray*}
O(t^3) dw\wedge  dw \wedge  d\bar\zeta' \wedge d\bar\zeta'&\lesssim& |t|^3( |t|\omega'(\zeta')+ |t|^{-1} \omega(w) ) \wedge ( \omega'(\zeta')+  \omega(w) ) \lesssim S_1\wedge S_3,\\
O(t^3) dw\wedge  dt \wedge  d\bar\zeta' \wedge d\bar\zeta'&\lesssim& |t|^3( |t|\omega'(\zeta')+ |t|^{-1} \omega(w) ) \wedge ( \omega'(\zeta')+  idt\wedge d\bar t ) \lesssim S_1\wedge S_3,\\
O(t^3) d\zeta'\wedge  d\zeta' \wedge  d\bar w \wedge d\bar w&\lesssim& |t|^3 (|t|\omega'(\zeta') +|t|^{-1}\omega(w))  \wedge ( \omega'(\zeta')+ \omega(w)) \lesssim S_1\wedge S_3,\\
 O(t^3)  d\zeta' \wedge d\zeta'\wedge d\bar w\wedge d\bar t&\lesssim& |t|^3( |t|\omega'(\zeta')+ |t|^{-1} \omega(w) ) \wedge (i dt\wedge d\bar t+  \omega'(\zeta'))\lesssim S_1\wedge S_3.
\end{eqnarray*}
In summary, we have shown
inequality  \eqref{e:admissible-estimates-p,q} for $(p,q)=(1,3).$

 \noindent {\bf Proof of  \eqref{e:admissible-estimates-p,q} for $(p,q)=(2,3).$} Using the above  estimates, we  see that
\begin{eqnarray*}
 \Gamma_2^{2,0}\wedge \Gamma_3^{0,2}&=& O(t^3) dw\wedge dw\wedge d\bar w\wedge d\bar w+O(t^3)  dw\wedge dw\wedge d\bar w\wedge d\bar t\\
 &+&O(t^4)  dw\wedge dw\wedge d\bar w\wedge d\bar \zeta'
 +O(t^4)  dw\wedge dw\wedge d\bar t\wedge d\bar \zeta'+ O(t^5)  dw\wedge dw\wedge d\bar \zeta'\wedge d\bar \zeta'
 \\
 &+&     O(t^2) dw \wedge dt\wedge d\bar w\wedge d\bar w+      O(t^2) dw \wedge dt\wedge d\bar w\wedge d\bar t\\
 &+& O(t^3)  dw \wedge dt\wedge  d\bar w\wedge d\bar \zeta' +O(t^3) dw\wedge dt \wedge d\bar t\wedge d\bar \zeta'+  O(t^4) 
 dw\wedge dt \wedge d\bar \zeta'\wedge d\bar \zeta'\\
 &+&O(t^3)dw\wedge d\zeta'\wedge d\bar w\wedge d\bar w+O(t^3)dw\wedge d\zeta'\wedge d\bar w\wedge d\bar t\\
 &+& O(t^4)dw\wedge d\zeta'\wedge d\bar w\wedge d\bar \zeta'+O(t^4)dw\wedge d\zeta'\wedge d\bar t\wedge d\bar \zeta'+O(t^5)dw\wedge d\zeta'\wedge d\bar \zeta'\wedge d\bar \zeta'\\
 &+& O(t^3)dt\wedge d\zeta'\wedge d\bar w\wedge d\bar w+O(t^3)dt\wedge d\zeta'\wedge d\bar w\wedge d\bar t\\
 &+& O(t^4)dt\wedge d\zeta'\wedge d\bar w\wedge d\bar \zeta'+O(t^4)dt\wedge d\zeta'\wedge d\bar t\wedge d\bar \zeta'+O(t^5)dt\wedge d\zeta'\wedge 
 d\bar \zeta'\wedge d\bar \zeta'\\
 &+& O(t^4)d\zeta'\wedge d\zeta'\wedge d\bar w\wedge d\bar w+O(t^4)d\zeta'\wedge d\zeta'\wedge d\bar w\wedge d\bar t\\
 &+& O(t^5)d\zeta'\wedge d\zeta'\wedge d\bar w\wedge d\bar \zeta'+O(t^5)d\zeta'\wedge d\zeta'\wedge d\bar t\wedge d\bar \zeta'+O(t^6)d\zeta'\wedge d\zeta'\wedge d\bar \zeta'\wedge d\bar \zeta'.
\end{eqnarray*}
Moreover, we also have that
\begin{eqnarray*}
  S_2\wedge S_3&\gtrsim & |t|^3\omega(w)^2+ |t|i dt\wedge d\bar t\wedge  \omega(w)+|t|^3\omega(w)\wedge \omega'(\zeta') \\
 &+& |t|^{3\over 2}i dt\wedge d\bar t\wedge  \omega'(\zeta')+|t|^{7\over 2} \omega'(\zeta')^2.
\end{eqnarray*}
Now we  treat the two terms  with exact order $O(t^2)$ on the RHS of the  expression for $\Gamma_2^{2,0}\wedge \Gamma_3^{0,2}.$
Applying the Cauchy--Schwarz inequality yields that
\begin{eqnarray*}
 O(t^2) dw \wedge dt\wedge d\bar w\wedge d\bar w&\lesssim& |t|^2\omega(w)\wedge(\omega(w)+idt\wedge d\bar t)=  |t|^2\omega(w)^2+ |t|^2i dt\wedge d\bar t\wedge  \omega(w)\lesssim S_2\wedge S_3,\\
 O(t^2)  dw \wedge dt\wedge d\bar w\wedge d\bar t&\lesssim& |t|^2i dt\wedge d\bar t\wedge  \omega(w)\leq
  |t|i dt\wedge d\bar t\wedge  \omega(w)\lesssim S_2\wedge S_3.
\end{eqnarray*}
 Next, we  treat all  terms $I$  with order at least $O(t^3)$ on the RHS of the  expression for $\Gamma_2^{2,0}\wedge \Gamma_3^{0,2}$ such that $I$  contains
neither the factor $d\zeta'\wedge d\zeta'$ nor the factor  $d\bar\zeta'\wedge d\bar\zeta'.$ There are 14  such terms.
Applying the Cauchy--Schwarz inequality yields that
\begin{equation*}
 I\lesssim |t|^{3}\omega(w)^2+ |t|^{3}i dt\wedge d\bar t\wedge  \omega(w)+|t|^{3}\omega(w)\wedge \omega'(\zeta') 
 + |t|^{3}i dt\wedge d\bar t\wedge  \omega'(\zeta')\lesssim S_2\wedge S_3.
\end{equation*}
Finally, we  treat   the 9 remaining  terms on the RHS of the  expression for $\Gamma_2^{2,0}\wedge \Gamma_3^{0,2}.$
Such a  term  $ I$ is  of order at least $O(t^4).$ 
Applying the Cauchy--Schwarz inequality yields that
\begin{equation*}
 I\lesssim |t|^4\omega(w)^2+ |t|^4i dt\wedge d\bar t\wedge  \omega(w)+|t|^4\omega(w)\wedge \omega'(\zeta') 
 + |t|^4i dt\wedge d\bar t\wedge  \omega'(\zeta')+ |t|^4 \omega'(\zeta')^2 \lesssim S_2\wedge S_3.
\end{equation*}
 In summary, we have shown
inequality  \eqref{e:admissible-estimates-p,q} for $(p,q)=(2,3).$   
\endproof

\subsection{Admissible  estimates for wedge-products}\label{SS:adm-est-wedge-prod}

Let $\U$ be an open neighborhood of $0$ in $\C^k.$
We use the local coordinates $y=(z,w)\in\C^{k-l}\times \C^l$ on $\U$
and recall  the notion of order   $\trianglelefteq  $
given in Definition
 \ref{D:precsim}.

\begin{lemma}\label{L:basic-positive-estimate}
 For every $1\leq j\leq q,$  $\Gamma_j$ and $S_j$ are  real currents of the same  bidegree $(p_j,p_j)$
 on $\U$
 such that $S_j$ is positive and $-S_j\leq  \Gamma_j\leq  S_j$ on $\U.$
 Then $$-c\cdot\,  S_1\wedge\ldots\wedge S_q\leq   \Gamma_1\wedge\ldots\wedge \Gamma_q \leq  c\cdot \,S_1\wedge\ldots\wedge S_q \quad\text{on}\quad \U .$$  
 Here $c$ is a constant that depends only on the dimension $k.$
\end{lemma}
\proof
We only need to prove the lemma for $q=2.$ The general case can be  proved by repeatedly applying  the case $q=2.$
Write
$$
\Gamma_1\wedge \Gamma_2=S_1\wedge S_2+(\Gamma_1-S_1)\wedge S_2+S_1\wedge (\Gamma_2-S_2)+(\Gamma_1-S_1)\wedge (\Gamma_2-S_2).
$$
Since for $j=1,2,$ we have  $-2S_j\leq  \Gamma_j-S_j\leq  0$ on $\U,$ it follows that
$$
-4 S_1\wedge S_2\leq \Gamma_1\wedge \Gamma_2\leq 5  S_1\wedge S_2.
$$
So  for $q=2$ the lemma  is true with $c=5.$
\endproof

\begin{lemma}\label{L:Hc}
  Let $ H_1,\ldots, H_q$ be  $q$ real $(1,1)$-forms in the class $\Hc=\Hc(\U)$  introduced  in Definition \ref{D:Hc}.
Then there is a constant $c>0$  such that
\begin{equation*}
\pm  H_1\wedge\ldots\wedge H_m\leq c \big( \sum_{j=0}^q  \pi^*\omega^j \wedge \hat\beta^{q-j} \big).  
\end{equation*}
\end{lemma}
\proof
By the Cauchy--Schwarz inequality, there is a constant $c>0$ such that
$\pm H_j\leq c (\pi^*\omega+\hat\beta)$  for all $1\leq  j\leq m.$
Using these inequalities and  applying  Lemma \ref{L:basic-positive-estimate}, the result follows.
\endproof

\begin{lemma}\label{L:Cauchy-Schwarz} {\rm (Cauchy--Schwarz inequality for wedge-products)}
 Let  $T$ be  a  positive   current of bidimension $(q,q)$  and  $\Gamma$ a  real  current of  bidimension $(q,q)$ on $\U$ such that
 $ -T\leq \Gamma\leq  T.$
 Let $R$ and $S$ be  continuous  $(q,0)$-forms on $\U.$ Then
 the following  inequalities hold:
 $$
 \big|\int_\U  R\wedge \overline S \wedge T\big |^2\leq \big( i^{q^2}\int_\U  R\wedge \overline R \wedge T\big )  \big( i^{q^2}\int_\U  S\wedge \overline S \wedge T\big ),
 $$
 $$
 \big|\int_\U  R\wedge \overline S \wedge \Gamma\big |^2\leq 9\big( i^{q^2}\int_\U  R\wedge \overline R \wedge T\big )  \big( i^{q^2}\int_\U  S\wedge 
 \overline S \wedge T\big ).
 $$
\end{lemma}
\proof
We  may assume that $T$ is a continuous  positive form.  The  general case  will follow by a  regularization  procedure. Let $\Leb$ be the  canonical
Lebesgue  measure in $\C^k.$  We can write
for $y\in \U,$
\begin{equation*}
  \big(i^{q^2} R\wedge \overline R \wedge T\big)(y)=\phi(y)\Leb(y),\quad  \big( i^{q^2}  S\wedge \overline S \wedge T\big)(y)=\psi(y)\Leb(y),\quad
   \big(i^{q^2} R\wedge \overline S \wedge T\big)(y)= f(y)\Leb(y),
\end{equation*}
where $\phi,$ $\psi$ and  $f$ are continuous  functions.
For every $y\in\U,$ consider  also the  quadratic  form   $g_y:\ \C\to\C $ defined  by 
$$ g_y(t)\Leb(y)=[ (R+tS)\wedge   (\overline R+\bar t\overline S)\wedge T ](y)\quad\text{for}\quad  t\in\C.$$
Since the current $T$ is positive, we  see that $g_y(t)\geq 0. $ Hence,  the discriminant  of $g_y$ is  $\leq 0,$ which implies  that  
  $|f(y)|^2\leq \phi(y)\psi(y)$ for $y\in  \U.$  So  by the Cauchy--Schwarz inequality, we  get
$$
\big|\int_\U  R\wedge \overline S \wedge T\big |=\int_{\U} |f(y)|\Leb(y)\leq  \int_\U\sqrt{\phi(y)\psi(y)}\Leb(y)\leq 
\big(\int_{\U} |\phi(y)|\Leb(y)\big )^{1\over 2}   \big (\int_{\U} |\psi(y)|\Leb(y)\big )^{1\over 2}.
$$
This  proves the  first inequality of the lemma.

Since  $-T\leq \Gamma\leq T,$ we infer that $\Gamma+T$ is a positive  current and $\Gamma+T\leq 2T.$ Consequently,  the  first inequality of the lemma implies that
$$
 \big|\int_\U  R\wedge \overline S \wedge (\Gamma+T)\big |\leq \big( i^{q^2}\int_\U  R\wedge \overline R \wedge (\Gamma+T)\big )^{1\over 2}  \big( i^{q^2}\int_\U  S\wedge \overline S \wedge (\Gamma+T)\big )^{1\over 2}.
 $$
  Since    $0\leq \Gamma+T\leq 2T,$ it follows that
 $$
 \big|\int_\U  R\wedge \overline S \wedge (\Gamma+T)\big |\leq 2\big( i^{q^2}\int_\U  R\wedge \overline R \wedge T\big )^{1\over 2}  \big( i^{q^2}\int_\U  S\wedge \overline S \wedge T\big )^{1\over 2}.
 $$
 Observe that
 $$
 \big|\int_\U  R\wedge \overline S \wedge \Gamma\big | \leq \big|\int_\U  R\wedge \overline S \wedge (\Gamma+T)\big | +\big|\int_\U  R\wedge \overline S \wedge T\big|
 $$
 This, combined  with the first inequality of the lemma and  the last inequality, implies the second  inequality of the lemma.
\endproof

 %
\begin{proposition} \label{P:spec-wedge}
  Let $T$ be a  positive current of bidimension $(q,q)$ on $\U.$
 Let   $S_1,\ldots,S_q$  and  $S'_1,\ldots, S'_q$ be $2q$ positive forms of bidegree $(1,1)$   on $\U.$
  Let $\Gamma_1,\ldots,\Gamma_q$ be  $q$  real $2$-forms on $\U$ and let $H_1,\ldots,H_q$ be  $q$  real $(1,1)$-forms  in the class $\Hc=\Hc(\U)$   such that
  \begin{eqnarray*}
  -S_j\leq \Gamma^\sharp_j+H_j&\leq& S_j\qquad\text{for}\qquad 1\leq j\leq q,\\
  (\Gamma_1,\ldots  ,\Gamma_q)&\trianglelefteq & (S'_1,\ldots,S'_q).
  \end{eqnarray*}
Then there is a constant $c>0$ which depends only on the  dimension $k$ such that 
   \begin{multline*}
  \big|\int_{\U}T\wedge \Gamma_1\wedge \ldots\wedge  \Gamma_q\big|^2  \leq c  \sum_{I,J}\sum_{j=0}^{|I|} \big(\int_{\U}T\wedge \pi^*\omega^j\wedge \hat\beta^{|I|-j} \wedge S_J\wedge S'_{(I\cup J)^\bfc}\big)\\
   \cdot    \big(\int_{\U}T\wedge \pi^*\omega^{|I|-j}\wedge \hat\beta^{j} \wedge S_J\wedge S'_{(I\cup J)^\bfc}\big).
   \end{multline*}
   Here,  the first  sum $\sum_{I,J}$ is  taken over all
   $I,J\subset \{1,\ldots,q\}$  such that  $H_j\not\equiv 0$ for $j\in I,$ and that  $I\cap J=\varnothing,$  and $|(I\cup J)^\bfc|$\footnote{For a  subset $I$ of a given set $K,$  $I^\bfc$ denotes the complement of $I$ in $K,$ that is,  $I^\bfc:=K\setminus I.$} is  even. 
\end{proposition}
 
\proof
Using  the notation  introduced in Definition
 \ref{D:precsim}, we get the  decomposition  $\Gamma_j=\Gamma_j^{1,1}+\Gamma_j^{0,2}+\Gamma_j^{2,0}$ for $1\leq j\leq q.$
A consideration of bidegree gives that
 $$
 T\wedge \Gamma_1\wedge \ldots\wedge  \Gamma_q=
 \sum_{K,K'}  T\wedge \bigwedge_{j\in K} \Gamma^{1,1}_j\wedge \bigwedge_{j\in K'} \Gamma^{2,0}_j\wedge \bigwedge_{j\in (K\cup K')^\bfc} \Gamma^{0,2}_j,
 $$
where the  sum $\sum_{K,K'}$ is taken over  all $K,K'\subset\{1,\ldots,q\}$  such that $K\cap K'=\varnothing$  and  $|K|+2|K'|=q$.
Using  the equality   $\Gamma^{1,1}_j=(\Gamma^{1,1}_j+H_j)-H_j$ for $j\in K,$  the above  expression  is equal to
\begin{eqnarray*}
 &&\sum_{K,K'}  T\wedge \bigwedge_{j\in K} ((\Gamma^{1,1}_j+H_j)-H_j)\wedge \bigwedge_{j\in K'} \Gamma^{2,0}_j\wedge \bigwedge_{j\in (K\cup K')^\bfc} \Gamma^{0,2}_j\\
 &=&\sum_{I,J, K'}(-1)^{|I|}T\wedge \bigwedge_{j\in J} (\Gamma^{1,1}_j+H_j)\wedge \bigwedge_{j\in I} H_j \wedge \bigwedge_{j\in K'} \Gamma^{2,0}_j\wedge \bigwedge_{j\in (I\cup J\cup K')^\bfc} \Gamma^{0,2}_j,
\end{eqnarray*}
where setting $K=I\cup J,$ the last sum is taken over all $I,J,K'\subset\{1,\ldots,q\}$ such that $I,J,K'$ are mutually disjoint and $|I|+|J|+2|K'|=q.$
Since $\pm\big(\Gamma^\sharp_j+H_j\big)\leq S_j$ for $1\leq j\leq q,$
we deduce from
Lemma \ref{L:basic-positive-estimate} that 
there is a constant $c>0$ such that
\begin{equation}\label{e:wedge-prod-H}-\pm  \bigwedge_{j\in J} (\Gamma^{1,1}_j+H_j) \leq  c\cdot \, S_J,\quad\text{where}\quad S_J:=\bigwedge_{j\in J}  S_j.\end{equation}
Moreover, since $(\Gamma_1,\ldots  ,\Gamma_q)\trianglelefteq  (S'_1,\ldots,S'_q),$
we infer from 
Definition 
 \ref{D:precsim-bis} that
 there is a constant $c>0$ such that   the  following   two inequalities  hold for all $1\leq j,j'\leq q$ 
 $$
  \Re\big[\Gamma^{2,0}_j\wedge \overline{\Gamma^{2,0}_{j'}}\big]\leq c S'_j\wedge  S'_{j'}\ \quad\text{and}\quad \Im\big[ \Gamma^{2,0}_j\wedge \overline{\Gamma^{2,0}_{j'}}\big]\leq c (S'_j\wedge  S'_{j'}).
 $$
 Consequently,  applying  
Lemma \ref{L:basic-positive-estimate} yields (see the notation in \eqref{e:wedge-prod-H}) that
\begin{equation*}
 \Re\big[\bigwedge_{j\in K'} \Gamma^{2,0}_j\wedge \bigwedge_{j\in (I\cup J\cup K')^\bfc} \Gamma^{0,2}_j\big]\lesssim   S'_{(I\cup J)^\bfc}\quad\text{and}\quad
\Im\big[\bigwedge_{j\in K'} \Gamma^{2,0}_j\wedge \bigwedge_{j\in (I\cup J\cup K')^\bfc} \Gamma^{0,2}_j\big]\lesssim   S'_{(I\cup J)^\bfc}
 \end{equation*}
 If $H_j\equiv 0$ for some $j\in I,$ then  clearly $\bigwedge_{j\in I} H_j=0.$
 So we only   consider  $I\subset \{1,\ldots, q\}$ such that $H_j\not\equiv 0$ for $j\in I.$
 By Lemma \ref{L:Hc} there is a constant $c>0$  such that
\begin{equation*}
\pm  \bigwedge_{j\in I} H_j \leq c \big( \sum_{j=0}^{|I|}  \pi^*\omega^j \wedge \hat\beta^{|I|-j} \big).  
\end{equation*}
 Combining  the last two inequalities and \eqref{e:wedge-prod-H}, we may apply the second inequality of Lemma \ref{L:Cauchy-Schwarz}.
 Consequently, there is a constant $c>0$ such that
\begin{multline*}
\big |\int_\U T\wedge \bigwedge_{j\in J} (\Gamma^{1,1}_j+H_j)\wedge \bigwedge_{j\in I} H_j \wedge \bigwedge_{j\in K'} \Gamma^{2,0}_j\wedge \bigwedge_{j\in (I\cup J\cup K')^\bfc} \Gamma^{0,2}_j\big|^2\\
\leq c  \sum_{j=0}^{|I|} \big(\int_{\U}T\wedge \pi^*\omega^j\wedge \hat\beta^{|I|-j} \wedge S_J\wedge S'_{(I\cup J)^\bfc}\big)
   \cdot    \big(\int_{\U}T\wedge \pi^*\omega^{|I|-j}\wedge \hat\beta^{j} \wedge S_J\wedge S'_{(I\cup J)^\bfc}\big).
\end{multline*}
 This implies
the result.
\endproof

\begin{lemma}\label{L:basic-positive-difference}
 Let $T$ be a  positive current of bidgree $(p,p)$ on $\bfU.$
 Let  $R_1,\ldots,R_{k-p}$  and   
 $S_1,\ldots,S_{k-p}$ and $S'_1,\ldots,S'_{k-p}$  be  positive $(1,1)$-currents on $\Tube(B,\bfr)\subset \E,$ and  for each $1\leq\ell\leq \ell_0$ 
 let $H_{\ell,1},\ldots,H_{\ell,k-p}$ be    real $(1,1)$-forms  in the class $\Hc$ on $\U_\ell$    such that we have
 $$
 \pm\big\lbrace  (\tilde\tau_\ell)^*[(\pi^*\theta_\ell)^{1\over k-p} R_j] -[(\pi^*\theta_\ell)^{1\over k-p}R_j]   -H_{\ell,j}\big\rbrace^\sharp\lesssim  S_j \quad\mbox{on $\U_\ell$  for  $1\leq \ell\leq \ell_0$ and $ 1\leq j\leq k-p;$}$$ and that on $\U_\ell$  for  $1\leq \ell\leq \ell_0,$ we have
 $$
 \Big((\tilde\tau_\ell)^*[(\pi^*\theta_\ell)^{1\over k-p} R_1] -[(\pi^*\theta_\ell)^{1\over k-p}R_1]  ,\ldots,
 (\tilde\tau_\ell)^*[(\pi^*\theta_\ell)^{1\over k-p} R_{k-p}] -[(\pi^*\theta_\ell)^{1\over k-p}R_{k-p}]\Big) \trianglelefteq \big(S'_1,\ldots, S'_{k-p}\big).
 $$
Then for $R:=R_1\wedge \ldots \wedge R_{k-p},$ we have 
   \begin{multline*}
\left|  \langle \tau_*T ,R\rangle  -\langle T^\hash ,R\rangle \right |^2 \leq  c\cdot\, \sum_{\ell=1}^{\ell_0}  \sum\limits_{I,J,K}\sum_{j=0}^{|I|} \big(\int_{\U_\ell}(\pi^*\theta_\ell)^{{|K|\over k-p}}(\tau_\ell)_* T \wedge R_K\wedge  \pi^*\omega^j\wedge \hat\beta^{|I|-j} \wedge S_J\wedge S'_{(I\cup J\cup K)^\bfc}\big)\\
   \cdot    \big(\int_{\U_\ell}(\pi^*\theta_\ell)^{{|K|\over k-p}}(\tau_\ell)_* T \wedge R_K \wedge \pi^*\omega^{|I|-j}\wedge \hat\beta^{j} \wedge S_J\wedge S'_{(I\cup J\cup K)^\bfc}\big).
   \end{multline*}
Here,  $T^\hash$ is  defined in \eqref{e:T-hash} and
\begin{itemize} 
\item[$\bullet$] $c$ is a constant that depends only on the dimension $k$ and  $\ell_0;$
\item[$\bullet$]  the second sum $\sum_{I,J,K}$ is  taken over all $I,J,K\subset \{1,\ldots,k-p\}$  such that    $H_j\not\equiv 0$ for $j\in I,$ and that $I,J,K$ are mutually disjoint,
   and $|(I\cup J\cup K)^\bfc|$ is  even, and  $K\not=\{1,\ldots,k-p\}.$
\end{itemize}
\end{lemma}
\proof Fix $\ell$ with $1\leq \ell\leq \ell_0.$
For  $1\leq j\leq k-p,$ consider 
\begin{equation}\label{e:Gamma_j} \Gamma_j:= (\tilde\tau_\ell)^*[(\pi^*\theta_\ell)^{1\over k-p} R_j]-[(\pi^*\theta_\ell)^{1\over k-p}R_j].
\end{equation}
By  hypothesis, we  get  that
\begin{equation}\label{e:basic-positive-difference(1)}
\begin{split}
 \pm\big(  \Gamma^\sharp_j+H_{\ell,j}\big)&\lesssim S_j    \quad\mbox{on $\U_\ell$  for   $ 1\leq j\leq k-p;$}\\
 \big( \Gamma_1,\ldots,\Gamma_{k-p})&\trianglelefteq (S'_1,\ldots, S'_{k-p})  \quad\mbox{on $\U_\ell.$ }
 \end{split}
 \end{equation}
  Therefore,  applying Lemma \ref{L:basic-T-hash} yields that
 \begin{eqnarray*}\langle \tau_*T,R\rangle -\langle T^\hash, R\rangle &=&\sum_{\ell=1}^{\ell_0}\langle(\tau_\ell)_* T, 
 (\tilde\tau_\ell)^*[(\pi^*\theta_\ell )R]-[(\pi^*\theta_\ell )R]\rangle \\
   &=&\sum_{\ell=1}^{\ell_0}\langle(\tau_\ell)_* T, \bigwedge_{j=1}^{k-p}  (\Gamma_j+(\pi^*\theta_\ell)^{1\over k-p}\cdot R_j)-  
   \bigwedge_{j=1}^{k-p} ( (\pi^*\theta_\ell)^{1\over k-p}\cdot R_j)  \rangle\\
   &=& \sum_{\ell=1}^{\ell_0}\sum\limits_K \langle(\tau_\ell)_* T,   (\pi^*\theta_\ell)^{{|K|\over k-p}}\cdot \big(    R_{K}\wedge \Gamma_{K^\bfc}\big)\rangle,
 \end{eqnarray*}
 where the  inner sum $\sum_K$ in the last line is taken over  $K\subsetneq\{1,\ldots,k-p\}.$ 
 So we have 
 $$
 \langle \tau_*T,R\rangle -\langle T^\hash, R\rangle=\sum_{\ell=1}^{\ell_0}\sum\limits_K \langle  (\pi^*\theta_\ell)^{{|K|\over k-p}}(\tau_\ell)_* T
 \wedge    R_{K}, \Gamma_{K^\bfc}\rangle.
 $$
Using \eqref{e:basic-positive-difference(1)} and applying Proposition \ref{P:spec-wedge} to the last line,
the desired  inequality follows.
\endproof

Let  $ T$ be a current  defined on $\bfU$  and $0\leq s<r\leq \bfr.$  Consider the currents  $T^{\#}_r$  and $T^{\#}_{s,r}$  defined  on $\U$ as  follows:
\begin{equation}\label{e:T-hash_r}
T^{\#}_r:=\sum_{\ell=1}^{\ell_0}  (\pi^*\theta_\ell)\cdot (\ind_{\Tube(B,r)}\circ \tilde \tau_\ell )\cdot(\tau_\ell)_*( T|_{\bfU_\ell})\quad\text{and}\quad
T^{\#}_{s,r}:=\sum_{\ell=1}^{\ell_0}  (\pi^*\theta_\ell)\cdot (\ind_{\Tube(B,s,r)}\circ \tilde \tau_\ell )\cdot(\tau_\ell)_*( T|_{\bfU_\ell}).
\end{equation}
The  following lemma permits us to  replace the  integral  $ \langle \tau_*T,\ind_{\Tube(B,r)}R\rangle$
(resp.    $\langle \tau_*T,\ind_{\Tube(B,s,r)}R\rangle$) by  a  simpler  one $\langle T^\hash_r, R\rangle$ (resp. $\langle T^\hash_{s,r}, R\rangle$).
\begin{lemma}\label{L:basic-T-hash-bis}
The following identities holds
 \begin{eqnarray*}\langle \tau_*T,\ind_{\Tube(B,r)}R\rangle -\langle T^\hash_r, R\rangle &=&\sum_{\ell=1}^{\ell_0}\langle(\tau_\ell)_* T,  (\ind_{\Tube(B,r)}\circ \tilde\tau_\ell)\cdot \big((\tilde\tau_\ell)^*((\pi^*\theta_\ell) R)-((\pi^*\theta_\ell) R)\big)\rangle ,\\
  \langle \tau_*T,\ind_{\Tube(B,s,r)}R\rangle -\langle T^\hash_{s,r}, R\rangle &=&\sum_{\ell=1}^{\ell_0}\langle(\tau_\ell)_* T,  (\ind_{\Tube(B,s,r)}\circ \tilde\tau_\ell)\cdot \big((\tilde\tau_\ell)^*((\pi^*\theta_\ell) R)-((\pi^*\theta_\ell) R)\big)\rangle .
 \end{eqnarray*}
\end{lemma}
\proof
We only give the proof of the  first identity since the proof of the  second one  is  similar.
Since  $\sum_{\ell=1}^{\ell_0}  \pi^*\theta_\ell  =1$  on an open  neighborhood  of $\pi^{-1}(\overline{\bfU\cap V})\subset \pi^{-1}(V),$  we have
$$\tau_*T=\sum_{\ell=1}^{\ell_0}  \pi^*\theta_\ell\cdot\tau_* T.$$
So using \eqref{e:T-hash_r} we get that
$$
\langle \tau_*T,\ind_{\Tube(B,r)}R\rangle -\langle T^\hash_r, R\rangle =\sum_{\ell=1}^{\ell_0}  \big( \langle \tau_*T, (\pi^*\theta_\ell )\ind_{\Tube(B,r)}R\rangle
- \langle (\tau_\ell)_*( T|_{\bfU_\ell}),(\ind_{\Tube(B,r)}\circ \tilde\tau_\ell) (\pi^*\theta_\ell)\cdot R\rangle \big ).
$$
Writing   $\tau_*T=(\tau\circ \tau_\ell^{-1})_*(\tau_\ell)_* T=(\tilde\tau_\ell)_*(\tau_\ell)_* T$ on $\U_\ell,$
we get that
\begin{eqnarray*}
&&\langle \tau_*T,\ind_{\Tube(B,r)}R\rangle -\langle T^\hash_r, R\rangle =\sum_{\ell=1}^{\ell_0}  \big( \langle (\tilde\tau_\ell)_*(\tau_\ell)_* T, (\pi^*\theta_\ell)\ind_{\Tube(B,r)}R\rangle \\
&&\qquad -  \langle (\tau_\ell)_*( T|_{\bfU_\ell}),(\ind_{\Tube(B,r)}\circ \tilde\tau_\ell) (\pi^*\theta_\ell)\cdot R\rangle \big ) \\
&=&\sum_{\ell=1}^{\ell_0} \big( \langle(\tau_\ell)_* T, (\ind_{\Tube(B,r)}\circ \tilde\tau_\ell) (\tilde\tau_\ell)^* [(\pi^*\theta_\ell)R]\rangle - \langle 
(\tau_\ell)_*( T),(\ind_{\Tube(B,r)}\circ \tilde\tau_\ell) [(\pi^*\theta_\ell) R]\rangle \big)\\
&=& \sum_{\ell=1}^{\ell_0}  \langle(\tau_\ell)_* T, (\ind_{\Tube(B,r)}\circ \tilde\tau_\ell)((\tilde\tau_\ell)^*[(\pi^*\theta_\ell) R]-[(\pi^*\theta_\ell) R])\rangle ,
\end{eqnarray*}
which implies the desired identity.
\endproof

\begin{lemma}\label{L:basic-positive-difference-bis}
 Let $T$ be a  positive current of bidgree $(p,p)$ on $\bfU.$
 Let  $R_1,\ldots,R_{k-p}$  and   
 $S_1,\ldots,S_{k-p}$ and $S'_1,\ldots,S'_{k-p}$  be  positive $(1,1)$-currents on $\Tube(B,\bfr)\subset \E,$ and  for each $1\leq\ell\leq \ell_0$ 
 let $H_{\ell,1},\ldots,H_{\ell,k-p}$ be    real $(1,1)$-forms  in the class $\Hc$ on $\U_\ell$    such that
 \begin{eqnarray*}
 \varphi^{1\over 2}R_j&\lesssim& S_j\quad  \text{and}\quad \varphi^{1\over 2}R_j\lesssim S'_j  \quad\mbox{on $\Tube(B,\bfr)$  for  $ 1\leq j\leq k-p;$}\\
 \pm[(\tilde \tau_\ell)^* R_j-R_j-H_{\ell,j}]^\sharp&\lesssim & S_j \quad\mbox{on $\U_\ell$  for  $1\leq \ell\leq \ell_0$ and $ 1\leq j\leq k-p;$}\\
 \big( (\tilde \tau_\ell)^* R_1-R_1,\ldots,    (\tilde \tau_\ell)^* R_{k-p}-R_{k-p}) &\trianglelefteq &\big(S'_1,\ldots, S'_{k-p}\big)\quad \mbox{on $\U_\ell$  for  $1\leq \ell\leq \ell_0.$}
 \end{eqnarray*}
Let $0<s<r\leq \bfr$ and set  $R:=R_1\wedge \ldots \wedge R_{k-p}.$ 
Suppose in addition that there are constants $0<c_5<1$ and $c_6>1$ and    positive $(1,1)$-forms $R'_1,\ldots,R'_{k-p}$  such that 
\begin{itemize}
 \item $R'_j\geq R_j$  for $1\leq j\leq k-p;$ 
\item if  $y\in  \U_\ell$  with  $0<\theta_\ell(y)< c_5,$  then  we may find  $1\leq \ell'\leq \ell_0$ 
and  an open neighborhood $\U_y$  of $y$ in  $\U$  such that for $x\in\U_y,$ we have that $\theta_{\ell'}(x)>  c_5$  and that $R_j(x)\leq  c_6(\tilde\tau_{\ell'}\circ \tilde\tau^{-1}_\ell)^* R'_j(x)$ and  that $S_j(x)\leq  c_6(\tilde\tau_{\ell'}\circ \tilde\tau^{-1}_\ell)^* S_j(x)$  and that
 $S'_j(x)\leq  c_6(\tilde\tau_{\ell'}\circ \tilde\tau^{-1}_\ell)^* S'_j(x).$   
\end{itemize}
Then there is a constant 
$c$ that depends on $c_5,c_6$  and $\ell_0$ such that
 \begin{multline*}
\left   |                \langle \tau_*T ,\ind_{\Tube(B,r)}R\rangle  -\langle T^\hash_r ,R\rangle \right |^2\\
\leq  c\cdot\, \sum_{\ell=1}^{\ell_0}  \sum\limits_{I,J,K}\sum_{j=0}^{|I|} \big( \int
 (\ind_{\Tube(B,r)}\circ \tilde\tau_\ell) (\pi^*\theta_\ell)(\tau_\ell)_* T 
\wedge R'_K\wedge  \pi^*\omega^j\wedge \hat\beta^{|I|-j} \wedge S_J\wedge S'_{(I\cup J\cup K)^\bfc}\big)\\
   \cdot    \big(  \int  (\ind_{\Tube(B,r)}\circ \tilde\tau_\ell) (\pi^*\theta_\ell)(\tau_\ell)_* T        \wedge R'_K \wedge \pi^*\omega^{|I|-j}\wedge \hat\beta^{j} \wedge S_J\wedge S'_{(I\cup J\cup K)^\bfc}\big).
   \end{multline*} 
\begin{multline*}
\left   |                \langle \tau_*T ,\ind_{\Tube(B,s,r)}R\rangle  -\langle T^\hash_{s,r} ,R\rangle \right |^2\\
\leq 
  c\cdot\, \sum_{\ell=1}^{\ell_0}  \sum\limits_{I,J,K}\sum_{j=0}^{|I|} \big( \int
 (\ind_{\Tube(B,s,r)}\circ \tilde\tau_\ell) (\pi^*\theta_\ell)(\tau_\ell)_* T 
\wedge R'_K\wedge  \pi^*\omega^j\wedge \hat\beta^{|I|-j} \wedge S_J\wedge S'_{(I\cup J\cup K)^\bfc}\big)\\
   \cdot    \big( \int   (\ind_{\Tube(B,s,r)}\circ \tilde\tau_\ell) (\pi^*\theta_\ell)(\tau_\ell)_* T        \wedge R'_K \wedge \pi^*\omega^{|I|-j}\wedge \hat\beta^{j} \wedge S_J\wedge S'_{(I\cup J\cup K)^\bfc}\big).
\end{multline*}
Here,   the  sum $\sum_{I,J,K}$ is  taken over all $I,J,K\subset \{1,\ldots,k-p\}$  such that    $H_j\not\equiv 0$ for $j\in I,$ and that $I,J,K$ are mutually disjoint,
   and $|(I\cup J\cup K)^\bfc|$ is  even, and  $K\not=\{1,\ldots,k-p\}.$
\end{lemma}
\begin{remark}\rm
 Lemma \ref{L:basic-positive-difference-bis} has the advantage over Lemma \ref{L:basic-positive-difference}
 in that the former gets rid of  the somehow cumbersome exponent  of the functions $(\pi^*\theta_\ell)^{|K|\over k-p}$ of the latter.
\end{remark}

\proof
We only give the proof of the first inequality, since the proof of the second one is similar.
For  $1\leq j\leq k-p,$ write $\Gamma_j:= (\tilde\tau_\ell)^*[(\pi^*\theta_\ell)^{1\over k-p} R_j]-[(\pi^*\theta_\ell)^{1\over k-p}R_j]\big).$   
Pick an arbitrary $\ell$ with $1\leq \ell\leq \ell_0$ and an arbitrary point $y\in\Tube(B,r)\cap\U.$  We argue using the three cases   as at the  end of the proof. Consequently, we may assume without loss of generality that 
$(\pi^*\theta_\ell)(y)\geq  c_5,$ where $0<c_5<1$ is  a constant.  Write $y=(z,w).$
Using that  $(\pi^*\theta_\ell)(y) -(\tilde \tau_\ell)^*(\pi^*\theta_\ell)(y)=O(z),$ we  deduce that
$$(\pi^*\theta_\ell)(y) -(\tilde \tau_\ell)^*(\pi^*\theta_\ell)(y)=O(\varphi^{1\over 2}).$$
Consequently, we   infer that
\begin{equation*}
\Gamma^\sharp_j -[(\tilde \tau_\ell)^* R_j-R_j]^\sharp=O(\varphi^{1\over 2}R_j)\lesssim S_j.
\end{equation*}
Moreover, using  the assumption   $R_j\leq  \varphi^{1\over 2} S'_j,$ we also  get that
\begin{equation*}
 \pm\Re[ \Gamma^{2,0}_p\wedge\Gamma^{0,2}_q]\leq \pm  \big[(\tilde \tau_\ell)^* R_p-R_p\big]^{2,0}\wedge \big[(\tilde \tau_\ell)^* R_q-R_q)\big]^{0,2}\lesssim (1+\varphi^{1\over 2})  S'_p\wedge  S'_q\lesssim   S'_p\wedge  S'_q.
\end{equation*}
Therefore, we have shown that   
\begin{equation}\label{e:basic-positive-difference-bis(1)}
\begin{split}
 \pm\big(  \Gamma^\sharp_j+H_{\ell,j}\big)&\leq S_j    \quad\mbox{on $\U_\ell$  for  $1\leq \ell\leq \ell_0$ and $ 1\leq j\leq k-p;$}\\
 \big( \Gamma_1,\ldots,\Gamma_{k-p})&\trianglelefteq (S'_1,\ldots, S'_{k-p}).
 \end{split}
 \end{equation}
  Therefore,  applying Lemma \ref{L:basic-T-hash-bis} and arguing as in the  proof of Lemma \ref{L:basic-positive-difference} yields that
 \begin{eqnarray*}&&\langle \tau_*T,\ind_{\Tube(B,r)}R\rangle -\langle T^\hash_r, R\rangle =\sum_{\ell=1}^{\ell_0}\langle(\tau_\ell)_* T, 
 (\ind_{\Tube(B,r)}\circ \tilde\tau_\ell) \big((\tilde\tau_\ell)^*[(\pi^*\theta_\ell )R]-[(\pi^*\theta_\ell )R]\big)\rangle \\
   &=&\sum_{\ell=1}^{\ell_0}\langle(\tau_\ell)_* T,(\ind_{\Tube(B,r)}\circ \tilde\tau_\ell) \big( \bigwedge_{j=1}^{k-p}  (\Gamma_j+(\pi^*\theta_\ell)^{1\over k-p}\cdot R_j)-  
   \bigwedge_{j=1}^{k-p} ( (\pi^*\theta_\ell)^{1\over k-p}\cdot R_j) \big) \rangle\\
   &=& \sum_{\ell=1}^{\ell_0}\sum\limits_K \langle(\tau_\ell)_* T,   (\ind_{\Tube(B,r)}\circ \tilde\tau_\ell) 
   \cdot(\pi^*\theta_\ell)^{{|K|\over k-p}}\cdot \big(    R_{K}\wedge \Gamma_{K^\bfc}\big)\rangle.
 \end{eqnarray*}
 where the  inner sum $\sum_K$ in the last line is taken over  $K\subsetneq\{1,\ldots,k-p\}.$ 
 So we have 
 $$
 \langle \tau_*T,\ind_{\Tube(B,r)}R\rangle -\langle T^\hash_r, R\rangle=\sum_{\ell=1}^{\ell_0}\sum\limits_K \langle  (\ind_{\Tube(B,r)}\circ \tilde\tau_\ell)  (\pi^*\theta_\ell)^{{|K|\over k-p}}(\tau_\ell)_* T
 \wedge    R_{K}, \Gamma_{K^\bfc}\rangle.
 $$
Using \eqref{e:basic-positive-difference-bis(1)} and applying Proposition \ref{P:spec-wedge} to the last line,
we  infer that
\begin{multline*}
\left   |                \langle \tau_*T ,\ind_{\Tube(B,r)}R\rangle  -\langle T^\hash_r ,R\rangle \right |^2\\
\leq  c\cdot\, \sum_{\ell=1}^{\ell_0}  \sum\limits_{I,J,K}\sum_{j=0}^{|I|} \big( \int
 (\ind_{\Tube(B,r)}\circ \tilde\tau_\ell) (\pi^*\theta_\ell)^{|K|\over k-p}(\tau_\ell)_* T 
\wedge R'_K\wedge  \pi^*\omega^j\wedge \hat\beta^{|I|-j} \wedge S_J\wedge S'_{(I\cup J\cup K)^\bfc}\big)\\
   \cdot    \big(  \int  (\ind_{\Tube(B,r)}\circ \tilde\tau_\ell) (\pi^*\theta_\ell)^{|K|\over k-p}(\tau_\ell)_* T        \wedge R'_K \wedge \pi^*\omega^{|I|-j}\wedge \hat\beta^{j} \wedge S_J\wedge S'_{(I\cup J\cup K)^\bfc}\big).
   \end{multline*} 
It remains  to  get rid of  the  undesired  exponent of $(\pi^*\theta_\ell)^{|K|\over k-p}.$  Pick an arbitrary $\ell$ with $1\leq \ell\leq \ell_0$ and an arbitrary
 point $y\in  \Tube(B,r)\cap \U.$ There are 3 cases to consider.
 
 \noindent {\bf Case 1:  $(\pi^*\theta_\ell)(y)=0.$} In this case $(\pi^*\theta_\ell)^{{|J^\bfc|\over k-p}}(y)=(\pi^*\theta_\ell)(y)=0.$

 \noindent {\bf Case 2:  
 $(\pi^*\theta_\ell)(y)\geq c_5.$} In this case  $(\pi^*\theta_\ell)^{{|J^\bfc|\over k-p}}(y)\approx (\pi^*\theta_\ell)(y)\approx 1.$

 \noindent {\bf Case 3:  
 $0<(\pi^*\theta_\ell)(y)< c_5.$} In this case by  the  assumption of  item (2),  
 we may find  $1\leq \ell'\leq \ell_0$ 
and  an open neighborhood $\U_y$  of $y$ in  $\U$  such that for $x\in\U_y,$ we have that $\theta_{\ell'}(x)>  c_5$  and that $R_j(x)\leq  c_6(\tilde\tau_{\ell'}\circ \tilde\tau^{-1}_\ell)^* R'_j(x)$ and  that $S_j(x)\leq  c_6(\tilde\tau_{\ell'}\circ \tilde\tau^{-1}_\ell)^* S_j(x)$  and that
 $S'_j(x)\leq  c_6(\tilde\tau_{\ell'}\circ \tilde\tau^{-1}_\ell)^* S'_j(x).$
  So   $(\pi^*\theta_{\ell'})^{{|J^\bfc|\over k-p}}(x)\approx (\pi^*\theta_{\ell'})(x).$
Let  $I,J,K\subset \{1,\ldots,k-p\}$  such that  $I,J,K$ are mutually disjoint,
   and $|(I\cup J\cup K)^\bfc|$ is  even and  $K\not=\{1,\ldots,k-p\}.$ 
   Then  we have  
\begin{eqnarray*} &&\int_{(\tau_{\ell'}\circ \tau_\ell^{-1})(\U_y)} (\tau_{\ell'})_* T \wedge R'_K\wedge  \pi^*\omega^j\wedge \hat\beta^{|I|-j} \wedge S_J\wedge S'_{(I\cup J\cup K)^\bfc}\\
&=& \int_{\U_y} (\tau_{\ell'}\circ \tau_\ell^{-1})^*\big[ (\tau_{\ell'})_* T \wedge R'_K\wedge  \pi^*\omega^j\wedge \hat\beta^{|I|-j} \wedge S_J\wedge S'_{(I\cup J\cup K)^\bfc}\big]\\
&=&\int_{\U_y}  (\tau_\ell)_* T \wedge (\tau_{\ell'}\circ \tau_\ell^{-1})^*\big[R'_K\wedge  \pi^*\omega^j\wedge \hat\beta^{|I|-j} \wedge S_J\wedge S'_{(I\cup J\cup K)^\bfc}\big]. \end{eqnarray*}
Since  $\tau_{\ell'}\circ \tau_\ell^{-1}$ is a holomorphic admissible map, we  see that
$$(\tau_{\ell'}\circ \tau_\ell^{-1})^*(\pi^*\omega+\hat\beta)\approx (\pi^*\omega+\hat\beta)\qquad\text{on}\qquad \U_y.$$ 
Using this and the  above   inequalities on $\U_y,$  
we may apply  Lemma \ref{L:basic-positive-estimate}. Consequently, there is a constant $c>0$  which depends only on $c_5,$ $c_6$ such that 
$$
c \int_{(\tau_{\ell'}\circ \tau_\ell^{-1})(\U_y)} (\tau_{\ell'})_* T \wedge R'_K\wedge  (\pi^*\omega+ \hat\beta)^{|I|} \wedge S_J\wedge S'_{(I\cup J\cup K)^\bfc}\geq \int_{\U_y}  (\tau_\ell)_* T \wedge R_K\wedge  (\pi^*\omega^j+ \hat\beta)^{|I|} \wedge S_J\wedge S'_{(I\cup J\cup K)^\bfc}.
$$
This completes the proof in  Case 3.

\endproof

\section{Positive currents and positive closed currents}\label{S:Positive-closed-currents} 
 
We keep the  Extended Standing Hypothesis introduced  in Subsection \ref{SS:Ex-Stand-Hyp}.

\subsection{Local and global mass indicators for positive currents}

We use the notation  introduced  at the beginning  of  Section \ref{S:Regularization}. 
 Following the  model  formulas \eqref{e:mass-indicators} and \eqref{e:T-hash},
 we introduce  the following  mass indicators   for a  positive  current $T$  of bidegree $(p,p)$ defined on $X.$
 For $0\leq j\leq k$ and $0\leq q\leq k-l$ and  $1\leq\ell\leq \ell_0,$ and  for $0<s<r\leq \bfr,$ 
 \begin{equation}\label{e:local-mass-indicators-bis}
 \begin{split}
 \Mc_j(T,r,\tau_\ell)&:={1\over  r^{2(k-p-j)}}\int (\ind_{\Tube(B,r)}\circ \tilde\tau_\ell)  (\pi^*\theta_\ell)\cdot (\tau_\ell)_*( T|_{\bfU_\ell})
 \wedge\pi^*\omega^j\wedge  \hat\beta^{k-p-j},\\
 \Kc_{j,q}(T,r,\tau_\ell)&:=\int(\ind_{\Tube(B,r)}\circ \tilde\tau_\ell)  (\pi^*\theta_\ell)\cdot (\tau_\ell)_*( T|_{\bfU_\ell})
 \wedge\pi^*\omega^j\wedge  \hat\beta^{k-p-q-j}\wedge \hat \alpha^q,\\
 \Kc_{j,q}(T,s,r,\tau_\ell)&:=\int (\ind_{\Tube(B,s,r)}\circ \tilde\tau_\ell) (\pi^*\theta_\ell)\cdot (\tau_\ell)_*( T|_{\bfU_\ell})
 \wedge\pi^*\omega^j\wedge  \hat\beta^{k-p-q-j}\wedge \hat \alpha^q.
 \end{split}
 \end{equation}
\begin{remark}\rm Recall  from Corollary \ref{C:Lelong-Jensen}
that $(\tau_\ell)_*( T|_{\bfU_\ell})
 \wedge\pi^*\omega^\upm$ is of full bidegree $(l,l)$ in $\{dw,d\bar w\}.$
Consequently,
by the bidegree reason, we  deduce  that  $\Mc_j(T,r,\tau_\ell),$ $\Kc_{j,q}(T,r,\tau)$ and $\Kc_{j,q}(T,s,r,\tau)$
are equal to $0$ provided that $j>\upm.$ 
\end{remark}
 
  We define   the following  global mass indicators.
 \begin{equation}\label{e:global-mass-indicators}
 \begin{split}
 \Mc_j(T,r)=\Mc_j(T,r,\Uc)&:=\sum_{\ell=1}^{\ell_0}\Mc_j(T,r,\tau_\ell) ,\\
 \Mc^\tot(T,r)=\Mc^\tot(T,r,\Uc)&:=\sum_{j=0}^\upm \Mc_j(T,r),\\
 \Kc_{j,q}(T,r)=    \Kc_{j,q}(T,r,\Uc) &:=\sum_{\ell=1}^{\ell_0}\Kc_{j,q}(T,r,\tau_\ell),\\
 \Kc_{j,q}(T,s,r)=    \Kc_{j,q}(T,s,r,\Uc) &:=\sum_{\ell=1}^{\ell_0}\Kc_{j,q}(T,s,r,\tau_\ell).
 \end{split}
 \end{equation}

 \begin{lemma}
 \label{L:Mc_j}
  \begin{equation*}
 \begin{split}
 \Mc_j(T,r)&= {1\over r^{2(k-p-j)}}\int  T^\hash_r\wedge\pi^*\omega^j\wedge  \hat\beta^{k-p-j} ,\\
 \Kc_{j,q}(T,r) &=     \int  T^\hash_r\wedge
 \pi^*\omega^j\wedge  \hat\beta^{k-p-q-j}\wedge \hat\alpha^q ,\\
 \Kc_{j,q}(T,s,r) &=\int  T^\hash_{s,r}
 \wedge\pi^*\omega^j\wedge  \hat\beta^{k-p-q-j}\wedge \hat\alpha^q.
 \end{split}
 \end{equation*}
  \end{lemma}
\proof
It follows  from \eqref{e:local-mass-indicators-bis}, \eqref{e:global-mass-indicators}  and \eqref{e:T-hash_r}.
\endproof

\begin{lemma}\label{L:comparison-Mc} For every constant $\rho >1$ there is a constant $c>0$ such that
for $0\leq j\leq k$ and for $0<r<s<\rho r\leq \bfr$ and for every  positive current $T$ of bidegree $(p,p)$ on $\bfU,$    we have
$\Mc_j(T,r)<c\Mc_j(T,s).$ In particular, it holds that $\Mc^\tot(T,r)<c\Mc^\tot(T,s).$
\end{lemma}
\proof We only prove  the  first  inequality  since by the definition of  $\Mc^\tot(T,r)$
the second  inequality  is a consequence of   the  first  one. 

Since $0<r<s\leq \bfr$ and $T$  is a positive currents and $\omega,$ $\hat\beta$ are positive forms, we have
$$
 \int (\ind_{\Tube(B,r)}\circ \tilde\tau_\ell)  (\pi^*\theta_\ell)\cdot (\tau_\ell)_*( T|_{\bfU_\ell})
 \wedge\pi^*\omega^j\wedge  \hat\beta^{k-p-j}\leq \int (\ind_{\Tube(B,s)}\circ \tilde\tau_\ell)  (\pi^*\theta_\ell)\cdot (\tau_\ell)_*( T|_{\bfU_\ell})
 \wedge\pi^*\omega^j\wedge  \hat\beta^{k-p-j} .
$$
Hence, we infer   from the definition   of  $\Mc_j(T,r,\tau_\ell)$ in \eqref{e:local-mass-indicators-bis}
$$
  r^{2(k-p-j)}\Mc_j(T,r,\tau_\ell)\leq   s^{2(k-p-j)}\Mc_j(T,s,\tau_\ell).
$$
This, combined with the definition of   $\Mc_j(T,r)$ in  \eqref{e:global-mass-indicators}, implies the first  inequality of the lemma for $c:= \rho^{2(k-p-j)}.$
\endproof
\subsection{Finiteness of the mass indicator $\Kc_{j,q}$ for K\"ahler  metrics}

Let $\omega$ be  a K\"ahler metric  on $V.$

\begin{lemma}\label{L:Lelong-smooth-forms}
 Let $T$ be a positive closed  $\Cc^1$-smooth form   on $\bfU.$
 Then  for every $\lowm \leq j\leq \upm,$  we have  $\nu_j(T,B,\tau)=0$ if $j\not=l-p$ and $\nu_j(T,B,\tau)\geq 0$ if $j=l-p.$ 
\end{lemma}
\proof
First consider the  case   $j\not=l-p.$  As  $\lowm \leq j\leq \upm,$ we have   $j>l-p ,$  and hence   $k-p-j<k-l.$
Then by Theorem  \ref{T:Lelong-Jensen-smooth-closed} (1), $\nu_j(T,B,\tau)=0.$

Now consider the case $j=l-p.$ So $j=\lowm.$  Since $\tau$ is  strongly  admissible  $d\tau|_{\overline B}$ is  $\C$-linear,
it follows from  the positivity of $T$ on $\bfU$ that  $(\tau_*T)|_{\overline B}$ is also  a positive form.
Hence, by Theorem  \ref{T:Lelong-Jensen-smooth-closed} (1), $\nu_j(T,B,\tau)\geq 0.$
\endproof

Let $\bfj=(j_1,j_2,j_3,j_4)$ with $j_1,j_3,j_4\in\N $ and $j_2\in{1\over 4}\N,$ and 
$k-p-j_1-j_3\geq 0.$ For $0<s<r\leq\bfr,$ and  for a real current $T$  on $\bfU,$  consider  
\begin{equation}\label{e:I_bfj}\begin{split}  I_{\bfj}(s,r)&:=\int_{\Tube(B,s,r)}\tau_*T\wedge \varphi^{j_2}(c_1-c_2\varphi)^{j_4}\hat\beta^{k-p-j_1-j_3}\wedge
(\pi^*\omega)^{j_3}\wedge \hat\alpha^{j_1},\\ 
I^\hash_{\bfj}(s,r)&:=\int_{\Tube(B,s,r)}T^\hash_{s,r}\wedge \varphi^{j_2}(c_1-c_2\varphi)^{j_4}\hat\beta^{k-p-j_1-j_3}\wedge
(\pi^*\omega)^{j_3}\wedge \hat\alpha^{j_1}.
\end{split}
\end{equation}
We define  $I_\bfj(r)$ and $I^\hash_\bfj(r)$ similarly  replacing the current $T^\hash_{s,r}$ (resp.  the domain of integration $\Tube(B,s,r)$) by
$T^\hash_r$  (resp. $\Tube(B,r)$).
\begin{remark}\label{R:I_bfj}
 \rm  Observe that $\Kc_{j,q}(T,r)=I^\hash_{(q,0,j,0)}(r)$ and $\Kc_{j,q}(T,s,r)=I^\hash_{(q,0,j,0)}(s,r).$
\end{remark}

\begin{lemma}\label{L:spec-wedge}  There is a constant $c$ independent of $T$ and $s,r$   such that the   following inequality holds
 \begin{equation*}
|I_\bfj(s,r)- I^\hash_\bfj(s,r)|^2 \leq c\big(\sum_{\bfj'} I^\hash_{\bfj'}(s,r)\big)\big ( \sum_{\bfj''} I^\hash_{\bfj''}(s,r) \big).  
\end{equation*}
Here, on the RHS:
\begin{itemize} \item[$\bullet$] the first sum  is taken over a finite number of multi-indices    $\bfj'=(j'_1,j'_2,j'_3,j'_4)$ as above  such that  $j'_1\leq  j_1$  and $j'_2\geq j_2;$ and either ($j'_3\leq j_3$) or ($j'_3>j_3$ and $j'_2\geq j_2+{1\over 2}$).
\item  the second sum   is taken over  a finite number of multi-indices $\bfj''=(j''_1,j''_2,j''_3,j''_4)$ as above   such that   either  ($j''_1< j_1$)
or ($j''_1=j_1$ and $j''_2\geq {1\over 4}+j_2$) or ($j''_1=j_1$ and $j''_3<j_3$).
\end{itemize}
\end{lemma}
\proof
By Propositions \ref{P:basic-admissible-estimates-I} and  \ref{P:basic-admissible-estimates-II},
there  are constants $c_3,c_4>0$ such that    $ c_3 \pi^*\omega +c_4\beta\geq 0$  on  $\pi^{-1}(V_0)\subset \E$  and that
for every $1\leq \ell \leq \ell_0,$
 the following inequalities hold  on $\U_\ell \cap \Tube(B,\bfr):$
\begin{equation}\label{e:admissible-estimates}
\begin{split}
|\tilde \tau_\ell^*(\varphi)-\varphi|&\leq c_3 \varphi^{3\over 2}\quad\text{and}\quad   |\tilde \tau_\ell^*(f)-f|\leq c_3 \varphi^{1\over 2},\\
  \pm\big(   \tilde \tau_\ell^*(\pi^*\omega)  -\pi^*\omega-H\big)^\sharp &\lesssim  c_3 \varphi^{1\over 2}  \pi^*\omega +c_4\varphi^{1\over 2} \beta,  
  \\
   \pm\big(   \tilde \tau_\ell^*(\hat\beta) -\hat\beta  \big)^\sharp &\lesssim c_3 \phi^{3\over 2}\cdot  \pi^*\omega +c_4\phi^{1\over 2}\cdot \hat\beta,\\
\pm\big(  \tilde \tau_\ell^*(\hat\alpha) -\hat\alpha \big)^\sharp &\lesssim  c_3  \phi^{3\over 2}\cdot\pi^*\omega +c_4\hat\beta +c_3\varphi^{1/2}\hat\alpha.
                \end{split}
\end{equation}
Here, on the first line  $f$ is an arbitrary  $\Cc^1$-smooth  function  on $\Tube(B,\bfr),$ and  on the second and third lines $H$ is some  form 
in the class $\Hc$ given in Definition  \ref{D:Hc}.
By  Theorem \ref{T:basic-admissible-estimates}, for every $1\leq \ell \leq \ell_0,$
 the following inequality hold  on $\U_\ell \cap \Tube(B,\bfr):$
\begin{equation}
 \label{e:admissible-estimates(2)}
 \begin{split}
  &\big \lbrace \big(\tilde\tau_\ell^*(\pi^*\omega)  -\pi^*\omega \big), \big( \tilde \tau_\ell^*(\hat \beta) -\hat\beta\big),
  \big(  \tilde \tau_\ell^*(\hat\alpha) -\hat\alpha \big) 
  \big\rbrace\\
  &\trianglelefteq  \big \lbrace \big(c_3\varphi^{1\over 2}  \pi^*\omega +c_4\varphi^{1\over 2} \beta\big), \big( c_3 \phi^{3\over 2}\cdot  \pi^*\omega +c_4\phi^{1\over 2}\cdot \hat\beta\big),  \big(   c_3 \pi^*\omega +c_4\hat\beta +c_3\varphi^{1/4}\hat\alpha\big)\big\rbrace.
 \end{split}
\end{equation}
Next, we will explain how  to apply Lemma \ref{L:basic-positive-difference-bis}. We come back the statement of this lemma. 
 Let $R_1,\ldots,R_{k-p}$ be the $k-p$ forms among $\{\pi^*\omega,\hat \beta,  \hat\alpha\}$  which appear in the integral of $ I_{\bfj}$ in 
 \eqref{e:I_bfj}.  So  setting  $R:=R_1\wedge \ldots \wedge R_{k-p}.$ we get $$I_{\bfj}(s,r)=\int_{\Tube(B,s,r)}\tau_*T\wedge \varphi^{j_2}(c_1-c_2\varphi)^{j_4} R.$$  
 Now  we define  $R'_1,\ldots,R'_{k-p}$ as  follows. If  $R_j=\pi^*\omega$ set $R'_j:= \pi^*\omega+\hat\beta,$
  otherwise   $R_j\in\{\hat\alpha,\hat\beta\}$ and set $R'_j:=R_j.$
 Let   $S_1,\ldots,S_{k-p}$ be  the  corresponding positive $(1,1)$-form associated to $R_1,\ldots,R_{k-p}$ respectively on the RHS of  each  line of \eqref{e:admissible-estimates}.
 Let $S'_1,\ldots,S'_{k-p}$  be  the  corresponding positive $(1,1)$-form associated to $R_1,\ldots,R_{k-p}$ respectively on the RHS  of \eqref{e:admissible-estimates(2)}.
 Let $H_1,\ldots,H_{k-p}$ be    the  corresponding  real $(1,1)$-forms  associated to $R_1,\ldots,R_{k-p}$ respectively on the LHS of  each  of the last three lines  of \eqref{e:admissible-estimates}.  Observe  that  $H_j=0$ for $S_j=\hat\alpha$ and also  for $S_j=\hat\beta$ (see the  last two lines  of  
 \eqref{e:admissible-estimates}).  We also  check easily that $\varphi^{1\over 2} R_j\lesssim S_j.$

 Let $f$ be  either the function  $\varphi$  or the  function $c_1-c_2\varphi.$

Fix  a  constant $0<c_5<1$  small  enough.  Let  $y\in  \U_\ell$  with  $0<\theta_\ell(y)< c_5.$ Since  $\Sigma_{1\leq\ell\leq\ell_0}\pi^*\theta_\ell=1$
on  $\Tube(B,\bfr),$   we may find  $1\leq \ell'\leq \ell_0$ 
and  a  small open neighborhood $\U_y$  of $y$ in  $\U$  such that for $x\in\U_y,$ we have that $\theta_{\ell'}(x)>  c_5.$
Moreover, since $\tilde\tau_{\ell'}\circ \tilde\tau^{-1}_\ell=\tau_{\ell'}\circ \tau^{-1}_\ell$ is a holomorphic admissible map, we can check using \eqref{e:admissible-estimates} that
there is a constant $c_6>0$ such that  for $x\in\U_y,$
 $R_j(x)\leq  c_6(\tilde\tau_{\ell'}\circ \tilde\tau^{-1}_\ell)^* R'_j(x)$ and  that $S_j(x)\leq  c_6(\tilde\tau_{\ell'}\circ \tilde\tau^{-1}_\ell)^* S_j(x)$  and that
 $S'_j(x)\leq  c_6(\tilde\tau_{\ell'}\circ \tilde\tau^{-1}_\ell)^* S'_j(x).$   
Hence, we are in the  position to apply Lemma \ref{L:basic-positive-difference-bis} (2).
There is a constant 
$c$ that depends on $c_5,c_6$  and $\ell_0$ such that
\begin{equation}\label{e:L:spec-wedge}
\left   |                \langle \tau_*T ,\ind_{\Tube(B,s,r)}R\rangle  -\langle T^\hash_{s,r} ,R\rangle \right |^2\\
\leq 
  c\cdot\, \sum_{\ell=1}^{\ell_0}  \sum\limits_{I,J,K}\sum_{j=0}^{|I|}
  \Ic_{j,I,J,K}.
  \end{equation}
Here, 
for $0\leq  j\leq |I|,$
  \begin{equation}\begin{split}\label{e:Ic_j,I,J,K}
  \Ic_{j,I,J,K}&:=\big( 
 \int(\ind_{\Tube(B,s,r)}\circ \tilde\tau_\ell) (\pi^*\theta_\ell)(\tau_\ell)_* T 
\wedge R'_K\wedge  \pi^*\omega^j\wedge \hat\beta^{|I|-j} \wedge S_J\wedge S'_{(I\cup J\cup K)^\bfc}\big)\\
 &  \cdot    \big( \int   (\ind_{\Tube(B,s,r)}\circ \tilde\tau_\ell) (\pi^*\theta_\ell)(\tau_\ell)_* T        \wedge R'_K \wedge \pi^*\omega^{|I|-j}\wedge \hat\beta^{j} \wedge S_J\wedge S'_{(I\cup J\cup K)^\bfc}\big),
\end{split}
\end{equation} 
and the  sum $\sum_{I,J,K}$ is  taken over all $I,J,K\subset \{1,\ldots,k-p\}$  such that $R_j=\pi^*\omega$ for $j\in I,$ and that  $I,J,K$ are mutually disjoint,
   and $|(I\cup J\cup K)^\bfc|$ is  even, and  $K\not=\{1,\ldots,k-p\}.$
 
 Pick a family $(I,J,K)$  as above.
 Observe that  the  above condition on $I,J,K$  implies that at least one of the three sets $I,J$ and $(I\cup J\cup K)^\bfc$ is  non-empty.

Consider the case  where $I=\varnothing.$ So  either $J\not=\varnothing$ or $(I\cup J\cup K)^\bfc\not=\varnothing.$
Since the  RHS in all inequalities of \eqref{e:admissible-estimates} and \eqref{e:admissible-estimates(2)}  either do not contain  any  term $\hat\alpha$  or contains  $\hat\alpha$ with coefficient at least $\varphi^{1\over 4},$
the  exponent of $\hat\alpha$  in each term  in
 $ \Ic_{j,I,J,K}$ given by  \eqref{e:Ic_j,I,J,K}  must be either  $< j_1$
or  is  equal to $j_1$ but $j_2$ increases by at least  ${1\over 4}.$  Moreover,  $\phi\omega\lesssim \hat\beta$ and $\phi\hat\alpha\lesssim \hat\beta.$ Therefore, we infer  that 
\begin{equation*}
  \Ic_{j,I,J,K}
=\big(I^\hash_{\bfj'}(s,r)\big)^2   ,
\end{equation*}
for some  $\bfj'=(j'_1,j'_2,j'_3,j'_4)$    with  $j'_1\leq  j_1$  and $j'_2\geq j_2,$ 
  and  either  ($j'_1< j_1$)
or ($j'_1=j_1$ and $j'_2\geq {1\over 4}+j_2$), and either ($j'_3\leq j_3$) or ($j'_3>j_3$ and $j'_2\geq j_2+{1\over 2}$).

Consider the case  where $I\not=\varnothing.$  We obtain that 
\begin{equation*}
  \Ic_{j,I,J,K}
=I^\hash_{\bfj'}(s,r) I^\hash_{\bfj''}(s,r) ,
\end{equation*}
for some  $\bfj'=(j'_1,j'_2,j'_3,j'_4)$ and  $\bfj''=(j''_1,j''_2,j''_3,j''_4)$  with  $j'_1\leq  j_1$  and $j'_2\geq j_2,$ and either ($j'_3\leq j_3$) or ($j'_3>j_3$ and $j'_2\geq j_2+{1\over 2}$),
  and    ($j''_1=j_1$ and $j''_3<j_3$).

Combining   both above cases and estimate
  \eqref{e:L:spec-wedge},
the result follows.
\endproof 

Fix  an open neighborhood $\bfW$ of $\partial B$ in $X$ with $\bfW\subset \bfU.$
Recall  the class $\widetilde\CL^{1,1}_p(\bfU,\bfW)$  given in Definition
 \ref{D:sup}. Recall  the intermediate Lelong means $\nu_{j,q}(T,B,r,\tau)$ introduced in Subsection \ref{SS:Intermediate-Lelong-means}.
 For $0<r\leq \bfr,$ define
 \begin{equation}\label{e:nu_tot}
  \nu_\tot(T,B,r,\tau):=\sum_{(j,q):\ 0\leq q\leq\min(k-p,k-l),\  0\leq  j\leq \min(\upm,k-p-q)}|\nu_{j,q}(T,B,r,\tau)|.
 \end{equation}

 \begin{theorem}\label{T:Lc-finite}
   There is a constants $c_7>0$ such that for every positive closed  current $T$ on $\bfU$ belonging to the class $\widetilde\CL^{1,1}_p(\bfU,\bfW),$  we have 
   $$  \Kc_{j,q}(T,r) \leq  c_7 \nu_\tot(T,B,r,\tau)$$
   for $0\leq q\leq k-l$ and $0\leq j\leq k-p-q.$
   In particular, by increasing $c_7$ if necessary, we  have that
   $\Kc_{j,q}(T,\bfr)<c_7.$
 \end{theorem}
\proof Since  the mass of $T$  on $\bfU$ is $\leq 1,$ there is a constant $c$ independent of $T$ such that  $0\leq \nu_\tot(T,B,\bfr,\tau)\leq c.$ Therefore,
the  second  assertion is an immediate consequence of the  first one.  So  we only need to prove  the first  assertion.
The proof is  divided into three steps.

\noindent {\bf Step 1:} {\it The case  $q=0.$}

In this case  there  is  no factor $\hat\alpha$  appearing in $\Kc_{j,0}(T,\bfr)$ Since  the  forms  $\omega$ and $\hat\beta$ are  positive smooth,     there is a constant $c_7$ such that for $0\leq r\leq \bfr,$
$$
\Kc_{j,0}(T,r)=\sum_{\ell=1}^{\ell_0}\int\limits_{(\Tube(B,r)\cap \U_\ell) \setminus V  } (\pi^*\theta_\ell)\cdot (\tau_\ell)_*( T|_{\bfU_\ell})\wedge \pi^*\omega^j\wedge \hat\beta^{k-p-j} \leq c_7  \nu_\tot(T,B,r,\tau)  .
$$
This proves the theorem  for $q=0.$

\noindent {\bf Step 2:} {\it The general strategy and  a useful estimate (inequality \eqref{e:P-Lc-bullet_finite(8)} below).}
 
 The general strategy is  to prove  the proposition by  increasing induction on $q$ 
 with $0\leq q\leq k-l.$ 
 In the proof $\bfr$ is  a  fixed   but sufficiently small  constant.
Fix  $0\leq q_0\leq k-l.$ 
Suppose   that the proposition is true  for all $q,j$ with $q<q_0.$ 
We need to show that   the proposition is  also true  for all $q,j$ with   $q\leq q_0.$ 
Let $0\leq j_0\leq \min(\upm,k-p-q_0).$   Set  $j'_0:=k-p-q_0-j_0\geq 0.$
Consider
\begin{equation}\label{e:P-Lc-bullet_finite(0)}
\begin{split}
 \Kc^{-}_{j,q}(T,s,r)&:=\sum\limits_{\text{either}\ (q'<q)\ \text{or} (q'=q\ \text{and}\ j'<j )}\Kc_{j',q'}(T,s,r),
 \quad
 \Kc^{+}_{j,q}(T,s,r):=\sum\limits_{j':\ j<j'\leq k-p-q }\Kc_{j',q}(T,s,r),
 \\
 \Kc_{q}(T,s,r)&:=\sum\limits_{q'\leq q}\Kc_{j,q'}(T,s,r).
\end{split}
\end{equation}
We  define   $\Kc^{-}_{j,q}(T,r),$ $ \Kc^{+}_{j,q}(T,r)$ and   $\Kc_{q}(T,r)$ similarly.

We may assume without loss of generality that $T$ is  $\Cc^1$-smooth
and  let $s,r\in  [0,\bfr]$ with $s<r.$   Note that
\begin{equation*}
 d [(\tau_*T)\wedge \pi^*\omega^{j_0}]\wedge\beta^{j'_0}
=[(\tau_*dT)\wedge \pi^*\omega^{j_0}]\wedge\beta^{j'_0}=0 ,
\end{equation*}
where the last equality holds  as $T$ is  closed. Therefore,
applying Theorem  \ref{T:Lelong-Jensen-smooth-closed}  to  $\tau_*T\wedge \pi^*(\omega^{j_0})\wedge \beta^{j'_0},$   we get that
\begin{equation}\label{e:P-Lc-bullet_finite(1)}
\begin{split}
{1\over  r^{2q_0}}\int_{\Tube(B,r)}\tau_*T\wedge \pi^*(\omega^{j_0})\wedge \beta^{k-p-j_0}
-{1\over  s^{2q_0}}\int_{\Tube(B,s)}\tau_*T\wedge \pi^*(\omega^{j_0})\wedge \beta^{k-p-j_0}\\
= \Vc\big(\tau_*T\wedge \pi^*(\omega^{j_0})\wedge \beta^{j'_0},s,r\big)+\int_{\Tube(B,s,r)}\tau_*T\wedge \pi^*(\omega^{j_0})\wedge \beta^{j'_0}\wedge\alpha^{q_0}.
\end{split}
\end{equation}
Moreover,  by Theorem \ref{T:vertical-boundary-closed}, we have the following estimate independently of $T:$
\begin{equation}\label{e:P-Lc-bullet_finite(2)}
\Vc\big(\tau_*T\wedge \pi^*(\omega^{j_0})\wedge \beta^{j'_0},s,r\big)=O(r).
\end{equation}
Therefore,
when  $s\to 0+,$  applying  Lemma \ref{L:Lelong-smooth-forms} and Theorem   \ref{T:Lelong-Jensen-smooth-closed} (1),  equality \eqref{e:P-Lc-bullet_finite(1)} becomes
\begin{equation}\label{e:P-Lc-bullet_finite(3)}
{1\over  r^{2q_0}}\int_{\Tube(B,r)}\tau_*T\wedge \pi^*(\omega^{j_0})\wedge \beta^{k-p-j_0}
 \geq O(r)+ \int_{\Tube(B,r)}\tau_*T\wedge \pi^*(\omega^{j_0})\wedge \beta^{j'_0}\wedge\alpha^{q_0}.
\end{equation}
Hence, for $\bfr$ small enough, there is a constant $c_7>0$ independent of $T$ such that for $0\leq r\leq \bfr,$
\begin{equation}\label{e:P-Lc-bullet_finite(4)}
\int_{\Tube(B,r)}\tau_*T\wedge \pi^*(\omega^{j_0})\wedge \beta^{j'_0}\wedge\alpha^{q_0}\leq  c_7r+\nu_{j_0,q_0}(T,B,r,\tau).
\end{equation}
In the remainder of Step 2, we will use \eqref{e:P-Lc-bullet_finite(4)} in order to establish an useful estimate. 

Recall from \eqref{e:hat-alpha} and \eqref{e:hat-beta} that
\begin{equation*} \alpha=\hat\alpha- c_1 \pi^*\omega-c_2\beta
=\hat\alpha -c_2\hat \beta +(c_2\varphi-c_1)\pi^*\omega
\quad\text{and}\quad  \beta=\hat\beta-c_1\varphi\cdot  \pi^*\omega.
 \end{equation*}
So we get  that  
\begin{eqnarray*} 
\beta^{j'_0}\wedge \alpha^{q_0}
&=&(\hat\beta-c_1\varphi\cdot  \pi^*\omega)^{j'_0}\wedge(\hat\alpha -c_2\hat \beta +(c_2\varphi-c_1)\pi^*\omega)^{q_0} \\
&=&
\hat\beta^{j'_0}\wedge\hat \alpha^{q_0}
+\sum_{j_1,j'_1,j''_1}^q {j'_0\choose  j'_1}{q_0 \choose j_1}{q_0-j_1 \choose j''_1}\\
&\qquad& \cdot\,\hat\beta^{j'_1+j''_1}(-c_1\varphi\pi^*\omega)^{j'_0-j'_1}\wedge 
((c_2\varphi-c_1)\pi^*\omega)^{q_0-j_1-j''_1}\wedge \hat\alpha^{j_1},
\end{eqnarray*}
where the last sum is  taken over all $(j_1,j'_1,j''_1)$   such that $0\leq j'_1\leq j'_0$ and $0\leq j_1,j''_1\leq q_0$ such that
$j_1+j''_1\leq  q_0$ and 
$(j'_1,j_1)\not=(j'_0,q_0).$  Using this 
and the first equality of  \eqref{e:I_bfj}, we rewrite the last integral on  the RHS of  \eqref{e:P-Lc-bullet_finite(1)} 
 as
\begin{equation}\label{e:P-Lc-bullet_finite(5)}
 \begin{split}
&\int_{\Tube(B,s,r)}\tau_*T\wedge \pi^*(\omega^{j_0})\wedge \beta^{j'_0}\wedge\alpha^{q_0}
  = I_{q_0,0,j_0,0}(T,s,r)\\&+\sum_{j_1,j'_1,j''_1} {j'_0\choose  j'_1}{q_0 \choose j_1}
  {q_0-j_1 \choose j''_1}(-c_1)^{j'_0-j'_1}(-1)^{q_0-j_1-j''_1}
 I_{j_1, j'_0-j'_1, q_0+j_0+j'_0-j_1-j'_1-j''_1,q_0-j_1-j''_1}(T,s,r).
\end{split}
\end{equation}
Let $s$ tend to $0+.$ Using \eqref{e:P-Lc-bullet_finite(2)} and \eqref{e:P-Lc-bullet_finite(4)} and increasing $c_7$ if necessary, we  deduce from the above  equality that
\begin{equation*}
 \begin{split}
  I_{q_0,0,j_0,0}(T,r)+\sum_{j_1,j'_1,j''_1} {j'_0\choose  j'_1}{q_0 \choose j_1}
  {q_0-j_1 \choose j''_1}(-c_1)^{j'_0-j'_1}(-1)^{q_0-j_1-j''_1}\\
 \cdot I_{j_1, j'_0-j'_1, q_0+j_0+j'_0-j_1-j'_1-j''_1,q_0-j_1-j''_1}(T,r)\leq  c_7r+\nu_{j_0,q_0}(T,B,r,\tau).
\end{split}
\end{equation*}
We rewrite  this  inequality as follows:
\begin{equation}\label{e:P-Lc-bullet_finite(6)}
 \Ic_1+\Ic_2+\Ic_3\leq c_7r+\nu_{j_0,q_0}(T,B,r,\tau),
\end{equation}
where
\begin{equation*}
 \begin{split}
 \Ic_1&:=I^\hash_{q_0,0,j_0,0}(T,r)+\sum_{j'_1,j''_1,j_1}{j'_0\choose  j'_1}{q_0 \choose j_1}
  {q_0-j_1 \choose j''_1}(-c_1)^{j'_0-j'_1}(-1)^{q_0-j_1-j''_1}\\
 &\cdot I^\hash_{j_1, j'_0-j'_1, q_0+j'_0-j_1-j''_1,q_0-j_1-j''_1}(T,r), \\
 \Ic_2&:= I_{q_0,0,j_0,0}(T,r)-I^\hash_{q_0,0,j_0,0}(T,r), \\
 \Ic_3&:=  \sum_{j'_1,j''_1,j_1} {j'_0\choose  j'_1}{q_0 \choose j_1}
  {q_0-j_1 \choose j''_1}(-c_1)^{j'_0-j'_1}(-1)^{q_0-j_1-j''_1}\\
  &\cdot\big(I_{j_1, j'_0-j'_1, q_0+j'_0-j_1-j''_1,q_0-j_1-j''_1}(T,r)-I^\hash_{j_1, j'_0-j'_1, q_0+j_0+j'_0-j_1-j'_1-j''_1,q_0-j_1-j''_1}(T,r)\big)  .
\end{split}
\end{equation*}
Consider an arbitrary term
$I^\hash_{j_1, j'_0-j'_1, q_0+j_0+j'_0-j_1-j'_1-j''_1,q_0-j_1-j''_1}(T,r)$ in the sum on the  expression of $\Ic_1.$
Observe that $q_0+j_0+j'_0-j_1-j'_1-j''_1=j_0+(j'_0-j'_1)+(q_0-j_1-j''_1)\geq j_0+0+0=j_0.$ Moreover, 
if the equality holds  then $j_1<q_0$  because $(j'_1,j_1)\not=(j'_0,q_0),$ and
hence  the term is $\lesssim \Kc_{q_0-1}(T,r). $
If the  equality does not hold then either $j'_0-j'_1>0$ or $q_0-j_1-j''_1 >0,$ and hence  the term
is either $\lesssim  r^2\Kc^+_{j_0,q_0}(T,r)$ or $\lesssim \Kc_{q_0-1}(T,r).$
In all cases, we get that
$$I^\hash_{j_1, j'_0-j'_1, q_0+j'_0-j_1-j''_1,q_0-j_1-j''_1}(T,r)\lesssim cr^2\Kc^+_{j_0,q_0}(T,r)+\Kc_{q_0-1}(T,r).$$
Consequently, we get that
\begin{equation}\label{e:P-Lc-bullet_finite(6bis)}
  |\Ic_1-I^\hash_{q_0,0,j_0,0}(T,r)|\leq  c r^2\Kc^+_{j_0,q_0}(T,r)+c\Kc_{q_0-1}(T,r). 
\end{equation}
Applying Lemma \ref{L:spec-wedge}  to  each  difference term  in $\Ic_2$ and $\Ic_3$ 
yields that
\begin{equation} \label{e:P-Lc-bullet_finite(7)}
|I_\bfi(r)- I^\hash_\bfi(r)|^2 \leq c\big(\sum_{\bfi'} I^\hash_{\bfi'}(r)\big)\big ( \sum_{\bfi''} I^\hash_{\bfi''}(r) \big).  
\end{equation}
Here, on the LHS  $\bfi=(i_1,i_2,i_3,i_4)$  is  either $(q_0,0,j_0,0)$ or $(j_1, j'_0-j'_1, q_0+j_0+j'_0-j_1-j'_1-j''_1,q_0-j_1-j''_1)$  with $j_1,j'_1,j''_1$ as above, 
and on the RHS:
\begin{itemize} \item[$\bullet$] the first sum  is taken over a finite number of multi-indices    $\bfi'=(i'_1,i'_2,i'_3,i'_4)$ as above  such that  $i'_1\leq  i_1$  and $i'_2\geq i_2;$ 
\item  the second sum   is taken over  a finite number of multi-indices $\bfi''=(i''_1,i''_2,i''_3,i''_4)$ as above   such that   either  ($i''_1< i_1$)
or ($i''_1=i_1$ and $i''_2\geq {1\over 4}+i_2$) or ($i''_1=i_1$ and $i''_3<i_3$).
\end{itemize}
 Observe that when $\bfr$ is  small  enough,  $c_1-c_2\varphi\approx  1$ and $\varphi \lesssim \bfr^2\ll 1$ on $\Tube(B,\bfr).$
 Therefore, $I_{i_1,i_2,i_3,i_4}(T,r)\leq c I_{i_1,0,i_3,0}(T,r)$ for a constant $c>0$ independent of $T$ and $0<r\leq \bfr.$
 Consequently, the first  sum on the RHS of \eqref{e:P-Lc-bullet_finite(7)} is  bounded from above by a constant times  $\Kc_{q_0}(T,r),$
 whereas  the second sum  is bounded from above by a constant times $\Kc^-_{j_0,q_0}(T,r)+ r^{1\over 2} \Kc_{q_0}(T,r).$
 In fact the factor  $r^{1\over 2}$ comes  from  $\varphi^{1\over 4}$ because $\varphi\lesssim r^2$ on $\Tube(B,r).$    
Consequently,   we infer from \eqref{e:P-Lc-bullet_finite(6)}--\eqref{e:P-Lc-bullet_finite(6bis)} that there is a constant $c>0$ such that
\begin{equation}
 \label{e:upper-bound-I-hash}
\begin{split}
I^\hash_{q_0,0,j_0,0}(T,r)&\leq  cr+ |\nu_{j_0,q_0}(T,B,r,\tau)|+ cr^2 \Kc^{+}_{j_0,q_0}(T,r)+c\Kc_{q_0-1}(T,r) \\ 
&+c\sqrt{\Kc_{q_0}(T,r)} \sqrt{ \Kc^-_{j_0,q_0}(T,r)+r^{1\over 2} \Kc_{q_0}(T,r)  }.
\end{split}
\end{equation}
Hence, 
\begin{equation}\label{e:P-Lc-bullet_finite(8)}
\begin{split}
I^\hash_{q_0,0,j_0,0}(T,r)&\leq c_7\big(r+ \Kc_{q_0-1}(T,r)+  |\nu_{j_0,q_0}(T,B,r,\tau)|+r^2\Kc^{+}_{j_0,q_0}(T,r)\\&
+r^{1\over 4} \Kc_{q_0}(T,r)+  \sqrt{\Kc_{q_0}(T,r)} \sqrt{ \Kc^-_{j_0,q_0}(T,r)}
\big).
\end{split}
\end{equation}
This is the  desired estimate of  Step 2.


\noindent {\bf Step 3:} {\it  End of the  proof.}
 
 Suppose that $\Kc_{q_0}(r)\leq  c_7 \nu_\tot(T,B,r,\tau)$  for $0<r\leq \bfr,$
 where  $\Kc_{q_0}(r)$ is  defined in  \eqref{e:P-Lc-bullet_finite(0)}.
 
 When $j_0=0,$ we see that $\Kc^-_{j_0,q_0}(T,r)=\Kc_{q_0-1}(T,r),$ and hence \eqref{e:P-Lc-bullet_finite(8)} becomes
 \begin{equation}\label{e:P-Lc-bullet_finite(9)}
 \begin{split}
I^\hash_{q_0,0,0,0}(T,r)&\leq  c_7\big(r+\Kc_{q_0-1}(T,r) +|\nu_{0,q_0}(T,B,r,\tau)|+r^2 \Kc^{+}_{0,q_0}(T,r)\\
&+r^{1\over 4} \Kc_{q_0}(T,r)+  \sqrt{\Kc_{q_0}(T,r)} \sqrt{ \Kc_{q_0-1}(T,r)}\big).
\end{split}
\end{equation}
Observe that    
\begin{equation}\label{e:P-Lc-bullet_finite(10)}
\Kc^{-}_{j,q_0}(T,r) = \Kc_{q_0-1}(T,r)+\sum_{m=0}^{j-1}I^\hash_{q_0,0,m,0}(T,r). 
\end{equation}
Consequently, applying  \eqref{e:P-Lc-bullet_finite(8)} for $j_0=1$ and  hence inserting \eqref{e:P-Lc-bullet_finite(9)} and    estimate
\eqref{e:P-Lc-bullet_finite(10)} for $j=1$
into the resulting  inequality, we get that
\begin{equation*}
\begin{split}
I^\hash_{q_0,0,1,0}(T,r)
&\leq c_7\big(r+ \Kc_{q_0-1}(T,r)+|\nu_{1,q_0}(T,B,r,\tau)|+r^2 \Kc^{+}_{0,q_0}(T,r) \\
&+r^{1\over 4} \Kc_{q_0}(T,r)+  \sqrt{\Kc_{q_0}(T,r)} \sqrt{ \Kc^-_{1,q_0}(T,r)}
\big)\\
&\leq  c_7  \Big(r+ \Kc_{q_0-1}(T,r)+|\nu_{1,q_0}(T,B,r,\tau)|+r^{1\over 8} \Kc_{q_0}(T,r)\\
&+  \big[(\Kc_{q_0}(T,r))^{1\over 4}+ r^{1\over 4}+  (\Kc_{q_0-1}(T,r))^{1\over 4} +|\nu_{0,q_0}(T,B,r,\tau)|^{1\over 4}+(r^2 \Kc^{+}_{0,q_0}(T,r)  )^{1\over 4}
\big]^4\\
&-\Kc_{q_0}(T,r)
\Big) .
\end{split}
\end{equation*}
Using that $\Kc^{+}_{0,q_0}(T,r)=\Kc_{1,q_0}(T,r)+\Kc^{+}_{1,q_0}(T,r),$ 
and by Remark \ref{R:I_bfj},
  $\Kc_{1,q}(T,r)=I^\hash_{(q,0,1,0)}(r),$ 
we deduce that
\begin{equation*}
\begin{split}
I^\hash_{q_0,0,1,0}(T,r)
&\leq c_7\Big(r+ \Kc_{q_0-1}(T,r)+|\nu_{0,q_0}(T,B,r,\tau)|+|\nu_{1,q_0}(T,B,r,\tau)| 
 +r^{1\over 8} \Kc_{q_0}(T,r)\\
 &+  \big[(\Kc_{q_0}(T,r))^{1\over 4}+ r^{1\over 4}+  (\Kc_{q_0-1}(T,r))^{1\over 4} +|\nu_{0,q_0}(T,B,r,\tau)|^{1\over 4}+|\nu_{1,q_0}(T,B,r,\tau)|^{1\over 4}\\&+(r^2 \Kc^{+}_{1,q_0}(T,r)  )^{1\over 4}
\big]^4
-\Kc_{q_0}(T,r)
\Big) 
.
\end{split}
\end{equation*}
Set  $m_0:=k-p-q_0.$ Note that $\Kc^{+}_{m_0,q_0}(T,r)=0.$ We continue  this  process for $1\leq j\leq m_0$ and  obtain  that
\begin{multline*}
I^\hash_{q_0,0,j,0}(T,r)
\leq  c_7  \Big[r+ \sum_{j=0}^{m_0}|\nu_{j,q_0}(T,B,r,\tau)|+r^{1\over 2^{j+2}} \Kc_{q_0}(T,r)\\
+  \Big(\big[\Kc_{q_0}(T,r))^{1\over 2^{j+1}}+ r^{1\over 2^{j+1}}+ (\Kc_{q_0-1}(T,r))^{1\over 2^{j+1}}+ \big(\sum_{j=0}^{m_0}|\nu_{j,q_0}(T,B,r,\tau)|\big)^{1\over 2^{j+1}}\big]^{2^{j+1}}
-\Kc_{q_0}(T,r)\Big)
\Big] .
\end{multline*}
Note that  
\begin{equation*}
\Kc_{q_0}(T,r) =\Kc_{q_0-1}(T,r)+\sum_{j=0}^{m_0} I^\hash_{q_0,0,j,0}(T,r).
\end{equation*}
This, combined with the previous  estimates, implies by increasing $c_7$ that 
\begin{multline*}
\Kc_{q_0}(T,r) 
\leq  c_7  \Big(r + \sum_{j=0}^{m_0}|\nu_{j,q_0}(T,B,r,\tau)|+r^{1\over 2^{m_0+2}} \Kc_{q_0}(T,r)\\
+  \big[(\Kc_{q_0}(T,r))^{1\over 2^{m_0+1}}+  r^{1\over 2^{m_0+1}}+(\Kc_{q_0-1}(T,r))^{1\over 2^{m_0+1}}+ \big(\sum_{j=0}^{m_0}|\nu_{j,q_0}(T,B,r,\tau)|\big)^{1\over 2^{m_0+1}}  \big]^{2^{m_0+1}}
-\Kc_{q_0}(T,r)\Big)
 .
\end{multline*}
Recall from the  assumption of Step 3 that $\Kc_{q_0-1}(T,r)<c_7.$
Introduce the positive variable $$t:={(\Kc_{q_0}(T,r))^{1\over 2^{m_0+1}}\over  \big(r+\Kc_{q_0-1}(T,r)+\sum_{j=0}^{m_0}|\nu_{j,q_0}(T,B,r,\tau)|\big)^{1\over 2^{m_0+1}}}.$$
Dividing the both side  of the inequality   by $r+\Kc_{q_0-1}(T,r)+\sum_{j=0}^{m_0}|\nu_{j,q_0}(T,B,r,\tau)|,$ we infer that
$P_r(t)\leq 0,$ where 
\begin{multline*}
 P_r(t):= \big(1-c_7 r^{1\over  2^{m_0+2} }\big)t^{2^{m_0+1}}\\-\sum_{m=0}^{2^{m_0+1}-1}
 \big[{r^{1\over 2^{m_0+1}}+(\Kc_{q_0-1}(T,r))^{1\over 2^{m_0+1}}+ \big(\sum_{j=0}^{m_0}|\nu_{j,q_0}(T,B,r,\tau)|\big)^{1\over 2^{m_0+1}}  \over \big(r+\Kc_{q_0-1}(T,r)+\sum_{j=0}^{m_0}|\nu_{j,q_0}(T,B,r,\tau)|\big)^{1\over 2^{m_0+1}} }   \big]^{2^{m_0+1}-m} t^m\\
 -c_7 \big[{r + \sum_{j=0}^{m_0}|\nu_{j,q_0}(T,B,r,\tau)|\over r+\Kc_{q_0-1}(T,r)+\sum_{j=0}^{m_0}|\nu_{j,q_0}(T,B,r,\tau)|  }\big].
\end{multline*}
Observe  that $P_r$ is a real polynomial of degree $ 2^{m_0+1}$  whose leading coefficient is
$1-c_7 r^{1\over  2^{m_0+2} }$ and other coefficients are constant (dependent on $r$),  but all these other coefficients are of  modulus  $\lesssim 1.$
When $r>0$ is  small enough, the leading coefficient ranges within the interval $ ({1\over  2},1],$ and $P_r(t)\leq 0.$
Consequently,  $t$ is  uniformly bounded independently of $r.$
This  proves that  $$\Kc_{q_0}(T,r)\leq c_7\big(r+\Kc_{q_0-1}(T,r)+\sum_{j=0}^{m_0}|\nu_{j,q_0}(T,B,r,\tau)|\big)$$ for some constant $c_7>0$ independent of  $T$ and $r.$  The conclusion of Step 3 follows.

The proof of the proposition is thereby completed.

\endproof

\begin{proposition}\label{P:Lc-finite-bis}
  For  $0<r_1<r_2\leq \bfr,$ there is  a constant $c_8>0$ such that for every $q\leq  \min(k-p,k-l)$ and 
    every positive closed  current $T$ on $\bfU$ belonging to the class $\widetilde\CL^{1,1}_p(\bfU,\bfW),$   we have the following estimate:
 $$|\kappa_{k-p-q}(T,{r_1\over \lambda},{r_2\over \lambda},\tau)|<c_8\sum\limits_{0\leq q'\leq q,\ 0\leq j'\leq \min(\upm,k-p-q')} \Kc_{j',q'}(T,{r_1\over \lambda},{r_2\over \lambda} )\qquad\text{for}\quad  \lambda>1.$$
 \end{proposition}
\proof  
Fix  $0\leq q_0\leq \min(k-p,k-l)$ and  set $j_0:=k-p-q_0.$
We will adapt some parts in the proof of Theorem \ref{T:Lc-finite} for $s:={r_1\over \lambda}$ and  $r:={r_2\over \lambda}.$
Note that in the present context  $j'_0=0.$
Arguing as  in the proof of \eqref{e:P-Lc-bullet_finite(5)}, we  get that
\begin{equation}\label{e:P-Lc-bullet_finite-bis(1)}
 \begin{split}
&\kappa_{k-p-q_0}(T,{r_1\over \lambda},{r_2\over \lambda},\tau)
  = I_{q_0,0,j_0,0}(T,{r_1\over \lambda},{r_2\over \lambda})\\&+\sum_{j''_1,j_1} {q_0 \choose j_1}
  {q_0-j_1 \choose j''_1}(-1)^{q_0-j_1-j''_1}
 I_{j_1, 0, q_0+j_0-j_1-j''_1,q_0-j_1-j''_1}(T,{r_1\over \lambda},{r_2\over \lambda}),
\end{split}
\end{equation}
where the  sum is  taken over all $(j_1,j''_1)$   such that  $0\leq j_1<q_0$ and $j''_1\leq q_0$ and
$j_1+j''_1\leq  q_0.$
Similarly  as in \eqref{e:P-Lc-bullet_finite(6)}, we rewrite  the expression on the RHS of \eqref{e:P-Lc-bullet_finite-bis(1)} as the sum $
 \Ic_1+\Ic_2+\Ic_3,
$
where
\begin{equation*}
 \begin{split}
 \Ic_1&:=I^\hash_{q_0,0,j_0,0}(T,{r_1\over \lambda},{r_2\over \lambda})+\sum_{j''_1,j_1}{q_0 \choose j_1}
  {q_0-j_1 \choose j''_1}(-1)^{q_0-j_1-j''_1}\\
 &\cdot I^\hash_{j_1, 0, q_0+j_0-j_1-j''_1,q_0-j_1-j''_1}(T,{r_1\over \lambda},{r_2\over \lambda}), \\
 \Ic_2&:= I_{q_0,0,j_0,0}(T,r)-I^\hash_{q_0,0,j_0,0}(T,{r_1\over \lambda},{r_2\over \lambda}), \\
 \Ic_3&:=  \sum_{j''_1,j_1} {q_0 \choose j_1}
  {q_0-j_1 \choose j''_1}(-1)^{q_0-j_1-j''_1}\\
  &\cdot\big(I_{j_1, 0, q_0+j_0-j_1-j''_1,q_0-j_1-j''_1}(T,{r_1\over \lambda},{r_2\over \lambda})-I^\hash_{j_1, 0, q_0+j_0-j_1-j''_1,q_0-j_1-j''_1}(T,{r_1\over \lambda},{r_2\over \lambda})\big)  .
\end{split}
\end{equation*}
Observe that  $\Ic_1$ is bounded from above  by a constant times $\sum_{0\leq q'\leq q,\ 0\leq j'\leq \min(\upm,k-p-q')} \Kc_{j',q'}(T,{r_1\over \lambda},{r_2\over \lambda} ).$
Applying Lemma \ref{L:spec-wedge}  to  each  difference term  in $\Ic_2$ and $\Ic_3$ 
 as  in the proof of \eqref{e:P-Lc-bullet_finite(7)} and the argument  which follows   \eqref{e:P-Lc-bullet_finite(7)} yields the  same
 estimate for   $\Ic_2$ and $\Ic_3.$ This, combined with \eqref{e:P-Lc-bullet_finite-bis(1)}, 
 gives the result.
\endproof
We conclude  this  subsection  with the following finiteness result of the mass indicators $\Mc_j.$  Its  proof will be postponed until Subsection 
\ref{SS:Other-charac-Lelong} below.
 \begin{proposition}\label{P:Lc-finite-closed}
 There is a   constant $c_9>0$ such that for  
  every positive closed  current $T$ on $\bfU$ belonging to the class $\widetilde\CL^{1,1}_p(\bfU,\bfW),$ 
   we have
   $\Mc_j(T,r)<c_{9}$ for $0\leq j\leq \upm$ and $0 <r \leq \bfr.$
 \end{proposition}

\subsection{Existence of Lelong numbers}

This  subsection  is  devoted to the proof of assertions (1)--(4) of  Theorem \ref{T:Lelong-closed-Kaehler}.

\proof[Proof of assertion (1) of Theorem \ref{T:Lelong-closed-Kaehler}]
 
First assume  that the current $T$ is a closed $\Cc^1$-smooth form. 
 Since $\omega$ is  K\"ahler,
 we have for $1\leq j\leq   \upm$ that
\begin{equation*}
d [(\tau_*T)\wedge \pi^*\omega^{j}]=d(\tau_*T)\wedge \pi^*\omega^{j}
=(\tau_*dT)\wedge \pi^*\omega^{j}=0 .
\end{equation*}
Applying Theorem  \ref{T:Lelong-Jensen-closed}  to  $\tau_*T\wedge \pi^*(\omega^{j})$  and using   the  above  equality,  we get  that
\begin{equation}\label{e:Lelong-closed-all-degrees-bis-asser(1)-Kaehler}
  \nu_{j}(T,B,r_2,\tau)- \nu_{j}(T,B,r_1,\tau)= \int_{\Tube(B,r_1,r_2)}\tau_*T\wedge \pi^*(\omega^{j})\wedge\alpha^{k-p-j}+\Vc(\tau_*T\wedge \pi^*(\omega^{j}),r_1,r_2).
\end{equation}
On  the other hand,  since $j\geq \lowm$ we  get that  $k-p-j\leq  k-l.$  Therefore, we can apply
Theorem \ref{T:vertical-boundary-closed} to the current $\tau_*T\wedge \pi^*(\omega^{j}),$ which gives  that  $\Vc(\tau_*T\wedge \pi^*(\omega^{j}),r_1,r_2)=O(r_2).$ This proves  assertion (1) in  the special case where $T$ is $\Cc^1$-smooth.

Now  we consider the general case where $T$ is a general positive closed  $(p,p)$-current such that $T=T^+-T^-,$ where $T^\pm$ are approximable
along $B\subset V$ by positive closed $\Cc^1$-smooth  $(p,p)$-forms $(T^\pm_n)$  with $\Cc^1$-control on boundary. 
So $T^+_n\to T^+$ and $T^-_n\to T^-$ as $n$ tends to infinity.  
By the  previous case applied  to $T^\pm_n,$ we get that
\begin{equation*}
\nu_j( T^\pm_n,B,r_2,\tau)- \nu_j(T^\pm_n,B,r_1,\tau)=\kappa_j(T^\pm_n,B,r_1,r_2,\tau) +O(r_2).
\end{equation*}
Letting  $n$ tend to infinity, we infer that
\begin{equation*}
\nu_j( T^\pm,B,r_2,\tau)- \nu_j(T^\pm,B,r_1,\tau)=\kappa_j(T^\pm,B,r_1,r_2,\tau)+O(r_2).
\end{equation*}
This implies assertion (1)  since $T=T^+-T^-.$
\endproof
\endproof

\proof[Proof of assertion (2) of Theorem \ref{T:Lelong-closed-Kaehler}]  Let $q:=k-p-j.$ Fix $r_1,\ r_2\in  (0,\bfr]$ with $r_1<r_2/2.$
Applying   Proposition  \ref{P:Lc-finite-bis}  yields that
 \begin{equation}\label{e:T:Lelong-closed(1)}|\kappa_{j}(T^\pm,{r_1\over \lambda} ,{r_2\over \lambda},\tau)|<c_8\sum_{0\leq q'\leq q,\ 0\leq j'\leq \min(\upm,k-p-q')} \Kc_{j',q'}(T,{r_1\over \lambda},{r_2\over \lambda} )\qquad\text{for}\quad  \lambda>1.
 \end{equation}
 On the other  hand,  since there is an  $M\in\N$ such that
 $$  1\leq \#\left\lbrace n\in\N:\  y\in  \Tube(B, {r_1\over 2^n} ,{r_2\over 2^n})\right\rbrace \leq M\quad\text{for}\quad y\in\Tube(B,\bfr),$$
 it follows that
 $$
 \sum_{n=1}^\infty \sum_{0\leq q'\leq q,\ 0\leq j'\leq \min(\upm,k-p-q')} \Kc_{j',q'}(T,{r_1\over 2^n},{r_2\over 2^n} )\leq M
  \sum_{0\leq q'\leq q,\ 0\leq j'\leq \upm} \Kc_{j',q'}(T,\bfr ).
 $$
 By Theorem \ref{T:Lc-finite} the  RHS is  finite.  Therefore, we infer from \eqref{e:T:Lelong-closed(1)} that
\begin{eqnarray*}
\sum_{n=1}^\infty|\kappa_{j}(T,{r_1\over 2^n},{r_2\over 2^n},\tau)|&\leq&  \sum_{n=0}^\infty|\kappa_{j}(T^+,{r_1\over 2^n},{r_2\over 2^n},\tau)| +\sum_{n=0}^\infty|\kappa_{j}(T^-,{r_1\over 2^n},{r_2\over 2^n},\tau)|\\
&\leq &   M c_8\sum_{0\leq q'\leq q,\ 0\leq j'\leq \min(\upm,k-p-q')} \Kc_{j',q'}(T,\bfr )
<\infty.
\end{eqnarray*}
Now  we apply  Lemma \ref{L:elementary} (2)  to functions $f^\pm$ and $\epsilon^\pm$ given by  
$$
f^\pm(r):=\nu(T^\pm,B,r,\tau)\qquad\text{and}\qquad  \epsilon^\pm_\lambda:=2c_8    
 \lambda^{-1}+c_8\sum_{0\leq q'\leq q,\ 0\leq j'\leq \min(\upm,k-p-q')} \Kc_{j',q'}(T,{r_1\over \lambda},{r_2\over \lambda} ).
$$
By  assertion  (1) and inequality \eqref{e:T:Lelong-closed(1)},  we have  by increasing the constant $c_8$  if necessary:
$$
|f^\pm({r_2\over \lambda})-f^\pm({r_1\over\lambda})|=|\kappa_j(T^\pm,{r_1\over \lambda},{r_2\over \lambda},\tau )+O(\lambda^{-1})|\leq \epsilon_\lambda.
$$
Hence, assertion (2) follows.
\endproof

\proof[Proof of assertion (3) of Theorem \ref{T:Lelong-closed-Kaehler}]
By \eqref{e:Lelong-log-bullet-numbers}  and   assertion (1), we have 
\begin{eqnarray*}\
  \kappa^\bullet_j(T,B,r,\tau)= \limsup\limits_{s\to0+}   \kappa_j(T,B,s,r,\tau)&=& \nu_j(T,B,r,\tau)-\liminf_{s\to0+} \nu_j(T,B,s,\tau)\\
 &=&  \nu_j(T,B,r,\tau)- \nu_j(T,B,\tau),
\end{eqnarray*} 
where the last equality holds by assertion (2).   Consequently, we  infer  from  assertion  (2) again that
\begin{equation*}
 \lim\limits_{r\to0+}\kappa^\bullet_j(T,B,r,\tau)= \lim\limits_{r\to0+}\nu_j(T,B,r,\tau)- \nu_j(T,B,\tau)=0.
\end{equation*}
\endproof

\proof[Proof of assertion (4) of Theorem \ref{T:Lelong-closed-Kaehler}]
It is  similar  to the  proof of assertion (4) of Theorem  \ref{T:top-Lelong-closed}.
\endproof

\subsection{Other characterizations of Lelong numbers and independence of admissible maps}
\label{SS:Other-charac-Lelong}

Consider the following mass indicators, for $\lowm \leq j\leq \upm:$ 
 Following   \eqref{e:mass-indicators} 
 we   define the mass indicators even when $\tau$ is not necessarily holomorphic
 \begin{equation}\label{e:other-mass-indicator} \hat\nu_j(T,r):={1\over  r^{2(k-p-j)}}\int\limits_{\Tube(B,r)}
 \tau_* T\wedge (\beta+c_1r^2\pi^*\omega)^{k-p-j}\wedge  \pi^*\omega^j.
 \end{equation}
 We also consider the following new mass indicators, where  $T^\hash$ and $T^\hash_r$ are given in \eqref{e:T-hash} and \eqref{e:T-hash_r}:
\begin{equation}\label{e:other-mass-indicators-bis}
 \begin{split}
 \widehat\Mc^\hash_j(T,r)&:={1\over   r^{2(k-p-j)}}\int\limits_{\Tube(B,r)}
 T^\hash \wedge (\beta +c_1r^2\pi^*\omega)^{k-p-j}\wedge  \pi^*\omega^j,\\
 \Mc^\hash_j(T,r)&:={1\over r^{2(k-p-j)}}\int
 T^\hash_r \wedge (\beta+c_1r^2\pi^*\omega)^{k-p-j}\wedge  \pi^*\omega^j
 .
 \end{split}
 \end{equation}


\begin{lemma} \label{L:Mc-hash-vs-nu} There is a  constant $c>0$ such that for every $\lowm\leq j\leq \upm$ and $0<r\leq\bfr:$
$$|\Mc^\hash_j(T,r)- \hat\nu_j(T,r)| 
\leq    cr\sum^{\upm}_{q=\lowm}  \Mc^\hash_q(T,r).$$
\end{lemma}

\proof
By Propositions \ref{P:basic-admissible-estimates-I}, 
there  are constants $c_3,c_4>0$ such that    $ c_3 r^2\pi^*\omega +c_4\beta\geq 0$  on  $\Tube(B,r)$  for $0<r\leq \bfr,$ and that
for every $1\leq \ell \leq \ell_0,$
 the following inequalities hold  on $\U_\ell \cap \Tube(B,r)$  for $0<r\leq \bfr:$
\begin{equation}\label{e:admissible-estimates-bis}
\begin{split}
  \pm\big(   \tilde \tau_\ell^*(\pi^*\omega)  -\pi^*\omega-H\big)^\sharp &\lesssim  c_3 r  \pi^*\omega +c_4r( \beta+c_1r^2\pi^*\omega),\\   
   \pm\big(   \tilde \tau_\ell^*(\beta+c_1r^2\pi^*\omega) - (\beta+c_1r^2\pi^*\omega) \big)^\sharp &\lesssim c_3 r^3  \pi^*\omega +c_4r (\beta+c_1r^2\pi^*\omega).
                \end{split}
\end{equation}
Here,  on the LHS of the  first line,  $H$ is some  form 
in the class $\Hc$ given in Definition  \ref{D:Hc}.

On the other hand, by Theorem \ref{T:basic-admissible-estimates}, 
for every $1\leq \ell \leq \ell_0,$
 the following inequalities hold  on $\U_\ell \cap \Tube(B,r)$  for $0<r\leq \bfr:$
 \begin{equation}\label{e:admissible-estimates-bis(2)}
\begin{split}
  &\left\lbrace \big(\tilde\tau_\ell^*(\pi^*\omega)  -\pi^*\omega\big),   \big(    \tilde \tau_\ell^*(\beta+c_1r^2\pi^*\omega) -(\beta+c_1r^2\pi^*\omega)  \big) \right\rbrace\\
  &\trianglelefteq    \left\lbrace  \big( c_3 r  \pi^*\omega +c_4r (\beta+c_1r^2\pi^*\omega)\big),\big(c_3 r^3  \pi^*\omega +c_4r (\beta+c_1r^2\pi^*\omega)\big)\right\rbrace.
  \end{split}
\end{equation}
Next, we will explain how  to apply Lemma \ref{L:basic-positive-difference-bis}. We come back the statement of this lemma. 
 Let $R_1,\ldots,R_{k-p}$ be the $k-p$ forms among $\{\pi^*\omega,\beta+c_1r^2\pi^*\omega\}$  which appear in the integral of $\hat\nu_j(T,r)$ in 
 \eqref{e:other-mass-indicator}.  So  setting  $R:=R_1\wedge \ldots \wedge R_{k-p}.$ we get $\hat\nu_j(T,r)={1\over  r^{2(k-p-j)}}\int_{\Tube(B,r)}\tau_*T\wedge R.$
 Now  we define  $R'_1,\ldots,R'_{k-p}$ as  follows. If  $R_j=\pi^*\omega$ set $R'_j:= \pi^*\omega+\beta,$
  otherwise   $R_j =\beta+c_1r^2\pi^*\omega$ and set $R'_j:=R_j.$
 Let   $S_1,\ldots,S_{k-p}$ be  the  corresponding positive $(1,1)$-form associated to $R_1,\ldots,R_{k-p}$ respectively on the RHS  of \eqref{e:admissible-estimates-bis}.
 Let $S'_1,\ldots,S'_{k-p}$  be  the  corresponding positive $(1,1)$-form associated to $R_1,\ldots,R_{k-p}$ respectively on the RHS   of \eqref{e:admissible-estimates-bis(2)}.
 Let $H_1,\ldots,H_{k-p}$ be    the  corresponding  real $(1,1)$-forms  associated to $R_1,\ldots,R_{k-p}$ respectively on the  LHS of  each  inequality of \eqref{e:admissible-estimates-bis}.

   Arguing  as  in the proof  of   Lemma  \ref{L:spec-wedge},
   we are in the  position to  apply Lemma \ref{L:basic-positive-difference-bis}.
There is a constant 
$c$ that depends on $c_5,c_6$  and $\ell_0$ such that
\begin{equation}\label{e:L:spec-wedge-bis(1)}
\left   |      {1\over  r^{2(k-p-j)}}          \langle \tau_*T ,\ind_{\Tube(B,r)}R\rangle  -\langle T^\hash_{r} ,R\rangle \right |^2\\
\leq 
  \big({1\over  r^{2(k-p-j)}}\big)^2 c\cdot\, \sum_{\ell=1}^{\ell_0}  \sum\limits_{I,J,K}\sum_{\iota=0}^{|I|}
  \Ic_{\iota,I,J,K}.
  \end{equation}
Here, 
for $0\leq  \iota\leq |I|,$ $\Ic_{\iota,I,J,K}$ is  given in 
  \eqref{e:Ic_j,I,J,K},
and the  sum $\sum_{I,J,K}$ is  taken over all $I,J,K\subset \{1,\ldots,k-p\}$  such that  $I,J,K$ are mutually disjoint,
   and $|(I\cup J\cup K)^\bfc|$ is  even, and  $K\not=\{1,\ldots,k-p\}.$ 
 
 Pick a family $(I,J,K)$  as above.
  As in the proof  of   Lemma  \ref{L:spec-wedge},  the  above condition on $I,J,K$  implies that at least one of the three sets $I,J$ and $(I\cup J\cup K)^\bfc$ is  non-empty.
 We rewrite \eqref{e:Ic_j,I,J,K}  as  
\begin{equation*}\begin{split}
  &{\Ic_{\iota,I,J,K}\over (r^{2(k-p-j)})^2}\\
  &=
 \big({1\over  r^{2(k-p-j)}}\int(\ind_{\Tube(B,r)}\circ \tilde\tau_\ell) (\pi^*\theta_\ell)(\tau_\ell)_* T 
\wedge R'_K\wedge  \pi^*\omega^\iota\wedge \hat\beta^{|I|-\iota} \wedge S_J\wedge S'_{(I\cup J\cup K)^\bfc}\big)\\
 &  \cdot    \big( {1\over  r^{2(k-p-j)}}\int   (\ind_{\Tube(B,r)}\circ \tilde\tau_\ell) (\pi^*\theta_\ell)(\tau_\ell)_* T        \wedge R'_K \wedge \pi^*\omega^{|I|-\iota}\wedge \hat\beta^{\iota} \wedge S_J\wedge S'_{(I\cup J\cup K)^\bfc}\big),
\end{split}
\end{equation*} 
 Consider the case where  either $J\not=\varnothing$ or $(I\cup J\cup K)^\bfc\not=\varnothing.$
Observe  that
the exponent of $r$ in the  coefficient of $\pi^*\omega$   (resp. in the coefficient of $\beta+c_1r^2\pi^*\omega$ on the LHS
of  all inequalities of \eqref{e:admissible-estimates-bis} does not exceed the corresponding  exponent of $r$ on the RHS    minus  1.
Therefore, we infer  from the above equality and  \eqref{e:other-mass-indicators-bis} that
\begin{equation}\label{e:L:spec-wedge-bis(2)}
\big({1\over  r^{2(k-p-j)}}\big)^2 \Ic_{\iota,I,J,K} 
\leq  cr^2\big(\sum_{q=\lowm}^\upm\Mc^\hash_q(T,r)   \big)^2 .
\end{equation}

Consider the case where   $I\not=\varnothing$ and   $J=\varnothing$ and $(I\cup J\cup K)^\bfc=\varnothing.$  The above equality becomes
 \begin{equation*}\begin{split}
  {\Ic_{\iota,I,J,K}\over r^{2|I|}(r^{2(k-p-j)})^2}
  &=
 \big({1\over  r^{2(k-p-j+|I|-\iota)}}\int(\ind_{\Tube(B,r)}\circ \tilde\tau_\ell) (\pi^*\theta_\ell)(\tau_\ell)_* T 
\wedge R'_K\wedge  \pi^*\omega^\iota\wedge \hat\beta^{|I|-\iota} \big)\\
 &  \cdot    \big( {1\over  r^{2(k-p-j+\iota)}}\int   (\ind_{\Tube(B,r)}\circ \tilde\tau_\ell) (\pi^*\theta_\ell)(\tau_\ell)_* T        \wedge R'_K \wedge \pi^*\omega^{|I|-\iota}\wedge \hat\beta^{\iota} \big).
\end{split}
\end{equation*} 
Consequently,
\begin{equation*}
\big({1\over  r^{2(k-p-j)}}\big)^2\Ic_{\iota,I,J,K}
\leq   cr^{2|I|}\sum_m \Mc^\hash_{m} (T,r)\Mc^\hash_{m+|I|-2\iota}(T,r).
\end{equation*}
Inserting this and \eqref{e:L:spec-wedge-bis(2)} in \eqref{e:L:spec-wedge-bis(1)}, the result follows.
\endproof


 \begin{proposition}\label{P:comparison-Mc}
 For $\lowm \leq j\leq \upm,$  we have that $$ \lim_{r\to 0+} \Mc^\hash_j(T,r)=\lim_{r\to 0+}  \hat\nu_j(T,r)=\sum\limits_{q=0}^{k-p-j}  {k-p-j\choose q} c^q_1\nu_{j+q}(T,B,\tau).$$
 \end{proposition}
\proof Using formula \eqref{e:other-mass-indicator}
 and arguing as in the proof of Lemma \ref{L:hat-nu}, we obtain  the following   identity which is similar to identity  \eqref{e:hat-nu-vs-nu}:
 \begin{equation*}
  \hat\nu_j(T,r)=\sum\limits_{q=0}^{k-p-j}  {k-p-j\choose q} c^q_1\nu_{j+q}(T,B,r,\tau).  
 \end{equation*}
 Next,  letting  $r$ tend to $0$
  in this identity, we infer
  from  Theorem  \ref{T:Lelong-closed-Kaehler} (2) that   
 \begin{equation}\label{e:Mc_j-vs-nu-bis}
  \lim_{r\to 0+} \hat\nu_j(T,r)=\sum\limits_{q=0}^{k-p-j}  {k-p-j\choose q} c^q_1\nu_{j+q}(T,B,\tau).  
 \end{equation}
 This proves the second identity of the  proposition.
 
 It remains  to show  the first  identity. 
 Applying  Lemma   \ref{L:Mc-hash-vs-nu}  yields that 
   there is a constant $c>0$ such that for $0<r\leq\bfr,$
 \begin{equation}\label{e:Mc_j-induction}
 \big|\sum_{j=\lowm}^\upm  \Mc^\hash_j(T,r)-\sum_{j=\lowm}^\upm\hat\nu_j(T,r)\big|\leq  c r \sum_{j=\lowm}^{\upm}\Mc^\hash_j(T,r).
 \end{equation}
This, combined  with   \eqref{e:Mc_j-vs-nu-bis}, implies that there is a constant  $c>0$ such  that
$$
 \sum_{j=\lowm}^{\upm}\Mc^\hash_j(T,r)\leq c\qquad\text{for}\qquad 0<r\leq\bfr.  
$$
Therefore, we infer from Lemma   \ref{L:Mc-hash-vs-nu} that
$
|\Mc^\hash_j(T,r)- \hat\nu_j(T,r)|\leq  cr
$ for $\lowm\leq j\leq \upm.$  Letting $r$ tend to $0,$  
the first  identity of the proposition  
 follows.
 \endproof
 
 \proof[Proof of assertion (5) of Theorem \ref{T:Lelong-closed-Kaehler}]
 Let $\tau$ and $\tau'$ be two   strongly admissible  maps.
 For $1\leq \ell\leq \ell_0$  we define  $\tilde\tau'_\ell:=\tau'\circ \tau^{-1}_\ell$  according to formula \eqref{e:tilde-tau_ell}.
 So   $\tilde\tau'_\ell$ is defined in the same was as  $\tilde\tau_\ell$
 using $\tau'$ instead of $\tau.$
 Similarly, we define $T^{'\hash}$ and $T^{'\hash}_r$ according to formulas \eqref{e:T-hash} and \eqref{e:T-hash_r}  by using $\tilde\tau'_\ell$ instead of $\tilde\tau_\ell.$ 
 Similarly, we  define
 $
 \widehat\Mc'^\hash_j(T,r)$ and 
 $\Mc'^\hash_j(T,r)$  according to formula \eqref{e:other-mass-indicators-bis} by using  $T^{'\hash}$ and $T^{'\hash}_r$  instead of $T^{\hash}$ and $T^{\hash}_r.$

 We need to show  that 
 \begin{equation}\label{e:tau-tau'}\nu_j(T,B,\tau)=\nu_j(T,B,\tau')\qquad\text{for}\qquad \lowm \leq j\leq \upm.
 \end{equation}
 By \eqref{e:admissible-estimates-bis} 
there  are constants $c_3,c_4>0$ such that    $ c_3 r^2\pi^*\omega +c_4\beta\geq 0$  on  $\Tube(B,r)$  for $0<r\leq \bfr,$ and that
for every $1\leq \ell \leq \ell_0,$
 the following inequalities hold  on $\U_\ell \cap \Tube(B,r)$  for $0<r\leq \bfr:$
\begin{equation}\label{e:admissible-estimates-bis-bis}
\begin{split}
  \pm\big(   \tilde \tau_\ell^*(\pi^*\omega)  -(\tilde \tau'_\ell)^*(\pi^*\omega)-H\big) &\lesssim  c_3 r  \pi^*\omega +c_4r( \beta+c_1r^2\pi^*\omega),\\   
\pm\big(   \tilde \tau_\ell^*(\beta+c_1r^2\pi^*\omega) - (\tilde \tau'_\ell)^*(\beta+c_1r^2\pi^*\omega)  \big) &\lesssim c_3 r^3  \pi^*\omega +c_4r (\beta+c_1r^2\pi^*\omega).
                \end{split}
\end{equation}
Here,    $H$ is some  form 
in the class $\Hc$ given in Definition  \ref{D:Hc}.

 By \eqref{e:admissible-estimates-bis(2)} for every $1\leq \ell \leq \ell_0,$
 the following inequality holds  on $\U_\ell \cap \Tube(B,r)$  for $0<r\leq \bfr:$
\begin{equation}\label{e:admissible-estimates-bis-bis(2)}
\begin{split}
& \left\lbrace \big(\tilde\tau_\ell^*(\pi^*\omega)  -(\tilde \tau'_\ell)^*(\pi^*\omega)\big),\big(    \tilde \tau_\ell^*(\beta+c_1r^2\pi^*\omega) -(\tilde \tau'_\ell)^*(\beta+c_1r^2\pi^*\omega)  \big)  \right\rbrace\\
&\trianglelefteq \left\lbrace \big(  c_3 r  \pi^*\omega +c_4r (\beta+c_1r^2\pi^*\omega)\big),\big(  c_3 r^3  \pi^*\omega +c_4r (\beta+c_1r^2\pi^*\omega)\big)\right\rbrace.
 \end{split}
 \end{equation}
Using \eqref{e:admissible-estimates-bis-bis}--\eqref{e:admissible-estimates-bis-bis(2)} and arguing as in the  proof of Lemma  \ref{L:Mc-hash-vs-nu} we can show that 
there is a  constant $c>0$ such that for every $\lowm\leq j\leq \upm$ and $0<r\leq\bfr:$
$$|\Mc^\hash_j(T,r)- {\Mc'}^\hash_j(T,r)|\leq   cr\sum_{j=\lowm}^{\upm}  \Mc^\hash_j(T,r).$$
Thus by  Proposition  \ref{P:comparison-Mc},  $|\Mc^\hash_j(T,r)- {\Mc'}^\hash_j(T,r)|\leq  cr.$ 
 So by this  proposition  again, we get that
\begin{equation*}
 \lim_{r\to 0}{\Mc'}^\hash_{j}(T,r)  =\lim_{r\to 0+}  \hat\nu_{j_0}(T,r)=\sum\limits_{q=0}^{k-p-j_0}  {k-p-j_0\choose q} \nu_{j_0+q}(T,B,\tau).  
 \end{equation*}
 Hence, for $\lowm\leq j\leq \upm,$ we have 
 $$
 \sum\limits_{q=0}^{k-p-j}  {k-p-j\choose q} \nu_{j+q}(T,B,\tau)=\sum\limits_{q=0}^{k-p-j}  {k-p-j\choose q} \nu_{j+q}(T,B,\tau').
 $$
 These  equalities imply \eqref{e:tau-tau'}.
 The proof is  thereby completed.
 \endproof


 \proof[Proof of assertion (6) of Theorem \ref{T:Lelong-closed-Kaehler}]
 Applying Proposition \ref{P:comparison-Mc}
 for $j=\upm$ yields  that 
 $$ \lim_{r\to 0+} \Mc^\hash_\upm(T,r)=\lim_{r\to 0+}  \hat\nu_\upm(T,r)=\nu_{\upm}(T,B,\tau).$$
 By \eqref{e:other-mass-indicators-bis},  $\Mc^\hash_\upm(T,r)\geq 0$  for $0<r\leq\bfr.$  Hence, $\nu_{\upm}(T,B,\tau)\geq 0.$
 \endproof

 \begin{proposition}\label{P:comparison-Mc-bis}
 For $\lowm \leq j\leq \upm,$  we have that $$ \lim_{r\to 0+} \widehat\Mc^\hash_j(T,r)=\lim_{r\to 0+} \Mc^\hash_j(T,r)  .$$
 \end{proposition}
 \proof
 Observe that there is  a constant $c>0$ such that  for every $1\leq \ell\leq \ell_0$ and $0<r\leq \bfr,$ we have
 \begin{equation}
  \label{e:difference-charac-functions}
  \big| (\ind_{\Tube (B,r)})(\pi^*\theta_\ell)-(\ind_{\Tube (B,r)}\circ \tilde\tau_\ell)(\pi^*\theta_\ell)\big|\leq c (\ind_{\Tube (B,r-cr^2,r+cr^2)})(\pi^*\theta_\ell).
 \end{equation}
 Indeed,  for $y=(z,w)\in \Tube(B,r),$  writing  $y'=(z',w')=\tilde\tau_\ell(z,w),$ we have 
 $\|z-z'\|=O(\|z\|^2)=O(r^2)$ since $\tilde\tau_\ell$ is  admissible, and this  estimate implies  \eqref{e:difference-charac-functions}.
 Hence,  we infer that
 \begin{multline*}
  \big| (\ind_{\Tube (B,r)})(\pi^*\theta_\ell)\cdot (\tau_\ell)_*(T|_{\bfU_\ell})-(\ind_{\Tube (B,r)}\circ \tilde\tau_\ell)(\pi^*\theta_\ell)\cdot (\tau_\ell)_*(T|_{\bfU_\ell})\big|\\
  \leq c (\ind_{\Tube (B,r-cr^2,r+cr^2)}\circ \tilde\tau_\ell)(\pi^*\theta_\ell)\cdot (\tau_\ell)_*(T|_{\bfU_\ell}).
 \end{multline*}
Putting this together  with \eqref{e:other-mass-indicators-bis}, we get  that
 \begin{equation*}
 |\widehat\Mc^\hash_j(T,r)-\Mc^\hash_j(T,r)|\leq{ ( r+cr^2)^{2(k-p-j)}
 \Mc^\hash_j(T,r+cr^2)    -  ( r-cr^2)^{2(k-p-j)}
 \Mc^\hash_j(T,r-cr^2)  \over r^{2(k-p-j)}}
 .
 \end{equation*}
 By  Proposition \ref{P:comparison-Mc-bis}, the RHS tends to $0$ as $r\to 0.$  Hence,  the result follows.
 \endproof
 
 \proof[End of the  proof of Proposition \ref{P:Lc-finite-closed}]
 It  follows  from  the  definition of $\Mc_j$ and $\widehat\Mc^\hash_j$  in 
 \eqref{e:local-mass-indicators-bis}, \eqref{e:global-mass-indicators}
  and \eqref{e:other-mass-indicators-bis}
 that  there is a constant $c>0$ such that for every $\lowm\leq j\leq\upm$ and $0<r\leq\bfr$ and  every positive closed  current $T$ on $\bfU$ belonging to the class $\widetilde\CL^{1,1}_p(\bfU,\bfW),$  we have
 $\Mc_j(T,r)\leq  c\Mc^\hash_j(T,r).$  
  By Proposition \ref{P:comparison-Mc-bis}
 there is a constant $c'$  such that  $  \widehat\Mc^\hash_j(T,r)<c'.$  Choosing $c_9:=cc',$
   we obtain $\Mc_j(T,r)<c_{9}$ as  desired. 
 \endproof
\subsection{Another variant of top Lelong number}

We use the notation  introduced  at the beginning  of  Section \ref{S:Regularization}. 
We use the local setting  introduced in Subsection  \ref{SS:Local-setting} for  each $\U_\ell\subset \E$  with  $1\leq \ell\leq \ell_0.$
Namely, 
we use the coordinates $(z,w)\in\C^{k-l}\times \C^l.$ 
We may assume that  $\U_\ell$ has the form $\U_\ell=\U'_\ell\times \U''_\ell,$ where $\U'_\ell$ (resp. $\U''_\ell)$ are open neighborhood of $0'$ in $\C^{k-l}$ of  ($0''$ in $\C^l$).
Let $V=\{z=0\}=\U''$ and let    $\bfr>0$  such that  $\{\|z\|<\bfr\}\times B\Subset \U_\ell.$ 
Consider  the trivial  vector bundle $\pi:\ \E \to  \U''_\ell$ with  $\E\simeq  \C^{k-l}\times \U''_\ell.$ For $\lambda\in\C^*,$  let $a_\lambda:\ \E\to \E$ be the multiplication by  $\lambda$
on fibers, that is, 
$a_\lambda(z,w):=(\lambda z,w)$ for $(z,w)\in \E.$

Consider the positive closed $(1,1)$-forms
$$\beta=\omega_z:=\ddc \|z\|^2\quad\text{and}\quad  \omega=\omega_w:=\ddc\|w\|^2\quad\text{and}\quad \alpha=\theta_z:=\ddc \log\|z\|^2.$$ 
Define  
\begin{equation}\label{e:top-Lelong-number-variant}
\begin{split}
\nu_\top^\Uc(T,B,r)&:= \sum_{\ell =1}^{\ell_0} \kappa_\top((\tau_\ell)_*T,\U_\ell,r),\\
\nu_\top^\Uc(T,B) &:=\lim\limits_{r\to 0+} \nu_\top^\Uc(T,B,r).
\end{split}
\end{equation}

\begin{proposition}
 There is a constant $c>1$ depending only on $\Uc$ such that  for every positive closed current $T$ on $X$  we have that
 \begin{equation*}
 c^{-1} \nu_\top(T,B,r)  \leq \nu_\top^\Uc(T,B,r)\leq  c\nu_\top(T,B,r).
 \end{equation*}
In particular,  letting $r$ tend to $0$ we obtain that 
$$c^{-1} \nu_\top(T,B)  \leq \nu_\top^\Uc(T,B)\leq  c\nu_\top(T,B).$$
\end{proposition}
\proof Since  the proof is  not difficult, we leave it   to the interested reader.   \endproof

The  drawback of  this notion of the top Lelong number is that  it is not intrinsic.

 \section{Quasi-positivity and quasi-monotonicity of the generalized  Lelong numbers}
 \label{S:quasi-monotone}
In this  section  we    establish   the quasi-positivity and quasi-monotonicity of the generalized  Lelong numbers
of  positive closed  currents.
\begin{definition}\label{D:Dc-Lelong}\rm
Let $\Dc$ be  a family of   real numbers    $\Dc=\{ d_{jq}\in\R:\ 0\leq j\leq k-p-q\text{ and}\quad 0\leq q\leq  k-l\}.$   
For $0\leq q\leq k-l,$ 
consider the number
\begin{equation}\label{e:Dc-nu_q}
 \nu^\Dc_q(T,B,r,\tau):=\sum_{(j',q'):\ j'\leq k-p-q',\ q'\leq q}^\upm d_{j'q'}\nu_{j',q'}(T,B,r,\tau).
\end{equation}
Consider also  the number    
\begin{equation}\label{e:nu-Dc_tot}
\nu^\Dc_\tot(T,B,r,\tau):=\sum_{(j,q):\ j\leq k-p-q,\  0\leq q\leq k-l}^\upm d_{jq}\nu_{j,q}(T,B,r,\tau).
\end{equation}
Note that  $\nu^\Dc_\tot(T,B,r,\tau)=\nu^\Dc_{k-l}(T,B,r,\tau).$ 
\end{definition}


\begin{proposition}\label{P:nu_tot-monotone} 
Let  $0< r_1<r_2\leq\bfr.$ Then there are  a  family  $\Dc=\{ d_{jq}\in\R:\ 0\leq j\leq k-p-q,\ 0\leq q\leq  k-l\}$ and   a constant $c>0 $ depending  on $r_1$ and $r_2$  such that for every positive closed   current $T$ on $\bfU$ belonging to the class $\widetilde\CL^{2,2}_p(\bfU,\bfW),$ 
the  following inequality  hold for $0\leq q\leq k-l:$
\begin{equation*}
  \nu^\Dc_q\big(T,B,{r_1\over \lambda},\tau\big)\leq  \nu^\Dc_q(T,B,{r_2\over \lambda},\tau\big) + {c\over \lambda}
  \quad\text{for}\quad \lambda\gg 1.
 \end{equation*}
Moreover, for  every  $\epsilon>0$  we can choose $\Dc$ such that  $d_{k-p-q,q} <\epsilon^q d_{k-p-q+1,q-1}$ for $1\leq q\leq k-l.$
\end{proposition}

\begin{lemma}\label{L:basic-difference-nu_j,q_estimate} Given $0<r_1<r_2\leq \bfr,$ 
there is  a constant $c>0$ such that  for every positive closed current $T\in  \widetilde\CL^{1,1}_p(\bfU,\bfW)$
and  $0\leq q\leq k-l$ and $0\leq j\leq \min(\upm,k-p-q),$ the following inequality holds:
 \begin{multline*}
\nu_{j,q}\big(T,B,{r_2\over \lambda},\tau\big)  -\nu_{j,q}\big(T,B,{r_1\over \lambda},\tau\big)\\
\geq  \Kc_{j,q}\big (T,{r_1\over\lambda},{r_2\over \lambda}  \big)-  c\lambda^{-1}-  c\lambda^{1\over 2} \Kc_{q}\big (T,{r_1\over\lambda},{r_2\over \lambda} \big )-c\Kc_{q-1} (T,{r_1\over\lambda},{r_2\over \lambda})  
-c\sqrt{\Kc_{q}(T,{r_1\over\lambda},{r_2\over \lambda})} \sqrt{ \Kc^-_{j,q}(T,{r_1\over\lambda},{r_2\over \lambda})  }.
\end{multline*}
\end{lemma}
\proof
Fix  $0\leq q_0\leq k-l.$   
Let $0\leq j_0\leq \min(\upm,k-p-q_0).$   Set  $j'_0:=k-p-q_0-j_0\geq 0.$
We may assume without loss of generality that $T$ is  $\Cc^1$-smooth
and  let $s,r\in  [0,\bfr]$ with $s<r.$   Since $T$ is  closed, it follows that 
\begin{equation*}
 d [(\tau_*T)\wedge \pi^*\omega^{j_0}]\wedge\beta^{j'_0}
=[(\tau_*dT)\wedge \pi^*\omega^{j_0}]\wedge\beta^{j'_0}=0 .
\end{equation*}
Therefore,
applying Theorem  \ref{T:Lelong-Jensen-smooth-closed}  to  $\tau_*T\wedge \pi^*(\omega^{j_0})\wedge \beta^{j'_0},$   we get that
\begin{equation}\label{e:P-Lc-bullet_finite(1)-mono}
\begin{split}
&{\lambda^{2q_0}\over  r_2^{2q_0}}\int_{\Tube(B,{r_2\over \lambda})}\tau_*T\wedge \pi^*(\omega^{j_0})\wedge \beta^{k-p-j_0}
-{\lambda^{2q_0}\over  r_1^{2q_0}}\int_{\Tube(B,{r_1\over\lambda})}\tau_*T\wedge \pi^*(\omega^{j_0})\wedge \beta^{k-p-j_0}\\
&= \Vc\big(\tau_*T\wedge \pi^*(\omega^{j_0})\wedge \beta^{j'_0},{r_1\over\lambda},{r_2\over \lambda}\big)+\int_{\Tube(B,{r_1\over\lambda},{r_2\over \lambda})}\tau_*T\wedge \pi^*(\omega^{j_0})\wedge \beta^{j'_0}\wedge\alpha^{q_0}.
\end{split}
\end{equation}
Moreover,  by Theorem \ref{T:vertical-boundary-closed}, we have the following estimate independently of $T:$
\begin{equation}\label{e:P-Lc-bullet_finite(2)-mono}
\Vc\big(\tau_*T\wedge \pi^*(\omega^{j_0})\wedge \beta^{j'_0},{r_1\over\lambda},{r_2\over \lambda}\big)=O(\lambda^{-1}).
\end{equation}
Therefore,  there is a constant $c>0$ independent of $T$ such that for $\lambda\geq 1,$
\begin{equation}\label{e:P-Lc-bullet_finite(4)-mono}
\big|\int_{\Tube(B,{r_1\over\lambda},{r_2\over \lambda})}\tau_*T\wedge \pi^*(\omega^{j_0})\wedge \beta^{j'_0}\wedge\alpha^{q_0}
-\big (\nu_{j_0,q_0}\big(T,B,{r_2\over \lambda},\tau\big)   -  \nu_{j_0,q_0}\big(T,B,{r_1\over \lambda},\tau\big) \big)  \big| \leq  c\lambda^{-1}.
\end{equation} 
 Arguing as  in the proof of \eqref{e:P-Lc-bullet_finite(5)}, we obtain  that
\begin{equation}\label{e:P-Lc-bullet_finite(5)-mono}
 \begin{split}
&\int_{\Tube(B,{r_1\over\lambda},{r_2\over \lambda})}\tau_*T\wedge \pi^*(\omega^{j_0})\wedge \beta^{j'_0}\wedge\alpha^{q_0}
  = I_{q_0,0,j_0,0}(T,{r_1\over\lambda},{r_2\over \lambda})\\&+\sum_{j_1,j'_1,j''_1} {j'_0\choose  j'_1}{q_0 \choose j_1}
  {q_0-j_1 \choose j''_1}(-c_1)^{j'_0-j'_1}(-1)^{q_0-j_1-j''_1}
 I_{j_1, j'_0-j'_1, q_0+j_0+j'_0-j_1-j'_1-j''_1,q_0-j_1-j''_1}(T,{r_1\over\lambda},{r_2\over \lambda}).
\end{split}
\end{equation}
 Using \eqref{e:P-Lc-bullet_finite(2)} and \eqref{e:P-Lc-bullet_finite(4)} and increasing $c$ if necessary, we  deduce from the above  equality that
\begin{equation*}
 \begin{split}
 \big| I_{q_0,0,j_0,0}(T,r)+\sum_{j_1,j'_1,j''_1} {j'_0\choose  j'_1}{q_0 \choose j_1}
  {q_0-j_1 \choose j''_1}(-c_1)^{j'_0-j'_1}(-1)^{q_0-j_1-j''_1}\\
 \cdot I_{j_1, j'_0-j'_1, q_0+j_0+j'_0-j_1-j'_1-j''_1,q_0-j_1-j''_1}(T,r) -\big(\nu_{j_0,q_0}\big(T,B,{r_2\over \lambda},\tau\big)  -\nu_{j_0,q_0}\big(T,B,{r_1\over \lambda},\tau\big)\big)\big|\leq  c\lambda^{-1}.
\end{split}
\end{equation*}
As in the  proof of \eqref{e:P-Lc-bullet_finite(6)} we rewrite  this  inequality as follows:
\begin{equation}\label{e:P-Lc-bullet_finite(6)-mono}
 \big|\Ic_1+\Ic_2+\Ic_3   -\big(\nu_{j_0,q_0}\big(T,B,{r_2\over \lambda},\tau\big)  -\nu_{j_0,q_0}\big(T,B,{r_1\over \lambda},\tau\big)\big)\big|\leq c\lambda^{-1},
\end{equation}
where
\begin{equation*}
 \begin{split}
 \Ic_1&:=I^\hash_{q_0,0,j_0,0}(T,{r_1\over\lambda},{r_2\over \lambda} )+\sum_{j'_1,j''_1,j_1}{j'_0\choose  j'_1}{q_0 \choose j_1}
  {q_0-j_1 \choose j''_1}(-c_1)^{j'_0-j'_1}(-1)^{q_0-j_1-j''_1}\\
 &\cdot I^\hash_{j_1, j'_0-j'_1, q_0+j'_0-j_1-j''_1,q_0-j_1-j''_1}(T,{r_1\over\lambda},{r_2\over \lambda}), \\
 \Ic_2&:= I_{q_0,0,j_0,0}(T,{r_1\over\lambda},{r_2\over \lambda})-I^\hash_{q_0,0,j_0,0}(T,{r_1\over\lambda},{r_2\over \lambda}), \\
 \Ic_3&:=  \sum_{j'_1,j''_1,j_1} {j'_0\choose  j'_1}{q_0 \choose j_1}
  {q_0-j_1 \choose j''_1}(-c_1)^{j'_0-j'_1}(-1)^{q_0-j_1-j''_1}\\
  &\cdot\big(I_{j_1, j'_0-j'_1, q_0+j'_0-j_1-j''_1,q_0-j_1-j''_1}(T,{r_1\over\lambda},{r_2\over \lambda})-I^\hash_{j_1, j'_0-j'_1, q_0+j_0+j'_0-j_1-j'_1-j''_1,q_0-j_1-j''_1}(T,{r_1\over\lambda},{r_2\over \lambda})\big)  .
\end{split}
\end{equation*}
Arguing as in the proof of \eqref{e:P-Lc-bullet_finite(6bis)} we can show that
\begin{equation}\label{e:P-Lc-bullet_finite(6bis)-mono}
  |\Ic_1-I^\hash_{q_0,0,j_0,0}(T,{r_1\over\lambda},{r_2\over \lambda})|\leq  c \lambda^{-2}\Kc^+_{j_0,q_0}(T,{r_1\over\lambda},{r_2\over \lambda})+c\Kc_{q_0-1}(T,{r_1\over\lambda},{r_2\over \lambda}). 
\end{equation}
Applying Lemma \ref{L:spec-wedge}  to  each  difference term  in $\Ic_2$ and $\Ic_3$ 
yields that
\begin{equation} \label{e:P-Lc-bullet_finite(7)-mono}
|I_\bfi({r_1\over\lambda},{r_2\over \lambda})- I^\hash_\bfi({r_1\over\lambda},{r_2\over \lambda})|^2 \leq c\big(\sum_{\bfi'} I^\hash_{\bfi'}({r_1\over\lambda},{r_2\over \lambda})\big)\big ( \sum_{\bfi''} I^\hash_{\bfi''}({r_1\over\lambda},{r_2\over \lambda}) \big).  
\end{equation}
Here, on the LHS  $\bfi=(i_1,i_2,i_3,i_4)$  is  either $(q_0,0,j_0,0)$ or $(j_1, j'_0-j'_1, q_0+j_0+j'_0-j_1-j'_1-j''_1,q_0-j_1-j''_1)$  with $j_1,j'_1,j''_1$ as above, 
and on the RHS $\bfi'$ and $\bfi''$ are  described by the  two  properties $\bullet$ which  follow \eqref{e:P-Lc-bullet_finite(7)}.

 Consequently, the first  sum on the RHS of \eqref{e:P-Lc-bullet_finite(7)} is  bounded from above by a constant times  $\Kc_{q_0}(T,{r_1\over\lambda},{r_2\over \lambda}),$
 whereas  the second sum  is bounded from above by a constant times $\Kc^-_{j_0,q_0}(T,{r_1\over\lambda},{r_2\over \lambda})+ \lambda^{-1\over 2} \Kc_{q_0}(T,{r_1\over\lambda},{r_2\over \lambda}).$
 In fact the factor  $\lambda^{-1\over 2}$ comes  from  $\varphi^{1\over 4}$ because $\varphi\lesssim \lambda^{-2}$ on $\Tube(B,{r_1\over\lambda},{r_2\over \lambda}).$    
Therefore,   we infer from \eqref{e:P-Lc-bullet_finite(6)-mono}--\eqref{e:P-Lc-bullet_finite(6bis)-mono} that there is a constant $c>0$ such that
\begin{equation*}
\begin{split}
&\big|I^\hash_{q_0,0,j_0,0}(T,{r_1\over\lambda},{r_2\over \lambda})  -\big(\nu_{j_0,q_0}\big(T,B,{r_2\over \lambda},\tau\big)  -\nu_{j_0,q_0}\big(T,B,{r_1\over \lambda},\tau\big)\big)\big|\\
&\leq  c\lambda^{-1}+  c\lambda^{-2} \Kc^{+}_{j_0,q_0}(T,{r_1\over\lambda},{r_2\over \lambda} )+c\Kc_{q_0-1}(T,{r_1\over\lambda},{r_2\over \lambda}) \\ 
&+c\sqrt{\Kc_{q_0}(T,{r_1\over\lambda},{r_2\over \lambda})} \sqrt{ \Kc^-_{j_0,q_0}(T,{r_1\over\lambda},{r_2\over \lambda})+\big ( {r_2\over \lambda}\big)^{1\over 2} \Kc_{q_0}(T,{r_1\over\lambda},{r_2\over \lambda})  }.
\end{split}
\end{equation*}
 Since $\Kc_{j_0,q_0}(T,{r_1\over\lambda},{r_2\over \lambda})=I^\hash_{q_0,0,j_0,0}(T,{r_1\over\lambda},{r_2\over \lambda})\geq 0,$ it follows  that
 \begin{multline*}
\nu_{j,q}\big(T,B,{r_2\over \lambda},\tau\big)  -\nu_{j,q}\big(T,B,{r_1\over \lambda},\tau\big)
\geq  \Kc_{j,q}\big(T,{r_1\over\lambda},{r_2\over \lambda}  \big)-  c\lambda^{-1}-  c\lambda^{-2} \Kc^{+}_{j,q}(T,{r_1\over\lambda},{r_2\over \lambda} )\\-c\Kc_{q-1}(T,{r_1\over\lambda},{r_2\over \lambda}) 
-c\sqrt{\Kc_{q}(T,{r_1\over\lambda},{r_2\over \lambda})} \sqrt{ \Kc^-_{j,q}(T,{r_1\over\lambda},{r_2\over \lambda})
+\big({r_2\over\lambda}\big)^{1\over 2} \Kc_{q}(T,{r_1\over\lambda},{r_2\over \lambda})  }.
\end{multline*}
 As  $\Kc^{+}_{j,q}(T,{r_1\over\lambda},{r_2\over \lambda} )
 \leq  \Kc_{q}(T,{r_1\over\lambda},{r_2\over \lambda} ),$
 the last inequality implies
 the  desired conclusion  of the lemma when  we choose the constant $c$ large enough.
\endproof

\begin{lemma}\label{L:mu}
 For every  $q\geq 1$  and $\mu>0,$ there are $(q+2)$ numbers $\lambda_j$  ($0\leq j\leq q$)  and $\mu_0$  such that
 $1=\lambda_0>\lambda_1>\ldots>\lambda_q>0$ and $\mu_0>0$  and that   for  $t_1,\ldots,t_q\geq 0$ with $t_1+\ldots+t_q\leq 1,$
  we have  $P_{\lambda,\mu} (t)\geq \mu_0,$ where 
 \begin{multline*}
  P_{\lambda,\mu} (t):= \lambda_0t^2_1+\lambda_1[(t_1+t_2)^2-t^2_1]+\ldots+\lambda_{q-1}[(t_1+\ldots+t_q)^2-(t_1+\ldots+t_{q-1})^2]\\
  +\lambda_q[1-(t_1+\ldots+t_q)^2]-2\mu\lambda_1t_1-2 \mu\lambda_2(t_1+t_2)-\ldots-2\mu\lambda_q(t_1+\ldots+t_q). 
 \end{multline*}
\end{lemma}

\proof
We have that
\begin{eqnarray*}
  P_{\lambda,\mu} (t)&=& \sum_{j=1}^q  \big[ (\lambda_{j-1}-\lambda_j)   (t_1+\ldots+t_j)^2-2\mu\lambda_j(t_1+\ldots+t_j)\big]+\lambda_q\\
  &=&\sum_{j=1}^q  \big[ (\lambda_{j-1}-\lambda_j)  \big[ (t_1+\ldots+t_j)-{\mu\lambda_j\over  \lambda_{j-1}-\lambda_j}\big]^2
  +\big[\lambda_q  - \sum_{j=1}^q  {\mu^2\lambda^2_j\over  \lambda_{j-1}-\lambda_j} \big].
\end{eqnarray*}
Therefore, we only need  to show that with a  suitable choice of $\lambda_j$  ($0\leq j\leq q$)  and $\mu_0,$
\begin{equation*}
 \lambda_q  - \sum_{j=1}^q  {\mu^2\lambda^2_j\over  \lambda_{j-1}-\lambda_j}>0.
\end{equation*}
Write $\lambda_j=k_j\lambda_q$ for  $1\leq j\leq q.$ So $k_j>0$ and $k_q=1.$ The last inequality  is reduced  to
\begin{equation*}
 \mu^{-2}>{\lambda_1^2\over 1- \lambda_1}  + \sum_{j=1}^{q-1}  {k^2_{j+1}\over  k_{j}-k_{j+1}}.
\end{equation*}
So we only need to choose   $k_1>\ldots> k_q=1$ and $ \lambda_1>0$ such that
\begin{equation}\label{e:choice-lambda}
  {k^2_{j+1}\over  k_{j}-k_{j+1}}<{1\over \mu^2 q}\qquad\text{and}\qquad {\lambda_1^2\over 1- \lambda_1}<{1\over \mu^2 q}.
\end{equation}
We first fix $k_{q-1}>0$ large  enough such that ${1\over  k_{q-1}-1}<{1\over \mu^2 q}.$
Suppose that  $k_j$ is   already fixed.
Next, we choose   $k_{j-1}>0$ large  enough such that  ${k^2_{j}\over  k_{j-1}-k_{j}}<{1\over \mu^2 q}.$
After  having  determined $k_j$  for $1\leq j\leq  q,$ it remains to choose $\lambda_1>0$ small enough such that
the second  estimate of \eqref{e:choice-lambda} is fulfilled.
\endproof

\begin{lemma} \label{L:nu_0}
Given  $0<r_1<r_2\leq \bfr$ 
and numbers  $d_{j0}>0$   for $0\leq  j\leq \upm, $  there is   a constant $c_0>0$ such that 
for every positive closed current $T\in  \widetilde\CL^{1,1}_p(\bfU,\bfW),$  
 the following inequality holds:
 \begin{equation*} 
  \nu^\Dc_0\big(T,B,{r_2\over \lambda},\tau\big)- \nu^\Dc_0\big(T,B,{r_1\over \lambda},\tau\big) \geq - {c_0\over \lambda}
  +c^{-1}_0\Kc_0\big(T,  {r_1\over \lambda},{r_2\over \lambda}      \big)\quad\text{for}\quad \lambda\gg 1,
 \end{equation*}
\end{lemma}
\proof
Following the model of \eqref{e:other-mass-indicator}  and \eqref{e:other-mass-indicators-bis},  consider, for $0\leq j\leq\upm$ and $0<s<r\leq \bfr:$
\begin{eqnarray*}
\hat \kappa_{j,0}(T,s,r)&:= &\int_{\Tube(B,s,r)} \tau_*T \wedge (\beta+c_1r^2\pi^*\omega)^{k-p-j}\wedge \pi^*\omega^j,\\
\hat\kappa^\hash_{j,0}(T,s,r)&:=&\int\limits_{\Tube(B,s,r)}
 T^\hash_{s,r} \wedge (\beta +c_1r^2\pi^*\omega)^{k-p-j}\wedge  \pi^*\omega^j,
 \end{eqnarray*}
 where   and $T^\hash_{s,r}$ is given in  \eqref{e:T-hash_r}.
  We adapt  the proof of   Lemma   \ref{L:Mc-hash-vs-nu} 
So similarly as  in  \eqref{e:Mc_j-induction},  
we may find a constant $c>0$ such that for $\lambda\geq 1,$
 \begin{equation*}
 \big|\sum_{j=0}^\upm  d_{j,0}\hat\kappa^\hash_{j,0}(T,{r_1\over\lambda},{r_2\over\lambda})-\sum_{j=0}^\upm d_{j,0}\hat\kappa_j(T,{r_1\over\lambda},{r_2\over\lambda})\big|\leq  c \lambda^{-1}
 \sum_{j=0}^\upm  d_{j,0}\hat\kappa^\hash_{j,0}(T,{r_1\over\lambda},{r_2\over\lambda}).
 \end{equation*}
On the other hand,
\begin{equation*}
 \nu^\Dc_0\big(T,B,{r_2\over \lambda},\tau\big)- \nu^\Dc_0\big(T,B,{r_1\over \lambda},\tau\big)
 =\sum_{j=0}^\upm d_{j,0}\kappa_j(T,B,{r_1\over\lambda},{r_2\over\lambda}).
\end{equation*}
By  Lemma   \ref{L:Mc-hash-vs-nu} again, there is a constant $c>0$  that
\begin{equation*}
 \big|\sum_{j=0}^\upm d_{j,0}\kappa_j(T,B,{r_1\over\lambda},{r_2\over\lambda},\tau)-\sum_{j=0}^\upm d_{j,0}\hat\kappa^\hash_j(T,B,{r_1\over\lambda},{r_2\over\lambda},\tau)\big|\\
 \leq  c\lambda^{-1}\sum_{j=0}^\upm  d_{j,0}\hat\kappa^\hash_{j,0}(T,{r_1\over\lambda},{r_2\over\lambda}).
\end{equation*}
Observe  that there is a constant  $c'>0$ independent of $T$  such that
$\lim_{s\to 0+}\hat\kappa^\hash_{j,0}(T,B,  s,\bfr,\tau)\leq c'.$ Moreover, 
$\Kc_0(T, {r_1\over\lambda},{r_2\over\lambda})\approx \sum_{j=0}^\upm d_{j,0}\kappa_j(T,{r_1\over\lambda},{r_2\over\lambda}).$
Combining all these estimates, the result follows.
\endproof

\proof[Proof of Proposition \ref{P:nu_tot-monotone}]
We   prove the following  assertion   by increasing induction on $0\leq q\leq k-l:$
\\{\it  There are  a  family  $\Dc_q=\{ d_{j'q'}\in\R^+:\ 0\leq j'\leq k-p- q' \text{ and}\quad q'\leq q\}$ and   a constant $c_q>1 $ depending  on $r_1$ and $r_2$  such that for every positive closed   current $T$ on $\bfU$ belonging to the class $\widetilde\CL^{2,2}_p(\bfU,\bfW),$ 
the  following inequality  hold for any family $\Dc$  which contains $\Dc_q$ and  for $0\leq q\leq \upm:$
\begin{equation}\label{e:monotone-closed-ind}
  \nu^\Dc_q\big(T,B,{r_2\over \lambda},\tau\big)- \nu^\Dc_q(T,B,{r_1\over \lambda},\tau\big) \geq - {c_q\over \lambda}
  +c^{-1}_q\Kc_q\big(T,  {r_1\over \lambda},{r_2\over \lambda}      \big)\quad\text{for}\quad \lambda\gg 1.
 \end{equation}}

Since $\Kc_q\big(T,  {r_1\over \lambda},{r_2\over \lambda}      \big)\geq 0,$  inequality \eqref{e:monotone-closed-ind}
implies the desired conclusion of the theorem. 

By  Lemma \ref{L:nu_0}, inequality \eqref{e:monotone-closed-ind} holds for $q=0.$
Suppose inequality \eqref{e:monotone-closed-ind} true for all $0\leq q<q_0$ with a given $0<q_0\leq k-l.$ We need to prove it for $q=q_0.$ More precisely,  we  need to find  the constants  $d_{j,q_0}>0$ for $0\leq j\leq q_0$ 
such that $\Dc_{q_0}:=\Dc_{q_0-1}\bigcup\{ d_{j,q_0}:\ 0\leq j\leq k-p-q_0\}$  satisfies inequality \eqref{e:monotone-closed-ind}
for $q=q_0.$ Write 
\begin{equation}\label{e:ind_nu-Dc_q_0}
 \begin{split}
  \nu^\Dc_{q_0}\big(T,B,{r_2\over \lambda},\tau\big)- \nu^\Dc_{q_0}(T,B,{r_1\over \lambda},\tau\big) &=\big(\nu^\Dc_{q_0-1}\big(T,B,{r_2\over \lambda},\tau\big)- \nu^\Dc_{q_0-1}(T,B,{r_1\over \lambda},\tau\big)\big)\\
  &+ \sum_{j=0}^{k-p-q_0}  d_{j,q_0} \big(\nu_{j,q_0}\big(T,B,{r_2\over \lambda},\tau\big)- \nu_{j,q_0}(T,B,{r_1\over \lambda},\tau\big)\big).
  \end{split}
\end{equation}
 By the  inductive  hypothesis  we get that
 \begin{equation*}
  \nu^\Dc_{q_0-1}\big(T,B,{r_2\over \lambda},\tau\big)- \nu^\Dc_{q_0-1}(T,B,{r_1\over \lambda},\tau\big) \geq - {c_{q_0-1}\over \lambda}
  +c^{-1}_{q_0-1}\Kc_{q_0-1}\big(T,  {r_1\over \lambda},{r_2\over \lambda}      \big)\quad\text{for}\quad \lambda\geq 1.
 \end{equation*}
 Let $\mu$ be the constant $c$ given  by  Lemma \ref{L:basic-difference-nu_j,q_estimate}.
 Applying Lemma  \ref{L:mu} to $q=k-p-q_0$ yields  the constants $1=\lambda_0 >\lambda_1>\ldots>\lambda_{q_0}>0.$
 Let $0<\theta\ll 1$ be a   small  enough number whose exact value will be determined later. 
 Choose  $d_{j,q}:= \theta\mu^{-1}c^{-1}_{q_0-1} \lambda_j.$ Applying  Lemma \ref{L:basic-difference-nu_j,q_estimate} 
there is  a constant $c>0$ such that  for every positive closed current $T\in  \widetilde\CL^{1,1}_p(\bfU,\bfW)$
and  $0\leq q\leq k-l$ and $0\leq j\leq \min(\upm,k-p-q),$ the following inequality holds:
 \begin{multline*}
\sum_{j=0}^{k-p-q_0} d_{j,q_0} \big(\nu_{j,q}\big(T,B,{r_2\over \lambda},\tau\big)  -\nu_{j,q}\big(T,B,{r_1\over \lambda},\tau\big)\big)\\
\geq  \theta\mu^{-1}c^{-1}_{q_0-1}\Big[ \sum_{j=0}^{k-p-q_0}   \lambda_j \Kc_{j,q}\big (T,{r_1\over\lambda},{r_2\over \lambda}  \big)
-2\mu\lambda_j \sqrt{\Kc_{q}(T,{r_1\over\lambda},{r_2\over \lambda})} \sqrt{ \Kc^-_{j,q}(T,{r_1\over\lambda},{r_2\over \lambda})  }\Big]\\
-   \theta c^{-1}_{q_0-1}(\sum_{j=0}^{k-p-q_0}\lambda_j) \big[\lambda^{-1}+\lambda^{1\over 2} \Kc_{q}\big (T,{r_1\over\lambda},{r_2\over \lambda} \big )+\Kc_{q-1} (T,{r_1\over\lambda},{r_2\over \lambda})\big]  
.
\end{multline*}
Note that
 $1<\sum_{j=0}^{k-p-q_0}\lambda_j<k-p-q_0+1\leq k+1.$ Inserting this  into equality \eqref{e:ind_nu-Dc_q_0} and using the  above inductive hypothesis, we obtain that
\begin{equation*} 
 \begin{split}
  &\nu^\Dc_{q_0}\big(T,B,{r_2\over \lambda},\tau\big)- \nu^\Dc_{q_0}(T,B,{r_1\over \lambda},\tau\big) \\
  &\geq \theta \mu^{-1}c^{-1}_{q_0-1}\Big[ \sum_{j=0}^{k-p-q_0}   \lambda_j \Kc_{j,q}\big (T,{r_1\over\lambda},{r_2\over \lambda}  \big)
-2\mu\lambda_j \sqrt{\Kc_{q}(T,{r_1\over\lambda},{r_2\over \lambda})} \sqrt{ \Kc^-_{j,q}(T,{r_1\over\lambda},{r_2\over \lambda})  }\Big]\\
&-   [c_{q_0-1}+\theta(k+1)c^{-1}_{q_0-1}] \lambda^{-1}-[\theta(q_0+1) c^{-1}_{q_0-1}]\lambda^{1\over 2} \Kc_{q_0}\big (T,{r_1\over\lambda},{r_2\over \lambda} \big ) \\
&+ [ (1    -  \theta(k+1) )c^{-1}_{q_0-1} ] \Kc_{q_0-1}\big (T,{r_1\over\lambda},{r_2\over \lambda} \big )
.
  \end{split}
\end{equation*}
Recall  that $\Kc^-_{j,q}=\Kc_{q-1}+\sum_{j'=0}^{j-1}\Kc_{j',q}.$
Applying  the elementary inequalities  for $a,b\geq 0:$
$$\sqrt{a+b}\leq \sqrt{a}+\sqrt{b}\quad\text{and}\quad 2\sqrt{ab}\leq  {\mu_0\over 2(k+1)}a +   { 2(k+1)\over \mu_0}b
$$
firstly to $a:=  \Kc_{q-1}$ and $b= \sum_{j'=0}^{j-1}\Kc_{j',q},$ and hence secondly to
$a:=\Kc_q$ and $b=\Kc_{q-1},$ we  infer that
\begin{equation*} 
 \begin{split}
  &\nu^\Dc_{q_0}\big(T,B,{r_2\over \lambda},\tau\big)- \nu^\Dc_{q_0}(T,B,{r_1\over \lambda},\tau\big) \\
  &\geq \theta \mu^{-1}c^{-1}_{q_0-1}\Big[ \sum_{j=0}^{k-p-q_0}   \lambda_j \Kc_{j,q}\big (T,{r_1\over\lambda},{r_2\over \lambda}  \big)
-2\mu\lambda_j \sqrt{\Kc_{q}(T,{r_1\over\lambda},{r_2\over \lambda})} \sqrt{\sum_{j'=0}^{j-1} \Kc_{j',q}(T,{r_1\over\lambda},{r_2\over \lambda})  }\Big]\\
&-   [c_{q_0-1}+\theta(k+1)c^{-1}_{q_0-1}] \lambda^{-1}-[\theta(k+1) c^{-1}_{q_0-1}\lambda^{1\over 2} + \theta \mu^{-1}c^{-1}_{q_0-1} {\mu_0\over 2}] \Kc_{q_0}\big (T,{r_1\over\lambda},{r_2\over \lambda} \big ) \\
&+ [ (1    -  \theta(k+1) -(k+1)^2\theta \mu_0^{-1}\mu^{-1})c^{-1}_{q_0-1} ] \Kc_{q_0-1}\big (T,{r_1\over\lambda},{r_2\over \lambda} \big )
.
  \end{split}
\end{equation*}
Define  $t_1,\ldots,t_{k-p-q_0}\geq  0$ as follows:
$$t^2_1:= { \Kc_{0,q}\big (T,{r_1\over\lambda},{r_2\over \lambda}  \big)\over  \Kc_{q}\big (T,{r_1\over\lambda},{r_2\over \lambda}  \big)},\qquad  (t_1+\ldots+t_j)^2 :={ \sum_{j'=0}^{j-1}\Kc_{j',q}\big (T,{r_1\over\lambda},{r_2\over \lambda}  \big)\over  \Kc_{q}\big (T,{r_1\over\lambda},{r_2\over \lambda}  \big)}\quad\text{for}\quad 1\leq j\leq q_0.$$
Recall that  
$  \Kc_{q_0}=\Kc_{q_0-1}+\sum_{j=0}^{k-p-q_0} \Kc_{j,q_0} .$ So  
$$1-(t_1+\ldots t_{k-p-q_0})^2={\Kc_{k-p-q_0,q_0}\over \Kc_{q_0}}+{\Kc_{q_0-1}\over \Kc_{q_0}}.$$
Using the quadratic polynomial $P_{\lambda,\mu}$ introduced in Lemma \ref{L:mu} and noting that $\lambda_{k-p-q_0}\leq 1,$
 we  may rewrite the above  inequality as
\begin{equation*} 
 \begin{split}
&\nu^\Dc_{q_0}\big(T,B,{r_2\over \lambda},\tau\big)- \nu^\Dc_{q_0}(T,B,{r_2\over \lambda},\tau\big) 
  \geq \theta \mu^{-1}c^{-1}_{q_0-1}\Kc_{q_0}(T,{r_1\over\lambda},{r_2\over \lambda})  P_{\lambda,\mu}(t)\\
 &  -   [c_{q_0-1}+\theta(k+1)c^{-1}_{q_0-1}] \lambda^{-1}-[\theta(k+1) c^{-1}_{q_0-1}\lambda^{1\over 2} + \theta \mu^{-1}c^{-1}_{q_0-1} {\mu_0\over 2}] \Kc_{q_0}\big (T,{r_1\over\lambda},{r_2\over \lambda} \big ) \\
&+ [ (1    -  \theta(k+1) -(k^2+2k+2)\theta \mu_0^{-1}\mu^{-1})c^{-1}_{q_0-1} ] \Kc_{q_0-1}\big (T,{r_1\over\lambda},{r_2\over \lambda} \big )
.
  \end{split}
\end{equation*}
Observe that there is $\theta_0>0$ small enough  such that   the coefficient of $\Kc_{q_0-1}$ is $> {1\over 2}c_{q_0-1}^{-1}$
for  $0\in[0,\theta_0).$ Moreover, there is  $\lambda_0\geq 1$  large enough such that
$(k+1) \lambda^{1\over 2} \leq  \mu^{-1} {\mu_0\over 4}$ for $\lambda\geq  \lambda_0.$
By Lemma  \ref{L:mu}, the expression on the RHS is bounded from below by 
$$
 \big[{\mu^{-1}c^{-1}_{q_0-1}\mu_0 \over 4}\big]    \Kc_{q}(T,{r_1\over\lambda},{r_2\over \lambda}) -   (c_{q_0-1}+\theta(k+1)c^{-1}_{q_0-1}) \lambda^{-1}.
$$
Choosing  $c_{q_0}$ such that $c_{q_0}>\max\big(c_{q_0-1}+\theta(k+1) c^{-1}_{q_0-1}, 4\mu c_{q_0-1}\mu^{-1}_0\big),$
we see that inequality  \eqref{e:monotone-closed-ind} holds for $q=q_0$ and $\lambda\gg 1.$

Since $\theta>0$ can be chosen  arbitrarily small,  we  can  choose  $d_{k-p-q_0,q_0} $ so that   $d_{k-p-q_0,q_0} <\epsilon^{q_0} d_{k-p-q_0+1,q_0-1}.$ 
\endproof

Here is the main result of this section.

\begin{theorem}\label{T:nu_tot-monotone} 
Let  $0< r_1<r_2\leq\bfr.$ Then there are  a  family  $\Dc=\{ d_{jq}\in\R^+_*:\ 0\leq j\leq k-p-q,\ 0\leq q\leq  k-l\}$ of positive  numbers and   a constant $c>0 $ depending  on $r_1$ and $r_2$  such that for every positive closed   current $T$ on $\bfU$ belonging to the class $\widetilde\CL^{2,2}_p(\bfU,\bfW),$ 
the  following inequality  hold for $0\leq q\leq \upm:$
\begin{equation}\label{e:monotone-closed}
  \nu^\Dc_q\big(T,B,{r_1\over \lambda},\tau\big)\leq  \nu^\Dc_q(T,B,{r_2\over \lambda},\tau\big) + {c\over \lambda}
  \quad\text{for}\quad \lambda\gg 1.
 \end{equation}
\begin{equation} \label{e:nu-DC-good}  
\nu^\Dc_\tot(T,B, r,\tau)\leq  c\Mc^\tot(T,r)\quad\text{and}\quad
c^{-1}\Mc^\tot(T,r)\leq \nu^\Dc_\tot(T,B, r,\tau)+ cr  
\quad\text{for}\quad 0<r\leq\bfr.
\end{equation}
\end{theorem}
\proof  Choose $\epsilon:=k^{-1}c_1^{-1}.$    Then  applying   Proposition  \ref{P:nu_tot-monotone}
we can choose $\Dc$ such that \eqref{e:monotone-closed} holds and that $d_{k-p-q,q} <\epsilon^q d_{k-p-q+1,q-1}$ for $1\leq q\leq k-l.$ It remains  to show   \eqref{e:nu-DC-good}.  The  first  inequality of  \eqref{e:nu-DC-good} is  easy.
So  we need to prove   the second  inequality of  \eqref{e:nu-DC-good}.

To this end, we  find constants $\mu_j>0$ for $\lowm \leq j\leq \upm$ independent of $T$ and $0<r\leq\bfr$     such that
\begin{equation}\label{e:d_vs_mu_j}
 \sum_{q=0}^{k-l}  d_{k-p-q,q}  \nu_{k-p-q,q} (T,B,r,\tau) =\sum_{j=\lowm}^\upm   \mu_j  \hat\nu_j(T,B,r,\tau)  .
\end{equation}
Indeed, by equality \eqref{e:hat-nu-vs-nu} 
\begin{equation*}
  \hat\nu_j(T,r)=\sum\limits_{q=0}^{\upm-j}  {k-p-j\choose q} c^q_1\nu_{j+q}(T,B,r,\tau).  
 \end{equation*}
 We insert this  equality into the RHS of  \eqref{e:d_vs_mu_j}
 Recall  that $ \nu_{k-p-q,q} (T,B,r,\tau)=\nu_{k-p-q}(T,B,r,\tau)$
on the LHS of   \eqref{e:d_vs_mu_j}. So   by equating   the coefficients   of $ \nu_{k-p-q}(T,B,r,\tau)$ in both sides of \eqref{e:d_vs_mu_j}
using  becomes the system of equations
\begin{equation}
 d_{k-p-q,q}  =\sum\limits_{j=\lowm}^{\min(\upm,k-p-q)} {k-p-j \choose k-p-q-j } \mu_j c^{k-p-q-j}_1.  
 \end{equation}
 We obtain a  triangular system  which  permits us  to calculate the $\mu_j$'s in terms of the $d_{k-p-q,q}$  uniquely.
 The  condition   $0<d_{k-p-q,q} <\epsilon^q d_{k-p-q+1,q-1}$ for $1\leq q\leq k-l$ allows us to show that  $\mu_j>0.$
 
 As in the proof of \eqref{e:Mc_j-induction}, we apply  Lemma   \ref{L:Mc-hash-vs-nu}. 
   So there is a constant $c>0$ such that for $0<r\leq\bfr,$
 \begin{equation*}
 \big|\sum_{j=\lowm}^\upm  \mu_j\Mc^\hash_j(T,r)-\sum_{j=\lowm}^\upm \mu_j\hat\nu_j(T,r)\big|\leq  c r \sum_{j=\lowm}^{\upm}\mu_j\Mc^\hash_j(T,r).
 \end{equation*}
 On the other hand, we infer from  \eqref{e:d_vs_mu_j} and \eqref{e:nu-Dc_tot} and \eqref{e:inter-Lelong-numbers} that there is a constant $c>0$ such that for $0<r\leq\bfr,$
 \begin{equation*} \big|\nu^\Dc_\tot(T,B,r,\tau)-\sum_{j=\lowm}^\upm \mu_j\hat\nu_j(T,r)\big|\leq  c r \sum_{j=\lowm}^{\upm}\mu_j\Mc^\hash_j(T,r).
 \end{equation*}
Therefore, we infer from Lemma   \ref{L:Mc-hash-vs-nu} that
$$
 1-cr\leq {\nu^\Dc_\tot(T,B,r,\tau)\over \sum_{j=\lowm}^\upm  \mu_j\Mc^\hash_j(T,r)  }\leq  1+cr.
$$  
Since  there is a constant $c$ such that $ c \sum_{j=\lowm}^\upm  \mu_j\Mc^\hash_j(T,r)\geq \Mc^\tot(T,r)$  for $0<r\leq\bfr,$
  the second inequality of \eqref{e:nu-DC-good} follows.
  \endproof

  \begin{corollary}\label{C:nu_tot-monotone}
   Let  $0< r_1<r_2\leq\bfr.$ Then there are  a  family  $\Dc=\{ d_{jq}\in\R^+_*:\ 0\leq j\leq k-p-q,\ 0\leq q\leq  k-l\}$ of positive  real numbers and   a constant $c>0 $ depending  on $r_1$ and $r_2$  such that for every positive closed   current $T$ on $\bfU$ belonging to the class $\widetilde\CL^{2,2}_p(\bfU,\bfW),$ 
the  following inequality  hold for $\lowm\leq j\leq \upm:$
\begin{equation*}
  \sum_{m=j}^\upm d_{m,k-p-m}\nu_j\big(T,B,{r_1\over \lambda},\tau\big)\leq  \sum_{m=j}^\upm d_{m,k-p-m} \nu_m(T,B,{r_2\over \lambda},\tau\big) + {c\over \lambda}
  \quad\text{for}\quad \lambda\gg 1.
 \end{equation*}
  \end{corollary}

\section{Positive plurisubharmonic currents  and holomorphic admissible maps}\label{S:psh-holomorphic}

In  this  section we deal with positive  plurisubharmonic  currents together with holomorphic admissible maps, and   we  prove  Theorem  \ref{T:top-Lelong-psh}  and  then  Theorem \ref{T:Lelong-psh-all-degrees}.
 This section  may be regarded   as  a preparation  for  the  proof of Theorems \ref{T:Lelong-psh}, where    the  general situation  with non-holomorphic  admissible maps  will be  investigated.

We  keep  the  global  setting  of Subsection  \ref{SS:Global-setting}
and  suppose   in addition that $T$ is a positive plurisubharmonic  on $X,$ $\tau$ is a holomorphic  admissible map,
and    $\omega$ is a  K\"ahler form on $V.$

\subsection{Top Lelong  number}

\proof[Proof of assertion (1) of Theorem \ref{T:top-Lelong-psh}]
Let $T^\pm_n$ be a  sequence of  approximating forms for $ T^\pm$ as an element of   $\SH^2_p(\overline B).$  
We may assume  that $T^\pm_n$ are in $\SH^2_p(\bfU,\bfW).$
Let $0<r_1<r_2\leq \bfr.$  Theorem \ref{T:Lelong-Jensen-smooth} applied  to  $T^\pm_n\wedge \pi^*(\omega^\upm)$  gives
\begin{multline*}
  \nu_\top(T^\pm_n,B,r_2,\tau)- \nu_\top(T^\pm_n,B,r_1,\tau)=  \int_{\Tube(B,r_1,r_2)}\tau_*T^\pm_n\wedge \pi^*(\omega^\upm)\wedge\alpha^{k-p-\upm}\\
  +\Vc(\tau_*T^\pm_n\wedge \pi^*(\omega^\upm),r_1,r_2)\\
 +  \int_{r_1}^{r_2} \big( {1\over t^{2(k-p-\upm)}}-{1\over r_2^{2(k-p-\upm)}}  \big)2tdt\int_{\Tube(B,t)} \ddc ( \tau_*T^\pm_n\wedge \pi^*\omega^\upm)\wedge \beta^{(k-p-\upm)-1} \\
 +  \big( {1\over r_1^{2(k-p-\upm)}}-{1\over r_2^{2(k-p-\upm)}}  \big) \int_{0}^{r_1}2tdt\int_{z\in \Tube(B,t)} \ddc (\tau_*T^\pm_n\wedge \pi^*\omega^\upm)\wedge \beta^{(k-p-\upm)-1}.
\end{multline*}
Since  $\tau$ is  holomorphic and      $\omega$ is a K\"ahler form on $B,$ 
it follows that
$$
 \ddc ( \tau_*T^\pm_n\wedge \pi^*\omega^\upm)=\tau_*\ddc T^\pm_n\wedge \pi^*\omega^\upm.
$$
Consider  a small neighborhood $V(x_0)$ of  an arbitrary  point $x_0\in  \Tube(B, r_0),$  where in a local chart $V(x_0)\simeq \D^l$ and  $\E|_{V(x_0)}\simeq \C^{k-l}\times \D^l.$
For $x\in \E|_{V(x_0)},$ write $x=(z,w).$  Since   $\upm=\min(l,k-p)$ and $T$ is  of bidegree $(p,p)$ and $\tau$ is  holomorphic,  we  see by the Fact  in the proof of Corollary \ref{C:Lelong-Jensen} that  $S:=\tau_*T^\pm_n\wedge \pi^*\omega^\upm$ and $\ddc S$ are full of bidegree $(l,l)$  in $dw,$ $d\bar w.$
Consequently,  we infer from \eqref{e:hat-alpha'} that
\begin{equation}\begin{split}\label{e:full-in-dw}
\tau_*T^\pm_n\wedge \pi^*(\omega^\upm)\wedge\alpha^{k-p-\upm}
&=\tau_*T^\pm_n\wedge \pi^*(\omega^\upm)\wedge(\hat\alpha')^{k-p-\upm},\\
\tau_*\ddc T^\pm_n\wedge \pi^*\omega^\upm\wedge\beta^{k-p-\upm-1}&=\tau_*\ddc T^\pm_n\wedge \pi^*\omega^\upm\wedge\hat\beta^{k-p-\upm-1}.
\end{split}
\end{equation}
Therefore, as in the proof of 
  Theorem \ref{T:top-Lelong-closed} (1), we deduce from  \eqref{e:Lelong-corona-numbers} that
\begin{equation*}
  \kappa_\top(T^\pm_n,B,r_1,r_2,\tau)= \int_{\Tube(B,r_1,r_2)}\tau_*T^\pm_n\wedge \pi^*(\omega^\upm)\wedge(\hat\alpha')^{k-p-\upm}. 
\end{equation*}
Moreover, by \eqref{e:Lelong-numbers} we also get that
\begin{equation*}
 \nu_\top(T^\pm_n,B,r,\tau)=   
 {1\over r^{2(k-p-\upm)}}\int_{\Tube(B,r)} (\tau_*T^\pm_n)\wedge \pi^*(\omega^\upm)
 \wedge \hat\beta^{k-p-\upm}.  
\end{equation*} 
Consider  the functions
\begin{eqnarray*}
 f^\pm_n(t)&:=&\int_{\Tube(B,t)}  \tau_*(\ddc T^\pm_n)\wedge (\pi^*\omega^\upm)\wedge \beta^{(k-p-\upm)-1},\\ 
 f^\pm(t)&:=&\int_{\Tube(B,t)}  \tau_*(\ddc T^\pm)\wedge (\pi^*\omega^\upm)\wedge \beta^{(k-p-\upm)-1}.
\end{eqnarray*}
By \eqref{e:full-in-dw}, we  get that 
\begin{eqnarray*}
 f^\pm_n(t)&:=&\int_{\Tube(B,t)}  \tau_*(\ddc T^\pm_n)\wedge (\pi^*\omega^\upm)\wedge \hat\beta^{(k-p-\upm)-1},\\ 
 f^\pm(t)&:=&\int_{\Tube(B,t)}  \tau_*(\ddc T^\pm)\wedge (\pi^*\omega^\upm)\wedge \hat\beta^{(k-p-\upm)-1}.
\end{eqnarray*}
So  $f^\pm_n$ and $f^\pm$  are  nonnegative-valued   functions on $(0,\bfr].$
Since  $T^\pm_n$  converge to $T^\pm$  weakly, we  infer that $f^\pm_n$ converge pointwise to $f^\pm$ on  $(0,\bfr]$  except for a countable  set (see  
\eqref{e;cut-off}--\eqref{e:continuity-cut-off}--\eqref{e:except-countable}).
We deduce  from \eqref{e:vertical-boundary-term-bis} and  the fact that   $T^\pm_n\wedge \pi^*\omega^\upm$ is  of full bidegree $(l,l)$  in $dw,$ $d\bar w$  that $\Vc( \tau_*T^\pm_n\wedge \pi^*(\omega^\upm),r_1,r_2)=0.$ 
 Combining  the  above  equalities, we get that 
\begin{equation}\label{-vs-nu}
\begin{split}
&\nu_\top(T^\pm_n,B,r_2,\tau)- \nu_\top(T^\pm_n,B,r_1,\tau)=  \kappa_\top(T^\pm_n,B,r_1,r_2,\tau)\\
  &+  \int_{r_1}^{r_2} \big( {1\over t^{2(k-p-\upm)}}-{1\over r_2^{2(k-p-\upm)}}  \big)2tf^\pm_n(t)dt +  \big( {1\over r_1^{2(k-p-\upm)}}-{1 \over r_2^{2(k-p-\upm)}}  \big) \int_{0}^{r_1} 2tf^\pm_n(t)dt. 
\end{split}
 \end{equation}
Observe that the non-negative functions $f^\pm_n(t),$ $f^\pm(t)$  are  increasing  in $t\in(0,\bfr].$  Moreover, since  $T^\pm_n$ are in $\SH^2_p(\bfU,\bfW)$ and 
$f^\pm(\bfr)<\infty$  and  $f^\pm_n(t)\to f^\pm(t)$ as $n\to\infty$ for all  $t\in(0,\bfr)$ except  for a countable set of values, we may find for every $\bfr'\in (0,\bfr),$ a constant $c=c(\bfr')>0$ such  
that  $f^\pm_n(t)\leq c$ for all $n\geq 1$ and  $t\in (0,\bfr').$
Consequently,  as $n$ tends to infinity,  Lebesgue  dominated  convergence yields that
\begin{equation}\label{e:top-Lelong-kappa-vs-nu}
\begin{split}
&\nu_\top(T,B,r_2,\tau)- \nu_\top(T,B,r_1,\tau)=  \kappa_\top(T,B,r_1,r_2,\tau)\\
  &+  \int_{r_1}^{r_2} \big( {1\over t^{2(k-p-\upm)}}-{1\over r_2^{2(k-p-\upm)}}  \big)2tdt\int_{\Tube(B,t)}  \tau_*(\ddc T)\wedge (\pi^*\omega^\upm)\wedge \hat\beta^{(k-p-\upm)-1} \\
 &+  \big( {1\over r_1^{2(k-p-\upm)}}-{1\over r_2^{2(k-p-\upm)}}  \big) \int_{0}^{r_1}2tdt\int_{z\in \Tube(B,t)} 
 \tau_*(\ddc T)\wedge (\pi^*\omega^\upm)\wedge \hat\beta^{(k-p-\upm)-1}. 
\end{split}
 \end{equation}
Since $T$ and $\ddc T$  are   positive  currents and $\omega,$ $\hat\alpha',$  $\hat\beta$ are positive forms
and the map $\tau$   is holomorphic, the second and third terms on the RHS  are $\geq 0.$
Hence, $\nu_\top(T,B,r_2,\tau)- \nu_\top(T,B,r_1,\tau)\geq \kappa_\top(T,B,r_1,r_2,\tau).$
By the same  positivity, we deduce from the above  expression of  $\kappa_\top(T,B,r_1,r_2,\tau)$ and $\nu_\top(T,B,r,\tau)$ that
they  are non-negative.
This  completes the proof of assertion (1).
\endproof 


\proof[Proof of  assertion (2) of Theorem \ref{T:top-Lelong-psh}]
Since  we know  by assertion (1) that the non-negative function $r\mapsto  \nu_\top(T,B,r,\tau)\geq 0$ is increasing for $r\in(0,\bfr],$
assertion (2) follows.
\endproof


\proof[Proof of  assertion (3) of Theorem \ref{T:top-Lelong-psh}]
By \eqref{e:Lelong-log-bullet-numbers}  and  the identity of assertion (1), we have 
\begin{eqnarray*}\
 0\leq\kappa^\bullet_\top(T,B,r,\tau)= \limsup\limits_{s\to0+}   \kappa_\top(T,B,s,r,\tau)&\leq& \nu_\top(T,B,r,\tau)-\liminf_{s\to0+} \nu_\top(T,B,s,\tau)\\
 &=&  \nu_\top(T,B,r,\tau)- \nu_\top(T,B,\tau),
\end{eqnarray*} 
where the last equality holds by assertion (2).   Consequently, we  infer  from  assertion  (2) again that
\begin{equation*}
0\leq \lim\limits_{r\to0+}\kappa^\bullet_\top(T,B,r,\tau)\leq \lim\limits_{r\to0+}\nu_\top(T,B,r,\tau)- \nu_\top(T,B,\tau)=0.
\end{equation*}
The  result follows.
\endproof


\proof[Proof of  assertion (4) of Theorem \ref{T:top-Lelong-psh}]
Applying Theorem  \ref{T:Lelong-Jensen-smooth}  to the  current  $\tau_*T_n\wedge \pi^*(\omega^\upm)$  and combining  the    equalities  before \eqref{e:top-Lelong-kappa-vs-nu}  in the proof of assertion (1), we get that 
\begin{equation*}
\begin{split}
&\nu_\top(T^\pm_n,B,r,\tau)- \lim_{s\to 0}\nu_\top(T^\pm_n,B,s,\tau)=  \kappa_\top(T^\pm_n,B,r,\tau)\\
  &+  \int_{0}^{r} \big( {1\over t^{2(k-p-\upm)}}-{1\over r^{2(k-p-\upm)}}  \big)2tdt\int_{\Tube(B,t)}  \tau_*(\ddc T^\pm_n)\wedge (\pi^*\omega^\upm)\wedge \hat\beta^{(k-p-\upm)-1}. 
\end{split}
 \end{equation*}
Since $T_n$ and $\ddc T_n$  are   positive  currents and $\omega,$ $\hat\alpha',$  $\hat\beta$ are positive forms
and the map $\tau$   is holomorphic, all the  terms on the LHS and   on the RHS  are $\geq 0.$
Hence, $$\nu_\top(T^\pm_n,B,r,\tau)\geq \int_{0}^{r} \big( {1\over t^{2(k-p-\upm)}}-{1\over r^{2(k-p-\upm)}}  \big)2tdt\int_{\Tube(B,t)}  \tau_*(\ddc T^\pm_n)\wedge (\pi^*\omega^\upm)\wedge \hat\beta^{(k-p-\upm)-1}.$$
On the  other hand,  since  $\|T^\pm_n\|_\bfU \to \|T^\pm\|_\bfU <\infty,$ we  see that there is a constant $c$ independent of $n$ and $0<r\leq\bfr$ such that
 \begin{equation}\label{e:ddc-Lelong-zero-holo-T_n}
 \begin{split}&\int_{0}^{r} \big( {1\over t^{2(k-p-\upm)}}-{1\over r^{2(k-p-\upm)}}  \big)2tdt\int_{\Tube(B,t)}  \tau_*(\ddc T^\pm_n)\wedge (\pi^*\omega^\upm)\wedge \hat\beta^{(k-p-\upm)-1}\\
 &\leq  \nu_\top(T^\pm_n,B,r,\tau)\leq  \nu_\top(T^\pm_n,B,\bfr,\tau)< c.
  \end{split}
  \end{equation}
Since  $\big( {1\over t^{2(k-p-\upm)}}-{1\over \bfr^{2(k-p-\upm)}}  \big)2t\geq 0$ and  the non-negative functions $f^\pm_n$ converge pointwise to $f^\pm$ on  $(0,\bfr]$  except for a countable  set , we infer from  Fatou's lemma  that 
\begin{equation}\label{e:ddc-Lelong-zero-holo-T}\int_{0}^{r} \big( {1\over t^{2(k-p-\upm)}}-{1\over r^{2(k-p-\upm)}}  \big)2tdt\int_{\Tube(B,t)}  \tau_*(\ddc T)\wedge (\pi^*\omega^\upm)\wedge \hat\beta^{(k-p-\upm)-1}<c.
 \end{equation}
By  Theorem  \ref{T:Lelong-closed-Kaehler} applied to the positive closed $(p+1,p+1)$-current $\ddc T,$ we have that
\begin{eqnarray*}
\int_{\Tube(B,t)}  \tau_*(\ddc T)\wedge (\pi^*\omega^\upm)\wedge \hat\beta^{(k-p-\upm)-1}&=&\int_{\Tube(B,t)}  \tau_*(\ddc T)\wedge (\pi^*\omega^\upm)\wedge \beta^{(k-p- \upm)-1}\\
&= & t^{2(k-p-\upm-1)}\nu_\top(\ddc T,B,t, \tau)\\
&\geq& t^{2(k-p-\upm-1)}\nu_\top(\ddc T,B, \tau).
\end{eqnarray*}
Inserting this  inequality into the LHS of \eqref{e:ddc-Lelong-zero-holo-T}, we deduce that
\begin{equation*}\big(\int_{0}^{r} \big( {1\over t^{2(k-p-\upm)}}-{1\over r^{2(k-p-\upm)}}  \big)2t^{2(k-p-\upm)-1}dt\big) \cdot \nu_\top(\ddc T,B, \tau)  <c.
 \end{equation*}
Choose $r:=\bfr.$ Since  the last integral is  equal to infinite,  it follows that  $\nu_\top(\ddc T,B, \tau)=0.$

\endproof

\proof[Proof of  assertion (5) of Theorem \ref{T:top-Lelong-psh}] The proof is  divided into two parts.

\noindent {\bf Proof of the interpretation of assertion (5)  in the sense  of 
Definition \ref{D:Lelong-log-numbers(2)}. }

Fix $0<r\leq\bfr$ and let  $0<\epsilon<r.$  Theorem \ref{T:Lelong-Jensen-eps} applied  to  $\tau_*T\wedge \pi^*(\omega^\upm)$  gives
\begin{equation}\label{e:T:top-Lelong-psh-(2)}\begin{split}
 &\qquad {1\over  (r^2+\epsilon^2)^{k-p-\upm}} \int_{\Tube(B,r)} \tau_*T\wedge \pi^*(\omega^\upm)\wedge \beta^{k-p-\upm} 
 = \Vc_\epsilon(\tau_*T\wedge \pi^*(\omega^\upm),r)\\&+   \int_{\Tube(B,r)} \tau_*T\wedge \pi^*(\omega^\upm)\wedge \alpha_\epsilon^{k-p-\upm}\\
 &+  \int_{0}^{r} \big( {1\over (t^2+\epsilon^2)^{k-p-\upm}}-{1\over (r^2+\epsilon^2)^{k-p-\upm}}  \big)2tdt\int_{\Tube(B,t)} \ddc   [ \tau_*T\wedge \pi^*(\omega^\upm) ]\wedge \beta^{k-p-\upm-1}. 
 \end{split}
\end{equation}  
Next, we let $\epsilon$ tend to $0.$  Then the LHS of \eqref{e:T:top-Lelong-psh-(2)}  tends to $\nu_\top(T,B,r,\tau).$  On the other hand,  we deduce  from \eqref{e:vertical-boundary-term-eps} and  the fact that   $\tau_*T\wedge \pi^*\omega^\upm$ is  of full bidegree $(l,l)$  in $dw,$ $d\bar w$  that $\Vc_\epsilon( \tau_*T\wedge \pi^*(\omega^\upm),r)=0.$ Moreover, using  the  functions $f^\pm$  introduced in the proof of assertion (1), the third term on the RHS of \eqref{e:T:top-Lelong-psh-(2)}  is  rewritten as
\begin{equation}\label{e:T:top-Lelong-psh-(3)}
 \int_{0}^{r} \big( {1\over (t^2+\epsilon^2)^{k-p-\upm}}-{1\over (r^2+\epsilon^2)^{k-p-\upm}}  \big)2t(f^+(t)-f^-(t))dt.
\end{equation}
Observe  that  for  $t\in (0,r],$  we have  as $\epsilon\searrow 0,$  
\begin{multline*}
0\leq  {1\over (t^2+\epsilon^2)^{k-p-\upm}}-{1\over (r^2+\epsilon^2)^{k-p-\upm}}\approx {(r^2-t^2)\over  (t^2+\epsilon^2)^{k-p-\upm}(r^2+\epsilon^2) } \nearrow {(r^2-t^2)\over  t^{2(k-p-\upm)}r^2 } \\
\approx  {1\over t^{2(k-p-\upm)}}-{1\over r^{2(k-p-\upm)}}.
\end{multline*}
So  using  that $f^\pm(t)\geq 0$ almost everywhere on $[0,r],$ an application of  Lebesgue's Monotone Convergence Theorem gives that the expression in \eqref{e:T:top-Lelong-psh-(3)} converges,  as $\epsilon\searrow 0,$  to  
\begin{equation}\label{e:T:top-Lelong-psh-(4)}
 \int_{0}^{r} \big( {1\over t^{2(k-p-\upm)}}-{1\over r^{2(k-p-\upm)}}  \big)2t(f^+(t)-f^-(t))dt.
\end{equation}
By \eqref{e:ddc-Lelong-zero-holo-T}, there is a constant $c>0$ independent of $T$ and $0<r\leq \bfr$  such that
\begin{equation}\label{e:T:top-Lelong-psh-(5)}
 \int_{0}^{r} \big( {1\over t^{2(k-p-\upm)}}-{1\over r^{2(k-p-\upm)}}  \big)2tf^\pm(t)dt\leq c.
\end{equation}
On the  other hand, 
\begin{equation}\label{e:T:top-Lelong-psh-(6)}
 \int_{0}^{r} {1\over r^{2(k-p-\upm)}} 2tf^\pm(t)dt= \int_{0}^{r} {1\over r^{2(k-p-\upm)}} 2t^{2(k-p-\upm)-1}\nu_\top(\ddc T^\pm,B,t,\tau)dt\to 0\,\, \text{as}\,\, r\to 0,
\end{equation}
because  $\nu_\top(\ddc T^\pm,B,\tau)=0$ by assertion (4).
This, combined with \eqref{e:T:top-Lelong-psh-(5)}, implies that by increasing  the constant $c,$
$$
 \int_{0}^{r} \big( {1\over t^{2(k-p-\upm)}} \big)2tf^\pm(t)dt\leq c\quad\text{for}\quad 0<r\leq\bfr.
 $$
 Since  $f^\pm(t)\geq 0,$  it follows that
 $$
 \lim_{r\to 0} \int_{0}^{r} \big( {1\over t^{2(k-p-\upm)}} \big)2tf^\pm(t)dt=0.
 $$
 This, coupled with  \eqref{e:T:top-Lelong-psh-(6)}, gives 
\begin{equation}\label{e:T:top-Lelong-psh-(7)}
 \lim_{r\to 0} \int_{0}^{r} \big( {1\over t^{2(k-p-\upm)}}-{1\over r^{2(k-p-\upm)}}  \big)2tf^\pm(t)dt =0
\end{equation}
Consequently,  by assertions (2) and (3),  the integral in  \eqref{e:T:top-Lelong-psh-(4)} is  bounded and it converges to $0$  as $r\to 0+.$ Putting this,   \eqref{e:T:top-Lelong-psh-(2)} and  \eqref{e:T:top-Lelong-psh-(3)} together,  
we obtain the  desired  interpretation  according to  Definition \ref{D:Lelong-log-numbers(2)}.

\noindent {\bf Proof of the interpretation of assertion (5)  in the sense  of 
Definition \ref{D:Lelong-log-numbers(1)}. }

Since $p>0$ and $l<k,$  it follows  from \eqref{e:m}  that $k-p-\upm<k-l.$
Therefore, we are in the position to apply  Theorem  \ref{T:Lelong-Jensen-smooth}  to the case where $q=k-p-\upm<k-l.$ 
Hence, we get that 
\begin{equation*}
\nu_\top( T^\pm_n,B,r,\tau)=\kappa_\top(T^\pm_n,B,r,\tau)+\Vc(\tau_*T^\pm_n\wedge \pi^*(\omega^\upm),r)
+\int_{0}^{r} \big( {1\over t^{2(k-p-\upm)}}-{1\over r^{2(k-p-\upm)}}  \big)2tf^\pm_n(t)dt.
\end{equation*}
Thus,  we  obtain 
\begin{eqnarray*}
\kappa_\top(T,B,r,\tau)&:=& \lim\limits_{n\to\infty} \kappa_\top(T^+_n-T^-_n,B,r,\tau)=\lim\limits_{n\to\infty}\kappa_\top(T^+_n,B,r,\tau)-\lim\limits_{n\to\infty}\kappa_\top(T^-_n,B,r,\tau)\\
&=&\lim\limits_{n\to\infty}\nu_\top(T^+_n,B,r,\tau)-\lim\limits_{n\to\infty}\nu_\top(T^-_n,B,r,\tau)\\
&-&\lim\limits_{n\to\infty}\int_{0}^{r} \big( {1\over t^{2(k-p-\upm)}}-{1\over r^{2(k-p-\upm)}}  \big)2t(f^+_n(t)-f^-_n(t))dt\\
&=&\nu_\top(T,B,r,\tau) -\lim\limits_{n\to\infty}\int_{0}^{r} \big( {1\over t^{2(k-p-\upm)}}-{1\over r^{2(k-p-\upm)}}  \big)2t(f^+_n(t)-f^-_n(t))dt.
\end{eqnarray*}
So the   interpretation  according to  Definition \ref{D:Lelong-log-numbers(1)} will hold if one can show that
\begin{equation}\label{e:T:top-Lelong-psh-(8)}
 \lim\limits_{n\to\infty}\int_{0}^{r} \big( {1\over t^{2(k-p-\upm)}}-{1\over r^{2(k-p-\upm)}}  \big)2t(f^+_n(t)-f^-_n(t))dt  \to 0\quad\text{as}\quad r\to 0.
\end{equation}
 Recall  from the proof of assertion (1) that   $f^\pm_n(t)\to f^\pm(t)$ as $n\to\infty$ for all  $t\in(0,\bfr)$ except  for a countable set of values, and that for every $\bfr'\in (0,\bfr),$ there is a constant $c=c(\bfr')>0$ such  
that  $0\leq f^\pm_n(t)\leq c$ for all $n\geq 1$ and  $t\in (0,\bfr').$
Consequently, 
\begin{equation*}
 \lim\limits_{n\to\infty}\int_{0}^{r} \big( {1\over t^{2(k-p-\upm)}}-{1\over r^{2(k-p-\upm)}}  \big)2t(f^+_n(t)-f^-_n(t))dt   =\int_{0}^{r} \big( {1\over t^{2(k-p-\upm)}}-{1\over r^{2(k-p-\upm)}}  \big)2t(f^+(t)-f^-(t))dt .
\end{equation*}
So the  desired estimate  \eqref{e:T:top-Lelong-psh-(8)} follows  immediately from inequality \eqref{e:T:top-Lelong-psh-(7)}.
\endproof

\proof[Proof of  assertion (6) of Theorem \ref{T:top-Lelong-psh}]
We argue as  in the  proof of assertion (5)  of  Theorem \ref{T:Lelong-closed-Kaehler}.
The present situation is even simpler  since $\tilde\tau'_\ell=\tau'\circ \tau_\ell^{-1}$ is  holomorphic admissible.
We leave the detals of the proof to the  interested  reader.
\endproof
 
\subsection{Other Lelong  numbers}
This  subsection is devoted to the proof of Theorem  \ref{T:Lelong-psh-all-degrees}.

\proof[Proof of  assertion (1) of Theorem \ref{T:Lelong-psh-all-degrees}]
Let $T^\pm_n$ be a  sequence of  approximating forms for $ T^\pm$ as an element of   $\SH^2_p(\overline B).$  
We may assume  that $T^\pm_n$ are in $\SH^2_p(\bfU,\bfW).$
Let $0<r_1<r_2\leq \bfr$  and $\lowm\leq j\leq \upm.$  Theorem \ref{T:Lelong-Jensen-smooth} applied  to  $T^\pm_n\wedge \pi^*\omega^j$  gives
\begin{multline*}
  \nu_j(T^\pm_n,B,r_2,\tau)- \nu_j(T^\pm_n,B,r_1,\tau)=  \int_{\Tube(B,r_1,r_2)}\tau_*T^\pm_n\wedge \pi^*(\omega^j)\wedge\alpha^{k-p-j}\\
  +\Vc(\tau_*T^j_n\wedge \pi^*(\omega^j),r_1,r_2)\\
 +  \int_{r_1}^{r_2} \big( {1\over t^{2(k-p-j)}}-{1\over r_2^{2(k-p-j)}}  \big)2tdt\int_{\Tube(B,t)} \ddc ( \tau_*T^\pm_n\wedge \pi^*\omega^j)\wedge \beta^{(k-p-j)-1} \\
 +  \big( {1\over r_1^{2(k-p-j)}}-{1\over r_2^{2(k-p-j)}}  \big) \int_{0}^{r_1}2tdt\int_{z\in \Tube(B,t)} \ddc (\tau_*T^\pm_n\wedge \pi^*\omega^j)\wedge \beta^{(k-p-j)-1}.
\end{multline*}
Since  $\tau$ is  holomorphic, and    $T,$  $\ddc T$  are   positive  currents, and   $\omega$ is a K\"ahler form on $B,$ and $\alpha,\ \beta$ are  positive closed form, 
it follows that $\tau_*T^\pm_n\wedge \pi^*(\omega^j)\wedge\alpha^{k-p-j}$ and 
$
\tau_*\ddc T^\pm_n\wedge \pi^*\omega^j\wedge\beta^{k-p-j-1}  $
are  positive currents.
Consider  the functions
\begin{eqnarray*}
 f^\pm_n(t)&:=&\int_{\Tube(B,t)}  \tau_*(\ddc T^\pm_n)\wedge (\pi^*\omega^j)\wedge \beta^{(k-p-j)-1},\\ 
 f^\pm(t)&:=&\int_{\Tube(B,t)}  \tau_*(\ddc T^\pm)\wedge (\pi^*\omega^j)\wedge \beta^{(k-p-j)-1}.
\end{eqnarray*}
So  $f^\pm_n$ and $f^\pm$  are  nonnegative-valued   functions on $(0,\bfr].$
Since  $T^\pm_n$  converge to $T^\pm$  weakly, we  infer that $f^\pm_n$ converge pointwise to $f^\pm$ on  $(0,\bfr]$  except for a countable  set (see  
\eqref{e;cut-off}--\eqref{e:continuity-cut-off}--\eqref{e:except-countable}).
By Theorem  \ref{T:vertical-boundary-terms}  we have    that $\Vc( \tau_*T^\pm_n\wedge \pi^*(\omega^\upm),r_1,r_2)=O(r_2).$ 
 Combining  the  above  equalities, we get that 
\begin{equation}
\begin{split}
&\nu_j(T^\pm_n,B,r_2,\tau)- \nu_j(T^\pm_n,B,r_1,\tau)=  O(r_2)+\kappa_j(T^\pm_n,B,r_1,r_2,\tau)\\
  &+  \int_{r_1}^{r_2} \big( {1\over t^{2(k-p-j)}}-{1\over r_2^{2(k-p-j)}}  \big)2tf^\pm_n(t)dt +  \big( {1\over r_1^{2(k-p-j)}}-{1 \over r_2^{2(k-p-j)}}  \big) \int_{0}^{r_1} 2tf^\pm_n(t)dt. 
\end{split}
 \end{equation}
Observe that the non-negative functions $f^\pm_n(t),$ $f^\pm(t)$  are  increasing  in $t\in(0,\bfr].$  Moreover, since  $T^\pm_n$ are in $\SH^2_p(\bfU,\bfW)$ and 
$f^\pm(\bfr)<\infty$  and  $f^\pm_n(t)\to f^\pm(t)$ as $n\to\infty$ for all  $t\in(0,\bfr)$ except  for a countable set of values, we may find for every $\bfr'\in (0,\bfr),$ a constant $c=c(\bfr')>0$ such  
that  $f^\pm_n(t)\leq c$ for all $n\geq 1$ and  $t\in (0,\bfr').$
Consequently,  as $n$ tends to infinity,  Lebesgue  dominated  convergence yields that
\begin{equation}
\begin{split}
&\nu_j(T,B,r_2,\tau)- \nu_j(T,B,r_1,\tau)=  O(r_2)+\kappa_j(T,B,r_1,r_2,\tau)\\
  &+  \int_{r_1}^{r_2} \big( {1\over t^{2(k-p-j)}}-{1\over r_2^{2(k-p-j)}}  \big)2tdt\int_{\Tube(B,t)}  \tau_*(\ddc T)\wedge (\pi^*\omega^j)\wedge\beta^{(k-p-j)-1} \\
 &+  \big( {1\over r_1^{2(k-p-j)}}-{1\over r_2^{2(k-p-j)}}  \big) \int_{0}^{r_1}2tdt\int_{z\in \Tube(B,t)} 
 \tau_*(\ddc T)\wedge (\pi^*\omega^j)\wedge \beta^{(k-p-j)-1}. 
\end{split}
 \end{equation}
By the  above positivity, the two last terms on the RHS  are $\geq 0,$ and hence  $$\nu_\top(T,B,r_2,\tau)- \nu_\top(T,B,r_1,\tau)\geq \kappa_\top(T,B,r_1,r_2,\tau)+O(r_2).$$
By the same  positivity, we deduce from the above  expression of  $\kappa_\top(T,B,r_1,r_2,\tau)$ and $\nu_\top(T,B,r,\tau)$ that
they  are non-negative.
This  completes the proof of assertion (1).
\endproof 

The proof  of the remaining  assertions  of  Theorem  \ref{T:Lelong-psh-all-degrees}
follow  along  almost the  same lines as those given in the  proof of Theorem   \ref{T:Lelong-psh}.
We only need to use   $j$ instead $\upm$ (resp.   $\alpha,$ $\beta$ instead of  $\hat\alpha',$ $\hat\beta$).

\section{Admissible estimates  for positive plurisubharmonic  currents}
\label{S:Admissible-psh}

In this  section we  develop  admissible  estimates   for positive plurisubharmonic  currents. 
These  estimates are more sophisticated than those  for positive  closed  currents since   the  curvature term comes into play  in the former estimates, whereas  this term  vanishes  automatically  in the latter ones. 

\subsection{Pointwise admissible  estimates}
 We keep the  Extended  Standing  Hypothesis  formulated in Subsection \ref{SS:Ex-Stand-Hyp}. Let $1\leq \ell\leq \ell_0$   and recall   that  $\tilde\tau_\ell (\U_\ell)=\tau(\bfU_\ell).$ 
 
\begin{lemma}\label{L:ddc-difference-function}
Let $1\leq \ell\leq \ell_0$ and   $f$  be  a smooth complex-valued function  defined  on $\tilde\tau_\ell(\U_\ell).$ 
Fix a holomorphic coordinate system $\zeta=(\zeta_1,\ldots, \zeta_k)$ of $\tilde\tau_\ell (\U_\ell).$ 
Set $s=(s_1,\ldots, s_k)=\tilde\tau_\ell.$ Then the  following  two identities hold
 \begin{equation*}\dbar [(\tilde\tau)^*f] -(\tilde\tau)^*[\dbar f]= \sum_{n=1}^k  {\partial f \over \partial\zeta_n}(s) 
 \overline\partial s_n ,
\end{equation*}
\begin{eqnarray*}
 -i\pi\big(\ddc [(\tilde\tau)^*f] -(\tilde\tau)^*[\ddc f]\big)&=& \sum_{m,n=1}^k  {\partial^2 f \over \partial \zeta_m\partial\zeta_n}(s)\partial  s_m\wedge 
 \overline\partial s_n+ \sum_{m,n=1}^k  {\partial^2 f\over \partial \overline\zeta_m\partial\overline\zeta_n}(s)
 \partial \overline s_m\wedge \overline \partial \overline s_n  \\
 &-&\sum_{m,n=1}^k  {\partial^2 f\over \partial \zeta_m\partial\overline\zeta_n}(s)(\partial s_m\wedge \partial \bar s_n +\dbar s_m\wedge \dbar \bar s_n)\\
 &+& \sum_{m=1}^k {\partial f\over \partial \zeta_m}(s)\ddbar s_m 
 +\sum_{m=1}^k {\partial f\over \partial\overline \zeta_m}(s)\ddbar\overline s_m .
 \end{eqnarray*}
\end{lemma}
\proof To prove the first identity,  observe that
$$(\tilde\tau_\ell)^*[\dbar f] = \sum_{n=1}^k  {\partial f \over \partial\bar \zeta_n}(s) \overline\partial  s_n.
$$
On the  other  hand,
$$
d [(\tilde\tau_\ell)^*f] = \sum_{n=1}^k  {\partial f \over \partial\zeta_n}(s) ds_n+  \sum_{n=1}^k  {\partial f \over \partial\bar\zeta_n}(s) d\bar s_n.
$$
Inserting the identities   $d s_m=\partial s_m+\dbar s_m$ and $d\overline s_n=\partial\overline s_n+\dbar \overline s_n$  to the  last  equality and  combining it with  the first  equality, the  first identity of the lemma follows.

 Recall that $\ddc={i\over\pi}\ddbar$ and we have 
\begin{equation*} 
\begin{split}
 \ddbar( f\circ \tilde\tau_\ell) &=\sum_{m,n=1}^k  {\partial^2 f \over \partial \zeta_m\partial\zeta_n}(s)\partial  s_m\wedge 
 \overline\partial s_n+ \sum_{m,n=1}^k  {\partial^2 f\over \partial \overline\zeta_m\partial\overline\zeta_n}(s)
 \partial \overline s_m\wedge \overline \partial \overline s_n  \\
 &+\sum_{m,n=1}^k  {\partial^2 f\over \partial\overline \zeta_m\partial\zeta_n}(s)\partial \overline s_m\wedge
 \overline\partial s_n +\sum_{m,n=1}^k  {\partial^2 f\over \partial \zeta_m\partial\overline\zeta_n}(s)\partial s_m\wedge
 \overline\partial \overline s_n\\
 &+ \sum_{m=1}^k {\partial f\over \partial \zeta_m}(s)\ddbar s_m 
 +\sum_{m=1}^k {\partial f\over \partial\overline \zeta_m}(s)\ddbar\overline s_m .
 \end{split}
 \end{equation*}
 On the other hand, we have that
 $$
 (\tilde\tau_\ell)^*(\ddbar f)=\sum_{m,n=1}^k  {\partial^2 f\over \partial \zeta_m\partial\overline\zeta_n}(s)d s_m\wedge
 d \overline s_n.
 $$
 Inserting the identities   $d s_m=\partial s_m+\dbar s_m$ and $d\overline s_n=\partial\overline s_n+\dbar \overline s_n$  to the  last  equality and  combining it with  the previous  one, the  second identity of the  lemma follows.
 \endproof
\begin{lemma}\label{L:ddc-difference-form}
Let $1\leq \ell\leq \ell_0$ and   $S$  be  a  $(p,q)$-smooth form on $\U_\ell.$
So in a  holomorphic coordinate system $y=(y_1,\ldots, y_k)$ of $\U_\ell,$ $S$ can be rewritten 
as
 $$ 
S:=\sum_{I,J\subset\{1,\ldots,k\}:\ |I|=p, |J|=q} S_{I,J} dy_I\wedge d\bar y_J,
$$
where $S_{I,J}$  are smooth  functions.
 Then the  following  two identities hold:
 \begin{eqnarray*}
 \dbar [(\tilde\tau_\ell)^*S] -(\tilde\tau_\ell)^*(\dbar S)&=&\sum_{I,J}\big(\dbar [(\tilde\tau_\ell)^* S_{I,J}]-(\tilde\tau_\ell)^*[ \dbar S_{I,J}]\big) \bigwedge _{\iota\in I}d[(\tilde\tau_\ell)^*y_\iota]\wedge \bigwedge_{j\in J} d[(\tilde\tau_\ell)^* \bar y_j]\\
 &+&i\pi\sum_{I,J}(\tilde\tau_\ell)^* (S_{I,J}) \wedge \bigwedge_{j\in J} d[(\tilde\tau_\ell)^* \bar y_j]\wedge\big( \sum_{\iota\in I}\pm   \ddc[(\tilde\tau_\ell)^*y_\iota]\wedge
 \bigwedge _{\iota'\in I\setminus \{\iota\}}d[(\tilde\tau_\ell)^*y_{\iota'}]
 \big)\\
 &+& i\pi\sum_{I,J}(\tilde\tau_\ell)^* (S_{I,J}) \wedge \bigwedge _{\iota\in I}d[(\tilde\tau_\ell)^*y_\iota]\wedge\big(\sum_{j\in J}  \ddc[(\tilde\tau_\ell)^* \bar y_j]\wedge \bigwedge_{j'\in J\setminus \{j\}} d[(\tilde\tau_\ell)^* \bar y_j]\big),
 \end{eqnarray*}
\begin{eqnarray*}
 \ddc [(\tilde\tau_\ell)^*S] -(\tilde\tau_\ell)^*(\ddc S)&=&\sum_{I,J}\big(\ddc [(\tilde\tau_\ell)^* S_{I,J}]-(\tilde\tau_\ell)^*[ \ddc S_{I,J}]\big) \bigwedge _{\iota\in I}d[(\tilde\tau_\ell)^*y_\iota]\wedge \bigwedge_{j\in J} d[(\tilde\tau_\ell)^* \bar y_j]\\
 &+&i\pi\sum_{I,J}(\tilde\tau_\ell)^* (dS_{I,J}) \wedge \bigwedge_{j\in J} d[(\tilde\tau_\ell)^* \bar y_j]\wedge\big( \sum_{\iota\in I}\pm   \ddc[(\tilde\tau_\ell)^*y_\iota]\wedge
 \bigwedge _{\iota'\in I\setminus \{\iota\}}d[(\tilde\tau_\ell)^*y_{\iota'}]
 \big)\\
 &+& i\pi\sum_{I,J}(\tilde\tau_\ell)^* (dS_{I,J}) \wedge \bigwedge _{\iota\in I}d[(\tilde\tau_\ell)^*y_\iota]\wedge\big(\sum_{j\in J}  \ddc[(\tilde\tau_\ell)^* \bar y_j]\wedge \bigwedge_{j'\in J\setminus \{j\}} d[(\tilde\tau_\ell)^* \bar y_j]\big).
 \end{eqnarray*}
\end{lemma}
\proof Since   $\ddc={-i\over \pi} \dbar d,$ it follows that
\begin{eqnarray*}
 \dbar [(\tilde\tau_\ell)^*S] &=&\sum_{I,J}\dbar [(\tilde\tau_\ell)^* S_{I,J}] \bigwedge _{\iota\in I}d[(\tilde\tau_\ell)^*y_\iota]\wedge \bigwedge_{j\in J} d[(\tilde\tau_\ell)^* \bar y_j]\\
 &+&i\pi\sum_{I,J}(\tilde\tau_\ell)^* (S_{I,J}) \wedge \bigwedge_{j\in J} d[(\tilde\tau_\ell)^* \bar y_j]\wedge\big( \sum_{\iota\in I}\pm   \ddc[(\tilde\tau_\ell)^*y_\iota]\wedge
 \bigwedge _{\iota'\in I\setminus \{\iota\}}d[(\tilde\tau_\ell)^*y_{\iota'}]
 \big)\\
 &+& i\pi\sum_{I,J}(\tilde\tau_\ell)^* (S_{I,J}) \wedge \bigwedge _{\iota\in I}d[(\tilde\tau_\ell)^*y_\iota]\wedge\big(\sum_{j\in J}  \ddc[(\tilde\tau_\ell)^* \bar y_j]\wedge \bigwedge_{j'\in J\setminus \{j\}} d[(\tilde\tau_\ell)^* \bar y_j]\big),
 \end{eqnarray*}
This, combined with the identity
\begin{equation*}
 (\tilde\tau_\ell)^*(\dbar S)=\sum_{I,J}(\tilde\tau_\ell)^*[ \dbar S_{I,J}] \bigwedge _{\iota\in I}d[(\tilde\tau_\ell)^*y_\iota]\wedge \bigwedge_{j\in J} d[(\tilde\tau_\ell)^* \bar y_j],
\end{equation*}
gives the first identity of the lemma.

To prove the  second identity  observe that
$$
d\big((\tilde \tau_\ell)^*S\big)=(\tilde\tau_\ell)^*( dS)=\sum_{I,J\subset\{1,\ldots,k\}:\ |I|=p,\ |J|=q}d((\tilde\tau_\ell)^* S_{I,J}) \bigwedge _{\iota\in I}d[(\tilde\tau_\ell)^*y_\iota]\wedge \bigwedge_{j\in J} d[(\tilde\tau_\ell)^* \bar y_j].
$$
Using   $\ddc={i\over \pi} \partial d,$ it follows that
\begin{eqnarray*}
 \ddc[ (\tau_\ell)^*S)]&=&\sum_{I,J}\ddc[(\tilde\tau_\ell)^* S_{I,J}]\wedge  \bigwedge _{\iota\in I}d[(\tilde\tau_\ell)^*y_\iota]\wedge \bigwedge_{j\in J} d[(\tilde\tau_\ell)^* \bar y_j]  \\
 &+&i\pi\sum_{I,J}(\tilde\tau_\ell)^* (dS_{I,J}) \wedge \bigwedge_{j\in J} d[(\tilde\tau_\ell)^* \bar y_j]\wedge\big( \sum_{\iota\in I}\pm   \ddc[(\tilde\tau_\ell)^*y_\iota]\wedge
 \bigwedge _{\iota'\in I\setminus \{\iota\}}d[(\tilde\tau_\ell)^*y_{\iota'}]
 \big)\\
 &+&i\pi \sum_{I,J}(\tilde\tau_\ell)^* (dS_{I,J}) \wedge \bigwedge _{\iota\in I}d[(\tilde\tau_\ell)^*y_\iota]\wedge\big(\sum_{j\in J}  \ddc[(\tilde\tau_\ell)^* \bar y_j]\wedge \bigwedge_{j'\in J\setminus \{j\}} d[(\tilde\tau_\ell)^* \bar y_j]\big).
 \end{eqnarray*}
 This, combined  with the  identity
 $$
 (\tilde\tau_\ell)^*(\ddc S)=\sum_{I,J}(\tilde\tau_\ell)^*(\ddc S_{I,J})\wedge  d[(\tilde\tau_\ell)^*y_I]\wedge d[(\tilde\tau_\ell)^* \bar y_J],
 $$
 implies the second identity of the lemma.
\endproof

\subsection{$m$-negligible test forms and basic  volume  estimate}
 
  Recall  from Subsection \ref{SS:Ex-Stand-Hyp} that  for every $1\leq\ell\leq\ell_0,$  there is a local  coordinate system $y=(z,w)$  on $\U_\ell$  with  $V\cap \U_\ell=\{z=0\}.$
  
 \begin{definition}
  \label{D:negligible}
  \rm  Let $S$ be  a differential  form (resp.  a  current) defined on $\Tube(B,r)\subset \E.$  for some $0<r\leq\bfr.$
   So we can write  in a  local representation of $S$ in coordinates $y=(z,w)\in\C^{k-l}\times\C^l: $
   \begin{equation}\label{e:S_IJKL} S=\sum_{M=(I,J;K,L)} S_Mdz_I\wedge d\bar z_J\wedge dw_K\wedge d\bar w_L,
   \end{equation}
where the $S_M=S_{I,J;K,L}(z,w)$  are  the component  functions  (resp.  component  distributions), and the sum is taken over   $M=( I,J;K,L)$ with $I,J\subset\{1,\ldots,k-l\}$ and $K,L\subset \{1,\ldots,l\}.$

For  $M=( I,J;K,L)$ as above, we also write $dy_M$ instead of $dz_I\wedge d\bar z_J\wedge dw_K\wedge d\bar w_L.$
  
    Given $0\leq  m\leq 2l,$ we say that a  differential form $S$ is  {\it $m$-weakly negligible}  if  in the  above  representation,
  for every $I,J,K,L,$ 
 it holds that  $S_{I,J;K,L}$ is  smooth  outside $V$ and   $S_{I,J;K,L}(z,w)=O(\|z\|^{|K|+|L|-m}).$

   Given $0\leq  m\leq 2l,$ we say that a bounded differential form $S$ is  {\it $m$-negligible}  if  in the  above  representation,
  for every $I,J,K,L,$ 
 it holds that  $S_{I,J;K,L}$ is  smooth  outside $V$ and   $S_{I,J;K,L}(z,w)=O(\|z\|^{\bfe(m,K,L)}),$
 where
 \begin{equation*}
  \bfe(m,K,L):=\max\big(0,|K|+|L|-m  \big)\in\N.
 \end{equation*}
 \end{definition}

 \begin{remark}\label{R:negligible}
  \rm   We keep  the  above  notation and let $0\leq  m\leq 2l-1.$   If $S$ is  $m$-weakly negligible then it is $(m+1)$-weakly negligible.  If $S$ is  $m$-negligible then it is $(m+1)$-negligible.
  If $S$ is  $m$-negligible then it is $m$-weakly negligible, but the converse statement is not true in general.
 \end{remark}

 \begin{definition}\rm
  \label{D:negligible-model}
  For $0 \leq m\leq 2l$  and $1\leq \ell\leq \ell_0,$ consider  the  $(k-p,k-p)$-smooth  form on $\U:$
  \begin{equation}\label{e:R_p,j}
  \begin{split}
 R^\star_{p,m}(y)&:=\sum_{q=0}^{\lfloor{ m\over 2}\rfloor}  (\pi^*\omega^q)(y)\wedge \hat\beta^{k-p-q}(y)+\sum_{q> {m\over 2}}^{k-p} \varphi(y)^{q-{m\over 2}} (\pi^*\omega^q)(y)\wedge \hat\beta^{k-p-q}(y) ,\,\, y\in\U;\\
 R^\dagger_{p,m}(y)&:=\sum_{q=0}^{k-p} \varphi(y)^{q-{m\over 2}} (\pi^*\omega^q)(y)\wedge \hat\beta^{k-p-q}(y) ,\,\, y\in\U.
 \end{split}
\end{equation}
 \end{definition}

Typical  negligible and  weakly  negligible   forms are provided  by the  following 
\begin{lemma}\label{L:ex-negligble} $R^\dagger_{p,m}$ are    $m$-weakly negligible and  $R^\star_{p,m}$ are  $m$-negligible.
\end{lemma}
\proof 
We  only give the proof of the   first assertion since the second one can be done similarly. Let $0\leq q\leq k-p.$
If $0\leq q\leq {m\over 2}$ set
$\Phi=\Phi_1:=\pi^*\theta_\ell \cdot \pi^*(\omega^q)\wedge \beta^{k-p-q}.$ Otherwise, set  $ \Phi=\Phi_2:=\varphi^{q-{m \over 2}}\pi^*\theta_\ell \cdot \pi^*(\omega^{q})\wedge \beta^{k-p-q}.$
By Definition
  \ref{D:negligible}, we only need to  show that $\Phi$  is  $m$-negligible.

We  check   $\Phi_1$ in which case    $q\leq {m\over 2}.$
Write $\Phi_1$ in  the form  \eqref{e:S_IJKL}  $\Phi_1=\sum_{M=(I,J;K,L)} S_Mdz_I\wedge d\bar z_J\wedge dw_K\wedge d\bar w_L.$
 Fix a  multi-index $M$ in the above sum and suppose without loss of generality that $|K|\geq |L|.$ 
 Since     $ \pi^*(\omega^q)$ in $\Phi_1$  gives  $q$ elements for $K$ and  also for $L,$  the other contribution for $K$ and $L$ come
 from  the factor $\beta^{k-p-q}.$ The latter contribution is  $|K|-q+|L|-q.$  Since  $q\leq {m\over 2},$ note that
 $
 |K|-q+|L|-q= |K|+|L|-2q\geq |K|+|L|-m.
 $
 Hence,  $|K|-q+|L|-q\geq  \bfe(m,K,L).$
On the other hand,
 the local expression of $\beta$ given in \eqref{e:alpha-beta-local} shows that
 each  coefficient of $\{dw,d\bar w\}$ in  $\hat\beta$ gives a factor of  order at least $\|z\|\approx \varphi^{1\over 2}.$ 
 Hence,  $\Phi_1$  is   $m$-negligible  according to 
 Definition
  \ref{D:negligible}.
  
We  check   $\Phi_2$ in which case    $q>{m\over 2}.$   Write $\Phi_2$ in  the form  \eqref{e:S_IJKL}  $\Phi_2=\sum_{M=(I,J;K,L)} S_Mdz_I\wedge d\bar z_J\wedge dw_K\wedge d\bar w_L.$
 Note that for every $M,$ we have $\min(|K|,|L|)\geq {m\over 2}$   because of the factor $ \pi^*(\omega^q)$ in $\Phi_2$  and   $q>{m\over 2}.$
Hence, $\bfe(m,K,L)=|K|+|L|-m.$
Using  this  and the local expression of $\beta$ given in \eqref{e:alpha-beta-local}, it can be checked that $\Phi_2$  is   $m$-negligible  according to 
 Definition
  \ref{D:negligible}.

\endproof

\begin{definition}\label{D:substraction}\rm
 Given  a multi-index $M=( I,J;K,L)$ with $I,J\subset\{1,\ldots,k-l\}$ and $K,L\subset \{1,\ldots,l\},$
 its {\it length}  $|M|$ is by  definition  $|M|:= |I|+|J|+|K|+|L|.$

  For two multi-indices of the same length $M=(I,J;K,L)$ and $M'=(I',J';K',L')$ with $| M|=|M'|,$ we define {\it the positive substraction  from  $M$ by $M'$} as  the  following nonnegative integer
  \begin{equation}\label{e:Delta}
   \Delta(M,M'):= \max\big( |I\setminus I'|+|J\setminus J'|+ |K\setminus K'|+|L\setminus L'|, 2(|K| +|L|  -|K'| -|L'|)   \big ).
  \end{equation}
  \end{definition}

 \begin{remark}\label{R:substraction} \rm Note that  $ \Delta(M,M')\geq 0$ and the   positive  substraction is  not symmetric, i.e., in general  $\Delta(M,M')\not= \Delta(M',M).$
 \end{remark}
 
 We collect here the basic properties of   the   positive  substraction.
 \begin{lemma}\label{L:comparison}
 \begin{enumerate}
 
 \item  $\Delta (M,M')=0$ if and only if $M=M'.$
 
 \item The triangle inequality holds:
 $\Delta(M,M'')\leq \Delta(M,M')+\Delta(M',M'').$ 
 \item 
  The following inequality  hold for $ M\not=M',$
 \begin{equation*}
  | K|+ |L|\leq  |K'|+ | L'|  +\Delta(M,M')-1.
 \end{equation*}
 \end{enumerate}
 \end{lemma}
\proof 
If $\Delta(M,M')=0,$ then   by  Definition \ref{D:substraction} $I\subset I',$ $J\subset J',$ $K\subset K',$ $L\subset L',$
and hence $ M\subset M',$  which  implies $M=M'$  because $|M|=|M'|.$
Conversely, if $M=M'$   we  see  by  Definition \ref{D:substraction} that $\Delta(M,M')=0.$  This proves assertion (1).

We make the following observation for an element  $j\in I\setminus I'':$  if $j\in I'$ then  $j\in I'\setminus I'',$ otherwise $j\not\in I'$ and hence
$j\in I\setminus I'.$  So $|I\setminus I'|+|I'\setminus  I''|\geq  |I\setminus I''|.$
Using  this and  similar inequalities for $J,K,L$ and the  equality
$$
(|K| +|L|  -|K'| -|L'|)+(|K'| +|L'|  -|K''| -|L''|)=
|K| +|L|  -|K''| -|L''|,$$
we infer from Definition \ref{D:substraction} that assertion (2) holds.

If  $| K|+ |L|<  |K'|+ | L'|,$ then assertion (3)  holds  because $\Delta(M,M')\geq 0.$
If $| K|+ |L|>  |K'|+ | L'|,$ then by  Definition \ref{D:substraction} 
$$\Delta(M,M')\geq  2(|K| +|L|  -|K'| -|L'|)\geq 1+(|K| +|L|  -|K'| -|L'|),$$
which implies assertion (3).
So to complete the proof of assertion (3), we need to treat the case where  $| K|+ |L|=  |K'|+ | L'|.$
In this  last case, assertion (3) becomes $\Delta(M,M')\geq 1$ for $M\not=M',$ which is  true by assertion (1).

\endproof
 \begin{definition}
  \label{D:negligible-bisbisbis}
  \rm  Let $S$ be   a current defined on $\Tube(B,\bfr)\subset \E.$  
   So we can write  in a  local representation of $S$ in coordinates $y=(z,w): $
   $$S=\sum_{M=(I,J;K,L)} S_Mdz_I\wedge d\bar z_J\wedge dw_K\wedge d\bar w_L,$$
where the $S_M=S_{I,J;K,L}(z,w)$  are  the component  functions, and the sum is taken over   $M=( I,J;K,L)$ with $I,J\subset\{1,\ldots,k-l\}$ and $K,L\subset \{1,\ldots,l\}.$ If moreover, $S$ is a  current of dimension $q$ then
we have $|M|=q.$
  
  
   We say that a current $R$ is  $\star$-negligible  (resp.   $\star$-fine) relative to $S$   if  in the  above  representation,  we have for each $M=(I,J;K,L),$
   \begin{equation*}
    R_M=\sum_{M'}  f_{M,M'} S_{M'},
   \end{equation*}
where  $f_{M,M'}$ is a  smooth functions   with $f_{M,M'}(z,w)=O(\|z\|^{\max(1,\Delta(M,M'))})$  (resp.    $f_{M,M'}(z,w)=O(\|z\|^{\Delta(M,M')})$). 
 \end{definition}

 \begin{lemma}\label{L:local-vs-global-forms}
  there is  a constant $c>0$ such that  for every $\lowm\leq j\leq \upm$ and every positive  $(p,p)$-current  on $\Tube (B,\bfr),$ we have
$$
T\wedge (\ddc \|w\|^2)^j(\ddc \|z\|^2)^{k-p-j}\leq c T\wedge \big(\sum_{q=j}^\upm (\pi^*\omega^q)\wedge (\ddc \beta)^{k-p-q}\big)\quad\text{on}\quad  \Tube (B,\bfr).
$$
 \end{lemma}
\proof
Using the expression of $\beta$ in \eqref{e:alpha-beta-local}, we see that
$$
(\ddc \|w\|^2)^j(\ddc \|z\|^2)^{k-p-j}\leq c \sum_{q=j}^l (\pi^*\omega^q)\wedge (\ddc \beta)^{k-p-q}.
$$
Hence, 
$$
T\wedge (\ddc \|w\|^2)^j(\ddc \|z\|^2)^{k-p-j}\leq c T\wedge \big(\sum_{q=j}^l (\pi^*\omega^q)\wedge (\ddc \beta)^{k-p-q}\big)\quad\text{on}\quad  \Tube (B,\bfr).
$$
Since $T$ is  of bidegree $(p,p),$ a degree consideration show that $T\wedge (\pi^*\omega^q)\wedge (\ddc \beta)^{k-p-q}=0$ for $q>\upm.$
The result follows.
\endproof
 
 \begin{lemma}\label{L:Demailly}
  Let  $S$ be a positive $(p,p)$-current on  $\U_\ell$  for some $1\leq \ell\leq \ell_0$ which has the  representation  according to
  Definition \ref{D:negligible} in coordinates $y=(z,w): $
   $$S=\sum_{M=(I,J;K,L)} S_Mdz_I\wedge d\bar z_J\wedge dw_K\wedge d\bar w_L=\sum_M S_Mdy_M,$$
where the $S_M=S_{I,J;K,L}(z,w)$  are  the component  distributions, and the sum is taken over   $M=( I,J;K,L)$ with $I,J\subset\{1,\ldots,k-l\}$ and $K,L\subset \{1,\ldots,l\}$  such that $|I|+|K|=|J|+|L|=p.$
Let  $M=(I,J;K,L)$ be a  multi-index as  above. Then for every $0<r\leq\bfr,$  the following assertions hold.
\begin{enumerate} \item  We have  $$
r^{-|K|-|L|}|S_M|\leq 2^{k-p}  \sum_{M'=(I',I';K',K')} r^{-2|K'|}  |S_{M'}|,
$$
where the sum on the RHS is  taken over all  $M'$ such that $I\cap J\subset I'\subset  I\cup J$ and $K\cap L\subset K'\subset K\cup L.$ 
Here $|S_{M'}|$ is  the absolute value of the measure $S_{M'}.$
\item There is a constant $c>0$ independent of $r$  such that
$${1\over r^{2(k-p-l)+|K|+|L|} }|\langle S_M dy_M\rangle_{\Tube(B,r)}|\leq  c \sum_{q=l-|K\cup L|}^{\upm} \nu_q(S, B,r,\id).$$

\end{enumerate}
 \end{lemma}
\proof
In order to  obtain  assertion (1),
we apply Proposition \ref{P:Demailly}  to the case where $\lambda_j:=1$ for $1\leq j\leq k-l,$   and $\lambda_j:=r^{-1}$  for $k-l+1\leq j\leq k.$ 

Applying  assertion (1) yields that
\begin{eqnarray*}
{1\over r^{2(k-p-l)+|K|+|L|} }|\langle S_M dy_M\rangle_{\Tube(B,r)}|\leq
2^{k-p}\sum_{M'=(I',I';K',K')}   {1\over r^{2(k-p-l+|K'|)} }|\langle S_{M'} dy_{M'}\rangle_{\Tube(B,r)}|.
\end{eqnarray*}
Consider  $M'=(I',I';K',K')$ and  set $j:=l-|K'|.$
By Lemma \ref{L:local-vs-global-forms}
we have that
$$  |S_M'|=
S\wedge (\ddc \|w\|^2)^j(\ddc \|z\|^2)^{k-p-j}\leq c S\wedge \big(\sum_{q=j}^\upm (\pi^*\omega^q)\wedge (\ddc \beta)^{k-p-q}\big)\quad\text{on}\quad  \Tube (B,\bfr).
$$
Consequently,  we  get that  
\begin{equation*}
 {1\over r^{2(k-p-l+|K'|)} }|\langle S_{M'} dy_{M'}\rangle_{\Tube(B,r)}|\leq  c\sum_{q=l-|K'|}^\upm\nu_{q} (S,B,r,\id),
\end{equation*}
and assertion (2)  follows.
\endproof
 
 \begin{proposition}\label{P:mass-fine-negligible}
  Let $T$ be  a  positive current   and   $\Phi$  a   real continuous  form of dimension $2p$  on $\Tube(B,\bfr).$
  Assume that  $R$  is  a  current on  $\Tube(B,\bfr)$ such that one of the following  conditions  is  satisfied:
  \begin{itemize}
  \item [(i)]  $R$  is   $\star$-fine  relative to $T$ and $\Phi$  is   $m$-negligible;
  \item[(ii)]  $R$  is   $\star$-negligible  relative to $T$ and $\Phi$  is   $(m+1)$-negligible;
  \end{itemize}
   Then there is  a constant $c=c_\Phi>0$ such that   for every $0<r\leq \bfr,$
\begin{equation*}
  {1\over r^{2(k-p)-m}} |\langle R,\Phi\rangle_{\Tube(B,r)}| \leq    c \sum_{q=\lowm}^{\upm} \nu_q(T,B, r,\id).
\end{equation*}
 \end{proposition}
  \proof We  divide the proof into two parts.
  
  \noindent {\bf Proof of Case (i):}   Since $R$ is    $\star$-fine relative to $T,$  by  Definition  \ref{D:negligible} we have, for each $M=(I,J;K,L),$
  the following representation
   \begin{equation}\label{e:R-sum-R_M}
    R_M=\sum_{M'}  f_{M,M'} T_{M'},
   \end{equation}
where  $f_{M,M'}$ is a  smooth functions   with    $f_{M,M'}(z,w)=O(\|z\|^{\Delta(M,M')}).$
Observe that
\begin{equation}\label{e:R-leq-sum-R_M}
  {1\over r^{2(k-p)-m}} |\langle R,\Phi\rangle_{\Tube(B,r)}|\leq  \sum_M {1\over r^{2(k-p)-m}} |\langle R_M dy_M,\Phi\rangle_{\Tube(B,r)}|.
\end{equation}

 \noindent{\bf  Subcase (i-1):}  $m\leq  2l-|K|-|L|.$
  
  Observe  that
  $
  \langle R_M dy_M,\Phi\rangle=\langle R_M dy_M,\Phi'\rangle,
  $
  where $\Phi'$ is the component of bidegree $M^\bfc=(I^\bfc,J^\bfc;K^\bfc,L^\bfc)$  of $\Phi.$
  Since $\Phi$ is $m$-negligible,  we deduce from Definition \ref{D:negligible} that $\Phi'(y)=O(\|z\|^{|K^\bfc|+|L^\bfc|-m})=  O(\|z\|^{2l-|K|-|L|-m}).$  Hence,   
    the RHS of \eqref{e:R-leq-sum-R_M} in this  subcase is  dominated by a constant times
 $$
   \sum_{M=(I,J;K,L)} {1\over r^{2(k-p)-m}} |\langle \|z\|^{2l-|K|-|L|-m}R_M dy_M\rangle_{\Tube(B,r)}|.
  $$  
  In order   to    majorize each term  in the sum on the RHS, fix a multi-index $M=(I,J;K,L).$  
  Since $2l-|K|-|L|-m\geq 0$ and $\|z\|\leq  r$ for $y=(z,w)\in\Tube(B,r),$
  it follows  that  each  term in the above sum is  majorized by
 $ {1\over r^{2(k-p-l)+|K|+|L|}} |\langle   R_M dy_M\rangle_{\Tube(B,r)}| .$
   By  \eqref{e:R-sum-R_M}, this expression is  dominated by 
   \begin{equation*}
    \sum_{M'=(I',J';K',L'):\ |I'|+|K'|=|J'|+|L'|=k-p}   {1\over r^{2(k-p-l)+|K|+|L|}} |\langle f_{M,M'} T_{M'}dy_M\rangle_{\Tube(B,r)}|.
  \end{equation*}
   By Definition \ref{D:substraction}
   we have $
  \max(0,| K|+ |L|-|K'|- | L'|)\leq    \Delta(M,M').$  Hence,  $f_{M,M'}(z,w)=O(\|z\|^{ \max(0,| K|+ |L|-|K'|- | L'|)}).$
  Therefore,   we infer that the term in the sum  of the last line is bounded from above by a constant times
    \begin{equation*}
      {1\over r^{2(k-p-l)+|K'|+|L'|}} |\langle  T_{M'}dy_{M'}\rangle_{\Tube(B,r)}|.
   \end{equation*}
    This integral  is, in turn,     bounded from above by a constant times
   $
   \sum_{q=\lowm}^{\upm} \nu_q(T,B, r,\id) 
   $ by Lemma \ref{L:Demailly} (2).  Hence, we obtain the desired conclusion in this  subcase.

 \noindent{\bf  Subcase (i-2):}  $m> 2l-|K|-|L|.$ 
 
  As in the previous subcase, observe  that
  $
  \langle R_M dy_M,\Phi\rangle=\langle R_M dy_M,\Phi'\rangle,
  $
  where $\Phi'$ is the component of bidegree $M^\bfc$  of $\Phi.$
  Since $\Phi$ is $m$-negligible and  $m>|K^\bfc|+|L^\bfc|,$ we deduce from Definition \ref{D:negligible} that  $\Phi'(y)=O(1).$  Hence,   
    the RHS of \eqref{e:R-leq-sum-R_M} in this  subcase is  dominated by a constant times
   $$
   \sum_{M=(I,J;K,L)} {1\over r^{2(k-p)-m}} |\langle R_M dy_M\rangle_{\Tube(B,r)}|.
  $$  
  In order   to    majorize each term  in the sum on the RHS, fix a multi-index $M=(I,J;K,L).$  
 We infer from   the assumption $m> 2l-|K|-|L|$ that the above term is   dominated by
 $ {1\over r^{2(k-p-l)+|K|+|L|}} |\langle   R_M dy_M\rangle_{\Tube(B,r)}| .$
 By  \eqref{e:R-sum-R_M}, this expression is  dominated by 
   \begin{equation*}
    \sum_{M'=(I',J';K',L'):\ |I'|+|K'|=|J'|+|L'|=k-p}   {1\over r^{2(k-p-l)+|K|+|L|}} |\langle f_{M,M'} T_{M'}dy_M\rangle_{\Tube(B,r)}|.
  \end{equation*}
    We conclude  the  proof of this  subcase  as in Subcase (i-1).

  \noindent {\bf Proof of Case (ii):}   Since $R$ is    $\star$-negligible relative to $T,$  by  Definition  \ref{D:negligible} we have, for each $M=(I,J;K,L),$
  the  representation
\eqref{e:R-sum-R_M},
where  $f_{M,M'}$ is a  smooth functions   with   
\begin{equation*}  f_{M,M'}(z,w)=\begin{cases} O(\|z\|^{\Delta(M,M')}),&\quad\text{if}\quad M'\not=M;\\
                                   O(\|z\|),&\quad\text{if}\quad M'=M.
                                 \end{cases}
\end{equation*}
On the  other hand, recall  from  Lemma \ref{L:comparison} (3)
that for $ M\not=M',$
we have $
  (| K|+ |L|)- (|K'|+ | L'|)\leq   \Delta(M,M')-1.$
  Using  the above two  inequalities, we argue as in the proof of  Case (i). Hence,  Case (ii) follows.
  \endproof

 \begin{proposition}\label{P:mass-fine-negligible-bis}
  Let $T$ be  a  positive current   and   $\Phi$  a   real continuous  form of dimension $2p$  on $\Tube(B,\bfr).$
  Assume that  $R$  is  a  current on  $\Tube(B,\bfr)$ such that one of the following  conditions  is  satisfied:
  \begin{itemize}
  \item [(i)]  $R$  is   $\star$-fine  relative to $T$ and $\Phi$  is   $m$-weakly negligible;
  \item[(ii)]  $R$  is   $\star$-negligible  relative to $T$ and $\Phi$  is   $(m+1)$-weakly negligible;
  \end{itemize}
   Then there is  a constant $c=c_\Phi>0$ such that   for every $0<r\leq \bfr,$
\begin{equation*}
  {1\over r^{2(k-p)-m}} |\langle R,\Phi\rangle_{\Tube(B,{r\over 2}, r)}| \leq    c\sum_{q=\lowm}^{\upm} \nu_q(T,B, r,\id).
\end{equation*}
 \end{proposition}
  \proof We only give the proof of case (i) since   case (ii) can be done 
    similarly.   
    
     Since $R$ is    $\star$-fine relative to $T,$    we have also , for each $M=(I,J;K,L),$
   representation \eqref{e:R-sum-R_M}.
Observe that instead of \eqref{e:R-leq-sum-R_M} we have
\begin{equation}\label{e:R-leq-sum-R_M-bis}
  {1\over r^{2(k-p)-m}} |\langle R,\Phi\rangle_{\Tube(B,{r\over 2}, r)}|\leq  \sum_M {1\over r^{2(k-p)-m}} |\langle R_M dy_M,\Phi\rangle_{\Tube(B,{r\over 2},r)}|.
\end{equation}
Note that
  $
  \langle R_M dy_M,\Phi\rangle=\langle R_M dy_M,\Phi'\rangle,
  $
  where $\Phi'$ is the component of bidegree $M^\bfc=(I^\bfc,J^\bfc;K^\bfc,L^\bfc)$  of $\Phi.$
  Since $\Phi$ is $m$-negligible,  we deduce from Definition \ref{D:negligible} that $\Phi'(y)=O(\|z\|^{|K^\bfc|+|L^\bfc|-m})=  O(\|z\|^{2l-|K|-|L|-m}).$  Hence,   
    the RHS of \eqref{e:R-leq-sum-R_M-bis}  is  dominated by a constant times
 $$
   \sum_{M=(I,J;K,L)} {1\over r^{2(k-p)-m}} |\langle \|z\|^{2l-|K|-|L|-m}R_M dy_M\rangle_{\Tube(B,{r\over 2},r)}|.
  $$  
  In order   to    majorize each term  in the sum on the RHS, fix a multi-index $M=(I,J;K,L).$  
  Since  ${r\over 2}\leq \|z\|\leq  r$ for $y=(z,w)\in\Tube(B,{r\over  2},r),$
  it follows  that  each  term in the above sum is  majorized by
 $ {1\over r^{2(k-p-l)+|K|+|L|}} |\langle   R_M dy_M\rangle_{\Tube(B,{r\over 2},r)}| .$
   By  \eqref{e:R-sum-R_M}, this expression is  dominated by 
   \begin{equation*}
    \sum_{M'=(I',J';K',L'):\ |I'|+|K'|=|J'|+|L'|=k-p}   {1\over r^{2(k-p-l)+|K|+|L|}} |\langle f_{M,M'} T_{M'}dy_M\rangle_{\Tube(B,{r\over 2},r)}|.
  \end{equation*}
   We  conclude the proof as  in the Subcase (i-1) of the proof of Proposition \ref{P:mass-fine-negligible}.

  \endproof

\subsection{Basic boundary formula}


For  every current $S$ of bidegree $(p,q)$ on $\E,$ we will  always fix a smooth approximating  $(p,q)$-forms  $(S_\epsilon)_{\epsilon>0}$  which can  be obtained  from  $S$  using   a standard  convolution locally and  patching  the local regularizations by   a partition  of unity.
Let $\Omega\Subset \Tube (B,\bfr)$  be an open set.  Suppose that $\|S\|(\partial \Omega)=0,  $ we get
\begin{equation}\label{e:limit_current}
 \lim\limits_{\epsilon\to 0} \int_\Omega S_\epsilon\wedge \phi =\int_\Omega S\wedge \phi\qquad\text{for}\qquad\phi\in\Cc^\infty(\Omega).
\end{equation}


\begin{definition}\label{D:boundary-value} \rm Let $\Sigma$ be  an open set of $\partial \Omega$  which is a $\Cc^1$-real  hypersurface in $\E,$ we define
\begin{equation*}\int_{\Sigma}S:= \lim_{\epsilon\to 0}\int_{\Sigma} S_\epsilon
\end{equation*}
provided that the limit exists and is finite.

In what follows, we  will use  $\Sigma:=\partial_\hor\Tube(B,r)$ for $0<r\leq\bfr.$
\end{definition}

Let $S$ be a  current of order $0.$ Then  the set
\begin{equation}
 \label{e:except-set}
 \Ec_S:=\left\lbrace r\in (0,\bfr]:  \|S\| (\partial_\hor\Tube(B,r))>0  \right\rbrace
\end{equation}
is at most countable.

\begin{proposition}\label{P:Stokes}
Fix  $\ell$  with $1\leq \ell\leq \ell_0$ and $r\in(0,\bfr].$    Set $\tilde\tau:=\tilde\tau_\ell$  and $\H:=\Tube (\widetilde V_\ell,r)\subset \E.$
Then, for every 
 every   current $S$ of bidimension $(q-1,q-1)$  defined on  $\U_\ell$ and every  smooth  form $\Phi$ of bidegree $(q,q)$  defined on $\tilde\tau(\H)$
 with $\pi(\supp(\Phi))\Subset \widetilde V_\ell,$  we have
\begin{eqnarray*}
 &&\langle  \ddc (\tilde\tau_* S) -\tilde\tau_*(\ddc  S),\Phi  \rangle_{\tilde \tau (\H)}
 =\langle  S,\tilde\tau^*(\ddc\Phi) -\ddc(\tilde\tau^*\Phi)\rangle_{\H} \\ &+&\big( \langle    S ,\dc (\tilde\tau^*\Phi)^\sharp -\tilde\tau^*(\dc\Phi) \rangle_{\partial \H} - \langle \tilde\tau^*[(\tilde\tau_*S)^\sharp]-S, \tilde\tau^*(\dc\Phi)\rangle_{\partial\H}\big) \\
 &-& {1\over2\pi i}\big( \langle  \tilde\tau^*[ (\tilde\tau_* S)^\sharp]  ,\tilde\tau^*(d\Phi) -d[(\tilde\tau^*\Phi)^\sharp] \rangle_{\partial \H} -\langle   S -\tilde\tau^*[(\tilde\tau_*S)^\sharp)],d[(\tilde\tau^*\Phi)^\sharp]  \rangle_{\partial \H}\big)\\
 &-&{1\over \pi i} \big(  \langle  \dbar (\tilde\tau_* S)^\sharp  ,\Phi  \rangle_{\partial[ \tilde \tau (\H)]} -\langle  \dbar  S ,(\tilde\tau^*\Phi)^\sharp  \rangle_{\partial \H} \big).
\end{eqnarray*}
Here,  we  have used the operator $\sharp$  introduced in Notation \ref{N:principal}.
\end{proposition}
\proof
Since  the current $S$ is of bidimension $(q-1,q-1)$  and the  smooth  form $\Phi$ is of bidegree $(q,q),$ the LHS  is  rewritten  as  follows:
\begin{eqnarray*}
 \langle  \ddc (\tilde\tau_* S) ,\Phi  \rangle_{\tilde \tau (\H)}-\langle \tilde\tau_*(\ddc  S), \Phi  \rangle_{\tilde \tau (\H)}
&=&\langle  \ddc (\tilde\tau_* S)^\sharp ,\Phi  \rangle_{\tilde \tau (\H)}-\langle \ddc  S ,\tilde\tau^*\Phi  \rangle_{ \H}\\
&=&\langle  \ddc (\tilde\tau_* S)^\sharp ,\Phi  \rangle_{\tilde \tau (\H)}-\langle \ddc  S ,(\tilde\tau^*\Phi)^\sharp  \rangle_{ \H}.
\end{eqnarray*}
By Stokes' theorem  (see e.g. \cite[Formula III.3.1]{Demailly}),  the last line  is  equal to
\begin{multline*}
 \big(\langle  (\tilde\tau_* S)^\sharp ,\ddc\Phi  \rangle_{\tilde \tau (\H)}+ \langle  \dc (\tilde\tau_* S)^\sharp ,\Phi  \rangle_{\partial[\tilde \tau (\H)]}
 - \langle  (\tilde\tau_* S)^\sharp ,\dc\Phi  \rangle_{\partial[\tilde \tau (\H)]}  \big)\\
 - \big(\langle  S ,\ddc (\tilde\tau^*\Phi)^\sharp  \rangle_{\H}+ \langle  \dc  S ,(\tilde\tau^*\Phi)^\sharp  \rangle_{\partial \H}
 - \langle    S ,\dc (\tilde\tau^*\Phi)^\sharp  \rangle_{\partial \H}\big).
\end{multline*}
Since we have  by a bidegree consideration
\begin{equation*} \langle  (\tilde\tau_* S)^\sharp ,\ddc\Phi  \rangle_{\tilde \tau (\H)}=\langle  \tilde\tau_* S ,\ddc\Phi  \rangle_{\tilde \tau (\H)}=\langle  S,\tilde\tau^*(\ddc\Phi)\rangle_{\H}\quad\text{and}\quad
\langle  S ,\ddc (\tilde\tau^*\Phi)^\sharp  \rangle_{\H}=\langle  S ,\ddc (\tilde\tau^*\Phi) \rangle_{\H},
\end{equation*}
 it follows that
 \begin{eqnarray*}
 &&\langle  \ddc (\tilde\tau_* S) -\tilde\tau_*(\ddc  S),\Phi  \rangle_{\tilde \tau (\H)}
 =\langle  S,\tilde\tau^*(\ddc\Phi) -\ddc(\tilde\tau^*\Phi)\rangle_{\H} \\ &+&\big(\langle    S ,\dc (\tilde\tau^*\Phi)^\sharp \rangle_{\partial \H} - \langle (\tilde\tau_*S)^\sharp, \dc\Phi\rangle_{\partial|\tilde\tau(\H)]}\big) 
 + \big(\langle  \dc (\tilde\tau_* S)^\sharp  ,\Phi  \rangle_{\partial[ \tilde \tau (\H)]} -\langle  \dc  S ,(\tilde\tau^*\Phi)^\sharp  \rangle_{\partial \H} \big)\\
 &=& I+II+III.
\end{eqnarray*}
 Using that $\tilde\tau$ is  diffeomorphic  from $\partial \H$ to $\partial|\tilde\tau(\H)],$ we have that
 \begin{equation*}
 II=\langle    S ,\dc (\tilde\tau^*\Phi)^\sharp -\tilde\tau^*(\dc\Phi) \rangle_{\partial \H} - \langle \tilde\tau^*[(\tilde\tau_*S)^\sharp]-S, \tilde\tau^*(\dc\Phi)\rangle_{\partial\H}.
 \end{equation*}
 Using  the  identity $\dc ={d\over 2\pi i} -{\dbar\over \pi i},$  we see that
 \begin{equation*}
  III={1\over 2\pi i}\big(\langle  d (\tilde\tau_* S)^\sharp  ,\Phi  \rangle_{\partial[ \tilde \tau (\H)]} -\langle  d S ,(\tilde\tau^*\Phi)^\sharp  \rangle_{\partial \H} \big)-{1\over \pi i}\big(\langle  \dbar (\tilde\tau_* S)^\sharp  ,\Phi  \rangle_{\partial[ \tilde \tau (\H)]} -\langle  \dbar  S ,(\tilde\tau^*\Phi)^\sharp  \rangle_{\partial \H} \big).
 \end{equation*}
By  Stokes' theorem applied to $\partial \H$  and $\partial[ \tilde \tau (\H)] $ and using the diffeomorphism $\tau$  again,
the  first  expression in parentheses is  equal to
\begin{eqnarray*}
 &&-\big(\langle  (\tilde\tau_* S)^\sharp  ,d\Phi  \rangle_{\partial[ \tilde \tau (\H)]} -\langle   S ,d[(\tilde\tau^*\Phi)^\sharp]  \rangle_{\partial \H}\big)\\
 &=&\langle  \tilde\tau^*[ (\tilde\tau_* S)^\sharp]  ,\tilde\tau^*(d\Phi) -d[(\tilde\tau^*\Phi)^\sharp] \rangle_{\partial \H} -\langle   S -\tilde\tau^*[(\tilde\tau_*S)^\sharp)],d[(\tilde\tau^*\Phi)^\sharp]  \rangle_{\partial \H}. 
\end{eqnarray*}
Using  the new  expressions for $II$ and $III,$
 we obtain the  desired formula.
\endproof


\subsection{Boundary differential  operators: First part}
\label{SS:Boundary-diff-oper-I}

Fix  a   smooth  increasing function $\chi:\ \R\to [0,1]$ 
which is  equal to $0$ on $(-\infty,-1]$ and is  equal to $1$ on $[-{1\over 2},\infty).$
For $0<r\leq\bfr$ and $0<\epsilon<r,$ set  $\chi_{r,\epsilon}(t):=     \chi({t-r\over \epsilon})$ for $t\in\R.$

\begin{definition}\label{D:diff-oper}\rm
Fix $1\leq \ell\leq\ell_0$ and $m\in\{0,1\}.$  Let $\Cc^\infty_\comp(\U_\ell)$ be  the space of 
smooth differential forms with compact support in $\U_\ell.$
Consider    the class  $\DO^m_\ell$ of  differential operators $D:\ \Cc^\infty_\comp(\U_\ell)\to  \Cc^\infty_\comp(\U_\ell)$ of  order $m$ on $\U_\ell$   whose  coefficients  are  the product of   the function $\theta_\ell$ and smooth forms
on $\Tube(B,\bfr).$    For a  current  $S$ of a  given degree $n$ on $\U_\ell,$ write
 $$S=\sum_{M=(I,J;K,L)} S_Mdz_I\wedge d\bar z_J\wedge dw_K\wedge d\bar w_L,$$
where the $S_M=S_{I,J;K,L}(z,w)$  are  the component  distributions, and the sum is taken over   $M=( I,J;K,L)$ with $I,J\subset\{1,\ldots,k-l\}$ and $K,L\subset \{1,\ldots,l\}$ such that $|M|=n.$
Here $|M|:=|I|+|J|+|K|+|L|.$

Consider  the   subclass  $\widehat\Dc^0_\ell\subset \DO^0_\ell$ consisting of  all $D\in\DO_\ell$  such that  for a current $S,$  by writing  $R:=DS,$  we have
\begin{equation}\label{e:diff-oper-0}
    R_M=\theta_\ell\sum_{M'}  f_{M,M'} S_{M'},
   \end{equation}
where  $f_{M,M'}$ is a  smooth functions   with     $f_{M,M'}(z,w)=O(\|z\|^{\Delta(M,M')}).$
A differential operator $D\in  \widehat\Dc^0_\ell$ is  said to be  {\it $\star$-fine of order $0$}.
If  moreover $f_{M,M'}(z,w)=O(\|z\|^{\max(1,\Delta(M,M'))})$ for all  $M,M'$  then  we say that $D$ is  a  {\it $\star$-negligible of order $0$}.
The set of all  $\star$-negligible differential  operators  $D$ of order $0$ is  denoted by $\Dc^0_\ell.$  So  $\Dc^0_\ell\subset\widehat \Dc^0_\ell.$

\end{definition}

\begin{definition}\rm
 \label{D:substraction-bis} Given two multi-indices  $M=(I_M,J_M), N=(I_N,J_N)\subset \{1,\ldots,k\}^2$ such that $|N|=|M|-1$  and an integer $j\in \{1,\ldots,k\} ,$  we define {\it the positive substraction by index $j$ from $M$  by $N$} is
 $$
 \Delta_j(M,N):=\min\limits_{P} \big(\delta_{j,P,M}+ \Delta(P,N)\big).
 $$
Here,
\begin{itemize}
 \item   the sum is taken is over  all multi-index $P=(I_P,J_P)\subset \{1,\ldots,k\}^2$ such that $P\subset M$ and $|P|=|M|-1;$
 \item $\Delta(P,N)$ is calculated by Definition \ref{D:substraction};
 \item $\delta_{j,P,M}=1$ if  we have  $j\in\{1,\ldots,k-l\}$
 and $M\setminus P\subset \{k-l+1,\ldots, k\}$   simultaneously. Otherwise, $\delta_{j,P,M}=0.$
\end{itemize}

\end{definition}

\begin{definition}\label{D:diff-oper-bis}\rm
Consider also  the   subclass  $\widehat\Dc^1_\ell$ consisting of  all $D\in\DO^1_\ell$   such that  for a current $S,$  by writing  $R:=DS,$  we have 
\begin{equation}\label{e:diff-oper-1}
    R_M=\theta_\ell\big(\sum_{M',j}  f_{M,M',j} {\partial S_{M'}\over \partial y_j} +  g_{M,M',j} {\partial S_{M'}\over \partial \bar y_j}\big),
   \end{equation}
    the sum being taken over all $M'$ with $|M'|=|M|-1=2k-q-1$ and $1\leq  j\leq k.$
Here  $f_{M,M',j},$ $g_{M,M',j}$ are  smooth functions  such that    $$f_{M,M',j}(z,w)=O(\|z\|^{\Delta_j(M,M')})\quad\text{ and}\quad g_{M,M',j}(z,w)=O(\|z\|^{\Delta_j(M,M')}).$$
A differential operator $D\in  \widehat\Dc^1_\ell$ is  said to be  {\it $\star$-fine of order $1$}.
If moreover  for every    $M,M',j$  we have
  $$f_{M,M',j}(z,w)=O(\|z\|^{\max(\Delta_j(M,M'),1)})\quad\text{  and}\quad g_{M,M',j}(z,w)=O(\|z\|^{\max(1,\Delta_j(M,M'))})$$  then  we say that $D$ is  a  {\it $\star$-negligible of order $1$}.
The set of all  $\star$-negligible  differential  operators  $D$ of order $1$ is  denoted by $\Dc^1_\ell.$  So  $\Dc^1_\ell\subset\widehat \Dc^1_\ell.$
\end{definition}
Let $D^\star$ be the  adjoint  operator  of $D,$ that is,  if $\Phi$ is a  smooth  form   compactly supported   in $\Tube(B,\bfr),$ then
\begin{equation}\label{e:adjoint}
 \langle   DS,\Phi\rangle=  \langle   S,D^\star\Phi\rangle.
\end{equation}

\begin{proposition}\label{P:boundary-vs-tube-eps}
Let $S$ be a  positive plurisubharmonic  current of  bidimension   $(q,q)$ on a neighborhood of $\Tube(B,\bfr)$ such that $S$ and $\ddc S$  such that $S$ is  $\Cc^1$-smooth   near $ \partial_\ver\Tube(B,\bfr).$  Let $\delta\in\{0,1\}.$ 
\begin{itemize}
 \item[(i-0)]  If  $D$ is  a  differential operator in the class $\widehat\Dc^0_\ell$ and   $\Phi$ is a  form of  degree $2q-1$  which is  $m$-negligible,
 then there are:
 \begin{itemize} 
 \item  [$\bullet$] a  bounded  form $S_0$ is   in a neighborhood of 
 $\partial_\ver\Tube(B,r)$  which depends only on $D$ and $S;$
  \item  [$\bullet$] 
  three differential operators $D_1,$ $D_2$ and $D_3$ in   the class $\widehat\Dc^0_\ell;$
  \item[$\bullet$] and three   forms
 $\Phi_1$ of degree $2q$  which is $(m+1)$-negligible,
 $\Phi_2$ of degree $2q$  which is $m$-negligible and 
 $\Phi_3$ of degree $(2q-1)$  which is  $m$-negligible;
 \end{itemize}
 such that  for every $0<r\leq\bfr,$  we have
 \begin{equation}\label{e:P:boundary-vs-tube-eps}
 \begin{split}
 \int_{\partial_\hor \Tube(B,r)} DS\wedge \Phi&= \int_{\partial_\ver\Tube(B,r)} S_0\wedge \Phi + \int_{\Tube(B,r)} D_1S\wedge \Phi_1 \\
 &+  {1\over r} \int_{\Tube(B,r)} D_2S\wedge \Phi_2+\lim\limits_{\epsilon\to 0+}\int_{\Tube(B,r-\epsilon,r)} D_3S\wedge d\chi_{r,\epsilon}\wedge \Phi_3,
 \end{split}
 \end{equation}
 
 \item[(ii-0)]  If  $D$ is  a  differential operator in the class $\Dc^0_\ell$ and   $\Phi$ is a smooth form of  degree $2q-1$  which is  $m$-negligible,  then the conclusion of assertion (i-0) also holds. Moreover, the three differential operators $D_0,$ $D_1$ and $D_2$ belong to    the class $\Dc^0_\ell.$
  
 \end{itemize}
\end{proposition}
\proof 

We only give   the proof of assertion  (i-0). Since  the proof  of assertion  (ii-0) is  similar, it is left  to the  interested  reader.
By \eqref{e:diff-oper-0} we may assume  without loss of generality that $D  S =\theta_\ell  f S_{M'}dy_M,$ where $ M,M'\subset\{1,\ldots,k\}$  are some  multi-indices, and $f$ is a  bounded form on $\Tube(B,\bfr)$ smooth out of $V$ and $f(z,w)=O(\|z\|^{\Delta(M,M')}).$ In what follows for $y\in\Tube(B,r)$ we write $y=(z,w)$
and note that  $\|y\|\approx\|z\|.$
Since $\|y\|=r$ for $y\in\partial_\hor \Tube(B,r),$ it follows that
\begin{equation*}
 \int_{\partial_\hor \Tube(B,r)} DS\wedge \Phi= \int_{\partial_\hor \Tube(B,r)}  {\|y\|^2\over r^2} (DS\wedge \Phi)(y)
 \end{equation*}
 So, by  Stokes' theorem, we have
 \begin{multline*}
 \int_{\partial_\hor \Tube(B,r)} DS\wedge \Phi=
 -\int_{\partial_\ver \Tube(B,r)}{\|y\|^2\over r^2}( DS\wedge \Phi)(y)+\int_{\Tube(B,r)} d( {\|y\|^2\over r^2} (DS\wedge \Phi)(y))\\
 =  -\int_{\partial_\ver \Tube(B,r)} {\|y\|^2\over r^2}(DS\wedge \Phi)(y)+\int_{\Tube(B,r)} d\big( {\|y\|^2\over r^2} \theta_\ell(y)f(y) S_{M'}(y)dy_M\wedge \Phi(y)\big).
 \end{multline*}
 The  first term  on the RHS  is  of the form  is  of the form  $\int_{\partial_\ver\Tube(B,r)}S_0\wedge \Phi_ ,$  where $S_0$ is the restriction of ${\|y\|^2\over r^2} DS$ to  $\partial_\ver\Tube(B,r).$ So $S_0$  is
 a  bounded form in a neighborhood of 
 $\partial_\ver\Tube(B,r)$ which depends only on $D$ and $S.$

The second term on the RHS  can be expanded into the expression
\begin{equation}\label{e:trick-with-y-over-r}
 \pm \int_{\Tube(B,r)}   {\|y\|^2\over r^2}\theta_\ell  f\Phi (d S_{M'}\wedge dy_M) \pm\int_{\Tube(B,r)} d(  {\|y\|^2\over r^2} \theta_\ell  f\Phi)\wedge  S_{M'}dy_M.
\end{equation}
Since 
\begin{equation*}
 1-\chi_{r,\epsilon}(y)=\begin{cases}
                         1,& \text{for},\,\, |y|\leq r-\epsilon;\\
                         0,&\text{for},\,\, |y|\geq r-{\epsilon\over 2};\\
                         \in [0,1],& \text{otherwise},
                        \end{cases}
\end{equation*}
and hence  $\lim_{\epsilon\to 0} 1-\chi_{r,\epsilon}(y)= 1$ for $y\in \Tube(B,r),$  the  first term of   expression \eqref{e:trick-with-y-over-r} can be rewritten as
\begin{multline*}
 \lim\limits_{\epsilon\to 0+} \int_{\Tube(B,r)}  (1-\chi_{r,\epsilon})  {\|y\|^2\over r^2}\theta_\ell f\Phi (d S_{M'}\wedge dy_M)=\lim\limits_{\epsilon\to 0+} \int_{\Tube(B,r)} d[ (1-\chi_{r,\epsilon}) {\|y\|^2\over r^2}\theta_\ell f\Phi] S_{M'}\wedge dy_M\\
 =-\lim\limits_{\epsilon\to 0+} \int_{\Tube(B,r)} d\chi_{r,\epsilon}\wedge  {\|y\|^2\over r^2} \theta_\ell f\Phi\wedge dy_M  \wedge S_{M'}+ 
 \lim\limits_{\epsilon\to 0+} \int_{\Tube(B,r)}  (1-\chi_{r,\epsilon}) d[ {\|y\|^2\over r^2}\theta_\ell f\Phi]\wedge  dy_M\wedge  S_{M'}.
 \end{multline*}
 Arguing as  in the analysis of  the second term of  expression  \eqref{e:trick-with-y-over-r} (see the paragraph below), 
 we see that the second  term of the last line can be written in the form
 \begin{equation*}
 \int_{\Tube(B,r)} D'_1S\wedge \Phi'_1+r^{-1}\int_{\Tube(B,r)} D'_2S\wedge \Phi'_2 .
 \end{equation*}
where $D'_1,D'_2$ are  differential operators   in the class $\widehat\Dc^0_\ell,$ and $\Phi'_1$ is   form of degree $2q$   which is  $(m+1)$-negligible, and $\Phi'_2$ is   form of degree $2q$   which is  $m$-negligible

The second   term  of  expression  \eqref{e:trick-with-y-over-r}
can be  rewritten as
\begin{eqnarray*}
&&\int_{\Tube(B,r)} d(  {\|y\|^2\over r^2}) \theta_\ell  f\Phi\wedge  S_{M'}dy_M\pm \int_{\Tube(B,r)}  {\|y\|^2\over r^2} d(\theta_\ell)  f\Phi\wedge  S_{M'}dy_M\\ &&\pm\int_{\Tube(B,r)}  {\|y\|^2\over r^2} \theta_\ell df\wedge\Phi\wedge  S_{M'}dy_M
\pm \int_{\Tube(B,r)} {\|y\|^2\over r^2} \theta_\ell  f(d\Phi)\wedge  S_{M'}dy_M\\
&=:& I_1+I_2+I_3+I_4.
\end{eqnarray*}
Observe that
$I_1$ is  of the form  $r^{-1}\int_{\Tube(B,r)} D_2S\wedge \Phi_2$ for  a  differential operator $D_0$  
in the class $\widehat\Dc^0_\ell$  and   form $\Phi_2$ of bidegree $2q$   which is  $m$-negligible.
Next, $I_2$ is  of the form  $\int_{\Tube(B,r)} D_1S\wedge \Phi_1$ for  a  differential operator $D_1$  
in the class $\widehat\Dc^0_\ell$  and   form $\Phi_1$ of bidegree $2q$   which is  $(m+1)$-negligible.
Since  $f(z,w)=O(\|z\|^{\Delta(M,M')}),$  it follows that  $\|y\|df(y)=O(\|z\|^{\Delta(M,M')}).$ This, combined with
the  inequality  $\|y\|<  r$ for $y\in \Tube(B,r),$ implies that 
$$
I_3=r^{-1}\int_{\Tube(B,r)}  {\|y\|\over r} \theta_\ell (\|y\|df(y)\wedge\Phi(y)\wedge  S_{M'}(y)dy_M=r^{-1}\int_{\Tube(B,r)} D_2S\wedge \Phi_2
$$
for  a  differential operator $D_2$  
in the class $\widehat\Dc^0_\ell$  and   form $\Phi_2$ of bidegree $2q$   which is  $m$-negligible.

Since  $\Phi$ is $m$-negligible,   we can check using  Definition \ref{D:negligible} that  $\|y\|d\Phi(y)$ is also $m$-negligible. This, combined with
the  inequality  $\|y\|<  r$ for $y\in \Tube(B,r),$ implies that 
$$
I_4=r^{-1}\int_{\Tube(B,r)}  {\|y\|\over r} \theta_\ell(y) f(y)\wedge(\|y\|d\Phi(y))\wedge  S_{M'}(y)dy_M=r^{-1}\int_{\Tube(B,r)} D_2S\wedge \Phi_2
$$
for  a  differential operator $D_2$  
in the class $\widehat\Dc^0_\ell$  and   form $\Phi_2$ of bidegree $2q$   which is  $m$-negligible.

Putting together the above estimates, the result follows.
\endproof

\begin{proposition}\label{P:boundary-vs-tube}
Let $S$ be a  positive   current of  bidimension   $(q,q)$ on a neighborhood of $\Tube(B,\bfr)$ such that $S$ and $\ddc S$  such that $S$ is  $\Cc^1$-smooth   near $ \partial_\ver\Tube(B,\bfr).$
Let $D$ be  a  differential operator  in  the class $\DO^0_\ell.$ 
 Let $\Phi$ be a $\Cc^2$-smooth  form $\Phi$ of degree $2q-1$ on  $\Tube(B,r).$  
 For  $0<s<r\leq\bfr,$  consider
 $$
 I_{s,r}:=\int_{s}^r\limsup\limits_{\epsilon\to 0+} \big|\int_{\Tube(B,t-\epsilon,t)} DS\wedge d\chi_{t,\epsilon}\wedge \Phi\big|dt.
 $$
 Suppose that  one of  the following two cases  happens:
 \begin{enumerate} 
 \item $D$ is   in  the class $\widehat\Dc^0_\ell$    and    $\Phi$   is  $m$-weakly negligible;  

 \item  If  $D$ is   in  the class $\Dc^0_\ell$    and    $\Phi$   is  $(m+1)$-weakly negligible. 
 \end{enumerate}
 Then there is a constant  $c>0$ independent of $s,\ r$ such that
 $$
 |I_{s,r}|\leq  c\int_{\Tube(B,s,r)} S\wedge   R^\dagger_{k-q,m} .
 $$
\end{proposition}
\proof
For  $0<t\leq \bfr,$ set
$$
 J_t:=\limsup\limits_{\epsilon\to 0+}\big|\int_{\Tube(B,t-\epsilon,t)} DS\wedge d\chi_{t,\epsilon}\wedge \Phi\big|.
 $$
So $I_{s,r}=\int_{s}^r J_tdt. $

\noindent  {\bf  Proof of Case  (1):} We consider two subcases.

\noindent {\bf Subcase (i):  $S$ is continuous.}

By \eqref{e:diff-oper-0} we may assume  without loss of generality that $D  S =\theta_\ell  f S_{M'}dy_M,$ where $ M,M'\subset\{1,\ldots,k\}$  are some  multi-indices, and $f$ is a  bounded form on $\Tube(B,\bfr)$ smooth out of $V$ and $f(z,w)=O(\|z\|^{\Delta(M,M')}).$
Since $S$ is a current of bidimension $(q,q),$ it follows that $dy_{M'}$ is also  of bidimension $(q,q).$

 Since  $d\chi_{t,\epsilon}(y)={1\over\epsilon}\chi'({\rho-t\over\epsilon})d\rho,$  where $\rho:=\|y\|,$ and by \eqref{e:varphi-new-exp} we have $\|y\|=\|A(w)z\|$ 
it follows that 
 $$
d\chi_{t,\epsilon}\wedge\Phi={1\over\epsilon}\chi'({\rho-t\over\epsilon}) \Psi,\qquad\text{where}\qquad \Psi(z,w):= d\| A(w) z\|\wedge \Phi(z,w).
 $$
  Recall  that $\Phi$ is  $m$-weakly negligible and $S$ is   continuous.
 Therefore, we infer from   the expressions  of $\Psi$  and of  $J_t$ that
 $$
 J_{s,r}= \int_{\Tube(B,s,r)} \theta_\ell |f| \cdot|S_{M'} dy_M\wedge  \Psi|.
 $$
 Let $\Phi_1$ be  the component of $\Psi$ correspponding  to  $dy_{M^\bfc},$
 where   for $M=(I,J;K,L),$  $M^\bfc$ denotes $(I^\bfc,J^\bfc;K^\bfc,L^\bfc).$
 Since $\Psi$ is $m$-weakly negligible, so is $\Phi_1.$  Write $M_1=(I_1,J_1;K_1,L_1):=M^\bfc.$
   By Definition \ref{D:negligible}, we may assume  without loss of generality that
 $$\Phi_1= g dz_{I_1}\wedge d\bar z_{J_1}\wedge dw_{K_1}\wedge d\bar w_{L_1}=gdy_{M_1},$$
where  $I_1,J_1\subset\{1,\ldots,k-l\}$ and $K_1,L_1\subset \{1,\ldots,l\}$ such that     $g(z,w)=O(\|z\|^{|K_1|+|L_1|-m})$. 
Since $dy_M\wedge \Psi=dy_M\wedge  \Psi_1,$ it follows that
$$
 J_{s,r}= \int_{\Tube(B,s,r)} \theta_\ell |f| \cdot|S_{M'} dy_M\wedge  \Psi_1|=  \int_{\Tube(B,s,r)} \theta_\ell |fg| \cdot|S_{M'} dy_M\wedge  dy_{M^\bfc}|.
 $$
 Next,  we find $M_2=(I_2,J_2;K_2,L_2)$ and $M_3=(I_3,J_3;K_3,L_3):=M_2^\bfc$  such that $|M_2|=|M|$ and  $ |K|+|L|= |K_2|+|L_2|$ and $ |K_3|+|L_3|= |K_1|+|L_1|$  and $dy_{M_2}$ is of bidegree $(q,q),$ that is $|I_2|+|K_2|=|J_2|+|L_2|=q.$ Indeed, it suffices to change some $dz_p$ (resp. $d\bar z_{p'}$) into $d\bar z_{p'}$ (resp. $dz_p$) and to change some $dw_q$ (resp. $d\bar w_{q'}$) into $d\bar w_{q'}$ (resp. $dw_q$). 
 So   $dy_M\wedge dy_{M_1}=\pm dy_{M_2}\wedge  dy_{M_3}.$ Consider  the  $(q,q)$-form  $\Psi_3:= g dy_{M_3}.$  We infer that
 $$
 J_{s,r}= \int_{\Tube(B,s,r)} \theta_\ell |fg| \cdot|S_{M'} dy_{M_2}\wedge  \Psi_3|.
 $$
 We also  deduce from the above  equalities and   Definition  \ref{D:substraction} that
 \begin{eqnarray*} 
 |K_1|+|L_1|-m&=& |K_3|+|L_3|-m,\\
\Delta(M,M')&\geq& \max\big(0,|K|+|L|-|K'|-|L'|\big)=
  \max\big(0,|K_2|+|L_2|-|K'|-|L'|\big).
  \end{eqnarray*}
 This,  combined with the last expression for $J_{s,r},$ implies that
 $$
 J_{s,r}= \int_{\Tube(B,s,r)} \theta_\ell \|z\|^{\max\big(0,|K_2|+|L_2|-|K'|-|L'|\big)}    \|z\|^{|K_3|+|L_3|-m}  \cdot|S_{M'} dy_{M_2}\wedge  dy_{M_3}|.
 $$
 Since  $|K_2|+|K_3|=|L_2|+|L_3|=l,$  it follows that
 $$
 J_{s,r}\leq \int_{\Tube(B,s,r)} \theta_\ell \|z\|^{2l-|K'|-|L'|-m}  \cdot|S_{M'} dy_{M_2}\wedge  dy_{M_3}|.
 $$
 By Lemma \ref{L:Demailly}  applied to the positive current  $S,$ we have that
 $$
\|z\|^{-|K'|-|L'|}|S_{M'}|\leq 2^{k-p}  \sum_{M''=(I'',I'';K'',K'')} \|z\|^{-2|K''|}  |S_{M''}|,
$$
where the sum on the RHS is  taken over all  $M''$ such that $I'\cap J'\subset I''\subset  I'\cup J'$ and $K'\cap L'\subset K''\subset K'\cup L'.$ 
Combining  the last two estimates, we   get that
$$
J_{s,r}\leq  2^{k-p}\int_{\Tube(B,s,r)} \theta_\ell        \sum_{M''=(I'',I'';K'',K'')} \|z\|^{2l-2|K''|-m}  |S_{M''}|                     d \Leb(y),
$$
where $d\Leb(y)$ is  the Lebesgue measure on  $\U_\ell.$ The integrand   on the RHS is  bounded from above by a constant times  $S\wedge  R^\dagger_{k-q,m}.$  Hence, there is a constant $c>0$ such that
 $
 J_{s,r}\leq  c\int_{\Tube(B,s,r)}  S\wedge  R^\dagger_{k-q,m}. 
 $
 This  completes  the proof of Subcase (i).
 
\noindent {\bf Subcase (ii):  $S$ is general.}
We leave it to the interested reader.
 
Case (1) is  thereby completed.

\noindent  {\bf  Proof of Case  (2):}

By \eqref{e:diff-oper-0} we may assume  without loss of generality that $D  S =\theta_\ell  f S_{M'}dy_M,$ where $ M,M'\subset\{1,\ldots,k\}$  are some  multi-indices, and $f$ is a  bounded form on $\Tube(B,\bfr)$ smooth out of $V$ and $f(z,w)=O(\|z\|^{\max(1,\Delta(M,M'))}).$
Since $S$ is a current of bidimension $(q,q),$ it follows that $dy_{M'}$ is also  of bidimension $(q,q).$
On the  other hand, recall  from  Lemma \ref{L:comparison} (3)
that for $ M\not=M',$
we have $
  (| K|+ |L|)- (|K'|+ | L'|)\leq   \Delta(M,M')-1.$
  So  $f(z,w)=O(\|z\|^{\max(1,1+(| K|+ |L|)- (|K'|+ | L'|)) }).$
  Using this, we argue as in the proof of  Case (1). Hence,  Case (2) follows.
\endproof
\endproof


\subsection{Boundary differentiel operators:  Second part}

Fix  $\ell$  with $1\leq \ell \leq\ell_0.$  
Fix  a local coordinate  system $y=(z,w)$ on $\U_\ell$  with  $V\cap \U_\ell=\{z=0\}.$
 Without loss of generality we   may assume   \eqref{e:max-coordinate}, that is,
$ 2|z_{k-l}| > \max\limits_{1\leq j\leq k-l}|z_j|.$ Recall  that  $y=(z,w).$ Write $y'=(z',w)\in\C^{k-l-1}\times\C^l=\C^{k-l-1},$ where $z=(z',z_{k-l}).$

We   introduce a new  coordinate  system $\tilde y=(y', u,  t)=\widetilde Y(y),$  where 
\begin{equation}\label{e:coor-tilde-y}
u=u(z_{k-l}):=|z_{k-l}|\big({\arg (z_{k-l})\over \pi}-1\big)\in[-|z_{k-l}|,|z_{k-l}|)\qquad\text{and}\qquad t:=\sqrt{\varphi (y)}=\|y\|\in[0,\infty),
\end{equation}
where $\arg(z_{k-l})\in[-\pi,\pi]$  is  the argument of $z_{k-l}\in\C^*.$
By \eqref{e:varphi-new-exp} we obtain  that
\begin{equation}\label{e:coor-tilde-y-t}
t= \| A(w)z\|\qquad\text{for}\qquad  z\in\C^{k-l},\ w\in \D^l.
\end{equation}
Using  this and  \eqref{e:coor-tilde-y},  a direct computation shows that
 \begin{equation}
  \label{e:coor-tilde-y-t-der}
  \begin{split}
   {\partial  u(z_{k-l})\over  \partial z_{k-l}} &= O(1),\qquad {\partial  u(z_{k-l})\over  \partial \bar z_{k-l}}=O(1),\\
   {\partial  u(z_{k-l})\over  \partial w} &= {\partial  u(z_{k-l})\over  \partial \bar w}=0\quad\text{and}\quad    {\partial  u(z_{k-l})\over  \partial z_p} =  {\partial  u(z_{k-l})\over  \partial \bar z_p}    =0\quad\text{for}\quad 1\leq p\leq k-l-1,\\
  {\partial  t(z,w)\over  \partial w}&=O(\|z\|)=O(t)\quad\text{and}\quad  {\partial  t(z,w)\over  \partial \bar w}=O(\|z\|)=O(t),\\
  {\partial  t(z,w)\over  \partial z} &=  {\partial  t(z,w)\over  \partial \bar z}=    O(1).
  \end{split}
 \end{equation}
 Let  $$\M:=\left\lbrace \tilde y=(\tilde y_1,\ldots,\tilde y_k)=(z',\tilde y_{k-l}, w) \in\  \D^{k-l-1}\times \D\times \D^l=\D^k:\,\, \tilde y_k=u+it\,\, \text{and}\,\,  \max_{1\leq j\leq k-l-1} |z_j| \leq  2|t| \right\rbrace.  
 $$
 For  $0<r\leq \bfr,$ let
 $$
 \M_r:=    \left\lbrace \tilde y=(\tilde y_1,\ldots,\tilde y_k)=(z',\tilde y_{k-l}, w) \in\  \D^{k-l-1}\times \D\times \D^l=\D^k:\quad  \tilde y_k=u+it\quad \text{and}\quad  t=r \right\rbrace          .$$
 Observe  that   $\M_r\subset (2r\D)^{k-l-1}\times  (r\D)\times  \D^{l}.$
 Write  $\tilde y:=\widetilde Y(y).$
 \begin{lemma}\label{L:widetilde_Y}
 By  using  a refinement  of     the family  $(\U_\ell)_{1\leq \ell\leq \ell_0}$ if necessary, 
  $\widetilde Y$ is a  smooth diffeomorphism on  each  $\U_\ell.$
 \end{lemma}
 \proof
 We only need  to  check that the  Jacobian  of $\widetilde Y$ is nonzero on each  $\U_\ell.$
 But this  follows from \eqref{e:coor-tilde-y-t-der}.
 \endproof
Write $y=Y(\tilde y).$
 By Lemma  \ref{L:widetilde_Y},   $Y$ is a smooth  diffeomorphism from $\U_\ell\setminus V$ onto  $\M$ with the inverse $\widetilde Y.$
As in Definition
  \ref{D:negligible},  we have a similar notion in  the new coordinate system $\tilde y=(z',u+it,w).$ 
 \begin{definition}
  \label{D:negligible-bisbis}
  \rm  Let $S$ be  a continuous differential  form (resp.  a  current) defined on $\M.$  
   So we can write  in a  local representation of $S$ in coordinates $\tilde y=(\tilde y_1,\ldots,\tilde y_k)= (z',u+it,w): $
   $$S=\sum_{M=(I,J;K,L)} S_Mdz_I\wedge d\bar z_J\wedge dw_K\wedge d\bar w_L.$$
Here $S_M=S_{I,J;K,L}(z,w)$  are  the component  functions  (resp.  component  distributions), and the sum is taken over   $M=( I,J;K,L)$ with $I,J\subset\{1,\ldots,k-l\}$ and $K,L\subset \{1,\ldots,l\}$  with  the following  convention 
\begin{equation}\label{e:convention-z}dz_{k-l}:=  d\tilde y_{k-l}=du+idt\qquad\text{and}\qquad d\bar z_{k-l}:=d\overline{\tilde y}_{k-l}=du-idt.
\end{equation}

  Let $0\leq  m\leq 2l.$ We say that $S$ is  {\it $m$-weakly negligible}  if  in the  above  representation,
 if for every $I,J,K,L$ with $|K|+|L|\geq  m,$ then    $S_{I,J;K,L}(\tilde y)=O(t^{|K|+|L|-m}).$ 
 \end{definition}

\begin{definition}\label{D:diff-oper-bisbis}\rm
Consider    the class  $\widetilde\DO^1_\ell$ of  differential operators $D:\ \Cc^\infty_\comp(\M)\to  \Cc^\infty_\comp(\M)$ of  order $1$ on $\M.$  For a  current  $S$ on $\M,$ write
 $$S=\sum_{M=(I,J;K,L)} S_Mdz_I\wedge d\bar z_J\wedge dw_K\wedge d\bar w_L,$$
where the $S_M=S_{I,J;K,L}(z,w)$  are  the component  distributions, and the sum is taken over   $M=( I,J;K,L)$ with $I,J\subset\{1,\ldots,k-l\}$ and $K,L\subset \{1,\ldots,l\}$ and  the convention 
\eqref{e:convention-z} is taken into account. For short  we also  write $dy_M$ instead of $dz_I\wedge d\bar z_J\wedge dw_K\wedge d\bar w_L.$

Consider also  the   subclass  $\widehat{\widetilde\Dc}^1_\ell$ consisting of  all $D\in\widetilde\DO^1_\ell$   with the following property:
There is an integer $n_D\geq 0$  such that 
given a current $S,$  by writing  $R:=DS,$  we have the following representation in the coordinates $\tilde y=(\tilde y_1,\ldots,\tilde y_k)$ as above: 
\begin{equation}\label{e:diff-oper-1-bis}
    R_M=\sum_{M',j,n}  \big(f_{M,M',j,n} {\partial S_{M'}\over \partial \tilde y_j} +  g_{M,M',j,n} {\partial S_{M'}\over \partial \overline{ \tilde y}_j}
    +h_{M,M',j,n} S_{M'}\big),
   \end{equation}
   with the  following two properties:
   \begin{enumerate}
   \item the sum is taken over all $M'$ with $|M'|=|M|-1=2k-q-1$ and $1\leq  j\leq k$ and $1\leq n\leq n_D.$
\item [$\widehat{(2)}$] $f_{M,M',j,},$ $g_{M,M',j,n},$ $h_{M,M',j,n}$  are  smooth functions  such that  
$$f_{M,M',j,n}(\tilde y)=O(t^{\Delta_j(M,M')}),\quad g_{M,M',j,n}(\tilde y)=O(t^{\Delta_j(M,M')})
\quad\text{ and}\quad h_{M,M',j,n}(\tilde y)=O(t^{\max(0,\Delta_j(M,M')-1)}).$$
\end{enumerate}
A differential operator $D\in  \widehat{\widetilde\Dc}^1_\ell$ is  said to be  {\it $\star$-fine of order $1$}.

   Consider  the following    property  (2) which is stronger than  propety  $\widehat{(2)}:$
 \begin{enumerate}
 \item  [(2)] $f_{M,M',j,},$ $g_{M,M',j,n},$ $h_{M,M',j,n}$  are  smooth functions  such that  
\begin{equation*}\begin{split}f_{M,M',j,n}(\tilde y )&=O(t^{\max(1,\Delta_j(M,M'))}),\quad g_{M,M',j,n}(\tilde y)=O(t^{\max(1,\Delta_j(M,M'))}),\\ h_{M,M',j,n}(\tilde y)&=O(t^{\max(0,\Delta_j(M,M')-1)})
                  \end{split}
\end{equation*}
 for all  $M,M',j.$
 \end{enumerate}
 If $D$ satisfies   both properties (1) and (2),  then  we say that $D$ is  a  {\it $\star$-negligible of order $1$}.
The set of all  $\star$-negligible  differential  operators  $D$ of order $1$ is  denoted by $\widetilde\Dc^1_\ell.$  So  $\widetilde\Dc^1_\ell\subset\widehat{\widetilde \Dc}^1_\ell.$
\end{definition}

\begin{lemma}\label{L:equiv-negli-forms}
 Let $\Phi$ be  a smooth  form  on $\Tube(B,\bfr).$  Then   
 $\Phi$ is $j$-negligible   if and only if  $\widetilde Y_*\Phi$ is  $j$-negligible, where $\widetilde Y$ is  the  diffeomorphism given in  \eqref{e:coor-tilde-y}.
\end{lemma}
\proof
By linearity it  suffices to show  the proposition for  the form $\Phi(y)=f(y) dy_I\wedge  d\bar y_J,$ 
where $f$ is  a  smooth function compactly supported in $\Tube(B,\bfr).$ 
Write  $I:=I'\cup\{ I''+(k-l)\}$ and  $J:=J'\cup \{J''+(k-l)\}$ for  $I',J'\subset \{1,\ldots,k-l\}$ and $I'',J''\subset \{1,\ldots,l\}.$
Here  $\{K+m\}:=\{j+m:\ j\in K\}$ for $K\subset \{1,\ldots,l\}$ and $0\leq m\leq k-l.$ 
We get that
\begin{equation*} \Phi:= f(z,w)dz_{I'}\wedge d\bar z_{J'}\wedge dw_{I''}\wedge d\bar w_{J''},
\end{equation*}
Note  that   by  \eqref{e:coor-tilde-y} and  by convention \ref{e:convention-z},
we have
\begin{equation}\label{e:widetilde-Y_z_j}
\begin{cases}  d(\widetilde Y_*z_j)=dz_j\quad{and}\quad  d(\widetilde Y_*\bar z_j)=d\bar z_j,&  \text{if}\quad j<k-l;\\
d(\widetilde Y_*z_{k-l})= d\tilde  y_k=du+idt\quad\text{ and}\quad  d(\widetilde Y_*\bar z_{k-l})=d\overline{\tilde y}_k=du-idt,&  \text{if}\quad j=k-l.
\end{cases} \end{equation}
Using this and  Definition \ref{D:negligible}, we  see that   $\Phi$ is $j$-negligible if and only if 
$$
\widetilde Y_*\Phi=f(Y(\tilde y))d(\widetilde Y_*z_{I'})\wedge d(\widetilde Y_*\bar z_{J'})\wedge d(\widetilde Y_*w_{I''})\wedge d(\widetilde Y_*
\bar w_{J''})
$$ is $j$-negligible.
\endproof

\begin{definition}\label{D:pull-push-opers}\rm  Let $D$  be  a  differential operator in the class  $\DO^1_\ell.$  Then we define
the differential operator $\widetilde Y_*D=Y^*D$ on $\M$ as  follows:
$$
\langle (\widetilde Y_*D)(\widetilde S),\widetilde\Phi\rangle_\M= \langle D(\widetilde Y^*(\widetilde  S)), \widetilde Y^*(\widetilde \Phi)\rangle_{\U_\ell} 
$$
 for all  current  $\widetilde S$ and  smooth  test forms $\widetilde \Phi$  on $\M.$ In other words,
 if $R:= D(\widetilde Y^*(\widetilde S)),$ then  $ \widetilde Y_*D(\widetilde S)=\widetilde Y_*(R).$
 
 Analogously, Let $\widetilde D$  be  a  differential operator in the class  $\widetilde\DO^1_\ell.$  Then we define
the differential operator $ Y_*\widetilde D=\widetilde Y^*\widetilde D$ on $\U_\ell$ as  follows:
$$
\langle ( Y_*\widetilde D)(S),\Phi\rangle_{\U_\ell}= \langle \widetilde D(Y^* S),  Y^*( \Phi)\rangle_{\M} 
$$
 for all  current  $ S$ and  smooth  test forms $\Phi$  on $\U_\ell.$ In other words,
 if $\widetilde R:=\widetilde D( Y^*( S)),$ then  $  (Y_*\widetilde  D)( S)= Y_*(\widetilde R).$
\end{definition}

\begin{lemma}\label{L:equiv-negli-opers}
  Let $D$  be  a  differential operator in the class  $\DO^1_\ell.$  Then
  \begin{enumerate} 
   \item $D\in  \widehat{\Dc}^1_\ell$ if and only if  $\widetilde Y_*D\in \widehat{\widetilde \Dc}^1_\ell.$
    \item $D\in  \Dc^1_\ell$ if and only if  $\widetilde Y_*D\in {\widetilde\Dc}^1_\ell.$
  \end{enumerate}
\end{lemma}
\proof
Let $\widetilde S$ be a  current   on  $\M.$ Set $S:= Y_*(\widetilde S)$ and $R:=  D(\widetilde Y^*(\widetilde S)).$
So  $\widetilde R=  (\widetilde Y_*D)(\widetilde S)= \widetilde Y_*(R).$ By \eqref{e:diff-oper-1} we have 
\begin{equation*}
    R_M=\sum_{M',j}  \big(f_{M,M',j} {\partial S_{M'}\over \partial  y_j} +  g_{M,M',j} {\partial S_{M'}\over \partial \bar  y_j}\big),
   \end{equation*}
   where  $f_{M,M',j},$ $g_{M,M',j}$ are  smooth functions   with    
   \begin{equation} \label{e:f,g_M,M'}
   \begin{cases} f_{M,M',j}(z,w)=O(\|z\|^{\Delta(M,M')})\quad\text{ and}\quad g_{M,M',j}(z,w)=O(\|z\|^{\Delta(M,M')}),  &  \text{if}\quad D\in  \widehat{\Dc}^1_\ell;\\
    f_{M,M',j}(z,w)=O(\|z\|^{\max(1,\Delta(M,M'))})\quad\text{ and}\quad g_{M,M',j}(z,w)=O(\|z\|^{\max(1,\Delta(M,M'))}),  &  \text{if}\quad D\in  \Dc^1_\ell. 
   \end{cases}
\end{equation}
   Therefore,
   \begin{equation}\label{e:equiv-negli-opers--widetilde-R}
\begin{split}    \widetilde R&=   \widetilde Y_*(R)
    =\sum_{M=(I,J;K,L)} \widetilde Y_*( R_M) d( \widetilde Y_*z_I)\wedge d(\widetilde Y_*\bar z_J)\wedge dw_K\wedge d\bar w_L\\
  &=  \sum_{M=(I,J;K,L)}  \big(\widetilde Y_*f_{M,M',j} \widetilde Y_*({\partial S_{M'}\over \partial  y_j}) +  \widetilde Y_*g_{M,M',j} \widetilde Y_*({\partial S_{M'}\over \partial \bar  y_j})\big)  d( \widetilde Y_*z_I)\wedge d(\widetilde Y_*\bar z_J)\wedge dw_K\wedge d\bar w_L . 
  \end{split}
  \end{equation}
  We  deduce from    \eqref{e:f,g_M,M'} and  \eqref{e:coor-tilde-y}--\eqref{e:coor-tilde-y-t} that
  \begin{equation}\label{e:f,g_M,M'-bis}
  \begin{cases} \widetilde Y_*f_{M,M',j}(z,w)=O(t^{\Delta(M,M')})\quad\text{ and}\quad \widetilde Y_*g_{M,M',j}(z,w)=O(t^{\Delta(M,M')}),  &  \text{if}\quad D\in  \widehat{\Dc}^1_\ell;\\
  \widetilde Y_*f_{M,M',j}(z,w)=O(t^{\max(1,\Delta(M,M'))})\quad\text{ and}\quad \widetilde Y_*g_{M,M',j}(z,w)=O(t^{\max(1,\Delta(M,M'))}), &  \text{if}\quad D\in  \Dc^1_\ell.
   \end{cases}
\end{equation}
  Moreover,  since  $\widetilde Y_*S=\widetilde S$  we have
  \begin{eqnarray*}
   \widetilde Y_*({\partial S_{M'}\over \partial  y_j})&= &{\partial \widetilde S_{M'}\over \partial  y_j}=\sum_{p=1}^k{\partial \widetilde S_{M'}\over \partial \tilde y_p}   {\partial \tilde y_p\over \partial  y_j}+\sum_{p=1}^k{\partial \widetilde S_{M'}\over \partial \overline{\tilde y}_p}   {\partial \overline{\tilde y}_p\over \partial  y_j}\\
  \widetilde Y_*({\partial S_{M'}\over \partial \bar y_j})&= &{\partial \widetilde S_{M'}\over \partial \bar y_j}=\sum_{p=1}^k{\partial \widetilde S_{M'}\over \partial \tilde y_p}   {\partial \tilde y_p\over \partial \bar y_j}+\sum_{p=1}^k{\partial \widetilde  S_{M'}\over \partial \overline{\tilde y}_p}   {\partial \overline{\tilde y}_p\over \partial \bar y_j} .
  \end{eqnarray*}
Here,
\begin{itemize}
 \item[$\bullet$] If  $j\not= k-l$ and $p\not=k-l,$ then ${\partial \tilde y_p\over \partial  \bar y_j}= {\partial \overline{\tilde y}_p\over \partial  y_j}=0$ and ${\partial \tilde y_p\over \partial  y_j}= {\partial \overline{\tilde y}_p\over \partial \bar y_j}=\delta_{jp},$ where $\delta_{jp}=1$ if  $j=p$ and   $\delta_{jp}=0$  otherwise.
 \item[$\bullet$] If  $p= k-l,$  then  we deduce from  \eqref{e:coor-tilde-y-t-der}  that
  $${\partial \tilde y_p\over \partial  y_j}={\partial (u+it)\over \partial y_j}=\begin{cases}
                                                                             O(t),& \text{if}\ \ k-l<j\leq k,\\
                                                                             O(1) ,& \text{if}\ \ j\leq k-l.
                                                                            \end{cases}
  $$
  Similar estimates hold for ${\partial \tilde y_p\over \partial \bar  y_j}.$  Moreover,  we have that 
  $${\partial\overline{ \tilde y}_p\over \partial  y_j}={\partial (u-it)\over \partial y_j}=\begin{cases}
                                                                             O(t),& \text{if}\ \ k-l<j\leq k,\\
                                                                             O(1) ,& \text{if}\ \ j\leq k-l.
                                                                            \end{cases}
  $$
  Similar estimates hold for ${\partial\overline{ \tilde y}_p\over \partial \bar  y_j}.$
   \item[$\bullet$] If  $j= k-l$ and $p\not=k-l,$ then ${\partial \tilde y_p\over \partial  y_j}= {\partial \overline{\tilde y}_p\over \partial \bar y_j}
   ={\partial \tilde y_p\over \partial \bar  y_j}= {\partial \overline{\tilde y}_p\over \partial  y_j}=0.$
\end{itemize}
Putting the  above estimates and equalities together  with  \eqref{e:f,g_M,M'-bis} and \eqref{e:widetilde-Y_z_j} into the RHS of \eqref{e:equiv-negli-opers--widetilde-R}, we  infer that
if $D\in  \widehat{\Dc}^1_\ell$ then  $\widetilde Y_*D\in \widehat{\widetilde \Dc}^1_\ell,$  and if 
     $D\in  \Dc^1_\ell$ then  $\widetilde Y_*D\in {\widetilde\Dc}^1_\ell.$
     
   The converse  implications can be proved  similarly.  We leave them to the interested reader.
\endproof

For  $\tilde y=(\tilde y_1,\ldots,\tilde y_k)\in \M\subset\C^k,$  write
$\tilde y_j=u_j+it_j,$ where $u_j,$ $t_j\in\R.$  Note that   for $j=k-l,$ $u_{k-l}=u$ and $t_{k-l}=t$ and  hence $\tilde y_{k-l}=u+it.$
For $0<r\leq \bfr$ consider the   real  hyperplane  $\M(r):=\{\tilde y\in\M:\   t=r\}.$ For $0<r_1<r_2\leq\bfr,$ consider the strip
$\M(r_1,r_2):=\{\tilde y\in\M:\  r_1< t<r_2\}.$
\begin{lemma}\label{L:integration-by-parts}
 Let  $S$ be  a  distribution on $\M(r)$ and  $f$ a smooth function on $\M_r.$
 Let $1\leq j\leq k.$ Let $d\Leb_{2k-1}$ be the  Lebesgue measure on $\M_r.$
 Then:
 \begin{enumerate}
  \item If $  j\not=k-l,$ then  
  \begin{eqnarray*} \int_{\M(r)}  {\partial  S\over  \partial  \tilde y_j} f  d\Leb_{2k-1}&=&-\int_{\M(r)}  S{\partial  f\over  \partial  \tilde y_j}   d\Leb_{2k-1}  +  \int_{\partial\M(r)}  ( Sf)(\tilde y) \big({\partial \over \partial \tilde y_j}\intprod \Leb_{2k-1}(\tilde y)\big),\\
  \int_{\M(r)}  {\partial  S\over  \partial \overline{ \tilde y}_j} f  d\Leb_{2k-1}&=&-\int_{\M(r)}  S{\partial  f\over  \partial \overline{ \tilde y}_j}   d\Leb_{2k-1}+\int_{\partial\M(r)}  ( Sf) (\tilde y)\big({\partial \over \partial \overline{\tilde y}_j}\intprod
  \Leb_{2k-1}(\tilde y)\big).
  \end{eqnarray*}
  \item  If  $j=k-l,$ then 
  \begin{eqnarray*}
   \int_{\M(r)}  {\partial  S\over  \partial  \tilde y_{k-l}} f  d\Leb_{2k-1}&=&-{1\over 2}\big( \int_{\M(r)}  S{\partial  f\over  \partial   u}   d\Leb_{2k-1} -i
    \int_{\M(r)}  {\partial  S\over  \partial   t} f  d\Leb_{2k-1}\big)\\
    &+&{1\over 2}\int_{\partial\M(r)}  ( Sf)(\tilde y) \big({\partial \over \partial \tilde u}\intprod \Leb_{2k-1}(\tilde y)\big)
    ,\\
     \int_{\M(r)}  {\partial  S\over  \partial  \overline{ \tilde y}_{k-l}} f  d\Leb_{2k-1}&=&-{1\over 2}\big(\int_{\M(r)}  S{\partial  f\over  \partial   u}   d\Leb_{2k-1} +i
    \int_{\M(r)}  {\partial  S\over  \partial   t} f  d\Leb_{2k-1}\big)\\
     &+&{1\over 2}\int_{\partial\M(r)}  ( Sf)(\tilde y) \big({\partial \over \partial \tilde u}\intprod \Leb_{2k-1}(\tilde y)\big)
    .
  \end{eqnarray*}
\item  Moreover,  for   $j=k-l$ and  $0<r_1<r_2\leq\bfr,$ we have 
\begin{equation*}
 \int_{r_1}^{r_2}\big(\int_{\M(t)}  {\partial  S\over  \partial   t} f  d\Leb_{2k-1}\big)dt=\int_{\M(r_1,r_2)}S {\partial  f\over  \partial   t}   d\Leb_{2k}+\int_{\M(r_2)}S   f   d\Leb_{2k-1}-\int_{ \M(r_1)}S   f   d\Leb_{2k-1}.
\end{equation*}
 \end{enumerate}

\end{lemma}

\begin{proposition}\label{P:boundary-vs-tube-eps-bis}
Let $S$ be a  positive plurisubharmonic  current of  bidimension   $(q,q)$ on a neighborhood of $\Tube(B,\bfr)$ such that $S$ and $\ddc S$  such that $S$ is  $\Cc^1$-smooth   near $ \partial_\ver\Tube(B,\bfr).$  
\begin{itemize}
 \item[(i-1)] Suppose that  $D$ is  a  differential operator in the class $\widehat\Dc^1_\ell$ and   $\Phi$ is a smooth form of  degree $2q-1$  which is  $m$-negligible.
 Then there are:
 \begin{itemize}
  \item [$\bullet$]   six differential operators $D_j$ for $1\leq j\leq 6$  in  the class $\widehat\Dc^0_\ell;$
  \item [$\bullet$]  A  bounded form $S_0$  of dimension $2q$   which is defined on a neighborhood of $\partial_\ver\Tube(B,\bfr);$
  \item[$\bullet$]  five  $2q$-forms $\Phi_1$   which is  $(m+2)$-negligible, $\Phi_2$   which is  $(m+1)$-negligible,
   $\Phi_3$  which  is  $m$-negligible,  $\Phi_4$  which  is  $(m+1)$-negligible,   and  $\Phi_5,$  $\Phi_6$    which  are both $m$-negligible;
 \end{itemize}
  such that if we set, for $  0<t\leq \bfr$:
  \begin{equation}\label{e:Ic_D(t)}
  \begin{split}
  \Ic_D(t)&:=\int_{\partial_\hor \Tube(B,t)} DS\wedge \Phi -  \int_{\partial_\ver\Tube(B,t)} S_0\wedge \Phi-  \int_{\Tube(B,t)} D_1S\wedge \Phi_1\\
  & -{1\over  t} \int_{\Tube(B,t)} D_2S\wedge \Phi_2-{1\over  t^2} \int_{\Tube(B,t)} D_3S\wedge \Phi_3, 
  \end{split}
  \end{equation}
then for every $0<r_1<r_2\leq\bfr$ and  every smooth  function $\chi$ on $(0,\bfr),$  we have
 \begin{equation}\label{e:Stokes-DS-order-1}
 \begin{split}
  \int_{r_1}^{r_2}\chi(t) \Ic_D(t)dt&= 
\int_{\Tube(B,r_1,r_2)}\chi(\|y\|) (D_4S\wedge \Phi_4)(y)+\int_{\Tube(B,r_1,r_2)}\chi'(\|y\|) (D_5S\wedge \Phi_5)(y)\\
&+ \int_{\partial_\hor\Tube(B,r_2) }\chi(r_2)( D_6S\wedge \Phi_6)(y)- \int_{\partial_\hor\Tube(B,r_1)} \chi(r_1)(D_6S\wedge \Phi_6)(y)  . 
 \end{split}
 \end{equation}

 \item[(ii-1)] Suppose  that  $D$ is  a  differential operator in the class $\Dc^1_\ell$ and   $\Phi$ is a smooth form of  degree $2q-1$    which is  $m$-negligible.
 Then  the  conclusion of assertion (i-1) holds. Moreover,  
 the six differential operators $D_j$ for $1\leq j\leq 6$  belong to   the class $\Dc^0_\ell.$
   
 \end{itemize}
 
\end{proposition}
\proof 
Consider the current  $\widetilde S$   on  $\M$ given by   $\widetilde S:= \widetilde Y_*( S)$ and  the  differential  operator 
$\widetilde D$    on  $\M$ given by   $\widetilde D:=  \widetilde Y_*D.$  Set $\widetilde R:=  \widetilde D(\widetilde S).$
By     Lemma \ref{L:equiv-negli-opers} and  \eqref{e:diff-oper-1} we have 
\begin{equation*}
   \widetilde R_M=\sum_{M',j}  \big(\tilde f_{M,M',j} {\partial  \widetilde S_{M'}\over \partial  \tilde y_j} + \tilde g_{M,M',j} {\partial \widetilde S_{M'}\over \partial \overline{ \tilde y}_j}\big),
   \end{equation*}
   where  $\tilde f_{M,M',j},$ $\tilde g_{M,M',j}$ are  smooth functions   satisfying      
     \begin{equation} \label{e:f,g_M,M'-tilde}
   \begin{cases} \tilde f_{M,M',j}(\tilde y)=O(t^{\Delta(M,M')})\quad\text{ and}\quad \tilde g_{M,M',j}(\tilde y)=O(t^{\Delta(M,M')}),  &  \text{if}\quad D\in  \widehat{\Dc}^1_\ell;\\
   \tilde f_{M,M',j}(\tilde y)=O(t^{\max(1,\Delta(M,M'))})\quad\text{ and}\quad \tilde g_{M,M',j}(\tilde y)=O(t^{\max(1,\Delta(M,M'))}),  &  \text{if}\quad D\in  \Dc^1_\ell. 
   \end{cases}
\end{equation}
 Consider the smooth form  $\widetilde \Phi:= \widetilde Y_*( \Phi)$ on $\M.$   Observe  that    
 \begin{equation*}
  \int_{\partial_\hor \Tube(B,r)\cap \U_\ell} DS\wedge \Phi=\int_{\M\cap  \{t=r\}}  \widetilde  D\widetilde S\wedge \widetilde \Phi
  =\sum_M \int_{\M(r)}  \widetilde R_M dy_M \wedge \widetilde \Phi\\
 \end{equation*}
This, combined  with  the above expression of   $\widetilde R_M,$  gives  that 
\begin{equation*}
   \int_{\partial_\hor \Tube(B,r)\cap \U_\ell} DS\wedge \Phi=\sum_{M,M',j} I_{M,M',j}(r),
   \end{equation*}
   where 
\begin{equation*}I_{M,M',j}(r):=\int_{\M(r)}    \big(\tilde f_{M,M',j} {\partial  \widetilde S_{M'}\over \partial  \tilde y_j} + \tilde g_{M,M',j} {\partial \widetilde S_{M'}\over \partial \overline{ \tilde y}_j}\big) d\tilde y_M \wedge \widetilde \Phi.
\end{equation*}
Write 
$$
\widetilde \Phi=\sum_M \widetilde \Phi(P) d\tilde y_P,
$$
where the sum is taken over all $P$  with $|P|=2q-1$ 
and   $ \widetilde \Phi(P)$ are  distributions coefficients.   Since  the above integrals are performed on $\M(r),$
we  see that if   $d\tilde y_M$ contains  $dt$ or if  $d\tilde y_P$ contains $dt$  or if $d\tilde y_M$ and 
$d\tilde y_P$ contains a common factor,
then the corresponding integral 
$$\int_{\M(r)}    \big(\tilde f_{M,M',j} {\partial  \widetilde S_{M'}\over \partial  \tilde y_j} + \tilde g_{M,M',j} {\partial \widetilde S_{M'}\over \partial \overline{ \tilde y}_j}\big) d\tilde y_M \wedge \widetilde \Phi(P) d\tilde y_P=0.$$
 So we only need to treat  every $M$ such that  $|M|=2k-2q$ and  that $d\tilde y_M$ does not contain $dt.$ For such  a multi-index $M$  there is a unique multi-index $P$ such that $|P|= 2q-1$ and  $$d\tilde y_M\wedge d\tilde y_P= ({\partial \over \partial t} \intprod   d\Leb_{2k}(\tilde y)\big)=du\wedge \prod_{q=1}^{l}idw_q\wedge d\bar w_q\wedge \prod_{p=1}^{k-l-1} idz_p\wedge d\bar z_p  .$$ Write $\widetilde \Phi_M:= \widetilde \Phi(P).$ Since   $\Phi$ is  $m$-negligible, so is  $ \Phi_M d\tilde y_P.$ 
We infer that
\begin{equation*}I_{M,M',j}(r):=\int_{\M(r)}    \big(\tilde f_{M,M',j} {\partial  \widetilde S_{M'}\over \partial  \tilde y_j} + \tilde g_{M,M',j} {\partial \widetilde S_{M'}\over \partial \overline{ \tilde y}_j}\big)  \widetilde \Phi_M d\Leb_{2k-1}(\tilde y).
\end{equation*}
Since $t=r$  for $\tilde y\in   \M(r),$ it  follows that
\begin{equation*}I_{M,M',j}(r):={1\over  r}\int_{\M(r)}    \big(t\tilde f_{M,M',j}  \widetilde \Phi_M{\partial  \widetilde S_{M'}\over \partial  \tilde y_j} + t\tilde g_{M,M',j}  \widetilde \Phi_M{\partial \widetilde S_{M'}\over \partial \overline{ \tilde y}_j}\big)  d\Leb_{2k-1}(\tilde y).
\end{equation*}

 To  handle  the   integral $I_{M,M',j},$ we  consider  two cases  according to the value  of $j.$
 
  \noindent  {\bf Case $j\not=k-l:$} Applying  Lemma \ref{L:integration-by-parts} (1) for a given $r\in (0,\bfr]$ yields that 
 \begin{equation*}
   I_{M,M',j}(r)={1\over r}\int_{\M(r)}    \big(\widetilde S_{M'}  {\partial( t\tilde f_{M,M',j}  \widetilde \Phi_M) \over \partial  \tilde y_j} +  \widetilde S_{M'} {\partial  (t\tilde g_{M,M',j}   \widetilde \Phi_M)\over \partial \overline{ \tilde y}_j}\big)d\Leb_{2k-1}(\tilde y) .
\end{equation*}
Let $Q\subset  M$  be such that $|Q|=|M|-1=2k-2q-1.$ So  $d\tilde y_{Q^\bfc} = d\tilde y_P\wedge dt \wedge   dx,$
where $dx\in\{dz_p,d\bar z_p,dw_q,d\bar w_q\}.$   Thus, we can write for  $\tilde y \in \M(r),$
\begin{eqnarray*}
 \widetilde S_{M'} {\partial( t\tilde f_{M,M',j}  \widetilde \Phi_M) \over \partial  \tilde y_j}d\Leb_{2k-1}(\tilde y)
 &= &
\widetilde S_{M'}\tilde f_{M,M',j} t{ \widetilde \Phi_M\over \partial  \tilde y_j}d\Leb_{2k-1}(\tilde y)
+\widetilde S_{M'} \big(t{\partial \tilde f_{M,M',j}\over \partial  \tilde y_j}  \big)  \widetilde \Phi_Md\Leb_{2k-1}(\tilde y)\\
&:= & E_1+E_2.
\end{eqnarray*}
To handle  $E_1,$  we rewrite it as follows:
$$
E_1:=\tilde f_{M,M',j}\widetilde S_{M'}d\tilde y_{Q} \wedge \widetilde{\widetilde\Phi }_M,
$$
where  $\widetilde{\widetilde\Phi }_M:= t {\partial  \widetilde \Phi_M\over \partial  \tilde y_j}d\tilde y_P\wedge   dx.$
Recall that  $ \Phi_M d\tilde y_P$ is  $m$-negligible. There are two subcases.
 
 \noindent  {\bf Subcase  $j\in\{1,\ldots,k-l\}$ and  $dx\in\{dw,d\bar w\}$:}  Using Definition \ref{D:negligible}, we check that
 $\widetilde{\widetilde\Phi }_M$ is $(m+1)$-negligible. Note that $\delta_{j,Q,M,M'}=1$ in the notation of Definition  \ref{D:substraction-bis}.
 By Definition  \ref{D:substraction-bis}, we can choose $Q\subset M$ such that $\Delta_j(M,M')=\Delta(Q,M') +1.$
Then we have
\begin{equation*}
 \tilde f_{M,M',j}\widetilde S_{M'}d\tilde y_{Q} \wedge \widetilde{\widetilde\Phi }_M=
 \tilde{\tilde f}_{Q,M'}\widetilde S_{M'}d\tilde y_{Q} \wedge \widetilde{\widetilde\Phi },
\end{equation*}
where $\tilde{\tilde f}_{Q,M'}(\tilde y):= t^{-1}\tilde f_{M,M',j} = O(t^{\Delta(Q,M')})$ and $\widetilde{\widetilde\Phi }:=t \widetilde{\widetilde\Phi }_M $ is $m$-negligible.

 \noindent  {\bf Subcase: the remaining  subcase. } Since we have either  $j\not\in\{1,\ldots,k-l\}$ or  $dx\not\in\{dw,d\bar w\}$,  using Definition \ref{D:negligible} we check that
 $\widetilde{\widetilde\Phi }_M$ is $m$-negligible. Note that $\delta_{j,Q,M,M'}=0$ in the notation of Definition  \ref{D:substraction-bis}.
 By Definition  \ref{D:substraction-bis}, we can choose $Q\subset M$ such that $\Delta_j(M,M')=\Delta(Q,M') .$
Then we have
\begin{eqnarray*}
 \tilde f_{M,M',j}\widetilde S_{M'}d\tilde y_{Q} \wedge \widetilde{\widetilde\Phi }_M=\tilde{\tilde f}_{Q,M'}\widetilde S_{M'}d\tilde y_{Q} \wedge \widetilde{\widetilde\Phi },
\end{eqnarray*}
where $\tilde{\tilde f}_{Q,M'}(\tilde y):= \tilde f_{M,M',j} = O(t^{\Delta(Q,M')})$ and $\widetilde{\widetilde\Phi }:= \widetilde{\widetilde\Phi }_M $ is $m$-negligible.

In summary, we have shown that 
\begin{equation*}
 E_1=
 \tilde{\tilde f}_{Q,M'}\widetilde S_{M'}d\tilde y_{Q} \wedge \widetilde{\widetilde\Phi },
\end{equation*}
where $\tilde{\tilde f}_{Q,M'}(\tilde y)  = O(t^{\Delta(Q,M')})$ and $\widetilde{\widetilde\Phi }  $ is $m$-negligible.
Similarly,   we can show that $E_2$ is also of this form.
Consequently,  there are a  differential operators $\widetilde D_{1j,M,M'}$ (depending on $M$ and $M'$)   in  the class $\widehat{\widetilde \Dc}^0_\ell$ 
and a test  form  $ \widetilde\Phi_{1,j,M,M'}$ of degree $2q-1$ which is  $m$-negligible such that 
$$
 I_{M,M',j}(r)={1\over r}\int_{\M(r)} D_{1,j,M,M'}S\wedge \Phi_{1,j,M,M'}.
$$
Therefore, there are a  differential operators $ D_{1j,M,M'}$ (depending on $M$ and $M'$)   in  the class $\widehat\Dc^0_\ell$ 
and a test  form  $ \Phi_{1,j,M,M'}$ of degree $2q-1$ which is  $m$-negligible such that 
$$
 I_{M,M',j}(r)={1\over r}\int_{\partial_\hor\Tube (B,r)} D_{1,j,M,M'}S\wedge \Phi_{1,j,M,M'}.
$$

\noindent  {\bf Case $j=k-l:$} By Lemma \ref{L:integration-by-parts} (2),  we have that  $I_{M,M',k-l}= {1\over 2}(-I^1_{M,M',k-l}+i I^2_{M,M',k-l}),$ where
\begin{eqnarray*} I^1_{M,M',k-l}(r)&:=&\int_{\M(r)}    \big(   \widetilde S_{M'} {\partial ( \tilde f_{M,M',j}  \widetilde \Phi_M)\over \partial  u} +  \widetilde S_{M'} {\partial (\tilde g_{M,M',j}   \widetilde \Phi_M) \over \partial u}\big)d\Leb_{2k-1}(\tilde y),\\
I^2_{M,M',k-l}(r)&:=&\int_{\M(r)}      \big(  \widetilde S_{M'} {\partial (\tilde f_{M,M',j}  \widetilde \Phi_M) \over \partial  t} -    \widetilde S_{M'} {\partial (\tilde g_{M,M',j}   \widetilde \Phi_M)\over \partial t}\big)d\Leb_{2k-1}(\tilde y).
\end{eqnarray*}
We handle  $I^1_{M,M',k-l}(r)$ in the same way as  for $ I_{M,M',j}(r)$ with $j\not=k-l.$
Therefore, there are a  differential operators $ D_{1,k-l,M,M'}$ (depending on $M$ and $M'$)   in  the class $\widehat\Dc^0_\ell$ 
and a test  form  $ \Phi_{1,k-l,M,M'}$ of degree $2q-1$ which is  $m$-negligible such that 
$$
 I^1_{M,M',k-l}(r)={1\over r}\int_{\partial_\hor\Tube (B,r)} D_{1,k-l,M,M'}S\wedge \Phi_{1,k-l,M,M'}.
$$
Set
\begin{equation}\label{e:Ic_D(t)-bis}
 \Ic_D(t):= \int_{\partial_\hor \Tube(B,t)} DS\wedge \Phi -\sum_{M,M'}\big(\sum_{j\not=k-l} I_{M,M',j}(t)+  I^1_{M,M',k-l}(t)\big).
\end{equation}
By  the above  discussion,  we see that
\begin{equation}\label{e:Ic_D(t)-bisbis}
  \Ic_D(t)=  \sum_{M,M'}   I^2_{M,M',k-l}(t).
\end{equation}
Moreover,  since  formula \eqref{e:Ic_D} yields that 
\begin{equation}\label{e:Ic_D}
  \int_{\partial_\hor \Tube(B,t)} DS\wedge \Phi- \Ic_D(t)=\sum_{M,M'}\big(\sum_{j\not=k-l} I_{M,M',j}(t)+  I^1_{M,M',k-l}(t)\big),
\end{equation}
we may find an integer $n\geq 1$ and $n$ differential operators $ D_{n}$    in  the class $\widehat\Dc^0_\ell$ 
and $n$ test  form  $ \Phi_{1,k-l,M,M'}$ of degree $2q-1$ which are  $m$-negligible such that 
$$ \int_{\partial_\hor \Tube(B,t)} DS\wedge \Phi- \Ic_D(t)={1\over r}\int_{\partial_\hor\Tube (B,r)} D_{n}S\wedge \Phi_n.
$$
Applying  Proposition \ref{P:boundary-vs-tube-eps}  to the  RHS  of the last line,  we obtain   identity \eqref{e:Ic_D(t)}. 

 By Lemma \ref{L:integration-by-parts} (3),  we obtain, 
for  $0<r_1<r_2\leq\bfr,$ that
\begin{multline*}
\int_{r_1}^{r_2}\chi(t) I^2_{M,M',j}(t)dt=\int_{\M(r_1,r_2)}    \big(\widetilde S_{M'}  { \partial(  \chi(t) \tilde f_{M,M',j}  \widetilde \Phi_M) \over \partial  t} - \widetilde S_{M'} {\partial  ( \chi(t) \tilde g_{M,M',j}   \widetilde \Phi_M)\over \partial t}\big)d\Leb_{2k}(\tilde y)\\
+  \int_{\M(r_2)}    \chi(r_2) \widetilde S_{M'} \widetilde \Phi_M( \tilde f_{M,M',j} -\tilde g_{M,M',j} )d\Leb_{2k-1}(\tilde y)
- \int_{\M(r_1)}   \chi(r_1) \widetilde S_{M'}  \widetilde \Phi_M( \tilde f_{M,M',j} -\tilde g_{M,M',j} )d\Leb_{2k-1}(\tilde y)\\
:=F_1+F_2.
\end{multline*}
Since we have
 \begin{equation*}  
\widetilde S_{M'}  { \partial(  \chi(t) \tilde f_{M,M',j}  \widetilde \Phi_M) \over \partial  t}=
\chi'(t)\widetilde S_{M'} (   \tilde f_{M,M',j}  \widetilde \Phi_M) +
 \chi(t)\widetilde S_{M'}  { \partial(  \tilde f_{M,M',j}  \widetilde \Phi_M) \over \partial  t},
 \end{equation*}
 we can find two differential operators $D_4$ and $D_5$ in the class of $\widehat{\Dc}^0_\ell$ and 
  two  $2q$-forms    $\Phi_4$  which  is  $(m+1)$-negligible   and  $\Phi_5$   which  is $m$-negligible;
  such that  
 \begin{equation*} 
  F_1= 
\int_{\Tube(B,r_1,r_2)}\chi(\|y\|) (D_4S\wedge \Phi_4)(y)+\int_{\Tube(B,r_1,r_2)}\chi'(\|y\|) (D_5S\wedge \Phi_5)(y).
 \end{equation*}
  We can check that there is a differential operators $D_6$   in  the class $\widehat\Dc^0_\ell$
  and a  $2q$-form    $\Phi_6$    which  is both $m$-negligible such that
 \begin{equation*}  F_2= \int_{\partial_\hor\Tube(B,r_2) }\chi(r_2)( D_6S\wedge \Phi_6)(y)- \int_{\partial_\hor\Tube(B,r_1)} \chi(r_1)(D_6S\wedge \Phi_6)(y)  . 
 \end{equation*}
 We obtain   identity \eqref{e:Stokes-DS-order-1}. The proof of  assertion (i-1) is thereby completed.
  \endproof

\begin{proposition}\label{P:basic_boundary-inequ}
Let $S$ be a   current in the class $\SH^{2,1}(\Tube(B,\bfr)).$ Let $D$ be  a  differential operator    and $\Phi$   a smooth form  on $\U_\ell.$ 

\begin{enumerate}
 \item    If  $D$ is  a  differential operator in the class $\widehat\Dc^0_\ell$ and   $\Phi$ is a  form of  degree $2q-1$  which is  $m$-negligible,
then  
for all  $r\in (0,\bfr]$ and  $s\in ({r\over 2},r),$ 
$$
{1\over  r^{2(k-p)-m} }\int_{r\over 2}^r\big|\int_{\partial_\hor \Tube(B,t)} DS\wedge \Phi\big|dt \leq 
c \sum_{j=\lowm}^\upm \nu_j( S,B,r,\id) , 
$$
where $c$ is a constant that depends only on  $\Phi,$ and  the LHS is  defined  using  Definition \ref{D:boundary-value}  for  all $r\in (0,\bfr]$ except   at most a  countable set.

\item   
If  $D$ is  a  differential operator in the class $\Dc^0_\ell$ and   $\Phi$ is a smooth form of  degree $2q-1$  which is  $m$-negligible,
then  
for all  $r\in (0,\bfr]$ and  $s\in ({r\over 2},r),$ 
$$
{1\over  r^{2(k-p)-m} }\int_{r\over 2}^r\big|\int_{\partial_\hor \Tube(B,t)} DS\wedge \Phi\big|dt \leq 
cr \sum_{j=\lowm}^\upm \nu_j( S,B,r,\id) , 
$$
where $c$ is a constant that depends only on $S$ and $\Phi,$ and  the LHS is  defined  using  Definition \ref{D:boundary-value}  for  all $r\in (0,\bfr]$ except   at most a  countable set.

\item   
If $D$ is  a  differential operator in the class $\widehat\Dc^1_\ell$ and   $\Phi$ is a smooth form of  degree $2q-1$  which is  $m$-negligible, then   the  function   $\Ic_D$ defined  by \eqref{e:Ic_D(t)} satisies the following inequality  for all $  0<r\leq \bfr:$
  \begin{equation}\label{e:inequal-order-1}
 {1\over  r^{2(k-p)-m} }\int_{r\over 2}^r \big| \int_{\partial_\hor \Tube(B,t)} DS\wedge \Phi -  \Ic_D(t)\big|dt \leq 
c\sum_{j=\lowm}^\upm \nu_j( S,B,r,\id) , 
 \end{equation}
 where $c$ is a constant that depends only on  $\Phi.$
 \item   If $D$ is  a  differential operator in the class $\Dc^1_\ell$ and   $\Phi$ is a smooth form of  degree $2q-1$    which is  $m$-negligible,  
 then   the  function   $\Ic_D$ defined  by \eqref{e:Ic_D(t)} satisies the following inequality  for all $  0<r\leq \bfr:$
  \begin{equation}\label{e:inequal-order-2}
 {1\over  r^{2(k-p)-m} }\int_{r\over 2}^r \big| \int_{\partial_\hor \Tube(B,t)} DS\wedge \Phi -  \Ic_D(t)\big|dt \leq 
cr\sum_{j=\lowm}^\upm \nu_j( S,B,r,\id) , 
 \end{equation}
 where $c$ is a constant that depends only on  $\Phi.$
 \end{enumerate}

\end{proposition}
\proof   

\noindent {\bf Proof of assertion (1):}
We are in  Case (i-0)  in Propositions \ref{P:boundary-vs-tube-eps}.

By  Propositions \ref{P:boundary-vs-tube-eps}  (i-0), there are:
\begin{itemize} 
 \item  [$\bullet$] a  bounded  form $S_0$ is   in a neighborhood of 
 $\partial_\ver\Tube(B,r)$  which depends only on $D$ and $S;$
  \item  [$\bullet$] 
  three differential operators $D_1,$ $D_2$ and $D_3$ in   the class $\widehat\Dc^0_\ell;$
  \item[$\bullet$] and three   forms
 $\Phi_1$ of degree $2q$  which is $(m+1)$-negligible,
 $\Phi_2$ of degree $2q$  which is $m$-negligible and 
 $\Phi_3$ of degree $(2q-1)$  which is  $m$-negligible;
 \end{itemize}
 such that  for every $0<s\leq\bfr,$  we have
 \begin{equation*} 
 \begin{split}
 \int_{\partial_\hor \Tube(B,s)} DS\wedge \Phi&= \int_{\partial_\ver\Tube(B,s)} S_0\wedge \Phi + \int_{\Tube(B,s)} D_1S\wedge \Phi_1 \\
 &+  {1\over r} \int_{\Tube(B,s)} D_2S\wedge \Phi_2+\lim\limits_{\epsilon\to 0+}\int_{\Tube(B,s-\epsilon,r)} D_3S\wedge d\chi_{s,\epsilon}\wedge \Phi_3,
 \end{split}
 \end{equation*}
 Integrating  both sides with respect to $s\in[{r\over 2},r]$ and  
 applying 
 Proposition \ref{P:boundary-vs-tube} to the  last term on the RHS, there is  a constant $c>0$ such that  for all $0<r\leq\bfr,$  we have  
 \begin{eqnarray*}
 \big|\int_{r\over 2}^r\big|\int_{\partial_\hor \Tube(B,s)} DS\wedge \Phi\big|ds &\leq& \int_{r\over 2}^r\big| \int_{\partial_\ver\Tube(B,s)} S_0\wedge \Phi\big| ds + \int_{r\over 2}^r\big| \int_{\Tube(B,s)} D_1S\wedge \Phi_1\big| ds\\
 &+&\int_{r\over 2}^r\big| {1\over r}\int_{\Tube(B,s)} D_2S\wedge \Phi_2\big| ds+ c\int_{\Tube(B,{r\over 2},r)} S\wedge  R^\dagger_{k-q,m}. 
 \end{eqnarray*}
 Since  $S_0$ and $\Phi$ are bounded in a neighborhood of $\partial_\ver\Tube(B,r),$
 the first term on the RHS is of order $O(r^{2k}).$
 Dividing both sides  by   $r^{2(k-p)-m},$ using the above description of  $D_j,$ $\Phi_j$ for $1\leq j\leq 3,$
 and then applying  Proposition \ref{P:mass-fine-negligible} (i)
 to  the second and third terms on the RHS, and then applying 
  Proposition \ref{P:mass-fine-negligible-bis} (i)   to the last on the RHS,
 the result follows.

\noindent {\bf Proof of assertion (2):} We are in  Case (ii-0)  in Propositions \ref{P:boundary-vs-tube-eps}.

We argue  as in the  proof of Case (i-0) using Propositions \ref{P:boundary-vs-tube-eps}  (ii-0), \ref{P:mass-fine-negligible} (ii) and  \ref{P:mass-fine-negligible-bis} (ii)  instead of Propositions \ref{P:boundary-vs-tube-eps}  (i-0), \ref{P:mass-fine-negligible} (i) and  \ref{P:mass-fine-negligible-bis} (i) respectively.  This completes  the proof of assertion (2).

\noindent {\bf Proof of assertion (3):} We are in Case (i-1)  in Propositions \ref{P:boundary-vs-tube-eps-bis}.

 By Proposition  \ref{P:boundary-vs-tube-eps-bis} (i-1) and formula \eqref{e:Ic_D(t)},
 there are:
 \begin{itemize}
  \item [$\bullet$]   five differential operators $D_j$ for $1\leq j\leq 5$  in  the class $\widehat\Dc^0_\ell;$
  \item [$\bullet$]  A  bounded form $S_0$  of dimension $2q$   which is defined on a neighborhood of $\partial_\ver\Tube(B,\bfr);$
  \item[$\bullet$]  five  $2q$-forms $\Phi_1$   which is  $(m+2)$-negligible,
   $\Phi_2$  which  is  $m$-negligible,  $\Phi_3$  which  is  $(m+1)$-negligible,   and  $\Phi_4,$  $\Phi_5$    which  are both $m$-negligible;
 \end{itemize}
  such that for $  0<t\leq \bfr$:
  \begin{equation*}
 \int_{\partial_\hor \Tube(B,t)} DS\wedge \Phi - \Ic_D(t)=  \int_{\partial_\ver\Tube(B,t)} S_0\wedge \Phi+
 \int_{\Tube(B,t)} D_1S\wedge \Phi_1 +{1\over  t^2} \int_{\Tube(B,t)} D_2S\wedge \Phi_2.  \end{equation*}
 Integrating  both sides with respect to $s\in[{r\over 2},r]$ yields that 
 \begin{multline*}
  \int_{r\over 2}^r\big|\int_{\partial_\hor \Tube(B,t)} DS\wedge \Phi - \Ic_D(t)\big|dt \leq
   \int_{r\over 2}^r\big|\int_{\partial_\ver\Tube(B,t)} S_0\wedge \Phi\big|dt\\
   +
 \int_{r\over 2}^r\big|\int_{\Tube(B,t)} D_1S\wedge \Phi_1 \big|dt+ \int_{r\over 2}^r\big| {1\over  t^2} \int_{\Tube(B,t)} D_2S\wedge \Phi_2\big|dt. 
 \end{multline*}
Arguing  as in 
the proof of assertion (1)  (without using Proposition \ref{P:mass-fine-negligible-bis} (i)), we can dominate  all terms on the RHS. This proves assertion (3).

\noindent {\bf Proof of assertion (4):}  We are in  Case (ii-1)  in Propositions \ref{P:boundary-vs-tube-eps-bis}.

 We argue as in the  proof of assertion (2). Indeed, we apply  Proposition  \ref{P:boundary-vs-tube-eps-bis} (ii-1)
 instead of  Proposition  \ref{P:boundary-vs-tube-eps-bis} (i-1)
 and use the proof of assertion (3) (instead of the  proof of assertion (1)).
\endproof
\subsection{Basic boundary estimates}

 We use the coordinate $y=(z,w)\in\C^{k-l}\times \C^l$ instead of  the homogeneous coordinates
 \eqref{e:homogeneous-coordinates}.
 As in  \eqref{e:w,zeta',t} we adopt  the  following notation for $n\in\N:$
 \begin{equation}\label{e:z,w,n}\begin{split} 
 O(t^n)dz\wedge d\bar z&:=\sum_{p,p'=1}^{k-l} O(t^n)dz_p\wedge d\bar z_{p'},\qquad O(t^n)dw\wedge d\bar w:=\sum_{q,q'=1}^{l}
  O(t^n)dw_q\wedge d\bar w_{q'},\\
 O(t^n)dz\wedge d\bar w &  :=\sum_{p=1}^{k-l}\sum_{q'=1}^l O(t^n)dz_p\wedge d\bar w_{q'},\qquad
   O(t^n)d\bar z\wedge d w   :=\sum_{p'=1}^{k-l}\sum_{q=1}^l O(t^n)d\bar z_p\wedge dw_{q}.
 \end{split}
 \end{equation}
Recall that  $\tau$ is  strongly  admissible  and  write  $\tau=(s_1,\ldots, s_k)$ in the  local coordinates $y=(z,w).$
Note  that $s_j= \tau^*z_j$  for $1\leq j\leq k-l$  and  $s_j= \tau^*w_{j-k+l}$  for  $k-l<j\leq k.$ 
In complement to  the first  collection of estimates  obtained in \eqref{e:diff-dz}--\eqref{e:diff-dw}, we infer from  Definition \ref{D:admissible-maps}   the  following  second collection dealing  with the Levi form  of the components of $\tau$  and their complex-conjugates:
 \begin{equation}\label{e:ddc-z}
  \ddc(\tau^*z_j)=  O(\|z\|) dz\wedge d\bar  z +      O(\|z\|) dz\wedge d\bar  w
 +  O(\|z\|^2) d\bar z\wedge d w +   O(\|z\|^2) dw\wedge d\bar  w,
 \end{equation}
 \begin{equation}\label{e:ddc-bar-z}
   \ddc(\tau^*\bar z_j)= O(\|z\|) dz\wedge d\bar  z  +      O(\|z\|^2) dz\wedge d\bar  w
 +  O(\|z\|) d\bar z\wedge d w +   O(\|z\|^2) dw\wedge d\bar  w,
 \end{equation}
 \begin{equation}\label{e:ddc-w}
 \ddc(\tau^* w_m )=     O(1) dz\wedge d\bar  z +     O(1) dz\wedge d\bar  w
 +  O(\|z\|) d\bar z\wedge d w +   O(\|z\|) dw\wedge d\bar  w,
 \end{equation}
 \begin{equation}\label{e:ddc-bar-w}
 \ddc(\tau^* \bar w_m )=    O(1) dz\wedge d\bar  z +       O(\|z\|) dz\wedge d\bar  w
 +   O(1) d\bar z\wedge d w 
 +    O(\|z\|) dw\wedge d\bar  w .
 \end{equation} 
 
 To   $I,J\subset \{1,\ldots,k\}$ we associate  
   $I',J'\subset \{1,\ldots,k-l\}$ and $I'',J''\subset \{1,\ldots,l\}$
such that  $I:=I'\cup\{ I''+(k-l)\}$ and  $J:=J'\cup \{J''+(k-l)\}.$
Here,  $\{K+p\}:=\{j+p:\ j\in K\}$ for $K\subset \{1,\ldots,l\}$ and $0\leq p\leq k-l.$ 
We can  write
$$
dy_I\wedge d\bar y_J=  dz_{I'}\wedge d\bar z_{J'}\wedge dw_{I''}\wedge d\bar w_{J''}.
$$
\begin{proposition}\label{P:negligible-forms} Let $0\leq m\leq 2l.$
Let $\Phi$ be a $\Cc^2$-smooth  $m$-negligible  form on $\Tube(B,\bfr)\subset \E.$
Then 
\begin{enumerate}
 \item  The forms  $\Phi^\sharp$ and $\tilde\tau^*\Phi$  are  $m$-negligible.
 \item If $\Phi$ is  of bidegree $(q,q)$ then  the form $\tilde\tau^*(\Phi)-[\tilde\tau^* (\Phi)]^\sharp$ is  $(m-1)$-negligible.
 
 \item  The forms $\dbar\Phi$  and $\partial\Phi$ and hence  $d\Phi$ are  $(m+1)$-negligible.
 
 \item The  forms $\dbar [(\tilde \tau_\ell)^*\Phi)] -(\tilde\tau_\ell)^*[\dbar \Phi]$ 
 and  $\partial [(\tilde \tau_\ell)^*\Phi)] -(\tilde\tau_\ell)^*[\partial \Phi]$ are   $m$-negligible.
\item The  form $\ddc [(\tilde \tau_\ell)^*\Phi)] -(\tilde\tau_\ell)^*[\ddc \Phi]$  is   $(m+1)$-negligible.
\end{enumerate}
\end{proposition}
\proof By linearity it  suffices to show  the proposition for  the form 
\begin{equation}\label{e:Phi-model} \Phi:=f(y)dy_I\wedge d\bar y_J= f(z,w)dz_{I'}\wedge d\bar z_{J'}\wedge dw_{I''}\wedge d\bar w_{J''},
\end{equation}
where $f$ is  a  smooth function compactly supported in $\Tube(B,\bfr).$ 
\\
\noindent  {\bf  Proof of assertion (1).}
By  Notation  \ref{N:principal}, $\Phi^\sharp$ given in \eqref{e:Phi-model} is equal to either  $\Phi$ or $0.$ 
Since $\Phi$ is  $m$-negligible, so  is $\Phi^\sharp.$  
  Moreover, using  \eqref{e:diff-dz}-\eqref{e:diff-dw} and \eqref{e:ddc-z}--\eqref{e:ddc-bar-w} 
we  can check  by Definition \ref{D:negligible} that 
\begin{equation*} \tilde\tau^*\Phi= (f \circ \tilde\tau)(z,w)d(\tilde\tau^*z_{I'})\wedge d(\tilde\tau^* \bar z_{J'})\wedge d(\tilde\tau^* w_{I''})\wedge d(\tilde\tau^*\bar w_{J''})
\end{equation*}
is  $m$-negligible.

\noindent  {\bf  Proof of assertion (2).}
Since the form  $\Phi$ is  $m$-negligible, we can  write 
$f(z,w)$  as  the sum of  finite functions of the  form  $  z^\bfm  \bar z^\bfn  g(z,w)$  with $|\bfm|+|\bfn|\geq  \max(0,|I''|+|J''|-m).$
Here, $z^\bfm:=  z_1^{m_1}\ldots z_{k-l}^{m_{k-l}}$ and $\bar z^\bfn:=  \bar z_1^{n_1}\ldots \bar z_{l}^{n_l}$ for  $\bfm:=(m_1,\ldots,m_{k-l})\in \N^{k-l}$
and $\bfn:=(n_1,\ldots,n_l)\in \N^l.$  Assume  without loss of generality that
$f(z,w)= z^\bfm  \bar z^\bfn  g(z,w).$

By  Definition \ref{D:Strongly-admissible-maps},  $\tilde\tau^*z-z=O(\|z\|^2)$ and  $\tilde\tau^*g-g=O(\|z\|).$  Therefore, we get that
\begin{eqnarray*}
\tilde\tau^*f(z,w) -f(z,w)&=&(z_1^{m_1}+O(\|z\|^2))\ldots (z_{k-l}^{m_{k-l}}+O(\|z\|^2))  (\bar z_1^{n_1}+O(\|z\|^2))\ldots (\bar z_{l}^{n_{l}}+O(\|z\|^2))\\&\cdot & (g(z,w)+O(\|z\|)  -  z^\bfm  \bar z^\bfn  g(z,w)\\
&=&  O(\|z\|^{|\bfm|+|\bfn|+1})=  O(\|z\|^{\max(0,|I''|+|J''|-(m-1))}).
\end{eqnarray*}
On  the  other hand,   we deduce  from  \eqref{e:diff-dz}-\eqref{e:diff-dw} that
the coefficients of  $d\bar z,$ $d\bar  w$  in   $d(\tilde\tau^* z_j)-dz_j$  and  in   $d(\tilde\tau^* w_m)-dw_m$ 
as well  as  the coefficients of  $d z,$ $ dw$  in   $d(\tilde\tau^* \bar z_j)-d\bar z_j$  and  in   $d(\tilde\tau^* \bar w_m)-d\bar w_m$  are of order $O(\|z\|).$ Using  this and  applying Lemma \ref{L:difference} (2), the result  follows.

\noindent  {\bf  Proof of assertion (3).}  We deduce  from the hypothesis $f(z,w)=O(\|z\|^{\max(0,|I''|+|J''|-m)})$ and the  equality
$$
\dbar\Phi=\dbar f \wedge dz_{I'}\wedge d\bar z_{J'}\wedge dw_{I''}\wedge d\bar w_{J''}\quad\text{and}\quad \partial\Phi=\partial f \wedge dz_{I'}\wedge d\bar z_{J'}\wedge dw_{I''}\wedge d\bar w_{J''}
$$
that both  forms $\dbar\Phi$  and $\partial\Phi$ are  $(m+1)$-negligible.

\noindent  {\bf  Proof of assertion (4).}
Consider the form $S=S_{I,J}dy_I\wedge d\bar y_J=\Phi.$
So $S_{I,J}=f$ and $dy_I\wedge d\bar y_J=dz_{I'}\wedge d\bar z_{J'}\wedge dw_{I''}\wedge d\bar w_{J''}.$
Applying the first   equality of  Lemma  \ref{L:ddc-difference-form} to $S$ yields that
\begin{eqnarray*}
 \dbar [(\tilde\tau_\ell)^*\Phi] -(\tilde\tau_\ell)^*(\dbar \Phi)&=&\big(\dbar [(\tilde\tau_\ell)^* f]-(\tilde\tau_\ell)^*[ \dbar f]\big) \bigwedge _{\iota\in I}d[(\tilde\tau_\ell)^*y_\iota]\wedge \bigwedge_{j\in J} d[(\tilde\tau_\ell)^* \bar y_j]\\
 &+&(\tilde\tau_\ell)^* (f) \wedge \bigwedge_{j\in J} d[(\tilde\tau_\ell)^* \bar y_j]\wedge\big( \sum_{\iota\in I}\pm   \ddc[(\tilde\tau_\ell)^*y_\iota]\wedge
 \bigwedge _{\iota'\in I\setminus \{\iota\}}d[(\tilde\tau_\ell)^*y_{\iota'}]
 \big)\\
 &+&(\tilde\tau_\ell)^* (f) \wedge \bigwedge _{\iota\in I}d[(\tilde\tau_\ell)^*y_\iota]\wedge\big(\sum_{j\in J}  \ddc[(\tilde\tau_\ell)^* \bar y_j]\wedge \bigwedge_{j'\in J\setminus \{j\}} d[(\tilde\tau_\ell)^* \bar y_j]\big).
 \end{eqnarray*}
Applying  Lemma \ref{L:ddc-difference-function}  to $f$ and  using  \eqref{e:diff-dz}-\eqref{e:diff-dw} and \eqref{e:ddc-z}--\eqref{e:ddc-bar-z},
\eqref{e:ddc-w}--\eqref{e:ddc-bar-w},
we see that the first term on the  RHS is $m$-negligible.
Using  \eqref{e:ddc-z}  and \eqref{e:ddc-w}, we see that the second sum  on the RHS is $m$-negligible.
Using  \eqref{e:ddc-bar-z}  and \eqref{e:ddc-bar-w}, we see that the third sum on the RHS  is $m$-negligible.
This  proves  the first  part of assertion (4).

Since $ d [(\tilde\tau_\ell)^*\Phi] -(\tilde\tau_\ell)^*(d\Phi)$ and $d=\partial+\dbar,$  it follows that
$$ \partial [(\tilde\tau_\ell)^*\Phi] -(\tilde\tau_\ell)^*(\partial \Phi) =  \dbar [(\tilde\tau_\ell)^*\Phi] -(\tilde\tau_\ell)^*(\dbar \Phi).$$
Hence, the second part of assertion (4) is a consequence of the first one.

\noindent  {\bf  Proof of assertion (5).} 
Consider the form $S=S_{I,J}dy_I\wedge d\bar y_J=\Phi.$
So $S_{I,J}=f$ and $dy_I\wedge d\bar y_J=dz_{I'}\wedge d\bar z_{J'}\wedge dw_{I''}\wedge d\bar w_{J''}.$
Applying the second  equality of  Lemma  \ref{L:ddc-difference-form} to $S$ yields that
\begin{eqnarray*}
 \ddc [(\tilde\tau_\ell)^*\Phi] -(\tilde\tau_\ell)^*(\ddc \Phi)&=&\big(\ddc [(\tilde\tau_\ell)^* f]-(\tilde\tau_\ell)^*[ \ddc f]\big) \bigwedge _{\iota\in I}d[(\tilde\tau_\ell)^*y_\iota]\wedge \bigwedge_{j\in J} d[(\tilde\tau_\ell)^* \bar y_j]\\
 &+&(\tilde\tau_\ell)^* (df) \wedge \bigwedge_{j\in J} d[(\tilde\tau_\ell)^* \bar y_j]\wedge\big( \sum_{\iota\in I}\pm   \ddc[(\tilde\tau_\ell)^*y_\iota]\wedge
 \bigwedge _{\iota'\in I\setminus \{\iota\}}d[(\tilde\tau_\ell)^*y_{\iota'}]
 \big)\\
 &+&(\tilde\tau_\ell)^* (df) \wedge \bigwedge _{\iota\in I}d[(\tilde\tau_\ell)^*y_\iota]\wedge\big(\sum_{j\in J}  \ddc[(\tilde\tau_\ell)^* \bar y_j]\wedge \bigwedge_{j'\in J\setminus \{j\}} d[(\tilde\tau_\ell)^* \bar y_j]\big).
 \end{eqnarray*}
Applying  Lemma \ref{L:ddc-difference-function}  to $f$ and  using  \eqref{e:diff-dz}-\eqref{e:diff-dw} and \eqref{e:ddc-z}--\eqref{e:ddc-bar-z},
\eqref{e:ddc-w}--\eqref{e:ddc-bar-w},
we see that the first term on the  RHS is $(m+1)$-negligible.
Using  \eqref{e:ddc-z}  and \eqref{e:ddc-w}, we see that the second sum  on the RHS is $(m+1)$-negligible.
Using  \eqref{e:ddc-bar-z}  and \eqref{e:ddc-bar-w}, we see that the third sum on the RHS  is $(m+1)$-negligible.
Hence, the result follows.
\endproof

\begin{corollary}\label{C:negligible-forms}  We keep  the  hypothesis  and notation of 
Proposition \ref{P:Stokes} and assume in addition that  the smooth test $(q,q)$-form $\Phi$ is $m$-negligible for some $0\leq m \leq 2l.$
Then  
\begin{enumerate}
 \item
The  forms 
 $\tilde\tau^*(\ddc\Phi) -\ddc(\tilde\tau^*\Phi)$ as well as 
 $d[(\tilde\tau^*\Phi)^\sharp]  ,$  $ \tilde\tau^*(\dc\Phi)$ are $(m+1)$-negligible.
 
 \item  The forms $\tilde\tau^*\Phi,$ $\dc (\tilde\tau^*\Phi)^\sharp -\tilde\tau^*(\dc\Phi), $
 $\tilde\tau^*(d\Phi) -d[(\tilde\tau^*\Phi)^\sharp] $  are $m$-negligible.
 
 \item   The form $\tilde\tau^*\Phi-(\tilde\tau^*\Phi)^\sharp$  is $(m-1)$-negligible.
 \end{enumerate}
\end{corollary}
\proof  It  follows from  Proposition \ref{P:negligible-forms}. 
\endproof

Recall here  Definitions \ref{D:diff-oper} and \ref{D:diff-oper-bis}.

\begin{proposition}\label{P:negligible-currents} 
Let $S$ be a  $(p,p)$-current on $\Tube (B,\bfr).$
Then 
\begin{enumerate}
 \item The operator $S\mapsto \dbar S $ is in the class $\widehat \Dc^1.$ 
 \item The operator $S\mapsto \tilde\tau^*[ (\tilde\tau_* S)^\sharp] $ is in the class $\widehat \Dc^0.$  
 and the  operator
 $S\mapsto \tilde\tau^*[(\tilde\tau_*S)^\sharp]-S $   is in the class $ \Dc^0.$

 \item The operator $S\mapsto \tilde\tau^*[ \dbar (\tilde\tau_* S)^\sharp ] -\dbar  S$ is in the class $ \Dc^1.$ 
\end{enumerate}

\end{proposition}
\proof[Proof of assertion (1) of Proposition \ref{P:negligible-currents}] Consider two multi-indices $M=(I,J;K,L)$ and $M'=(I',J';K',L')$ with
$|M|=2p+1$ and $|M'|=2p.$  In the  representation \eqref{e:diff-oper-1} we can write
\begin{equation*}
    (\dbar S)_M=\theta_\ell\big(\sum_{M',j}    g_{M,M',j} {\partial S_{M'}\over \partial \bar y_j}\big),
   \end{equation*}
   the sum being taken over all $M,M',j$ such that  $|I'|+|K'|=|J'|+|L'|=p$ and  $j\not \in K'\cup (k-l+L')$ and $|M|=|M'|+1=2p+1.$ Moreover,    
 $g_{M,M',j}=\pm 1$ if   $M=M'\cup \bar j$  and  $g_{M,M',j}=0$ otherwise.  Here,
 \begin{equation*}  M'\cup \bar j:=\begin{cases}  (I, J\cup \{j\};K,L),&  \text{if}\quad 1\leq j\leq k-l;\\
                                    (I, J;K,L\cup \{j-(k-l)\}),&  \text{if}\quad k-l+1\leq j\leq k.
                                   \end{cases}
    \end{equation*}
 So   when  $g_{M,M',j}\not=0,$ as $M=M'\cup \bar j$ we choose simply $P:=M',$  and hence $\Delta(P,M')=\Delta(M',M')=0.$   By  Definition  \ref{D:substraction-bis}, we can check that $\delta_{j,P,M}=0.$
 Hence,   $S\mapsto \dbar S $ is in the class $\widehat \Dc^1.$
\endproof

  Prior to the proof of the remaining assertions of Proposition \ref{P:negligible-currents},  we state  the  following (First) transfer rule for 
  the pull-back and pushforward $\tilde\tau^*$ and $\tilde\tau_*.$ Indeed,
  according to   \eqref{e:diff-dz} and \eqref{e:diff-dw}  we obtain the following  table:
  
  \begin{tabular}{|l|l|l|l|l|}
  \hline
  \multicolumn{5}{|c|}{ First transfer rule for $\tilde\tau^*$ (the same rule also holds for $\tilde\tau_*$)} \\
  \hline
  Source  & Target   && Source & Target  \\ 
  \hline
  \multirow{2}{*}{$d(\tilde\tau^*z_j)-dz_j$} &$ O(\|z\|) dz_p$ && {$d(\tilde\tau^*\bar z_j)-d\bar z_j$}     &  $O(\|z\|) d\bar z_p$    \\
     & $ O(\|z\|^2) \{d\bar z_p,dw_q,d\bar w_q\} $  &&            &   $ O(\|z\|^2) \{d z_p,dw_q,d\bar w_q\} $                   \\
    \hline
     \multirow{2}{*}{$d(\tilde\tau^*w_m)-dw_m$}  &$ O(1) dz_p$ && {$d(\tilde\tau^*\bar w_m)-d\bar w_m$}     &  $O(1) d\bar z_p$    \\
    & $ O(\|z\|) \{d\bar z_p,dw_q,d\bar w_q\} $  &&             &   $ O(\|z\|) \{d z_p,dw_q,d\bar w_q\} $                   \\
    \hline
    \end{tabular}
    
    \medskip
 We interpret  the table as follows.  The  term in each  source  column  is  replaced by one  of the terms proposed in the corresponding  target  column.
 Here  $\{\text{item} \ 1,   \text{item} \ 2, \text{item} \ 3   \}$ in the second line of each  target  case  means that we can choose  one the three proposed items. 
 
 The following  result is needed.
 \begin{lemma}\label{L:tau-negigible}
   Let
 $$S=\sum_M  S_M(y) dy_M=  \sum_MS_M(z,w)dz_I\wedge d\bar z_J\wedge dw_K\wedge d\bar w_L$$
 be a  $p$-current on $\Tube(B,\bfr),$ where the  $S_M$ are  distributions and  the sum is taken over all $M$  with 
 $|M|=2p.$
 Then
 \begin{enumerate}
  \item $(\tilde\tau^*S)_M=\sum_{M'} f_{M,M'}(\tilde\tau^*S_{M'}),$ where the sum is taken over all
   $M'$  with 
 $|M'|=2p$  and $f_{M,M'}$ is a smooth function with  $f_{M,M'}(z,w)=O(\|z\|^{\Delta(M,M')})$ for $M'\not=M$ and 
  $f_{M,M}(z,w)=1 +O(\|z\|).$
  \item   $(\tilde\tau_*S)_M=\sum_{M'} g_{M,M'}(\tilde\tau_*S_{M'}),$  where the multi-indices $M'$ and the functions
  $g_{M,M'}$ have the same property as in assertion (1). 
 \end{enumerate}
 \end{lemma}
\proof
We only give the proof of assertion (1) since  the same proof also  works  for assertion making the obviously necessary changes.
Write 
\begin{multline*}
 \tilde\tau^*S-\sum_{M'=(I',J';K',L')}  (\tilde\tau^*S_{M'})(y) dy_{M'}=\sum_{M'} (\tilde\tau^*S)_{M'}(z,w)\big[d  (\tilde\tau^*z_{I'})\wedge d(\tilde\tau^*\bar z_{J'})\wedge d(\tilde\tau^*w_{K'})\wedge d(\tilde\tau^*\bar w_{L'})\\
 -dz_{I'}\wedge d\bar z_{J'}\wedge dw_{K'}\wedge d\bar w_{L'}\big],
\end{multline*}
where the sums are taken over all multi-indices $M'=(I',J';K',L')$ with $|M'|=2p.$
Applying Lemma \ref{L:difference} (2) to each term  on  brackets  yields that
\begin{eqnarray*}
&& \big[d  (\tilde\tau^*z_{I'})\wedge d(\tilde\tau^*\bar z_{J'})\wedge d(\tilde\tau^*w_{K'})\wedge d(\tilde\tau^*\bar w_{L'})-
 dz_{I'}\wedge d\bar z_{J'}\wedge dw_{K'}\wedge d\bar w_{L'}\big]\\&=&  \sum_{M''=(I'',J'';K'',L'')} (dz_{I'})_{I''}\wedge (d\bar z_{J'})_{J''}\wedge (dw_{K'})_{K''}\wedge (d\bar w_{L'})_{L''},
\end{eqnarray*}
where  the last sum  is taken over all $M''=(I'',J'';K'',L'')$  such that $I''\subset  I',$ $J''\subset J',$ $K''\subset K'$ and $L''\subset L'$ and there is  at least one
 nonempty set among four sets $I'',$ $J'',$ $K'',$ $L''.$ 
 Here,  
 \begin{eqnarray*}(dz_{I'})_{I''}&=&\big(\bigwedge_{j\in  I'\setminus I''}  dz_j\big)\wedge  \big(\bigwedge_{j\in  I''}  [d(\tau^*z_j)-dz_j]\big),\\
 (d\bar z_{J'})_{J''}&=&\big(\bigwedge_{j\in  J'\setminus J''}  d \bar z_j\big)\wedge  \big(\bigwedge_{j\in  J''}  [d(\tau^*\bar z_j)-d\bar z_j]\big),
\end{eqnarray*}
and we have similar  definitions for    $(dw_{K'})_{K''}$  and $ (d\bar w_{L'})_{L''},$ namely,
\begin{eqnarray*}
(dw_{K'})_{K''}&=&\big(\bigwedge_{m\in  K'\setminus K''}  dw_m\big)\wedge  \big(\bigwedge_{m\in  K''}  [d(\tau^*w_m)-dw_m]\big),\\
 (d\bar w_{L'})_{L''}&=&\big(\bigwedge_{m\in  L'\setminus L''}  d \bar w_m\big)\wedge  \big(\bigwedge_{m\in  L''}  [d(\tau^*\bar w_m)-d\bar w_m]\big).
 \end{eqnarray*}
Next, we replace each term  in brackets $[d(\tau^*z_j)-dz_j],$  $[ d(\tau^*\bar z_j)-d\bar z_j],$ $[d(\tau^*w_m)-dw_m],$  $[ d(\tau^*\bar w_m)-d\bar w_m]$
by one  of  its  four possible corresponding target terms  in the above table, and  we  expand out  all possible  combinations.
Let $\exponent[\text{source term}]$  be the exponent of $\|z\|$  of   four possible corresponding target terms.
So each  $\exponent[\text{source term}]$ has two  possible values, it is a multi-valued function.  For each possible combination  we fix only one  among these two values, and hence
when either  fixing a possible combination or fixing  the  target item,  $\exponent$ becomes a univalued function. 
Consequently,  we  obtain
$$
(\tilde\tau^*S)=\sum_{M=(I,J;K,L)} (\tilde\tau^*S)_M dy_M=  \sum_{M=(I,J;K,L)} (\tilde\tau^*S)_M dy_M
$$
 where the sum is taken over all
   $M=(I,J;K,L)$  with 
 $|M|=p.$  Observe  that 
\begin{equation}\label{e:f_M,M'}(\tilde\tau^*S)_M=\sum_{M'} f_{M,M'}(\tilde\tau^*S_{M'}),
\end{equation}
 where the sum is taken over all
   $M'$  with 
 $|M'|=2p$  and the  functions $f_{M,M'}$'s are some smooth functions
 satisfying the  growth control $f_{M,M'}(z,w)=O(\|z\|^{\delta(M,M')}),$ where the $\delta(M,M')$'s  are some nonnegative integers.
 Observe  also that functions $f_{M,M'}$'s are uniquely determined  by the  relation
 \begin{equation}\label{e:f_M,M'-bis}
 \tilde\tau^*dy_{M'}=\sum_{M} f_{M,M'}dy_M,
\end{equation}
 To complete the proof of assertion (1) we need to show that  \begin{equation}\label{e:Delta-vs-delta}
 \delta(M,M')\geq \Delta(M,M')\quad\text{and}\quad f_{M,M}(z,w)=1+O(\|z\|) .
                                                               \end{equation}

 To this  end  fix multi-indices $M,$ $M'$ with $|M|=|M'|=2p.$  It follows  from the above expressions  that
 \begin{equation}\label{e:inequ-delta}
  \delta(M,M')\geq  \min_{M''=(I'',J'';K'',L'')}\delta(M,M',M''),
  \end{equation}
  the minimum being taken over all $M''=(I'',J'';K'',L'')$  such that $I''\subset  I',$ $J''\subset J',$ $K''\subset K'$ and $L''\subset L'$ and that there is  at least one
 nonempty set among four sets $I'',$ $J'',$ $K'',$ $L''.$ Here,
 \begin{equation} \label{e:delta-bis}
  \begin{split}
  \delta(M,M',M'')&:= \sum_{j\in I''} \exponent[ d(\tau^* z_j)-d z_j]
  +  \sum_{j\in J''} \exponent[ d(\tau^*\bar z_j)-d\bar z_j]\\
  &+ \sum_{m\in K''} \exponent[ d(\tau^* w_m)-d w_m]
  +  \sum_{m\in L''} \exponent[ d(\tau^*\bar w_m)-d\bar w_m].
 \end{split}\end{equation}
 Fix  such a  set $M''.$ In order to  show \eqref{e:Delta-vs-delta},  we only need to prove  that  
 \begin{equation}\label{e:Delta-vs-delta-bis} \delta(M,M',M'')\geq \Delta(M,M')\quad\text{and}\quad
 f_{M,M}(z,w)=1+O(\|z\|).
 \end{equation}
 If  $ |K|+|L| \geq |K'|+| L'|,$ then according to the above  table,
 there are at least  $(|K|+|L| - |K'|-| L'|)$ items
  $dw_m$ or $d\bar w_m$ such that each  of them  is the target  item of a source term  $d(\tau^* z_j)-d z_j$  for some $j\in I'\setminus I$
  or   the target  item of a source term  $d(\tau^* \bar z_j)-d \bar z_j$  for some $j\in J'\setminus J.$
  Since  in this case  $\exponent[ d(\tau^* z_j)-d z_j]\geq 2$ and  $\exponent[ d(\tau^* \bar z_j)-d \bar z_j]\geq 2$, it follows that
  \begin{equation*}
   \delta(M,M',M'')\geq  \sum_{j\in I''} \exponent[ d(\tau^* z_j)-d z_j]
  +  \sum_{j\in J''} \exponent[ d(\tau^*\bar z_j)-d\bar z_j]\geq  2(|K|+|L| - |K'|-| L'|).
  \end{equation*}
  Next, we divide $\delta(M,M',M'')$ into two disjoint parts:
  \begin{equation*}
   \delta(M,M',M'')=\delta_z(M,M',M'')+\delta_w(M,M',M''),
  \end{equation*}
where $\delta_z(M,M',M'')$ (resp. $\delta_w(M,M',M'')$) is  the sum of all exponents 
  in \eqref{e:delta-bis}  such that the target  items are either  $dz_j$ or $d\bar z_j$ (resp.  the target  items are either  $dw_m$ or $d\bar w_m$).
  
  According to the above table, we obtain the following:
  
  \noindent { \bf Fact.}{
  The only case when  
the $\exponent (\text{source item}, \text{target item})= 0$ is   either  (the source item is $d(\tilde\tau^*w_m)-dw_m$ and the target item is $O(1)dz_p$) or
the complex-conjugate situation, thats is,                (the source item is $d(\tilde\tau^*\bar w_m)-d\bar w_m$ and the target item is $O(1)d\bar z_p$).
  }
  
If  $m\in  K\setminus K'$ (resp. $m\in L\setminus L'$), then  using the above fact, $dw_m$ (resp. $d\bar w_m$)  should be 
the target  item  of a  source item  whose corresponding  $\exponent (\text{source}) \geq  1.$ 
Hence, 
\begin{equation*}
 \delta_w(M,M',M'')\geq  |K\setminus K'|+|L\setminus L'|.
\end{equation*}
Analogously, if  $j\in  I\setminus I'$ (resp. $j\in J\setminus J'$), then using the above fact, $dz_j$ (resp. $d\bar z_j$)  should be  
the target  item of a  source item  whose corresponding  $\exponent (\text{source}) \geq  1.$ Hence, 
\begin{equation*}
 \delta_z(M,M',M'')\geq  |I\setminus I'|+|K\setminus K'|.
\end{equation*}
  Putting together the last three estimates and the last equality  on exponents and using the formula of $\Delta(M,M')$ given in  \eqref{e:Delta}, the first inequality of \eqref{e:Delta-vs-delta} 
  follows.
  
   To  complete the proof of assertion (1), it remains  to show that   $f_{M,M}(z,w)=1+O(\|z\|).$ Fix a combination and we need to show that
   the exponent  of $\|z\|$ in  $f_{M,M}(z,w)-1$ is $\geq 1,$ that is, $\delta(M,M,M'')\geq 1.$
   Suppose   in order to reach a contradiction that  $\delta(M,M,M'')=0 .$ 
    We deduce    \eqref{e:delta-bis} that all terms on the RHS  are equal to $0.$
   Using the above fact,  we deduce  that $I''=\varnothing,$ $J''=\varnothing.$  So  either $K''\not=\varnothing$ or $L''\not=\varnothing.$
   Suppose without loss of generality that $K''\not=\varnothing.$   Using  the above fact and the fact that $\exponent[ d(\tau^* w_m)-d w_m]=0$
   for all $m\in K''$ and $\exponent[ d(\tau^* \bar w_m)-d \bar w_m]=0$ for all $m\in L'',$ we  deduce that $K=K'\setminus K''$ and $I=I'\cup K''.$
   This contradicts the assumption $M=M'.$

  The proof of assertion (1) of the lemma is thereby completed.

  
 
 \endproof
 
 Now  we arrive at 
 \proof[Proof of assertion (2) of Proposition  \ref{P:negligible-currents}]
Consider the  $2p$-current $R:= (\tilde\tau_* S)- (\tilde\tau_* S)^\sharp.$ Write 
 $R=\sum_{M=(I,J;K,L):\ |M|=2p} R_Md y_M,$
 where the $R_M$'s  are distributions.
 Applying   Lemma \ref{L:tau-negigible} (2) yields that 
  $R_M=0$ unless $|I|+|K|=|J|+|L|=p$ and 
 $R_M=\sum_{M'} g_{M,M'}(\tilde\tau_*S_{M'}),$ where the sum is taken over all
   $M'$  with 
 $|M'|=2p$  and $g_{M,M'}$ is a smooth function with  $g_{M,M'}(z,w)=O(\|z\|^{\Delta(M,M')})$ for $M'\not=M$ and 
  $g_{M,M}(z,w)=1 +O(\|z\|).$
  Next, applying   Lemma \ref{L:tau-negigible} (1) to $R$ yields that 
  \begin{equation*}
   (\tilde\tau^*R)_N=  \sum_{M'}  h_{N,M'} S_{M'},\quad\text{and}\quad   h_{N,M'}:=\sum_{M} f_{N,M}(z,w)g_{M,M'}(z,w) ,
  \end{equation*}
  for each multi-index $N$ with $|N|=2p$ and  the first (resp. second) sum is taken over all
   $M'$ (resp.    $M=(I,J;K,L)$)  with 
 $|M'|=2p$  (resp.  $|I|+|K|=|J|+|L|=p$). 
  Since 
  $$\tilde\tau^*R:=\tilde\tau^* [(\tilde\tau_* S)- (\tilde\tau_* S)^\sharp]=S- \tilde\tau^* [(\tilde\tau_* S)^\sharp], $$
  and  by Lemma  \ref{L:comparison} (2), $\Delta(N,M)+\Delta(M,M')\geq \Delta(N,M'),$  it follows  that
  $f_{N,M}(z,w)g_{M,M'}(z,w)=O(\|z\|^{\Delta(N,M')}),$ and hence  $h_{N,M'}(z,w)=O(\|z\|^{\Delta(N,M')}).$
  
  Consider the  case $N=M'.$  By  Lemma \ref{L:tau-negigible} (1) and (2) we  see that either  $M=N=M'$  and $f_{N,M}(z,w)=1+O(\|z\| $ and  $g_{M,M'}(z,w)=1+O(\|z\|),$
  or  $M\not=N$ and  $f_{N,M}(z,w)=O(\|z\| )$ and  $g_{M,M'}(z,w)=O(\|z\|),$
  Consequently,  $h_{N,N}(z,w)=1+O(\|z\|).$  
  
\endproof

Prior to the  proof of assertion (3) of Proposition \ref{P:negligible-currents}, the following intermediate results are   needed.

\begin{lemma}\label{L:assert(3)-P:negligible-currents} 
Consider the  $(2p+1)$-current $R:=\dbar [(\tilde\tau_*S)^\sharp] -(\tilde\tau)_*(\dbar S).$ 
Writing  
 \begin{equation*}
R=\sum_{M=(I,J;K,L):\ |M|=2p+1}  R_Mdy_M,
\end{equation*}
where $R_M$'s  are distributions. Then  the following representation holds
\begin{equation}
\label{e:diff-oper-1-bisbis}
    R_M=\theta_\ell \sum_{M',j,n}  \big(f_{M,M',j,n} {\partial (\tilde\tau_*S_{M'})\over \partial \tilde y_j} +  g_{M,M',j,n} {\partial (\tilde\tau_*S_{M'})\over \partial \overline{ \tilde y}_j}
    +h_{M,M',j,n}( \tilde\tau_*S_{M'})\big),
\end{equation}
  where  the sum is taken over all $M'$ with $|M'|=|M|-1=2p$ and $1\leq  j\leq k$ and $1\leq n\leq n_0,$ and $n_0$ is a positive integer.  
Here  $f_{M,M',j,n},$ $g_{M,M',j,n},$ $h_{M,M',j,n}$  are  smooth functions  such that 
 \begin{equation*}  \begin{split}  f_{M,M',j,n}(\tilde y )&=O(t^{\max(1,\Delta_j(M,M'))}),\quad g_{M,M',j,n}(\tilde y)=O(t^{\max(1,\Delta_j(M,M'))}),\\ h_{M,M',j,n}(\tilde y)&=O(t^{\max(0,\Delta_j(M,M')-1)})
                    \end{split}
                    \end{equation*}
 for all  $M,M',j.$  
 \end{lemma}
 \begin{remark}
  \rm It is  worthy to compare  the  conclusion of Lemma \ref{L:assert(3)-P:negligible-currents} with the class $\Dc^1$  given un  Definition \ref{D:diff-oper-bis}.
 \end{remark}

\proof
As  in the  proof of assertion (2) of Proposition \ref{P:negligible-currents}, consider the  $2p$-current $R:= (\tilde\tau_* S)- (\tilde\tau_* S)^\sharp.$ Writing
 $R=\sum_{M'=(I',J';K',L'):\ |M'|=2p} R_{M'}d y_{M'},$
 where the $R_{M'}$'s  are distributions,
 we know  that 
  $R_{M'}=0$ unless $|I'|+|K'|=|J'|+|L'|=p$ and 
 $R_{M'}=\sum_{M''} \lambda_{M',M''}(\tilde\tau_*S_{M''}),$ where the sum is taken over all
   $M''$  with 
 $|M''|=2p$  and $\lambda_{M',M''}$ is a smooth function with  $\lambda_{M',M''}(z,w)=O(\|z\|^{\Delta(M',M'')})$ for $M''\not=M'$ and 
  $\lambda_{M',M''}(z,w)=1 +O(\|z\|).$
  Observe that
  \begin{equation*}
   \dbar   (R_{M'}d y_{M'})= \sum_{M'',\ 1\leq j\leq k} {\partial  \lambda_{M',M''}\over \partial\overline{\tilde y}_j }  (\tilde\tau_*S_{M''})d \overline{\tilde y_j}\wedge d y_{M'}+ \sum_{M'',\ 1\leq j\leq k} \lambda_{M',M''}   {\partial (\tilde\tau_*S_{M''})
   \over \partial\overline{\tilde y}_j } d \overline{\tilde y}_j\wedge d y_{M'}.
  \end{equation*}
If $j\in J'\cup ( L'+(k-l))$ then $d\overline{\tilde y_j}\wedge d y_{M'}=0$ and there is  nothing to do.
Otherwise,  we  set $M:=M'\cup\{j'\}$ if $1\leq j'\leq k-l$ and $M:=M'\cup\{j'-(k-l)\}$ if $k-l< j'\leq k.$
Moreover,  we set $$f_{M,M',j,n}:=0,\quad g_{M,M',j,n}:=  \lambda_{M',M''}  \quad\text{and}\quad h_{M,M',j,n}:={\partial  \lambda_{M',M''}\over \partial\overline{\tilde y}_j } .$$
We can check  that  $\dbar R=  \sum_{M'}\dbar  (R_{M'}d y_{M'})$ has  the desired property  stated in \eqref{e:diff-oper-1-bisbis}.  
  Therefore, writing  
 \begin{equation}
  \dbar [(\tilde\tau_*S)^\sharp] -(\tilde\tau)_*(\dbar S)=[\dbar (\tilde\tau_*S) -(\tilde\tau)_*(\dbar S)]  -\dbar R,
 \end{equation}
it remains  to us  to show that $\dbar (\tilde\tau_*S) -(\tilde\tau)_*(\dbar S)$ also possesses  the desired property  stated in \eqref{e:diff-oper-1-bisbis}.

 Write
 $S=\sum_{M'=(I',J';K',L')} S_{M'}d y_{M}',$
 where the $R_M$'s  are distributions and  the sum is taken over all $M'=(I',J';K',L')$ with $|I'|+|K'|=|J'|+|L'|=p.$
By Lemma \ref{L:ddc-difference-form}, we can write
 \begin{equation}
 \dbar [(\tilde\tau)_*S] -(\tilde\tau)_*(\dbar S)=T_0+T_1,
 \end{equation}
 where 
 \begin{equation}\label{e:T_1-negligible}
 T_1:=\sum_{M'}\big(\dbar [\tilde\tau_* S_{M'}]-(\tilde\tau)_*[ \dbar S_{M'}]\big) \bigwedge _{j\in I'}d(\tilde\tau_*z_j)\wedge 
 \bigwedge _{j\in J'}d(\tilde\tau_*z_j)\wedge\bigwedge _{m\in K'}d(\tilde\tau_*w_m)\wedge\bigwedge _{m\in L'}d(\tilde\tau_*w_m),
 \end{equation}
 and  
  \begin{equation}\label{e:T_0-negligible}
T_0:=\sum_{M'=(I',J';K',L')} (\tilde\tau_* S_{M'}) R_{M'}.
\end{equation}
 Here    $R_{M'}$ is  the smooth form, which is  not necessarily of bidegree $(|I'|+|K'|,|J'|+|L'|),$   given by
 \begin{eqnarray*}
  R_{M'}&:=&  \big( \sum_{j\in I'}\pm   \ddc(\tilde\tau_*z_j)\wedge
 \bigwedge_{\iota\in I'\setminus \{j\}}d(\tilde\tau_*z_\iota)
 \big)   \bigwedge_{j\in J'} d(\tilde\tau_* \bar z_j)  
 \wedge  \bigwedge_{m\in K'} d(\tilde\tau_*  w_m)\wedge \bigwedge_{m\in L'} d(\tilde\tau_* \bar w_m)          \\
 &+ & \bigwedge_{j\in I'} d(\tilde\tau_*  z_j) \big( \sum_{j\in J'}\pm   \ddc(\tilde\tau_*\bar z_j)\wedge
 \bigwedge_{\iota\in J'\setminus \{j\}}d(\tilde\tau_*\bar z_\iota)
 \big)    
 \wedge  \bigwedge_{m\in K'} d(\tilde\tau_*  w_m) \wedge \bigwedge_{m\in L'} d(\tilde\tau_* \bar w_m)          \\
 &+&  \bigwedge_{j\in I'} d(\tilde\tau_*  z_j)  \wedge  \bigwedge_{j\in J'} d(\tilde\tau_*  \bar z_j)
 \big( \sum_{m\in K'}\pm   \ddc(\tilde\tau_* w_m)\wedge
 \bigwedge_{\iota\in K'\setminus \{m\}}d(\tilde\tau_* w_m)
 \big)    
  \wedge \bigwedge_{m\in L'} d(\tilde\tau_* \bar w_m)          \\
  &+&  \bigwedge_{j\in I'} d(\tilde\tau_*  z_j)  \wedge  \bigwedge_{j\in J'} d(\tilde\tau_*  \bar z_j)
  \wedge \bigwedge_{m\in K'} d(\tilde\tau_*  w_m)
 \big( \sum_{m\in L'}\pm   \ddc(\tilde\tau_* \bar w_m)\wedge
 \bigwedge_{\iota\in L'\setminus \{m\}}d(\tilde\tau_* \bar w_m)
 \big)    
  .
 \end{eqnarray*} 
 We  only need to show  that both $T_0$ and $T_1$ possess  the desired property  stated in \eqref{e:diff-oper-1-bisbis}.
This is  the content of   Lemmas \ref{L:T_1-negligible} and \ref{L:T_1-negligible} below. Modulo these lemmas, the proof is  thereby completed.
\endproof

\begin{lemma}\label{L:T_1-negligible}
Let $T_1$ be  the $(2p+1)$-current  given by \eqref{L:T_1-negligible} and  write $T_1=\sum_M  (T_1)_M dy_M,$ where
$(T_1)_M$  are  distributions and the sum is taken over all multi-indices $M$ with $|M|=2p+1.$
 Then  the following representation holds
\begin{equation}
\label{e:diff-oper-1-bisbis-T_1}
    (T_1)_M=\sum_{M',j,n}  \big(f_{M,M',j,n} {\partial (\tilde\tau_*S_{M'})\over \partial \tilde y_j} +  g_{M,M',j,n} {\partial (\tilde\tau_*S_{M'})\over \partial \overline{ \tilde y}_j}
\big),
\end{equation}
  where  the sum is taken over all $M'$ with $|M'|=|M|-1=2p$ and $1\leq  j\leq k$ and $1\leq n\leq n_0,$ and $n_0$ is a positive integer.  
Here  $f_{M,M',j,n}$   are  smooth functions  such that  
 $   f_{M,M',j,n}(\tilde y )=O(t^{\max(1,\Delta_j(M,M'))}) $
 for all  $M,M',j.$  
\end{lemma}
\proof
 We fix  $M'=(I',J';K',L')$ with $|I'|+|K'|=|J'|+|L'|=p$ and  prove this fact for each term in the sum \eqref{e:T_1-negligible}
 \begin{equation}\label{e:T_1-negligible-bis}
 \big(\dbar [\tilde\tau_* S_{M'}]-(\tilde\tau)_*[ \dbar S_{M'}]\big) \bigwedge _{j\in I'}d(\tilde\tau_*z_j)\wedge 
 \bigwedge _{j\in J'}d(\tilde\tau_*z_j)\wedge\bigwedge _{m\in K'}d(\tilde\tau_*w_m)\wedge\bigwedge _{m\in L'}d(\tilde\tau_*w_m).
 \end{equation}
 Applying  Lemma \ref{L:ddc-difference-function} yields that
 \begin{equation}\label{e:dbar-tau-commute}\dbar [\tilde\tau_*S_{M'}] -(\tilde\tau)_*[\dbar S_{M'}]= \sum_{j=1}^{k-l}  {\partial (\tilde\tau_*S) \over \partial  z_j}\dbar (\tilde\tau_*   z_j) 
 + \sum_{m=1}^{l}  {\partial (\tilde\tau_*S )\over \partial  w_m}\dbar (\tilde\tau_*   w_m) .
\end{equation}
By Definition \ref{D:Strongly-admissible-maps}, we know that for $1\leq j\leq k-l$ and $1\leq m\leq l,$
\begin{equation}\label{e:dbar-tau}
 \dbar (\tilde\tau_*   z_j) =O(\|z\|^2)\quad\text{and}\quad  \dbar (\tilde\tau_*   w_m) =O(\|z\|).
\end{equation}
Using  \eqref{e:dbar-tau-commute} and  \eqref{e:dbar-tau}, we can check  that
\begin{equation*}
 \big(\dbar [\tilde\tau_* S_{M'}]-(\tilde\tau)_*[ \dbar S_{M'}]\big) \bigwedge _{j\in I'}dz_j\wedge 
 \bigwedge _{j\in J'}d\bar z_j\wedge\bigwedge _{m\in K'}dw_m\wedge\bigwedge _{m\in L'}d\bar w_m
\end{equation*}
possess  the desired property  stated in \eqref{e:diff-oper-1-bisbis-T_1}. Hence, in order to prove
this  property  for the current given in  \eqref{e:T_1-negligible-bis}, we are reduced to proving this for the current
\begin{equation}\label{e:T_1-negligible-bisbis}
\begin{split}
 T_{1,M'}&:=\big(\dbar [\tilde\tau_* S_{M'}]-(\tilde\tau)_*[ \dbar S_{M'}]\big) \wedge \big[\bigwedge _{j\in I'}d(\tilde\tau_*z_j) \wedge 
 \bigwedge _{j\in J'}d(\tilde\tau_*z_j)\wedge\bigwedge _{m\in K'}d(\tilde\tau_*w_m)\wedge\bigwedge _{m\in L'}d(\tilde\tau_*w_m)\\
 &-  \bigwedge _{j\in I'}dz_j\wedge 
 \bigwedge _{j\in J'}d\bar z_j\wedge\bigwedge _{m\in K'}dw_m\wedge\bigwedge _{m\in L'}d\bar w_m \big].
 \end{split}
 \end{equation}
 Arguing as  in the proof of \eqref{e:f_M,M'}--\eqref{e:f_M,M'-bis} and \eqref{e:Delta-vs-delta}, we  infer that
 \begin{equation*}T_{1,M'}=\big(\dbar [\tilde\tau_* S_{M'}]-(\tilde\tau)_*[ \dbar S_{M'}]\big)\wedge \big(\sum_{M} f_{M,M'}dy_M\big),
\end{equation*}
 where the sum is taken over all
   $M$  with 
 $|M|=2p$  and the  functions $f_{M,M'}$'s are some smooth functions
 satisfying the  growth control $f_{M,M'}(z,w)=O(\|z\|^{\Delta(M,M')})$ and   $f_{M',M'}(z,w)=O(\|z\|) .$

 For  $\mathcal M:=M\cup\{\bar j\}$ 
 set $P:=M.$ 
 Using the  last  growth control of $f_{M,M'}$ and  the last expression  of  $T_{1,M'},$ we can check that $T_{1,M'}$  possesses  the desired property  stated in \eqref{e:diff-oper-1-bisbis-T_1}
 with the multi-index $\mathcal M$ instead of $M.$
\endproof

According to   \eqref{e:ddc-z}--\eqref{e:ddc-bar-z} and \eqref{e:ddc-w}--\eqref{e:ddc-bar-w}  we obtain the following  table:
  
  \medskip
  
  \begin{tabular}{|l|l|l|l|l|}
  \hline
  \multicolumn{5}{|c|}{Second transfer rule for $\tilde\tau^*$ (the same rule also holds for $\tilde\tau_*$)} \\
  \hline
  Source  & Target   && Source & Target  \\ 
  \hline
  \multirow{2}{*}{$\ddc(\tilde\tau^*z_j)$} &$ [O(\|z\|) \{dz_p\}]\wedge \{d\bar z_{p'}, d\bar w_{q'}\}$ &&  {$\ddc(\tilde\tau^*\bar z_j)$} 
  &   $ [O(\|z\|) \{d\bar z_p\}]\wedge \{dz_{p'}, d w_{q'}\}$   \\
     & $ [O(\|z\|^2) \{d\bar z_p,d\bar w_q\}]\wedge \{dw_{q'}\} $  &&            &   $[ O(\|z\|^2) \{d z_p,d w_q\}]\wedge \{d\bar w_{q'}\} $ \\                 
    \hline
    \multirow{2}{*}{$\ddc(\tilde\tau^*w_m)$} &$[ O(1) \{dz_p\}]\wedge \{d\bar z_{p'}, d\bar w_{q'}\}$ &&{  $\ddc(\tilde\tau^*\bar w_m)$}     &   $ [O(1) 
    \{d\bar z_p\}]\wedge \{dz_{p'}, d w_{q'}\}$   \\
     & $ [O(\|z\|) \{d\bar z_p,d\bar w_q\}]\wedge \{dw_{q'}\} $  &&            &   $[ O(\|z\|) \{d z_p,d w_q\}]\wedge \{d\bar w_{q'}\} $                  
    \\
    \hline
    \end{tabular}
    
    \medskip

\begin{lemma}
 \label{L:T_0-negligible}
Let $T_0$ be  the $(2p+1)$-current  given by \eqref{L:T_0-negligible} and  write $T_0=\sum_M  (T_0)_M dy_M,$ where
$(T_0)_M$  are  distributions and the sum is taken over all multi-indices $M$ with $|M|=2p+1.$
 Then  the following representation holds
\begin{equation}
\label{e:diff-oper-1-bisbis-T_0}
    R_M=\sum_{M',j,n}   
    h_{M,M',j,n}( \tilde\tau_*S_{M'}) ,
\end{equation}
  where  the sum is taken over all $M'$ with $|M'|=|M|-1=2p$ and $1\leq  j\leq k$ and $1\leq n\leq n_0,$ and $n_0$ is a positive integer.  
Here    $h_{M,M',j,n}$  are  smooth functions  such that 
 $   h_{M,M',j,n}(\tilde y)=O(t^{\Delta_j(M,M')})$
 for all  $M,M',j.$  
\end{lemma}
\proof
We  interpret  the  above table  as  follows. The term in  each  source column is  replaced by one of the terms  proposed in the  corresponding  target column. Each term in a target column has the form $[\text{Factor 1}]\wedge \text{Factor 2},$
where  $\text{Factor 1}$ is of the form $$O(\|z\|^s) \ \text{times one of the item in the first }\ \{...\}$$ for some $0\leq s\leq 2,$
and   $\text{Factor 2}$ is one of the items in  the second $\{...\}.$ 

The key observation is  that
\begin{itemize}
 \item  each  Factor 1  of $\ddc(\tilde\tau^*z_j)$  (resp. of   $\ddc(\tilde\tau^*\bar z_j)$) is similar  to  the corresponding  target item
 of $d(\tilde\tau^*z_j)-dz_j$  (resp.   of $d(\tilde\tau^*\bar z_j)-d\bar z_j$);
 \item each  Factor 1  of $\ddc(\tilde\tau^*w_m)$  (resp. of   $\ddc(\tilde\tau^*\bar w_m)$) is similar  to  the corresponding  target item
 of $d(\tilde\tau^*w_m)-dw_m$  (resp.   of $d(\tilde\tau^*\bar w_m)-d\bar w_m$).
\end{itemize}
    Next,  we  argue as in the proof of   Lemma \ref{L:tau-negigible} but using the  pushforward $\tau_*$ instead of the pull-back $\tau^*$
    and  using the second transfer rule instead of the first  one and  using the above key observation. Consequently, we obtain
    the representation \eqref{e:diff-oper-1-bisbis-T_0} with $dy_j$ (or $d\bar y_j$) is  the $\text{Factor 2}$ of a target item and  $M=M'\cup\{j\}.$
\endproof

 We still need  some  properties of   the   positive  substraction.
 \begin{lemma}\label{L:comparison-bis}
 \begin{enumerate}
 \item  Let $M=(I_M,J_M;K_M,L_M)$ and $N=(I_N,J_N;K_N,L_N)$ be two multi-indices with $|M|=|N|$
 and  let $j'\in N.$ 
 Set $\widetilde N:=N\setminus \{j'\}.$
 If  $j'\in  M$ then we set $\widetilde M:=M\setminus \{j'\}.$
Otherwise  if $K_M\cup L_M\not=\varnothing$ then   then we set $\widetilde M:=M\setminus \{j''\}$
 for some $j''\in K_M\cup L_M.$
 Otherwise,  if $K_M= L_M=\varnothing$ then 
 we choose  an arbitrary  subset $\widetilde M\subset M$ of length $|M|-1.$
 Then 
 $ \Delta(\widetilde N,\widetilde M)\leq \Delta(N,M).$
 Moreover, in the last two cases (that is, if $j'\not\in N$), we have 
 $ \Delta(\widetilde N,\widetilde M)\leq \Delta(N,M)-1.$
 
 \item For  $j\in \{1,\ldots, k\}\cup  \{\bar 1,\ldots, \bar k\}$ and for multi-indices $M,$ $N,$ $Q$ with
 $|M|=|N|=|Q|+1,$ the following  triangle inequality holds:
 $\Delta_j(M,Q)\leq \Delta(N,M)+\Delta_j(N,Q).$ 
 \end{enumerate}
 \end{lemma}
\proof
We prove  assertion (1)
considering  each of the  third  proposed  cases
by using  Definition  \ref{D:substraction}.

In the  first case (that is,   $j'\in  N$), 
we see  that
\begin{eqnarray*}
 |I_N\setminus I_M|+|J_N\setminus J_M|+|K_N\setminus K_M|+|L_N\setminus L_M|&=&  |I_{\widetilde N}\setminus I_{\widetilde M}|+|J_{\widetilde N} \setminus J_{\widetilde M}|+|K_{\widetilde N}\setminus K_{\widetilde M}|+|L_{\widetilde N}\setminus L_{\widetilde M}|,\\
 |K_N|+|L_N| -|K_M|-|L_M|&=&  |K_{\widetilde N}|+|L_{\widetilde N}| -|K_{\widetilde M}|-|L_{\widetilde M}|.
\end{eqnarray*}
So by  Definition  \ref{D:substraction},
$ \Delta(\widetilde N,\widetilde M)= \Delta(N,M).$

In the second case (that is,   $j'\not\in  N$ and  $K_N\cup L_N\not=\varnothing$), 
we see  that
\begin{equation*}
1\geq \big( |I_N\setminus I_M|+|J_N\setminus J_M|+|K_N\setminus K_M|+|L_N\setminus L_M|\big)-\big(  |I_{\widetilde N}\setminus I_{\widetilde M}|+|J_{\widetilde N} \setminus J_{\widetilde M}|+|K_{\widetilde N}\setminus K_{\widetilde M}|+|L_{\widetilde N}\setminus L_{\widetilde M}|\big) \geq 0.
\end{equation*}
On the other hand, we also  have
\begin{equation*}
 |K_N|+|L_N| -|K_M|-|L_M|\geq  |K_{\widetilde N}|+|L_{\widetilde N}| -|K_{\widetilde M}|-|L_{\widetilde M}|.
\end{equation*}
This,  coupled with the previous  estimate, implies assertion (1) in  the second case.

In the last case (that is,   $j'\not\in  N$ and  $K_N= L_N=\varnothing$),  we have $K_{\widetilde N}= L_{\widetilde N}=\varnothing,$ and hence
we see  that
\begin{equation*}
1\geq \big( |I_N\setminus I_M|+|J_N\setminus J_M|\big)-\big(  |I_{\widetilde N}\setminus I_{\widetilde M}|+|J_{\widetilde N} \setminus J_{\widetilde M}|\big) \geq 0.
\end{equation*}
On the other hand, since  $-|K_M|-|L_M|\leq 0$ and $ -|K_{\widetilde M}|-|L_{\widetilde M}|,$
we deduce from Definition \ref{D:substraction} that 
\begin{equation*}
  \Delta(N,M)= |I_N\setminus I_M|+|J_N\setminus J_M|\qquad\text{and}\qquad \Delta(\widetilde N,\widetilde M)=  |I_{\widetilde N}\setminus I_{\widetilde M}|+|J_{\widetilde N} \setminus J_{\widetilde M}|.
\end{equation*}
This,  coupled with the previous  estimate, implies assertion (1) in  the last case.

We turn to  the proof of  assertion (2). By Definition  \ref{D:substraction-bis}, let $\widetilde N\subset N$ such that $|\widetilde N|=|N|-1$ and  that 
$$
 \Delta_j(N,Q)= \delta_{j,\widetilde N,N}+ \Delta(\widetilde N,Q)\quad\text{and}\quad \delta_{j,\widetilde N,N}\in\{0,1\}.
 $$
 Set $j':=N\setminus  \widetilde N$ and  define  $\widetilde M$ as in assertion (1).
By assertion (1) we have   that $\Delta(\widetilde N,\widetilde M)\leq \Delta(N,M).$ 
On the  other hand,  by  
 Lemma  \ref{L:comparison} we get $\Delta(\widetilde N,\widetilde M)+\Delta(\widetilde N,Q)\geq \Delta(\widetilde M,Q).$ So
 $$ \Delta(\widetilde N,Q)+\Delta(N,M)  \geq  \Delta(\widetilde N,\widetilde M)+\Delta(\widetilde N,Q)\geq \Delta(\widetilde M,Q) .$$
 By Definition \ref{D:substraction-bis}, we have  
 $$
 \Delta_j(M,Q)\leq \delta_{j,\widetilde M,M}+ \Delta(\widetilde M,Q)\quad\text{and}\quad \delta_{j,\widetilde M,M}\in\{0,1\}.
 $$
 If $\delta_{j,\widetilde N,N}\geq \delta_{j,\widetilde M,M},$ then putting  the last  three estimates 
 together assertion (2) follows.

 It remains  to consider the case where  $\delta_{j,\widetilde M,M}=1$  and $\delta_{j,\widetilde N,N}=0.$
 If  we were in the  first case of assertion (1), then $M\setminus  \widetilde M=N\setminus \widetilde N,   $ and hence $\delta_{j,\widetilde M,M}$ would be equal to  $\delta_{j,\widetilde N,N}.$ Since this is  not the case, we are in the context of the  last two cases of assertion (1).
 So by this assertion, we have   that $\Delta(\widetilde N,\widetilde M)\leq \Delta(N,M)-1.$ 
On the  other hand,  by  
 Lemma  \ref{L:comparison} we get $\Delta(\widetilde N,\widetilde M)+\Delta(\widetilde N,Q)\geq \Delta(\widetilde M,Q).$ So
 $$ \Delta(\widetilde N,Q)+\Delta(N,M)  \geq  \Delta(\widetilde N,\widetilde M)+\Delta(\widetilde N,Q)+1\geq \Delta(\widetilde M,Q) +1.$$
 This implies  assertion (2).
\endproof

 Now  we arrive at 
 \proof[Proof of assertion (3) of Proposition  \ref{P:negligible-currents}]
Consider  the $(2p+1)$-current  $R:= \dbar[(\tilde\tau_*S)^\sharp] -(\tilde\tau)_*(\dbar S).$ 
 Observe that
  $$  \tilde\tau^*R:=\tilde\tau^*\big [  \dbar [(\tilde\tau_*S)^\sharp] -(\tilde\tau)_*(\dbar S)    \big]=\tilde\tau^*[\dbar((\tilde\tau_*S)^\sharp)]
  -\dbar S .$$
 So we need  to show that the operator $S\mapsto  \tilde\tau^*R$ is  in the class $\Dc^1.$ 
Writing  
 \begin{equation*}
R=\sum_{M=(I,J;K,L):\ |M|=2p+1}  R_Mdy_M,
\end{equation*}
where $R_M$'s  are distributions. Then  by Lemma \ref{L:assert(3)-P:negligible-currents},  the following representation holds
\begin{equation}
\label{e:diff-oper-1-bisbisbis}
    R_M=\sum_{M',j,n}  \big(f_{M,M',j,n} {\partial (\tilde\tau_*S_{M'})\over \partial \tilde y_j} +  g_{M,M',j,n} {\partial (\tilde\tau_*S_{M'})\over \partial \overline{ \tilde y}_j}
    +h_{M,M',j,n}( \tilde\tau_*S_{M'})\big),
\end{equation}
  where  the sum is taken over all $M'$ with $|M'|=|M|-1=2p$ and $1\leq  j\leq k$ and $1\leq n\leq n_0,$ and $n_0$ is a positive integer.  
Here  $f_{M,M',j,n},$ $g_{M,M',j,n},$ $h_{M,M',j,n}$  are  smooth functions  such that 
 \begin{equation}\label{e:fgh_M,M',j,n}  \begin{split}  f_{M,M',j,n}(\tilde y )&=O(t^{\max(1,\Delta_j(M,M'))}),\quad g_{M,M',j,n}(\tilde y)=O(t^{\max(1,\Delta_j(M,M'))}),\\ h_{M,M',j,n}(\tilde y)&=O(t^{\max(0,\Delta_j(M,M')-1)})
                    \end{split}
                    \end{equation}
 for all  $M,M',j.$  
  Next, applying   Lemma \ref{L:tau-negigible} (1) to $R$ yields that 
   $$(\tilde\tau^*R)_N=\sum_{M} F_{N,M}(\tilde\tau^*R_{M}),$$
   where the sum is taken over all
   $M$  with 
 $|M|=2p+1$  and $F_{N,M}$ is a smooth function with  
 \begin{equation}\label{e:F_N,M}F_{N,M}(z,w)=O(\|z\|^{\Delta(N,M)})\quad\text{ for}\quad N\not=M\quad\text{and}\quad 
  F_{M,M}(z,w)=1 +O(\|z\|).
  \end{equation}
  On the  other hand, there are smooth functions $\alpha_{jq}(z,w),$ $\beta_{jq}(z,w),$  and $\gamma_{jq}(z,w),$ $\delta_{jq}(z,w)$  such that 
  \begin{eqnarray*}
   \tilde\tau^*\big({\partial (\tilde\tau_*S_{M'})\over \partial \tilde y_j}\big)&=&\sum_{q=1}^k\big(\alpha_{jq}(z,w) {\partial S_{M'}\over \partial \tilde y_q}
   +\beta_{jq}(z,w)  {\partial S_{M'}\over \partial \overline{ \tilde y}_q}\big),\\
   \tilde\tau^*\big({\partial (\tilde\tau_*S_{M'})\over \partial \overline{\tilde y}_j}\big)&=&\sum_{q=1}^k \big(\gamma_{jq}(z,w) {\partial S_{M'}\over \partial \tilde y_q}
   +\delta_{jq}(z,w)  {\partial S_{M'}\over \partial \overline{ \tilde y}_q}\big).
  \end{eqnarray*}
  This, combined with \eqref{e:diff-oper-1-bisbisbis}, yields that
  \begin{eqnarray*}
   (\tilde\tau^*R)_N&=&  \sum_{M',j,n}  \big(\tilde f_{N,M',j,n} {\partial S_{M'}\over \partial \tilde y_j} +  \tilde g_{N,M',j,n} {\partial S_{M'}\over \partial \overline{ \tilde y}_j}\big)\\
   &+&  \sum_{M',j,n}  \big(\tilde{\tilde f}_{N,M',j,n} {\partial S_{M'}\over \partial \tilde y_j} +  \tilde{\tilde g}_{N,M',j,n} {\partial S_{M'}\over \partial \overline{ \tilde y}_j}\big)
    +\sum_{M',j,n} \tilde h_{N,M',j,n}S_{M'}.
  \end{eqnarray*}
  Here
   we have 
  \begin{eqnarray*} 
  \tilde f_{N,M',j,n}&:=&\sum_{M} F_{N,M}(z,w) f_{M,M',j,n}    (z,w)\alpha_j(z,w),\quad \tilde g_{N,M',j,n}:=\sum_{M} F_{N,M}(z,w) f_{M,M',j,n}    (z,w)\beta_j(z,w),\\
  \tilde{\tilde f}_{N,M',j,n}&:=&\sum_{M} F_{N,M}(z,w) g_{M,M',j,n}    (z,w)\gamma_j(z,w),\quad \tilde{\tilde g}_{N,M',j,n}:=\sum_{M} F_{N,M}(z,w) g_{M,M',j,n}    (z,w)\delta_j(z,w),\\
  \tilde h_{N,M',j,n}&:=&\sum_{M} F_{N,M}(z,w) h_{M,M',j,n}    (z,w),
  \end{eqnarray*}
  and  the  functions $\alpha_j,$ $\beta_j,$ $\gamma_j$ and $\delta_j$ are  given by 
  \begin{equation*}
   \alpha_j:=\sum_{q=1}^k  \alpha_{qj},\quad  \beta_j:=\sum_{q=1}^k  \beta_{qj},\quad
    \gamma_j:=\sum_{q=1}^k  \gamma_{qj},\quad  \delta_j:=\sum_{q=1}^k  \delta_{qj}.
  \end{equation*}

   Combining  this together with \eqref{e:fgh_M,M',j,n} and \eqref{e:F_N,M} and applying Lemma \ref{L:comparison-bis}, we  see that
  $ \tilde f_{N,M',j,n},$ $\tilde{\tilde f}_{N,M',j,n},$ $\tilde g_{N,M',j,n},$  $\tilde{\tilde g}_{N,M',j,n},$ $\tilde h_{N,M',j,n}$
  have the  desired asymptotic property.
\endproof

\begin{lemma}\label{L:Stokes} Fix  $\ell$  with $1\leq \ell\leq \ell_0$ and $r\in(0,\bfr].$    Set $\tilde\tau:=\tilde\tau_\ell$  and $\H_r:=\Tube (\widetilde V_\ell,r)\subset \E.$
Then, for every 
 every   current $S$ of bidimension $(q-1,q-1)$  defined on  $\U_\ell$ and every  smooth  form $\Phi$ of bidegree $(q,q)$  defined on $\tilde\tau(\H_r)$
 with $\pi(\supp(\Phi))\Subset \widetilde V_\ell.$ 
 Then \begin{equation*}\langle  \dbar (\tilde\tau_* S)^\sharp  ,\Phi  \rangle_{\partial[ \tilde \tau (\H_r)]}
 =\langle  \tilde\tau^*[ \dbar (\tilde\tau_* S)^\sharp] ,(\tilde\tau^*\Phi) \rangle_{\partial \H_r} .
 \end{equation*}
\end{lemma}
\proof
 Consider the  canonical  injections
 $\iota_{\partial \H}:\  \partial \H\hookrightarrow \E$ and $\iota_{\partial[ \tilde \tau (\H)]}:\ \partial[ \tilde \tau (\H)]\hookrightarrow \E.$
 Since   $ \iota_{\partial[ \tilde \tau (\H)]}\circ (\tilde\tau|_{\partial \H})= \tilde\tau\circ \iota_{\partial \H},$  it follows that
 $$   (\tilde\tau|_{\partial \H})^*\circ \iota^*_{\partial[ \tilde \tau (\H)]}=\iota^*_{\partial \H}\circ \tilde\tau^*.$$
Therefore,  we obtain that 
\begin{eqnarray*}
\langle  \dbar (\tilde\tau_* S)^\sharp  ,\Phi  \rangle_{\partial[ \tilde \tau (\H)]}
&=& \int_{\partial[ \tilde \tau (\H)]}  \iota^*_{\partial[ \tilde \tau (\H)]}  \dbar (\tilde\tau_* S)^\sharp  \wedge \Phi \\
&=&\int_{\partial\H}  (\tilde\tau|_{\partial \H})^*\big( \iota^*_{\partial[ \tilde \tau (\H)]} [ \dbar (\tilde\tau_* S)^\sharp  \wedge \Phi]\big)\\
&=&\int_{\partial\H}  \iota^*_{\partial \H}\big( \tilde\tau^*[ \dbar (\tilde\tau_* S)^\sharp  \wedge \Phi]\big)\\
&=& \int_{\partial\H}  \iota^*_{\partial \H}\big( \tilde\tau^*[ \dbar (\tilde\tau_* S)^\sharp]  \wedge \tilde\tau^*\Phi\big)\\
&=&\langle  \tilde\tau^*[ \dbar (\tilde\tau_* S)^\sharp] ,(\tilde\tau^*\Phi) \rangle_{\partial \H} .
\end{eqnarray*}
\endproof

\begin{proposition}\label{P:basic-bdr-estimates}
Fix  $\ell$  with $1\leq \ell\leq \ell_0$ and   set $\tilde\tau:=\tilde\tau_\ell.$ For $r\in(0,\bfr],$   set $\H_r:=\Tube (\widetilde V_\ell,r)\subset \E.$
Let $S$ be a   current in the class $\SH^{2,1}(\Tube(B,\bfr)).$  
Let $\Phi$ be  the  product of $\theta_\ell$ and a  smooth $(q,q)$-form on $\Tube(B,\bfr)$ 
which is  $2j$-negligible. 
 Then there are 
 \begin{itemize}
 \item [$\bullet$]  two  functions $\Ic_1,\ \Ic_2 :\ (0,\bfr]\to\R;$ 
  \item [$\bullet$]   three differential operators $D_{10},$ $D_{11},$  $D_{12}$ in  the class $\widehat\Dc^0_\ell;$
   and three differential operators $D_{20},$ $D_{21},$  $D_{22}$   in  the class $\Dc^0_\ell;$
  \item[$\bullet$]  three smooth $2q$-forms $\Phi_{10}$   which is  $(2j-1)$-negligible,
   $\Phi_{11}$  which  is  $2j$-negligible,  $\Phi_{12}$  which  is  $(2j-1)$-negligible;
   and three smooth $2q$-forms $\Phi_{20}$   which is  $2j$-negligible,
   $\Phi_{21}$  which  is  $(2j+1)$-negligible,  $\Phi_{22}$  which  is  $2j$-negligible;
 \end{itemize}
  such that 
 every $0<r_1<r_2\leq\bfr$ and  every smooth  function $\chi$ on $(0,\bfr),$  we have for $\nu\in\{1,2\},$
 \begin{equation}\label{e:Stokes-ddc-difference}
 \begin{split}
  \int_{r_1}^{r_2}\chi(t) \Ic_\nu(t)dt&= 
\int_{\Tube(B,r_1,r_2)}\chi(\|y\|) (D_{\nu 1}S\wedge \Phi_{\nu1})(y)+\int_{\Tube(B,r_1,r_2)}\chi'(\|y\|) (D_{\nu2}S\wedge \Phi_{\nu2})(y)\\
&+ \int_{\partial_\hor\Tube(B,r_2) }\chi(r_2)( D_{\nu0}S\wedge \Phi_{\nu0})(y)- \int_{\partial_\hor\Tube(B,r_1)} \chi(r_1)(D_{\nu 0}S\wedge \Phi_{\nu0})(y)  ,
 \end{split}
 \end{equation}
 and  that the  following inequality holds 
 for all $  0<r\leq \bfr:$
  \begin{equation}\label{e:inequal-ddc-difference}
 {1\over  r^{2(k-p-j)} } \int_{r\over 2}^r\big| \langle  \ddc (\tilde\tau_* S) -\tilde\tau_*(\ddc  S),\Phi  \rangle_{\tilde \tau (\H_t)}
 -  \Ic_1(t) -\Ic_2(t)\big|dt
  \leq 
\sum_{m=\lowm}^\upm \nu_m( S,B,r,\id).
 \end{equation}
\end{proposition}
\proof

 By  Proposition \ref{P:Stokes} we have
\begin{equation}\label{e:sum-ddc-difference}
 \langle  \ddc (\tilde\tau_* S) -\tilde\tau_*(\ddc  S),\Phi  \rangle_{\tilde \tau (\H_t)}=I_1(t)+I_2(t)-I_3(t)+{1\over2\pi i}I_4(t)
 -{1\over2\pi i}I_5(t)-{1\over \pi i} \tilde I_6(t),
 \end{equation}
 where 
 \begin{eqnarray*}
 I_1(t)&:=&\langle  S,\tilde\tau^*(\ddc\Phi) -\ddc(\tilde\tau^*\Phi)\rangle_{\H_t}, \\ 
 I_2(t)&:=& \langle    S ,\dc (\tilde\tau^*\Phi)^\sharp -\tilde\tau^*(\dc\Phi) \rangle_{\partial \H_t},\\
 I_3(t)&:=& \langle \tilde\tau^*[(\tilde\tau_*S)^\sharp]-S, \tilde\tau^*(\dc\Phi)\rangle_{\partial\H_t}, \\
 I_4(t)&:=&  \langle  \tilde\tau^*[ (\tilde\tau_* S)^\sharp]  ,\tilde\tau^*(d\Phi) -d[(\tilde\tau^*\Phi)^\sharp] \rangle_{\partial \H_t},\\
 I_5(t)&:=&\langle   S -\tilde\tau^*[(\tilde\tau_*S)^\sharp)],d[(\tilde\tau^*\Phi)^\sharp]  \rangle_{\partial \H_t},\\
 \tilde I_6(t)&:=& \langle  \dbar (\tilde\tau_* S)^\sharp  ,\Phi  \rangle_{\partial[ \tilde \tau (\H_t)]} -\langle  \dbar  S ,(\tilde\tau^*\Phi)^\sharp  \rangle_{\partial \H_t} .
\end{eqnarray*}
By Lemma \ref{L:Stokes}, we have 
\begin{eqnarray*}  
\tilde I_6(t)
 &=&\langle  \tilde\tau^*[ \dbar (\tilde\tau_* S)^\sharp] ,(\tilde\tau^*\Phi) \rangle_{\partial \H_r} -\langle  \dbar  S ,(\tilde\tau^*\Phi)^\sharp  \rangle_{\partial \H_t} \\
 &=&\langle  \tilde\tau^*[ \dbar (\tilde\tau_* S)^\sharp] - \dbar  S,(\tilde\tau^*\Phi) \rangle_{\partial \H_r} +\langle  \dbar  S ,(\tilde\tau^*\Phi)-(\tilde\tau^*\Phi)^\sharp  \rangle_{\partial \H_t}\\
 &=& I_6(t)+I_7(t).
 \end{eqnarray*}
This, coupled with  \eqref{e:sum-ddc-difference}, implies the following reduction.
To  prove  the  proposition, we only need to show for $1\leq n\leq 7,$  that
there are 
 \begin{itemize}
 \item [$\bullet$]  two  functions $\Ic^{(n)}_1,\ \Ic^{(n)}_2 :\ (0,\bfr]\to\R;$ 
  \item [$\bullet$]   three differential operators $D^{(n)}_{10},$ $D^{(n)}_{11},$  $D^{(n)}_{12}$ in  the class $\widehat\Dc^0_\ell;$
   and three differential operators $D^{(n)}_{20},$ $D^{(n)}_{21},$  $D^{(n)}_{22}$   in  the class $\Dc^0_\ell;$
  \item[$\bullet$]  three smooth $2q$-forms $\Phi^{(n)}_{10}$   which is  $(2j-1)$-negligible,
   $\Phi^{(n)}_{11}$  which  is  $2j$-negligible,  $\Phi^{(n)}_{12}$  which  is  $(2j-1)$-negligible;
   and three smooth $2q$-forms $\Phi^{(n)}_{20}$   which is  $2j$-negligible,
   $\Phi^{(n)}_{21}$  which  is  $(2j+1)$-negligible,  $\Phi^{(n)}_{22}$  which  is  $2j$-negligible;
 \end{itemize}
  such that 
 every $0<r_1<r_2\leq\bfr$ and  every smooth  function $\chi$ on $(0,\bfr),$  we have for $\nu\in\{1,2\},$
 \begin{equation}\label{e:Stokes-ddc-difference-n}
 \begin{split}
  \int_{r_1}^{r_2}\chi(t) \Ic_\nu(t)dt&= 
\int_{\Tube(B,r_1,r_2)}\chi(\|y\|) (D^{(n)}_{\nu1}S\wedge \Phi^{(n)}_{\nu1})(y)+\int_{\Tube(B,r_1,r_2)}\chi'(\|y\|) (D^{(n)}_{\nu2}S\wedge \Phi^{(n)}_{\nu2})(y)\\
&+ \int_{\partial_\hor\Tube(B,r_2) }\chi(r_2)( D^{(n)}_{\nu0}S\wedge \Phi^{(n)}_{\nu0})(y)- \int_{\partial_\hor\Tube(B,r_1)} \chi(r_1)(D^{(n)}_{\nu0}S\wedge \Phi^{(n)}_{\nu0})(y)  ,
 \end{split}
 \end{equation}
 and  that the  following inequality holds 
 for all $  0<r\leq \bfr:$
  \begin{equation}\label{e:inequal-ddc-difference-n}
 {1\over  r^{2(k-p-j)} } \int_{r\over 2}^r\big|  I_n(t)
 -  \Ic^{(n)}_1(t) -\Ic^{(n)}_2(t)\big|dt
  \leq 
\sum_{m=\lowm}^\upm \nu_m( S,B,r,\id).
 \end{equation}
 Indeed, it suffices to  consider the functions $\Ic_\nu: \ (0,\bfr]\to\R$ for $\nu=1,2,$ defined by 
 \begin{equation*}
  \Ic_\nu(t):= \Ic^{(1)}_\nu(t)+\Ic^{(2)}_\nu(t)-\Ic^{(3)}_\nu(t)+{1\over2\pi i}\Ic^{(4)}_\nu(t)
 -{1\over2\pi i}\Ic^{(5)}_\nu(t)-{1\over \pi i}  \Ic^{(6)}_\nu(t) -{1\over \pi i}  \Ic^{(7)}_\nu(t)\quad\text{for}\quad t\in(0,\bfr].
 \end{equation*}
Then  we see that  the equality \eqref{e:Stokes-ddc-difference-n}  (resp. the inequality   \eqref{e:inequal-ddc-difference-n})  follows  from combining  the equalities  \eqref{e:Stokes-ddc-difference-n}  (resp. the inequalities  \eqref{e:inequal-ddc-difference-n})
for $1\leq n\leq 7.$
\endproof

\section{Positive  plurisubharmonic  currents and quasi-monotonicity of the Lelong numbers}
\label{S:Positive-psh-currents}

\subsection{Preliminary  estimates }
\label{SS:Preliminaries-Section:psh-and-quasi-mono}

Let $T$ be  positive plurisubharmonic  current $T$ of bidegree $(p,p)$ on $\bfU.$ 
Consider the integers 
\begin{equation}\label{e:m+}
 \lowm^+:=\max(0,l-p-1)\qquad\text{and}\qquad  \upm^+:=\min(l,k-p-1).
\end{equation}
In other words,    $\lowm^+,$ $\upm^+$ are    associated  to the  $(p+1,p+1)$-current
$\ddc T$ in the same  way  as    $\lowm,$ $\upm$ are    associated  to the  $(p,p)$-current
$T$ in formula \eqref{e:m}. 

Following  the model of \eqref{e:global-mass-indicators}, consider the following   mass indicators, for $0<r\leq\bfr,$
\begin{equation}\label{e:global-mass-indicators-bis}
 \begin{split}
 \Mc^\tot( T,r)&:=\sum_{j=0}^{\upm} \Mc_j( T,r),\quad    \Mc^\tot( \ddc T,r):=\sum_{j=0}^{\upm^+} \Mc_j( \ddc T,r) \\
 \Nc(T,r)&:= \Mc^\tot(T,r)+\Mc^\tot(\ddc T,r)=\sum_{j=0}^{\upm} \Mc_j( T,r)+\sum_{j=0}^{\upm^+} \Mc_j(\ddc T,r).
 \end{split}
 \end{equation}
where  
the $\Mc_j$'s   are defined  in \eqref{e:global-mass-indicators}.

In this section 
following   Definition \ref{D:sup}, we  introduce the following class of currents.


\begin{definition}\label{D:sup-bis}\rm
Fix an open neighborhood $\bfU$ of $\overline B$ and an open neighborhood $\bfW$ of $\partial B$ in $X$ with $\bfW\subset \bfU.$
Let $\widetilde\SH^{3,3}_p(\bfU,\bfW)$ be the  set of all $T\in \SH^{3,3}_p(\bfU,\bfW)$  whose  a sequence of approximating  forms $(T_n)_{n=1}^\infty$
satisfies the following   condition:
 \begin{equation}\label{e:unit-SH-3,3} \|T_n\|_{\bfU}\leq  1 \quad\text{and}\quad \| \ddc T_n\|_{\bfU}\leq 1  \quad\text{and}\quad  \| T_n\|_{\Cc^3(\bfW)}\leq 1.\end{equation}

Given   a  class of currents  $\Fc$ and  a mass indicator $\Mc(T)$   for all currents $T\in\Fc,$  
 We denote  by $\sup_{T\in\Fc}\Mc(T)$  the supremum of $\Mc(T)$  when  $T$ is taken over $\Fc.$
\end{definition}

Recall some notation from  the Extended Standing  Hypothesis in Subsection \ref{SS:Ex-Stand-Hyp}.
Consider  a  strongly admissible map $\tau:\ \bfU\to\tau(\bfU)$  along $B,$  with $\bfU$ a neighborhood of $\overline B$ in $X.$
By shrinking $\bfU$ if necessary, we may   fix a finite collection $\Uc=(\bfU_\ell,\tau_\ell)_{1\leq \ell\leq \ell_0} ,$
 of holomorphic admissible maps for $\bfU.$  More precisely,  
we fix  a  finite cover of $\overline \bfU$ by open subsets $\bfU_\ell,$ $1\leq \ell\leq  \ell_0,$ of $X$
such that   there is   a holomorphic coordinate system on $\overline \bfU_\ell$ in $X$ and $\bfU_\ell$  is  biholomorphic  to  $\U_\ell:=\tau_\ell(\bfU_\ell)\subset \E$
by a  holomorphic admissible map $\tau_\ell.$ By  choosing $\bfr>0$ small enough, we may assume  without loss of generality that $\overline\Tube(B,\bfr)\Subset \U:=\bigcup_{\ell=1}^{\ell_0} \U_\ell.$ 
Choose   a partition of unity  $(\theta_\ell)_{1\leq \ell\leq \ell_0}$ subordinate   to the open cover  $(\bfU_\ell\cap V)_{1\leq \ell\leq \ell_0}$   of $\overline{\bfU\cap V}$  in $V$  such that $\sum_{1\leq \ell\leq \ell_0}  \theta_\ell=1$ on an open neighborhood of $\overline {\bfU\cap V} \subset V.$ We suppose   without loss of generality  that there are open  subsets
$\widetilde V_\ell\subset V$  for $1\leq \ell\leq \ell_0$ such that
\begin{equation*}
\supp(\theta_l)\subset  \widetilde V_\ell\Subset \bfU_\ell\cap V \quad\text{and}\quad 
\tau( \widetilde V_\ell)\Subset  \U_\ell\quad\text{and}\quad \pi^{-1}(\supp(\theta_\ell))\cap \U\subset \U_\ell.
\end{equation*}
 For $1\leq \ell\leq \ell_0$ set 
 \begin{equation*}
 \tilde \tau_\ell:=\tau\circ\tau_\ell^{-1}.
 \end{equation*}
  So  $\tilde \tau_\ell$ defines a map  from  $\U_\ell\subset \E$  onto  $\tau(\bfU_\ell)\subset \E.$
  We may suppose  that  for every $1\leq\ell\leq\ell_0,$  there is a local  coordinate system $y=(z,w)$  on $\U_\ell$  with  $V\cap \U_\ell=\{z=0\}.$

Fix an integer $j$ with $\lowm\leq j\leq \upm.$ 
Consider  the  forms on $\bfU$:
\begin{equation}\label{e:partition-can-forms-j,l}
\Phi:= \pi^*(\omega^j)\wedge \beta^{k-p-j-1} \quad\text{and}\quad \Phi^{(\ell)}:= (\pi^*\theta_\ell)\cdot\pi^*(\omega^j)\wedge \beta^{k-p-j-1}\quad\text{for}\quad
1\leq\ell\leq \ell_0.
\end{equation}
So we have
\begin{equation}\label{e:sum-forms-j,l}
\Phi=\sum_{\ell=1}^{\ell_0} \Phi^{(\ell)}\qquad\text{on}\qquad \bfU.
\end{equation}
For $\ell$  with $1\leq \ell\leq \ell_0$ and   set $\tilde\tau:=\tilde\tau_\ell.$
For $r\in(0,\bfr],$   set $\H_r:=\Tube (\widetilde V_\ell,r)\subset \E.$

Let $T$ be a positive  plurisubharmonic   current on $\bfU$ in the class  $\widetilde\SH^{3,3}_p(\bfU,\bfW).$ 
Consider the  current
\begin{equation}\label{e:S_ell}
 S^{(\ell)}:=(\tau_\ell)_*(T|_{\bfU_\ell})   .
\end{equation}
By \eqref{e:T-hash} we get that
\begin{equation}
 T^\hash =\sum_{\ell=1}^{\ell_0} (\pi^*\theta_\ell)\cdot S^{(\ell)}.
\end{equation}
Note that the current $S^{(\ell)}$ is  positive plurisubharmonic on $\H_\bfr.$ Moreover, by Lemma
\ref{L:ex-negligble},   $\Phi^{(\ell)}$  is a $2j$-negligible smooth form. 
By Proposition \ref{P:basic-bdr-estimates},
there are 
 \begin{itemize}
 \item [$\bullet$]  two  functions $\Ic^{(\ell)}_1,\ \Ic^{(\ell)}_2 :\ (0,\bfr]\to\R;$ 
  \item [$\bullet$]   three differential operators $D^{(\ell)}_{10},$ $D^{(\ell)}_{11},$  $D^{(\ell)}_{12}$ in  the class $\widehat\Dc^0_\ell;$
   and three differential operators $D^{(\ell)}_{20},$ $D^{(\ell)}_{21},$  $D^{(\ell)}_{22}$   in  the class $\Dc^0_\ell;$
  \item[$\bullet$]  three smooth $2q$-forms $\Phi^{(\ell)}_{10}$   which is  $(2j-1)$-negligible,
   $\Phi^{(\ell)}_{11}$  which  is  $2j$-negligible,  $\Phi^{(\ell)}_{12}$  which  is  $(2j-1)$-negligible;
   and three smooth $2q$-forms $\Phi^{(\ell)}_{20}$   which is  $2j$-negligible,
   $\Phi^{(\ell)}_{21}$  which  is  $(2j+1)$-negligible,  $\Phi^{(\ell)}_{22}$  which  is  $2j$-negligible;
 \end{itemize}
  such that 
 every $0<r_1<r_2\leq\bfr$ and  every smooth  function $\chi$ on $(0,\bfr],$  we have for $\nu\in\{1,2\},$
 \begin{equation}\label{e:Stokes-ddc-difference-bis}
 \begin{split}
  \int_{r_1}^{r_2}\chi(t) \Ic^{(\ell)}_\nu(t)dt&= 
\int_{\Tube(B,r_1,r_2)}\chi(\|y\|) (D^{(\ell)}_{\nu 1}S^{(\ell)}\wedge \Phi^{(\ell)}_{\nu 1})(y)+\int_{\Tube(B,r_1,r_2)}\chi'(\|y\|) (D^{(\ell)}_{ \nu 2}S^{(\ell)}\wedge \Phi^{(\ell)}_{\nu 2})(y)\\
&+ \int_{\partial_\hor\Tube(B,r_2) }\chi(r_2)( D^{(\ell)}_{\nu 0}S^{(\ell)}\wedge \Phi^{(\ell)}_{\nu 0})(y)- \int_{\partial_\hor\Tube(B,r_1)} \chi(r_1)(D^{(\ell)}_{\nu 0}S^{(\ell)}\wedge \Phi^{(\ell)}_{\nu 0})(y)  ,
 \end{split}
 \end{equation}
 and  that the  following inequality holds 
 for all $  0<t\leq \bfr:$
  \begin{equation}\label{e:inequal-ddc-difference-bis}
 {1\over  r^{2(k-p-j)} } \int_{r\over 2}^r\big| \langle  \ddc [(\tilde\tau_\ell)_* S^{(\ell)}] -(\tilde\tau_\ell)_*(\ddc  S^{(\ell)}),\Phi^{(\ell)}  \rangle_{\tilde \tau (\H_t)}
 -  \Ic^{(\ell)}_1(t) -\Ic^{(\ell)}_2(t)\big|dt
  \leq 
\sum_{m=\lowm}^\upm \nu_m( S^{(\ell)},B,r,\id).
 \end{equation} 
 The following  auxiliary  results are needed.
 \begin{lemma}\label{L:partition}
  The  following   equalities hold:
  \begin{eqnarray*}
     (\tilde\tau_\ell)_* S^{(\ell)}&=&\tau_*T \qquad \text{and}\quad (\tilde\tau_\ell)_*(\ddc  S^{(\ell)})=\tau_*(\ddc T)\quad\text{on}\quad \bfU_\ell,\\
     \sum_{\ell=1}^{\ell_0}  \ddc[ (\tilde\tau_\ell)_* S^{(\ell)}]\wedge \Phi^{(\ell)}  &=&\ddc (\tau_*T)\wedge \Phi\quad\text{and}\quad
      \sum_{\ell=1}^{\ell_0} (\tilde\tau_\ell)_*(  \ddc S^{(\ell)})\wedge \Phi^{(\ell)}  = \tau_*(\ddc T)\wedge \Phi
      \quad\text{on}\quad \bfU.
  \end{eqnarray*}
 \end{lemma}
\proof Since 
$ \tilde \tau_\ell\circ\tau_\ell=\tau$  on $\bfU_\ell,$  the  first  equality follow from \eqref{e:S_ell}.
Since $\tau_\ell$ is  holomorphic,  we infer from 
\eqref{e:S_ell} that $\ddc S^{(\ell)}=(\tau_\ell)_*(\ddc T|_{\bfU_\ell})   .$ The second equality can be proved as the  first one.

We deduce from the first  equality that  
 $\ddc[ (\tilde\tau_\ell)_* S^{(\ell)}]\wedge \Phi^{(\ell)}  =\ddc (\tau_*T)\wedge \Phi^{(\ell)}$  on $\bfU.$
 Summing this equality over $1\leq \ell\leq  \ell_0,$ we obtain  the  third equality.

We deduce from the second   equality that $(\tilde\tau_\ell)_*(  \ddc S^{(\ell)})\wedge \Phi^{(\ell)}  = \tau_*(\ddc T)\wedge \Phi^{(\ell)}$ on $\bfU.$  Summing this equality over $1\leq \ell\leq \ell_0,$ we obtain  the  last equality.
 \endproof

 \begin{lemma}\label{L:Ic-0}
  Under  the above  hypotheses and  notations,  there is a constant $c$ independent of $T$  such that for $\nu\in \{1,2\}$ and for all $1\leq \ell\leq \ell_0$ and for all $0<r\leq \bfr:$
  \begin{equation*}
   {1\over  r^{2(k-p-j)} } \int_{r\over 2}^r\big|  \int_{\partial_\hor\Tube(B,t) }( D^{(\ell)}_{\nu0}S^{(\ell)}\wedge \Phi^{(\ell)}_{\nu0}) \big|dt\leq  cr^2\Mc^\tot(T,r).
  \end{equation*}
 \end{lemma}
\proof Using the above-mentioned  property of  the operators $D_{\nu 0}$ and the smooth forms $\Phi_{\nu 0},$ 
we are able  to  apply  Proposition \ref{P:basic_boundary-inequ}
for $m=2j+1.$ This, combined  with \eqref{e:global-mass-indicators} and  \eqref{e:global-mass-indicators-bis}, implies the result.
\endproof
 Consider  two functions $\chi_1,\chi_2: (0,r]\to \R^+$ defined by    
  \begin{equation}\label{e:chi_1-chi_2}\chi_1(t):= {t\over  r^{2(k-p-j)} }\quad\text{and}\quad
  \chi_2(t):={1\over  t^{2(k-p-j)-1} } \quad\text{for}\quad t\in (0,r].
  \end{equation}
  \begin{lemma}\label{L:Ic-1-and-2}
Under  the above  hypotheses and  notations, let  $0<r\leq \bfr.$
  Then  there is a constant $c$ independent of $T$ and $r$ such that for $\nu\in \{1,2\}$ and for all $1\leq \ell\leq \ell_0$ and for
  all $0<s<r:$
  \begin{eqnarray*}
  \big |\int_{\Tube(B,s,r)}\chi(\|y\|) (D^{(\ell)}_{\nu 1}S^{(\ell)}\wedge \Phi^{(\ell)}_{\nu1})(y)\big|  &\leq & c\sum_{n=0}^\infty {r\over 2^n}\Mc^\tot(T,{r\over 2^n}),\\
  \big|  \int_{\Tube(B,s,r)}\chi'(\|y\|) (D^{(\ell)}_{\nu2}S\wedge \Phi^{(\ell)}_{\nu2})(y) \big|&\leq & c\sum_{n=0}^\infty {r\over 2^n}\Mc^\tot(T,{r\over 2^n}).
  \end{eqnarray*}
  Here $\chi$ is  either the function $\chi_1$ or the function $\chi_2$ given in \eqref{e:chi_1-chi_2}.
 \end{lemma}
\proof  There is $N\in\N$  such that $s':=2^{-N}r$ satisfies ${s'\over 2}\leq s\leq  s'.$ 
Observe that
\begin{eqnarray*}
 \big |\int_{\Tube(B,s,r)}\chi(\|y\|) (D^{(\ell)}_{\nu 1}S^{(\ell)}\wedge \Phi^{(\ell)}_{\nu1})(y)\big| 
 &\leq &  \sum_{n=0}^{N} |\int_{\Tube(B,{r\over 2^{n+1}},{r\over 2^n})}\chi(\|y\|) (D^{(\ell)}_{\nu 1}S^{(\ell)}\wedge \Phi^{(\ell)}_{\nu1})(y)\big|\\
 &+& |\int_{\Tube(B,s,s')}\chi(\|y\|) (D^{(\ell)}_{\nu 1}S^{(\ell)}\wedge \Phi^{(\ell)}_{\nu1})(y)\big|
\end{eqnarray*}
By  Proposition \ref{P:mass-fine-negligible},
\begin{eqnarray*}
|\int_{\Tube(B,{r\over 2^{n+1}},{r\over 2^n})}\chi(\|y\|) (D^{(\ell)}_{\nu 1}S^{(\ell)}\wedge \Phi^{(\ell)}_{\nu1})(y)\big|
&\leq&c{r\over 2^n}\Mc^\tot(T,{r\over 2^n}),\\
|\int_{\Tube(B,s,s')}\chi(\|y\|) (D^{(\ell)}_{\nu 1}S^{(\ell)}\wedge \Phi^{(\ell)}_{\nu1})(y)\big|
&\leq& c s\Mc^\tot(T,s).
\end{eqnarray*}
Combining this estimates,  the  first  inequality of the lemma follows. 

The second  inequality can be proved in the  same way.
\endproof

 \subsection{Quasi-positivity and quasi-monotonicity of the  Lelong numbers  and finiteness of the mass indicators $\Mc_j$}
 \label{SS:quasi-monotone-psh}

\begin{lemma}\label{L:ddc-difference-lambda} For  all $r_1,r_2\in  (0,\bfr]$ with  $r_1<r_2,$  
there is a constant $c>0$  such that  for every
$j$ with $\lowm\leq j\leq \upm,$  and every $m$ with $0\leq m\leq j,$ and every positive plurisubharmonic  current $T$   in the class  $\widetilde\SH^{3,3}_p(\bfU,\bfW),$
  the following two inequalities hold for all $\lambda\geq 1:$
\begin{multline*}
\big| \int_{r_1}^{r_2} \big( {1\over t^{2(k-p-j)}}-{1\over r_2^{2(k-p-j)}}  \big)2tdt\int_{\Tube(B,t)} (A_{\lambda})_*(\ddc (\tau_*T)-\tau_*(\ddc T))\wedge \pi^*(\omega^{j-m})\wedge \beta^{k-p-j+m-1} \big| \\
 \leq c\sum_{n=0}^\infty {1\over (2^n\lambda)^{2m+1}}\Mc^\tot(  T,{r_2\over 2^n\lambda} ),
 \end{multline*}
 \begin{multline*}
  \big( {1\over r_1^{2(k-p-j)}}-{1\over r_2^{2(k-p-j)}}  \big) 
 \big|\int_{0}^{r_1}2tdt \int_{\Tube(B,t)} (A_{\lambda})_*(\ddc (\tau_*T)-\tau_*(\ddc T))\wedge\pi^*(\omega^{j-m})\wedge \beta^{k-p-j+m-1}\big|\\
  \leq c\sum_{n=0}^\infty{1\over (2^n \lambda)^{2m+1}}\Mc^\tot( T,{r_1\over 2^n\lambda} ).
\end{multline*} 
 \end{lemma}
 \proof
  We only give the proof of the  first  inequality since the second one can be done similarly.
In fact,  the first  inequality will  follow if one can  show that for $i\in\{1,2\}$ and for $\chi_i$  defined in \eqref{e:chi_1-chi_2},  \begin{multline*}
\big| \int_{r_1\over \lambda }^{r_2\over\lambda} \chi_i(t)dt\int_{\Tube(B,t)} (\ddc (\tau_*T)-\tau_*(\ddc T))\wedge \pi^*(\omega^{j-m})\wedge \beta^{k-p-j+m-1} \big| \\
 \leq {c\over\lambda} \sum_{n=0}^\infty {1\over 2^{n(2m+1)}}\Mc^\tot(  T,{r_2\over 2^n\lambda} ).
 \end{multline*}
  Combining \eqref{e:inequal-ddc-difference-bis} and \eqref{e:Stokes-ddc-difference-bis} and  Lemmas \ref{L:partition}, \eqref{L:Ic-0}
   and \eqref{L:Ic-1-and-2}, the  last inequality  follows.
 \endproof

\begin{lemma}\label{L:ddc-difference-m-1-lambda} 
   For  all $r_1,r_2\in  (0,\bfr]$ with  $r_1<r_2,$  there is a constant $c>0$  such that  for every
$j$ with $\lowm\leq j\leq \upm,$ and every $m$ with $1\leq m\leq j,$  and every positive plurisubharmonic  current $T$   in the class  $\widetilde\SH^{3,3}_p(\bfU,\bfW),$
 the  following inequality  holds for  for every $\lambda\geq 1$:
\begin{multline*}
 \Big|\int_{r_1}^{r_2} \big( {1\over t^{2(k-p-j)}}-{1\over r_2^{2(k-p-j)}}  \big)2tdt\int_{\Tube(B,t)} (A_{\lambda})_*(\tau_*(\ddc T))\wedge \pi^*(\omega^{j-m})\wedge \beta^{k-p-j+m-1}\Big| \\
 \leq {c\over \lambda^{2m}}\Mc^\tot(\ddc  T,{r_2\over \lambda} ),
 \end{multline*}
 \begin{multline*} \Big|\big( {1\over r_1^{2(k-p-j)}}-{1\over r_2^{2(k-p-j)}}  \big) 
  \int_{0}^{r_1}2tdt\int_{\Tube(B,t)} (A_{\lambda})_*(\tau_*(\ddc T))\wedge \pi^*(\omega^{j-m})\wedge \beta^{k-p-j+m-1} \Big|\\
  \leq {c\over\lambda^{2m}}\Mc^\tot(\ddc  T,{r_1\over \lambda} ).
\end{multline*} 
\end{lemma}
\proof We argue as  in the proof of Lemma
 \ref{L:Mc-hash-vs-nu}.
\endproof

\begin{corollary}\label{C:ddc-estimate-m-1-lambda}  
 For  $r_1,r_2\in  (0,\bfr]$ with  $r_1<r_2,$ and for $\lowm \leq j\leq \upm$   and   for   $m$ with $0\leq m\leq j,$ there is a constant $c>0$  such that  the following two inequalities hold $\lambda\geq 1:$
\begin{multline*}
\Big |\int_{r_1}^{r_2} \big( {1\over t^{2(k-p-j)}}-{1\over r_2^{2(k-p-j)}}  \big)2tdt\int_{\Tube(B,t)} (A_{\lambda_n})_*(\ddc (\tau_*T))\wedge \pi^*(\omega^{j-m})\wedge \beta^{k-p-j+m-1} \Big | \\
 \leq c\sum_{n=0}^\infty {1\over (2^n\lambda)^{2m+1}}\Mc^\tot(  T,{r_2\over 2^n\lambda} )+{c\over \lambda^{2m}}\Mc^\tot( \ddc T,{r_2\over \lambda} ),
 \end{multline*}
 \begin{multline*}
  \Big |\big( {1\over r_1^{2(k-p-j)}}-{1\over r_2^{2(k-p-j)}}  \big) 
 \int_{0}^{r_1}2tdt \int_{\Tube(B,t)} (A_{\lambda_n})_*(\ddc (\tau_*T))\wedge \pi^*(\omega^{j-m})\wedge \beta^{k-p-j+m-1}\Big |\\
  \leq c\sum_{n=0}^\infty {1\over (2^n\lambda)^{2m+1}}\Mc^\tot(  T,{r_1\over 2^n\lambda} )+{c\over \lambda^{2m}}\Mc^\tot( \ddc T,{r_1\over \lambda} ).
\end{multline*} 
\end{corollary}
\proof
It follows from a combination of Lemmas \ref{L:ddc-difference-lambda} and \ref{L:ddc-difference-m-1-lambda}. 
\endproof

\begin{proposition}\label{P:nu_tot-monotone-psh} 
Let  $0< r_1<r_2\leq\bfr.$ Then there are  a  family  $\Dc=\{ d_{jq}\in\R:\ 0\leq j\leq k-p-q,\ 0\leq q\leq  k-l\}$ and   a constant $c>0 $ depending  on $r_1$ and $r_2$  such that for every positive closed   current $T$ on $\bfU$ belonging to the class $\widetilde\SH^{3,3}_p(\bfU,\bfW),$ 
the  following inequality  hold for $0\leq q\leq k-l:$
\begin{equation*}
  \nu^\Dc_q\big(T,B,{r_1\over \lambda},\tau\big)\leq  \nu^\Dc_q(T,B,{r_2\over \lambda},\tau\big) + {c\over \lambda}+ c\lambda^{-1}\sum_{n=0}^\infty 
  2^{-n}\Mc^\tot(T,{r_2\over 2^n\lambda})+c\Mc^\tot(\ddc T,{r_2\over \lambda})
  \quad\text{for}\quad \lambda\gg 1.
 \end{equation*}
Moreover, for  every  $\epsilon>0$  we can choose $\Dc$ such that  $d_{k-p-q,q} <\epsilon^q d_{k-p-q+1,q-1}$ for $1\leq q\leq k-l.$
\end{proposition}

\proof
 We argue as in the   Proposition \ref{P:nu_tot-monotone}
using  Lemma \ref{L:basic-difference-nu_j,q_estimate-psh} below instead of Lemma 
\ref{L:basic-difference-nu_j,q_estimate}.
\endproof

\begin{lemma}\label{L:basic-difference-nu_j,q_estimate-psh} Given $0<r_1<r_2\leq \bfr,$ 
there is  a constant $c>0$ such that  for every positive closed current $T\in  \widetilde\SH^{3,3}_p(\bfU,\bfW)$
and  $0\leq q\leq k-l$ and $0\leq j\leq \min(\upm,k-p-q),$ the following inequality holds:
 \begin{multline*}
\nu_{j,q}\big(T,B,{r_2\over \lambda},\tau\big)  -\nu_{j,q}\big(T,B,{r_1\over \lambda},\tau\big)
\geq  \Kc_{j,q}\big (T,{r_1\over\lambda},{r_2\over \lambda}  \big)-  c\lambda^{-1} -  c\lambda^{-1}\sum_{n=0}^\infty 2^{-n}\Mc^\tot(T,{r_2\over \lambda})\\-c\Mc^\tot(\ddc T,{r_2\over \lambda})
   -  c\lambda^{1\over 2} \Kc_{q}\big (T,{r_1\over\lambda},{r_2\over \lambda} \big )-c\Kc_{q-1} (T,{r_1\over\lambda},{r_2\over \lambda})  
-c\sqrt{\Kc_{q}(T,{r_1\over\lambda},{r_2\over \lambda})} \sqrt{ \Kc^-_{j,q}(T,{r_1\over\lambda},{r_2\over \lambda})  }.
\end{multline*}
\end{lemma}
\proof
Fix  $0\leq q_0\leq k-l.$   
Let $0\leq j_0\leq \min(\upm,k-p-q_0).$   Set  $j'_0:=k-p-q_0-j_0\geq 0.$
We may assume without loss of generality that $T$ is  $\Cc^3$-smooth.   
Since the $(1,1)$-smooth forms $\omega$ and $\beta$ are closed, it follows that
$$
\ddc[\tau_*T\wedge \pi^*(\omega^{j_0})\wedge \beta^{j'_0}]=\ddc(\tau_*T)\wedge \pi^*(\omega^{j_0})\wedge \beta^{j'_0}.
$$
Applying Theorem  \ref{T:Lelong-Jensen-smooth}  to  $\tau_*T\wedge \pi^*(\omega^{j_0})\wedge \beta^{j'_0},$   we get that
\begin{equation*}
\begin{split}
&{\lambda^{2q_0}\over  r_2^{2q_0}}\int_{\Tube(B,{r_2\over \lambda})}\tau_*T\wedge \pi^*(\omega^{j_0})\wedge \beta^{k-p-j_0}
-{\lambda^{2q_0}\over  r_1^{2q_0}}\int_{\Tube(B,{r_1\over\lambda})}\tau_*T\wedge \pi^*(\omega^{j_0})\wedge \beta^{k-p-j_0}\\
&= \Vc\big(\tau_*T\wedge \pi^*(\omega^{j_0})\wedge \beta^{j'_0},{r_1\over\lambda},{r_2\over \lambda}\big)+\int_{\Tube(B,{r_1\over\lambda},{r_2\over \lambda})}\tau_*T\wedge \pi^*(\omega^{j_0})\wedge \beta^{j'_0}\wedge\alpha^{q_0}\\
&+  \int_{r_1\over\lambda}^{r_2\over\lambda} \big( {1\over t^{2q_0}}-{\lambda^{2q_0}\over r_2^{2q_0}}  \big)2tdt\int_{\Tube(B,t)} \ddc (\tau_*T)\wedge \pi^*(\omega^{j_0})\wedge \beta^{q_0+j'_0-1}  \\
&+\big( {\lambda^{2q_0}\over r_1^{2q_0}}-{\lambda^{2q_0}\over r_2^{2q_0}}  \big) 
 \int_{0}^{r_1\over \lambda }2tdt \int_{\Tube(B,t)} \ddc (\tau_*T)\wedge \pi^*(\omega^{j_0})\wedge \beta^{q_0+j'_0-1} .
\end{split}
\end{equation*}
By Corollary \ref{C:ddc-estimate-m-1-lambda},
the  last  two double  integrals  are  of order smaller than  
$$
{c \lambda^{-1}}\sum_{n=0}^\infty 2^{-n}\Mc^\tot( T,{r_2\over 2^n\lambda} )+c\Mc^\tot( \ddc T,{r_2\over \lambda} )
.$$ Moreover,  by Theorem \ref{T:vertical-boundary-terms}, we have the following estimate independently of $T:$
\begin{equation*}
\Vc\big(\tau_*T\wedge \pi^*(\omega^{j_0})\wedge \beta^{j'_0},{r_1\over\lambda},{r_2\over \lambda}\big)=O(\lambda^{-1}).
\end{equation*}
Therefore,  there is a constant $c>0$ independent of $T$ such that for $\lambda\geq 1,$
\begin{multline*}
\big|\int_{\Tube(B,{r_1\over\lambda},{r_2\over \lambda})}\tau_*T\wedge \pi^*(\omega^{j_0})\wedge \beta^{j'_0}\wedge\alpha^{q_0}
-\big (\nu_{j_0,q_0}\big(T,B,{r_2\over \lambda},\tau\big)   -  \nu_{j_0,q_0}\big(T,B,{r_1\over \lambda},\tau\big) \big)  \big| \\
\leq  c\lambda^{-1}+ {c \lambda^{-1}}\sum_{n=0}^\infty  2^{-n}\Mc^\tot( T,{r_2\over 2^n \lambda} )+c\Mc^\tot( \ddc T,{r_2\over \lambda} ).
\end{multline*} 
 Arguing as  in the proof of \eqref{e:P-Lc-bullet_finite(5)}, we obtain  that
\begin{equation*}
 \begin{split}
&\int_{\Tube(B,{r_1\over\lambda},{r_2\over \lambda})}\tau_*T\wedge \pi^*(\omega^{j_0})\wedge \beta^{j'_0}\wedge\alpha^{q_0}
  = I_{q_0,0,j_0,0}(T,{r_1\over\lambda},{r_2\over \lambda})\\&+\sum_{j_1,j'_1,j''_1} {j'_0\choose  j'_1}{q_0 \choose j_1}
  {q_0-j_1 \choose j''_1}(-c_1)^{j'_0-j'_1}(-1)^{q_0-j_1-j''_1}
 I_{j_1, j'_0-j'_1, q_0+j_0+j'_0-j_1-j'_1-j''_1,q_0-j_1-j''_1}(T,{r_1\over\lambda},{r_2\over \lambda}).
\end{split}
\end{equation*}
 Using \eqref{e:P-Lc-bullet_finite(2)} and \eqref{e:P-Lc-bullet_finite(4)} and increasing $c$ if necessary, we  deduce from the above  equality that
\begin{equation*}
 \begin{split}
 \big| I_{q_0,0,j_0,0}(T,r)+\sum_{j_1,j'_1,j''_1} {j'_0\choose  j'_1}{q_0 \choose j_1}
  {q_0-j_1 \choose j''_1}(-c_1)^{j'_0-j'_1}(-1)^{q_0-j_1-j''_1}\\
 \cdot I_{j_1, j'_0-j'_1, q_0+j_0+j'_0-j_1-j'_1-j''_1,q_0-j_1-j''_1}(T,r) -\big(\nu_{j_0,q_0}\big(T,B,{r_2\over \lambda},\tau\big)  -\nu_{j_0,q_0}\big(T,B,{r_1\over \lambda},\tau\big)\big)\big|\\
 \leq  c\lambda^{-1}+ {c \lambda^{-1}}\sum_{n=0}^\infty 2^{-n}\Mc^\tot( T,{r_2\over 2^n \lambda} )+c\Mc^\tot( \ddc T,{r_2\over \lambda} ).
\end{split}
\end{equation*}
As in the  proof of \eqref{e:P-Lc-bullet_finite(6)} we rewrite  this  inequality as follows:
\begin{equation}\label{e:P-Lc-bullet_finite(6)-mono-psh}
\begin{split}
 &\big|\Ic_1+\Ic_2+\Ic_3   -\big(\nu_{j_0,q_0}\big(T,B,{r_2\over \lambda},\tau\big)  -\nu_{j_0,q_0}\big(T,B,{r_1\over \lambda},\tau\big)\big)\big|\\
 &\leq c\lambda^{-1}+ {c \lambda^{-1}}\sum_{n=0}^\infty 2^{-n}\Mc^\tot( T,{r_2\over 2^n\lambda} )+c\Mc^\tot( \ddc T,{r_2\over \lambda} ),
 \end{split}
\end{equation}
where $\Ic_1,$ $\Ic_2$ and $\Ic_3$ are given in \eqref{e:P-Lc-bullet_finite(6)-mono}.
 
Repeating the argument   from  \eqref{e:P-Lc-bullet_finite(6bis)-mono} to  the end  of the  proof of Lemma 
\ref{L:basic-difference-nu_j,q_estimate-psh}, the  result  follows.
\endproof

Here is the main result of this section.

\begin{theorem}\label{T:nu_tot-monotone-psh} 
Let  $0< r_1<r_2\leq\bfr.$ Then there are  a  family  $\Dc=\{ d_{jq}\in\R:\ 0\leq j\leq k-p-q,\ 0\leq q\leq  k-l\}$ and   a constant $c>0 $ depending  on $r_1$ and $r_2$  such that for every positive plurisubharmonic  current $T$ on $\bfU$ belonging to the class $\widetilde\SH^{3,3}_p(\bfU,\bfW),$ 
the  following inequality  hold for $0\leq q\leq \upm:$
\begin{equation}\label{e:monotone-psh}
  \nu^\Dc_q\big(T,B,{r_1\over \lambda},\tau\big)\leq  \nu^\Dc_q(T,B,{r_2\over \lambda},\tau\big) + {c\over \lambda}+ {c\over \lambda}\sum_{n=0}^\infty \Mc^\tot(T,{r_2\over 2^n \lambda})+c\Mc^\tot(\ddc T,{r_2\over \lambda})
  \quad\text{for}\quad \lambda\gg 1. 
 \end{equation}
 Moreover,  the  following  two inequalities also hold:
\begin{equation} \label{e:nu-DC-good-psh}  \begin{split}
\nu^\Dc_\tot(T,B, r,\tau)&\leq  cr+c\sum_{n=0}^\infty 2^{-n}\Mc^\tot(T,r),\\
\Mc^\tot(T,r)&\leq cr+ c\nu^\Dc_\tot(T,B, r,\tau)
\quad\text{for}\quad 0<r\leq\bfr.\end{split}
\end{equation}
\end{theorem}
\proof We proceed as in the  proof of Theorem \ref{T:nu_tot-monotone}
using  Proposition \ref{P:nu_tot-monotone-psh}  instead of  Proposition \ref{P:nu_tot-monotone}. 
\endproof
\begin{corollary}\label{C:nu_tot-monotone-psh} 
 Let  $0< r_1<r_2\leq\bfr.$ Then there  is   a constant $c>0 $ depending  on $r_1$ and $r_2$  such that for every positive plurisubharmonic  current $T$ on $\bfU$ belonging to the class $\widetilde\SH^{3,3}_p(\bfU,\bfW),$ 
 and  every  $0<r\leq \bfr,$ 
we have 
\begin{equation*}
  \nu^\Dc_\tot(T,B,{r_1\over \lambda },\tau)\leq (1+c\lambda^{-1})  \nu^\Dc_\tot(T,B,{r_2\over \lambda },\tau)  +c\lambda^{-1}+
  c\sum_{n=1}^\infty {1\over 2^n\lambda}\Mc^\tot(T, {r_2\over  2^n\lambda} )+c \Mc^\tot(\ddc T,{r_2\over \lambda}).
 \end{equation*}
\end{corollary}
\proof  Applying  inequality \eqref{e:monotone-psh} to  $ q:=\upm$ yields
\begin{equation*}
  \nu^\Dc_\tot\big(T,B,{r_1\over \lambda},\tau\big)\leq  \nu^\Dc_\tot(T,B,{r_2\over \lambda},\tau\big) + {c\over \lambda}+ {c\over \lambda}\sum_{n=0}^\infty 2^{-n}  \Mc^\tot(T,{r_2\over 2^n\lambda})+c\Mc^\tot(\ddc T,{r_2\over \lambda})
  \quad\text{for}\quad \lambda\gg 1. 
 \end{equation*}
This,  combined with the second inequality of \eqref{e:nu-DC-good-psh}, gives the  result.
\endproof

\begin{proposition}\label{P:nu_tot-sum-psh}  There is  a constant $c_{10}>0$ such that for every positive plurisubharmonic   current $T$ on $\bfU$ belonging to the class $\widetilde\SH^{3,3}_p(\bfU,\bfW),$ 
and  every  $0<r\leq \bfr,$ 
we have 
\begin{equation*}
  \sum_{n=0}^\infty {1\over 2^n}\Mc^\tot(T,{r\over 2^n})\leq c \Nc^\tot(T,r) +cr.
 \end{equation*}
\end{proposition}
\proof Consider  $r_1:={\bfr\over 2}$ and $r_2=\bfr.$ 
Applying the  second inequality of \eqref{e:nu-DC-good-psh} yields that
\begin{eqnarray*}
   \sum_{n=0}^\infty {1\over 2^n}\Mc^\tot(T,{r\over 2^n})\leq cr
   + c \sum_{n=0}^\infty {1\over 2^n}\Mc^\tot(\ddc T,{r\over 2^n})+c \sum_{n=0}^\infty {1\over 2^n}\nu^\Dc_\tot(T,B, {r\over 2^n},\tau).
\end{eqnarray*}
Since $\ddc T\in\CL_p^{1,1}(\bfU,\bfW),$  it follows   from  Theorem   \ref{T:nu_tot-monotone} that
\begin{eqnarray*}
  \sum_{n=0}^\infty {1\over 2^n}\Mc^\tot(\ddc T,{r\over 2^n})&\leq& cr+c  \sum_{n=0}^\infty {1\over 2^n}\nu^\Dc_\tot(\ddc T,B,{r\over 2^n},\tau)\\
  &\leq & 2cr+c  \sum_{n=0}^\infty {1\over 2^n}\nu^\Dc_\tot(\ddc T,B,r,\tau)\leq c'r +c\Mc^\tot(\ddc T,r).
\end{eqnarray*}
On the other hand,
applying   Corollary  \ref{C:nu_tot-monotone-psh} for $r_1={r_2\over 2}$  and using   the elementary inequality $1+t\leq e^t$ for $t\geq 0,$
yield   that  
\begin{equation*}
  \nu^\Dc_\tot(T,B,{r_2\over 2\lambda },\tau)\leq   e^{c\lambda^{-1}}\nu^\Dc_\tot(T,B,{r_2\over \lambda },\tau)  +c\lambda^{-1}+
  c\sum_{m=1}^\infty {1\over 2^m\lambda}\Mc^\tot(T, {r_2\over  2^m\lambda} )+c \Mc^\tot(\ddc T,{r_2\over \lambda}).
 \end{equation*}
For each $n\geq 0$ we apply  this   inequality  for $\lambda$ such that ${r_2\over \lambda}={r\over  2^n}.$
Consequently, we get a constant $c>0$ independent of $T$ and $n,r$  such that
\begin{eqnarray*}
  \nu^\Dc_\tot(T,B, {r\over 2^n },\tau)&\leq &c \nu^\Dc_\tot(T,B,r,\tau) +cr\sum_{m=1}^\infty \min{(m,n)} 2^{-m}\nu_\tot^\Dc(T,B,  {r\over 2^m },\tau)  +cr\\
  &+&c \sum_{m=1}^n\Mc^\tot(\ddc T,{r\over 2^m}).
  \end{eqnarray*}
Since  $\nu^\Dc_\tot(T,B,r,\tau)\lesssim \Mc^\tot(T,r),$ there is a constant $c'>0$  independent of $T$ and $n,r$  such that
  \begin{eqnarray*}
   \nu^\Dc_\tot(T,B, {r\over 2^n },\tau)&\leq &  c'\big(r+\Mc^\tot(T,r)+ \Mc^\tot(\ddc T,r)+r\sum_{m=1}^\infty \min{(m,n)}2^{-m}\nu_\tot^\Dc(T,B,  {r\over 2^m },\tau)\\
   &+& \sum_{m=1}^n\Mc^\tot(\ddc T,{r\over 2^m})\big).
 \end{eqnarray*}
Hence,  we infer that
\begin{multline*}
 \sum_{n=0}^\infty {1\over 2^n}\nu^\Dc_\tot(T,B, {r\over 2^n},\tau)\leq 
 c'r\sum_{n=0}^\infty\sum_{m=1}^\infty  \min{(m,n)}2^{-m-n}\nu_\tot^\Dc(T,B,  {r\over 2^m },\tau) +c'r+c'\Mc^\tot(T,r)\\
 + c'\sum_{n=0}^\infty\sum_{m=1}^n 2^{-n}\Mc^\tot(\ddc T,{r\over 2^m}).
\end{multline*}
Since  $\sum_{n=0}^\infty \min{(m,n)}2^{-n}\leq m2^{-m+1}+\sum_{n=0}^\infty n2^{-n}<  1+\sum_{n=0}^\infty<\infty,$  the first double sum on the  RHS 
is dominated  by a constant times the LHS. On othe other hand,  the second double sum on the  RHS 
is dominated  by a constant times $\sum_{m=0}^\infty 2^{-m} \Mc^\tot(\ddc T,{r\over 2^m}),$
which is, by Theorem \ref{T:nu_tot-monotone},  bounded by a constant times $ r+\Mc^\tot(\ddc T,r)  . $   Taking into account the factor $r$  in front of  this  double sum, we   get for a constant $c''\gg 1$  that
$$
\sum_{n=0}^\infty {1\over 2^n}\nu^\Dc_\tot(T,B, {r\over 2^n},\tau)\leq c''r+c''\Nc^\tot(T,r).
$$
This  completes the proof.
\endproof
 \begin{proposition}\label{P:Lc-finite-psh}
 There is a   constant $c_{11}>0$ such that for  
  every positive plurisubharmonic  current $T$  belonging to the class $\widetilde\SH^{3,3}_p(\bfU,\bfW),$ 
   we have
   $\Mc_j(T,r)<c_{11}$ for $0\leq j\leq \upm$ and $0 <r \leq \bfr.$
 \end{proposition}
\proof  We apply  Theorem  \ref{T:nu_tot-monotone-psh}  to $r_1:={\bfr\over 2}$ and $r_2=\bfr.$  \endproof
 We close  the section with the  following synthesis.
 \begin{corollary}
  \label{C::nu_tot-monotone-psh-bis}
Let  $0< r_1<r_2\leq\bfr.$ Then there are  a  family  $\Dc=\{ d_{jq}\in\R:\ 0\leq j\leq k-p-q,\ 0\leq q\leq  k-l\}$ and   a constant $c>0 $ depending  on $r_1$ and $r_2$  such that for every positive plurisubharmonic  current $T$ on $\bfU$ belonging to the class $\widetilde\SH^{3,3}_p(\bfU,\bfW),$ 
the  following inequality  hold for $0\leq q\leq \upm:$
\begin{equation*} \begin{split}
  \nu^\Dc_q\big(T,B,{r_1\over \lambda},\tau\big)&\leq  \nu^\Dc_q(T,B,{r_2\over \lambda},\tau\big) + {c\over \lambda}  +c\Nc^\tot( T,{r_2\over \lambda})
  \quad\text{for}\quad \lambda\gg 1,\\
\nu^\Dc_\tot(T,B, r,\tau)&\leq  cr+c\Nc^\tot(T,r)\quad\text{for}\quad 0<r\leq\bfr.
  \end{split}
 \end{equation*}
 Moreover,  the  following   inequality also holds for $ \lambda\gg 1:$
 \begin{equation*}
  \nu^\Dc_\tot(T,B,{r_1\over \lambda },\tau)\leq (1+c\lambda^{-1})  \nu^\Dc_\tot(T,B,{r_2\over \lambda },\tau)  +c\lambda^{-1}+
  c \Nc^\tot(T,{r_2\over \lambda}).
 \end{equation*}
 \end{corollary}
 \proof
Using  
Proposition \ref{P:nu_tot-sum-psh}, the first and the  second inequalities (resp. the third one)  
follow  from
Theorem \ref{T:nu_tot-monotone-psh}  (resp. Corollary \ref{C:nu_tot-monotone-psh}).
\endproof

\section{Positive plurisubharmonic currents and finiteness of the mass indicators $\Kc_{j,q}$ and $\Lc_{j,q}$}
\label{S:finitness-Kc-Lc}

 \subsection{Preliminary estimates}
 \label{SS:preli-estimates-psh}
\begin{lemma}\label{L:ddc-difference} There is a constant $c>0$  such that  for every
$j$ with $\lowm\leq j\leq \upm,$  and every $m$ with $0\leq m\leq j,$ and every positive plurisubharmonic  current $T$   in the class  $\widetilde\SH^{3,3}_p(\bfU,\bfW),$
 there exists  a function  $(0,\bfr]\ni r\to \tilde r$ (depending on $T$)  with ${r\over 2}\leq \tilde r\leq r$
such  that the  following  two inequalities hold for  $0<s<r\leq \bfr:$
 \begin{multline*} \big|\int_{\tilde s}^{\tilde r} \big({1\over t^{2(k-p-j)}} -{1\over r^{2(k-p-j)}}\big)2tdt  \int_{\Tube(B,t)}  \big( \ddc (\tau_*T)- \tau_*(\ddc T)\big)\wedge \pi^*(\omega^{j-m})\wedge \beta^{k-p-j+m-1}\big|\\
 \leq cr^{2m+1}\Mc^\tot(T,r),
 \end{multline*}

\end{lemma}
\proof 
By Lemma  \ref{L:Ic-0}, there is  a function  $(0,\bfr]\ni r\to \tilde r$
with the  following  two properties.
\begin{itemize} \item  ${r\over 2}\leq \tilde r\leq r;$
 \item   there is a constant $c$ independent of $T$  such that for $\nu\in \{1,2\}$ and for all $1\leq \ell\leq \ell_0$ and for all $0<r\leq \bfr:$
  \begin{equation*}
   {1\over  r^{2(k-p-j)} } \big|  \int_{\partial_\hor\Tube(B,t) }( D^{(\ell)}_{\nu0}S^{(\ell)}\wedge \Phi^{(\ell)}_{\nu 0}) \big|\leq  cr^{2m+1}\Mc^\tot(T,r).
  \end{equation*}
\end{itemize}
Now let $0<s<r\leq \bfr.$ The two $\bullet$ imply that for all $t$ with $0<t<r,$ 
\begin{equation*}\big|\int_{\partial_\hor\Tube(B,\tilde t) }\chi(\tilde t)( D^{(\ell)}_{\nu 0}S^{(\ell)}\wedge \Phi^{(\ell)}_{\nu 0})(y) \big|\leq  cr^{2m+2}\Mc^\tot(T,r),
\end{equation*}
where  $\chi$ is  either the function $\chi_1$ or the function $\chi_2$ given in \eqref{e:chi_1-chi_2}.
Using  this, we apply  Lemma \ref{L:Ic-1-and-2}  to equality \eqref{e:Stokes-ddc-difference-bis} for $r_1:=\tilde s$ and $r_2:=\tilde r.$
Hence, we get from \eqref{e:Stokes-ddc-difference-bis} and Lemma \ref{L:Ic-1-and-2} that
\begin{equation*}
  \big|\int_{\tilde s}^{\tilde r}\chi(t) \Ic^{(\ell)}_\nu(t)dt\big|\leq  cr^{2m+1}\Mc^\tot(T,r).
 \end{equation*}
 On the other hand, applying  \eqref{e:inequal-ddc-difference-bis} to ${r\over 2^n}$ $(n\in\N$) instead of $r$ and summing  the obtained inequalities  yields that
\begin{equation*}
  \int_0^r  \chi(t)\big| \langle  \ddc [(\tilde\tau_\ell)_* S^{(\ell)}] -(\tilde\tau_\ell)_*(\ddc  S^{(\ell)}),\Phi^{(\ell)}  \rangle_{\tilde \tau_\ell (\H_t)}
 -  \Ic^{(\ell)}_1(t) -\Ic^{(\ell)}_2(t)\big|dt
  \leq 
cr^{2m+1}\Mc^\tot(T,r) .
 \end{equation*}
 This, combined with the previous inequality, implies that
 \begin{equation*}
  \big|\int_{\tilde s}^{\tilde r} \chi(t)\cdot \langle  \ddc [(\tilde\tau_\ell)_* S^{(\ell)}] -(\tilde\tau_\ell)_*(\ddc  S^{(\ell)}),\Phi^{(\ell)}  \rangle_{\tilde \tau_\ell (\H_t)} \big|\leq  cr^{2m+1}\Mc^\tot(T,r).
 \end{equation*}
 Summing  this inequality   for $1\leq \ell\leq \ell_0,$ we get that
 \begin{equation*}
  \big|\int_{\tilde s}^{\tilde r} \chi(t)\cdot \sum_{\ell=1}^{\ell_0}\langle  \ddc [(\tilde\tau_\ell)_* S^{(\ell)}],  \Phi^{(\ell)}  \rangle_{\tilde \tau_\ell (\H_t)}     -\int_{\tilde s}^{\tilde r}\chi(t)\cdot \sum_{\ell=1}^{\ell_0}\langle(\tilde\tau_\ell)_*(\ddc  S^{(\ell)}),  \Phi^{(\ell)}  \rangle_{\tilde \tau (\H_t)}\big|\leq  cr^{2m+1}\Mc^\tot(T,r).
 \end{equation*}
  By the two last equalities of Lemma \ref{L:partition},   the last  inequality is  rewritten as
  \begin{equation*}  |I_\chi|\leq  cr^{2m+1}\Mc^\tot(T,r),
  \end{equation*}
where
  \begin{equation*}
  I_\chi:=\int_{\tilde s}^{\tilde r} \chi(t)\cdot \langle  \ddc (\tau_* T),  \Phi\rangle_{\Tube(B,t)}     -\int_{\tilde s}^{\tilde r}\chi(t)\cdot \langle \tau_*(\ddc  T),  \Phi  \rangle_{\Tube(B,t)}.
 \end{equation*}
 Since  the integral in the   inequality of the lemma is equal to $I_{\chi_2}-I_{\chi_1},$ 
  the lemma follows from the last  estimate.

\endproof

\begin{lemma}\label{L:ddc-difference-m-1} 
There is a constant $c>0$  such that  for every
$j$ with $\lowm\leq j\leq \upm,$ and every $m$ with $1\leq m\leq j,$  and every positive plurisubharmonic  current $T$   in the class  $\widetilde\SH^{3,3}_p(\bfU,\bfW),$
 the  following properties  hold for every $0<s<r\leq \bfr:$
 \begin{multline*} \Big|\int_{\tilde s}^{\tilde r} \big({1\over t^{2(k-p-j)}} -{1\over r^{2(k-p-j)}}\big)2tdt  \int_{\Tube(B,t)}   \tau_*(\ddc T)\wedge \pi^*(\omega^{j-m})\wedge \beta^{k-p-j+m-1} \Big|\\
\leq cr^{2m}\Mc^\tot(  \ddc T,r ).
 \end{multline*}

\end{lemma}

\proof  
Consider the positive  closed $(p+1,p+1)$-current $S:=\ddc T$ on $\bfU.$
Observe that $T$  belongs to the class $\widetilde\CL^{1,1}_p(\bfU,\bfW).$
Applying Proposition  \ref{P:Lc-finite-closed} yields that
   $\Mc_j(S,r)<c_{9}$ for $\lowm\leq j\leq \upm$ and $0 <r \leq \bfr.$
   Using  formula \eqref{e:Lelong-numbers} we rewrite the integral in assertion (1) as
 $$\int_0^r \big({1\over t^{2(k-p-j)}} -{1\over r^{2(k-p-j)}}\big)2t^{2(k-p-j)-1} \nu(\ddc T,B,t,\tau)dt.$$
It  follows from the above   inequality that $| \nu(\ddc T,B,t,\tau)|\leq c_9.$ Therefore, we infer that the absolute value of the above integral is bounded from above by
$$
c_9\int_0^r \big({2\over t} +{1\over r^{2(k-p-j)}}\big)2t^{2(k-p-j)-1} <c ,$$
where $c>0$ is a constant independent of $T$ and  $r.$
\endproof

As  an immediate consequence of  Lemmas \ref{L:ddc-difference} and \ref{L:ddc-difference-m-1}, we obtain  the following result.

\begin{corollary}\label{C:ddc-estimate-m-1} There is a constant $c>0$  such that  for every
$j$ with $\lowm\leq j\leq \upm,$ and every $m$ with $1\leq m\leq j,$ and every positive plurisubharmonic  current $T$   in the class  $\widetilde\SH^{3,3}_p(\bfU,\bfW),$
 the  following inequality  holds:
 \begin{multline*} \Big |\int_{\tilde s}^{\tilde r} \big({1\over t^{2(k-p-j)}} -{1\over r^{2(k-p-j)}}\big)2tdt  \int_{\Tube(B,t)}   \ddc (\tau_*T)\wedge \pi^*(\omega^{j-m})\wedge \beta^{k-p-j+m-1}\Big|\\
 \leq cr^{2m}\Nc(  T,r ).
 \end{multline*}

\end{corollary}


Fix  an open neighborhood $\bfW$ of $\partial B$ in $X$ with $\bfW\subset \bfU.$
Recall  the class $\widetilde\SH^{3,3}_p(\bfU,\bfW)$  given in Definition
 \ref{D:sup-bis}.
 
   For $0<r\leq \bfr$ and $0 \leq q\leq k-l$ and $0\leq j\leq k-p-q,$ consider  following  global mass indicator
  \begin{equation}\label{e:mass-indicators-log}
 \Lc_{j,q}(T,r) :=     \int_0^r {2dt\over t^{2q-1}}    \big( \int_{\Tube(B,t)}(\ddc T)^\hash\wedge
 \pi^*\omega^j\wedge  (\beta+c_1t^2\pi^*\omega)^{k-p-j-1}\big).
 \end{equation}
 Since  $\beta+c_1t^2\pi^*\omega$ is a positive form on $\Tube(B,t),$ it follows that $\Lc_{j,q}(T,r)\geq 0.$
  \begin{lemma}\label{L:Lc-vs-analog}
  There is a constant  $c>0$ such that for every positive plurisubharmonic   current $T$ on $\bfU$ belonging to the class $\widetilde\SH^{3,3}_p(\bfU,\bfW),$  and  every $0<r\leq \bfr$ and $j,q\geq 0$ as  above, we have
  \begin{equation*}
   \big|  \int_0^r {2dt\over t^{2q-1}}   \big( \int_{\Tube(B,t)}\tau_*(\ddc T)\wedge
 \pi^*\omega^j\wedge  (\beta+c_1t^2\pi^*\omega)^{k-p-j-1}\big)-\Lc_{j,q}(T,r) \big|\leq cr\Mc^\tot(\ddc T,r).
   \end{equation*}
  
  \end{lemma}
 
 \proof  We argue as  in the proof of Lemma
 \ref{L:Mc-hash-vs-nu}.   \endproof

  \begin{lemma}\label{L:Lc-vs-analog-bis}
  There is a constant  $c>0$ such that for every positive plurisubharmonic   current $T$ on $\bfU$ belonging to the class $\widetilde\SH^{3,3}_p(\bfU,\bfW),$  and  every $0<r\leq \bfr$ and $j,q\geq 0$ as  above, we have
  \begin{multline*}
   \big|  \int_0^r {2dt\over t^{2q-1}}   \big( \int_{\Tube(B,t)}\tau_*(\ddc T)\wedge
 \pi^*\omega^j\wedge  (\beta+c_1t^2\pi^*\omega)^{k-p-j-1}\big)\\
 - \int_0^r {2dt\over t^{2q-1}}   \big( \int_{\Tube(B,t)}\ddc (\tau_*T)\wedge
 \pi^*\omega^j\wedge  (\beta+c_1t^2\pi^*\omega)^{k-p-j-1}\big)
 \big|\leq cr\Nc( T,r).
   \end{multline*}
  
  \end{lemma}
 
 \proof 
First we argue as  in the  proof of Lemma \ref{L:ddc-difference-lambda} in order to show that the  expression on the LHS is
dominated by a constant times  $ r\sum_{n=0}^\infty \Mc^\tot(T, {r\over 2^n})+r^2+r\Mc^\tot(\ddc T,r) .$ Second   we apply Proposition 
\ref{P:nu_tot-sum-psh} in order to see that the last expression is  in turn dominated by  a constant times  $r\Nc(T,r).$
 \endproof

  \begin{lemma}\label{L:Lc-vs-analog-bis'}
  There is a constant  $c>0$ such that for every positive plurisubharmonic   current $T$ on $\bfU$ belonging to the class $\widetilde\SH^{3,3}_p(\bfU,\bfW),$  and  every $0<r\leq \bfr$ and $j,q\geq 0$ as  above, we have
  \begin{multline*}
   \big|  \int_0^r {2dt\over t^{2q-1}}   \big( \int_{\Tube(B,t)}\tau_*(\ddc T)\wedge
 \pi^*\omega^j\wedge  \beta^{k-p-j-1}\big)\\
 - \int_0^r {2dt\over t^{2q-1}}   \big( \int_{\Tube(B,t)}(\ddc T)^\hash\wedge
 \pi^*\omega^j\wedge  \beta^{k-p-j-1}\big)
 \big|\leq cr\Mc^\tot( \ddc T,r).
   \end{multline*}
  
  \end{lemma}
 
 \proof  We argue as  in the proof of Lemma
 \ref{L:Mc-hash-vs-nu}.   \endproof
 
 \begin{lemma}\label{L:Lc-vs-analog-bisbis}  There is a constant  $c>0$ such that for every $j$ with $\lowm^+\leq j\leq \upm^+,$ and for every positive plurisubharmonic   current $T$ in the class $\widetilde\SH^{3,3}_p(\bfU,\bfW),$  
  and  for every $0<r\leq \bfr,$ we have
  \begin{multline*}
   \big|\Lc_{j,q}(T,r)-  \int_0^{ r} {2dt\over t^{2q-1}} \big(  \int_{\Tube(B,t)}  
   (\ddc T)^\hash\wedge \pi^*(\omega^{j})\wedge \beta^{k-p-j-1}\big)  \big|\\
   \leq 
   c\sum_{j'=1}^{\upm-j}  \Lc_{j+j',q-j'}(T,r)
   +
   c\Nc(T,r).
  \end{multline*}
 \end{lemma}
\proof  Write
$$
  (\beta+c_1\pi^*\omega)^{k-p-j-1}-\beta^{k-p-j-1}=
  \sum_{j'=1}^{k-p-j-1} (-1)^{j'+1} {k-p-j-1\choose j'} c_1^{j'} (\beta+c_1\pi^*\omega)^{k-p-j-j'-1} (\pi^*\omega)^{j'}.$$
  Using this and \eqref{e:mass-indicators-log}, we argue as  in the proof of Lemma
 \ref{L:Mc-hash-vs-nu}. 
\endproof

  \begin{lemma}\label{L:Lc-vs-analog-bis-new}
  There is a constant  $c>0$ such that for every positive plurisubharmonic   current $T$ on $\bfU$ belonging to the class $\widetilde\SH^{3,3}_p(\bfU,\bfW),$  and  every $0<r\leq \bfr$ and $j,q\geq 0$ as  above, we have
  \begin{multline*}
   \big|  \int_0^r \big({1\over t^{2q}} -  {1\over r^{2q}}\big)2tdt   \big( \int_{\Tube(B,t)}\tau_*(\ddc T)\wedge
 \pi^*\omega^j\wedge  \beta^{k-p-j-1}\big)\\
 - \int_0^r  \big({1\over t^{2q}} -  {1\over r^{2q}}\big)2tdt \big( \int_{\Tube(B,t)}\ddc (\tau_*T)\wedge
 \pi^*\omega^j\wedge  \beta^{k-p-j-1}\big)
 \big|\leq cr\Nc( T,r).
   \end{multline*}
  
  \end{lemma}
 
 \proof We argue as in the proof of 
  Lemma \ref{L:Lc-vs-analog-bis}.  \endproof

  \begin{lemma}\label{L:Lc-vs-analog-bis'-new}
  There is a constant  $c>0$ such that for every positive plurisubharmonic   current $T$ on $\bfU$ belonging to the class $\widetilde\SH^{3,3}_p(\bfU,\bfW),$  and  every $0<r\leq \bfr$ and $j,q\geq 0$ as  above, we have
  \begin{multline*}
   \big|  \int_0^r   \big({1\over t^{2q}} -  {1\over r^{2q}}\big)2tdt  \big( \int_{\Tube(B,t)}\tau_*(\ddc T)\wedge
 \pi^*\omega^j\wedge  \beta^{k-p-j-1}\big)\\
 - \int_0^r  \big({1\over t^{2q}} -  {1\over r^{2q}}\big)2tdt  \big( \int_{\Tube(B,t)}(\ddc T)^\hash\wedge
 \pi^*\omega^j\wedge  \beta^{k-p-j-1}\big)
 \big|\leq cr\Mc^\tot( \ddc T,r).
   \end{multline*}
  
  \end{lemma}
 
 \proof We argue as in the proof of
  Lemma \ref{L:Lc-vs-analog-bis'}.   \endproof
 
 \begin{lemma}\label{L:Lc-vs-analog-bisbis-new}  There is a constant  $c>0$ such that for every $j$ with $\lowm^+\leq j\leq \upm^+,$ and for every positive plurisubharmonic   current $T$ in the class $\widetilde\SH^{3,3}_p(\bfU,\bfW),$  
  and  for every $0<r\leq \bfr,$ we have
  \begin{multline*}
   \big|\Lc_{j,q}(T,r)-  \int_0^{ r} \big({1\over t^{2q}} -  {1\over r^{2q}}\big)2tdt  \big(  \int_{\Tube(B,t)}  
   (\ddc T)^\hash\wedge \pi^*(\omega^{j})\wedge \beta^{k-p-j-1}\big)  \big|\\
   \leq 
   c\sum_{j'=1}^{\upm-j}  \Lc_{j+j',q-j'}(T,r)
   +
   c\Nc(T,r).
  \end{multline*}
 \end{lemma}
\proof  We argue as in the proof of 
 Lemma \ref{L:Lc-vs-analog-bisbis}. 
 Since  the proof is  not difficult, we leave it   to the interested reader.   
\endproof

\subsection{Finiteness of the mass indicators $\Kc_{j,q} $ and   $\Lc_{j,q} $ }
 
\begin{lemma}\label{L:Lelong-smooth-forms-psh}
 Let $T$ be a positive plurisubharmonic   $\Cc^2$-smooth $(p,p)$-form   on $\bfU.$
 Then  for every $\lowm \leq j\leq \upm,$  we have  $\nu_j(T,B,\tau)=0$ if $j\not=l-p$ and $\nu_j(T,B,\tau)\geq 0$ if $j=l-p.$ 
\end{lemma}
\proof
First consider the  case   $j\not=l-p.$  As  $\lowm \leq j\leq \upm,$ we have   $j>l-p ,$  and hence   $k-p-j<k-l.$
Then by Theorem  \ref{T:Lelong-Jensen-smooth} (1), $\nu_j(T,B,\tau)=0.$

Now consider the case $j=l-p.$ So $j=\lowm.$  Since $\tau$ is  strongly  admissible  $d\tau|_{\overline B}$ is  $\C$-linear,
it follows from  the positivity of $T$ on $\bfU$ that  $(\tau_*T)|_{\overline B}$ is also a positive form.
Hence, by Theorem  \ref{T:Lelong-Jensen-smooth} (1)  again, $\nu_j(T,B,\tau)\geq 0.$
\endproof

 \begin{theorem}\label{T:Lc-finite-psh}
   There is  a constant $c_{10}>0$ such that for every positive plurisubharmonic   current $T$ on $\bfU$ belonging to the class $\widetilde\SH^{3,3}_p(\bfU,\bfW),$  we have 
   \begin{equation}\label{e:Lc-finite-psh-main}\begin{split}
    \Kc_{j,q}(T,r)&\leq  c_{10} (\nu_\tot(T,B,r,\tau)+\nu_\tot(\ddc T,B,r,\tau)),\\
    \Lc_{j,q}(T,r)&\leq  c_{10} (\nu_\tot(T,B,r,\tau)+\nu_\tot(\ddc T,B,r,\tau))
   \end{split}
   \end{equation}
   for $0\leq q\leq  k-l$ and   $0\leq   j\leq k-p-q$ and for $0<r\leq \bfr.$ 
   In particular,
   $$\Kc_{j,q}(T,\bfr)<c_{10}\quad\text{and}\quad  \Lc_{j,q}(T,\bfr)<c_{10} .
   $$
   \end{theorem}
\proof

Since  the masses of $T$ and of $\ddc T$ on $\bfU$ is $\leq 1,$ there is a constant $c$ independent of $T$ such that  $0\leq \nu_\tot(T,B,\bfr,\tau)\leq c$ and  $0\leq \nu_\tot(\ddc T,B,\bfr,\tau)\leq c.$ Therefore,
the last  two  inequalities  follow from  the  first ones.  So  we only need to prove  the first  two inequalities
\eqref{e:Lc-finite-psh-main}.
The proof of \eqref{e:Lc-finite-psh-main} is  divided into three steps.
We indicate how to adapt the proof of Theorem  \ref{T:Lc-finite} in the  present context.
The proof is also divided into three steps.

\noindent {\bf Step 1:} {\it The case  $q=0.$}

In this case  there  is  no factor $\hat\alpha$  appearing in $\Kc_{j,0}(T,\bfr)$ Since  the  forms  $\omega$ and $\hat\beta$ are  positive smooth  and  the mass of $T$  on $\bfU$ is $\leq 1,$   there is a constant $c_{10}$ such that 
$$
\Kc_{j,0}(T,\bfr)=\sum_{\ell=1}^{\ell_0}\int\limits_{(\Tube(B,\bfr)\cap \U_\ell) \setminus V  } (\pi^*\theta_\ell)\cdot (\tau_\ell)_*( T|_{\bfU_\ell})\wedge \pi^*\omega^j\wedge \hat\beta^{k-p-j} <c\|T\|_{\bfU}<c_{10}.
$$
Similarly, we obtain 
$$
 \Lc_{j,0}(T,r) :=     \int_0^r 2tdt    \big( \int_{\Tube(B,t)}(\ddc T)^\hash\wedge
 \pi^*\omega^j\wedge  (\beta+c_1t^2\pi^*\omega)^{k-p-j-1}\big)<c\|\ddc T\|_{\bfU}<c_{10}.
 $$
This proves the theorem for $q=0,$ and hence concludes Step 1.

 The general strategy is  to prove  the proposition by  increasing induction on $q$ 
 with $0\leq q\leq k-l.$ But the  induction  procedure  is  more  complicated  than that of Theorem  \ref{T:Lc-finite} since  
  a  double  induction is   needed  in the  present context. 
 In the proof $\bfr$ is  a  fixed   but sufficiently small  constant.
Fix  $0\leq q_0\leq k-l.$ 
Suppose   that \eqref{e:Lc-finite-psh-main} is true  for all $q,j$ with $q<q_0.$ 
We need to show that   it is  also true  for all $q,j$ with   $q\leq q_0.$ 

Recall  from \eqref{e:P-Lc-bullet_finite(0)} the   mass indicators $ \Kc^{\pm}_{j,q}(T,s,r)$ and $\Kc_{q}(T,s,r).$
We also introduce the following new mass indicators:
\begin{equation}\label{e:P-Lc-bullet_finite(0)-psh-bullet}
 \Kc^\bullet_{q}(T,s,r):=\sum\limits^{\bullet}_{q}\Kc_{j',q'}(T,s,r),
\end{equation}
where $\sum\limits^{\bullet}_{ q}$ means that the sum is taken over all $(j',q')$ such that either ($q'<q$ and $0\leq j'\leq k-p-q'$) or ($q'=q$ and $0\leq j'<k-p-q'$).  
So we have 
\begin{equation}\label{e:P-Lc-bullet_finite(0)-psh}
 \Kc_{q}(T,s,r):=\Kc_{q-1}(T,s,r)+\sum\limits_{j=0}^{k-p-q}\Kc_{j,q}(T,s,r)= \Kc^\bullet_{q}(T,s,r)+\Kc_{k-p-q,q}(T,s,r).
\end{equation}
We may assume without loss of generality that $T$ is  $\Cc^3$-smooth
and  let $s,r\in  [0,\bfr)$ with $s<r.$

{\it  Set $m_0:= k-p-q,$  $m_1:=m_0-1.$ 
In the first induction we will prove that there is a constant $c_{10}$ independent of $T$ and $r$
such that 
\begin{equation}\label{e:first-ind}
\Kc^\bullet_{q}(T,r)
\leq c_{10}\Nc^\bullet_q(T,r)\qquad\text{and}\qquad \Lc_{j,q}(T,r)\leq c_{10}\Nc^\bullet_q(T,r),
\end{equation}
for every $0\leq q\leq k-l,$ $j\geq 0$  with $j\leq m_1,$  and for every $0<r\leq \bfr.$
Here,
\begin{eqnarray*}
\Nc^\bullet_q(T,r)&:=&r+\Kc_{q-1}(T,r)+ \Lc_{q-1}(T,r)+\Nc(T,r)+\sum_{j=0}^{m_1} |\nu_{j,q}(T,r)|,\\
\Lc_{q}(T,r)&:=&\sum_{j,q':\  q'\leq q\quad\text{and}\quad j+q'\leq k-p} \Lc_{j,q'}(T,r).
\end{eqnarray*}

The proof of \eqref{e:first-ind} will be completed in Steps 2 and 3 below.
}

\noindent {\bf Step 2:} {\it  Let $q_0:=q$ and define $m_0$ and $m_1$ as  above using $q_0$ instead $q.$ There is a constant $c_{10}>0$ such that for every  $j_0,q_0\geq 0$ with $j_0\leq  m_1$ and   every $0<r\leq \bfr,$
\begin{equation}\label{e:P-Lc-bullet_finite(8)-psh}
\begin{split}
I^\hash_{q_0,0,j_0,0}(T,r)&\leq  c_{10}\big(|\nu_{j_0,q_0}(T,B,r,\tau)|+r +r^{1\over 4}\Nc(T,r)+r^{1\over 4}\Kc^{+}_{j_0,q_0}(T,r) +r^{1\over 4} \Kc_{q_0}^\bullet(T,r)\\
&+  \sqrt{\Kc^\bullet_{q_0}(T,r)} \sqrt{ \Kc^-_{j_0,q_0}(T,r)  }\big).
\end{split}
\end{equation}
where the expression on the LHS is  given by \eqref{e:I_bfj}  (see also Remark \ref{R:I_bfj}).} 
 
Let $0\leq j_0\leq \min(\upm, k-p-q_0).$ Set 
 $j'_0:=k-p-q_0-j_0\geq 0$ and $m_0:=k-p-q_0.$
Suppose that $j'_0\geq 1.$

Note that
\begin{equation*}
 \ddc [(\tau_*T)\wedge \pi^*\omega^{j_0}\wedge\beta^{j'_0}]
=\ddc (\tau_*T)\wedge \pi^*\omega^{j_0}\wedge\beta^{j'_0} .
\end{equation*}
Applying Theorem  \ref{T:Lelong-Jensen-smooth}  to  $\tau_*T\wedge \pi^*(\omega^{j_0})\wedge \beta^{j'_0}$  and using  the  above  equality, we get for $0<r\leq\bfr$ that
\begin{equation}\label{e:P-Lc-bullet_finite(1)-psh}
\begin{split}
&{1\over  r^{2q_0}}\int_{\Tube(B,r)}\tau_*T\wedge \pi^*(\omega^{j_0})\wedge \beta^{k-p-j_0}
-\lim_{s\to 0+}{1\over  s^{2q_0}}\int_{\Tube(B,s)}\tau_*T\wedge \pi^*(\omega^{j_0})\wedge \beta^{k-p-j_0}\\
&= \int_{\Tube(B,r)}\tau_*T\wedge \pi^*(\omega^{j_0})\wedge \beta^{j'_0}\wedge\alpha^{q_0}+\Vc\big(\tau_*T\wedge \pi^*(\omega^{j_0})\wedge \beta^{j'_0},r\big) \\
&+\int_0^r \big({1\over t^{2q_0}} -{1\over r^{2q_0}}\big)2tdt   \int_{\Tube(B,t)}   \ddc (\tau_* T)\wedge \pi^*(\omega^{j_0})\wedge \beta^{q_0+j'_0-1}.  
\end{split}
\end{equation}
If  $j'_0\geq 1,$
then by  Corollary \ref{C:ddc-estimate-m-1} there is a  constant $c$ independent of $T$ and $r$   such that
\begin{equation}\label{e:P-Lc-bullet_finite(1)-bis-psh}\Big|\int_0^r \big({1\over t^{2q_0}} -{1\over r^{2q_0}}\big)2tdt   \int_{\Tube(B,t)}   \ddc (\tau_* T)\wedge \pi^*(\omega^{j_0})\wedge \beta^{q_0+j'_0-1}\Big|\leq cr^{2j'_0}\Nc(T,r).
 \end{equation}
Moreover, if  $j_0\geq 0,$ then  by Theorem \ref{T:vertical-boundary-terms} we have the following estimate independently of $T:$
\begin{equation}\label{e:P-Lc-bullet_finite(2)-psh}
\Vc\big(\tau_*T\wedge \pi^*(\omega^{j_0})\wedge \beta^{j'_0},s,r\big)=O(r).
\end{equation}
Therefore,
when  $s\to 0+,$  applying  Lemma \ref{L:Lelong-smooth-forms-psh} and Theorem   \ref{T:Lelong-Jensen-smooth} (1)  yields that
\begin{equation}\label{e:P-Lc-bullet_finite(3)-psh}
 \begin{split}  &{1\over  r^{2q_0}}\int_{\Tube(B,r)}\tau_*T\wedge \pi^*(\omega^{j_0})\wedge \beta^{k-p-j_0}
 -\lim_{s\to 0}{1\over s^{2q_0}} \int_{\Tube(B,s)}\tau_*T\wedge \pi^*(\omega^{j_0})\wedge \beta^{k-p-j_0}\\
 &=\nu_{j_0,q_0}(T,B,r,\tau)- \nu_{j_0,q_0}(T,B,\tau)\leq   \nu_{j_0,q_0}(T,B,r,\tau) .
 \end{split}
\end{equation}
 Thus, we deduce from  \eqref{e:P-Lc-bullet_finite(1)-psh}--\eqref{e:P-Lc-bullet_finite(3)-psh} that
\begin{equation}\label{e:P-Lc-bullet_finite(4)-psh}
\int_{\Tube(B,\bfr)}\tau_*T\wedge \pi^*(\omega^{j_0})\wedge \beta^{j'_0}\wedge\alpha^{q_0} \leq \nu_{j_0,q_0}(T,B,r,\tau)+  cr+cr^{2j'_0}\Nc(T,r).
\end{equation}
Arguing as in the proof of Theorem \ref{T:Lc-finite}, we  obtain
the  following  equality
\begin{equation*}
 \begin{split}
&\int_{\Tube(B,r)}\tau_*T\wedge \pi^*(\omega^{j_0})\wedge \beta^{j'_0}\wedge\alpha^{q_0}
   = I_{q_0,0,j_0,0}(T,r)\\&+\sum_{j'_1,j''_1,j_1} {j'_0\choose  j'_1}{q_0 \choose j_1}
  {q_0-j_1 \choose j''_1}(-c_1)^{j'_0-j'_1}(-1)^{q_0-j_1-j''_1}
 I_{j_1, j'_0-j'_1, q_0+j_0+j'_0-j_1-j'_1-j''_1,q_0-j_1-j''_1}(T,r).
\end{split}
\end{equation*}
    Observe  that  RHS can be rewritten as  the  sum
        $\Ic_1+\Ic_2+\Ic_3,$  where  $\Ic_j$ for $1\leq j\leq 3$    were defined in \eqref{e:P-Lc-bullet_finite(6)}.
Combining inequalities \eqref{e:P-Lc-bullet_finite(1)-bis-psh}, \eqref{e:P-Lc-bullet_finite(2)-psh}, \eqref{e:P-Lc-bullet_finite(3)-psh} and 
\eqref{e:P-Lc-bullet_finite(4)-psh} and increasing $c$ if necessary, we  deduce from the above  equality that
\begin{equation}\label{e:P-Lc-bullet_finite(6)-psh}
 \Ic_1+\Ic_2+\Ic_3\leq cr+cr^{2j'_0}\Nc(T,r)+\nu_{j_0,q_0}(T,B,r,\tau).
\end{equation}
Applying Lemma \ref{L:spec-wedge}  to  each  difference term  in $\Ic_2$ and $\Ic_3$ 
yields that
\begin{equation} \label{e:P-Lc-bullet_finite(7)-psh}
|I_\bfi(r)- I^\hash_\bfi(r)|^2 \leq c\big(\sum_{\bfi'} I^\hash_{\bfi'}(r)\big)\big ( \sum_{\bfi''} I^\hash_{\bfi''}(r) \big).  
\end{equation}
 Here, on the LHS  $\bfi=(i_1,i_2,i_3,i_4)$  is  either $(q_0,0,j_0,0)$ or $(j_1, j'_0-j'_1, q_0+j_0+j'_0-j_1-j'_1-j''_1,q_0-j_1-j''_1)$  with $j_1,j'_1,j''_1$ as above, 
and on the RHS:
\begin{itemize} \item[$\bullet$] the first sum  is taken over a finite number of multi-indices    $\bfi'=(i'_1,i'_2,i'_3,i'_4)$ as above  such that  $i'_1\leq  i_1$  and $i'_2\geq i_2;$ and either ($i'_3\leq i_3$) or ($i'_3>i_3$ and $i'_2\geq i_2+{1\over 2}$);   
\item  the second sum   is taken over  a finite number of multi-indices $\bfi''=(i''_1,i''_2,i''_3,i''_4)$ as above   such that   either  ($i''_1< i_1$)
or ($i''_1=i_1$ and $i''_2\geq {1\over 4}+i_2$) or ($i''_1=i_1$ and $i''_3<i_3$).
\end{itemize}
Using that $j'_0\geq 1$  and  arguing as  in the  proof of \eqref{e:P-Lc-bullet_finite(7)}--\eqref{e:P-Lc-bullet_finite(8)}, we  see that
 the first  sum on the RHS of \eqref{e:P-Lc-bullet_finite(7)-psh} is  bounded from above by a constant times  $\Kc^-_{j_0,q_0}(T,r)+ \Kc_{j_0,q_0}(T,r)+r^{1\over 2} \Kc^+_{j_0,q_0}(T,r),$
 whereas  the second sum  is bounded from above by a constant times $\Kc^-_{j_0,q_0}(T,r)+r^{1\over 2} \Kc_{j_0,q_0}(T,r)+ r^{1\over 2} \Kc^+_{j_0,q_0}(T,r).$
 Consequently,   we infer from \eqref{e:P-Lc-bullet_finite(6)-psh}, \eqref{e:P-Lc-bullet_finite(7)-psh} and \eqref{e:P-Lc-bullet_finite(6bis)} that there is a constant $c>0$ such that
\begin{equation}\label{e:P-Lc-bullet_finite(9)-psh}
 \begin{split}
&I^\hash_{q_0,0,j_0,0}(T,r)\leq  cr+ |\nu_{j_0,q_0}(T,r)|+cr^2\Nc(T,r)+ \\ 
&+c\sqrt{\Kc^-_{j_0,q_0}(T,r)+ \Kc_{j_0,q_0}(T,r)+r^{1\over 2} \Kc^+_{j_0,q_0}(T,r)} \sqrt{ \Kc^-_{j_0,q_0}(T,r)+r^{1\over 2} \Kc_{j_0,q_0}(T,r)+r^{1\over 2} \Kc^+_{j_0,q_0}(T,r)  }.
\end{split}
\end{equation}
Since 
$I^\hash_{q_0,0,j_0,0}(T,r) = \Kc_{j_0,q_0}(T,r)$ by Remark \ref{R:I_bfj}, it follows that
there is a constant $c>0$ such that
\begin{equation}\label{e:P-Lc-bullet_finite(9)-bis-psh}
 \begin{split}
I^\hash_{q_0,0,j_0,0}(T,r)&\leq  cr+ c|\nu_{j_0,q_0}(T,r)|+cr^2\Nc(T,r)+ \\ 
&+c\sqrt{\Kc^-_{j_0,q_0}(T,r)+ \Kc_{j_0,q_0}(T,r)+r^{1\over 2} \Kc^+_{j_0,q_0}(T,r)} \sqrt{ \Kc^-_{j_0,q_0}(T,r)+r^{1\over 2} \Kc^+_{j_0,q_0}(T,r)  }.
\end{split}
\end{equation}
When $j_0=m_0$ we deduce using  $\Kc^+_{m_0,q_0}(T,r)=0$   that 
\begin{eqnarray*}
 I^\hash_{q_0,0,m_0,0}(T,r)=\Kc_{m_0,q_0}(T,r)&\leq&  cr+ c|\nu_{m_0,q_0}(T,r)|+cr^2\Nc(T,r) \\
&+&c\sqrt{\Kc^\bullet_{q_0}(T,r)+ \Kc_{m_0,q_0}(T,r)} \sqrt{ \Kc^\bullet_{q_0}(T,r)  }.
\end{eqnarray*}
This  implies that
\begin{equation}\label{e:Kc_m_0,q_0}
 \Kc_{m_0,q_0}(T,r)\leq  c \Kc^\bullet_{q_0}(T,r)+  c  \big (r+ |\nu_{m_0,q_0}(T,r)|+r^2\Nc(T,r)\big) .
\end{equation}
 Hence, for $j'_0\geq 1$ we obtain    
 \begin{eqnarray*}\Kc^-_{j_0,q_0}(T,r)+ \Kc_{j_0,q_0}(T,r)+r^{1\over 2} \Kc^+_{j_0,q_0}(T,r)&\lesssim & \Kc^\bullet_{q_0}(T,r)+
 r^{1\over 2}\Kc_{m_0,q_0}(T,r)\\
 &\lesssim &\Kc^\bullet_{q_0}(T,r)+ r+ |\nu_{m_0,q_0}(T,r)|+r^2\Nc(T,r).
 \end{eqnarray*}
 Putting this  together with  the easily obtained inequality
 $$  \sqrt{\Kc^\bullet_{q_0}(T,r)}\sqrt{ r^{1\over 2} \Kc^+_{j_0,q_0}(T,r)}
 \leq  r^{1\over 4} \Kc^\bullet_{q_0}(T,r)+ r^{1\over 4}  \Kc^+_{j_0,q_0}(T,r) ,$$ 
 \eqref{e:P-Lc-bullet_finite(8)-psh} follows. 
This is the  desired estimate of  Step 2.

\noindent {\bf Step 3:} {\it  
End of the proof of
\eqref{e:first-ind}.}

Recall that $m_1:=k-p-q_0-1$ and that  by \eqref{e:Kc_m_0,q_0}  we have
 \begin{equation*}\Kc^+_{m_1,q_0}(T,r)=\Kc_{m_0,q_0}(T,r)
  \leq  c \Kc^\bullet_{q_0}(T,r)+  c  \big (r+ |\nu_{m_0,q_0}(T,r)|+r^2\Nc(T,r)\big) .
\end{equation*}
Using this and estimate \eqref{e:P-Lc-bullet_finite(8)-psh} and  arguing as in Step 3 of the proof of Theorem \ref{T:Lc-finite}, we  can prove for $1\leq j\leq m_1$ that
\begin{multline*}
I^\hash_{q_0,0,j,0}(T,r)
\leq  c  \Big (\Kc_{q_0-1}(T,r)+r +r^{1\over 2^{j+2}}\Nc(T,r)+\sum_{j=0}^{m_1} |\nu_{j,q_0}(T,r)|+ r\nu_{m_0,q_0}(T,r)+r^{1\over 2^{j+2}} \Kc^\bullet_{q_0}(T,r)+  \\
\big[(\Kc^\bullet_{q_0}(T,r))^{1\over 2^{j+1}}+  (\Kc^-_{q_0-1}(T,r))^{1\over 2^{j+1}}+   \big( \sum_{j=0}^{m_1} |\nu_{j,q_0}(T,r)|^{1\over 2^{j+1}}   \big)\big]^{2^{j+1}}
-\Kc^\bullet_{q_0}(T,r)
\Big ).
\end{multline*}
Using  \eqref{e:P-Lc-bullet_finite(0)-psh-bullet}
 and \eqref{e:P-Lc-bullet_finite(0)-psh}, we see that
\begin{equation*}
\Kc^\bullet_{q_0}(T,r) =\Kc_{q_0-1}(T,r)+\sum_{j=0}^{m_1} I^\hash_{q_0,0,j,0}(T,r).
\end{equation*}
This, combined with the previous  estimate, implies by increasing $c_{10}$ that 
\begin{multline*}
\Kc^\bullet_{q_0}(T,r) 
\leq  c_{10}  \Big(\Kc_{q_0-1}(T,r)+r+ r^{1\over 2^{m_1+2}}\Nc(T,r)+\sum_{j=0}^{m_1} |\nu_{j,q_0}(T,r)|+ r\nu_{m_0,q_0}(T,r) +r^{1\over 2^{m_1+2}} \Kc^\bullet_{q_0}(T,r)\\
+  \big[(\Kc^\bullet_{q_0}(T,r))^{1\over 2^{m_1+1}}+  (\Kc_{q_0-1}(T,r))^{1\over 2^{m_1+1}}+\big( \sum_{j=0}^{m_1} |\nu_{j,q_0}(T,r)|^{1\over 2^{m_1+1}}\big)\big]^{2^{m_1+1}}
-\Kc^\bullet_{q_0}(T,r)\Big)
 .
\end{multline*}
Using the last   estimate and using $m_1$ instead of $m_0,$ we argue as in  the end of Step 3 of  the proof of Theorem
\ref{T:Lc-finite}.   Hence, the first  inequality of \eqref{e:first-ind} follows.

Combining  \eqref{e:P-Lc-bullet_finite(1)-bis-psh} and Lemmas  \ref{L:Lc-vs-analog-bis-new}, \ref{L:Lc-vs-analog-bis'-new} and  \ref{L:Lc-vs-analog-bisbis-new}, 
\begin{equation*}
\Lc_{j_0,q_0}(T,r)\leq c \Nc(T,r)+c\sum_{j'=1}^{\upm-j_0}  \Lc_{j_0+j',q_0-j'}(T,r).
\end{equation*}
This,  coupled   with the  inequality $\Nc^\bullet_q(T,r)\geq \Lc_{q-1}(T,r)$, implies    the   second  inequality of \eqref{e:first-ind}.
The conclusion of Step 3 is thereby completed.

{\it 
Now it remains to treat the case  where $j=m_0:=k-p-q,$ that is,  there is a constant $c_{10}$ independent of $T$
such that 
\begin{equation}\label{e:second-ind}
\Kc_{q}(T,r)
\leq c_{10}\Nc_q(T,r)\qquad\text{and}\qquad \Lc_{j,q}(T,r)\leq c_{10}\Nc_q(T,r),
\end{equation}
for every $0\leq q\leq \min(k-l,k-p).$ Here,
\begin{eqnarray*}
\Nc_0(T,r)&:=&r+\Nc(T,r)+\sum_{j=0}^{k-p} |\nu_{j,0}(T,r)|,\\
\Nc_q(T,r)&:=&\Nc_{q-1}(T,r)+\sum_{j=0}^{m_0} |\nu_{j,q}(T,r)|\quad\text{for}\quad q\geq 1.
\end{eqnarray*}
The proof of \eqref{e:second-ind}  will be completed in Steps 4 and 5 below.
By Steps  2 and 3,  inequality  \eqref{e:second-ind} is reduced to proving that
\begin{equation}\label{e:second-ind-bis}
\Kc_{m_0,q}(T,r)
\leq c_{10}\Nc_{q}(T,r)\qquad\text{and}\qquad \Lc_{m_0,q}(T,r)\leq c_{10}\Nc_{q}(T,r).
\end{equation}
}

\noindent {\bf Step 4: }{\it  Inequality \eqref{e:second-ind-bis}
 holds 
for every $0\leq q<  k-p-\upm.$ 
}
    
 Since $k-p-q>\upm$ and $\pi^*\theta_\ell\cdot (\tau_\ell)_* (T|_{\bfU_\ell})\wedge \pi^*\omega^\upm$ is  full in 
 bidegree $(dw,d\bar w),$ it follows that  $\pi^*\theta_\ell\cdot (\tau_\ell)_* (T|_{\bfU_\ell})\wedge\pi^*\omega^{k-p-q}=0,$ and hence $ \Kc_{m_0,q}(T,r)=0$ and  $ \Lc_{m_0,q}(T,r)=0.$ 
 So  \eqref{e:second-ind-bis} is trivially  fulfilled  in this  case.
 Step 4  follows.

\noindent {\bf Step 5: }{\it  Inequality \eqref{e:second-ind}
 holds 
for every $k-p-\upm\leq q\leq  k-p-\lowm.$ 
}

We make the second induction on $q.$ Suppose  inequality \eqref{e:second-ind}
 holds 
for every $q$ with $0\leq q<q_0,$  where $q_0$ is a  given  integer with $k-p-\upm\leq q_0\leq  k-p-\lowm.$
We need to show that  \eqref{e:second-ind} also
 holds for $q_0.$  Set $j_0:=m_0=k-p-q_0.$

By   Lemma \ref{L:Lc-vs-analog},
we have that 
\begin{equation*}
\int_{0}^{ r} \big({1\over t^{2q_0}} -{1\over r^{2q_0}}\big)2tdt  \int_{\Tube(B,t)}   \tau_*(\ddc T)\wedge \pi^*(\omega^{j_0})\wedge (\beta+c_1t^2\pi^*\omega)^{q_0-1}\\
 \geq -cr\Mc^\tot(  \ddc T,r ).
\end{equation*}
Using the identity
 $$  \pi^*(\omega^{j_0})\wedge (\beta+c_1t^2\pi^*\omega)^{q_0-1}=\sum_{j=0}^{q_0-1}  {q_0-1\choose j}  c_1^jt^{2j} \pi^*(\omega^{j_0+j})\wedge \beta^{q_0-1-j},$$
it follows  that
$I_1+I_2\geq -cr\Mc^\tot(  \ddc T,r),$ where 
\begin{eqnarray*}
 I_1&:=&
 \int_0^r\big({1\over t^{2q_0}} -{1\over r^{2q_0}}\big)2tdt  \int_{\Tube(B,t)}   \tau_*(\ddc T)\wedge \pi^*(\omega^{j_0})\wedge \beta^{q_0-1},\\
 I_2&:=&\sum_{j=1}^{q_0-1}  {q_0-1\choose j}c_1^j\int_{0}^{ r} \big({1\over t^{2q_0}} -{1\over r^{2q_0}}\big)2t^{2j+1}dt  \int_{\Tube(B,t)}   \tau_*(\ddc T)\wedge \pi^*(\omega^{j_0+j})\wedge \beta^{q_0-1-j}.
 \end{eqnarray*}
  Write  each double integral of $I_2$  as  follows:
  \begin{multline*}
   \int_{0}^{ r} {2dt\over t^{2(q_0-j)-1}}  \int_{\Tube(B,t)}   \tau_*(\ddc T)\wedge \pi^*(\omega^{j_0+j})\wedge \beta^{q_0-1-j}\\
 - \int_{0}^{ r} {2t^{2j+1}dt\over r^{2q_0}}  \int_{\Tube(B,t)}   \tau_*(\ddc T)\wedge \pi^*(\omega^{j_0+j})\wedge \beta^{q_0-1-j}
 =
 \end{multline*}
Combining  Lemma \ref{L:Lc-vs-analog-bis'} and Lemma  \ref{L:Lc-vs-analog-bisbis}, the absolute value of the first integral  is  bounded by  a constant times $\Lc_{j_0+j,q_0-j}(T,r)+\Nc_{q_0-j}(T,r).$
  Moreover,  the second integral is equal to
  \begin{equation*}
   {1\over r^{2q_0}}\int_0^r t^{2q_0-1}\nu(\ddc T,B,t,\tau) dt,
  \end{equation*}
which is  bounded in absolute value  by a constant times $\Mc^\tot(\ddc T,r).$
  
 Since  $q:=q_0-j<q_0$ for $j\geq 1,$ we can apply  the induction hypothesis of Step 5 in order  to  conclude that
 $I_2\leq c\Nc_{q_0-1}(T,r).$  
 Hence,  we  can find  a constant $c>0$ independent of $T$ and $0<r\leq \bfr$ such that
\begin{equation}\label{e:P-Lc-bullet_finite(1)-bisbis-psh} \int_0^r\big({1\over t^{2q_0}} -{1\over r^{2q_0}}\big)2tdt  \int_{\Tube(B,t)}   \tau_*(\ddc T)\wedge \pi^*(\omega^{j_0})\wedge \beta^{q_0-1}
 \geq -c\Mc^\tot(  \ddc T,r) - c\Nc_{q_0-1}(T,r).
 \end{equation} 

 Now  we repeat the argument which has been used from \eqref{e:P-Lc-bullet_finite(1)-psh} to the end of Step  2.
 In the present context  $j'_0=0.$ Note that $\Kc^{-}_{j_0,q_0}(T,r)= \Kc^\bullet_{q_0}(T,r).$
 We use   inequality \eqref{e:P-Lc-bullet_finite(1)-bisbis-psh} instead of 
  \eqref{e:P-Lc-bullet_finite(1)-bis-psh}.
As $\Kc^{+}_{j_0,q_0}(T,r)=0$ and $\Kc^-_{j_0,q_0}(T,r)=\Kc^\bullet_{q_0}(T,r)$ we deduce  from \eqref{e:P-Lc-bullet_finite(9)-psh}
    that there is  a constant $c>0$ such that for   every $0<r\leq \bfr,$
\begin{equation*}
\begin{split}
I^\hash_{q_0,0,j_0,0}(T,r)&\leq   cr+ |\nu_{j_0,q_0}(T,r)|+c\Kc_{q_0-1}(T,r)+c\Nc(T,r)+ \\ 
&+c\sqrt{\Kc^\bullet_{q_0}(T,r)+ \Kc_{j_0,q_0}(T,r)} \sqrt{ \Kc^\bullet_{q_0}(T,r) + (\Kc_{j_0,q_0}(T,r))^{1\over 2} }
\end{split}
\end{equation*}
So  we  infer that
\begin{equation*}
I^\hash_{q_0,0,j_0,0}(T,r)=\Kc_{j_0,q_0}(T,r)\leq   cr+ |\nu_{j_0,q_0}(T,r)| +c\Nc(T,r)+c \Kc^\bullet_{q_0}(T,r)
.
\end{equation*}
Using  \eqref{e:P-Lc-bullet_finite(0)-psh-bullet}
 and \eqref{e:P-Lc-bullet_finite(0)-psh}, we see that
\begin{equation*}
\Kc_{q_0}(T,r) =\Kc^\bullet_{q_0}(T,r)+I^\hash_{q_0,0,j_0,0}(T,r)=\Kc^\bullet_{q_0}(T,r)+ \Kc_{j_0,q_0}(T,r).
\end{equation*}
This, combined with the previous  estimate, 
implies  that
$$
\Kc_{q_0}(T,r)\leq   cr+ |\nu_{j_0,q_0}(T,r)|+(c+1)\Kc^\bullet_{q_0}(T,r)+c\Nc(T,r). 
$$
This, coupled  with  the first inequality in  \eqref{e:first-ind}, 
 gives  the first  inequality of \eqref{e:second-ind}.

We  turn  to the proof of  the second  inequality of \eqref{e:second-ind}.  Using \eqref{e:P-Lc-bullet_finite(1)-psh} and \eqref{e:P-Lc-bullet_finite(2)-psh} and  \eqref{e:P-Lc-bullet_finite(3)-psh} for $j_0=m_0=k-p-q_0,$  we get that 
\begin{eqnarray*}
\int_0^r \big({1\over t^{2q_0}} -{1\over r^{2q_0}}\big)2tdt   \int_{\Tube(B,t)}   \ddc (\tau_* T)\wedge \pi^*(\omega^{j_0})\wedge \beta^{q_0-1} &\leq & c  \Kc_{q_0}(T,r) + \nu_{m_0,q_0}(T,B,r,\tau)+cr\\
&\leq &   c_{10}\Nc_{q_0}(T,r),
 \end{eqnarray*}
 where the  last  estimate  holds by the first  inequality of \eqref{e:second-ind}
 and  $c_{10}$ is  a constant large enough independent of $T$ and $r.$
Using  this and applying Lemma \ref{L:Lc-vs-analog-bis}, 
we get that
\begin{eqnarray*}
\int_0^r \big({1\over t^{2q_0}} -{1\over r^{2q_0}}\big)2tdt   \int_{\Tube(B,t)}  ( \ddc  T)^\hash \wedge \pi^*(\omega^{j_0})\wedge \beta^{q_0-1}\leq &   c \Nc(T,r).
 \end{eqnarray*}
Using this and  applying  Lemma \ref{L:Lc-vs-analog-bisbis}  and applying  the second  inequality of \eqref{e:second-ind}
for $q<q_0$ (the inductive hypothesis), 
we get  the second  inequality \eqref{e:second-ind}.
This proves Step 5, and  the proof of the theorem is thereby completed.

  \endproof

\begin{corollary}\label{C:ddc-positive-finite}
There is a constant $c_{11}>0$  such that for every  positive plurisubharmonic   current $T$ in the class $\widetilde\SH^{3,3}_p(\bfU,\bfW),$    and every $q,j$ with $0\leq q\leq  \min(k-l,k-p-1)$ and  $0\leq j\leq k-p-q-1,$ we have  
$$
 \int_0^\bfr{dt\over t^{2q-1}}  \big( \int_{\Tube(B,t)}(\ddc T)^\hash\wedge
 \pi^*\omega^j\wedge  \hat\beta^{k-p-j-1}\big)
<c_{11}.
$$ 
\end{corollary}

\proof Since $\hat\beta\leq c(\beta+c_1t^2\pi^*\omega)$ on $\Tube(B,t)$ for a constant $c$ independent of $t,$
the  desired estimate follows  immediately from    the  inequality $\Lc_{j,q}(T,\bfr)<c_{10}$  obtained in Theorem \ref{T:Lc-finite-psh}.

\endproof
\begin{theorem}\label{T:vanishing-Lelong}
 For every positive plurisubharmonic   current $T$ such that $T=T^+-T^-$ on an open neighborhood of $\overline B$ in $X$  with $T^\pm$  in the class $\SH^{3,3}_p(B),$    and every $\lowm^+\leq j\leq \upm^+,$
 we have  $\nu_j(\ddc T,B,\tau)=0.$
\end{theorem}
\proof
Suppose  that there is  an index $j$ with $\lowm^+\leq j\leq \upm^+$ such that $ \nu_j(\ddc T,B,t,\tau) \not=0.$ We may assume  without loss of generality that  $T$ is in the class $\widetilde\SH^{3,3}_p(\bfU,\bfW),$
 By 
 Lemmas  \ref{L:Lc-vs-analog-bis'} and  \ref{L:Lc-vs-analog-bisbis},   there is a  constant  $c>0$ such that  for every $0<r\leq\bfr,$
  \begin{equation*}
   \big|  \int_0^{ r} {2dt\over t^{2(k-p-j)-1}} \big(  \int_{\Tube(B,t)}  \tau_*(\ddc T)\wedge \pi^*(\omega^{j})\wedge \beta^{k-p-j-1}\big) \big|\leq c \Lc_{k-p-j}(T,r)+c\Nc(T,r).
  \end{equation*}
  By   Theorem \ref{T:Lc-finite-psh}, the absolute value of the  expression on the RHS is   bounded from above by a constant $c'$ independent of $T$ and $r.$ In particular, 
  the absolute value of the expression on the LHS is $\leq c'.$
  We  rewrite this  inequality as:
 \begin{equation*}
   \big|  \int_0^{r} {\nu_j(\ddc T,B,t,\tau)  dt\over t}  \big|\leq c'.
  \end{equation*}
   Since $\lim_{t\to 0}\nu_j(\ddc T,B,t,\tau)=\nu_j(\ddc T, B,\tau)\not=0,$ it follows that   for $r>0$  small enough,
   $\big|  \int_0^{ r}  { dt\over t}  \big|<\infty.$ This is  a contradiction.
\endproof

 \begin{proposition}\label{P:Lc-finite-psh-bis}
  For  $0<r_1<r_2\leq \bfr,$ there is  a constant $c_{11}>0$ such that for every $q\leq  \min(k-p,k-l)$ and 
    every positive plurisubharmonic   current $T$ in the class $\widetilde\SH^{3,3}_p(\bfU,\bfW),$   we have the following estimate:
 $$|\kappa_{k-p-q}(T,{r_1\over \lambda},{r_2\over \lambda},\tau)|<c_{11}\sum\limits_{0\leq q'\leq q,\ 0\leq j'\leq \min(\upm,k-p-q')} \Kc_{j',q'}(T,{r_1\over \lambda},{r_2\over \lambda} )\qquad\text{for}\quad  \lambda>1.$$
 \end{proposition}
\proof  
 It follows along the same lines as those of  the proof of Proposition  \ref{P:Lc-finite-bis}.
 
 \endproof

\subsection{End of the proof for positive  plurisubharmonic  currents}

This  subsection  is  devoted to the proof  of  Theorem \ref{T:Lelong-psh}.

\proof[Proof of assertion (1) of Theorem \ref{T:Lelong-psh}]
Fix $r_1,\ r_2\in  (0,\bfr]$ with $r_1<r_2.$ 
 Fix  $j$ with $0\leq j\leq \upm$ and  let $\lambda >1.$  
Applying Theorem  \ref{T:Lelong-Jensen}  to  $(A_\lambda)_*\tau_*T\wedge \pi^*(\omega^m),$  we get that
\begin{multline*}
  \nu_j(T,B,{r_2\over \lambda},\tau)- \nu_j(T,B,{r_1\over \lambda},\tau)-  \kappa_j(  T,B,{r_1\over \lambda}, {r_2\over \lambda},\tau   )\\
 =  \int_{r_1\over\lambda }^{r_2\over \lambda} \big( {1\over t^{2(k-p-j)}}-{\lambda^{2(k-p-j)} \over r_2^{2(k-p-j)}}  \big)2tdt\int_{\Tube(B,t)}  \ddc [\tau_*T\wedge \pi^*(\omega^j)\wedge \beta^{k-p-j-1}] \\
 +  \big( {\lambda^{2(k-p-j)}\over r_1^{2(k-p-j)}}-{\lambda^{2(k-p-j)}\over r_2^{2(k-p-j)}}  \big) \int_{0}^{r_1\over\lambda }2tdt\int_{z\in \Tube(B,t)} \ddc [\tau_*T\wedge \pi^*(\omega^j)\wedge \beta^{k-p-j-1}].
\end{multline*}
By Corollary \ref{C:ddc-estimate-m-1-lambda} with $m=0,$ the  two terms on the RHS are of  modulus  smaller  than  a constant times  
  $ \lambda^{-1}\Mc^\tot( T,{r_2\over \lambda} )+\Mc^\tot( \ddc T,{r_2\over \lambda} ).$
 For $\lambda\geq 1$ set  
 $$
\epsilon_\lambda:= |\kappa_j(  T,B,{r_1\over \lambda}, {r_2\over \lambda},\tau   )| + \lambda^{-1}\Mc^\tot(T, {r_2\over \lambda})+\Mc^\tot( \ddc T,{r_2\over \lambda} ).
$$
  Consequently,  we infer that
 \begin{equation*}
\big| \nu_j(T,B,{r_2\over \lambda},\tau)- \nu_j(T,B,{r_1\over \lambda},\tau) \big|\leq \epsilon_\lambda.
\end{equation*}
 We need to show  that $\sum_{n=0}^\infty\epsilon_{2^n\lambda}<\infty $ for $\lambda\geq 1.$ 
Applying   Proposition  \ref{P:Lc-finite-psh-bis}  yields that
 $$|\kappa_{j}(T,{r_1\over \lambda},{r_2\over \lambda},\tau)|<{c_{10}\over \lambda}+c_{10}\sum_{q'\leq q} \Kc_{j',q}(T,{r_1\over \lambda},{r_2\over \lambda})\quad\text{for}\quad 0<s<r<\bfr.$$
 Since Theorem  \ref{T:Lc-finite-psh}  says that
 $\Kc_{j',q'}(T,\bfr)<c_9$ for $0\leq q'\leq k-l$ and $0\leq j'\leq k-p-q',$  we infer that
$$
 \sum_{n=0}^\infty|\kappa_{j}(T,{r_1\over 2^n},{r_2\over 2^n},\tau)|\leq  c \Kc_{k-p-j}(T,\bfr)<\infty.
$$
Moreover,   there is a constant $c>0$ independent  of  $T$   such that  obtain   that
$$
 \sum_{n=0}^\infty\Mc^\tot( \ddc T,{r_2\over 2^n} )\leq c\Lc_{j,k-p-j}(T,\bfr)
$$
 By  Theorem  \ref{T:Lc-finite-psh}, the RHS is  finite.
 Next,   
by Proposition \ref{P:nu_tot-sum-psh},  there is  a constant $c>0$ such that 
\begin{equation*}
  \sum_{n=0}^\infty {1\over 2^n}\Mc^\tot(T,{r_2\over 2^n})\leq cr_2^{-1} \Nc^\tot(T,r_2) +c<\infty.
 \end{equation*} 
Combining   together the last inequalities,   we have shown  that $\sum_{n=0}^\infty\epsilon_{2^n\lambda}<\infty. $ 
Applying  Lemma \ref{L:elementary} (1) yields that
 $\lim_{r\to 0+} \nu_j(T,B,r,\tau) \in \R ,$
  and assertion (1) follows.
\endproof

\proof[Proof of assertion (2) of Theorem \ref{T:Lelong-psh}]
We need to show that
\begin{equation*}
 \lim\limits_{r\to 0}   \sup_{ s\in(0,r)}|\kappa_j(T,B,s,r)|=0.
\end{equation*}
The above limit does not exceed  
\begin{equation*}
 \lim\limits_{r\to 0}\sup_{ s\in(0,r)}|\kappa_j(T,B,{s\over 2},r)|+  \lim\limits_{s\to 0}|\kappa_j(T,B,{s\over 2},s)|.
\end{equation*}
Therefore, we  are reduced  to proving  that
\begin{equation}\label{e:reduced-kappa-bullet}
 \lim\limits_{r\to 0}  \sup_{ s\in (0,{r\over 2}]}|\kappa_j(T,B,s,r)|=0.
\end{equation}
Using  $ 0<s\leq {r\over 2},$ we  argue as in the proof of 
 Proposition \ref{P:Lc-finite-psh-bis}. Consequently, we  get  the following estimate:
 $$|\kappa_{j}(T,s,r,\tau)|<c_{11}\sum\limits_{0\leq q'\leq k-p-j,\ 0\leq j'\leq \min(\upm,k-p-q')} \Kc_{j',q'}(T,{s\over 2},2r )\qquad\text{for}\quad  
 0<r<{\bfr\over 2}.$$
 On the  other hand,  since Theorem  \ref{T:Lc-finite-psh}  says that
 $\Kc_{j',q'}(T,\bfr)<c_9$ for $0\leq q'\leq k-l$ and $0\leq j'\leq k-p-q',$  we infer that
$$  \lim\limits_{r\to 0}\sup_{s\in (0,{r\over 2}]}\Kc_{j',q'}(T,{s\over 2},2r )=0.$$
This,  combined  with the above upper-bound for $|\kappa_{j}(T,s,r,\tau)|,$ gives  the  desired  estimate \eqref{e:reduced-kappa-bullet}.
\endproof

\proof[Proof of assertion (3) of Theorem \ref{T:Lelong-psh}]

\noindent {\bf Proof of the interpretation of assertion (3)  in the sense  of 
Definition \ref{D:Lelong-log-numbers(2)}. }

Fix  an index $j$ with $\lowm\leq j\leq\upm.$ Fix $0<r\leq\bfr$ and let  $0<\epsilon<r.$  Theorem \ref{T:Lelong-Jensen-eps} applied  to  $\tau_*T\wedge \pi^*(\omega^j)$  gives
\begin{equation}\label{e:T:Lelong-psh-(2)}\begin{split}
 &\qquad {1\over  (r^2+\epsilon^2)^{k-p-j}} \int_{\Tube(B,r)} \tau_*T\wedge \pi^*(\omega^j)\wedge \beta^{k-p-j} 
 = \Vc_\epsilon(\tau_*T\wedge \pi^*(\omega^j),r)\\&+   \int_{\Tube(B,r)} \tau_*T\wedge \pi^*(\omega^j)\wedge \alpha_\epsilon^{k-p-j}\\
 &+  \int_{0}^{r} \big( {1\over (t^2+\epsilon^2)^{k-p-j}}-{1\over (r^2+\epsilon^2)^{k-p-j}}  \big)2tdt\int_{\Tube(B,t)} \ddc   [ \tau_*T\wedge \pi^*(\omega^j) ]\wedge \beta^{k-p-j-1}. 
 \end{split}
\end{equation}  
Next, we let $\epsilon$ tend to $0.$  Then the LHS of \eqref{e:T:Lelong-psh-(2)}  tends to $\nu_j(T,B,r,\tau).$  On the other hand,  we deduce  from Theorem \ref{T:vertical-boundary-eps} that
$\Vc_\epsilon( \tau_*T\wedge \pi^*(\omega^j),r)=O(r).$ Moreover, the third term on the RHS of \label{e:T:Lelong-psh-(1)}  is  rewritten as
\begin{equation}\label{e:T:Lelong-psh-(3)}
 \int_{0}^{r} \big( {1\over (t^2+\epsilon^2)^{k-p-j}}-{1\over (r^2+\epsilon^2)^{k-p-j}}  \big)2t(f^+(t)-f^-(t))dt,
\end{equation}
where
\begin{equation*}
 f^\pm(t):=\int_{\Tube(B,t)}  \tau_*(\ddc T^\pm)\wedge (\pi^*\omega^j)\wedge \beta^{(k-p-j)-1}.
\end{equation*}
Combining a variant of  Lemma  \ref{L:Lc-vs-analog-bis'-new} and
 Lemma \ref{L:Lc-vs-analog-bisbis-new}, 
  there is a constant  $c>0$ such that  
  \begin{equation*}
       \int_{0}^{r} \big( {1\over t^{2(k-p-j)}}-{1\over r^{2(k-p-j)}}  \big)2t|f^\pm(t)|dt 
   \leq 
   c\sum_{j'=0}^{\upm-j}  \Lc_{j+j',q-j'}(T,r)
   + 
   c\Nc(T,r).
  \end{equation*}
  By  Theorem \ref{T:Lc-finite-psh}, the  RHS is  bounded by a constant $c_{10}.$
    So there is a constant $c>0$ independent of $T$ and $0<r\leq \bfr$  such that
\begin{equation}\label{e:T:Lelong-psh-(5)}
 \int_{0}^{r} \big( {1\over t^{2(k-p-j)}}-{1\over r^{2(k-p-j)}}  \big)2t|f^\pm(t)|dt\leq c.
\end{equation}
Observe  that  for  $t\in [0,r],$  we have  as $\epsilon\searrow 0,$  
\begin{multline*}
0\leq  {1\over (t^2+\epsilon^2)^{k-p-j}}-{1\over (r^2+\epsilon^2)^{k-p-j}}\approx {(r^2-t^2)\over  (t^2+\epsilon^2)^{k-p-j}(r^2+\epsilon^2) } \nearrow {(r^2-t^2)\over  t^{2(k-p-j)}r^2 } \\
\approx  {1\over t^{2(k-p-j)}}-{1\over r^{2(k-p-j)}}.
\end{multline*}
An application of  Lebesgue  Dominated Convergence  yields that the expression in \eqref{e:T:Lelong-psh-(3)} converges,  as $\epsilon\searrow 0,$  to  
\begin{equation}\label{e:T:Lelong-psh-(4)}
 \int_{0}^{r} \big( {1\over t^{2(k-p-j)}}-{1\over r^{2(k-p-j)}}  \big)2t(f^+(t)-f^-(t))dt.
\end{equation}
On the  other hand, 
\begin{equation}\label{e:T:Lelong-psh-(6)}
 \int_{0}^{r} {1\over r^{2(k-p-j)}} 2t|f^\pm(t)|dt= \int_{0}^{r} {1\over r^{2(k-p-j)}} 2t^{2(k-p-j)-1}|\nu_j(\ddc T^\pm,B,t,\tau)|dt\to 0,
\end{equation}
because  $\nu_j(\ddc T^\pm,B,\tau)=0$ by assertion (4).
This, combined with \eqref{e:T:Lelong-psh-(5)}, implies that  
 $$
 \lim_{r\to 0} \int_{0}^{r} \big( {1\over t^{2(k-p-j)}} \big)2t|f^\pm(t)|dt=0.
 $$
 This, coupled with  \eqref{e:T:Lelong-psh-(6)}, gives 
\begin{equation}\label{e:T:Lelong-psh-(7)}
 \lim_{r\to 0} \int_{0}^{r} \big( {1\over t^{2(k-p-j)}}-{1\over r^{2(k-p-j)}}  \big)2t|f^\pm(t)|dt =0
\end{equation}
Consequently,  by assertions (2) and (3),  the integral in  \eqref{e:T:Lelong-psh-(4)} is  bounded and it converges to $0$  as $r\to 0+.$ Putting this,   \eqref{e:T:Lelong-psh-(2)} and  \eqref{e:T:Lelong-psh-(3)} together,  
we obtain the  desired  interpretation  according to  Definition \ref{D:Lelong-log-numbers(2)}.

\noindent {\bf Proof of the interpretation of assertion (3)  in the sense  of 
Definition \ref{D:Lelong-log-numbers(1)}. }

Since $j>l-p$   it follows  from \eqref{e:m}  that $k-p-j<k-l.$
Therefore, we are in the position to apply  Theorem  \ref{T:Lelong-Jensen-smooth}  to the case where $q:=k-p-j<k-l.$ 
Hence, we get that 
\begin{equation*}
\nu_j( T^\pm_n,B,r,\tau)=\kappa_j(T^\pm_n,B,r,\tau)+\Vc(\tau_*T^\pm_n\wedge \pi^*(\omega^j),r)
+\int_{0}^{r} \big( {1\over t^{2(k-p-j)}}-{1\over r^{2(k-p-j)}}  \big)2tf^\pm_n(t)dt.
\end{equation*}
Thus,  we  obtain 
\begin{eqnarray*}
\kappa_j(T,B,r,\tau)&:=& \lim\limits_{n\to\infty} \kappa_j(T^+_n-T^-_n,B,r,\tau)=\lim\limits_{n\to\infty}\kappa_j(T^+_n,B,r,\tau)-\lim\limits_{n\to\infty}\kappa_j(T^-_n,B,r,\tau)\\
&=&\lim\limits_{n\to\infty}\nu_j(T^+_n,B,r,\tau)-\lim\limits_{n\to\infty}\nu_j(T^-_n,B,r,\tau)\\
&-&\lim\limits_{n\to\infty}\int_{0}^{r} \big( {1\over t^{2(k-p-j)}}-{1\over r^{2(k-p-j)}}  \big)2t(f^+_n(t)-f^-_n(t))dt\\
&=&\nu_j(T,B,r,\tau) -\lim\limits_{n\to\infty}\int_{0}^{r} \big( {1\over t^{2(k-p-j)}}-{1\over r^{2(k-p-j)}}  \big)2t(f^+_n(t)-f^-_n(t))dt,
\end{eqnarray*}
where 
$$
f^\pm_n(t):=\int_{\Tube(B,t)}  \tau_*(\ddc T^\pm_n)\wedge (\pi^*\omega^j)\wedge \beta^{(k-p-j)-1}.
$$
So the   interpretation  according to  Definition \ref{D:Lelong-log-numbers(1)} will hold if one can show that
\begin{equation}\label{e:T:top-Lelong-psh-(8)-bis}
 \lim\limits_{n\to\infty}\int_{0}^{r} \big( {1\over t^{2(k-p-j)}}-{1\over r^{2(k-p-j)}}  \big)2t(f^+_n(t)-f^-_n(t))dt  \to 0\quad\text{as}\quad r\to 0.
\end{equation}
  It is  not  difficult  to see that   $f^\pm_n(t)\to f^\pm(t)$ as $n\to\infty$ for all  $t\in(0,\bfr)$ except  for a countable set of values. Moreover, since $f^\pm_n(t)= t^{2(k-p-j-1)} \nu(\ddc,B,t,\tau),$
  we infer from
  Proposition \ref{P:Lc-finite-closed} applied  to $\ddc T\in \widetilde\CL^{1,1}_p(\bfU,\bfW)$ that for every $\bfr'\in (0,\bfr),$ there is a constant $c=c(\bfr')>0$ such  
that  $0\leq |f^\pm_n(t)|\leq ct^{2(k-p-j-1)}$ for all $n\geq 1$ and  $t\in (0,\bfr').$
Consequently, 
\begin{equation*}
 \lim\limits_{n\to\infty}\int_{0}^{r} \big( {1\over t^{2(k-p-j)}}-{1\over r^{2(k-p-j)}}  \big)2t(f^+_n(t)-f^-_n(t))dt   =\int_{0}^{r} \big( {1\over t^{2(k-p-j)}}-{1\over r^{2(k-p-j)}}  \big)2t(f^+(t)-f^-(t))dt .
\end{equation*}
So the  desired estimate  \eqref{e:T:top-Lelong-psh-(8)-bis} follows  immediately from inequality \eqref{e:T:Lelong-psh-(7)}.
\endproof

\proof[Proof of assertions (4) and (5) of Theorem \ref{T:Lelong-psh}]
Using assertion (1) we can  show that
all  the results   established in  
Subsection 
\ref{SS:Other-charac-Lelong} still hold when $T$ is a current in $ \widetilde\SH^{3,3}_p(\bfU,\bfW).$
In particular,  arguing as in the proof of assertions (5) and (6) of Theorem  \ref{T:Lelong-closed-Kaehler},
we  obtain   assertions (5) and (4) of  Theorem \ref{T:Lelong-psh}.
\endproof

\proof[Proof of assertions (6) of Theorem \ref{T:Lelong-psh}]
Notice  that   the assumption  $T^\pm\in \SH^{3,3}_p(\bfU,\bfW)$ is only  necessary  
to infer that $\ddc T\in \CL^{1,1}_p(\bfU,\bfW).$
When  $T=T^+-T^-$ for some   positive  pluriharmonic $(p,p)$-currents $T^\pm\in \PH^{2,2}_p(\bfU,\bfW),$
we have $\ddc T^\pm=0,$ and  in particular $\ddc T^\pm\in \CL^{1,1}_p(\bfU,\bfW).$ Hence,  all the  above assertions  still hold.  
\endproof

\section{
Non-K\"ahler metrics}
\label{S:non-Kaehler-metrics}

In this  section we  study 
positive closed  currents along a  submanifold   endowed with a
non-K\"ahler metric.

\subsection{Preliminary  estimates }


The following   result  is  the  analogue of   Proposition \ref{P:Stokes}
in this  section. 

\begin{proposition}\label{P:Stokes-non-Kaehler}
 Fix  $\ell$  with $1\leq \ell\leq \ell_0$ and $r\in(0,\bfr].$    Set $\tilde\tau:=\tilde\tau_\ell$  and $\H:=\Tube (\widetilde V_\ell,r)\subset \E.$
Then, for every 
 every   current $S$ of bidimension $(q-1,q-1)$  defined on  $\U_\ell$ and every  smooth  form $\Phi$ of bidegree $(q,q)$  defined on $\tilde\tau(\H)$
 with $\pi(\supp(\Phi))\Subset \widetilde V_\ell,$  we have
\begin{equation*}
 \langle  \dbar (\tilde\tau_* S) -   \tilde\tau_* (\dbar S),\Phi\rangle_{\tilde \tau (\H)}
   = \big \langle  \tilde\tau^*[ (\tilde\tau_* S)^\sharp]  ,\tilde\tau^*\Phi \big \rangle_{\partial \H}
 - \big \langle  \tilde\tau^*[ (\tilde\tau_* S)^\sharp]  ,\tilde\tau^*(d\Phi) \big \rangle_{ \H}
 -\langle \dbar S ,\tilde\tau^*\Phi\rangle_{ \H}.
\end{equation*}
\end{proposition}
\proof On the one hand, we have  
$$
 \langle   \tilde\tau_* (\dbar S),\Phi\rangle_{\tilde \tau (\H)}=\langle \dbar S ,\tilde\tau^*\Phi\rangle_{ \H}.
$$
On the  other hand, by a bidegree  consideration we write
\begin{equation*}
 \langle  \dbar (\tilde\tau_* S) ,\Phi\rangle_{\tilde \tau (\H)}=\langle  \dbar (\tilde\tau_* S)^\sharp ,\Phi\rangle_{\tilde \tau (\H)}=\langle  d (\tilde\tau_* S)^\sharp ,\Phi\rangle_{\tilde \tau (\H)}.
 \end{equation*}
 By Stokes' theorem,   the last expression is  equal to
 \begin{equation*}
  \langle   (\tilde\tau_* S)^\sharp ,\Phi\rangle_{\partial [\tilde \tau (\H)]}- \langle   (\tilde\tau_* S)^\sharp ,d\Phi\rangle_{\tilde \tau (\H)},
 \end{equation*}
 which is, by  coming back to $\H$ and $\partial H$ via  $\tilde\tau,$ equal to 
 $$
\big \langle  \tilde\tau^*[ (\tilde\tau_* S)^\sharp]  ,\tilde\tau^*\Phi \big \rangle_{\partial \H}
 - \big \langle  \tilde\tau^*[ (\tilde\tau_* S)^\sharp]  ,\tilde\tau^*(d\Phi) \big \rangle_{ \H}.
 $$ Hence, the result follows.
\endproof

The following   result  is  the  analogue of   Proposition \ref{P:basic-bdr-estimates}
in this  section.

\begin{proposition}\label{P:basic-bdr-estimates-non-Kaehler}
Fix  $\ell$  with $1\leq \ell\leq \ell_0$ and   set $\tilde\tau:=\tilde\tau_\ell.$ For $r\in(0,\bfr],$   set $\H_r:=\Tube (\widetilde V_\ell,r)\subset \E.$
Let $S$ be a   positive  current of bidimension $(q,q)$ such that $\dbar S$ is  a current of order $0.$
Let $\Phi$ be  the  product of $\theta_\ell$ and a  smooth $(q,q-1)$-form on $\Tube(B,\bfr)$ 
which is  $(2j+1)$-negligible. 
 Then there are 
 \begin{itemize}
 \item [$\bullet$]  two  functions $\Ic_1,\ \Ic_2 :\ (0,\bfr]\to\R;$ 
  \item [$\bullet$]   three differential operators $D_{10},$ $D_{11},$  $D_{12}$ in  the class $\widehat\Dc^0_\ell;$
   and three differential operators $D_{20},$ $D_{21},$  $D_{22}$   in  the class $\Dc^0_\ell;$
  \item[$\bullet$]  three smooth $2q$-forms $\Phi_{10}$   which is  $(2j-1)$-negligible,
   $\Phi_{11}$  which  is  $2j$-negligible,  $\Phi_{12}$  which  is  $(2j-1)$-negligible;
   and three smooth $2q$-forms $\Phi_{20}$   which is  $2j$-negligible,
   $\Phi_{21}$  which  is  $(2j+1)$-negligible,  $\Phi_{22}$  which  is  $2j$-negligible;
 \end{itemize}
  such that 
 every $0<r_1<r_2\leq\bfr$ and  every smooth  function $\chi$ on $(0,\bfr),$  we have for $\nu\in\{1,2\},$
 \begin{equation}\label{e:Stokes-ddc-difference-non-Kaehler}
 \begin{split}
  \int_{r_1}^{r_2}\chi(t) \Ic_\nu(t)dt&= 
\int_{\Tube(B,r_1,r_2)}\chi(\|y\|) (D_{\nu 1}S\wedge \Phi_{\nu1})(y)+\int_{\Tube(B,r_1,r_2)}\chi'(\|y\|) (D_{\nu2}S\wedge \Phi_{\nu2})(y)\\
&+ \int_{\partial_\hor\Tube(B,r_2) }\chi(r_2)( D_{\nu0}S\wedge \Phi_{\nu0})(y)- \int_{\partial_\hor\Tube(B,r_1)} \chi(r_1)(D_{\nu 0}S\wedge \Phi_{\nu0})(y)  ,
 \end{split}
 \end{equation}
 and  that the  following inequality holds 
 for all $  0<r\leq \bfr:$
  \begin{equation}\label{e:inequal-ddc-difference-non-Kaehler}
 {1\over  r^{2(k-p-j)} } \int_{r\over 2}^r\big| \langle  \dbar (\tilde\tau_* S) -\tilde\tau_*(\dbar  S),\Phi  \rangle_{\tilde \tau (\H_t)}
 -  \Ic_1(t) -\Ic_2(t)\big|dt
  \leq 
\sum_{m=\lowm}^\upm \nu_m( S,B,r,\id).
 \end{equation}
\end{proposition}
\proof
We  argue as in  the proof of Proposition \ref{P:basic-bdr-estimates} using  
 Proposition \ref{P:Stokes-non-Kaehler} instead of Proposition \ref{P:Stokes}.  
\endproof

  As in  Subsection
\ref{SS:Preliminaries-Section:psh-and-quasi-mono} we
 recall some notation from  the Extended Standing  Hypothesis in Subsection \ref{SS:Ex-Stand-Hyp}
  
Let $\omega$ be  a Hermitian  metric  on $V$ such that $\ddc \omega^j=0$ for all $1\leq j\leq \upm-1.$  
Fix an integer $j$ with $\lowm\leq j\leq \upm.$ 
Consider  the  forms on $\bfU$:
\begin{equation}\label{e:partition-can-forms-j,l}
\Phi:= \pi^*(  \partial (\omega^j))\wedge \beta^{k-p-j-1} \quad\text{and}\quad \Phi^{(\ell)}:= (\pi^*\theta_\ell)\cdot\pi^*(\partial (\omega^j))\wedge \beta^{k-p-j-1}\quad\text{for}\quad
1\leq\ell\leq \ell_0.
\end{equation}
So we have
\begin{equation}\label{e:sum-forms-j,l}
\Phi=\sum_{\ell=1}^{\ell_0} \Phi^{(\ell)}\qquad\text{on}\qquad \bfU.
\end{equation}
For $\ell$  with $1\leq \ell\leq \ell_0$ and   set $\tilde\tau:=\tilde\tau_\ell.$
For $r\in(0,\bfr],$   set $\H_r:=\Tube (\widetilde V_\ell,r)\subset \E.$

Let $T$ be a positive  closed   current on $\bfU$ in the class  $\widetilde\CL^{2,2}_p(\bfU,\bfW).$ 
Consider the  current
\begin{equation}\label{e:S_ell-non-Kaehler}
 S^{(\ell)}:=(\tau_\ell)_*(T|_{\bfU_\ell})   .
\end{equation}
By \eqref{e:T-hash} we get that
\begin{equation}
 T^\hash =\sum_{\ell=1}^{\ell_0} (\pi^*\theta_\ell)\cdot S^{(\ell)}.
\end{equation}
Note that the current $S^{(\ell)}$ is  positive plurisubharmonic on $\H_\bfr.$ Moreover, by Lemma
\ref{L:ex-negligble},   $\Phi^{(\ell)}$  is a $(2j+1)$-negligible smooth form. 
By Proposition \ref{P:basic-bdr-estimates-non-Kaehler},
there are 
 \begin{itemize}
 \item [$\bullet$]  two  functions $\Ic^{(\ell)}_1,\ \Ic^{(\ell)}_2 :\ (0,\bfr]\to\R;$ 
  \item [$\bullet$]   three differential operators $D^{(\ell)}_{10},$ $D^{(\ell)}_{11},$  $D^{(\ell)}_{12}$ in  the class $\widehat\Dc^0_\ell;$
   and three differential operators $D^{(\ell)}_{20},$ $D^{(\ell)}_{21},$  $D^{(\ell)}_{22}$   in  the class $\Dc^0_\ell;$
  \item[$\bullet$]  three smooth $2q$-forms $\Phi^{(\ell)}_{10}$   which is  $(2j-1)$-negligible,
   $\Phi^{(\ell)}_{11}$  which  is  $2j$-negligible,  $\Phi^{(\ell)}_{12}$  which  is  $(2j-1)$-negligible;
   and three smooth $2q$-forms $\Phi^{(\ell)}_{20}$   which is  $2j$-negligible,
   $\Phi^{(\ell)}_{21}$  which  is  $(2j+1)$-negligible,  $\Phi^{(\ell)}_{22}$  which  is  $2j$-negligible;
 \end{itemize}
  such that 
 every $0<r_1<r_2\leq\bfr$ and  every smooth  function $\chi$ on $(0,\bfr],$  we have for $\nu\in\{1,2\},$
 \begin{equation}\label{e:Stokes-ddc-difference-bis-non-Kaehler}
 \begin{split}
  \int_{r_1}^{r_2}\chi(t) \Ic^{(\ell)}_\nu(t)dt&= 
\int_{\Tube(B,r_1,r_2)}\chi(\|y\|) (D^{(\ell)}_{\nu 1}S^{(\ell)}\wedge \Phi^{(\ell)}_{\nu 1})(y)+\int_{\Tube(B,r_1,r_2)}\chi'(\|y\|) (D^{(\ell)}_{ \nu 2}S^{(\ell)}\wedge \Phi^{(\ell)}_{\nu 2})(y)\\
&+ \int_{\partial_\hor\Tube(B,r_2) }\chi(r_2)( D^{(\ell)}_{\nu 0}S^{(\ell)}\wedge \Phi^{(\ell)}_{\nu 0})(y)- \int_{\partial_\hor\Tube(B,r_1)} \chi(r_1)(D^{(\ell)}_{\nu 0}S^{(\ell)}\wedge \Phi^{(\ell)}_{\nu 0})(y)  ,
 \end{split}
 \end{equation}
 and  that the  following inequality holds 
 for all $  0<t\leq \bfr:$
  \begin{equation}\label{e:inequal-ddc-difference-bis-non-Kaehler}
 {1\over  r^{2(k-p-j)} } \int_{r\over 2}^r\big| \langle  \dbar [(\tilde\tau_\ell)_* S^{(\ell)}] -(\tilde\tau_\ell)_*(\dbar  S^{(\ell)}),\Phi^{(\ell)}  \rangle_{\tilde \tau (\H_t)}
 -  \Ic^{(\ell)}_1(t) -\Ic^{(\ell)}_2(t)\big|dt
  \leq 
\sum_{m=\lowm}^\upm \nu_m( S^{(\ell)},B,r,\id).
 \end{equation} 
 Note that by \eqref{e:S_ell-non-Kaehler},
 $\dbar S^{(\ell)}=0$    since $T$ is a  closed  $(p,p)$-current.

 The following  auxiliary  results are needed.
 \begin{lemma}\label{L:partition-non-Kaehler}
  The  following   equalities hold:
  \begin{eqnarray*}
     (\tilde\tau_\ell)_* S^{(\ell)}&=&\tau_*T \qquad \text{and}\quad (\tilde\tau_\ell)_*(\dbar S^{(\ell)})=\tau_*(\dbar T)\quad\text{on}\quad \bfU_\ell,\\
     \sum_{\ell=1}^{\ell_0}  \dbar[ (\tilde\tau_\ell)_* S^{(\ell)}]\wedge \Phi^{(\ell)}  &=&\dbar (\tau_*T)\wedge \Phi\quad\text{and}\quad
      \sum_{\ell=1}^{\ell_0} (\tilde\tau_\ell)_*(  \dbar S^{(\ell)})\wedge \Phi^{(\ell)}  = \tau_*(\dbar T)\wedge \Phi
      \quad\text{on}\quad \bfU.
  \end{eqnarray*}
 \end{lemma}
\proof 
We argue as in the proof of Lemma \ref{L:partition}.  
 \endproof

 \begin{lemma}\label{L:Ic-0-non-Kaehler}
  Under  the above  hypotheses and  notations,  there is a constant $c$ independent of $T$  such that for $\nu\in \{1,2\}$ and for all $1\leq \ell\leq \ell_0$ and for all $0<r\leq \bfr:$
  \begin{equation*}
   {1\over  r^{2(k-p-j)} } \int_{r\over 2}^r\big|  \int_{\partial_\hor\Tube(B,t) }( D^{(\ell)}_{\nu0}S^{(\ell)}\wedge \Phi^{(\ell)}_{\nu0}) \big|dt\leq  cr^2\Mc^\tot(T,r).
  \end{equation*}
 \end{lemma}
\proof  We argue as in the proof of Lemma \ref{L:Ic-0}.
\endproof 
  \begin{lemma}\label{L:Ic-1-and-2-non-Kaehler}
Under  the above  hypotheses and  notations, let  $0<r\leq \bfr.$
  Then  there is a constant $c$ independent of $T$ and $r$ such that for $\nu\in \{1,2\}$ and for all $1\leq \ell\leq \ell_0$ and for
  all $0<s<r:$
  \begin{eqnarray*}
  \big |\int_{\Tube(B,s,r)}\chi(\|y\|) (D^{(\ell)}_{\nu 1}S^{(\ell)}\wedge \Phi^{(\ell)}_{\nu1})(y)\big|  &\leq & c\sum_{n=0}^\infty {r\over 2^n}\Mc^\tot(T,{r\over 2^n}),\\
  \big|  \int_{\Tube(B,s,r)}\chi'(\|y\|) (D^{(\ell)}_{\nu2}S\wedge \Phi^{(\ell)}_{\nu2})(y) \big|&\leq & c\sum_{n=0}^\infty {r\over 2^n}\Mc^\tot(T,{r\over 2^n}).
  \end{eqnarray*}
  Here $\chi$ is  either the function $\chi_1$ or the function $\chi_2$ given in \eqref{e:chi_1-chi_2}.
 \end{lemma}
\proof   We argue as in the proof of Lemma \ref{L:Ic-1-and-2}.
\endproof

\subsection{Finiteness of the mass indicator $\Kc_{j,q}$ }

Let $\omega$ be  a Hermitian  metric  on $V$ such that $\ddc \omega^j=0$ for all $1\leq j\leq \upm-1.$

Fix  an open neighborhood $\bfW$ of $\partial B$ in $X$ with $\bfW\subset \bfU.$
Recall  the class $\widetilde\CL^{2,2}_p(\bfU,\bfW)$  given in Definition
 \ref{D:sup}.
  
  The following result  states    the main  difference  with   the siuation where $\omega$ is  K\"ahler.
  \begin{lemma}\label{L:non-Kaehler}
  Let  $T$ be   a  closed  $(p,p)$-current on $\bfU.$
  Then
  \begin{equation*}
 \ddc [(\tau_*T)\wedge \pi^*\omega^{j_0} ]
= (\dbar \tau_*T)\wedge \pi^*(\partial \omega^{j_0})=  (\dbar (\tau_*T)-\tau_*(\dbar T))\wedge \pi^*(\partial \omega^{j_0}).
\end{equation*}
  \end{lemma}
  \proof
  Since $T$ is closed and  $\omega$ is  $\ddc$-closed, it follows that $d T=\dbar T=0.$ A straightforward calculation  gives the  desired identity.
  \endproof

\begin{lemma}\label{L:ddc-difference-lambda-non-Kaehler} For  all $r_1,r_2\in  (0,\bfr]$ with  $r_1<r_2,$  
there is a constant $c>0$  such that  for every
$j$ with $\lowm\leq j\leq \upm,$  and every $m$ with $0\leq m\leq j,$ and every positive plurisubharmonic  current $T$   in the class  $\widetilde\CL^{2,2}_p(\bfU,\bfW),$
  the following two inequalities hold for all $\lambda\geq 1:$
\begin{multline*}
\big| \int_{r_1}^{r_2} \big( {1\over t^{2(k-p-j)}}-{1\over r_2^{2(k-p-j)}}  \big)2tdt\int_{\Tube(B,t)} (A_{\lambda})_*\big(\ddc [(\tau_*T)\wedge \pi^*(\omega^{j-m})]\wedge \beta^{k-p-j+m-1}\big) \big| \\
 \leq {c\over \lambda^{2m+1}}\Mc^\tot(  T,{r_2\over \lambda_n} ),
 \end{multline*}
 \begin{multline*}
  \big( {1\over r_1^{2(k-p-j)}}-{1\over r_2^{2(k-p-j)}}  \big) 
 \big|\int_{0}^{r_1}2tdt \int_{\Tube(B,t)} (A_{\lambda})_*\big(\ddc [(\tau_*T)\wedge\pi^*(\omega^{j-m})]\wedge \beta^{k-p-j+m-1}\big)\big|\\
  \leq {c\over \lambda^{2m+1}_n}\Mc^\tot( T,{r_1\over \lambda} ).
\end{multline*} 
 \end{lemma}
 \proof
 We only give the proof of the  first  inequality since  the second one can be  obtained   similarly.
 By Lemma \ref{L:non-Kaehler}, the first inequality is  equivalent to
 \begin{multline*}
\big| \int_{r_1}^{r_2} \big( {1\over t^{2(k-p-j)}}-{1\over r_2^{2(k-p-j)}}  \big)2tdt\int_{\Tube(B,t)} (A_{\lambda})_*\big([\dbar (\tau_*T)-\tau_*(\dbar T)]\wedge \partial [\pi^*(\omega^{j-m})]\wedge \beta^{k-p-j+m-1}\big) \big| \\
 \leq {c\over \lambda^{2m+1}}\Mc^\tot(  T,{r_2\over \lambda_n} ),
 \end{multline*}
 But  this inequality is true by applying Proposition \ref{P:basic-bdr-estimates-non-Kaehler}  and
Theorem \ref{T:nu_tot-monotone}.
 \endproof

\begin{lemma}\label{L:basic-difference-nu_j,q_estimate-non-Kaehler} Given $0<r_1<r_2\leq \bfr,$ 
there is  a constant $c>0$ such that  for every positive closed current $T\in  \widetilde\SH^{3,3}_p(\bfU,\bfW)$
and  $0\leq q\leq k-l$ and $0\leq j\leq \min(\upm,k-p-q),$ the following inequality holds:
 \begin{multline*}
\nu_{j,q}\big(T,B,{r_2\over \lambda},\tau\big)  -\nu_{j,q}\big(T,B,{r_1\over \lambda},\tau\big)
\geq  \Kc_{j,q}\big (T,{r_1\over\lambda},{r_2\over \lambda}  \big)-  c\lambda^{-1} -  c\lambda^{-1}\Mc^\tot(T,{r_2\over \lambda})
   \\-  c\lambda^{1\over 2} \Kc_{q}\big (T,{r_1\over\lambda},{r_2\over \lambda} \big )-c\Kc_{q-1} (T,{r_1\over\lambda},{r_2\over \lambda})  
-c\sqrt{\Kc_{q}(T,{r_1\over\lambda},{r_2\over \lambda})} \sqrt{ \Kc^-_{j,q}(T,{r_1\over\lambda},{r_2\over \lambda})  }.
\end{multline*}
\end{lemma}
\proof
Fix  $0\leq q_0\leq k-l.$   
Let $0\leq j_0\leq \min(\upm,k-p-q_0).$   Set  $j'_0:=k-p-q_0-j_0\geq 0.$
We may assume without loss of generality that $T$ is  $\Cc^3$-smooth.   
Applying Theorem  \ref{T:Lelong-Jensen-smooth}  to  $\tau_*T\wedge \pi^*(\omega^{j_0})\wedge \beta^{j'_0}$
and noting that $\beta$ is  closed, we get that
\begin{equation*}
\begin{split}
&{\lambda^{2q_0}\over  r_2^{2q_0}}\int_{\Tube(B,{r_2\over \lambda})}\tau_*T\wedge \pi^*(\omega^{j_0})\wedge \beta^{k-p-j_0}
-{\lambda^{2q_0}\over  r_1^{2q_0}}\int_{\Tube(B,{r_1\over\lambda})}\tau_*T\wedge \pi^*(\omega^{j_0})\wedge \beta^{k-p-j_0}\\
&= \Vc\big(\tau_*T\wedge \pi^*(\omega^{j_0})\wedge \beta^{j'_0},{r_1\over\lambda},{r_2\over \lambda}\big)+\int_{\Tube(B,{r_1\over\lambda},{r_2\over \lambda})}\tau_*T\wedge \pi^*(\omega^{j_0})\wedge \beta^{j'_0}\wedge\alpha^{q_0}\\
&+  \int_{r_1\over\lambda}^{r_2\over\lambda} \big( {1\over t^{2q_0}}-{\lambda^{2q_0}\over r_2^{2q_0}}  \big)2tdt\int_{\Tube(B,t)} \ddc [(\tau_*T)\wedge \pi^*(\omega^{j_0})]\wedge \beta^{q_0+j'_0-1}  \\
&+\big( {\lambda^{2q_0}\over r_1^{2q_0}}-{\lambda^{2q_0}\over r_2^{2q_0}}  \big) 
 \int_{0}^{r_1\over \lambda }2tdt \int_{\Tube(B,t)} \ddc[ (\tau_*T)\wedge \pi^*(\omega^{j_0})]\wedge \beta^{q_0+j'_0-1} .
\end{split}
\end{equation*}
By Lemma \ref{L:ddc-difference-lambda-non-Kaehler},
the  last  two double  integrals  are  of order smaller than  
$
{c \lambda^{-1}}\Mc^\tot( T,{r_2\over \lambda} ) 
.$ Moreover,  by Theorem \ref{T:vertical-boundary-terms}, we have the following estimate independently of $T:$
\begin{equation*}
\Vc\big(\tau_*T\wedge \pi^*(\omega^{j_0})\wedge \beta^{j'_0},{r_1\over\lambda},{r_2\over \lambda}\big)=O(\lambda^{-1}).
\end{equation*}
Therefore,  there is a constant $c>0$ independent of $T$ such that for $\lambda\geq 1,$
\begin{multline*}
\big|\int_{\Tube(B,{r_1\over\lambda},{r_2\over \lambda})}\tau_*T\wedge \pi^*(\omega^{j_0})\wedge \beta^{j'_0}\wedge\alpha^{q_0}
-\big (\nu_{j_0,q_0}\big(T,B,{r_2\over \lambda},\tau\big)   -  \nu_{j_0,q_0}\big(T,B,{r_1\over \lambda},\tau\big) \big)  \big| \\
\leq  c\lambda^{-1}+ {c \lambda^{-1}}\Mc^\tot( T,{r_2\over \lambda} ).
\end{multline*} 
The remainder of the  proof follows along the same lines as those given in the proof of  Lemma 
\ref{L:basic-difference-nu_j,q_estimate-psh}.
\endproof

Theorem \ref{T:nu_tot-monotone} is still valid in this  more general context. For the  reader convenience,
we  record  here the  new  statement
  
\begin{theorem}\label{T:nu_tot-monotone-non-Kaehler} 
Let  $0< r_1<r_2\leq\bfr.$ Then there are  a  family  $\Dc=\{ d_{jq}\in\R:\ 0\leq j\leq k-p-q,\ 0\leq q\leq  k-l\}$ and   a constant $c>0 $ depending  on $r_1$ and $r_2$  such that for every positive closed   current $T$ on $\bfU$ belonging to the class $\widetilde\CL^{2,2}_p(\bfU,\bfW),$ 
the  following inequality  hold for $0\leq q\leq \upm:$
\begin{equation*}
  \nu^\Dc_q\big(T,B,{r_1\over \lambda},\tau\big)\leq  \nu^\Dc_q(T,B,{r_2\over \lambda},\tau\big) + {c\over \lambda}
  \quad\text{for}\quad \lambda\gg 1.
 \end{equation*}
\begin{equation*}   
\nu^\Dc_\tot(T,B, r,\tau)\leq  c\Mc^\tot(T,r)\quad\text{and}\quad
c^{-1}\Mc^\tot(T,r)\leq \nu^\Dc_\tot(T,B, r,\tau)+ cr  
\quad\text{for}\quad 0<r\leq\bfr.
\end{equation*}
\end{theorem}

  The main technical  result of this  section is  the following
  
 \begin{theorem}\label{T:Lc-finite-non-Kaehler}
   There is  a constant $c_7>0$ such that for every positive closed  current $T$ on $\bfU$ belonging to the class $\widetilde\CL^{2,2}_p(\bfU,\bfW),$  we have 
   \begin{equation}  \label{e:Lc-finite-non-Kaehler-main}
   \Kc_{j,q}(T,r) \leq  c_7 \nu^\tot(T,B,r,\tau)
   \end{equation}
   for $0\leq q\leq k-l$ and $0\leq j\leq k-p-q.$
   Here  $\nu^\tot(T,B,r,\tau)$ is  defined by \eqref{e:nu_tot}.
   In particular,
   $\Kc_{j,q}(T,\bfr)<c_7.$
 \end{theorem}
\proof   
We indicate how to adapt the proof of Theorem  \ref{T:Lc-finite-psh} in the  present context.
The proof is also divided into three steps.

\noindent {\bf Step 1:} {\it The case  $q=0.$}

This step   is similar to that of the proof of Theorem  \ref{T:Lc-finite-psh}.

 The general strategy is  to prove  the proposition by  increasing induction on $q$ 
 with $0\leq q\leq k-l.$ But the  induction  procedure  is  somehow  simpler   than that of Theorem  \ref{T:Lc-finite-psh}.
 In the proof $\bfr$ is  a  fixed   but sufficiently small  constant.
Fix  $0\leq q_0\leq k-l.$ 
Suppose   that \eqref{e:Lc-finite-psh-main} is true  for all $q,j$ with $q<q_0.$ 
We need to show that   it is  also true  for all $q,j$ with   $q\leq q_0.$ 
We may assume without loss of generality that $T$ is  a $\Cc^2$-smooth closed $(p,p)$-form
and  let $s,r\in  [0,\bfr)$ with $s<r.$

{\it  Set $m_0:= k-p-q,$  $m_1:=m_0-1.$ 
In the first induction we will prove that there is a constant $c_{10}$ independent of $T$ and $r$
such that 
\begin{equation}\label{e:first-ind-non-Kaehler}
\Kc^\bullet_{q}(T,r)
\leq c_{10}\Nc^\bullet_q(T,r)\qquad\text{and}\qquad \Lc_{j,q}(T,r)\leq c_{10}\Nc^\bullet_q(T,r),
\end{equation}
for every $0\leq q\leq k-l,$ $j\geq 0$  with $j\leq m_1,$  and for every $0<r\leq \bfr.$
Here,
\begin{eqnarray*}
\Nc^\bullet_q(T,r)&:=&r+\Kc_{q-1}(T,r)+ \Lc_{q-1}(T,r)+\sum_{j=0}^{m_1} |\nu_{j,q}(T,r)|,\\
\Lc_{q}(T,r)&:=&\sum_{j,q':\  q'\leq q\quad\text{and}\quad j+q'\leq k-p} \Lc_{j,q'}(T,r).
\end{eqnarray*}

The proof of \eqref{e:first-ind} will be completed in Steps 2 and 3 below.
}

\noindent {\bf Step 2:} {\it  Let $q_0:=q$ and define $m_0$ and $m_1$ as  above using $q_0$ instead $q.$ There is a constant $c_{10}>0$ such that for every  $j_0,q_0\geq 0$ with $j_0\leq  m_1$ and   every $0<r\leq \bfr,$
\begin{equation}\label{e:P-Lc-bullet_finite(8)-non-Kaehler}
\begin{split}
I^\hash_{q_0,0,j_0,0}(T,r)&\leq  c_{10}\big(|\nu_{j_0,q_0}(T,B,r,\tau)|+r +r^{1\over 4}\Mc^\tot(T,r)+r^{1\over 4}\Kc^{+}_{j_0,q_0}(T,r) +r^{1\over 4} \Kc_{q_0}^\bullet(T,r)\\
&+  \sqrt{\Kc^\bullet_{q_0}(T,r)} \sqrt{ \Kc^-_{j_0,q_0}(T,r)  }\big).
\end{split}
\end{equation}
where the expression on the LHS is  given by \eqref{e:I_bfj}  (see also Remark \ref{R:I_bfj}).} 
 
Let $0\leq j_0\leq \min(\upm, k-p-q_0).$ Set 
 $j'_0:=k-p-q_0-j_0\geq 0$ and $m_0:=k-p-q_0.$
Suppose that $j'_0\geq 1.$

By Lemma \ref{L:non-Kaehler} we have that
\begin{equation*}
 \ddc [(\tau_*T)\wedge \pi^*\omega^{j_0}\wedge \beta^{j'_0}]
= (\dbar \tau_*T)\wedge \pi^*(\partial \omega^{j_0})\wedge \beta^{j'_0}=  (\dbar (\tau_*T)-\tau_*(\dbar T))\wedge \pi^*(\partial \omega^{j_0})\wedge \beta^{j'_0}.
\end{equation*}
Applying Theorem  \ref{T:Lelong-Jensen-smooth}  to  $\tau_*T\wedge \pi^*(\omega^{j_0}) $  and using  the  above  equality, we  argue as in the  end of  Step 2  of the proof  of Theorem  \ref{T:Lc-finite-psh}.

\noindent {\bf Step 3:} {\it  
End of the proof of
\eqref{e:first-ind-non-Kaehler}.}
We  argue as in  Step 3  of the proof  of Theorem  \ref{T:Lc-finite-psh}.

{\it 
Now it remains to treat the case  where $j=m_0:=k-p-q,$ that is,  there is a constant $c_{10}$ independent of $T$
such that 
\begin{equation}\label{e:second-ind-non-Kaehler}
\Kc_{q}(T,r)
\leq c_{10}\Mc^\tot(T,r)\qquad\text{and}\qquad \Lc_{j,q}(T,r)\leq c_{10}\Mc^\tot(T,r),
\end{equation}
for every $0\leq q\leq \min(k-l,k-p).$ 
The proof of \eqref{e:second-ind-non-Kaehler}  will be completed in Steps 4 and 5 below.
By Steps  2 and 3,  inequality  \eqref{e:second-ind-non-Kaehler} is reduced to proving that
\begin{equation}\label{e:second-ind-bis-non-Kaehler}
\Kc_{m_0,q}(T,r)
\leq c_{10}\Mc^\tot(T,r)\qquad\text{and}\qquad \Lc_{m_0,q}(T,r)\leq c_{10}\Mc^\tot(T,r).
\end{equation}
}

\noindent {\bf Step 4: }{\it  Inequality \eqref{e:second-ind-bis-non-Kaehler}
 holds 
for every $0\leq q<  k-p-\upm.$ 
}
     
   We  argue as in the  end of  Step 4  of the proof  of Theorem  \ref{T:Lc-finite-psh}. 

\noindent {\bf Step 5: }{\it  Inequality \eqref{e:second-ind-non-Kaehler}
 holds 
for every $k-p-\upm\leq q\leq  k-p-\lowm.$ 
}
 
 We  argue as in the  end of  Step 5  of the proof  of Theorem  \ref{T:Lc-finite-psh}.

\endproof

\begin{proposition}\label{P:Lc-finite-bis-non-Kaehler}
  For  $0<r_1<r_2\leq \bfr,$ there is  a constant $c_8>0$ such that for every $q\leq  \min(k-p,k-l)$ and 
    every positive closed  current $T$ on $\bfU$ belonging to the class $\widetilde\CL^{2,2}_p(\bfU,\bfW),$   we have the following estimate:
 $$|\kappa_{k-p-q}(T,{r_1\over \lambda},{r_2\over \lambda},\tau)|<c_8\sum\limits_{0\leq q'\leq q,\ 0\leq j'\leq \min(\upm,k-p-q')} \Kc_{j',q'}(T,{r_1\over \lambda},{r_2\over \lambda} )\qquad\text{for}\quad  \lambda>1.$$
 \end{proposition}
\proof
We argue as in the proof of Proposition \ref{P:Lc-finite-bis}  making the necessary changes.
 \endproof
We conclude  this  subsection  with the following finiteness result of the mass indicators $\Mc_j.$  
 \begin{proposition}\label{P:Lc-finite-closed-non-Kaehler}
 There is a   constant $c_9>0$ such that for  
  every positive closed  current $T$ on $U$ belonging to the class $\widetilde\CL^{2,2}_p(\bfU,\bfW),$ 
   we have
   $\Mc_j(T,r)<c_{9}$ for $0\leq j\leq \upm$ and $0 <r \leq \bfr.$
 \end{proposition}
\proof Since  the proof is  not difficult, we leave it   to the interested reader.   \endproof

\subsection{End of the proof for positive closed currents}
This  subsection is  devoted to the  proof of Theorem \ref{T:Lelong-closed}  using  
Theorem \ref{T:Lc-finite-non-Kaehler}
and Proposition \ref{P:Lc-finite-bis-non-Kaehler}.

\proof[Proof of assertion (1) of  Theorem \ref{T:Lelong-closed}]
We may assume without loss of generality that $T$ is  $\Cc^2$-smooth.   
Applying Theorem  \ref{T:Lelong-Jensen-smooth}  to  $\tau_*T\wedge \pi^*(\omega^{j})$
and noting that $\beta$ is  closed, we get that
\begin{equation*}
\begin{split}
&{1\over  r_2^{2(k-p-j)}}\int_{\Tube(B,r_2)}\tau_*T\wedge \pi^*(\omega^{j})\wedge \beta^{k-p-j}
-{1\over  r_1^{2(k-p-j)}}\int_{\Tube(B,r_1)}\tau_*T\wedge \pi^*(\omega^{j})\wedge \beta^{k-p-j}\\
&= \Vc\big(\tau_*T\wedge \pi^*(\omega^{j}),r_1,r_2\big)+\int_{\Tube(B,r_1,r_2)}\tau_*T\wedge \pi^*(\omega^{j})\wedge\alpha^{k-p-j}\\
&+  \int_{r_1}^{r_2} \big( {1\over t^{2(k-p-j)}}-{1\over r_2^{2(k-p-j)}}  \big)2tdt\int_{\Tube(B,t)} \ddc [(\tau_*T)\wedge \pi^*(\omega^{j})]\wedge \beta^{k-p-j-1}  \\
&+\big( {1\over r_1^{2(k-p-j)}}-{1\over r_2^{2(k-p-j)}}  \big) 
 \int_{0}^{r_1}2tdt \int_{\Tube(B,t)} \ddc[ (\tau_*T)\wedge \pi^*(\omega^{j})]\wedge \beta^{k-p-j-1} .
\end{split}
\end{equation*}
By Lemma \ref{L:ddc-difference-lambda-non-Kaehler},
the  last  two double  integrals  are  of order smaller than  
$
 c r_2\Mc^\tot( T,r_2 ) 
.$ Moreover,  by Theorem \ref{T:vertical-boundary-terms}, we have the following estimate independently of $T:$
\begin{equation*}
\Vc\big(\tau_*T\wedge \pi^*(\omega^{j}),r_1,r_2\big)=O(r_2).
\end{equation*}
Therefore,  there is a constant $c>0$ independent of $T$ such that  
\begin{equation*}
\big|   \nu_{j}(T,B,r_2,\tau)-   \nu_{j}(T,B,r_1,\tau)-  \kappa_{j}(T,B,r_1,r_2,\tau)\big| 
\leq  cr_2+ cr_2\Mc^\tot( T,r_2 ).
\end{equation*} 
This, combined  with   Proposition \ref{P:Lc-finite-closed-non-Kaehler}, implies the  result.
\endproof
\proof[Proof of assertions (2)--(6) of  Theorem \ref{T:Lelong-closed}]
It follows along  the same lines as those given in the proof of Theorem \ref{T:Lelong-closed-Kaehler}.

\section{Existence of tangent currents}
\label{S:Existence-tangent-currents}
Recall the  Standing Hypothesis from Subsection \ref{SS:Global-setting}.
The  main purpose of this section is  to prove  the existence of tangent  currents 
in the following  three cases:
positive closed currents, positive harmonic  currents  and plurisubharmonic   currents.

\subsection{Positive closed  currents}

The main goal of this  subsection is  to prove  the following 
\begin{theorem}\label{T:Tangent-closed}
   We  keep the     Standing Hypothesis.  
Suppose  that    the $(p,p)$-current   $T$  is  positive closed  and $T=T^+-T^-$ on an open neighborhood of $\overline B$ in $X$  with
$T^\pm$  in the class $\CL_p^{1,1} (B).$  Suppose in addition   that   $\omega$ is  a K\"ahler form on $V$ and  that there is  at least  one strongly  admissible map
along $B.$
  Then the following  assertions hold: 
  \begin{enumerate}

 \item   Consider   a collection of 
admissible maps $\tau_\ell :\ U_\ell \to \U_\ell:=\tau_\ell(U_\ell)\subset \E $  along $B\cap U_\ell$  for $\ell$  in an index set  $ L$
with  $\overline B\subset \bigcup_{\ell\in L} U_\ell.$
 Then, for every $\ell\in L,$   the  family of currents  $T_{\lambda,\ell}:=(A_\lambda)_*(\tau_\ell)_*(T)$    with $\lambda\in\C^*$  which are  defined on  $\pi^{-1}(B\cap U_\ell)\subset \E$  is  relatively compact.
 In particular, if  $L$  is  at most  countable, then  for every sequence
 $(\lambda_n ) \subset \C^*$ converging to $\infty,$ we can extract a subsequence $(\lambda_{m_n}) \subset \C^*$  such that
 the tangent current to $T$ along $B$  associated   to  the sequence $(\lambda_{m_n})$ and  the collection  $(\tau_\ell)_{\ell\in L}$ 
 in the  sense of Definition \ref{D:tangent-currents} exists.
 
 \item  Consider   another  collection of 
admissible maps $\tau'_{\ell'} :\ U'_{\ell'} \to \U'_{\ell'}:=\tau'_{\ell'}(U'_{\ell'})\subset \E $  along $B\cap U'_{\ell'}$  for $\ell'$  in an index set  $ L'$ with  $\overline B\subset \bigcup_{\ell'\in L'} U'_{\ell'}.$ As in assertion (1)  consider the family of currents $T'_{\lambda,\ell'}:=(A_\lambda)_*(\tau'_{\ell'})_*(T)$    with $\lambda\in\C^*$  which are  defined on  $\pi^{-1}(B\cap U'_{\ell'})\subset \E.$ 
 Then  the  family of currents  $T_{\lambda,\ell} -T'_{\lambda,\ell'}$ converge  weakly to  $0$ on  $\pi^{-1}(B\cap U_\ell\cap U'_{\ell'})\subset \E$
 as  $\lambda$ tends to infinity.  

 \item Let  $T_\infty$ be  the tangent current  to $T$ along $B$ associated, by Definition \ref{D:tangent-currents},   to a sequence
 $(\lambda_n ) \subset \C^*$ converging to $\infty$  and the collection of 
admissible maps $(\tau_\ell)_{\ell\in L}$ in assertion (2).
Then  $T_\infty$ is  also the tangent current  to $T$ along $B$ associated, by Definition \ref{D:tangent-currents},   to the same sequence $(\lambda_n ) $   and the collection of 
admissible maps $(\tau'_{\ell'})_{\ell'\in L'}$ in assertion (3).  
 
\end{enumerate}
 \end{theorem}

Prior to the proof of this theorem some  auxiliary  results are needed.
Fix  a  holomorphic admissible map  $\tau:\ U^0\to  \E$ along $V\cap U^0,$ where $U^0$ is a  small  open subset of $U$  with $U^0\cap V\not=\varnothing.$
We  use  the notation introduced  in Subsection \ref{SS:Local-setting} and identify $U^0,$ via  a local holomorphic chart, with
   the unit polydisc  of $\C^k.$     
We use the holomorphic coordinate system  $y=(z,w)\in\C^{k-l}\times \C^l$ and write $U^0=U^0_z\times U^0_w,$ 
where $U^0_z$ (resp. $U^0_w)$ is the unit polydisc of  $\C^{k-l}$ of  (resp. of  $\C^l$). We may assume that 
 $V\cap U^0=\{z=0\}=\{0_z\}\times U^0_w.$ 
Consider  the trivial  vector bundles $\pi^\dagger:\ U^0_z\times \C^l \to  U^0_z$  and  $\pi:\ \E\to U^0_w$ with  $\E\simeq  \C^{k-l}\times U^0_w.$ For $\lambda\in\C^*,$  let $a_\lambda:\ \E\to \E$ be the multiplication by  $\lambda$
on fibers, that is, 
$a_\lambda(z,w):=(\lambda z,w)$ for $(z,w)\in \E.$
Consider the positive closed $(1,1)$-forms
$$\omega_z:=\ddc \|z\|^2\quad\text{and}\quad  \omega_w:=\ddc\|w\|^2.$$  
Recall from \eqref{e:m} the two nonnegative integers  $\lowm:=\max(0,l-p)$ and $\upm:=\min(l,k-p).$

\begin{lemma}\label{L:mass-estimates}
\begin{enumerate}
\item For every current $T$ of bidegree $(p,p)$ with  measure coefficients of bounded mass  and for $0\leq j\leq l,$  we have 
$$  (a_\lambda)_*\big(\pi^\dagger_* (T\wedge  \omega^j_w)\big) =  \pi^\dagger_* \big( (a_\lambda)_*(T\wedge  \omega^j_w )\big).$$

\item There is a  constant $c$ which depends uniquely on the dimension $k$ such that for  all positive  $(p,p)$-current $T$ on $U^0,$
$$\|(a_\lambda)_*T\|_{U^0}\leq  c \sum_{j=\lowm}^{\upm} \|  (a_\lambda)_*\big(\pi^\dagger_* (T\wedge  \omega^j_w)\big)\|_{U^0_z}.
 $$
 \end{enumerate}
\end{lemma}

\proof
{\bf Proof of assertion (1).}  We only need to prove the assertion for  $T$ of the form
$$T=\sum_{K,L} T_{I,J;K,L}dz_I\wedge d\bar z_J\wedge dw_K\wedge d\bar w_L,$$ where $ I,J\subset\{1,\ldots,k-l\}$ are   fixed, and the sum is taken over $K,L\subset \{1,\ldots,l\}$
with $|K|+|I|=|L|+|J|=p.$
Consider two cases.

\noindent {\bf  Case  $|I|=|J|=p-l+j:$ }  By a consideration of bidegree we have that
\begin{eqnarray*}(a_\lambda)_*\big(\pi^\dagger_*( T\wedge \omega^j_w)\big)&=&(a_\lambda)_*\Big(\sum_K \big( \int_{w\in\C^l}T_{I,J; K,K}(z,w)\wedge
dw_K\wedge d\bar w_K\wedge \omega^j_w\big) dz_I\wedge d\bar z_J\Big)\\
&=&|\lambda|^{-2(p-l+j)}\sum_{K} \big( \int_{w\in\C^l} T_{I,J;K,K} (z,w)
\omega^l_w\big)dz_I\wedge d\bar z_J,
\end{eqnarray*}
where the sums are taken  over all $K\subset\{1,\ldots,l\}$ and $|K|=l-j.$ 

On the  other hand,
\begin{eqnarray*}
 \pi^\dagger_* \big( (a_\lambda)_*(T\wedge  \omega^j_w )\big)&=&\lambda^{-|I|} \bar{\lambda}^{-|J|} \pi^\dagger_{*}\Big(\sum_K
T_{I,J; K,K}(z,w) dz_I\wedge d\bar z_J\wedge dw_K\wedge d\bar w_K\wedge \omega^j_w\Big)\\
 &=&|\lambda|^{-2(p-l+j)}\sum_{K} \big( \int_{w\in\C^l} T_{I,J;K,K} (z,w)
\omega^l_w\big)dz_I\wedge d\bar z_J,
\end{eqnarray*}
where the sums are taken  over all $K\subset\{1,\ldots,l\}$ and $|K|=l-j.$  So assertion (1)  is true in this case.

\noindent{\bf  Case either $|I|\not=p-l+j$ or $|J|\not=p-l+j$:}  By a consideration of bidegree we see that
$$
(a_\lambda)_*\big(\pi^\dagger_* (T\wedge  \omega^j_w)\big) = 0\quad\text{and}\quad \pi^\dagger_* \big( (a_\lambda)_*(T\wedge  \omega^j_w )\big)=0.
$$
Hence,  assertion (1)  follows.

\noindent {\bf Proof of assertion (2).}
As an immediate consequence of assertion (1), we may replace $(a_\lambda)_*T$ by $T,$ that is, we may assume that $\lambda=1.$
By Proposition \ref{P:Demailly}, we only need to  prove that for every fixed multi-index $I\subset\{1,\ldots, k-l\}$ and  every fixed multi-index $K\subset  \{1,\ldots,l\}$ with $|I|+|K|=p,$ 
\begin{equation}\label{e:local-mass-est}
\begin{split}
\|T_{I,I;K,K}\|_{U^0}&\leq c  \|  \pi^\dagger_* (T\wedge  \omega^j_w)\|_{U^0_z},\quad\text{where}\quad  j=l-|K|,\\
\|  \pi^\dagger_* (T\wedge  \omega^j_w)\|_{U^0_z}&=0\quad\text{for}\quad  j\not\in[\lowm,\upm].
\end{split}
\end{equation} 
To  prove  the inequality of \eqref{e:local-mass-est}, observe that
\begin{eqnarray*}
\|  \pi^\dagger_* (T\wedge  \omega^j_w)\|_{U^0}&\gtrsim& \|  \pi^\dagger_* (T_{I,I,K,K}dz_I\wedge d\bar z_I\wedge dw_K\wedge d\bar w_K\wedge  \omega^j_w)\|_{U^0_z}\\
&=& \|\big(\int_{w\in\C^l} T_{I,I;K,K} (z,w)
\omega^l_w\big)dz_I\wedge d\bar z_I\|_{U^0_z}\\
&=& \int_{z\in U^0}\big(\int_{w\in\C^l} T_{I,I;K,K} (z,w)
\omega^l_w\big)\omega_z^{k-l}=\|T_{I,I;K,K}\|_{U^0}.
\end{eqnarray*}
It remains to us  to prove the equality  of \eqref{e:local-mass-est}.
Since $T$ is of bidegree $(p,p),$   $T\wedge \omega_w^j$ is of bidegree $(p+j,p+j)$ and  hence it is  zero
if $p+j>k.$  Moreover,  $T\wedge \omega_w^j$ is zero if $j>l$  as  $\omega_w^{l+1}=0.$
So   $T\wedge \omega_w^j$ is zero if $j>\upm.$

On the other hand, $\pi^\dagger_* (T\wedge  \omega^j_w)=0$ if  $T\wedge \omega^j_w$ is  not of full bidegree $(l,l)$ in $\{dw,d\bar w\}$
and this is  the case if $p+j<l.$  So  $\pi^\dagger_* (T\wedge  \omega^j_w)=0$ for $j<\lowm.$
This completes the proof of the equality  of \eqref{e:local-mass-est}.
\endproof
   
  \proof[Proof of assertion (1) of Theorem  \ref{T:Tangent-closed}]
 We  fix  an $\ell\in L$  and  write $U'$ (resp.  $\tau$) (resp.  $T_\lambda$)  instead  of  $U_\ell$ (resp. $\tau_\ell$)  (resp.  $T_{\lambda,\ell}$).
 To  prove the compacness of  the family
 $(T_\lambda)_{\lambda\in\C^*},$
 we only need to show that the masses of  the currents of this family are locally uniformly bounded on $\pi^{-1}(V\cap U^0)\subset \E.$ 
Fix  an arbitrary $r_0>0,$  we need to show that   there is a constant $c_0$ such that
\begin{equation}\label{e:T_lambda-bounded-mass}|\langle  (A_\lambda)_* (\tau_*T),\Phi\rangle|\leq c_0
\end{equation}
 for every  continuous test form $\Phi$ supported in $\Tube(B,r_0)$
with $\|\Phi\|_{\Cc^0}\leq 1.$ 
Observe  that    for a given $\lambda_0>0,$ we can find  $c_0$ such that \eqref{e:T_lambda-bounded-mass} holds for 
$\lambda\in\C^*$ with $|\lambda|\leq \lambda_0.$  Therefore, we may assume  without loss of generality that $r_0\leq {1\over 2}\bfr,$ and  we only  need to prove  \eqref{e:T_lambda-bounded-mass} for   
$\lambda\in\C^*$ with $|\lambda|\leq 1.$

By Proposition \ref{P:Lc-finite-closed} and  Lemma \ref{L:Mc_j}, we have
 for $0\leq j\leq \upm,$  
 \begin{equation} \label{e:Mc-finite-local}   \sup\limits_{r\in (0,\bfr]} {1\over r^{2(k-p-j)}}\int_{\|z\|<r,\ w\in U^0_w} T^\hash \wedge \omega_w^j \wedge \omega_z^{k-p-j}<c_{10},
 \end{equation}
 where the  positive  current $T^\hash$ is  defined in \eqref{e:T-hash}.
Setting $r:={r_0\over |\lambda|},$ we infer from the above  inequality that
$$\sup_{\lambda\in\C:\  |\lambda|\geq 1} \|  (a_\lambda)_*\big(\pi^\dagger_* (T^\hash\wedge  \omega^j_w)\big)\|_{U^0}<\infty.$$
This, combined with  Lemma \ref{L:mass-estimates} (2), implies that   
\begin{equation}\label{e:T-hash-lambda-finite}\sup\limits_{|\lambda|\geq 1} \|(a_\lambda)_*(T^\hash)\|_{U_0}<\infty,
 \end{equation}
where $U_0:=\{(z,w)\in U^0:\ \|z\|<2r_0\}.$
By  \eqref{e:T-hash_r}, we  infer from  \eqref{e:T-hash-lambda-finite} that 
\begin{equation*}\sup\limits_{|\lambda|\geq 1} \|(a_\lambda)_*(T^\hash_{r_0\over \lambda})\|<\infty.
\end{equation*}
Hence, $|\langle T^\hash_{r_0\over |\lambda|}, (a_\lambda)^*\Phi\rangle|\lesssim c_0$ independent of $\Phi$ as  above  and of $\lambda.$ 
Applying  Lemma  \ref{L:basic-T-hash-bis} yields for $r:={r_0\over |\lambda|}$ that
 \begin{equation}\label{e:tangent-current-closed}\begin{split} &\langle \tau_*T,\ind_{\Tube(B,{r_0\over |\lambda|})}(a_\lambda)^*\Phi\rangle -\langle T^\hash_{r_0\over |\lambda|}, (a_\lambda)^*\Phi\rangle \\
 &=\sum_{\ell=1}^{\ell_0}\langle(\tau_\ell)_* T,  (\ind_{\Tube(B,{r_0\over |\lambda|})}\circ \tilde\tau_\ell)\cdot \big((\tilde\tau_\ell)^*((\pi^*\theta_\ell) (a_\lambda)^*\Phi)-((\pi^*\theta_\ell) (a_\lambda)^*\Phi)\big)\rangle ,
\end{split}
 \end{equation}
 where  we recall  from \eqref{e:tilde-tau_ell}  that $\tilde\tau_\ell:=\tau\circ \tau_\ell^{-1}.$
On the one hand, the RHS is rewritten as
$$
\sum_{\ell=1}^{\ell_0}\langle(a_\lambda)_*(\tau_\ell)_* T,  (\ind_{\Tube(B,{r_0\over |\lambda|})}\circ \tilde\tau_\ell\circ a_{1\over \lambda})\cdot \big((\tilde\tau_\ell)^*((\pi^*\theta_\ell) (\Phi)-((\pi^*\theta_\ell) (\Phi)\big)\rangle .
$$
Observe that if  $ (\ind_{\Tube(B,{r_0\over |\lambda|})}\circ \tilde\tau_\ell\circ a_{1\over \lambda})(y)=1$ then $y\in U_0.$ 
Moreover, the  $\Cc^0(U_0)$-norm of the  test form  $(\tilde\tau_\ell)^*(\pi^*\theta_\ell) (\Phi)-(\pi^*\theta_\ell) (\Phi)$
is  $\leq  c_0$ independent of $\Phi$ as  above. Therefore,
using   \eqref{e:T-hash-lambda-finite}  we see easily that the modulus of   the last expression is $ \lesssim c_0$ independent of $\Phi$ as  above  and of $\lambda\in\C^*$ with $|\lambda|\leq 1.$
Hence,  \eqref{e:T_lambda-bounded-mass} follows.
\endproof


To prove assertion (2) of    Theorem  \ref{T:Tangent-closed}, the  following  result is needed.
\begin{lemma}\label{L:negligibility-tangent-currents}
Let $\Phi$ be a $\Cc^1$-smooth test form with $\|\Phi\|_{\Cc^1(U^0)}\leq 1$ as in the proof of assertion  (1) of    Theorem  \ref{T:Tangent-closed}. 
 For  every $1\leq\ell\leq  \ell_0$ and for every $\lambda\in\C$ with $|\lambda|\geq 1,$ there are  $N$ continuous  functions $\psi_{j,\lambda}$ defined
 on $\Tube(B',{r_0\over |\lambda|})$ 
 and  $N$ continuous test forms $\Psi_{j,\lambda}$  defined on $U^0$  
 such that $\|\psi_{j,\lambda}\|_{\Cc^0(\Tube(B',{r_0\over |\lambda|}))}\leq c_0|\lambda|^{-1}$ and  $\|\Psi_{j,\lambda}\|_{\Cc^0(U_0)}\leq c_0$  and  that
$$ (\tilde\tau_\ell)^*((\pi^*\theta_\ell) (a_\lambda)^*\Phi)-((\pi^*\theta_\ell) (a_\lambda)^*\Phi)=\sum_{j=1}^N\psi_{j,\lambda} \cdot ( (a_\lambda)^*\Psi_{j,\lambda}).$$
Here, $c_0>0$ and $N\in\N$  are  constants independent of $\Phi$ and $\lambda.$
\end{lemma}
\proof
For simplicity write  write  $\tilde\tau$ (resp. $\Phi$) instead of $\tilde\tau_\ell$ (resp.  $(\pi^*\theta_\ell)\Phi$).
We need to show  that  $\tilde\tau^*(a_\lambda^*)(\Phi)-(a_\lambda^*)(\Phi)$ has the  desired form.
In order  to obtain this   result, we  study the action of $ \tilde\tau^*(a_\lambda^*)$ and  that of $(a_\lambda^*)$
on $\Cc^1$-smooth functions and on linear $1$-forms.  The form $\Phi$  is  built using  these functions and $1$-forms.

Let $f$ be a $\Cc^1$-smooth function with compact  support in $U'.$  For  $(z,w)\in \Tube(B',{r_0\over |\lambda|}),$  write $(z',w'):=\tilde\tau_\ell(z,w).$
Then we have that
$$
\tilde\tau^*_\ell(a_\lambda^*)f(z,w))-(a_\lambda^*)f(z,w)=f(\lambda  z',w') -f(\lambda z,w).
$$
Since $\tilde\tau$ is   admissible, it follows from Definition \ref{D:admissible-maps} that  
$$\|(\lambda  z',w') -(\lambda z,w)\|=|\lambda|\|z'-z\|+\|w'-w\|=|\lambda|O(\|z\|^2)+ O(\|z\|) =O(\lambda^{-1}).$$
The $\Cc^1$-smoothness of $f$ implies that  the above expression is uniformly bounded by a constant times $|\lambda|^{-1}.$

Consider now  the  forms $
\tilde\tau^*_\ell(a_\lambda^*)dw_q-(a_\lambda^*)dw_q$ and  $
\tilde\tau^*_\ell(a_\lambda^*)d\bar w_q-(a_\lambda^*)d\bar w_q$ for $1\leq q\leq l.$ We only discuss  the first form; the other  form can  be treated similarly.
Since  $(a_\lambda^*)dw_q=dw_q$ and  $\tilde\tau_\ell$ is   admissible, it follows from Definition \ref{D:admissible-maps} that  
$$
\tilde\tau^*_\ell(a_\lambda^*)dw_q-(a_\lambda^*)dw_q=\sum_{q'=1}^l \big[O(\|z\|)(a_\lambda^*) dw_{q'} +O(\|z\|)(a_\lambda^*) d\bar w_{q'}\big ]
+\sum_{p=1}^{k-l} \big[ O(\lambda^{-1})(a_\lambda^*) dz_p +O(\lambda^{-1})(a_\lambda^*) d\bar z_{p}\big].
$$
The LHS has the desired form  because $\|z\|=O(\lambda^{-1}).$

Consider now  the  forms $
\tilde\tau^*(a_\lambda^*)dz_p-(a_\lambda^*)dz_p$ and  $
\tilde\tau^*(a_\lambda^*)d\bar z_p-(a_\lambda^*)d\bar z_p$ for $1\leq p\leq k- l.$ We only discuss  the first form; the other  form can  be treated similarly.
Since  $(a_\lambda^*)dz_p=\lambda dz_p$ and  $\tilde\tau_\ell$ is   admissible, it follows from Definition \ref{D:admissible-maps} that  
$$
\tilde\tau^*(a_\lambda^*)dz_p-(a_\lambda^*)dz_p=\sum_{q=1}^l  O(\lambda \|z\|^2)(a_\lambda^*) dw_{q} +O(\lambda\|z\|^2)(a_\lambda^*) d\bar w_{q}
+\sum_{p'=1}^{k-l} \big( O(\|z\|)(a_\lambda^*) dz_{p'} +O(\|z\|)(a_\lambda^*) d\bar z_{p'}\big).
$$
The LHS has the desired form  because $O(\lambda \|z\|^2)=O(\lambda^{-1})$ and $O(\|z\|)=O(\lambda^{-1}).$
The proof is thereby completed.
\endproof

  \proof[Proof of assertion (2) of Theorem  \ref{T:Tangent-closed}]
 We  fix  an $\ell\in L$ (resp. an $\ell'\in L'$) and  write $\tau,$ $\tau'$) (resp.  $T_\lambda,$  $T'_\lambda$)  instead  of  $\tau_\ell,$ $\tau'_{\ell'}$)  (resp.  $T_{\lambda,\ell},$     $T'_{\lambda,\ell'}$).
 We also fix a connected component $U^0$ of $U_\ell\cap U'_{\ell'}$ and  a compact subset $K\Subset  U^0.$ 
 To  prove assertion (2), we need  to show that  for every continuous  test form $\Phi$ supported in $K,$
\begin{equation}\label{e:T_lambda-cv}\lim_{\lambda\to\infty}\langle  T_\lambda-T'_\lambda ,\Phi\rangle=0.
\end{equation}
Since  we know  by  assertion (1) that  the masses of $T_\lambda$ and $T'_\lambda$  are  unifomly bounded on compact subsets of $\pi^{-1}(U^0)$
independently of $\lambda\in\C^*,$  we may assume that $\Phi$ is of class $\Cc^1$ with support in $\Tube(B,r_0)$ for some $r_0>0.$

Let  $(\tilde\tau'_\ell)_{1\leq \ell\leq \ell_0}$ be the family which is obtained from  $\tau'$ in exactly the same way as 
 $(\tilde\tau_\ell)_{1\leq \ell\leq \ell_0}$ associated to $\tau.$  Using \eqref{e:tangent-current-closed} for  $\tau$ and $\tau',$
 we get that \begin{eqnarray*}&&\langle  T_\lambda-T'_\lambda ,\Phi\rangle=\langle \tau_*T,\ind_{\Tube(B,{r_0\over |\lambda|})}(a_\lambda)^*\Phi\rangle -
 \langle \tau'_*T,\ind_{\Tube(B,{r_0\over |\lambda|})}(a_\lambda)^*\Phi\rangle\\
 &=&\sum_{\ell=1}^{\ell_0}\langle(\tau_\ell)_* T,  (\ind_{\Tube(B,{r_0\over |\lambda|})}\circ \tilde\tau_\ell)\cdot \big((\tilde\tau_\ell)^*((\pi^*\theta_\ell) (a_\lambda)^*\Phi)-((\pi^*\theta_\ell) (a_\lambda)^*\Phi)\big)\rangle\\
 &-& \sum_{\ell=1}^{\ell_0}\langle(\tau'_\ell)_* T,  (\ind_{\Tube(B,{r_0\over |\lambda|})}\circ \tilde\tau'_\ell)\cdot \big((\tilde\tau'_\ell)^*((\pi^*\theta_\ell) (a_\lambda)^*\Phi)-((\pi^*\theta_\ell) (a_\lambda)^*\Phi)\big)\rangle
 .
\end{eqnarray*}
The assertion will follow  if one can show that both terms on the RHS tends to $0$ as $\lambda$ tends to infinity.
We will prove this  for the  first term since  the proof for the second one is  similar.
Applying Lemma \ref{L:negligibility-tangent-currents}, the first term  is  equal to
$$
\sum_{\ell=1}^{\ell_0}\sum_{j=1}^N\langle(\tau_\ell)_* T,  (\ind_{\Tube(B,{r_0\over |\lambda|})}\circ \tilde\tau_\ell)\cdot\psi_{j,\lambda}\cdot  (a_\lambda^*) \Psi_{j,\lambda}\rangle,
$$
where, for $1\leq j\leq N,$  $\Psi_{j,\lambda}$ is a continuous test form with $\|\Psi_{j,\lambda}\|_{\Cc^0(U^0)}\leq c_0$ and $\psi_{j,\lambda}$ is a  continuous  function  defined
 with
  $\|\psi_{j,\lambda}\|_{\Cc^0(\Tube(B',{r_0\over |\lambda|}))}\leq c_0|\lambda|^{-1}.$
Using  this  and assertion (1) and inequality \eqref{e:T-hash-lambda-finite}, we see that   the above expression  tends to  $0$ as $\lambda$ tends to infinity.
This completes the proof of assertion (2).
\endproof


  \proof[Proof of assertion (3) of Theorem  \ref{T:Tangent-closed}]
  Pick arbitrary    $\ell\in L$  and $\ell'\in L'$ such that $B^0:=B\cap U_\ell\cap U'_{\ell'}\not=\varnothing.$
  We  only need to show that  $T_\infty=\lim_{\lambda\to\infty} T'_{\lambda,\ell'}$ on  $\pi^{-1}(B^0)\subset \E.$
  On the one hand,  we know  by  the hypothesis that $T_\infty=\lim_{\lambda\to\infty} T_{\lambda,\ell}$ on  $\pi^{-1}(B^0)\subset \E.$
  On the other hand, by assertion (2), 
   the  family of currents  $T_{\lambda,\ell} -T'_{\lambda,\ell'}$ converge  weakly to  $0$ on  $\pi^{-1}(B^0)\subset \E$
 as  $\lambda$ tends to infinity.  
 Hence, the result follows. \endproof
  
  
  \begin{remark} \rm \label{R:Tangent-closed}  Theorem \ref{T:Tangent-closed} still holds
  if 
    $\omega$ is a Hermitian metric on $V$  such that $ \ddc\omega^j=0$ on $V$ for all  $1\leq j\leq   \upm-1.$ However, we need a stronger assumption on $T,$ namely,
   the $(p,p)$-current   $T$  is  positive closed  and $T=T^+-T^-$ on an open neighborhood of $\overline B$ in $X$  with
$T^\pm$  in the class $\CL_p^{2,2} (B).$  To see this,  we  apply  Proposition  \ref{P:Lc-finite-closed-non-Kaehler} instead of  
 Proposition \ref{P:Lc-finite-closed} in order to obtain  inequality 
 \eqref{e:Mc-finite-local}.  The rest of the proof follows along  the same lines as those given in the proof
 of  Theorem \ref{T:Tangent-closed}.
  \end{remark}
  
\subsection{Positive  pluriharmonic  currents and positive plurisubharmonic  currents}
Now  we are ready to state  and  prove   the  existence  of tangent currents for positive plurisubharmonic  currents.

 \begin{theorem}\label{T:Tangent-psh}
We  keep the     Standing Hypothesis. Suppose that $\omega$ is  K\"ahler  and the $(p,p)$-current   $T$  is  positive plurisubharmonic    and $T=T^+-T^-$ on an open neighborhood of $\overline B$ in $X$  with
$T^\pm$  in the class $\SH_p^{3,3} (B).$  Suppose in addition   that
   there is  at least  one strongly  admissible map
along $B.$
  Then the same  assertions (1)--(3) as  those of Theorem  \ref{T:Tangent-closed} hold. 
 \end{theorem}
 \proof
 
By Proposition \ref{P:Lc-finite-psh} and  Lemma \ref{L:Mc_j}, we have
 for $\lowm\leq j\leq \upm,$  
 \begin{equation} \label{e:Mc-finite-local-psh}   \sup\limits_{r\in (0,\bfr]} {1\over r^{2(k-p-j)}}\int_{\|z\|<r,\ w\in U^0_w} T^\hash \wedge \omega_w^j \wedge \omega_z^{k-p-j}<c_{12},
 \end{equation}
 Using this  instead of  \eqref{e:Mc-finite-local}, we argue as in the proof of Theorem \ref{T:Tangent-closed}.
 \endproof
  Similarly,  we also obtain  the  existence  of tangent currents for positive pluriharmonic  currents.
 \begin{theorem}\label{T:Tangent-ph}
We  keep the     Standing Hypothesis. Suppose that $\omega$ is  K\"ahler  and the $(p,p)$-current   $T$  is  positive plurisubharmonic    and $T=T^+-T^-$ on an open neighborhood of $\overline B$ in $X$  with
$T^\pm$  in the class $\PH_p^{2,2} (B).$  Suppose in addition   that
   there is  at least  one strongly  admissible map
along $B.$
  Then the same  assertions (1)--(3) as  those of Theorem  \ref{T:Tangent-closed} hold. 
 \end{theorem}
 \proof  As in Definition \ref{D:sup}, we have  the following
  \begin{definition}\label{D:sup_PH}  \rm
Fix an open neighborhood $\bfU$ of $\overline B$ and an open neighborhood $\bfW$ of $\partial B$ in $X$ with $\bfW\subset \bfU.$
Let $\widetilde\PH^{2,2}_p(\bfU,\bfW)$ be the  set of all $T\in \PH^{2,2}_p(\bfU,\bfW)$  whose  a sequence of approximating  forms $(T_n)_{n=1}^\infty$
satisfies the following   condition:
 \begin{equation}\label{e:unit-PH-2,2} \|T_n\|_{\bfU}\leq  1\qquad\text{and}\qquad  \| T_n\|_{\Cc^2(\bfW)}\leq 1.\end{equation}
 \end{definition}
 By  Theorem \ref{T:Lelong-psh}   (6), we can obtain the following   
 result which is  the  analogue of Proposition \ref{P:Lc-finite-psh}   for positive pluriharmonic  currents.
 \begin{proposition}\label{P:Lc-finite-ph}
  There is a   constant $c_{11}>0$ such that for  
  every positive pluriharmonic  current $T$  belonging to the class $\widetilde\PH^{2,2}_p(\bfU,\bfW),$ 
   we have
   $\Mc_j(T,r)<c_{11}$ for $0\leq j\leq \upm$ and $0 <r \leq \bfr.$
 \end{proposition}
By Proposition \ref{P:Lc-finite-ph} and  Lemma \ref{L:Mc_j}, we also  obtain
 inequality \eqref{e:Mc-finite-local-psh} for $\lowm\leq j\leq \upm.$   
 Using this  instead of  \eqref{e:Mc-finite-local}, we argue as in the proof of Theorem \ref{T:Tangent-closed}.

 \endproof
\section{$V$-conic and pluriharmonicity  of tangent currents}
\label{S:V-conic-and-plurihar}

Recall the  Standing Hypothesis from Subsection \ref{SS:Global-setting}.
The  main purpose of this section is  to  establish some basic properties of tangent  currents  in   three  families of currents: the positive closed currents, the positive pluriharmonic  currents 
and the positive plurisubharmonic  currents.
\subsection{Positive  closed  currents}

\begin{theorem}\label{T:conic-closed}
   We  keep the     Standing Hypothesis. Suppose that $\omega$ is  a K\"ahler on $V.$ 
   Suppose in addition that    the current   $T$  is  positive closed  and $T=T^+-T^-$ on an open neighborhood  of $\overline B$ in $X$
with $T^\pm$ in the class $\CL_p^{1,1}( B).$
   Let $T_\infty$ be a tangent current to $T$ along $B$  given by  Theorem  \ref{T:Tangent-closed} (3).
  Then  $T_\infty$ is a $V$-conic  positive  closed  $(p,p)$-current on $\pi^{-1}(B)\subset \E.$
 \end{theorem}
 
 
\proof[Proof that  $T_\infty$ is positive closed]
Consider the covering family of  holomorphic admissible maps $\Uc=(\bfU_\ell,\tau_\ell)_{1\leq \ell\leq \ell_0}$ for $B$  introduced in  Subsection 
\ref{SS:Ex-Stand-Hyp}.  By  Theorem  \ref{T:Tangent-closed} (3),
$T_\infty$ is   the tangent current  to $T$ along $B$ associated, by Definition \ref{D:tangent-currents},   to a sequence
 $(\lambda_n ) \subset \C^*$ converging to $\infty$  and to the family $\Uc.$ More precisely, for $1\leq \ell\leq \ell_0,$ we have 
 $$
 T_\infty=\lim_{n\to\infty} T_{\lambda_n,\tau_\ell}\qquad\text{on}\qquad  \pi^{-1}(B\cap \bfU_\ell).  
 $$
 Since  $T$ is  positive closed and  $\tau_\ell$ is  holomorphic, we infer from the formula
 $T_{\lambda_n,\tau_\ell}=(A_{\lambda_n})_*((\tau_\ell)_*T)$  
that $T_{\lambda_n,\tau_\ell}$ is  positive closed.  Hence, the above limit implies that  $T_\infty$ is  also positive closed.
 \endproof


\proof[Proof that  $T_\infty$ is $V$-conic]
Let  $0<r_1<r_2\leq \bfr$  and $\lambda\geq 1$ and $\lowm\leq j\leq \upm.$
Applying  Theorem  \ref{T:Lelong-Jensen-closed}  to $(A_{\lambda_n})_*(\tau_*T)$ yields that  
\begin{multline*}
 \nu_j((A_{\lambda_n})_*(\tau_*T),B,r_2,\id)-\nu_j((A_{\lambda_n})_*(\tau_*T),B,r_1,\id)\\
 =\kappa_j((A_{\lambda_n})_*(\tau_*T),B,r_1,r_2,\id)+
 \Vc((A_{\lambda_n})_*(\tau_*T),r_1,r_2).
\end{multline*}
 By Proposition \ref{P:scale}, the LHS is   equal to
 \begin{equation*}
 \nu_j(T,B,{r_2\over |\lambda_n|},\tau)-\nu_j(T,B,{r_1\over |\lambda_n|},\tau).
 \end{equation*}
This quantity converges, by   Theorem  \ref{T:Lelong-closed} applied to $T,$
to  $\nu_j(T_\infty,B,\id)- \nu_j(T_\infty,B,\id)=0.$
 On the  other hand, by  Theorem \ref{T:vertical-boundary-closed},  $\Vc((A_{\lambda_n})_*(\tau_*T),r_1,r_2)=O(\lambda^{-1}).$
 Moreover, as $n$ tends to infinity, $ (A_{\lambda_n})_*(\tau_*T)$ tends to $T_\infty,$ we have that
 \begin{equation*}
  \kappa_j((A_{\lambda_n})_*(\tau_*T),B,r_1,r_2,\id)\to  \kappa_j(T_\infty,B,r_1,r_2,\id).
 \end{equation*}
Consequently, we infer that   $\kappa_j(T_\infty,B,r_1,r_2,\id)=0$ for $\lowm\leq j\leq \upm.$ So 
   \begin{equation}\label{e:kappa-T_infty-equal-zero}
    \int_{\Tube (B,r_1,r_2)} T_\infty\wedge \alpha^{k-p-j} \wedge \pi^*\omega^j=0 \quad\text{for}\quad \lowm\leq j\leq \upm\quad\text{and}\quad 0<r_1<r_2\leq \bfr.
   \end{equation}
  Since $T\wedge \pi^*\omega^\upm$  is  of full bidegree $(l,l)$ in $\{dw,d\bar w\}$ by the fact in Corollary \ref{C:Lelong-Jensen}, we deduce from \eqref{e:kappa-T_infty-equal-zero} that for all $0<r_1<r_2\leq \bfr$ and for $\lowm\leq j \leq k:$
   \begin{equation*}
    \int_{\Tube(B,r_1,r_2)} T_\infty\wedge (\alpha+c_1 \pi^*\omega)^{k-p-j} \wedge \pi^*\omega^j=0.
   \end{equation*}
   Recall from  \eqref{e:hat-alpha'} that  $\hat\alpha'=\alpha+c_1\pi^*\omega.$
By  \eqref{e:hat-alpha'-vs-alpha_ver},  $
\hat\alpha'\geq  c^{-1}_1\alpha_\ver\geq 0.$  Moreover, $T_\infty$ is   a positive current.
So by letting $r_1\to 0$ and $r_2\to\bfr,$  we  get that  for $\lowm\leq j\leq k,$
\begin{equation*} T_\infty\wedge (\pi^*\omega)^j\wedge \alpha_\ver^{k-p-j}=0\qquad\text{on}\qquad \Tube(B,\bfr)\setminus B.
\end{equation*}
Since  $\alpha_\vert^{k-l}=0,$ it follows that    $\alpha_\ver^{k-p-j}=0$ for $0\leq j<\lowm.$ This,  combined with the last equality, implies that 
 for $0\leq j\leq k,$
\begin{equation}\label{e:T_infty-proj_iden} T_\infty\wedge (\pi^*\omega)^j\wedge \alpha_\ver^{k-p-j}=0\qquad\text{on}\qquad \Tube(B,\bfr)\setminus B.
\end{equation}
We place ourselves on  an open set of $\C^{k-l}$ defined by $z_{k-l}\not=0.$
We   may assume without loss of generality that
$2|z_{k-l}| > \max\limits_{1\leq j\leq k-l}|z_j|$
and use the projective coordinates  introduced  in \eqref{e:homogeneous-coordinates}:
$$
\zeta_1:={z_1\over z_{k-l}},\ldots, \zeta_{k-l-1}:={z_{k-l-1}\over z_{k-l}},\quad \zeta_{k-l}=z_{k-l}.
$$
In the coordinates  $\zeta=(\zeta_1,\ldots,\zeta_{k-l})=(\zeta',\zeta_{k-l}),$ the form $\omega_\FS([z])$  can be  rewritten as  $$\ddc \log{ (1+|\zeta_1|^2+\cdots+|\zeta_{k-l-1}|^2)}.$$
We recall  from \eqref{e:omega_FS-vs-omega'} that
\begin{equation}\label{e:omega_FS-vs-omega'-bis}
 \omega'(\zeta')\approx \alpha_\ver,
\end{equation}
since  both  of them  are  equivalent to $\omega_\FS([z]).$
Here $\omega'(\zeta'):=\ddc  (|\zeta_1|^2+\cdots+|\zeta_{k-l-1}|^2).$
Let $\Theta_{I,J;K,L}$ be  the coefficients of a current $\Theta$  in the coordinates $(\zeta,w)$ 
 according to  Definition
  \ref{D:negligible}.
We have that for $0\leq j\leq m,$
$$
\Theta\wedge\omega_w^j\wedge   {\omega}'(\zeta')^{k-p-j}=\sum_{I\ni   k-l} \Theta_{I,I;K,K}  \Leb(\zeta,w),
$$
where
\begin{equation*}
  \Leb(\zeta,w):=    ( i d\zeta_1\wedge d\bar \zeta_1)\wedge \ldots\wedge    ( i d\zeta_{k-l}\wedge d\bar \zeta_{k-l}) 
\wedge ( i dw_1\wedge d\bar w_1)\wedge \ldots\wedge    ( i dw_{l}\wedge d\bar w_{l}) .\end{equation*}

Now  set $\Theta:=T_\infty.$ Combining  equality  \eqref{e:T_infty-proj_iden} and estimate \eqref{e:omega_FS-vs-omega'-bis} gives that
$\Theta_{I,I;K,K}=0$ for $I,K$  with $ k-l\in I.$ Using this, Proposition \ref{P:Demailly} applied with  $\lambda_{k-l}\geq 0$ arbitrary  and $\lambda_j=1$ for $j\not=k-l,$ yields that
$\Theta_{I,J;K,K}=0$ if  
both $I$ and $J$ contains the element $k-l.$
The same argument also shows that
$$
\lambda_{k-l}|\Theta_{I,J;K,K}|\leq 2^{k-p} \sum_{M:\ k-l\not\in M} \Theta_{M,M}
$$ 
if either  $I$ or $J$   contains $k-l.$ Letting  $\lambda_{k-l}\to\infty$ we infer that $\Theta_{I,J;K,K}=0$  
in this case.

The same  argument also   shows that $\Theta_{I,J;K,L}=0$
if  either  $I$ or $J$   contains $k-l.$
Since $T_\infty$ is  closed, we infer from  the last  equality that
$${\partial \Theta_{I,J;K,L}\over \partial \zeta_{k-l}}={\partial \Theta_{I,J;K,L}\over \partial \overline\zeta_{k-l}}=0\quad\text{for all}\quad I,J,K,L.$$
So $T_\infty$ depends only on the variables $\zeta'$ and $w.$
As  the projection $\Pi:\  \C^{k-l}\setminus \{0\})\times  \C^l\to  \P^{k-l-1}\times \C^l$ may be rewritten as
$(\zeta,w)\mapsto (\zeta',w),$ we see that $T_\infty|_{ \C^{k-l}\setminus \{0\})\times  \C^l}$ is the preimage by $\Pi$   of a positive closed current on
 $\T_\infty$  living on $\P^{k-l-1}\times \C^l.$ Hence, $T_\infty$ is $V$-conic.
\endproof

\begin{remark} \rm \label{R:conic-closed}  Theorem \ref{T:conic-closed} still holds
  if 
    $\omega$ is a Hermitian metric on $V$  such that $ \ddc\omega^j=0$ on $V$ for all  $1\leq j\leq   \upm-1.$ However, we need a stronger assumption on $T,$ namely,
   the $(p,p)$-current   $T$  is  positive closed  and $T=T^+-T^-$ on an open neighborhood of $\overline B$ in $X$  with
$T^\pm$  in the class $\CL_p^{2,2} (B).$  To see this,  we  apply  Theorem  \ref{T:Lelong-closed} instead of  
Theorem \ref{T:Lelong-closed-Kaehler}.  The rest of the proof follows along  the same lines as those given in the proof
 of  Theorem \ref{T:conic-closed}.
  \end{remark}

\subsection{Positive  pluriharmonic currents and positive plurisubharmonic  currents}
Now  we are ready to state  and  prove   some  basic properties of tangent currents for positive plurisubharmonic  currents.

 \begin{theorem}\label{T:conic-psh}
We  keep the     Standing Hypothesis. Suppose that $\omega$ is  K\"ahler and  the $(p,p)$-current   $T$  is  positive plurisubharmonic    and $T=T^+-T^-$ on an open neighborhood of $\overline B$ in $X$  with
$T^\pm$  in the class $\SH_p^{3,3} (B).$
  Suppose in addition    that there is  at least  one strongly  admissible map
along $B.$  Let $T_\infty$ be a tangent current to $T$ along $B$  given by  Theorem \ref{T:Lelong-psh}.
  Then the following    assertions  hold:
  \begin{enumerate}
   \item  $T_\infty$ is a   positive  plurisubharmonic  $(p,p)$-current on $\pi^{-1}(B)\subset \E.$
   
   \item  $T_\infty$ is  partially pluriharmonic  in the  sense that the  current $T_\infty\wedge \pi^*(\omega^\lowm)$ is  pluriharmonic.
   
   \item   $T_\infty$ is  partially  $V$-conic in the  sense that the  current $T_\infty\wedge \pi^*(\omega^\lowm)$ is  $V$-conic.
 
  \end{enumerate}

 \end{theorem}

\proof[Proof  of assertion (1)]
We  keep the notation introduced  in the above proof of Theorem \ref{T:conic-closed}. 
Since  $T$ is  positive plurisubharmonic  and  $\tau_\ell$ is  holomorphic, we infer from the formula
 $T_{\lambda_n,\tau_\ell}=(A_{\lambda_n})_*((\tau_\ell)_*T)$  
that $T_{\lambda_n,\tau_\ell}$ is  also positive  plurisubharmonic.  Hence, the  limit $T_\infty=\lim_{n\to\infty}  T_{\lambda_n,\tau_\ell}$ on $\pi^{-1}(B\cap \bfU_\ell)$  implies that  
$T_\infty$ is  also positive  plurisubharmonic.
\endproof


\proof[Proof of assertion (2)]
 By  Theorem  \ref{T:Tangent-closed} (3),
$T_\infty$ is   the tangent current  to $T$ along $B$ associated, by Definition \ref{D:tangent-currents},   to a sequence
 $(\lambda_n ) \subset \C^*$ converging to $\infty$  and to the family $\Uc.$
Fix $r_1,\ r_2\in  (0,\bfr)$ with $r_1<r_2.$  Let $\lambda\in\R$ with $\lambda \geq 1.$  

For every $j$ with  $\lowm\leq j\leq \upm,$  applying Theorem  \ref{T:Lelong-Jensen}  to  $(A_{\lambda_n})_*(\tau_*T)\wedge \pi^*(\omega^j)$ yields that
\begin{multline*}
  \nu_j(T,B,{r_2\over |\lambda_n|},\tau)- \nu_j(T,B,{r_1\over |\lambda_n|},\tau)=  \Vc\big((A_{\lambda_n})_*(\tau_*T)\wedge \pi^*(\omega^j),r_1,r_2\big)\\ +\int_{\Tube(B,r_1,r_2)}(A_{\lambda_n})_*(\tau_*T)\wedge \pi^*(\omega^{j})\wedge\alpha^{k-p-j}\\ 
   +  \int_{r_1}^{r_2} \big( {1\over t^{2(k-p-j)}}-{1\over r_2^{2(k-p-j)}}  \big)2tdt\int_{\Tube(B,t)}  \ddc (A_{\lambda_n})_*(\tau_*T)\wedge \pi^*(\omega^j)\wedge \beta^{k-p-j-1} \\
  +  \big( {1\over r_1^{2(k-p-j)}}-{1\over r_2^{2(k-p-j)}}  \big) \int_{0}^{r_1}2tdt\int_{z\in \Tube(B,t)} \ddc (A_{\lambda_n})_*(\tau_*T)\wedge \pi^*(\omega^j)\wedge \beta^{k-p-j-1}.
\end{multline*}
We let $n$ tend to infinity.  The LHS  tends to $0$   since by
Theorem \ref{T:Lelong-psh} (1), $\lim_{n\to\infty}\nu_j(T,B,{r\over |\lambda_n|},\tau)=\nu_j(T,B,\tau)\in\R$ for $0<r\leq\bfr.$
By  Theorem \ref{T:vertical-boundary-terms},  $\Vc\big((A_{\lambda_n})_*(\tau_*T)\wedge \pi^*(\omega^j),r_1,r_2\big)\to 0$ as $n\to\infty.$
Therefore, we obtain that 
\begin{equation}\label{e:sum_integrals_equal_zero}\begin{split}
  0&=   \int_{\Tube(B,r_1,r_2)}T_\infty\wedge \pi^*(\omega^{j})\wedge\alpha^{k-p-j}\\
  & +  \int_{r_1}^{r_2} \big( {1\over t^{2(k-p-j)}}-{1\over r_2^{2(k-p-j)}}  \big)2tdt\int_{\Tube(B,t)}  \ddc T_\infty\wedge \pi^*(\omega^j)\wedge \beta^{k-p-j-1} \\
  &+  \big( {1\over r_1^{2(k-p-j)}}-{1\over r_2^{2(k-p-j)}}  \big) \int_{0}^{r_1}2tdt\int_{z\in \Tube(B,t)} \ddc T_\infty\wedge \pi^*(\omega^j)\wedge \beta^{k-p-j-1}.
  \end{split}
\end{equation}
Next,  we  argue as in the proof of assertion (1) of Theorem \ref{T:top-Lelong-psh}.
Consider  a small neighborhood $V(x_0)$ of  an arbitrary  point $x_0\in  \Tube(B, r_0),$  where in a local chart $V(x_0)\simeq \D^l$ and  $\E|_{V(x_0)}\simeq \C^{k-l}\times \D^l.$
For $x\in \E|_{V(x_0)},$ write $x=(z,w).$  Since   $\upm=\min(l,k-p)$ and $T_\infty$ is  of bidegree $(p,p)$   we  see that  $T_\infty\wedge \pi^*\omega^\upm$ is  of full bidegree $(l,l)$  in $dw,$ $d\bar w.$
Consequently, we infer from \eqref{e:hat-alpha'} that
\begin{eqnarray*}
T_\infty\wedge \pi^*(\omega^m)\wedge\alpha^{k-p-\upm}
&=&T_\infty\wedge \pi^*(\omega^\upm)\wedge(\hat\alpha')^{k-p-\upm},\\
\ddc T_\infty\wedge \pi^*\omega^\upm\wedge\beta^{k-p-\upm}&=&\ddc T_\infty\wedge \pi^*\omega^\upm\wedge\hat\beta^{k-p-\upm}.
\end{eqnarray*}
This,  combined  with \eqref{e:sum_integrals_equal_zero} for $j:=\upm,$  implies  that 
\begin{eqnarray*}
  0&&= \int_{\Tube(B,r_1,r_2)}T_\infty\wedge \pi^*(\omega^\upm)\wedge(\hat\alpha')^{k-p-\upm}\\
  &+&  \int_{r_1}^{r_2} \big( {1\over t^{2(k-p-\upm)}}-{1\over r_2^{2(k-p-\upm)}}  \big)2tdt\int_{\Tube(B,t)}  \ddc T_\infty\wedge (\pi^*\omega^\upm)\wedge \hat\beta^{(k-p-\upm)-1} \\
 &+&  \big( {1\over r_1^{2(k-p-\upm)}}-{1\over r_2^{2(k-p-\upm)}}  \big) \int_{0}^{r_1}2tdt\int_{z\in \Tube(B,t)} 
 \ddc T_\infty\wedge (\pi^*\omega^\upm)\wedge \hat\beta^{(k-p-\upm)-1}. 
\end{eqnarray*}
Since we have just shown that $T_\infty$ is  positive  plurisubharmonic,  both $T_\infty$ and $\ddc T_\infty$  are   positive  currents.
Moreover,  $\omega,$ $\hat\alpha',$  $\hat\beta$ are positive forms.  Consequently,  all integrals of  the RHS of the last line  are $\geq 0.$
On the ther hand, their sum  is  equal to $0.$  So  all integrals  are  $0,$ that is,
\begin{equation*}
  \int_{z\in \Tube(B,r_2)} \ddc T_\infty\wedge \pi^*(\omega^{\upm})\wedge \hat\beta^{k-p-\upm-1}=0\quad\text{and}\quad  \int_{\Tube(B,r_1,r_2)}T_\infty\wedge \pi^*(\omega^\upm)\wedge\alpha^{k-p-\upm}=0.
\end{equation*}
Note that  $\hat\beta$  and $\pi^*\omega$ are  smooth  strictly positive $(1,1)$ forms on $\Tube(B,\bfr),$ and that for every   smooth   positive $(1,1)$ form $H$ on $\Tube(B,\bfr),$ we can find a constant $c>0$ such that $H\leq c (\hat\beta+\pi^*\omega)$ on $\Tube(B,\bfr).$
Since $0<r_1<r_2\leq\bfr$ are arbitrarily chosen, 
 we infer that  the following  equality holds for all $j$ with $\upm\leq j\leq k$:
 \begin{equation}\label{e:fact-partial-plurihar}
  \ddc T_\infty\wedge\pi^*(\omega^j)=0\quad\text{and}\quad T_\infty\wedge \pi^*(\omega^j)\wedge\alpha^{k-p-j}=0 \quad\text{on}\quad\Tube(B,\bfr)\quad\text{for}\quad \lowm\leq j\leq k. 
 \end{equation}
Suppose  that \eqref{e:fact-partial-plurihar} holds for all $j$  with $j_0<j\leq\upm,$ where $j_0$ is a given integer
with $\lowm\leq j_0<\upm.$  We need to  prove \eqref{e:fact-partial-plurihar}  for  $j=j_0.$ 

 Using  \eqref{e:fact-partial-plurihar} for all $j$  with $j_0<j\leq k,$ we  infer from \eqref{e:hat-alpha'}  that
\begin{eqnarray*}
T_\infty\wedge \pi^*(\omega^{j_0})\wedge\alpha^{k-p-j_0}
&=&T_\infty\wedge \pi^*(\omega^{j_0})\wedge(\hat\alpha')^{k-p-j_0},\\
\ddc T_\infty\wedge \pi^*\omega^{j_0}\wedge\beta^{k-p-j_0}&=&\ddc T_\infty\wedge \pi^*\omega^{j_0}\wedge\hat\beta^{k-p-j_0}.
\end{eqnarray*}
This,  combined  with \eqref{e:sum_integrals_equal_zero} for $j:=j_0,$  implies  that 
\begin{eqnarray*}
  0&&= \int_{\Tube(B,r_1,r_2)}T_\infty\wedge \pi^*(\omega^{j_0})\wedge(\hat\alpha')^{k-p-{j_0}}\\
  &+&  \int_{r_1}^{r_2} \big( {1\over t^{2(k-p-j_0)}}-{1\over r_2^{2(k-p-j_0)}}  \big)2tdt\int_{\Tube(B,t)}  \ddc T_\infty\wedge (\pi^*\omega^{j_0})\wedge \hat\beta^{(k-p-j_0)-1} \\
 &+&  \big( {1\over r_1^{2(k-p-j_0)}}-{1\over r_2^{2(k-p-j_0)}}  \big) \int_{0}^{r_1}2tdt\int_{z\in \Tube(B,t)} 
 \ddc T_\infty\wedge (\pi^*\omega^{j_0})\wedge \hat\beta^{(k-p-j_0)-1}. 
\end{eqnarray*}
We repeat the above argument   using that   both $T_\infty$ and $\ddc T_\infty$  are   positive  currents and that $\pi^*\omega,$  $\hat\alpha',$  $\hat\beta$ are positive forms.  Consequently,    all integrals on the RHS  are  $0.$
Therefore, \eqref{e:fact-partial-plurihar} holds for  $j=j_0.$
Hence,  the proof of \eqref{e:fact-partial-plurihar} is  completed.
In particular,  $\ddc T_\infty\wedge\pi^*(\omega^\lowm)=0$ on $\Tube(B,\bfr).$
Since we will prove shortly below that  $T_\infty\wedge \pi^*(\omega^\lowm)$ is $V$-conic,
it follows that   $\ddc T_\infty\wedge\pi^*(\omega^j)=0$ on $\pi^{-1}(B)\subset \E.$
\endproof


\proof[Proof of assertion (3)]   Recall from  \eqref{e:hat-alpha'} that  $\hat\alpha'=\alpha+c_1\pi^*\omega,$ and from
  \eqref{e:hat-alpha'-vs-alpha_ver} that  $
\hat\alpha'\geq  c^{-1}_1\alpha_\ver\geq 0.$  Moreover, $T_\infty$ is   a positive current.
 Therefore, we infer from the  second identity of \eqref{e:fact-partial-plurihar} that
  \begin{equation}\label{e:T_infty-proj_iden-bis} T_\infty\wedge (\pi^*\omega)^j\wedge \alpha_\ver^{k-p-j}=0\quad\text{on}\quad \Tube(B,\bfr)\quad\text{for}\quad \lowm\leq j\leq k.
\end{equation}
Consider  the  positive  pluriharmonic  current $\Theta:=T_\infty\wedge\pi^*(\omega^\lowm).$
We need  to show  that $\Theta$ is  $V$-conic.
Let $\Theta_{I,J;K,L}$ be  the coefficients of the current  in the coordinates $(\zeta,w)$ 
 according to  Definition
  \ref{D:negligible}.  
Using \eqref{e:T_infty-proj_iden-bis}  we  argue   as in the proof  that $T_\infty$ is  conic when $T$ is a positive closed  current. Therefore, we can show that $\Theta_{I,J;K,L}=0$
if  either  $I$ or $J$   contains $k-l.$
Since  $\Theta$ is  $\ddc$-closed, we infer that
$${\partial^2 \Theta_{I,J;K,L}\over \partial \zeta_{k-l}\partial \overline\zeta_{k-l}}=0\quad\text{for all}\quad I,J,K,L.$$
So for fixed $(\zeta',w),$ $\Theta_{I,J;K,L}(\zeta,w)$ are  harmonic functions of  $\zeta_{k-l}.$  

We  choose  a  basis of $\Lambda^{k-p-\lowm,k-p-\lowm}(\C^k)$ consisting of  strictly  positive constant forms $\gamma_1,\ldots,\gamma_N,$
where $N:={k\choose p+\lowm}^2.$
Since  $\Theta$ is  positive, we can  write  $\Theta\wedge \gamma_j=f_j(\zeta,w)\Leb(\zeta,w)$  for $1\leq j\leq N,$  where $f_j$ is a non-negative function.
The harmonicity of the functions  $\Theta_{I,J;K,L}$ with respect to $\zeta_{k-l}$ and the  constant forms $\gamma_1,\ldots,\gamma_N,$
imply that $f_j$ is also   harmonic  functions with respect to $\zeta_{k-l}$.
So  the $f_j$'s are positive harmonic  functions with respect to $\zeta_{k-l}$.  On the other hand, positive harmonic  functions on the complex line
are necessarily constant.  Therefore,  we infer that  the $f_j$'s  are functions  depending only on  $\zeta'$ and $w.$
So  $\Theta_{I,J;K,L}$
depends only on the variables $\zeta'$ and $w.$ 
As  the projection $\Pi:\  \C^{k-l}\setminus \{0\})\times  \C^l\to  \P^{k-l-1}\times \C^l$ may be rewritten as
$(\zeta,w)\mapsto (\zeta',w),$ we see that $\Theta|_{ \C^{k-l}\setminus \{0\})\times  \C^l}$ is the preimage by $\Pi$   of a positive pluriharmonic 
current $\Theta_\infty$ on
 $\P^{k-l-1}\times \C^l.$ Hence, $\Theta$ is $V$-conic. 
\endproof
 
 
 We end the section with   some  basic properties of tangent currents for positive pluriharmonic  currents.
 
 \begin{theorem}\label{T:conic-ph}
We  keep the     Standing Hypothesis. Suppose that $\omega$ is  K\"ahler and  the $(p,p)$-current   $T$  is  positive plurisubharmonic    and $T=T^+-T^-$ on an open neighborhood of $\overline B$ in $X$  with
$T^\pm$  in the class $\PH_p^{2,2} (B).$
  Suppose in addition    that there is  at least  one strongly  admissible map
along $B.$  Let $T_\infty$ be a tangent current to $T$ along $B$  given by  Theorem \ref{T:Lelong-psh}.
  Then   $T_\infty$ is also  $V$-conic positive pluriharmonic.

 \end{theorem}

\proof
We  keep the notation introduced  in the above proof of Theorem \ref{T:conic-closed}. 
Since  $T$ is   positive  pluriharmonic  and  $\tau_\ell$ is  holomorphic, we infer from the formula
 $T_{\lambda_n,\tau_\ell}=(A_{\lambda_n})_*((\tau_\ell)_*T)$  
that $T_{\lambda_n,\tau_\ell}$ is  also positive  pluriharmonic.  Hence, the  limit $T_\infty=\lim_{n\to\infty}  T_{\lambda_n,\tau_\ell}$ on $\pi^{-1}(B\cap \bfU_\ell)$  implies that  
$T_\infty$ is  also positive  pluriharmonic. 

It  remains  to show that $T_\infty$ is  $V$-conic.
   For simplicity  write  $\Theta:=T_\infty.$
So $\Theta$ is a  positive  pluriharmonic  current.
We need  to show  that $\Theta$ is  $V$-conic. Let $\Theta_{I,J;K,L}$ be  the coefficients of the current  in the coordinates $(\zeta,w)$ 
 according to  Definition
  \ref{D:negligible}. 
Arguing as  in  the proof of assertion (3) of Theorem \ref{T:conic-psh} and
using \eqref{e:T_infty-proj_iden-bis}, we can show that $\Theta_{I,J;K,L}=0$
if  either  $I$ or $J$   contains $k-l.$
Since  $\Theta$ is  $\ddc$-closed, we infer that
$${\partial^2 \Theta_{I,J;K,L}\over \partial \zeta_{k-l}\partial \overline\zeta_{k-l}}=0\quad\text{for all}\quad I,J,K,L.$$
So for fixed $(\zeta',w),$ $\Theta_{I,J;K,L}(\zeta,w)$ are  harmonic functions of  $\zeta_{k-l}.$

We  choose  a  basis of $\Lambda^{k-p,k-p}(\C^k)$ consisting of  strictly  positive constant forms $\tilde\gamma_1,\ldots,\tilde\gamma_{\tilde N},$
where $\tilde N:={k\choose p}^2.$
Since  $\Theta$ is  positive, we can  write  $\Theta\wedge \tilde\gamma_j=\tilde f_j(\zeta,w)\Leb(\zeta,w)$  for $1\leq j\leq \tilde N,$  where $\tilde f_j$ is a non-negative function.
The harmonicity of the functions  $\Theta_{I,J;K,L}$ with respect to $\zeta_{k-l}$ and the  constant forms $\tilde\gamma_1,\ldots,\tilde\gamma_{\tilde N},$
imply that $\tilde f_j$ is also   harmonic  functions with respect to $\zeta_{k-l}$.
So  the $\tilde f_j$'s are positive harmonic  functions with respect to $\zeta_{k-l}$.
The rest of the proof follows along the same lines as in the proof of assertion (3)  of Theorem \ref{T:conic-psh}.
\endproof


\part{Geometric  characterizations of the generalized  Lelong numbers}

\section{Grassmannian  bundles}
\label{S:Grassmannian}

\subsection{Grassmannian bundles and canonical projections}
The  following construction  which is  analog of   the blow-up is  necessary  in order  to
obtain  a geometric characterization of  the generalized Lelong numbers.
We  will  use  the notation  introduced  in   Sections  \ref{S:preparatory_results} and \ref{S:Lelong}.
So we  will keep the {\bf Standing Hypothesis}. In particular, let 
 $\omega$ be as usual a Hermitian  form on $V$
and  let  $\pi:\ \E\to  V$ be   the    normal bundle to $V$ in $X.$ 
For every $x\in V,$ $\E_x$ is as usual the fiber of $\E$ over $x$  which  is isomorphic to $\C^{k-l}.$
We  identify $x\in V$ with  the  vector  zero $0_x$  of $\E_x$ so that $V$ is canonically identified to the zero section $x\mapsto 0_x$ of $\E.$
Given a  $\C$-vector space $F$ of dimension $q$ and  an integer $j$ with  $1\leq j\leq q,$ let $\G_j(F)$  (resp. $\G_{j,q}$)
be the Grassmannian of all $j$-dimensional linear subspaces $H$ of  $F$ (reps.  of  $\C^q$).


In what follows  $j$ is  an integer with $1\leq j\leq k-l.$

Let $\pi_j:\ \G_j(\E)\to V$ be  the  holomorphic   bundle  which is obtained  from the vector bundle $\pi:\ \E\to V$ by   
 taking  the Grassmannian of all $j$-dimensional linear subspaces of  each fiber of $\E.$
 So, for $x\in V,$ 
the fiber  $\G_j(\E)_x$   of $\G_j(\E)$ over $x$ is  simply $\G_j(\E_x),$ the Grassmannian of all $j$-dimensional linear subspaces of $\E_x.$

Let  $\X_j=\X_j(\E)$ be  the holomorphic subbundle of the fibre product  (or equivalently, the Whitney sum) $\E\oplus \G_j(\E)$  whose  fiber  over  
every $x\in V$ is  given  by
$$
\X_j(x):=\left\lbrace  (y,H)\in \E_x\times \G_j(\E_x ):\quad  y\in H  \right\rbrace\subset \E_x\times \G_j(\E_x )=(\E\oplus \G_j(\E))_x.
$$
There are two natural holomorphic  bundles associated to each $\X_j$   corresponding  to the  projection  on the first  
 factor $\Pi_j:\ \X_j \to  \E$  (resp.  the  projection  on the  second factor   $\Pr_j:\ \X_j \to  \G_j(\E)$).
 
Consider the  holomorphic  bundle    corresponding  to the  projection  on the first  
 factor $\Pi_j:\ \X_j \to  \E.$  For  every $x\in V,$  
let   $\Pi_{j,x}$ be the restriction of  $\Pi_j$ to  $\X_j(x).$  So
we have  a holomorphic  bundle  $\Pi_{j,x}:\X_j(x)\to \E_x.$ Moreover, for every $y\in \E_x\setminus \{0_x\},$  the fiber  of $y$ is 
$$
\Pi_{j,x}^{-1}(y)=\left\lbrace (y,H):\  H\in \G_j(\E_x )\quad\text{and}\quad y\in H   \right\rbrace\simeq \G_{j-1,k-l-1}.
$$
Therefore, if $y\in \E_x\setminus \{0\},$ $\dim  \Pi_{j,x}^{-1}(y)=(j-1)(k-l-j).$ The exceptional  fiber of the bundle  $\Pi_{j,x}:\X_j(x)\to \E_x$
is 
$
\Pi_{j,x}^{-1}(0_x)\simeq \G_{j,k-l}$ has dimension  $j(k-l-j).$
Since  $\dim \X_j(x)=k-l+(j-1)(k-l-j),$ we  obtain $\codim  \Pi_{j,x}^{-1}(0_x)=j.$ 
  
  We have the  following   expression  for 
$\Pi_j:  \X_j  \to \E$:
$$\Pi_j(y,H)=\Pi_{j,\pi_j(y)}(y,H)\qquad\text{for}\qquad  (y,H)\in \X_j.$$
Since  $ \Pi_j^{-1}(V)=\{  \Pi_{j,x}^{-1}(0_x):\ x\in V\},$  we see that the complex  manifold $\Pi_j^{-1}(V)$ which is  a  holomorphic  bundle  over $V$ whose fibers
are $ \G_{j,k-l}$ has codimension $j$ in $\X_j.$
Observe that the restriction of $\Pi_j$  to $\X_j\setminus  \Pi_j^{-1}(V)$:
$$
\Pi_j^\bullet :\   \X_j\setminus  \Pi_j^{-1}(V)\to  \E\setminus V
$$
is  a bundle map  whose  fiber  over  $y\in \E\setminus V$ is
$$          (\Pi_j^\bullet)^{-1}(y)        :=               \Pi_{j,\pi(y)}^{-1}(y)$$
which is  of  dimension $(j-1)(k-l-j).$   We have  $\codim \Pi_j^{-1}(V)=j$ in $\X_j$ and 
$$\dim \X_j(x)=k+(j-1)(k-l-j).$$

Consider the  holomorphic  bundle    corresponding  to the  projection  on the second  
 factor $\Pr_j:  \X_j  \to \G_j(\E). $  For  every $x\in V,$  
let   $\Pr_{j,x}$ be the restriction of  $\Pr_j$ to  $\X_j(x).$  So
we have  a holomorphic  bundle  $\Pr_{j,x}:\X_j(x)\to \G_j(\E_x).$ Moreover, for every $H\in\G_j(\E_x),$  the fiber  of $H$ is
$$
\Pr\nolimits_j^{-1}(H)=\Pr\nolimits_{j,x}^{-1}(H)=\left\lbrace (y,H):\  y\in H   \right\rbrace\simeq H.
$$

\subsection{Canonical vertical forms}
Recall  that  $\varphi$ is the function  given by  \eqref{e:varphi}. 
Fix a point $x\in V$ and a Hermitian metric $h_x$ on $\E_x.$
There exists a  canonical  K\"ahler  form  $\Upsilon_{j,x}$  on $\G_j(\E_x)$  with respect to  $h_x.$
This  is  the  unique positive $(1,1)$-form  which is invariant under the action of the unitary group $\U(k-l)$ of of degree  $k-l$ (with respect to  $h_x$) and  which is  so normalized that
$\int_{\G_j(\E_x)}\Upsilon_{j,x}^{j(k-l-j)}=1.$  
 Consider the following  form on $\E_x:$  \begin{equation}
      \label{e:alpha_v,x}
      \alpha_{\ver,x}:=\ddcv \log\varphi\qquad\text{on}\qquad \E_x,
     \end{equation}
where  $\ddcv$ is  the operator defined in \eqref{e:ddc-ver}.  In other words, $\alpha_{\ver,x}$ is  just the  restriction of $\alpha_\ver$  (given in \eqref{e:alpha-beta-ver})  to $\E_x.$
We record the following important identity: 
\begin{lemma}
 \label{L:Siu} {\rm (Siu \cite{Siu})}
  For every $0\leq q\leq k-l-j,$  there is a constant $\gamma_{q,j}>0$ such that for every $x\in V,$
 \begin{equation*}
 (\Pi^\bullet_{j,x})_*(\Pr\nolimits_{j,x}^\bullet)^*  (\Upsilon_{j,x}^{q+(j-1)(k-l-j)})=\gamma_{q,j} \alpha_{\ver,x}^q.
 \end{equation*} 
\end{lemma}

Now  we construct a  form   $\Upsilon_j$ on $\G_j(\E)$  which is    the aggregate of the forms $\{\Upsilon_{j,x}\},$ where $x\in V.$
The holomorphic  bundle  $\pi_j:\ \G_j(\E)\to V$  allows  us  to  obtain  the  following canonical decomposition 
of the  holomorphic tangent  bundle  of $\G_j(\E)$ into  two parts: the horizontal  part $\Tan_\hor(\G_j(\E))$  and the vertical part
 $\Tan_\ver(\G_j(\E)).$  More specifically, for every $H\in  \G_j(\E),$  let $x=\pi_j(H)$ and  write
 $$\Tan_\hor(\G_j(\E))_H:= \Tan_x(V)\qquad\text{and}\qquad  \Tan_\ver(\G_j(\E))_H:=\Tan_H(\G_j(\E_x)).$$ 
Then  we have
$$
\Tan\G_j(\E)_H=\Tan_\hor(\G_j(\E))_H\oplus\Tan_\ver(\G_j(\E))_H,
$$
and   
$$
\Tan\G_j(\E)_H\otimes \C= \Tan\G_j(\E)^{1,0}_H\oplus \Tan\G_j(\E)^{0,1}_H,       
$$
where we have 
 \begin{eqnarray*} \Tan\G_j(\E)^{1,0}_H&:=& \Tan_\hor(\G_j(\E))_H^{1,0}\oplus\Tan_\ver(\G_j(\E))^{1,0}_H,\\
 \Tan\G_j(\E)^{0,1}_H&:=& \Tan_\hor(\G_j(\E))_H^{0,1}\oplus\Tan_\ver(\G_j(\E))^{0,1}_H.
 \end{eqnarray*}
This induces a dual decomposition for the cotangent  bundles
 $$
\Cotan\G_j(\E)_H=\Cotan_\hor(\G_j(\E))_H\oplus\Cotan_\ver(\G_j(\E))_H,
$$
and   
$$
\Cotan\G_j(\E)_H\otimes \C= \Cotan\G_j(\E)^{1,0}_H\oplus \Cotan\G_j(\E)^{0,1}_H,       
$$
where we have 
 \begin{eqnarray*} \Cotan\G_j(\E)^{1,0}_H&:=& \Cotan_\hor(\G_j(\E))_H^{1,0}\oplus\Cotan_\ver(\G_j(\E))^{1,0}_H,\\
 \Cotan\G_j(\E)^{0,1}_H&:= &\Cotan_\hor(\G_j(\E))_H^{0,1}\oplus\Cotan_\ver(\G_j(\E))^{0,1}_H.
 \end{eqnarray*}
 These  decompositions also induce the decomposition of the complex of $m$-form   on $\G_j(\E)$ into forms of type $(p,q)$ with $p+q=m$
 and into 
 a  canonical composition of vertical and horizontal forms:
 $$
 \bigwedge^{m} \Cotan\G_j(\E)\otimes \C= \bigoplus_{p+q=m} \Cotan \G_j(\E)^{p,q},
 $$
 where  the bundle  $ \Cotan \G_j(\E)^{p,q}$ is  equal to
 \begin{multline*}
 \bigwedge^p \Cotan\G_j(\E)^{1,0}\otimes \bigwedge^q \Cotan\G_j(\E)^{0,1}
 =\bigoplus_{p'+p''=p,q'+q''=q} \bigwedge^{p'} \Cotan_\hor \G_j(\E)^{1,0}\\
 \otimes       \bigwedge^{p''} \Cotan_\ver\G_j(\E)^{1,0}
 \bigotimes\bigwedge^{q'} \Cotan_\hor \G_j(\E)^{0,1}\otimes       \bigwedge^{q''} \Cotan_\ver\G_j(\E)^{0,1}.
 \end{multline*}
 Consider  a smooth test form $\Psi$  compactly  supported in $\G_j(\E).$
 Let  $\widetilde\Psi$ be the sum  of all components corresponding to $p'=l,q'=l$ in the above  decomposition.  So $\widetilde\Psi$ is uniquely determined 
 by $\Psi.$
 Let $\widehat\Psi$ be the unique  form in $\bigwedge^{q'} \Cotan_\ver \G_j(\E)^{1,0}\otimes       \bigwedge^{q''} \Cotan_\ver\G_j(\E)^{0,1}$
 such that $\widetilde\Psi=(\pi_j^*\omega)^l\cdot  \widehat\Psi.$ 
 So  $\widehat\Psi$ is uniquely determined 
 by $\Psi.$
 
For every $x\in V$ let $\iota_x$ the canonical injection $\iota_x:\  \E_x\hookrightarrow  \E.$ It induces  canonically the  injection (still denoted by) $\iota_x:\  \G_j(\E_x)\hookrightarrow  \G_j(\E).$   We  are in the position to define the form $\Upsilon_j$ on $\G_j(\E)$ as  follows:
 \begin{equation}\label{e:Upsilon_j}
  \langle \Upsilon_j,\Phi\rangle:= \int_{x\in V} \langle \Upsilon_{j,x},  \iota_x^*\widehat \Phi\rangle \omega^l(x).
 \end{equation}
Then  $\Upsilon_j$  is  a uniquely-defined    positive smooth $(1,1)$-form  on  $\G_j(\E).$ Note that   $d_\ver  \Upsilon_j=0$ (see \eqref{e:ddc-ver} for the  definition of $d_\ver$) but  $d\Upsilon_j$ may not be $0,$ 
in other words, $\Upsilon_j$ is vertically closed 
but it is not necessarily closed.
Moreover,   it satisfies the  identity
$$
\iota_x^*\Upsilon_j =\Upsilon_{j,x}\qquad\text{on}\qquad  \G_j(\E_x).
$$


Consider the  following canonical decomposition 
of the  holomorphic tangent  bundle  of $\pi:\ \E\to  V$ into  two parts: the horizontal  part $\Tan_\hor(\E)$  and the vertical part
 $\Tan_\ver(\E).$  More specifically, for every $y\in  \E,$  let $x=\pi(y)$ and  write
 $$\Tan_\hor(\E)_y:= \Tan_x(V)\qquad\text{and}\qquad  \Tan_\ver(\E)_y:=\E_x.$$ 
Then  we have
$$
\Tan(\E)_y=\Tan_\hor(\E)_y\oplus\Tan_\ver(\E)_y,
$$
and   
$$
\Tan(\E)_y\otimes \C= \Tan(\E)^{1,0}_y\oplus \Tan(\E)^{0,1}_y,       
$$
where we have 
 \begin{eqnarray*} \Tan(\E)^{1,0}_y&:=& \Tan_\hor(\E)_y^{1,0}\oplus\Tan_\ver(\E)^{1,0}_y,\\
 \Tan(\E)^{0,1}_y&:=& \Tan_\hor(\E)_y^{0,1}\oplus\Tan_\ver(\E)^{0,1}_y.
 \end{eqnarray*}
This induces a dual decomposition for the cotangent  bundles
 $$
\Cotan(\E)_y=\Cotan_\hor(\E)_y\oplus\Cotan_\ver(\E)_y,
$$
and   
$$
\Cotan(\E)_y\otimes \C= \Cotan(\E)^{1,0}_y\oplus \Cotan(\E)^{0,1}_y,       
$$
where we have 
 \begin{eqnarray*} \Cotan(\E)^{1,0}_y&:=& \Cotan_\hor(\E)_y^{1,0}\oplus\Cotan_\ver(\E)^{1,0}_y,\\
 \Cotan(\E)^{0,1}_y&:= &\Cotan_\hor(\E)_y^{0,1}\oplus\Cotan_\ver(\E)^{0,1}_y.
 \end{eqnarray*}
 These  decompositions also induce the decomposition of the complex of $m$-form   on $\E$ into forms of type $(p,q)$ with $p+q=m$
 and into 
 a  canonical composition of vertical and horizontal forms:
 $$
 \bigwedge^{m} \Cotan(\E)\otimes \C= \bigoplus_{p+q=m} \Cotan (\E)^{p,q},
 $$
 where  the bundle  $ \Cotan (\E)^{p,q}$ is  equal to
 \begin{multline*}
 \bigwedge^p \Cotan(\E)^{1,0}\otimes \bigwedge^q \Cotan(\E)^{0,1}
 =\bigoplus_{p'+p''=p,q'+q''=q} \bigwedge^{p'} \Cotan_\hor (\E)^{1,0}\\
 \otimes       \bigwedge^{p''} \Cotan_\ver(\E)^{1,0}
 \bigotimes\bigwedge^{q'} \Cotan_\hor (\E)^{0,1}\otimes       \bigwedge^{q''} \Cotan_\ver(\E)^{0,1}.
 \end{multline*}
Consider  a smooth test form $\Phi$  compactly  supported in $\E.$
 Let  $\widetilde\Phi$ be the sum  of all components corresponding to $p'=l,q'=l$ in the above  decomposition.  So $\widetilde\Phi$ is uniquely determined 
 by $\Phi.$
 Let $\widehat\Phi$ be the unique  form in $\bigwedge^{q'} \Cotan_\ver (\E)^{1,0}\otimes       \bigwedge^{q''} \Cotan_\ver(\E)^{0,1}$
 such that $\widetilde\Phi=(\pi_j^*\omega)^l\cdot  \widehat\Phi.$ 
 So  $\widehat\Phi$ is uniquely determined 
 by $\Phi.$

 Consider the  form
 $\Psi:=  (\Pr\nolimits_j)_\diamond (\Pi_j)^\diamond\Phi$ on $\G_j(\E).$
 Observe  that
 \begin{equation}\label{e:Phi-Psi}
  \widehat\Psi= (\Pr\nolimits_j)_\diamond (\Pi_j)^\diamond(\widehat\Phi).
 \end{equation}

Let $
\Pr_j^\bullet$  (resp.    $
\Pr\nolimits_{j,x}^\bullet$ for  each  $x\in V$)   be the restriction of $\Pr_j$ to $\X_j\setminus  \Pi_j^{-1}(V)$  (resp.   the restriction of $\Pr_{j,x}$ to $\X_j(x)\setminus  \Pi_{j,x}^{-1}(0_x)$).  So we obtain the holomorphic  bundles
$$
\Pr\nolimits_j^\bullet :\   \X_j\setminus  \Pi_j^{-1}(V)\to  \G_j(\E) \qquad\text{and}\quad
\Pr\nolimits_{j,x}^\bullet :\   \X_j(x)\setminus  \Pi_{j,x}^{-1}(0_x)\to  \G_j(\E_x)
.$$
Moreover, for every $H\in\G_j(\E_x),$  the fiber  of $H$ is
$$
(\Pr\nolimits^\bullet_j)^{-1}(H)=(\Pr\nolimits^\bullet_{j,x})^{-1}(H)=\left\lbrace (y,H):\  y\in H\setminus \{0\}   \right\rbrace\simeq H\setminus \{0\}.
$$

 Recall from  \eqref{e:alpha-beta-ver} the following  form on $\E:$  \begin{equation}
      \label{e:alpha_v}
      \alpha_\ver:=\ddcv \log\varphi.
     \end{equation}
So $\alpha_\ver$ is a smooth positive $(1,1)$-form, it is  vertically  closed  but it is  not necessarily  closed.

\begin{notation}\rm
 \label{N:Diamond}
 For  $1\leq j\leq k-l,$  we  set
 \begin{eqnarray*}
  (\Pi_j)_\diamond&:=&(\Pi^\bullet_j)_*\qquad\text{and}\qquad (\Pi_j)^\diamond:=(\Pi^\bullet_j)^*,\\
 (\Pr\nolimits_j)_\diamond  &:= &(\Pr^\bullet\nolimits_j)_*\qquad\text{and}\qquad (\Pr\nolimits_j)^\diamond:= (\Pr^\bullet\nolimits_j)^*.
 \end{eqnarray*}
\end{notation}

 \begin{lemma}\label{L:Upsilon_j}
 For every $0\leq q\leq k-l-j,$ there is a constant $\gamma_{q,j}>0$ such that 
 \begin{equation*}
 (\Pi_j)_\diamond(\Pr\nolimits_j)^\diamond  (\Upsilon_j^{q+(j-1)(k-l-j)})=\gamma_{q,j} \alpha_\ver^q.
 \end{equation*} 
 \end{lemma}
 \proof
 Let $\Phi$ be a  smooth compactly  supported test form of bidegree $(k-q,  k-q)$ on $\E.$  Consider the  form
 $\Psi:=  (\Pr\nolimits_j)_\diamond (\Pi_j)^\diamond\Phi$ on $\G_j(\E).$
 We need to show that
 $$
  \langle \Upsilon_j^{q+(j-1)(k-l-j)},\Psi \rangle=\gamma_j\langle \alpha_\ver^q,\Phi\rangle.
 $$
 By \eqref{e:Upsilon_j} the LHS is  equal to
 $$
  \int_{x\in V} \langle \Upsilon_{j,x}^{q+(j-1)(k-l-j)}, \iota_x\widehat \Psi\rangle\omega^l(x).
 $$
 By  Lemma \ref{L:Siu}  and equality \eqref{e:Phi-Psi}, this expression is equal to
 $$
  \gamma_j\int_{x\in V} \langle \alpha_{\ver,x}^q ,\iota_x \widehat\Phi\rangle\omega^l(x).
 $$
 By \eqref{e:alpha_v,x} and \eqref{e:alpha_v},  the last expression is equal to  $\gamma_j\langle\alpha_\ver^q, \Phi\rangle.$ 
 \endproof

Consider  the  following natural positive smooth form  $(1,1)$-form on $\X_j$:
\begin{equation}\label{e:omega_j}
\omega_j:=\Pi_j^* (c_1\pi^*\omega+\beta_\ver) +(\Pr\nolimits_j)^*\Upsilon_j.
\end{equation}
Here, we recall  from  \eqref{e:alpha-beta-ver} that
\begin{equation}
 \label{e:beta_v}
 \beta_\ver:=\ddcv  \varphi.
\end{equation}



 \section{$\C$-flatness and extension currents}
 \label{S:Flatness}

 \subsection{$\C$-normal currents, Federer-type $\C$-flatness Theorem }
 
 Recall some definitions and  results  of Bassanelli  \cite{Bassanelli}   (see also  Sibony \cite{Sibony} for related  notions and  results
 on pluripositive  currents). 
 \begin{definition}\label{D:flat-normal} \rm 
 Let $T$ be  a current    on an open set  $\Omega$ in a complex manifold of dimension $k.$
 We say that $T$ is $\C$-flat if there exist     currents $S, $  $G$ and $H$ on $\Omega$   with coefficients in $L^1_\loc(\Omega)$ such that
 \begin{equation*}
  T=S+\partial G+\dbar H\qquad \text{on}\qquad \Omega.
 \end{equation*}
We say that   $T$  is $\C$-normal   if  $T$ and $\ddc T$  have measure coefficients
 \end{definition}
By  \cite[Theorem 1.18]{Bassanelli},  $\C$-normal currents  are $\C$-flat.

Let $F$ be a closed  subset of $\Omega.$  If  $T$ is a  current on $\Omega\setminus F$ with locally finite mass across $F,$
then there exists a  current, denoted by $T_\bullet,$  which is  the trivial  extension  of $T$ to $\Omega.$
More precisely,  $T_\bullet$ coincides with $T$ on $\Omega\setminus F$ and  $T_\bullet$ has no mass on $F.$

\begin{proposition} \label{P:Bassanelli} {\rm  (see  \cite[Proposition  1.22 and Lemma 1.11]{Bassanelli})} 
\begin{enumerate}
\item If $T$ is $\C$-flat current  on $\Omega\setminus F$  with  locally finite mass across $F,$ then      $T_\bullet$ is $\C$-flat on $\Omega.$
\item If $R$ is a $\C$-flat current with  measure coefficients on $\Omega,$ then $\ind_F R$ is $\C$-flat.
\end{enumerate}
\end{proposition}
The  following  Federer-type $\C$-flatness theorem is very useful.
\begin{theorem}\label{T:Bassanelli}
 Let $Z$ be an analytic subset of $\Omega$ and let $T$ be a $\C$-flat positive current of bidimension $(p,p)$ on $\Omega,$  supported in $Z.$
 Then there is a unique current $S$ of bidimension $(p,p)$ on $Z$ such that $T=\iota_*S,$ where $\iota:\ Z\hookrightarrow \Omega$
 is the canonical inclusion.
\end{theorem}
\proof  All  assertions except  the positivity of $S$  have been proved in   \cite[Theorem  1.24]{Bassanelli}. 
But the positivity of $S$ follows easily from  that of $T.$
\endproof
Let $T$ be  positive $(p,p)$-current which is  $\C$-flat on $X.$ Let $Z$ be an analytic subset of $X$ of pure codimension $p.$
The current $R:= T|_{X\setminus Z}$   is $\C$-flat. As $T\geq  R\geq 0,$  $R$ has locally finite mass across $Z,$  and hence we get that
$$
T=\ind_Z T+R_\bullet.
$$
Since  $\ind_Z T$ is  $\C$-flat  and positive,  by Theorem \ref{T:Bassanelli},  we get
$$
\ind_Z T=f[Z]
$$
for a suitable   function $f\in L^1_\loc(Z).$  Therefore,  it follows that
$$
T=f[Z]+R_\bullet.
$$
Moreover, if $T$ is  positive  plurisubharmonic,  Bassanelli proves the following
\begin{theorem}\label{T:Bassanelli-psh}{\rm  (\cite[Theorem 4.10]{Bassanelli})}
 If $T$ is a  positive plurisubharmonic current of bidegree $(p,p)$ on $X$ and $Z$  is an analytic  subset of $X$ of pure codimension
 $p,$ then there exists a weakly plurisubharmonic  function $f:\ Z\to  \R,$ $f\geq 0,$   such that $ \ind_Z T=f[Z].$
\end{theorem}

 \subsection{Extension  currents}

 \begin{lemma}\label{L:Siu-inequ}
  Let $\Omega,\Omega'$ be   open subsets of $\U\subset \E$  with $\Omega\Subset \Omega'\Subset \U.$ Let $1\leq p\leq k-l.$ Then  there is a constant $c>0$  such that
  for every positive  smooth  $(p,p)$-form $R$ on $\Omega\setminus V$  and every integer $j$ with $0\leq j\leq k-l,$  we have
  $$
  \int_{\Pi_{j}^{-1}(\Omega\setminus V)}\big((\Pi_{j})^\diamond R\big)\wedge\omega_{j}^{\dim \X_j-p}\leq c\sum_{0\leq q\leq k-l-j}
  \int_{\Omega'\setminus V}  R\wedge \alpha_\ver^{q}\wedge  (c_1\pi^*\omega + \beta_\ver)^{k-p-q}.$$
 \end{lemma}
\proof
By formula \eqref{e:omega_j}, the  expression on the LHS is  equal to
\begin{equation*}
\sum_{q=-(j-1)(k-l-j)}^{k-p}  {\dim \X_j-p \choose k-p-q  }\int_{\Pi_{j}^{-1}(\Omega\setminus V)}\big((\Pi_{j})^\diamond R\big)\wedge\big( \Pi_j^* (c_1\pi^*\omega+\beta_\ver)^{k-p-q}\big) \wedge \big((\Pr\nolimits_j)^\diamond \Upsilon_j\big)^{q+(j-1)(k-l-j)}.
\end{equation*}
Since the  fiber of $\Pi_j^\bullet$ is of dimension $(j-1)(k-l-j),$  all the  integrals corresponding to  $q<0$ on the  RHS  vanish.
On the other hand,  all the  integrals corresponding to  $q>k-l-j$ on the  RHS  vanish because the $(1,1)$-form $\Upsilon_{j,x}$ lives on $\G_j(\E_x)$ whose dimension is $j(k-l-j).$

Applying  Lemma \ref{L:Upsilon_j} for $0\leq q\leq  k-l-j,$ the  last expression is  equal to
\begin{equation*}
\sum_{q=0}^{k-l-j}  {\dim \X_j-p \choose k-p-q  } \int_{\Omega\setminus V}  R\wedge \alpha_\ver^{q}\wedge  (c_1\pi^*\omega + \beta_\ver)^{k-p-q}.
\end{equation*}
The  result follows.
\endproof
 
 \begin{proposition}\label{P:existence}
 Let $\Omega$ be an  open subset of $\U\subset \E.$  
  Let $(R_n)$  be a sequence of  positive smooth  $(p,p)$-forms on $\Omega\setminus V.$
  Let $j$ be an integer with $1\leq j\leq k-l.$
  Assume that   
  $$
  \sup\limits_{n\in\N}\int_{\Omega\setminus V}  R_n\wedge (\hat \alpha')^{q}\wedge  \pi^*(\omega^{m})\wedge \beta_\ver^{k-p-q-m}<\infty 
  $$
  for every   $0\leq q\leq k-l-j$  and $0\leq m\leq k-p-q.$
  Then there exist  currents  $\widehat R^{(j,q)},$ $ R^{(j,q)}$  for $0\leq q\leq k-l-j$ on  $\Omega$  and  $\widetilde R^{(j)}$ on $\Pi_{j}^{-1}(\Omega)$ such that, for a suitable subsequence $(R_{N_n})$ the following properties hold:
  \begin{enumerate}
  \item $\lim_{n\to \infty} \big(R_{N_n}\wedge (\hat\alpha')^q\big)_\bullet =\widehat R^{(j,q)}$  weakly on $\Omega.$
   \item $\lim_{n\to \infty} \big(R_{N_n}\wedge \alpha_\ver^q\big)_\bullet =R^{(j,q)}$  weakly on $\Omega.$
   \item  $\lim_{n\to \infty} \big(\Pi_j^\diamond R_{N_n}\big)_\bullet =\widetilde R^{(j)}$  weakly on $\Pi_j^{-1}(\Omega).$
  \end{enumerate}

 \end{proposition}
\proof
The  assumption  implies that  for every $0\leq q\leq k-l-j,$
$$
  \sup\limits_{n\in\N}\int_{\Omega\setminus V}  R_n\wedge  (\hat\alpha')^q\wedge (c_1 \pi^*(\omega)+ \beta_\ver)^{k-p-q}<\infty. 
  $$
  Since  $c_1\pi^*(\omega)+ \beta_\ver $ is a smooth  strictly positive $(1,1)$-form on $\U,$  assertion (1) follows from  Lemma  \ref{L:subsequence}.
  
  Using  inequality  \eqref{e:hat-alpha'-vs-alpha_ver}, the  assumption  implies that  for every $0\leq q\leq k-l-j,$
$$
  \sup\limits_{n\in\N}\int_{\Omega\setminus V}  R_n\wedge \alpha_\ver^{q}\wedge (c_1 \pi^*(\omega)+ \beta_\ver)^{k-p-q}<\infty. 
  $$
  Therefore, arguing  as  in the  proof of assertion (1),  assertion (2)  follows.
  
  Applying Lemma \ref{L:Siu-inequ} to  each $R_n$ yields that
  $$
   \sup\limits_{n\in\N}\int_{\Pi_{m}^{-1}(\Omega\setminus V)}\big( \Pr\nolimits_{m}^\diamond R\big)\wedge\omega_{m}^{\dim \X_{m}-p} <\infty.$$
  So assertion (2) follows from  Lemma  \ref{L:subsequence}.
\endproof

\begin{proposition} \label{P:existence-finite}
 Let $T$ be  a positive plurisubharmonic current in the  class $\SH^{3,3}_p(B)$ with a sequence of approximating forms $(T_n)_{n=1}^\infty.$
 Let $1\leq \ell\leq \ell_0.$ Then:
 \begin{enumerate} \item The assumption of Proposition \ref{P:existence} is  satisfied for $R_n:=(\tau_\ell)_*T_n$ and $\Omega:=\U_\ell:=\tau_\ell(\bfU_\ell)\subset \E.$
\item There exist  currents  $\widehat R^{(j,q)}_{[\ell]},$ $ R^{(j,q)}_{[\ell]}$  for $0\leq q\leq k-l-j$ on  $\Omega$  and  $\widetilde R^{(j)}_{[\ell]}$ on $\Pi_{j}^{-1}(\U_\ell)$ such that, for a suitable subsequence $(R_{N_n})$ the following properties hold:
  \begin{enumerate}
  \item $\lim_{n\to \infty} \big(R_{N_n}\wedge (\hat\alpha')^q\big)_\bullet =\widehat R^{(j,q)}_{[\ell]}$  weakly on $\U_\ell.$
   \item $\lim_{n\to \infty} \big(R_{N_n}\wedge \alpha_\ver^q\big)_\bullet =R^{(j,q)}_{[\ell]} $ weakly on $\U_\ell.$
   \item  $\lim_{n\to \infty} \big(\Pi_j^\diamond R_{N_n}\big)_\bullet =\widetilde R^{(j)}_{[\ell]}$  weakly on $\Pi_j^{-1}(\U_\ell).$
  \end{enumerate}

 \end{enumerate}
 \end{proposition}
\proof
We may assume without loss of generality that $T$ is in the class $\widetilde\SH^{3,3}_p(\bfU,\bfW).$
By Theorem \ref{T:Lc-finite-psh}, $\Kc_{j,q}(R_n,\bfr)\leq c_{10}.$
Using this and \eqref{e:hat-alpha-vs-alpha_ver} and \eqref{e:hat-beta-vs-beta_ver} and \eqref{e:T-hash}, the first assertion  follows.

Using the  first assertion, the second one is a consequence of  Proposition \ref{P:existence}.
\endproof

\begin{definition}
 \rm  \label{D:Cut-Off} 
 Let $S$ be  a real current  defined on an open set $\Omega$ and $V$ an  analytic subset of $\Omega.$
 We  say that $S$ enjoy the cut-off property through $V$ in  $\Omega$ if  the  following decomposition  holds
 $$  S=\ind_VS+ (S|_{\Omega\setminus V})_\bullet,$$
 where  $(S|_{\Omega\setminus V})_\bullet$  is the trivial  extension through $V$ to $\Omega$ of the  current $S|_{\Omega\setminus V}$, which is the restriction
 of $S$ to $\Omega\setminus V.$ 
\end{definition}

\begin{lemma}\label{L:Cut-Off}
Let $S$ be  a real current  defined on an open set $\Omega$ and $V$ an  analytic subset of $\Omega.$
Suppose that for every $x\in\Omega$ there is a $\Cc^2$-diffeomorphism $\tau_x: \ U_x\to W_x,$ where $U_x,$ $W_x$ are  open  neighborhood  of $x$
in $\Omega$ with the  following  properties:
\begin{itemize}
\item $\tau_x$ is  admissible along $V\cap U_x;$
 \item 
$(\tau_x)_*(S|_{U_x})$ is $\C$-flat positive  current on $W_x;$
\item  $\tau_x|_{V\cap U_x}$ is the identity.
\end{itemize}
 Then  $S$ enjoy the cut-off property through $V$ in  $\Omega$ 
\end{lemma}
\proof
Since  the  problem is local and $\tau_x|_{V\cap U_x}$ is the identity,  we may  work  locally with $(\tau_x)_*(S|_{U_x})$ instead of $S|_{U_x}.$
Therefore, we may assume without loss of generality that  $S$ is a $\C$-flat positive current.
The result follows then  from   Definition \ref{D:Cut-Off} and  Proposition \ref{P:Bassanelli} and Theorem \ref{T:Bassanelli}.
\endproof

\begin{lemma}\label{L:closed-ph-psh-smooth} Let $\lowm\leq j\leq \upm$ and $1\leq \ell\leq \ell_0.$ Let $T$  be a  real current on $\bfU_\ell.$ Let $\tau$ be a  holomorphic   admissible map  from $\bfU_\ell$ onto $\tau(\bfU_\ell).$
  Consider the current  
  $ R:=\tau_* (T) \wedge  (\hat\alpha')^{k-p-j}$  on $\tau(\bfU_\ell)\setminus  V.$ 
  Then the following  assertions hold:
  \begin{enumerate}  \item If $T$  is a positive closed  $\Cc^1$-smooth current on $\bfU_\ell$  then $R_\bullet$ is a positive closed current  on $\tau(\bfU_\ell).$
   \item  If $T$  is a positive pluriharmonic $\Cc^2$-smooth current on $\bfU_\ell$  then $R_\bullet$ is a positive pluriharmonic current  on $\tau(\bfU_\ell).$
   \item  If $T$  is a positive plurisubharmonic $\Cc^2$-smooth current on $\bfU_\ell$  then $R_\bullet$ is a positive plurisubharmonic current  on $\tau(\bfU_\ell).$
  \end{enumerate}
 
\end{lemma}
 
\proof
First  we prove assertion (2). The proof of  assertion (1) is  similar.
Let $R$ be a  smooth  differential form   compactly supported in $\U \cap \pi^{-1}(B).$   We have
\begin{eqnarray*}
 &&\langle\ddc \big(\tau_*T\wedge  \alpha^{n-p-j}\big)_\bullet,  R\rangle=\lim_{\epsilon\to 0} \int_{\U\setminus  \overline{\Tube} (B,\epsilon) } \tau_*T\wedge  \alpha^{n-p-j}\wedge \ddc R\\
 &=&\lim_{\epsilon\to 0}\big(  \int_{\partial_\hor\Tube (B,\epsilon) }- \tau_*T\wedge  \alpha^{n-p-j}\wedge i\dbar R -
  \int_{\partial_\hor\Tube (B,\epsilon) }i\partial(\tau_* T)\wedge  \alpha^{n-p-j}\wedge R\big),
\end{eqnarray*}
where the  second equality  holds  because $\ddc(\tau_*T)= \tau_* (\ddc T)=0.$
By  Lemma \ref{L:Lelong-Jensen_1},  $j^*_\epsilon (\alpha)={1\over \epsilon^2} j^*_\epsilon (\beta).$
Moreover,   $j^*_\epsilon (\alpha^{k-l})=0$  since  $\partial_\hor \Tube(B,\epsilon)$ has real dimension $2(k-l)-1.$
Using  these two equalities and the smoothness of $T$ and $R,$ we  can check that both integral  in the last line  are of order $O(\epsilon).$
Letting $\epsilon\to 0,$ we infer that   $\langle\ddc \big(\tau_*T\wedge  \alpha^{n-p-j}\big)_\bullet,  R\rangle=0$ as  desired.

To  prove  assertion (3), we pick  a  positive  smooth  differential form  $R$ compactly supported in $\U \cap \pi^{-1}(B).$
Since $T$ is   plurisubharmonic, the current   $\ddc \big(\tau_*T\wedge  \alpha^{n-p-j}\big)=\tau_*(\ddc T)\wedge  \alpha^{n-p-j}\big)$
is positive. Therefore, we get that
 $$\int_{\U\setminus  } \tau_*T\wedge  \alpha^{n-p-j}\wedge \ddc R\geq 0.$$
 On the other hand,  using the smoothness of $T$ and $R$ and applying  Lemma \ref{L:vertical-boundary-terms}, we can show that 
 $$
 \lim_{\epsilon\to 0} \int_{\Tube (B,\epsilon) } \tau_*T\wedge  \alpha^{n-p-j}\wedge \ddc R=0.
 $$
 This, combined with the previous inequality, implies that 
 \begin{equation*}
 \langle\ddc \big(\tau_*T\wedge  \alpha^{n-p-j}\big)_\bullet,  R\rangle=\lim_{\epsilon\to 0} \int_{\U\setminus  \overline{\Tube} (B,\epsilon) } \tau_*T\wedge  \alpha^{n-p-j}\wedge \ddc R\geq  0.
 \end{equation*}
\endproof

  \begin{corollary}\label{C:Cut-off}
We keep  the hypothesis and  the  conclusion of Proposition \ref{P:existence-finite} and
 let $1\leq \ell\leq \ell_0.$ Then the  currents  $\widehat R^{(j,q)}_{[\ell]},$ $ R^{(j,q)}_{[\ell]}$  enjoy the cut-off property through $V$ in  $\\U_\ell$  and the  current $\widetilde R^{(j)}_{[\ell]}$ enjoys the cut-off property through $\Pi_j^{-1}(V)$ in  $\Pi_{j}^{-1}(\U_\ell).$ 
  \end{corollary}

  \proof
  Combining Lemma \ref{L:closed-ph-psh-smooth} and
Proposition \ref{P:existence-finite}, we  see that the   currents  $\widehat R^{(j,q)}_{[\ell]},$ $ R^{(j,q)}_{[\ell]}$  are positive plurisubharmonic  on  $\U_\ell$  and  the  current $\widetilde R^{(j)}_{[\ell]}$ is positive plurisunharmmonic  on $\Pi_{j}^{-1}(\U_\ell).$
Hence, by  Theorem \ref{T:Bassanelli-psh}, these currents  enjoy  the corresponding  cut-off property. 
  \endproof

\section{Geometric  characterizations for positive closed  and  positive  pluriharmonic currents with holomorphic admissible maps}
\label{S:Charac-pos-closed-and-plurihar-I}

This  section is   devoted to geometric  characterizations  of the  generalized  Lelong  numbers 
for positive closed  currents   and   positive pluriharmonic  currents  with  holomorphic  admissible maps. 

For  $j$ with $\lowm \leq j\leq \upm,$ define
\begin{equation}\label{e:hatj}{\hat j}:=j+p-l.\end{equation}
Note  that ${\hat j}\in [0,k-l]$ and  ${\hat j}+(k-p-j)=k-l.$

\begin{proposition}\label{P:Cut-Off-hol}
Let $T$  be a  current in the class $ \SH^{2,2}_p(B)$ (resp.  $ \PH^{2,2}_p(B)$,  resp.  $\CL^{1,1}_p(B)$) introduced in Definition  \ref{D:classes}
with   an approximating sequence of $(T_n)_{n=1}^\infty.$
Let $\tau $ be a    holomorphic  admissible map along $B.$ 
Consider the  real currents $R_n:=\tau_*(T_n)$ on $\U$  for $n\geq 1.$
Then  the following   assertions hold:

\begin{enumerate}

\item   There exist  currents   $\widehat R^{(\hat j,k-p-j)},$ $R^{(\hat j,k-p-j)}$   on  $\U$  and  $\widetilde R^{(\hat j)}$ on $\Pi_{j}^{-1}(\U)$ such that, for a suitable subsequence $(R_{N_n})$ of  the  sequence $\big(R_n\big)_{n=1}^\infty,$ the following properties hold:
  \begin{enumerate}
   \item $\lim_{n\to \infty} \big(R_{N_n}\wedge(\hat \alpha')^{k-p- j}\big)_\bullet = \widehat R^{(\hat j,k-p-j)} $  weakly on $\U;$
   \item $\lim_{n\to \infty} \big(R_{N_n}\wedge \alpha_\ver^{k-p- j}\big)_\bullet = R^{(\hat j,k-p-j)} $  weakly on $\U;$
   \item  $\lim_{n\to \infty} \big(\Pi_{\hat j}^\diamond R_{N_n}\big)_\bullet =\widetilde R^{(\hat j)}$  weakly on $\Pi_{\hat j}^{-1}(\U).$
  \end{enumerate}

 
\item  Set   $T^{(j)} :=R^{(\hat j,k-p-j)}$ and  $\widehat T^{(j)} :=\widehat R^{(\hat j,k-p-j)}$ 
and    $\widetilde T^{(j)} :=\widetilde R^{(\hat j)}.$  
Then   $T^{(j)} $ and $\widehat T^{(j)} $ enjoy the cut-off property  through $V$ in $\U,$  and  
$ \widetilde T^{(j)} $  enjoys the cut-off property  through $\Pi_{{\hat j}}^{-1}(V)$ in $\Pi_{{\hat j}}^{-1}(\U).$

\item 
There exist    positive currents $f_j,$ $\hat f_j$  of bidegree $(l-j,l-j)$ on $B$ and  a  positive $(l-j,l-j)$-current  $\tilde f_j$ on $\Pi_{\hat j}^{-1}(B)$  such that 
$$
\ind_{B}( T^{(j)})=(\iota_{B,\E})_*(f_j)\quad\text{and}\quad \ind_{B}( \widehat T^{(j)})=(\iota_{B,\E})_*(\widehat f_j)  \quad\text{and}\quad  \ind_{\Pi_{\hat j}^{-1}(B)}( \widetilde T^{(j)})= 
(  \iota_{\Pi_{\hat j}^{-1}(B),\X_{\hat j}})_*  (\tilde f_j). 
$$
\item If $j=l$  then $\hat f_j=f_j$ and   $f_j$ and $\tilde f_j$ are positive  plurisubharmonic  functions.
Moreover,  if   $T$ belongs to $ \CL^{2,2}_p(B),$ then  $f_j,$ $\hat f_j$ and $\tilde f_j$ are non-negative constant.

\end{enumerate}
\end{proposition}
\proof
{\bf Proof of assertion (1).}   By  Proposition  \ref{P:existence-finite},
   The assumption of  Proposition \ref{P:existence} is  satisfied
for $ R_n:=(\tau_\ell)_*(T_n) $ and $\Omega:=\U_\ell.$ 
Note that  $ k-l-\hat j=k-p-j.$ Consequently, 
  the  sequences $R_n\wedge (\hat\alpha')^{ k-p-j},$ $R_n\wedge \alpha_\ver^{ k-p-j}$  and  $(\Pi_{{\hat j}})^\diamond(R_n)$ are relatively compact in the  weak-$\star$ topology.
  we obtain  by Proposition \ref{P:existence} the  existence of the  desired  currents. This  completes the proof of assertion (1).

By Lemma  \ref{L:closed-ph-psh-smooth},   $\hat T^{(j)},$
and    $\widetilde T^{(j)} $  are  positive plurisubharmonic  currents. Hence, by Theorem \ref{T:Bassanelli}, $\hat T^{(j)}$
 enjoys the cut-off property  through $V$ in $\U,$  and  
$ \widetilde T^{(j)} $  enjoys the cut-off property  through $\Pi_{{\hat j}}^{-1}(V)$ in $\Pi_{{\hat j}}^{-1}(\U).$
Using  identity \eqref{e:hat-alpha'} we can express   $T^{(j)} $ as  a linear combinations with  real coefficients of   $\widehat T^{(j+m)}\pi^*\omega^m$
for $0\leq m\leq \upm-j.$  Hence,   $ T^{(j)}$ also
 enjoys the cut-off property  through $V$ in $\U.$  We  obtain  the  desired conclusion of assertions (2) and  (3).
 
 Assertion (4)  follows from Theorem \ref{T:Bassanelli-psh}.

\endproof

 \begin{proposition}\label{P:Cut-Off-hol-rep}
 We keep the assumption, notation and  conclusion of Proposition  \ref{P:Cut-Off-hol}.
 Then  the following  assertions hold.
 \begin{enumerate}
  \item For  every subdomain  $D\Subset B,$ we have 
  \begin{equation*}
  \int_{D} f_j \wedge \pi^*\omega^j= 
  \int_{(\Pi_{\hat j})^{-1}(D)}  \tilde f_j\wedge\Upsilon_{\hat j}^{\dim\X_{\hat j}- p-j}\wedge \Pi_{\hat j}^*(\pi^*(\omega^j)).
 \end{equation*}
\item  When  $j=l$ (so $j=l=\upm),$ then   
$f_{\upm}$ and  $\tilde f_{\upm}$   are functions related by
$$
f_{\upm}(x)=\int_{\Pi_p^{-1}(0_x)} \tilde f_\upm \Upsilon_p^{p(k-l-p)}
$$
for Lebesgue almost every $x\in B.$
\end{enumerate}
 \end{proposition}
 
 \proof

 By assertion (4) of  Proposition \ref{P:Cut-Off}, we can write
 \begin{equation*}
 T^{(j)}=(\iota_{V,\E})_*(f_j)+P\qquad\text{and}\qquad   \widetilde T^{(j)}= ( \iota_{\Pi_{\hat j}^{-1}(V),\X_{\hat j}  })_*(\tilde f_j)+Q, 
 \end{equation*}
where $P$ and $Q$ are positive  currents whose masses vanish  on  $V$ and $\Pi_{\hat j}^{-1}(V)$ respectively.

Let $D\Subset B$ be a  subdomain.  By Lemma \ref{L:Siu}, we have  that
\begin{eqnarray*}
 &&\int_{D} f_j \wedge \pi^*\omega^j+\int_{\Tube(D,r)}  P\wedge  \pi^*\omega^j=  \int_{\Tube(D,r)}T^{(j)} \wedge \pi^*\omega^j\\
 &=& \lim_{n\to\infty} \int_{(\Pi^\bullet_{\hat j})^{-1}(\Tube(D,0,r))}\Pi_{\hat j}^\diamond (\tau_* (T_n))\wedge  (\Pr\nolimits_{\hat j})^\diamond (\Upsilon_{\hat j}^{\dim{\X_{\hat j}}- p-j })\wedge (\Pi^*_{\hat j})(\pi^*(\omega^j))\\
 &=&   \lim_{n\to\infty} \int_{(\Pi_{\hat j})^{-1}(\Tube(D,r))} \widetilde T^{(j)}\wedge(\Pr\nolimits_{\hat j})^*(\Upsilon_{\hat j}^{\dim\X_{\hat j}- p-j})\wedge \Pi_{\hat j}^*(\pi^*(\omega^j))
 \\
 &=& \int_{(\Pi_{\hat j})^{-1}(D)}  \tilde f_j\wedge\Upsilon_{\hat j}^{\dim\X_{\hat j}- p-j}\wedge  \Pi_{\hat j}^*(\pi^*(\omega^j))
 +\int_{(\Pi_{\hat j})^{-1}(\Tube(D,r))}  Q\wedge(\Pr\nolimits_{\hat j})^*(\Upsilon_{\hat j}^{\dim\X_{\hat j}- p-j})\wedge \Pi_{\hat j}^*(\pi^*(\omega^j)).
\end{eqnarray*}
Observe  that
\begin{equation*}
 \int_{\Tube(D,r)}  P\wedge  \pi^*\omega^j\leq \|P\|  \big( \Tube(D,r) \big)
 \end{equation*}
 and 
 \begin{equation*}
 \int_{(\Pi_{\hat j})^{-1}(\Tube(D,r))}  Q\wedge(\Pr\nolimits_{\hat j})^*(\Upsilon_{\hat j}^{\dim\X_{\hat j}- p-j})\wedge \Pi_{\hat j}^*(\pi^*(\omega^j))\leq \|Q\|  \big( (\Pi_{\hat j})^{-1}(\Tube(D,r))\big).
 \end{equation*}
 Moreover, both  RHSs tend to $0$ as $r\to 0.$  Thus,
 \begin{equation*}
  \int_{D} f_j \wedge \pi^*\omega^j=\lim_{r\to 0}\int_{\Tube(D,r)}T^{(j)} \wedge  \pi^*\omega^j=
  \int_{(\Pi_{\hat j})^{-1}(D)}  \tilde f_j\wedge\Upsilon_{\hat j}^{\dim\X_{\hat j}- p-j}\wedge  \Pi_{\hat j}^*(\pi^*(\omega^j)).
 \end{equation*}
  This proves  assertion (1).   
 
To prove assertion (2) observe that when $j=l$ both $f_{\upm}$ and  $\tilde f_{\upm}$   are functions.
Applying    assertion (1) to  $D:=\B(x,r),$ the ball  with center $x$ and radius $r$ for all $r>0$ small  enough,  
we get the  desired identity.
 \endproof

\begin{definition}\rm
 \label{D:hor-constant-metric}
We   say that a metric  $h$ on $\E$ is  horizontally constant  if $\alpha=\alpha_\ver$ and  $\beta=\beta_\ver,$ in other words,  if
the  horizontal parts   $\alpha-\alpha_\ver$ (resp.  $\beta-\beta_\ver$) of $\alpha$  (resp. $\beta$) vanish simultaneously.
 \end{definition}

 \begin{theorem}\label{T:charac-closed-all-degrees}
   We  keep the    Standing Hypothesis. Suppose that $ \ddc\omega^j=0$ on $B$ for all  $1\leq j\leq   \upm-1.$
Suppose  that    the current   $T$  is  positive closed  and $T=T^+-T^-$ on an open  neighborhood of $\overline B$ in $X$ with $T^\pm$ in the class $\CL_p^{2,2}(B).$
Suppose in addition  that  the metric $h$ of $\E$  is horizontally constant and 
  there is  a holomorphic admissible map $\tau$ for $B.$  Then, 
  \begin{enumerate}
   
   \item  for  $\lowm\leq j\leq \upm,$ we have  
  $$    \nu_j(T,B,h) =\int_{B} f_j \wedge \omega^j=
  \int_{(\Pi_{\hat j})^{-1}(B)}  \tilde f_j\wedge\Upsilon_{\hat j}^{\dim\X_{\hat j}- p-j}\wedge (\Pi_{\hat j}^*)(\pi^*(\omega^j)),
  $$
  where the  currents $f_j$ and $\tilde f_j$ given by  Proposition  \ref{P:Cut-Off} (4) are positive plurisubharmonic
  and   $\hat j$ is  given by \eqref{e:hatj}.
  \item   If  moreover $\omega$ is  K\"ahler, then    the above assertion   still holds if 
    $T^\pm$ in the class $\CL_p^{1,1}(B).$
\end{enumerate}
 \end{theorem}
\proof
Let  $\lowm\leq j\leq \upm,$
By  Proposition \ref{P:existence}, there exists a sequence of smooth forms $(T_{N_n})_{n=1}^\infty$ such that
$$\lim_{n\to\infty} \big(\tau_*(T_{N_n})\wedge \alpha_\ver^{ k-p-j}\big)_\bullet=T^{(j)}\qquad \text{weakly on }\qquad \U. $$
Clearly, $T^{(j)}$ is  a  current of order $0.$
We  will check that  $T^{(j)}$ is  closed. Let $\Phi$ be  a  smooth  form compactly supported in $\U.$
Since $T_{N_n}$ is  closed, it follows  that $d (\tau_*T_{N_n})= \tau_* (dT_{N_n})=0.$ We  also  have $\alpha_\ver=\alpha$  because
the metric $h$  is constant.  Hence,
\begin{eqnarray*}
 d\big(\tau_*(T_{N_n})\wedge \alpha_\ver^{ k-p-j}\big)_\bullet(\Phi)&=& \lim_{r\to 0}  \int_{\U\setminus  \Tube(B,r)} \tau_*(T_{N_n})\wedge \alpha_\ver^{ k-p-j}
 \wedge d\Phi\\
 &=&\lim_{r\to 0}  \int_{\U\cap\partial \Tube(B,r)}\tau_*(T_{N_n})\wedge \alpha_\ver^{ k-p-j}\wedge \Phi  =0,
\end{eqnarray*}
where the last  equality  holds  by Lemma \ref{L:vertical-boundary-terms} (2) since  $k-p-j\leq k-l.$
Hence,  $T^{(j)}$ is  a closed current of order $0.$
So  it is  also $\C$-normal  and $f_j$ is   also a closed current on $B.$
Since $T_{N_n}$  are positive,  we  see that $f_j$ is also  a positive  current.
Similarly,  we also see that  $\tilde f_j$ is a  positive closed current  on  $\Pi_{\hat j}^{-1}(B).$

Applying Theorem  \ref{T:Lelong-Jensen-smooth} to $  \tau_*(T_{N_n})\wedge   \pi^*(\omega^j)$ yields that
\begin{equation*}
  \int_{\Tube(B,r)} \tau_*(T_{N_n})\wedge \alpha^{ k-p-j}\wedge    \pi^*(\omega^j)={1\over r^{2(k-p-j)}} \int_{\Tube(B,r)} \tau_*(T_{N_n})\wedge \beta^{ k-p-j}\wedge \pi^*(\omega^j)+\Vc(\tau_*(T_{N_n})\wedge   \pi^*(\omega^j),r).
\end{equation*}
By Theorem  \ref{T:vertical-boundary-terms} $\Vc(\tau_*(T_{N_n})\wedge   \pi^*(\omega^j),r)=O(r)$ as $r\to 0.$
This, combined with   Proposition \ref{P:Cut-Off} (7),  implies that
as $n$ tends to infinity
\begin{equation*}
 {1\over r^{2(k-p-j)}} \int_{\Tube(B,r)} \tau_*(T)\wedge \beta^{ k-p-j}\wedge \pi^*(\omega^j)  = \int_{B} f_j\wedge\pi^*(\omega^j)+   O(r).
\end{equation*}
Taking again $r\to 0,$ we obtain assertion (1).

Since such a  form $  \tau_*(T_{N_n})\wedge   \pi^*(\omega^j)$ is  $d$-closed,
assertion (2) can be proved  in the same way  as  using   Theorem  \ref{T:vertical-boundary-closed} instead of Theorem  \ref{T:vertical-boundary-terms}.
\endproof

 When  the metric $h$  of  the normal bundle $\E$ is  not constant, we only obtain a geometric characterization for the top-Lelong number.
 
 \begin{theorem}\label{T:charac-closed-top-hol}
   We  keep the    Standing Hypothesis. Suppose that $ \ddc\omega^j=0$ on $B$ for all  $1\leq j\leq   \upm-1.$
Suppose  that    the current   $T$  is  positive closed  and $T=T^+-T^-$  on an open  neighborhood of $\overline B$ in $X$ with $T^\pm$ in the class $\CL_p^{2}(B).$
Suppose in addition  that   
  there is  a holomorphic admissible map $\tau$ for $B.$  Then  one  and  only one of the  following assertion holds:
  \begin{enumerate}

   \item If $\upm=k-p,$ then  
   $\nu_\upm(T,B,\tau)$ is  simply the mass of the  measure $T\wedge \pi^*(\omega^\upm)$ on $B.$
   \item If   $\upm\not=k-p,$ then  $\upm=l$ and  the function  $\tilde f_\upm$  given by  Proposition  \ref{P:Cut-Off} (4) is constant
on fibers  of $ \Pi_p,$  that is,  we have $f_\upm\circ \Pi_p=\tilde  f_\upm,$  
and
   we have  
  $$    \nu_\upm(T,B,\tau) =\int_{B} f_\upm \omega^l. 
  $$
  \item   If  moreover $\omega$ is  K\"ahler, then    the above two assertions   still hold if 
    $T^\pm$ belong to the class $\CL_p^{1}(B).$
\end{enumerate}
 \end{theorem}
 \proof
 To prove the    assertion (1), 
 observe that  when $\upm=k-p,$  we  have
 $$
 \nu_\upm(T,B)= \lim_{r\to 0} \int_{\Tube(B,r)} \tau_*(T)\wedge\pi^*(\omega^{k-p})= (T\wedge  \tau^*\pi^*(\omega^{k-p}))(B).
 $$
 Since $\tau$ is  an  admissible map,  it follows from   Proposition \ref{P:basic-admissible-estimates-I} (2) that 
 $\tau^*(\pi^*(\omega^{k-p}))=\omega^{k-p}+ O(\|z\|)\omega^{k-p}+ O(1)  dz_j +O(1) d\bar z_j.$
 Consequently,    we  infer that
 $$
 \nu_\upm(T,B)= \lim_{r\to 0} \int_{\Tube(B,r)} \tau_*(T)\wedge\pi^*(\omega^{k-p})= (T\wedge  \omega^{k-p})(B)
 +O(r) \sum_{j=\lowm}^\upm \nu_j(T,B).
 $$
 Hence, assertion (1) follows.
 
 To prove     assertion (2), 
 observe that  when $\upm=l,$  we  have
 $$\beta^{k-p-l}\wedge (\pi^*\omega^l)=\beta_\ver^{k-p-l}\wedge (\pi^*\omega^l)\qquad\text{and}\qquad 
 \alpha^{k-p-l}\wedge (\pi^*\omega^l)=\alpha^{k-p-l}_\ver\wedge (\pi^*\omega^l).$$
 Using this,  we follow  along the  same lines as  those  given in the  proof of  Theorem  \ref{T:charac-closed-all-degrees}
 for $j:=\upm.$

 Assertion  (3) can be proved in the same way as we did  for assertions (1) and (2) using  that  such a  form $  \tau_*(T_{N_n})\wedge   \pi^*(\omega^j)$ is  $d$-closed.
 \endproof

\begin{theorem}\label{T:charac-ph-hol}
 We  keep the    Standing Hypothesis. Suppose that $\omega$ is  K\"ahler on $B.$ 
Suppose  that    the current   $T$  is  positive pluriharmonic  and $T=T^+-T^-$ on an open  neighborhood of $\overline B$ in $X$ with $T^\pm$ in the class $\PH_p^{2,2}(B).$
Suppose in addition  that  the metric $h$ of $\E$  is horizontally constant and 
  there is  a holomorphic admissible map $\tau$ for $B.$  Then, 
     for  $\lowm\leq j\leq \upm,$ we have  
  $$    \nu_j(T,B,h) =\int_{B} f_j \wedge \pi^*\omega^j=
  \int_{(\Pi_{\hat j})^{-1}(B)}  \tilde f_j\wedge\Upsilon_{\hat j}^{\dim\X_{\hat j}- p-j}\wedge (\Pi_{\hat j}^*)(\pi^*(\omega^j)),
  $$
  where the  currents $f_j$ and $\tilde f_j$ given by  Proposition  \ref{P:Cut-Off} (4) are positive plurisubharmonic
  and   $\hat j$ is  given by \eqref{e:hatj}.
\end{theorem}

 When  the metric $h$  of  the normal bundle $\E$ is  not constant, we only obtain a geometric characterization for the top-Lelong number.
 
 \begin{theorem}\label{T:charac-ph-top-hol}
   We  keep the    Standing Hypothesis.  Suppose that $\omega$ is  K\"ahler on $B.$ 
Suppose  that    the current   $T$  is  positive pluriharmonic  and $T=T^+-T^-$  on an open  neighborhood of $\overline B$ in $X$ with $T^\pm$ in the class $\PH_p^{2}(B).$
Suppose in addition  that   
  there is  a holomorphic  admissible map $\tau$ for $B.$  Then  one  and  only one of the  following assertion holds:
  \begin{enumerate}

   \item If $\upm=k-p,$ then  
   $\nu_\upm(T,B,\tau)$ is  simply the mass of the  measure $T\wedge \pi^*(\omega^\upm)$ on $B.$
   \item  If   $\upm\not=k-p,$ then  $\upm=l$ and  the function  $\tilde f_\upm$  given by  Proposition  \ref{P:Cut-Off} (4) is constant
on fibers  of $ \Pi_p,$  that is,  we have $f_\upm\circ \Pi_p=\tilde  f_\upm,$  
and
   we have  
  $$    \nu_\upm(T,B,\tau) =\int_{B} f_\upm \omega^l. $$
   
\end{enumerate}
 \end{theorem}
  
\proof
 Since  the proof is  not difficult, we leave it   to the interested reader.   
\endproof
 \begin{remark}
  \rm
  We  regard  the  above results from  a geometric viewpoint for the top case where $p\leq k-l$ and hence $\upm=l.$
  By   Theorem  \ref{T:charac-ph-top-hol}  (resp.  Theorem \ref{T:charac-closed-top-hol}), we have
  $$\lim_{n\to\infty}  \big ((\Pr\nolimits_p^\diamond) (\tau_*T_{N_n})\big)_\bullet=\widetilde T^{(\upm)}\qquad\text{weakly on}\qquad  \Pr\nolimits_p^{-1}(\U).$$
  Since  $T_n$ is  smooth,  $\big ((\Pr_p^\diamond) (\tau_*T_n)\big)_\bullet= (\Pr_p^*) (\tau_*T_n).$ Hence,  $\widetilde T^{(\upm)}$
  is positive pluriharmonic (resp.  positive closed).  Observe that
  $$
  \ind_{\Pi^{-1}_p(B)}\widetilde T^{(\upm)}=\tilde f_\upm [\Pi^{-1}_p(B)] \quad\text{and}\quad f_\upm\circ \Pi_p=\tilde  f_\upm,
  $$
  where  $\tilde f_\upm$  and $f_\upm$  are positive  weakly plurisubharmonic  functions. Moreover, both  functions are non-negative constant if $T$ is either 
  in $\CL^{2}(B)$  or in  $\PH^2(B,\comp).$
  Indeed,  positive weakly plurisubharmonic  function on a compact manifold is necessarily constant.

  So  we have
  $$
  \widetilde T^{(\upm)}=\tilde f_\upm [\Pi^{-1}_p(B)]+ \big( \widetilde T^{(\upm)}|_{\Pi^{-1}_p(\U\setminus V) }\big)_\bullet\quad\text{and}\quad
   \nu_\upm(T,B,\tau) =\langle  f_\upm [B], \omega^l \rangle.$$
  This  means that the positive  plurisubharmonic  function $f$ which gives the top Lelong number
of $T$  along $B$ is simply  the  density  of the mass of  $\widetilde T^{(\upm)}$  over $\Pi^{-1}_p(B).$ 
But
  $$\widetilde T^{(\upm)}|_{\Pi^{-1}_p(\U\setminus V) }=\lim_{n\to\infty}  (\Pr\nolimits_p^\diamond) (\tau_*T_{N_n}).
  $$
In other words, $\widetilde T^{(\upm)}$ is completely determined  by $T$ and $\tilde f_\upm,$  that is, it is  independent of the choice of the  approximating  sequence. Therefore,  we can define  $\Pi^*_p(T):=\widetilde T^{(\upm)}.$
\end{remark}

\begin{remark}\rm
Our results follow the model  of Siu \cite{Siu}. Indeed,  suppose that $V=\{x\}$ is  a  single point $T$ is  a positive closed  $(p,p)$-current defined
  on an open neighborhood $U$ of $x.$  Writing $\widetilde T$
  instead of $\widetilde T^{(0)},$
   Siu's result  and  our Theorem  \ref{T:charac-closed-top-hol}  say that
  $$
  \nu(T,x)=\|\widetilde T\|(\G_p(\C^k))\quad\text{and}\quad\widetilde T=\Pi^\diamond_p(T|_{U\setminus \{x\}})+\nu(T,x)[\G_p(\C^k)].
  $$
  In particular, $\widetilde T$ is independent of the choice of  approximating  forms $(T_n)_{n=1}^\infty$  for $T.$
  So
  we can define 
  $$\Pi_p^* T :=\widetilde T=\Pi^\diamond_p(T|_{U\setminus \{x\}})+\nu(T,x)[\G_p(\C^k)].$$
 \end{remark}

\section{Strongly admissible maps and geometric  characterizations for positive closed  and  positive  pluriharmonic currents}
\label{S:Charac-pos-closed-and-plurihar-II}

\subsection{Cut-off  along  $V$ on $\E$}
When the admissible map $\tau$ is  not necessarily  holomorphic, we have the following analogous result but only for the top degree
$j=\upm.$

 \begin{proposition}\label{P:Cut-Off}
 Suppose that $p<k-l.$
Let $T$  be a  current in the class $ \SH^{3,3}_p(B)$ (resp.  $ \PH^{2,2}_p(B)$, resp. $\CL^{1,1}(B)$) introduced in Definition  \ref{D:classes}
with   an approximating sequence of $(T_n)_{n=1}^\infty.$
Let $\tau $ be a   strongly  admissible map along $B.$ 
Consider the  real currents $R_n:=\tau_*(T_n)$ on $\U$  for $n\geq 1.$
Then  the following   assertions hold:

\begin{enumerate}
 
\item   The  sequences  $R_n\wedge (\hat \alpha')^{ k-p-l}$ and  $R_n\wedge \alpha_\ver^{ k-p-l}$    are relatively compact in the  weak-$\star$ topology on $\U.$

\item   There exist  currents   $\widehat R^{(p,k-p-l)}$ and $R^{(p,k-p-l)}$   on  $\U$   such that, for a suitable subsequence $(R_{N_n})$ of  the  sequence $\big(R_n\big)_{n=1}^\infty,$ the following properties hold:
  \begin{enumerate}
   \item $\lim_{n\to \infty} \big(R_{N_n}\wedge(\hat \alpha')^{k-p- l}\big)_\bullet = \widehat R^{(p,k-p-l)} $  weakly on $\U;$
   \item $\lim_{n\to \infty} \big(R_{N_n}\wedge \alpha_\ver^{k-p- l}\big)_\bullet = R^{(p,k-p-l)} $  weakly on $\U.$
    
  \end{enumerate}

\item  Set  $T^{(l)}:= R^{(p,k-p-l)}$  
and      $\widehat T^{(l)}:= \widehat R^{(p,k-p-l)}. $ 
Then   $T^{(l)} $ and  $\widehat T^{(l)}$  enjoys the cut-off property  through $V$ in $\U.$  v

\item For $1\leq \ell\leq\ell_0,$ set $\widehat R_{[\ell]}:=\widehat R^{(p,k-p-l)}_{[\ell]}$  and     $  R_{[\ell]}:=R^{(p,k-p-l)}_{[\ell]},$   
  where the current on the  RHS  is  defined by Corollary  \ref{C:Cut-off}.
Then 
  $$
\ind_{B}( T^{(l)})=  \ind_{B}(  \widehat R_{[\ell]})= \ind_{B}( R_{[\ell]})   . 
$$ 
\item 
There exists   a positive plurisubharmonic  function $f$   on $B$    such that 
$
\ind_{B}( T^{(l)})=(\iota_{B,\E})_*(f) . 
$
Moreover,  if   $T$ belongs to $ \CL^{1,1}_p(B),$ then  $f$ is a non-negative constant.
\end{enumerate}
\end{proposition}
\proof \noindent
{\bf Proof of assertion (1).} Since $(\U_\ell)_{\ell=1}^{\ell_0}$ is an open cover of $\U,$ we  only need to show that  the sequences   $R_n\wedge (\hat \alpha')^{ k-p-l}$ and $R_n\wedge \alpha_\ver^{ k-p-l}$  are  relatively compact in the  weak-$\star$ topology on $\U_\ell$
for a given $\ell.$  We compose $\tilde\tau_\ell^*$ to both sequences and note that $(\tilde\tau_\ell)_*\circ (\tau_\ell)_*(T_n)=\tau_*T_n=R_n,$ we are reduced to proving the relative  compactness of 
 the sequences   $ (\tau_\ell)_*(T_n)\wedge (\tilde\tau^*_\ell (\hat \alpha'))^{ k-p-l}$ and $ (\tau_\ell)_*(T_n)\wedge (\tilde\tau^*_\ell(\alpha_\ver))^{ k-p-l}.$ 
 
 \begin{lemma}\label{L:difference-alpha}
For every $1\leq \ell\leq \ell_0,$  the forms
 $(\tilde\tau^*_\ell (\hat \alpha'))^{ k-p-l}- (\hat \alpha')^{ k-p-l}$   and 
$(\tilde\tau^*_\ell (\alpha_\ver))^{ k-p-l}-  \alpha_\ver^{ k-p-l}$) can be rewritten  as a finite sum
$\sum_N  f_{[\ell],N}  R_{[\ell],N},$  where   
the $f_{[\ell],N}$'s  are continuous   forms and   the $R_{[\ell],N}$'s are some real $(m,m)$-forms  with $0\leq m\leq  k-p-l$ such that
$$
\pm R_{[\ell],N}\lesssim  \sum_{q=0}^{\min(m,k-p-l-1)}(\pi^*\omega+\hat\beta)^{m-q} \wedge (\hat \alpha')^{q}.
$$
\end{lemma}
\proof  It follows from Theorem \ref{T:basic-admissible-estimates}.
\endproof

Using Lemma \ref{L:difference-alpha}, we  see that the  desired  compactness  will follow if one can show that 
\begin{equation}\label{e:rel-compact-R_n}
\sup_{n\in\N}\int_{\U}  (T_n)^\hash\wedge \hat\alpha^q\wedge \pi^*(\omega^{q'})\wedge \hat \beta^{k-p-q-q'}<\infty
\end{equation}
for every $0\leq q\leq k-l$ and $0\leq q'\leq k-p-q.$  But the last inequality holds by  
arguing  as in the  proof of Proposition \ref{P:existence-finite}.

 \noindent{\bf Proof of assertion (2).}
It is  an immediate consequence of assertion (1).

 \noindent{\bf Proof of assertion (3).} Since the  problem is local and $\U:=\bigcup_{\ell=1}^{\ell_0}\U_\ell,$ we will prove  the assertion on each $\U_\ell.$
By Lemma  \ref{L:Cut-Off},  we need  to prove the cut-off property for $(\tilde\tau_\ell)^*(T^{(l)})$  
and   $(\tilde\tau_\ell)^*(\widehat T^{(l)}). $ By  assertion (2), we have
 \begin{enumerate}
   \item $\lim_{n\to \infty} \big((\tau_\ell)_*T_{N_n}\wedge(\tilde\tau_\ell)^*(\hat \alpha'^{k-p- l})\big)_\bullet =(\tilde\tau_\ell)_*(\widehat T^{(l)}) $  weakly on $\U_\ell;$
   \item $\lim_{n\to \infty} \big((\tau_\ell)_*R_{N_n}\wedge (\tilde\tau_\ell)^*(\alpha_\ver^{k-p- l})\big)_\bullet = (\tilde\tau_\ell)^*(\widetilde T^{(l)})$  weakly on $\U_\ell.$
  \end{enumerate}
By   \eqref{e:rel-compact-R_n}  we see  that $ (\tau_\ell)_*(T_n)\wedge (\hat\alpha')^q\wedge \pi^*(\omega^{q'})\wedge (c_1\pi^*\omega+ \beta)^{q''}$  are positive plurisubharmonic  currents of uniform  bounded mass. Hence, by Theorem \ref{T:Bassanelli-psh} 
$$
R_{q,q',q''}:=\lim_{n\to\infty}(\tau_\ell)_*(T_{N_n})\wedge (\hat\alpha')^q\wedge \pi^*(\omega^{q'})\wedge (c_1\pi^*\omega+ \beta)^{q''}
$$ possess the cut-off property. 
On  the  other hand, arguing as in the proof of assertion (1), we see that the measure coefficient of  the  three currents $(\tilde\tau_\ell)_*(T^{(l)}), $  
and    $(\tilde\tau_\ell)_*(\widetilde T^{(l)})$  and $(\tilde\tau_\ell)_*(\widehat T^{(l)})$ are dominated  by a combination (with smooth functions)
of the  coefficients of the currents $R_{q,q',q''}.$  Hence, the  former currents inherit the cut-off property from the latter ones.

\noindent{\bf Proof of assertion (4).} Fix an $\ell$ with $1\leq \ell\leq \ell_0.$
By  assertion (3), there are real  functions $f,$ $ f_{[\ell]},$ $\hat f_{[\ell]}$  defined on $\U_\ell,$ 
such that 
\begin{equation*}
\ind_{B}( T^{(l)})=(\iota_{B,\E})_*(f)\quad\text{and}\quad 
\ind_{B}(  R_{[\ell]})=(\iota_{B,\E})_*( f_{[\ell]})\quad\text{and}\quad 
\ind_{B}( \widehat R_{[\ell]})=(\iota_{B,\E})_*(\hat f_{[\ell]}).
\end{equation*} 
We  need to show that
\begin{equation}\label{e:two-equalities-on-fs}
 f=\hat f_{[\ell]}=f_{[\ell]}.
\end{equation}
Let  $g$ be  a real smooth  test  function in $\Cc^\infty_0(\U_\ell).$
By  assertion (2), we get that
 \begin{eqnarray*}
 \int_{\U_\ell} (\hat f_{[\ell]} g)\cdot \omega^l&=&  \langle  \ind_{B}( \widehat R_{[\ell]}),g\omega^l\rangle=\lim\limits_{t\to 0}\Big(\lim_{n\to \infty} \int_{\Tube(B,t)}  (\tau_\ell)_*(T_{N_n})\wedge   g\cdot(\hat \alpha')^{k-p- l}\wedge  \pi^*(\omega^l) \Big),\\
 \int_{\U_\ell} ( f_{[\ell]} g)\cdot \omega^l&=&  \langle  \ind_{B}( R_{[\ell]}),g\omega^l\rangle=\lim\limits_{t\to 0}\Big(\lim_{n\to \infty} \int_{\Tube(B,t)}   (\tau_\ell)_*(T_{N_n})\wedge    g\cdot\alpha_\ver^{k-p- l}\wedge  \pi^*(\omega^l) \Big).
\end{eqnarray*}
On the  one hand, we see easily that
$$ \alpha_\ver^{k-p- l}\wedge  \pi^*(\omega^l)= (\hat \alpha')^{k-p- l}\wedge  \pi^*(\omega^l)=\alpha^{k-p- l}\wedge  \pi^*(\omega^l),
$$
Consequently,  we obtain the identity $f_{[\ell]}=\hat f_{[\ell]}.$

Next, by  assertion (2), we get that
 \begin{eqnarray*}  \int_{\U_\ell} (fg)\cdot \omega^l&=&  \langle  \ind_{B}( T^{(l)}),g\omega^l\rangle=\lim\limits_{t\to 0}\Big(\lim_{n\to \infty} \int_{\Tube(B,t)}   g\cdot \tau_*(T_{N_n})\wedge \alpha_\ver^{k-p- l}\wedge  \pi^*(\omega^l) \Big),\\
 &=&  \lim\limits_{t\to 0}\Big(\lim_{n\to \infty} \int_{\Tube(B,t)}   \tau_*(T_{N_n})\wedge  g\cdot \alpha^{k-p- l}\wedge  \pi^*(\omega^l) \Big).
 \end{eqnarray*}
 By  Lemma \ref{L:basic-positive-difference-bis},  we have
 \begin{eqnarray*}
   \int_{\Tube(B,t)}   \tau_*(T_{N_n})\wedge  g\cdot \alpha^{k-p- l}\wedge  \pi^*(\omega^l)-  \left\langle (T_{N_n})^\hash_t,  g\cdot \alpha^{k-p- l}\wedge  \pi^*(\omega^l)\right\rangle\leq \sum.
 \end{eqnarray*}
 Applying Proposition \ref{P:Lc-finite-psh} and Lemma \ref{L:mass-tends-to-zero} to $T_{N_n}$, the  sum $\sum$ on the RHS is
 uniformly of order $O(t)$ independently of $n.$
 Hence,  taking  $t\to 0,$ and using formula  $(T_{N_n})^\hash_t$  given in \eqref{e:T-hash_r}, we get that
 $$
 \int_{\U_\ell} (fg)\cdot \omega^l-\sum_{\ell=1}^{\ell_0} (\pi^*\theta_\ell) \cdot (\ind_{\Tube(B,t)\circ \tilde\tau_\ell})\int_{\U_\ell} (\hat f_{[\ell]} g)\cdot \omega^l
 =0. $$
  Since   this  is true for all test functions $g,$  we obtain  
  $$f=\sum_{\ell=1}^{\ell_0}\theta_\ell \hat f_{[\ell]}.$$
  In the remainder  of the  proof, we will show that 
  \begin{equation}  \label{e:equality-on-fs} \hat f_{[\ell]}= \hat f_{[\ell']}\qquad\text{ on}\qquad B\cap \U_\ell\cap\U_{\ell'}.
  \end{equation}
  Taking  for granted  this equality for the moment, we infer from  the previous equality  and  the  identity $\sum _{\ell=1}^{\ell_0}\theta_\ell =1$
  that   $f=\hat f_{[\ell]}$ on $\U_\ell$  and    equalities \eqref{e:two-equalities-on-fs}  follow.

  To  finish the proof of assertion (4), it remains  to establish \eqref{e:equality-on-fs}. Let $g$ be a test function in  the class 
  $\Cc^\infty_0(\U_\ell\cap \U_{\ell'}).$ Write
  \begin{eqnarray*}
 \int_{\U_\ell} (\hat f_{[\ell]} g)\cdot \omega^l&=&  \langle  \ind_{B}( \widehat R_{[\ell]}),g\omega^l\rangle=\lim\limits_{t\to 0}\Big(\lim_{n\to \infty} \int_{\Tube(B,t)}  (\tau_\ell)_*(T_{N_n})\wedge   g\cdot \alpha^{k-p- l}\wedge  \pi^*(\omega^l) \Big),\\
 &=&  \lim\limits_{t\to 0}\Big(\lim_{n\to \infty} \int_{\sigma^{-1}(\Tube(B,t))}   (\tau_{\ell'})_*(T_{N_n})\wedge  \sigma^* \big( g\cdot\alpha^{k-p- l}\wedge  \pi^*(\omega^l) \big)\Big),
\end{eqnarray*}
where $\sigma:= \tau_\ell\circ\tau_{\ell'}^{-1}.$
Observe  that $$ \sigma^* \big( g\cdot\alpha^{k-p- l}\wedge  \pi^*(\omega^l) \big) -  \big( g\cdot\alpha^{k-p- l}\wedge  \pi^*(\omega^l) \big)$$ is $(2j-1)$-negligible. Hence,  by Proposition \ref{P:mass-fine-negligible} and Proposition \ref{P:Lc-finite-psh}, we get the estimate independently of $n:$
\begin{eqnarray*}
 \int_{\sigma^{-1}(\Tube(B,t))}   (\tau_{\ell'})_*(T_{N_n})\wedge  \sigma^* \big( g\cdot\alpha^{k-p- l}\wedge  \pi^*(\omega^l) \big)
- \int_{\Tube(B,t)}   (\tau_{\ell'})_*(T_{N_n})\wedge   \big( g\cdot\alpha^{k-p- l}\wedge  \pi^*(\omega^l) \big)=O(t).
\end{eqnarray*}
Therefore,  it follows that 
$$\int_{U_\ell\cap \U_{\ell'}} (\hat f_{[\ell]} g)\cdot \omega^l=\int_{U_\ell\cap \U_{\ell'}} (\hat f_{[\ell']} g)\cdot \omega^l.$$
This proves  \eqref{e:equality-on-fs}.
\endproof

\begin{lemma}\label{L:mass-tends-to-zero}
 We  keep the hypothesis of 
 Proposition \ref{P:Cut-Off}. Then for every $\lowm\leq j\leq\upm$ and $0\leq q< k-p-j,$
 \begin{equation*}
  \sup_{r\in(0,\bfr]}{1\over r^{2(k-p-j-q)} }\Big(\sup_{T\in\Mc} \int_{\Tube(B,r)}  (T^\hash)\wedge  \hat\alpha^q\wedge \pi^*(\omega^j)\wedge \hat\beta^{k-p-j-q}\Big) <\infty.
 \end{equation*}
 Here $\Mc$ stands for one  of the following  classes $\widetilde \SH^{3,3}_p(\bfU,\bfW),$ $\widetilde \PH^{2,2}_p(\bfU,\bfW),$
 and $\widetilde \CL^{1,1}_p(\bfU,\bfW).$
\end{lemma}
\proof
We only give the proof for the case where $\Mc:=\widetilde \SH^{3,3}_p(\bfU,\bfW).$ The remaining two cases can be treated similarly.
  By Proposition \ref{P:Lc-finite-psh}, $\Mc_{j+q}(T,r)<c_{11}.$
This, combined  with inequality \eqref{e:hat-alpha-vs-hat-beta} $\varphi\hat\alpha\leq c_3\hat\beta,$ implies the desired conclusion.
\endproof

\subsection{Effect of  strongly admissible maps on the cut-off  along  the  exceptional fiber}
\label{SS:Effect-strongly-adm-maps}

Throughout the  subsection we  always assume that  $1\leq p<k-l.$
Consider the projection  $\Pi_p:\  \X_p\to \E.$

Recall the homogeneous  coordinates introduced in  \eqref{e:homogeneous-coordinates}.
We place ourselves on  an open set of $\C^{k-l}$ defined by $z_{k-l}\not=0.$
We   may assume without loss of generality as in \eqref{e:max-coordinate} that
\begin{equation*} 2|z_{k-l}| > \max\limits_{1\leq j\leq k-l}|z_j|.
\end{equation*}
and use the projective coordinates
\begin{equation*} 
\zeta_1:={z_1\over z_{k-l}},\ldots, \zeta_{k-l-1}:={z_{k-l-1}\over z_{k-l}},\quad \zeta_{k-l}=z_{k-l}=t.
\end{equation*}
In the coordinates  $\zeta=(\zeta_1,\ldots,\zeta_{k-l})=(\zeta',\zeta_{k-l})=(\zeta',t),$ the form $\omega_\FS([z])$  can be  rewritten as  $$\ddc \log{ (1+|\zeta_1|^2+\cdots+|\zeta_{k-l-1}|^2)},$$
and a direct computation shows that 
$$
\omega_\FS([z])\approx  (1+\|\zeta'\|^2)^{-2}\omega'(\zeta'),\quad\text{where}\quad \omega'(\zeta'):=\ddc (|\zeta_1|^2+\cdots+|\zeta_{k-l-1}|^2).
$$
Throughout the  subsection we  always assume that  $1\leq p<k-l.$
Consider the projection  $\Pi_p:\  \X_p\to \E.$

Fix  a point $x\in V.$ We add to   the  coordinates $z=(z_1,\ldots,z_{k-l})$ the  coordinates $w=(w_1,\ldots,w_k)$ so that
$(z,w)$ is a local coordinate  around $x.$
Let  $H$ be an element of $\G_p(\E_x).$ Then  $H$ is a $p$-linear subspace of $E_x.$
We may assume without loss of generality that  $H_0:=H\cap\{z_1=0\}$ is  a linear subspace of dimension $p-1.$
So $H_0$ defines an element in  $\G_{p-1}(\E_x).$
We may assume without loss of generality that
$$
H_0:=  \left\lbrace z_1=\cdots= z_{k-l-p+1}=0\right\rbrace.
$$
For $z=(z_1,\ldots,z_{k-l}),$  write  $z^{('p)}=(z_1,\ldots,z_{k-l-p+1})\in\C^{k-l-p+1}.$
If  $z^{('p)}\not=0,$ let  $[z^{('p)}]$ be the image of $z^{('p)}$ by the  canonical projection $\C^{k-l-p+1}\setminus\{0\}\to \P^{k-l-p}.$
Consider \begin{equation}\label{e:X_p,H_0} \X_{p,H_0}:=  \{(z,H)\in \X_p:\  H_0\subset  H \}  \end{equation} and let 
$
\Pi_{p,H_0}$ be the restriction of $\Pi_p$ on    $ \X_{p,H_0}.$ Observe that    $H$ defines an element $[z^{('p)}]\in \P^{k-l-p}.$
We see that       $\X_{p,H_0}$ is the closure of  $\X'_{p,H_0}$ in $\C^{k-l-p+1}\times \P^{k-l-p},$ where  
\begin{equation} \label{e:blow-up-H_0}
 \X'_{p,H_0}\simeq  \{(z^{('p)},[z^{('p)}]):\  z^{('p)}\in\C^{k-l-p+1}\setminus\{0\}\}\quad\text{and}\quad
\Pi_{p,H_0}(z,H)=z^{('p)}.
\end{equation}
Consequently, we obtain  the model of blow-up  at the origin  in $\C^{k-l-p+1}.$
We place ourselves on the chart $\{\zeta^{('p)}\in \D^{k-l-p+1}:\ 2|\zeta_1|\geq  |\zeta_j|\quad\text{for}\quad 1< j\leq k-l-p\}.$
On this chart,  $\Pi_{p,H_0}(z,H)$ reads as $(\zeta_1,\zeta_1\zeta_2, \ldots,\zeta_1\zeta_{k-l-p+1}).$
\begin{lemma}\label{L:pushforward_Pi_p}
  Let $S$ be a continuous  real form of bidimension $(m,n)$ on $\X_p$  where $0\leq m,n\leq \dim  \X_p.$ 
  
  \begin{enumerate}
   \item  $(\Pi_p)_\diamond S$ is  a continuous  form of bidimension  $(m,n)$  with respect to the  homogeneous coordinates $(\zeta',\zeta_{k-l},w)$ on $\E\setminus V$
   and the bidegree  of $(d\zeta',d\overline\zeta')$ of each component of $(\Pi_p)_\diamond S$ is  $\leq (k-p-l,k-p-l).$
   
   \item  When $m=n=p,$  there is a constant $c=c_S>0$ such that  
 $$\pm (\Pi_p)_\diamond S\leq  c \sum_{0\leq j\leq l,\  0\leq q\leq k-l-p}  \pi^*(\omega^j)\wedge \alpha_\ver^q\wedge \beta_\ver^{k-p-j-q}.$$
  \end{enumerate}

\end{lemma}
\proof
Since  the proof is  not difficult, we leave it   to the interested reader.   
\endproof

\begin{lemma}\label{L:pushforward_Pi_p-bis}
 Let $S$ be a continuous  real form of dimension $2p$ on $\X_p. $ 
 Then  there is  a constant $c=c_S>0$  such that   for $0<r\leq \bfr,$ the  following  inequality holds
 
  \begin{equation*} \pm H\leq  cr \sum_{0\leq j\leq l,\  0\leq q\leq k-l-p}  \pi^*(\omega^j)\wedge \alpha_\ver^q\wedge \beta_\ver^{k-p-j-q}\qquad\text{on}\qquad  \Tube(B,r),\\
  \end{equation*}
 where $H$ is  either  $R$ or $R'$ with
 \begin{equation*} R:= [\tilde\tau_\ell^* ((\Pi_p)_\diamond (S-S^\sharp))]^\sharp\qquad\text{and}\qquad
   R':= [(\tilde\tau_\ell^*-\id) ((\Pi_p)_\diamond (S))]^\sharp
 \end{equation*}

\end{lemma}
\proof
 Since  the proof is  not difficult, we leave it   to the interested reader.   
\endproof

 \begin{proposition}\label{P:Cut-Off-bis}
 Suppose that $1\leq p<k-l.$
Let $T$  be a  current in the class $ \SH^{3,3}_p(B)$ (resp.  $ \PH^{2,2}_p(B)$, resp. $\CL^{1,1}(B)$) introduced in Definition  \ref{D:classes}
with   an approximating sequence of $(T_n)_{n=1}^\infty.$
Let $\tau $ be a   strongly  admissible map along $B.$ 
Consider the  real currents $R_n:=\tau_*(T_n)$ on $\U$  for $n\geq 1.$
Then  the following   assertions hold:

\begin{enumerate}
 
\item   The  sequence   $(\Pi_{p})^\diamond(R_n)$ is relatively compact in the  weak-$\star$ topology.

\item   There exist  currents   $\widetilde R^{(p)}$ on $\Pi_{p}^{-1}(\U)$ such that, for a suitable subsequence $(R_{N_n})$ of  the  sequence $\big(R_n\big)_{n=1}^\infty,$   $\lim_{n\to \infty} \big(\Pi_{p}^\diamond R_{N_n}\big)_\bullet =:\widetilde T^{(l)}$  weakly on $\Pi_{p}^{-1}(\U).$

\item  For $1\leq \ell\leq\ell_0,$ set  
$\widetilde R_{[\ell]}:=\widetilde R^{(p)}_{[\ell]},$
  where the current on the  RHS  is  defined by Corollary  \ref{C:Cut-off}.
Then 
  $   \ind_{\Pi_{p}^{-1}(B)}( \widetilde T^{(l)})= 
\ind_{\Pi_{p}^{-1}(B)}( \widetilde R_{[\ell]}). 
$

\item   
$ \widetilde T^{(l)}$  enjoys the cut-off property  through $\Pi_{p}^{-1}(V)$ in $\Pi_{p}^{-1}(\U).$

\item 
There exist   a  positive  plurisubharmonic function  $\tilde f$ on $\Pi_{p}^{-1}(B)$  such that 
$   \ind_{\Pi_{p}^{-1}(B)}( \widetilde T^{(l)})= 
(  \iota_{\Pi_{p}^{-1}(B),\X_{\hat j}})_*  (\tilde f). 
$ 
Moreover,  if   $T$ belongs to $ \CL^{1,1}_p(B),$ then   $\tilde f$ is a  non-negative constant.
\end{enumerate}
\end{proposition}
\proof \noindent{\bf Proof  of assertion (1).}
Pick $1\leq \ell\leq\ell_0.$
Let  $S$ be a continuous  test form of dimension $2p$ on $\X_p$ which is compactly supported on $\Pi_p^{-1}(\U_\ell).$
Write
\begin{eqnarray*}
\langle (\Pi_{p})^\diamond(R_n), S\rangle&=&\langle\tau_*(T_n), (\Pi_{p})_\diamond(S)\rangle=\langle (\tau_\ell)_*( T_n), \tilde\tau_\ell^*((\Pi_{p})_\diamond(S))\rangle\\
&=& \langle (\tau_\ell)_*(T_n), (\Pi_{p})_\diamond(S^\sharp)\rangle+\langle (\tau_\ell)_*(T_n), (\tilde\tau^*_\ell-\id)((\Pi_{p})_\diamond(S^\sharp))\rangle\\
&+&\langle (\tau_\ell)_*( T_n), \tilde\tau^*_\ell((\Pi_{p})_\diamond(S-S^\sharp))\rangle.
\end{eqnarray*}
Applying Lemma 
\ref{L:pushforward_Pi_p} to  $S^\sharp$ yields that
 $$ |\langle (\tau_\ell)_*(T_n), (\Pi_{p})_\diamond(S^\sharp)\rangle|\leq  c \sum_{0\leq j\leq l,\  0\leq q\leq k-l-p} \int  (\tau_\ell)_*(T_n) \wedge  \pi^*(\omega^j)\wedge \alpha_\ver^q\wedge \beta_\ver^{k-p-j-q}.$$
By Proposition  \ref{P:existence-finite} and  \ref{P:existence}, the  RHS is  uniformly bounded independent of $ n.$
So is  $ |\langle (\tau_\ell)_*(T_n), (\Pi_{p})_\diamond(S^\sharp)\rangle|.$

Set $R':= (\tilde\tau^*_\ell-\id)((\Pi_{p})_\diamond(S^\sharp)).$
Applying Lemma \ref{L:tilde-tau-dzeta_j} and Lemma \ref{L:difference} yields that
\begin{equation*}
 \langle (\tau_\ell)_*(T_n), R'\rangle= \langle (\tau_\ell)_*(T_n), (R')^\sharp\rangle \leq  c \sum_{0\leq j\leq l,\  0\leq q\leq k-l-p} \int  (\tau_\ell)_*(T_n) \wedge  \pi^*(\omega^j)\wedge \alpha_\ver^q\wedge \beta_\ver^{k-p-j-q}
\end{equation*}
By Proposition  \ref{P:existence-finite} and  \ref{P:existence}, the  RHS is  uniformly bounded independent of $ n.$
So is  $\big|\langle (\tau_\ell)_*(T_n), (\tilde\tau^*_\ell-\id)((\Pi_{p})_\diamond(S^\sharp))\rangle\big|.$

Set $R:=[\tilde\tau^*_\ell((\Pi_{p})_\diamond(S-S^\sharp))]^\sharp.$
Applying Lemma \ref{L:pushforward_Pi_p-bis}  yields that
\begin{equation*}
 \langle (\tau_\ell)_*( T_n), \tilde\tau^*_\ell((\Pi_{p})_\diamond(S-S^\sharp))\rangle = \langle (\tau_\ell)_*(T_n), R\rangle \leq  c r\sum_{0\leq j\leq l,\  0\leq q\leq k-l-p} \int  (\tau_\ell)_*(T_n) \wedge  \pi^*(\omega^j)\wedge \alpha_\ver^q\wedge \beta_\ver^{k-p-j-q}
\end{equation*}
By Proposition  \ref{P:existence-finite} and  \ref{P:existence}, the  RHS is  uniformly bounded by $cr$  independent of $ n.$
So is $\big|\langle (\tau_\ell)_*( T_n), \tilde\tau^*_\ell((\Pi_{p})_\diamond(S-S^\sharp))\rangle\big|.$

Putting together the above  three estimates, we get  $| \langle (\Pi_{p})^\diamond(R_n), S\rangle|\leq  c\|S\|_{\Cc^0},$
for a constant $c>0$ independent of $S.$  This proves assertion (1).

\noindent{\bf Proof  of assertion (2).} It is  an immediate consequence of assertin (1).


\noindent{\bf Proof  of assertion (3).} Fix $1\leq \ell\leq\ell_0.$
Let  $S$ be a continuous  test form of dimension $2p$ on $\X_p$ which is compactly supported on $\Pi_p^{-1}(\U_\ell).$
Let $0<r\leq\bfr.$
Write
\begin{eqnarray*}
\langle \widetilde T^{(l)}-  \widetilde R_{[\ell]} ,S \rangle_{\Tube(B,r)}&=&\lim_{n\to\infty} \langle (\Pi_{p})^\diamond(\tau_*T_n  -(\tau_\ell)_* T_n), S\rangle_{\Tube(B,r)}\\
&=&  \langle (\tau_\ell)_*(T_n), (\tilde\tau^*_\ell-\id)((\Pi_{p})_\diamond(S))\rangle_{\Tube(B,r)}.
\end{eqnarray*}
Applying Lemma \ref{L:pushforward_Pi_p-bis} to the  expression on the RHS yields that
\begin{equation*}
 \langle \widetilde T^{(l)}-  \widetilde R_{[\ell]} ,S \rangle_{\Tube(B,r)} \leq  c r\sum_{0\leq j\leq l,\  0\leq q\leq k-l-p} \int  (\tau_\ell)_*(T_n) \wedge  \pi^*(\omega^j)\wedge \alpha_\ver^q\wedge \beta_\ver^{k-p-j-q}
\end{equation*}
By Proposition  \ref{P:existence-finite} and  \ref{P:existence}, the  RHS is  uniformly bounded by $cr$  independent of $ n.$
So   $\lim_{r\to 0} \langle \widetilde T^{(l)}-  \widetilde R_{[\ell]} ,S \rangle_{\Tube(B,r)}=0.$  This proves assertion (3).

\noindent{\bf Proof  of assertion (4).}  It follows  from combining   assertion (3) and Proposition \ref{P:Cut-Off-hol} (2).

\noindent{\bf Proof  of assertion (5).}  It follows  from combining   assertion (4) and Proposition \ref{P:Cut-Off-hol} (4).
\endproof

\subsection{Geometric characterizations}

 \begin{proposition}\label{P:Cut-Off-rep}
 We keep the assumption, notation and  conclusion of Proposition  \ref{P:Cut-Off-bis}. Suppose   that $p<k-l.$
 Then      
$f_{l}$ and  $\tilde f_{l}$   are functions related by
$$
f_{l}(x)=\int_{\Pi_p^{-1}(0_x)} \tilde f_l \Upsilon_p^{p(k-l-p)}
$$
for Lebesgue almost every $x\in B.$
 \end{proposition}
 
 \proof

 By assertions (5) of Propositions \ref{P:Cut-Off} and \ref{P:Cut-Off-bis}, we can write
 \begin{equation*}
 T^{(l)}=(\iota_{V,\E})_*(f_l)+P\qquad\text{and}\qquad   \widetilde T^{(l)}= ( \iota_{\Pi_{p}^{-1}(V),\X_p  })_*(\tilde f_l)+Q, 
 \end{equation*}
where $P$ and $Q$ are positive  currents whose masses vanish  on  $V$ and $\Pi_{\hat j}^{-1}(V)$ respectively.

Let $D\Subset B\cap \bfU_\ell$ be a  subdomain  for some $1\leq \ell\leq \ell_0.$ 
By Proposition \ref{P:Cut-Off} (4),
\begin{equation*}
 \int_{D} f_l\wedge \pi^*\omega^j  =\lim_{r\to 0} \int_{\Tube(D,r)} R_{[\ell]}  \wedge \pi^*\omega^l.
 \end{equation*}
By Lemma \ref{L:Upsilon_j}, we have  that
\begin{eqnarray*}
 &&\int_{\Tube(D,r)} R_{[\ell]}  \wedge \pi^*\omega^l\\
 &=& \lim_{n\to\infty} \int_{(\Pi^\bullet_{p})^{-1}(\Tube(D,0,r))}\Pi_{p}^\diamond ((\tau_\ell)_* (T_n))\wedge  (\Pr\nolimits_p)^\diamond (\Upsilon_p^{\dim{\X_p}- p-l })\wedge (\Pi^*_p)(\pi^*(\omega^l))\\
 &=&   \lim_{n\to\infty} \int_{(\Pi_p)^{-1}(\Tube(D,r))} \widetilde R_{[\ell]}\wedge(\Pr\nolimits_p)^*(\Upsilon_p^{\dim\X_p- p-l})\wedge \Pi_p^*(\pi^*(\omega^l))
 \\
 &=& \int_{(\Pi_p)^{-1}(D)}  \tilde f_l\wedge\Upsilon_p^{\dim\X_p- p-l}\wedge  \Pi_p^*(\pi^*(\omega^l))
 +\int_{(\Pi_p)^{-1}(\Tube(D,r))}  Q\wedge(\Pr\nolimits_p)^*(\Upsilon_p^{\dim\X_p- p-l})\wedge \Pi_p^*(\pi^*(\omega^l)),
\end{eqnarray*}
where for the last equality we apply  Proposition \ref{P:Cut-Off-bis} (3).
 Since 
 \begin{equation*}
 \int_{(\Pi_p)^{-1}(\Tube(D,r))}  Q\wedge(\Pr\nolimits_p)^*(\Upsilon_p^{\dim\X_p- p-l})\wedge \Pi_p^*(\pi^*(\omega^l))\leq \|Q\|  \big( (\Pi_p)^{-1}(\Tube(D,r))\big)\to 0 \quad\text{as}\quad r\to 0,
 \end{equation*}
 we infer that
 \begin{equation*}
  \int_{D} f_l \wedge \pi^*\omega^l=\lim_{r\to 0} \int_{\Tube(D,r)} R_{[\ell]}  \wedge \pi^*\omega^l=
  \int_{(\Pi_p)^{-1}(D)}  \tilde f_p\wedge\Upsilon_p^{\dim\X_p- p-l}\wedge  \Pi_p^*(\pi^*(\omega^l)).
 \end{equation*}
 Hence, 
 \begin{equation*}
  \int_{D} f_l \wedge \pi^*\omega^l= 
  \int_{(\Pi_p)^{-1}(D)}  \tilde f_p\wedge\Upsilon_p^{\dim\X_p- p-l}\wedge  \Pi_p^*(\pi^*(\omega^l)).
 \end{equation*}
 By a routine partition of unity we can show that the above  equality holds for arbitrary domains $D\Subset B.$
 Finally,
applying    this  equality to  $D:=\B(x,r),$ the ball  with center $x$ and radius $r$ for all $r>0$ small  enough,  
we get the  desired identity.
 \endproof

 \begin{theorem}\label{T:charac-closed-top}
   We  keep the    Standing Hypothesis. Suppose that $ \ddc\omega^j=0$ on $B$ for all  $1\leq j\leq   \upm-1.$
Suppose  that    the current   $T$  is  positive closed  and $T=T^+-T^-$  on an open  neighborhood of $\overline B$ in $X$ with $T^\pm$ in the class $\CL_p^{2}(B).$
Suppose in addition  that   
  there is  a strongly  admissible map $\tau$ for $B.$  Then  one  and  only one of the  following assertion holds:
  \begin{enumerate}

   \item If $\upm=k-p,$ then  
   $\nu_\upm(T,B,\tau)$ is  simply the mass of the  measure $T\wedge \pi^*(\omega^\upm)$ on $B.$
   \item If   $\upm\not=k-p,$ then  $\upm=l$ and  the function  $\tilde f_\upm$  given by  Proposition  \ref{P:Cut-Off} (4) is constant
on fibers  of $ \Pi_p,$  that is,  we have $f_\upm\circ \Pi_p=\tilde  f_\upm,$  
and
   we have  
  $$    \nu_\upm(T,B,\tau) =\int_{B} f_\upm \omega^l. 
  $$
  \item   If  moreover $\omega$ is  K\"ahler, then    the above two assertions   still hold if 
    $T^\pm$ belong to the class $\CL_p^{1}(B).$
\end{enumerate}
 \end{theorem}
 \proof
  Since  the proof is  not difficult, we leave it   to the interested reader.   
 \endproof

 \begin{theorem}\label{T:charac-ph-top}
   We  keep the    Standing Hypothesis.  Suppose that $\omega$ is  K\"ahler on $B.$ 
Suppose  that    the current   $T$  is  positive pluriharmonic  and $T=T^+-T^-$  on an open  neighborhood of $\overline B$ in $X$ with $T^\pm$ in the class $\PH_p^{2}(B).$
Suppose in addition  that   
  there is  a strongly  admissible map $\tau$ for $B.$  Then  one  and  only one of the  following assertion holds:
  \begin{enumerate}

   \item If $\upm=k-p,$ then  
   $\nu_\upm(T,B,\tau)$ is  simply the mass of the  measure $T\wedge \pi^*(\omega^\upm)$ on $B.$
   \item  If   $\upm\not=k-p,$ then  $\upm=l$ and  the function  $\tilde f_\upm$  given by  Proposition  \ref{P:Cut-Off} (4) is constant
on fibers  of $ \Pi_p,$  that is,  we have $f_\upm\circ \Pi_p=\tilde  f_\upm,$  
and
   we have  
  $$    \nu_\upm(T,B,\tau) =\int_{B} f_\upm \omega^l. $$
   
\end{enumerate}
 \end{theorem}
  
\proof
 Since  the proof is  not difficult, we leave it   to the interested reader.   
\endproof

 \section{Geometric characterizations for positive plurisubharmonic  currents with holomorphic admissible maps}
 \label{S:Charac-psh-hol}
  
  In this  section we  assume  that  $\tau$ is  a holomorphic  admissible map and $p<k-l.$  The latter   assumption is  equivalent to   $\upm=l$ and $\upm\not=k-p.$

\subsection{Mass estimates}

    Let $T$ be  a positive plurisubharmonic  current  in the class $\SH^{2}_p(B)$ and $(T_n)_{n=1}^\infty$ a sequence of approximating forms
   for $T.$    So  $(\ddc T_n)_{n=1}^\infty$ a sequence of approximating forms
   for $\ddc T$   in the class $\CL^{0}_{p+1}(B).$
   
     Let    $(R_n)_{n=1}^\infty$ be  a  sequence of currents defined on $\bfU.$   For  an integer $0\leq j\leq k-l$ and a subsequence $(N_n)_{n=1}^\infty\subset \N,$ we  denote by $R^{(j)}$   the current
     $$(R)^{(j)} := \lim_{n\to\infty} \big( \tau_*R_{N_n} \wedge   \alpha^{j}\big)_\bullet,$$
   provided  that the limit (of course depending on  the choice of  the subsequence $(N_n)_{n=1}^\infty$) exists in the sense  of currents. Here $(\cdot_\bullet)$ denotes the  trivial  extension across $V$ in $\E.$
   In what  follows,  we apply this  notation for two sequences  $R_n:=T_n$ and $R_n:=\ddc T_n.$
  \begin{lemma}\label{L:mass-ddc-T-near-B-hol}
 There is a  subsequence $(T_{N_n})_{n=1}^\infty$
   such that  for every $ \lowm \leq j\leq\upm,$ and $1\leq m\leq  k-j-p$ and $0\leq  q<m$ and  every $0<r\leq\bfr,$ the integral
   \begin{equation*}
   \int_0^r  {2tdt\over  t^{2(q+1)}}\int_{\Tube(B,t)}(\ddc T)^{(m-q-1)} \wedge \beta^{k-j-p-m+q}\wedge \pi^*(\omega^j) 
   \end{equation*}
   is  finite non-negative.
 
  \end{lemma}
\proof
Since $\tau$ is  holomorphic, 
 Corollary \ref{C:Lelong-Jensen} applied   to $ \big( \tau_*T_{N_n} \wedge   \alpha^{m-q-1}\wedge \beta^{k-j-p-m+q+1}\wedge \pi^*(\omega^j) ,$ yields that
\begin{multline*}
 \int_{\Tube(B,r)} ( \tau_*T_{N_n} \wedge   \alpha^m )\wedge \beta^{k-j-p-m}\wedge  \pi^*(\omega^j)={1\over r^{2(q+1)}}\int_{\Tube(B,r)} 
 ( \tau_*T_{N_n} \wedge   \alpha^{m-q-1})\wedge \beta^{k-j-p-m+q+1}\wedge  \pi^*(\omega^j)\\
 -\lim_{n\to\infty} \int_0^r  {2tdt\over  t^{2(q+1)}}\int_{\Tube(B,t)}\tau_*(\ddc T_{N_n}) \wedge  \alpha^{m-q-1}\wedge \beta^{k-j-p-m+q}\wedge \pi^*(\omega^j)\\
 +{1\over r^{2(q+1)}}\int_0^r 2tdt\int_{\Tube(B,t)} \tau_*(\ddc T_{N_n})\wedge \alpha^{m-q-1}\wedge \beta^{k-j-p-m+q}\wedge  \pi^*(\omega^j).
\end{multline*}
We pass the limit of the above  equality for $n\to\infty.$
Observe  that the  LHS tends to   
$$
\int_{\Tube(B,r)} T^{(m)}\wedge \beta^{k-j-p-m}\wedge  \pi^*(\omega^j)=\int_{\Tube(B,r)} T^{(m)}\wedge \hat\beta^{k-j-p-m}\wedge  \pi^*(\omega^j),
$$
which is    finite non-negative by Proposition \ref{P:existence-finite}.
The first   integral on the RHS  tends to
$$
{1\over r^{2(q+1)}}\int_{\Tube(B,r)} T^{(m-q-1)}\wedge \beta^{k-j-p-m+q+1}\wedge  \pi^*(\omega^j)={1\over r^{2(q+1)}}\int_{\Tube(B,r)} T^{(m-q-1)}\wedge \hat\beta^{k-j-p-m+q+1}\wedge  \pi^*(\omega^j)
$$
which is finite non-negative by Proposition \ref{P:existence-finite}. 

To treat the second  and third integrals  on the RHS,  consider the  function $G_n:\ (0,\bfr]\to\R$ defined by
$$
G_n(t):=\int_{\Tube(B,t)}\tau_*(\ddc T_{N_n}) \wedge  \alpha^{m-q-1}\wedge \beta^{k-j-p-m+q}\wedge \pi^*(\omega^j).
$$
Since  we can write
$$
G_n(t)=\int_{\Tube(B,t)}\tau_*(\ddc T_{N_n}) \wedge  (\hat\alpha')^{m-q-1}\wedge \hat\beta^{k-j-p-m+q}\wedge \pi^*(\omega^j),
$$
the function $G_n$ is increasing  non-negative valued and the following limit hold for every $t$ except at most  a countable set  
$$
\lim_{n\to\infty}G_n(t)=\int_{\Tube(B,r)} T^{(m-q-1)}\wedge \beta^{k-j-p-m+q+1}\wedge  \pi^*(\omega^j).
$$ By Proposition \ref{P:existence-finite}, the sequence $G_n(\bfr)$  is  bounded.
Therefore, by Lebesgue dominated  convergence,
the third   integral on the RHS  tends to
\begin{multline*}
{1\over r^{2(q+1)}}\int_0^r 2tdt\int_{\Tube(B,t)} (\ddc T)^{(m-q-1)}\wedge \beta^{k-j-p-m+q}\wedge  \pi^*(\omega^j)\\=
{1\over r^{2(q+1)}}\int_0^r 2tdt\int_{\Tube(B,t)} (\ddc T)^{(m-q-1)}\wedge \hat\beta^{k-j-p-m+q}\wedge  \pi^*(\omega^j)
\end{multline*}
which is finite non-negative. 

On the other hand, the above  discussion shows that  the  second  integral on the RHS is uniformly bounded indepentdent of $n$ since other integrals are so.
By Fatou lemma, we infer that  $\int_0^r  {2tdt\over  t^{2(q+1)}}\liminf_{n\to\infty} G_n(t)<\infty.$
This proves the lemma.


\endproof

Now  we come back Lemma \ref{L:mass-ddc-T-near-B-hol} for $q=0.$

\begin{lemma}\label{L:S(j)-widetilde-S-hol}
 There exist   currents $S^{(0)},\ldots, S^{(k-l-p-1)}$ on $\U$ and
 a current $\widetilde S$ on $\Pi_{p+1}^{-1}(\U)$ such that  for a  suitable  subsequence $(T_{N_n})_{n=1}^\infty,$  we have
 \begin{eqnarray*}
  \lim_{N\to\infty} \big(  -\log\varphi \cdot \tau_* (\ddc T_{N_n})\wedge \alpha^m    \big)_\bullet&=&S^{(m)}\quad\text{for}\quad m=0,\ldots, k-l-p-1,\\
  \lim_{N\to\infty} \big[\Pi^\diamond_{p+1}\big(  -\log\varphi\cdot \tau_* ( \ddc T_{N_n})   \big)\big]_\bullet&=&\widetilde S.
 \end{eqnarray*}
Moreover, for $\lowm \leq j\leq \upm$ and $m=1,\ldots, k-l-p-1,$
\begin{multline*}
   \lim_{n\to\infty}\int_0^r  {2tdt\over  t^{2}}\int_{\Tube(B,t)}\tau_*(\ddc T_{N_n}) \wedge  \alpha^{m-1}\wedge \beta^{k-j-p-m}\wedge \pi^*(\omega^j)\\
   =\log r\int_{\Tube(B,r)}(\ddc T)^{(m-1)} \wedge \beta^{k-j-p-m}\wedge \pi^*(\omega^j)+  
   \int_{\Tube(B,r)} S^{(m-1)} \wedge   \beta^{k-j-p-m}\wedge \pi^*(\omega^j).
   \end{multline*}

\end{lemma}
\proof   By  Fubini's theorem we have 
\begin{equation}\label{e:L:S(j)-widetilde-S-hol}
\begin{split}
 &\int_0^r  {2tdt\over  t^{2}}\int_{\Tube(B,t)}\tau_*(\ddc T_{N_n}) \wedge  \alpha^{m-1}\wedge \beta^{k-j-p-m}\wedge \pi^*(\omega^j)\\
 &=\int_{y\in \Tube(B,t)}\big(\int_{|y|}^r  {2tdt\over  t^{2}}\big)\tau_*(\ddc T_{N_n}) \wedge  \alpha^{m-1}\wedge \beta^{k-j-p-m}\wedge \pi^*(\omega^j)\\
 &= \log{ r}  \int_{\Tube(B,r)}\tau_*(\ddc T_{N_n}) \wedge  \alpha^{m-1}\wedge \beta^{k-j-p-m}\wedge \pi^*(\omega^j)\\
 &+  \int_{\Tube(B,r)}(-\log {\varphi})\tau_*(\ddc T_{N_n}) \wedge  \alpha^{m-1}\wedge \beta^{k-j-p-m}\wedge \pi^*(\omega^j).
 \end{split}
 \end{equation}
  By Lemma \ref{L:mass-ddc-T-near-B-hol}, the LHS  converges as $n\to\infty.$
By  Proposition \ref{P:existence-finite}  applied to $\ddc T_{N_n},$  the  integral in the third line of \eqref{e:L:S(j)-widetilde-S-hol} also  converges to
$$
\log{ r}  \int_{\Tube(B,r)}\ddc \big(T_{N_n}^{(m-1)}\big)\wedge \beta^{k-j-p-m}\wedge \pi^*(\omega^j)
$$
for every $r$  except  at most a  countable  set.  Observe  that the  following  integral is  a finite  linear combination with  real coefficients  of the integral on the last line   of \eqref{e:L:S(j)-widetilde-S-hol}:
$$
 \int_{\Tube(B,r)}(-\log {\varphi})\tau_*(\ddc T_{N_n}) \wedge (\hat \alpha')^{m-1}\wedge \hat\beta^{k-j-p-m}\wedge \pi^*(\omega^j)\geq 0.
$$
Therefore, we infer that the  latter integral  is uniformly bounded, that is
$$
\sup_{n\geq 1}\int_{\Tube(B,\bfr)}(-\log {\varphi})\tau_*(\ddc T_{N_n}) \wedge  \alpha^{m-1}\wedge \beta^{k-j-p-m}\wedge \pi^*(\omega^j)<\infty.
$$
This, combined with Propositions \ref{P:existence} and  \ref{P:existence-finite}, implies the existence of the currents  
  $S^{(0)},\ldots, S^{(k-l-p-1)}$  and $\widetilde S.$

  Now  we prove that these  currents enjoy the cut-off property through $B$ in $\Tube(B,\bfr)$.  Consider the  following  currents on $  \U\setminus V$:
  \begin{equation*}
  R_n:=-(\log {\varphi})\tau_*(\ddc T_{N_n}) \wedge  \alpha^{m}.
  \end{equation*}
  Arguing as in the proof of Theorem \ref{T:charac-closed-all-degrees},
  we  see that $\ddc  (R_n)_\bullet=(\ddc  R_n)_\bullet.$
So we  infer that
  \begin{equation}\label{e:ddc-R_n} \begin{split} \ddc  R_n&= -\ddc (\tau_*T_{N_n}) \wedge  \alpha^{m+1}= - \tau_*(\ddc T_{N_n}) \wedge (\hat \alpha'-c_1\pi^*\omega)^{m+1} \\
  &=  -\sum_{j=0}^{m+1} {m+1\choose j}  c_1^j \tau_*(\ddc T_{N_n}) \wedge (\hat \alpha')^{m+1-j}\wedge \pi^*(\omega^j).
  \end{split}\end{equation}
  Since  by passing to a  subsequence if necessary, the weak limit  $\lim_{n\to\infty}\tau_*(\ddc T_{N_n}) \wedge (\hat \alpha')^{m+1-j}\wedge \pi^*(\omega^j)$ is  a positive closed current,  we infer that  both  $S^{(m)}$ and  $\ddc S^{(m)}$ are  currents of order $0$.  Hence, by Theorem \ref{T:Bassanelli}, $S^{(0)},\ldots, S^{(k-l-p-1)}$ enjoy the cut-off property through $B$ in $\Tube(B,\bfr).$
 
  Let
  $$
  S_n:=\Pi_{p+1}^\diamond (R_n)= \Pi_{p+1}^\diamond\big( (-\log {\varphi})\ddc (\tau_*T_{N_n}) \wedge  \alpha^{m}\big) \qquad\text{on}\qquad  \Pi_{p+1}^{-1}(\U\setminus V).
  $$
  We get $\ddc S_n=\Pi_{p+1}^\diamond (\ddc R_n)$ on $ \Pi_{p+1}^{-1}(\U\setminus V).$
  
Let us  show that  
\begin{equation} \label{e:ddc-S_n-bullet}\ddc  (S_n)_\bullet=(\ddc  S_n)_\bullet.
\end{equation}
To this end let $\Phi$ be a smooth test form  compactly supported in   $\Pi_{p+1}^{-1}(\U).$  We have
\begin{eqnarray*}
\langle \ddc  (S_n)_\bullet-(\ddc  S_n)_\bullet, \Phi\rangle&=&\lim_{r\to 0}\int_{|z\circ \Pi_{p+1}|>r} S_n\wedge \ddc \Phi-\ddc S_n\wedge \Phi\\
&=&\lim_{r\to 0}\big( -\int_{|z\circ \Pi_{p+1}|=r} S_n\wedge \dc \Phi-\dc S_n\wedge \Phi\big)\\
&=&\lim_{r\to 0}\big( \log r\int_{\varphi\circ \Pi_{p+1}=r^2} \Pi_{p+1}^\diamond(\ddc (\tau_*T_{N_n}))\wedge \alpha^m)\wedge \dc \Phi\\
&+&{1\over r^2}  \int_{\varphi\circ \Pi_{p+1}=r^2} (\dc\varphi \circ \Pi_{p+1})\wedge  \Pi_{p+1}^\diamond(\ddc (\tau_*T_{N_n}))  \wedge \Phi\big).
\end{eqnarray*}
  Since the  forms involved in the last integrals are  all  $\Cc^3$-smooth  in   $\Pi_{p+1}^{-1}(\U),$ the last limit is  equal to $0.$
  This  proves \eqref{e:ddc-S_n-bullet}.
  
  Next, we insert    the  expression of $\ddc R_n$ given  in  \eqref{e:ddc-R_n} into the equality $\ddc S_n=\Pi_{p+1}^\diamond (\ddc R_n)$ on $ \Pi_{p+1}^{-1}(\U\setminus V)$ and 
  use \eqref{e:ddc-S_n-bullet}.  Since  by passing to a  subsequence if necessary, the weak limit 
  $\lim_{n\to\infty}\Pi_{p+1}^\diamond\big(\tau_*(\ddc T_{N_n}) \wedge (\hat \alpha')^{m+1-j}\wedge \pi^*(\omega^j)\big)$ is  a positive closed current,, we see that  both  $\widetilde S$ and  $\ddc \widetilde S$ are  currents of order $0$.  Hence, by Theorem \ref{T:Bassanelli}, $\widetilde S$ enjoy the cut-off property through $B$ in $\Tube(B,\bfr).$
 
 Taking the limit in \eqref{e:L:S(j)-widetilde-S-hol} for $n\to\infty,$ we get the last identity of the lemma.
  \endproof
Recall  from the  above discussion  that
$$
T^{(1)}:=\lim_{n\to\infty} \big(\tau_* (T_{N_n})\wedge \alpha \big)_\bullet\qquad\text{and}\qquad
S^{(0)}:=  \lim_{n\to\infty} \big((-\log\varphi) \cdot \tau_* (\ddc T_{N_n})\big)_\bullet
$$
 for a  suitable  subsequence $(T_{N_n})_{n=1}^\infty.$

\begin{lemma}\label{L:T(1)-S(0)-hol}
 The following identity holds
 \begin{equation*}
  \nu(T,B,\tau)=\nu(T^{(1)}+S^{(0)}, B,\id).
 \end{equation*}
\end{lemma}

\proof
Applying    Theorem  \ref{T:Lelong-Jensen-smooth} and  Corollary \ref{C:Lelong-Jensen} to  the current $\tau_*T_n\wedge\beta^{k-l-p-1}\wedge \pi^*(\omega^l)$ and  for $q=1$ yields for $0<r\leq\bfr$ that
\begin{equation}\label{T(1)-S(0)-hol-Lelong-Jensen}
 \begin{split}
 \int_{\Tube(B,r)}\tau_*T_n\wedge\alpha \wedge\beta^{k-l-p-1}\wedge \pi^*(\omega^l)&={1\over r^2}
\int_{\Tube(B,r)}\tau_*T_n\wedge\beta^{k-l-p}\wedge \pi^*(\omega^l)\\
 &-\int_0^r{2tdt\over t^2}\int_{\Tube(B,t)}\ddc (\tau_*T_n)\wedge\beta^{k-l-p-1}\wedge \pi^*(\omega^l)\\
 &+ {1\over r^2}\int_0^r2tdt\int_{\Tube(B,t)}\ddc (\tau_* T_n)\wedge\beta^{k-l-p-1}\wedge \pi^*(\omega^l).
 \end{split}
\end{equation}
Since $\tau$ is  holomorphic, we have   $\ddc(\tau_*T_n)=\tau_*(\ddc T_n)$ on  the  RHS of \eqref{T(1)-S(0)-hol-Lelong-Jensen}.
Observe  that  by  Propositions  \ref{P:existence} and \ref{P:existence-finite}, the  LHS of \eqref{T(1)-S(0)-hol-Lelong-Jensen} converges  as  $n\to\infty$  to
$$
\int_{\Tube(B,r)} T^{(1)} \wedge\beta^{k-l-p-1}\wedge \pi^*(\omega^l)
$$
for every $r\in (0,\bfr]$ except at most a countable values of $r.$  On the other hand,  the first integral on the  RHS of \eqref{T(1)-S(0)-hol-Lelong-Jensen} converges  as  $n\to\infty$  to
$$
{1\over r^2}
\int_{\Tube(B,r)}\tau_*T\wedge\beta^{k-l-p}\wedge \pi^*(\omega^l).
$$
Since 
$$\int_{\Tube(B,t)}\tau_*(\ddc T_n)\wedge\beta^{k-l-p-1}\wedge \pi^*(\omega^l)=\int_{\Tube(B,t)}\tau_*(\ddc T_n)\wedge\hat\beta^{k-l-p-1}\wedge \pi^*(\omega^l)$$
is a nonnegative increasing function of $t\in (0,\bfr],$  we infer from the  dominated convergence  theorem that
the limit of the last integral on the RHS of \eqref{T(1)-S(0)-hol-Lelong-Jensen} is
$$
{1\over r^2}\int_0^r2tdt\int_{\Tube(B,t)}\tau_*(\ddc T)\wedge\beta^{k-l-p-1}\wedge \pi^*(\omega^l).
$$
By Lemma  \ref{L:S(j)-widetilde-S-hol} for $m=1,$ the  second integral on the  RHS of \eqref{T(1)-S(0)-hol-Lelong-Jensen} converges  as  $n\to\infty$  to
$$
-\log r\int_{\Tube(B,r)}\tau_*(\ddc T)\wedge\beta^{k-l-p-1}\wedge \pi^*(\omega^l)-\int_{\Tube(B,r)}S^{(0)}\wedge\beta^{k-l-p-1}\wedge \pi^*(\omega^l).
$$
Summing up, we obtain  that
\begin{equation}\label{e:T(1)-S(0)-hol}
\begin{split}
&{1\over r^{2(k-l-p-1)} } \int_{\Tube(B,r)}(T^{(1)}+S^{(0)})\wedge\beta^{k-l-p-1}\wedge \pi^*(\omega^l)\\
&={1\over r^{2(k-l-p)} } \int_{\Tube(B,r)}\tau_*T\wedge\beta^{k-l-p}\wedge \pi^*(\omega^l)\\
&-{\log r\over r^{2(k-l-p-1)} } \int_{\Tube(B,r)}\tau_*(\ddc T)\wedge\beta^{k-l-p-1}\wedge \pi^*(\omega^l)\\
&+{1\over r^{2(k-l-p)} }\int_0^r  2tdt \int_{\Tube(B,t)}\tau_*(\ddc T)\wedge\beta^{k-l-p-1}\wedge \pi^*(\omega^l).
\end{split}
\end{equation}
We will show that all  terms of  \eqref{e:T(1)-S(0)-hol}  converge as $r\to 0$ and that  the last two terms on the RHS of \eqref{e:T(1)-S(0)-hol}  converge to $0$ as $r\to 0.$
Consider
$$
G(t):=\int_{\Tube(B,t)} \tau_*(\ddc T)\wedge  \beta^{k-l-p-1}\wedge  \pi^*(\omega^l). 
$$
 Lemma \ref{L:mass-ddc-T-near-B-hol} with $m=k-l-p$ and $q=m-1$ implies that the function
 $(0,\bfr]\ni t\mapsto  {2t\over t^{2(k-l-p)}}G(t)$ is  integrable. Hence,
 $$
  {1\over r^{2(k-l-p)}}\int_0^r 2tG(t)dt\leq  \int_0^r{2t\over t^{2(k-l-p)}}G(t)
 $$
 and the RHS converges to $0$ as $r\to 0.$ So the last integral of the RHS of \eqref{e:T(1)-S(0)-hol} converges   to $0.$

 On the other hand, observe that
 \begin{eqnarray*}
 T^{(1)} +S^{(0)} &=&\lim_{n\to\infty} \big( \tau_*T_{N_n}\wedge\alpha-(\log\varphi)\tau_*(\ddc T_{N_n}) \big)_\bullet \, ,\\
 \ddc\big( \tau_*T_{N_n}\wedge\alpha-(\log\varphi)\tau_*(\ddc T_{N_n}) \big)&= &\tau_*(\ddc T_{N_n})\wedge\alpha-\ddc (\log\varphi)\wedge \tau_*(\ddc T_{N_n})= 0.
 \end{eqnarray*}
 Consequently, arguing as in the  proof of  Theorem \ref{T:charac-closed-all-degrees}   we can show that
 $$
 \ddc\big( \tau_*T_{N_n}\wedge\alpha\big)_\bullet-\ddc\big((\log\varphi)\ddc (\tau_*T_{N_n}) \big)_\bullet=0.
 $$
 Therefore, by passing $n\to\infty$ we  infer that $T^{(1)} +S^{(0)}$ is a pluriharmonic $(p+1,p+1)$-current. 
 
 Moreover,
 write 
 $$
 \tau_*T_{N_n}\wedge\alpha-(\log\varphi)\tau_*(\ddc T_{N_n})=\big(\tau_*T_{N_n}\wedge\hat\alpha'-(\log\varphi)\ddc (\tau_*T_{N_n})\big)
 -c_1\tau_*T_{N_n}\wedge\pi^*\omega.
 $$
 Since  $T_{N_n}$  is a positive  plurisubharmonic $\Cc^3$-smooth forms and $\hat\alpha'$, $\pi^*\omega$  are  positive  smooth $(1,1)$-forms,  we can check   that  both   forms on the RHS are positive  plurisubharmonic.  By passing $n\to\infty$ we  see that  $T^{(1)} +S^{(0)}$
 is  the  difference of  two  plurisubharmonic $(p+1,p+1)$-current. Therefore, by Theorem \ref{T:Lelong-psh},
  the LHS of \eqref{e:T(1)-S(0)-hol}  converges  to $\nu( T^{(1)} +S^{(0)},B,\id)$ as $r\to 0$ and the limit  is finite.
  
  Next, by Theorem \ref{T:Lelong-psh}  the first term on the LHS of \eqref{e:T(1)-S(0)-hol}  converges  to $\nu(T,B,\tau)$ as $r\to 0.$

 Therefore, all terms in \eqref{e:T(1)-S(0)-hol} (except the second one on the RHS) converge as $r\to 0.$
 Hence,  the second term on the RHS  also  converges and its limit is  finite, 
 in other  word,
$\lim_{r\to 0} {-\log r\over  r^{2(k-l-p)}}G(r)\in\R.$ This, coupled with
the finiteness of  $\int_0^r{2t\over t^{2(k-l-p)}}G(t)dt,$ implies that $\lim_{r\to 0} {-\log r\over  r^{2(k-l-p-1)}}G(r)=0.$
In summary, we have  shown  that on the RHS of \eqref{e:T(1)-S(0)-hol}, as $r\to 0$  the first integral converges to  $\nu(T,B,\tau)$ and 
the   last two integrals   converges  to $0.$ So the integral on  the LHS of  \eqref{e:T(1)-S(0)-hol} also converges to   $\nu(T,B,\tau).$ This completes 
 the proof.
 
\endproof




\begin{lemma}\label{L:T_n-change-alpha-to-beta-hol}
 For every $0<r\leq \bfr,$ we have that
 \begin{multline*}
  \int_{\Tube(B,r)} \big( \tau_*T_{n}\wedge\alpha-(\log\varphi)\ddc (\tau_*T_{n})\big)\wedge \alpha^{k-l-p-1}\wedge \pi^*(\omega^l)\\
  ={1\over r^{2(k-l-p-1)}}\int_{\Tube(B,r)} \big( \tau_*T_{n}\wedge\alpha-(\log\varphi)\ddc (\tau_*T_{n})\big)\wedge \beta^{k-l-p-1}\wedge \pi^*(\omega^l)
 \end{multline*}
\end{lemma}
\proof
Set
$$  S_n:=   \tau_*T_{n}\wedge\alpha-(\log\varphi)\ddc (\tau_*T_{n})\wedge \pi^*(\omega^l).$$ We know that  
$\ddc S_n=0.$  Applying  Theorem \ref{T:Lelong-Jensen-compact-support}  and Corollary \ref{C:Lelong-Jensen} to $S_n$  and $q=k-l-p-1,$  we  obtain  for $0<s<r\leq\bfr$ that
\begin{equation}\label{e:T_n-change-alpha-to-beta-hol}
\begin{split}
 &\int_{\Tube(B,s,r)}  S_n\wedge \alpha^{k-l-p-1}
 \\
  &= {1\over r^{2(k-l-p-1)}}\int_{\Tube(B,r)}S_n\wedge \beta^{k-l-p-1}
  -{1\over s^{2(k-l-p-1)}}\int_{\Tube(B,s)} S_n\wedge \beta^{k-l-p-1}.
\end{split}\end{equation}
Using the  expression of $S_n$ and  the $\Cc^2$-smoothness of $T_n,$ we  can show by Lemma \ref{L:vertical-boundary-terms} that
$$
\lim_{s\to 0}{1\over s^{2(k-l-p-1)}}\int_{\Tube(B,s)} S_n\wedge \beta^{k-l-p-1}=0\quad\text{and}\quad \lim_{s\to 0}\int_{\Tube(B,s)}  S_n\wedge \alpha^{k-l-p-1}=0.
$$

Letting $s\to 0$ in \eqref{e:T_n-change-alpha-to-beta-hol} and using the above  discussion, the result  follows.
\endproof

\subsection{Geometric characterizations}

The following central result of the section provides a geometric characterization of the  top Lelong number in the case of holomorphic admissible maps. 

\begin{theorem}\label{T:charac-Lelong-for-psh-hol}
 Let  $T$ be  a positive plurisubharmonic  current  in the class $\SH^{2}_p(B)$ with $(T_n)_{n=1}^\infty$ a  sequence of approximating forms.
 Then  there exists  a subsequence $(T_{N_n})_{n=1}^\infty$ and  an open  neighborhood $\U'$ of $\overline B$ in $\E$  
 with $\U'\subset  \U$ such  that the  following  properties  holds.
 \begin{enumerate}
 \item  The following  currents  are well-defined:
 \begin{eqnarray*}
  \widetilde T&:=&\lim_{n\to\infty} \Pi^*_p  (\tau_*T_{N_n})\qquad\text{on}\qquad \Pi_p^{-1}(\U')\\
  \widetilde S&:=&\lim_{n\to\infty}\Big( \Pi^\diamond_{p+1}  \big((-\log{\varphi})\cdot \tau_*(\ddc T_{N_n})\big)\Big)_\bullet\qquad\text{on}\qquad \Pi_{p+1}^{-1}(\U').
 \end{eqnarray*}
\item
There exist two  functions $f,g\in L^1_\loc(B)$   such that
$$
\ind_{\Pi_p^{-1}(B)}=(f\circ \Pi_p)[ \Pi_p^{-1}(B)]\qquad\text{and}\qquad \ind_{\Pi_{p+1}^{-1}(B)}=(g\circ \Pi_p)[ \Pi_{p+1}^{-1}(B)].
$$
Moreover, both   function  $f$ and $f+g$ are  non-negative and  $f$ is    plurisubharmonic on $B$  and  $f+g$ is the difference of two  plurisubharmonic functions on $B$ and 
\begin{equation*}
 \nu(T,B,\tau)=\int_B (f+g)\omega^l.
\end{equation*}
\end{enumerate}
\end{theorem}
\proof 
Let  $\U'$ be an open  neighborhood of $\overline B$ in $\E$   such that 
 with $\U'\subset  \U$ and that
all  currents $T_{n}$'s are defined on $\U'.$ 

The existence of $\widetilde T$ follows from    Proposition \ref{P:Cut-Off-hol}.  The  existence of $ \widetilde S$ is  a consequence of Lemma  \ref{L:S(j)-widetilde-S-hol}. This completes the proof of assertion (1).

  We turn to the  proof of  assertion (2). Propositions \ref{P:existence-finite}  applied  to  the sequence $(T_{n})_{n=1}^\infty$ yields that
  the   sequence $ ( T_n\wedge \alpha\big)_{n=1}^\infty$
satisfies   the assumption  of  Proposition
\ref{P:existence-finite}.  Consequently,    we can extract  a  subsequence $(T_{N_n})_{n=1}^\infty$
  such  that
the  current
$$
\widehat T:=\lim_{n\to\infty}  \big(  \Pi^\diamond_{p+1}(\tau_*T_{N_n}\wedge \alpha)  \big)_\bullet
$$
is  well-defined in  $ \Pi_{p+1}^{-1}(\U').$
Moreover, by  Lemma \ref{L:S(j)-widetilde-S-hol}, the current $\widetilde S$ in the  statement of assertion (1) is  well-defined.

Now  we will show that
\begin{equation}\label{e:charac-Lelong-for-psh-hol(1)}
 \ddc  \big(  \Pi^\diamond_{p+1}(\tau_*T_{N_n}\wedge \alpha)  \big)_\bullet=   \big[\ddc  \big(  \Pi^\diamond_{p+1}(\tau_*T_{N_n}\wedge \alpha)\big)\big]_\bullet.
\end{equation}
To this end let $\Phi$ be a smooth  test form compactly  supprted in  $\Pi_{p+1}^{-1}(\U').$  We have  that
\begin{eqnarray*}
 &&\left\langle \ddc  \big(  \Pi^\diamond_{p+1}(\tau_*T_{N_n}\wedge \alpha)  \big)_\bullet-   \big[\ddc  \big(  \Pi^\diamond_{p+1}(\tau_*T_{N_n}\wedge \alpha)\big)\big]_\bullet,   \Phi\right\rangle\\
 &=&\lim_{r\to 0}  \big[\int_{\varphi\circ \Pi_{p+1}>r^2}    \big(  \Pi^\diamond_{p+1}(\tau_*T_{N_n}\wedge \alpha)\wedge \ddc\Phi-\ddc  \big(  \Pi^\diamond_{p+1}(\tau_*T_{N_n}\wedge \alpha)\wedge \Phi\big]\\
 &=&\lim_{r\to 0} \big [  \int_{\varphi\circ \Pi_{p+1}=r^2}     \Pi^\diamond_{p+1}(\tau_*T_{N_n}\wedge \alpha)\wedge \dc\Phi-\dc  \big(  \Pi^\diamond_{p+1}(\tau_*T_{N_n}\wedge \alpha)\wedge \Phi \big].
\end{eqnarray*}
Applying Lemma \ref{L:Lelong-Jensen_1}, the last line is  equal to 
$$
\lim_{r\to 0} {1\over r^2}\big [  \int_{\varphi\circ \Pi_{p+1}=r^2}     \Pi^\diamond_{p+1}(\tau_*T_{N_n}\wedge \beta)\wedge \dc\Phi-\dc  \big(  \Pi^\diamond_{p+1}(\tau_*T_{N_n}\wedge \beta)\wedge \Phi \big].
$$
The last  limit is  equal to $0$  as  $T_{N_n}$  is a $\Cc^3$-smooth form. This  proves \eqref{e:charac-Lelong-for-psh-hol(1)}.

Using \eqref{e:charac-Lelong-for-psh-hol(1)} we can  show that $\ddc \widehat T$ is a current of order $0.$ Hence,  $\widehat  T$ is  $\C$-normal.

On the other hand, we can show that
$$
\ddc  \widetilde S=  \lim_{n\to\infty} \Pi^\diamond_{p+1}  \big(-\ddc (\tau_*T_{N_n})\wedge \alpha \big)_\bullet\qquad\text{on}\qquad \Pi_{p+1}^{-1}(\U').$$
This implies that $ \ddc\widetilde S$  is  a current of order $0.$  So  $\widetilde S$ is $\C$-normal.

Summing up, we have shown that $\widehat T$ and $\widetilde S$  are both $\C$-normal.
By  Proposition 
\ref{P:Cut-Off}, there  exist  non-negative functions $\hat f,$ $\tilde g\in L^1_\loc(\Pi_{p+1}^{-1}(B))$ 
such that
\begin{equation}\label{e:charac-Lelong-for-psh-hol(2)}
 \ind_{\Pi_{p+1}^{-1}(B)}\widehat T=\hat f [\Pi_{p+1}^{-1}(B)]\qquad\text{and}\qquad \ind_{\Pi_{p+1}^{-1}(B)}\widetilde  S=\tilde g [\Pi_{p+1}^{-1}(B)].
\end{equation}
Moreover,  by  Theorem  \ref{T:Bassanelli-psh} $\hat f$ and $\hat f+\hat g$ are  positive plurisubharmonic. Therefore,  they are constant on fibers.
So there are   functions $ f_0,$ $ g\in L^1_\loc(B)$  such that
\begin{equation}\label{e:hatf-f_tildeg-g}\hat f=f_0\circ \Pi_{p+1}\qquad \text{and}\qquad   \tilde  g=g\circ \Pi_{p+1}\qquad \text{on}\qquad\Pi_{p+1}^{-1}(B).
\end{equation}
By  Lemma \ref{L:T(1)-S(0)-hol},  we have 
\begin{multline*}
 \nu(T,B,\tau)=\nu(T^{(1)}+S^{(0)},B,\tau)\\
 =\lim_{r\to 0}\lim_{n\to\infty} {1\over  r^{2(k-l-p-1)}}\int_{\Tube(B,r)} \big(\tau_*T_{N_n}\wedge\alpha-(\log\varphi)\tau_*(\ddc T_{N_n}) \big)\wedge\beta^{k-l-p-1}\wedge \pi^*(\omega^l).
\end{multline*}
By  Lemma 
\ref{L:T_n-change-alpha-to-beta-hol}, the last line  is  equal to
\begin{equation*}
\lim_{r\to 0}\lim_{n\to\infty} \int_{\Tube(B,r)} \big(\tau_*T_{N_n}\wedge\alpha-(\log\varphi)\tau_*(\ddc T_{N_n}) \big)\wedge\alpha^{k-l-p-1}\wedge \pi^*(\omega^l).
\end{equation*}
Since $\alpha^{k-l-p-1}\wedge \pi^*(\omega^l)=\alpha_\ver^{k-l-p-1}\wedge \pi^*(\omega^l),$ the last line  is  equal  to 
\begin{equation*}
\lim_{r\to 0}\lim_{n\to\infty} \int_{\Tube(B,r)} \big(\tau_*T_{N_n}\wedge\alpha-(\log\varphi)\tau_*(\ddc T_{N_n}) \big)\wedge\alpha_\ver^{k-l-p-1}\wedge \pi^*(\omega^l)
\end{equation*}
 By assertion (1), this is  equal to
\begin{equation*}
\lim_{r\to 0}  \int_{\Pi_{p+1}^{-1}(\Tube(B,r))} (\widehat T+\widetilde S )\wedge\Pi_{p+1}^\diamond\big (\alpha_\ver^{k-l-p-1}\wedge \pi^*(\omega^l)\big).
\end{equation*}
By  Lemma \ref{L:Upsilon_j} and equality \eqref{e:hatf-f_tildeg-g},  the last  expression is  equal to
\begin{equation*}
 \lim_{r\to 0}  \int_{\Pi_{p+1}^{-1}(\Tube(B,r))} (\widehat T+\widetilde S )\wedge\Pr\nolimits_{p+1}^\diamond(\Upsilon^{(p+1)(k-l-p-1)})\wedge \Pi_{p+1}^\diamond\big(\pi^*(\omega^l)\big)
 =\int_B (f_0+g)\omega^l,
\end{equation*}
where the  equality  follows from \eqref{e:charac-Lelong-for-psh-hol(2)} and  Proposition \ref{P:Cut-Off}.

On the other  hand,
by Proposition \ref{P:Cut-Off-hol} for $j=l$ and hence $\hat j=p$, there is a function  $f\in L^1_\loc(B)$   such that 
$$\ind_{\Pi_p^{-1}(B)}\widetilde T=(f\circ \Pi_p)[ \Pi_p^{-1}(B)]$$
and that for every  $\Cc^2$-piecewise smooth subdomain  $D\subset B,$ 
\begin{equation*}
 \lim_{r\to 0}\lim_{n\to\infty} \int_{\Tube(D,r)} \tau_*T_{N_n}\wedge\alpha^{k-l-p}\wedge \pi^*(\omega^l)=\|\widetilde T\|(\Pi^{-1}_p(D))
 =\int_D f\omega^l.
\end{equation*}
Observe that the expression on the LHS is  also equal to $\|\widehat T\|(\Pi^{-1}_{p+1}(D))=\int_D f_0\omega^l.$
So   $\int_D (f_0-f)\omega^l=0.$ Since this equality holds for every  $\Cc^2$-piecewise smooth subdomain  $D\subset B,$ we  infer that $f_0=f.$
The proof of assertion (2) is thereby completed.
\endproof

\begin{example}
 \rm    We place ourselves in the  setting of   Subsection \ref{SS:Local-setting}.
 Suppose  that    $p<  k-l $ and consider  the current  $T:=\alpha^p$ on $U$ and a domain  with $\Cc^2$-piecewise smooth boundary
 $B\subset U''.$ Let  $\lowm\leq j\leq \upm.$

 If $j<k-p,$ we have  $\alpha^p\wedge \omega_z^{k-p-j}=\|z\|^{-2p} \omega_z^{k-j},
 $ and hence
 \begin{eqnarray*}
 \nu(T,B,\id)&=&\lim_{r\to 0}  {1\over r^{2(k-p-j)}}\int_{\|z\|<r,\ w\in B} T\wedge \omega_w^j \wedge \omega_z^{k-p-j}\\
 &=&\lim_{r\to 0}  \int_{w\in B}\big( {1\over r^{2(k-p-j)}}\int_{\|z\|<r}\|z\|^{-2p} \omega_z^{k-j}\big)\wedge \omega_w^j\\
 &=& \int_B\omega_w^j.
 \end{eqnarray*}
 Hence,  if $j<k-p,$  we have
 $
 \nu(T,B,\id)=
               1$ for $j=l$ and  
              $
 \nu(T,B,\id)= 0 $ otherwise.

 If  $j=k-p,$  we can  show that $\nu(T,B,\id)=0$ using that $j>l.$
 
 Summing  up,  the only  nonzero Lelong number   is $\nu_\top(T,B,\id).$
 
 For every $n\geq 1$ pick $\varphi_n\in\Cc^\infty_0({1\over 2n},\infty)$ with $\varphi_n\geq 0$ and $\int \varphi_n(t)dt=1.$
 Consider the functions
 $$
 \psi_n:={\varphi_n\over  -t^p\log t}\qquad\text{and}\qquad \chi_n(t):=\int_0^t{1\over s} \big( \int_0^s\psi_n(r)dr \big)ds.
 $$
 So $\chi_n$ is a smooth non-negative function   with $\supp(\chi_n)\subset ({1\over 2n},\infty).$
 Consider the  smooth $(p,p)$-form  $T_n:=\chi_n(\|z\|^2) T.$  We have 
 $$
 \ddc T_n= \big(\chi''_n(\|z\|^2)+{\chi'_n(\|z\|^2)\over \|z\|^2}\big) \omega_z^{p+1}=\psi_n(\|z\|^2) \omega_z^{p+1} \geq 0.
 $$
 So $T_n$ is smooth plurisubharmonic.
 

 We have for a suitable  subsequence $(T_{N_n})_{n=1}^\infty,$
 \begin{eqnarray*}
\int_B f\omega^l&=&\lim_{r\to 0}\lim_{n\to\infty} \int_{\Tube(B,r)} T_{N_n}\wedge\alpha \wedge\alpha^{k-l-p-1}\wedge \pi^*(\omega^l_w)\\
&=& \lim_{r\to 0}\lim_{n\to\infty} \int_{\|z\|<r,\ w\in B }  \alpha^{k-l}\wedge \pi^*(\omega^l_w)
=0,
\end{eqnarray*}
where the last equality holds  because $\alpha^{k-l}=0.$
 So $f=0$ almost everywhere on $B.$

 Similarly,
 \begin{eqnarray*}
\int_B g\omega^l\lim_{r\to 0}&=&\lim_{n\to\infty} \int_{\Tube(B,r)} -(\log\varphi)\ddc (T_{N_n}) \wedge\alpha^{k-l-p-1}\wedge \pi^*(\omega^l_w)\\
&=&\lim_{n\to\infty} \int_{\|z\|<r,\ w\in B} -(\log{\|z\|})\psi_n(\|z\|^2)\omega_z^{p+1}\wedge\alpha^{k-l-p-1}\wedge \pi^*(\omega^l_w)\\
&=&\big(\int_{B}\omega_w^l\big)\big(\lim_{n\to\infty} \int_{\|z\|<r} -(\log{\|z\|})\|z\|^{-2(k-l-p-1)} \psi_n(\|z\|^2)\omega_z^{k-l}\big).
\end{eqnarray*}
Using  polar coordinates, the second integral of the last line is  equal to
$$
\int_0^r\varphi_n(\rho^2) \rho d\rho=1.
$$
\end{example}

\begin{remark}
 \rm The  above  example  shows that  the decomposition 
 \begin{equation*}
\nu(T,B,\tau)=\int_B (f+g)\omega^l
\end{equation*}
 depends on the choice of the  approximating sequences. 
 
 This  example  also shows that in general (when $T$ is a positive plurisubharmonic  current),
 the term  $\int_B g\omega^l$ expressing the mass of $\widetilde S$ over $B$ is  necessary.
\end{remark}

 \section{Geometric characterizations for positive plurisubharmonic  currents with strongly admissible maps}
 \label{S:Charac-psh}
  
  As in the  previous  section,  we  assume in this  section that $\upm=l$ and $\upm\not=k-p.$ This  assumption is  equivalent to  $p<k-l.$

\subsection{Mass estimates}

    Let $T$ be  a positive plurisubharmonic  current  in the class $\SH^{3,3}_p(B)$ and $(T_n)_{n=1}^\infty$ a sequence of approximating forms
   for $T.$    So  $(\ddc T_n)_{n=1}^\infty$ a sequence of approximating forms
   for $\ddc T$   in the class $\CL^{1,1}_{p+1}(B).$
  \begin{lemma}\label{L:mass-ddc-T-near-B}
  We have 
 \begin{equation*}
\sup_{n\geq 1}\sum_{j=\lowm}^\upm \int_0^\bfr {2tdt\over t^{2(k-p-j)}}\int_{\Tube(B,t)} (\ddc T_n)^\hash \wedge  \hat\beta^{k-j-p-1}\wedge  \pi^*(\omega^j)<\infty. 
\end{equation*}
  \end{lemma}
 \proof
 It follows from Corollary \ref{C:ddc-positive-finite}.
 \endproof
 
The following analogous version of  Lemma \ref{L:S(j)-widetilde-S-hol} still holds in the context of strongly  admissible maps.
 
\begin{lemma}\label{L:S(j)-widetilde-S}
 There exist   currents $S^{(0)},\ldots, S^{(k-l-p-1)}$ on $\U$ and
 a current $\widetilde S$ on $\Pi_{p+1}^{-1}(\U)$ such that  for a  suitable  subsequence $(T_{N_n})_{n=1}^\infty,$  we have
 \begin{eqnarray*}
  \lim_{N\to\infty} \big(  -\log\varphi \cdot \tau_* (\ddc T_{N_n})\wedge \alpha^m    \big)_\bullet&=&S^{(m)}\quad\text{for}\quad m=0,\ldots, k-l-p-1,\\
  \lim_{N\to\infty} \big[\Pi^\diamond_{p+1}\big(  -\log\varphi\cdot \tau_* ( \ddc T_{N_n})   \big)\big]_\bullet&=&\widetilde S.
 \end{eqnarray*}
Moreover, for  $\lowm\leq j\leq \upm$ and  $ m=1,\ldots, k-l-p-1,$
\begin{multline*}
   \lim_{n\to\infty}\int_0^r  {2tdt\over  t^{2}}\int_{\Tube(B,t)}\tau_*(\ddc T_{N_n}) \wedge  \alpha^{m-1}\wedge \beta^{k-j-p-m}\wedge \pi^*(\omega^j)\\
   =\log r\int_{\Tube(B,r)}(\ddc T)^{(m-1)} \wedge \beta^{k-j-p-m}\wedge \pi^*(\omega^j)+  
   \int_{\Tube(B,r)} S^{(m-1)} \wedge   \beta^{k-j-p-m}\wedge \pi^*(\omega^j).
   \end{multline*}

\end{lemma}
\proof 
  By  Fubini's theorem we have 
\begin{equation}\label{e:L:S(j)-widetilde-S}
\begin{split}
 &\int_0^r  {2tdt\over  t^{2}}\int_{\Tube(B,t)}\tau_*(\ddc T_{N_n}) \wedge  \alpha^{m-1}\wedge \beta^{k-j-p-m}\wedge \pi^*(\omega^j)\\
 &=\int_{y\in \Tube(B,t)}\big(\int_{|y|}^r  {2tdt\over  t^{2}}\big)\tau_*(\ddc T_{N_n}) \wedge  \alpha^{m-1}\wedge \beta^{k-j-p-m}\wedge \pi^*(\omega^j)\\
 &= \log{ r}  \int_{\Tube(B,r)}\tau_*(\ddc T_{N_n}) \wedge  \alpha^{m-1}\wedge \beta^{k-j-p-m}\wedge \pi^*(\omega^j)\\
 &+  \int_{\Tube(B,r)}(-\log {\varphi})\tau_*(\ddc T_{N_n}) \wedge  \alpha^{m-1}\wedge \beta^{k-j-p-m}\wedge \pi^*(\omega^j).
 \end{split}
 \end{equation}
  By Lemma \ref{L:mass-ddc-T-near-B}, the LHS  converges as $n\to\infty.$
By  Proposition \ref{P:existence-finite}  applied to $\ddc T_{N_n},$  the  integral in the third line of \eqref{e:L:S(j)-widetilde-S} also  converges to
$$
\log{ r}  \int_{\Tube(B,r)}\ddc \big(T_{N_n}^{(m-1)}\big)\wedge \beta^{k-j-p-m}\wedge \pi^*(\omega^j)
$$
for every $r$  except  at most a  countable  set.  Observe  that the  following  integral is  a finite  linear combination with  real coefficients times  a power of $\varphi$  of the integral on the last line   of \eqref{e:L:S(j)-widetilde-S}:
$$
 \int_{\Tube(B,r)}(-\log {\varphi})\tau_*(\ddc T_{N_n}) \wedge (\hat \alpha')^{m-1}\wedge \hat\beta^{k-j-p-m}\wedge \pi^*(\omega^j).
$$

Let $\bfj=(j_1,j_2,j_3,j_4)\in\N^4$ with 
$k-p-j_1-j_3\geq 0$ and $j_4\in\{0,1\}.$  For $0<r\leq\bfr,$ and  for a real current $T$  on $\bfU,$  consider  
\begin{equation}\label{e:I_bfj-bis}\begin{split}  I_{\bfj}(s,r)&:=\int_{\Tube(B,s,r)}\tau_*(\ddc T)\wedge \varphi^{j_2}(-\log\varphi)^{j_4}\hat\beta^{k-p-j_1-j_3}\wedge
(\pi^*\omega)^{j_3}\wedge \hat\alpha^{j_1},\\ 
I^\hash_{\bfj}(s,r)&:=\int_{\Tube(B,s,r)}(\ddc T)^\hash_{r}\wedge \varphi^{j_2}(-\log \varphi)^{j_4}\hat\beta^{k-p-j_1-j_3}\wedge
(\pi^*\omega)^{j_3}\wedge \hat\alpha^{j_1}.
\end{split}
\end{equation}
 
\begin{lemma}\label{L:spec-wedge-bis}  There is a constant $c$ independent of $T$ and $r$   such that the   following inequality holds
 \begin{equation*}
|I_\bfj(r)- I^\hash_\bfj(r)|^2 \leq c\big(\sum_{\bfj'} I^\hash_{\bfj'}(r)\big)\big ( \sum_{\bfj''} I^\hash_{\bfj''}(r) \big).  
\end{equation*}
Here, on the RHS:
\begin{itemize} \item[$\bullet$] the first sum  is taken over a finite number of multi-indices    $\bfj'=(j'_1,j'_2,j'_3,j'_4)$ as above  such that  $j'_1\leq  j_1$  and $j'_2\geq j_2$ and $j'_4\leq  j_4;$ 
\item  the second sum   is taken over  a finite number of multi-indices $\bfj''=(j''_1,j''_2,j''_3,j''_4)$ as above   such that $j''_4\leq j_4$ and that   either  ($j''_1< j_1$)
or ($j''_1=j_1$ and $j''_2\geq {1\over 4}+j_2$) or ($j''_1=j_1$ and $j''_3<j_3$).
\end{itemize}
\end{lemma}
\proof
Since  the proof is very similar to  that of  Lemma \ref{L:spec-wedge}, it is therefore left to the interested reader.
The only new thing is  the following estimate: there is a constant $c_3>0$ such that for every $1\leq \ell\leq \ell_0,$ 
\begin{equation*}
 |\tilde\tau_\ell^*(\log\varphi) - \log\varphi |\leq  c_3\varphi^{1\over 2}\qquad\text{on}\qquad \U_\ell\cap \Tube(B,\bfr).
\end{equation*}
\endproof
Applying  Lemma  \ref{L:spec-wedge-bis} yields that
\begin{multline*}
 \big |\int_{\Tube(B,\bfr)}(-\log {\varphi}) \tau_*(\ddc T_{N_n}) \wedge  (\hat\alpha')^{m-1}\wedge \hat\beta^{k-j-p-m}\wedge \pi^*(\omega^j)\\
 -\int_{\Tube(B,\bfr)}(-\log {\varphi}) (\ddc T_{N_n})^\hash \wedge  (\hat\alpha')^{m-1}\wedge \hat\beta^{k-j-p-m}\wedge \pi^*(\omega^j)\big|^2 \leq c\big(\sum_{\bfj'} I^\hash_{\bfj'}(r)\big)\big ( \sum_{\bfj''} I^\hash_{\bfj''}(r) \big)
\end{multline*}
Applying Lemma \ref{L:mass-ddc-T-near-B} yields   a constant $c>0$ independent of $r$ and $n$  such that 
\begin{eqnarray*}
 \sum_{\bfj'} I^\hash_{\bfj'}(r)&\leq& c\sum_{m,j} \int_{\Tube(B,\bfr)}(-\log {\varphi}) (\ddc T_{N_n})^\hash \wedge  (\hat\alpha')^{m-1}\wedge \hat\beta^{k-j-p-m}\wedge \pi^*(\omega^j),\\
  \sum_{\bfj''} I^\hash_{\bfj''}(r) &\lesssim& cr.
\end{eqnarray*}
Therefore, we infer that  for $\lowm\leq j\leq \upm,$
\begin{equation}\label{e:finiteness-log-ddc}
\sup_{n\geq 1}\int_{\Tube(B,\bfr)}(-\log {\varphi}) (\ddc T_{N_n})^\hash \wedge  (\hat\alpha')^{m-1}\wedge \hat\beta^{k-j-p-m}\wedge \pi^*(\omega^j)<\infty.
\end{equation}
This, combined with Propositions \ref{P:existence} and  \ref{P:existence-finite}, implies the existence of the currents  
  $S^{(0)},\ldots, S^{(k-l-p-1)}$  and $\widetilde S.$
 \endproof
 
\begin{lemma}\label{L:T(1)-S(0)}
 There is a  sequence  $(r_N)_{N=1}^\infty\subset (0,\bfr)$ with  $r_N \searrow 0$ as $N\to\infty$  such that the following identity holds
 \begin{equation*}
  \nu(T,B,\tau)=\lim_{N\to\infty} \nu(T^{(1)}+S^{(0)}, B,r_N,\id).
 \end{equation*}
\end{lemma}

\proof
 Applying    Theorem  \ref{T:Lelong-Jensen-smooth} and  Corollary \ref{C:Lelong-Jensen} to 
 the current $\tau_*T_n\wedge\beta^{k-l-p-1}\wedge \pi^*(\omega^l)$ and  for $q=1$  and  for $0<r\leq\bfr$  as in  \eqref{T(1)-S(0)-hol-Lelong-Jensen}, we can rewrite \eqref{T(1)-S(0)-hol-Lelong-Jensen} as 
\begin{equation}\label{T(1)-S(0)-Lelong-Jensen}
 \begin{split}
& \int_{\Tube(B,r)}\tau_*T_n\wedge\alpha \wedge\beta^{k-l-p-1}\wedge \pi^*(\omega^l)={1\over r^2}
\int_{\Tube(B,r)}\tau_*T_n\wedge\beta^{k-l-p}\wedge \pi^*(\omega^l)\\
 &-\int_0^r{2tdt\over t^2}\int_{\Tube(B,t)}\tau_*(\ddc T_n)\wedge\beta^{k-l-p-1}\wedge \pi^*(\omega^l)\\
 &+ {1\over r^2}\int_0^r2tdt\int_{\Tube(B,t)}\tau_*(\ddc T_n)\wedge\beta^{k-l-p-1}\wedge \pi^*(\omega^l)\\
 &-\int_0^r{2tdt\over t^2}\int_{\Tube(B,t)}\big(\ddc (\tau_*T_n)-\tau_*(\ddc T_n)\big)\wedge\beta^{k-l-p-1}\wedge \pi^*(\omega^l)\\
  &+ {1\over r^2}\int_0^r2tdt\int_{\Tube(B,t)}\big(\ddc (\tau_*T_n)-\tau_*(\ddc T_n)\big)\wedge\beta^{k-l-p-1}\wedge \pi^*(\omega^l).
 \end{split}
\end{equation}
Observe  that  by  Propositions  \ref{P:existence} and \ref{P:existence-finite}, the  LHS of \eqref{T(1)-S(0)-hol-Lelong-Jensen} converges  as  $n\to\infty$  to
$$
\int_{\Tube(B,r)} T^{(1)} \wedge\beta^{k-l-p-1}\wedge \pi^*(\omega^l)
$$
for every $r\in (0,\bfr]$ except at most a countable values of $r.$  On the other hand,  the first integral on the  RHS of \eqref{T(1)-S(0)-Lelong-Jensen} converges  as  $n\to\infty$  to
$$
{1\over r^2}
\int_{\Tube(B,r)}\tau_*T\wedge\beta^{k-l-p}\wedge \pi^*(\omega^l).
$$
Since 
$$\int_{\Tube(B,t)}(\ddc T_n)^\hash\wedge\beta^{k-l-p-1}\wedge \pi^*(\omega^l)=\int_{\Tube(B,t)}(\ddc T_n)^\hash\wedge\hat\beta^{k-l-p-1}\wedge \pi^*(\omega^l),$$
We  may apply  Lemma  \ref{L:spec-wedge} to the  RHS.  Consequently,  we infer from the  dominated convergence  theorem that
the limit of the third integral on the RHS of \eqref{T(1)-S(0)-Lelong-Jensen} is
$$
{1\over r^2}\int_0^r2tdt\int_{\Tube(B,t)}\tau_*(\ddc T)\wedge\beta^{k-l-p-1}\wedge \pi^*(\omega^l).
$$
By Lemma  \ref{L:S(j)-widetilde-S}, the  second integral on the  RHS of \eqref{T(1)-S(0)-Lelong-Jensen} converges  as  $n\to\infty$  to
$$
-\log r\int_{\Tube(B,r)}\tau_*(\ddc T)\wedge\beta^{k-l-p-1}\wedge \pi^*(\omega^l)-\int_{\Tube(B,r)}S^{(0)}\wedge\beta^{k-l-p-1}\wedge \pi^*(\omega^l).
$$
Summing up, we obtain  that
\begin{equation}\label{e:T(1)-S(0)}
\begin{split}
&{1\over r^{2(k-l-p-1)} } \int_{\Tube(B,r)}(T^{(1)}+S^{(0)})\wedge\beta^{k-l-p-1}\wedge \pi^*(\omega^l)\\
&={1\over r^{2(k-l-p)} } \int_{\Tube(B,r)}\tau_*T\wedge\beta^{k-l-p}\wedge \pi^*(\omega^l)\\
&-{\log r\over r^{2(k-l-p-1)} } \int_{\Tube(B,r)}\tau_*(\ddc T)\wedge\beta^{k-l-p-1}\wedge \pi^*(\omega^l)\\
&+{1\over r^{2(k-l-p)} }\int_0^r  2tdt \int_{\Tube(B,t)}\tau_*(\ddc T)\wedge\beta^{k-l-p-1}\wedge \pi^*(\omega^l)\\
&-\int_0^r{2tdt\over t^2}\int_{\Tube(B,t)}\big(\ddc (\tau_*T_n)-\tau_*(\ddc T_n)\big)\wedge\beta^{k-l-p-1}\wedge \pi^*(\omega^l)\\
  &+ {1\over r^2}\int_0^r2tdt\int_{\Tube(B,t)}\big(\ddc (\tau_*T_n)-\tau_*(\ddc T_n)\big)\wedge\beta^{k-l-p-1}\wedge \pi^*(\omega^l)\\
  &\equiv\sum_{j=1}^5 I_j(r).
\end{split}
\end{equation}
Clearly, by Theorem  \ref{T:Lelong-psh}  $\lim_{r\to 0} I_1(r)=\nu_l(T,B,\tau).$
We will show that there is  a  decreasing sequence $(r_N)_{n=1}^\infty\searrow 0$  through which
all  terms of  \eqref{e:T(1)-S(0)-hol}  converge  and that   $I_j(r_N)$  with $2\leq j\leq 5$  on the RHS of \eqref{e:T(1)-S(0)}  converge to $0$ as $N\to \infty.$

Consider
$$
G(t):=\sum_{j=\lowm}^\upm t^{-2(j-l)}\int_{\Tube(B,t)} (\ddc T)^\hash \wedge  \hat\beta^{k-j-p-1}\wedge  \pi^*(\omega^j). 
$$
 Lemma \ref{L:mass-ddc-T-near-B}  implies that the function
 $(0,\bfr]\ni t\mapsto  {2t\over t^{2(k-l-p)}}G(t)$ is  integrable. Hence,
 $$
  {1\over r^{2(k-l-p)}}\int_0^r 2tG(t)dt\leq  \int_0^r{2t\over t^{2(k-l-p)}}G(t)
 $$
 and the RHS converges to $0$ as $r\to 0.$ So  $\lim_{r\to 0} I_3(r)=0.$

 By Lemma \ref{L:ddc-difference},  for every $0<\epsilon <1$ there is a subset  $\bfI_\epsilon\subset(0,\bfr)$  such that $|\bfI_\epsilon\cap (r/2,r)|\geq (1-\epsilon){r\over 2}$ and
 $\lim_{r\to 0,\ r\in L_\epsilon} I_j(r)=0$ for $j\in\{4,5\}.$

    
Fix $0<\epsilon<1.$
Let  $\rho:=\liminf_{r\to 0,\ r \in\bfI_\epsilon} {-\log r\over  r^{2(k-l-p-1)}}G(r).$
We will  show that  $\rho=0.$ Suppose in order to reach a contradiction that $\rho>0.$
The finiteness of  $\int_0^r{2t\over t^{2(k-l-p)}}G(t)dt$ implies that $\int_{r\in\bfI_\epsilon} {  dr\over r |\log{  r}|} <\infty.$
On the  other hand, a straightforward
computation shows that there is a constant $c_\epsilon >0$ such that
$$
c_\epsilon \int_{r\in\bfI_\epsilon} {  dr\over r \log  r} \geq \int_{0}^\bfr {  dr\over r |\log{  r}|}=\infty. 
$$
This is a contradiction.

Since $\rho=0$ there is a  sequence  $(r_N)_{N=1}^\infty\subset \bfI_\epsilon \searrow 0$ such that   
$$\lim_{r_N\to 0} {-\log r_N\over  r_N^{2(k-l-p-1)}}G(r_N)= \lim_{r_N\to 0}I_2(t_N) =0.$$

In summary, we have  shown  that  the RHS of \eqref{e:T(1)-S(0)}, as $r=r_N$ and $N\to \infty$    converges  to    $\nu(T,B,\tau).$
So the integral on  the LHS of  \eqref{e:T(1)-S(0)} also converges to   $\nu(T,B,\tau).$ 
This completes 
 the proof. 
  \endproof
\subsection{Geometric characterizations}
Let  $T$ be  a positive plurisubharmonic  current  in the class $\SH^{3,3}_p(B)$ with $(T_n)_{n=1}^\infty$ a  sequence of approximating forms.
\begin{proposition}\label{P:widetilde-T-S}
 There exists  a subsequence $(T_{N_n})_{n=1}^\infty$ and  an open  neighborhood $\U'$ of $\overline B$ in $\E$  
 with $\U'\subset  \U$ such  that the  following  properties  hold:
 \begin{enumerate}
 \item  The following  currents  are well-defined:
 \begin{eqnarray*}
  \widetilde T&:=&\lim_{n\to\infty} \Pi^*_p  (\tau_*T_{N_n})\qquad\text{on}\qquad \Pi_p^{-1}(\U')\\
  \widetilde S&:=&\lim_{n\to\infty}\Big( \Pi^\diamond_{p+1}  \big((-\log{\varphi})\cdot \tau_*(\ddc T_{N_n})\big)\Big)_\bullet\qquad\text{on}\qquad \Pi_{p+1}^{-1}(\U').
 \end{eqnarray*}
 
 \item   For all $1\leq \ell\leq \ell_0,$  the following  currents  are well-defined:
 \begin{eqnarray*}
  \widetilde T_\ell&:= & \lim_{n\to\infty} \Pi^*_p  ((\tau_\ell)_*T_{N_n})\qquad\text{on}\qquad \Pi_p^{-1}(\U_\ell),\\
    \widetilde S_\ell &:=&  \lim_{n\to\infty}\Big( \Pi^\diamond_{p+1}  \big((-\log{\varphi})\cdot (\tau_\ell)_*(\ddc T_{N_n})\big)\Big)_\bullet\qquad\text{on}\qquad \Pi_{p+1}^{-1}(\U_\ell).
 \end{eqnarray*}
\item For all $1\leq \ell\leq \ell_0,$   the  current $\widetilde T_\ell$ (resp. $\widetilde S_\ell $)  enjoys  the cut-off property through $\Pi_p^{-1}(B)$ in $\Pi_p^{-1}(\U_\ell)$
(resp.  through $\Pi_{p+1}^{-1}(B)$ in $\Pi_{p+1}^{-1}(\U_\ell)$).

\item  The  current $\widetilde T$ (resp. $\widetilde S $)  enjoys  the cut-off property through $\Pi_p^{-1}(B)$ in $\Pi_p^{-1}(\U)$
(resp.  through $\Pi_{p+1}^{-1}(B)$ in $\Pi_{p+1}^{-1}(\U)$).  Moreover,  for all $1\leq \ell\leq \ell_0,$  the following  equalities hold   
\begin{equation*}
 \ind_{\Pi_p^{-1}(B\cap \U_\ell)} \widetilde T=  \ind_{\Pi_p^{-1}(B\cap \U_\ell)}  \widetilde T_\ell\quad\text{and}\quad
  \ind_{\Pi_{p+1}^{-1}(B\cap \U_\ell)}  \widetilde S=  \ind_{\Pi_{p+1}^{-1}(B\cap \U_\ell)} \widetilde S_\ell.
 \end{equation*}
 
 \end{enumerate}
 
\end{proposition}
\proof All the  assertions  for the currents $\widetilde T$ and $\widetilde T_\ell$ have been proved  in Proposition \ref{P:Cut-Off-bis}. So  we only need to prove these assertions  for the currents $\widetilde S$ and $\widetilde S_\ell.$

\noindent{\bf Proof  of assertion (1).} It follows from  inequality  \eqref{e:finiteness-log-ddc}.

\noindent{\bf Proof  of assertion (2).} It follows from  inequality  \eqref{e:finiteness-log-ddc}.

\noindent{\bf Proof  of assertion (3).}  For $n\geq 1$  consider the $(p+1,p+1)$-form
\begin{equation*}S_{\ell, n}:=(-\log{\varphi})\cdot (\tau_\ell)_*(\ddc T_{N_n}).
\end{equation*}
Since  $0<\varphi <1$ and  $\ddc T_n\geq  0,$ we  see that $S_{\ell, n}\geq 0.$
Arguing as in the proof of Theorem \ref{T:charac-closed-all-degrees},
  we  see that $\ddc  (S_{\ell,n})_\bullet=(\ddc  S_{\ell,n})_\bullet.$
So we  infer as in \eqref{e:ddc-R_n} that
  \begin{equation*}
  \begin{split} \ddc \Pi^\diamond_{p+1}( S_{\ell,n})&= -\Pi^\diamond_{p+1}\big(\ddc ((\tau_\ell)_*T_{N_n}) \wedge  \alpha\big)= - \Pi^\diamond_{p+1}\big((\tau_\ell)_*(\ddc T_{N_n}) \wedge (\hat \alpha'-c_1\pi^*\omega)\big) \\
  &=  -\Pi^\diamond_{p+1}\big((\tau_\ell)_*(\ddc T_{N_n}) \wedge (\hat \alpha') \big)+   c_1 \Pi^\diamond_{p+1}\big(\ddc ((\tau_\ell)_*T_{N_n}) \wedge \pi^*(\omega)\big).
  \end{split}\end{equation*}
  Since  by passing to a  subsequence if necessary, the weak limits of both terms on the last line are positive closed currents,  we infer that  both  $\widetilde S_\ell$ and  $\ddc \widetilde S_\ell$ are  currents of order $0$.  Hence, by Theorem \ref{T:Bassanelli}, $ \widetilde S_\ell$ enjoy the cut-off property  through $\Pi_{p+1}^{-1}(B)$ in $\Pi_{p+1}^{-1}(\U_\ell).$

\noindent{\bf Proof  of assertion (4).} 
 Fix $1\leq \ell\leq\ell_0.$
Let  $\Phi$ be a continuous  test form of dimension $2p+2$ on $\X_{p+1}$ which is compactly supported on $\Pi_{p+1}^{-1}(\U_\ell).$
Let $0<r\leq\bfr.$
Write
\begin{eqnarray*}
\langle \widetilde S-  \widetilde S_\ell ,\Phi \rangle_{\Tube(B,r)}&=&\lim_{n\to\infty} \langle (\Pi_{p+1})^\diamond((-\log{\varphi})(\tau_*(\ddc T_n)  -(\tau_\ell)_* (\ddc T_n))), \Phi\rangle_{\Tube(B,r)}\\
&=&  \langle (\tau_\ell)_*(\ddc T_n), (-\log{\varphi})(\tilde\tau^*_\ell-\id)((\Pi_{p+1})_\diamond(\Phi))\rangle_{\Tube(B,r)}=I_r.
\end{eqnarray*}
Applying Lemma \ref{L:pushforward_Pi_p-bis} to $I_r$ and using  the inequality $|\log\varphi|\varphi^{1\over 2}\lesssim \varphi^{1\over 4}$ we infer    that
\begin{equation*}
 I_r \leq  c r^{1\over 2}\sum_{0\leq j\leq l,\  0\leq q\leq k-l-p} \int  (\tau_\ell)_*(T_n) \wedge  \pi^*(\omega^j)\wedge \alpha_\ver^q\wedge \beta_\ver^{k-p-j-q}
\end{equation*} 
By Proposition  \ref{P:existence-finite} and  \ref{P:existence}, $I_r$  uniformly bounded by $cr$  independent of $ n.$
So   $\lim_{r\to 0} \langle  \widetilde S-  \widetilde S_\ell  ,\Phi \rangle_{\Tube(B,r)}=0.$  This proves assertion (4).
\endproof

\begin{proposition}\label{P:widehat-T}
    There exists  a subsequence $(T_{N_n})_{n=1}^\infty$ and  an open  neighborhood $\U'$ of $\overline B$ in $\E$  
 with $\U'\subset  \U$ such  that the  following  properties  hold:
\begin{enumerate}
 \item  The following  currents  are well-defined:
\begin{equation*}
\widehat T:=\lim_{n\to\infty}  \big(  \Pi^\diamond_{p+1}(\tau_*T_{N_n}\wedge \alpha)  \big)_\bullet\qquad\text{on}\qquad  \Pi_{p+1}^{-1}(\U').
\end{equation*}

 \item   For all $1\leq \ell\leq \ell_0,$  the following  currents  are well-defined:
  \begin{equation*}
\widehat T_\ell:=\lim_{n\to\infty}  \big(  \Pi^\diamond_{p+1}((\tau_\ell)_*T_{N_n}\wedge \alpha)  \big)_\bullet\qquad\text{on}\qquad  \Pi_{p+1}^{-1}(\U_\ell).
\end{equation*}
 
 \item For all $1\leq \ell\leq \ell_0,$   the  current $\widehat T_\ell$    enjoys  the cut-off property     through $\Pi_{p+1}^{-1}(B)$ in $\Pi_{p+1}^{-1}(\U_\ell).$

 \item  The  current $\widehat T$    enjoys  the cut-off property   through $\Pi_{p+1}^{-1}(B)$ in $\Pi_{p+1}^{-1}(\U).$  Moreover,  for all $1\leq \ell\leq \ell_0,$  the following  equalities hold   
\begin{equation*}
  \ind_{\Pi_{p+1}^{-1}(B\cap \U_\ell)}  \widehat T=  \ind_{\Pi_{p+1}^{-1}(B\cap \U_\ell)} \widehat T_\ell.
 \end{equation*}

 \end{enumerate}
\end{proposition}
\proof
 \noindent{\bf Proof  of assertion (1).}  We argue as in the proof of Proposition \ref{P:Cut-Off-bis} (1).
Pick $1\leq \ell\leq\ell_0.$
Let  $\Phi$ be a continuous  test form of dimension $2p+2$ on $\X_p$ which is compactly supported on $\Pi_{p+1}^{-1}(\U_\ell).$
Write
\begin{eqnarray*}
\langle (\Pi_{p+1})^\diamond(\tau_*T_{N_n}\wedge \alpha), \Phi\rangle&=&\langle\tau_*(T_{N_n}), \alpha\wedge (\Pi_{p+1})_\diamond(\Phi)\rangle=\langle (\tau_\ell)_*( T_n), \tilde\tau_\ell^*(\alpha\wedge (\Pi_{p+1})_\diamond(\Phi))\rangle\\
&=& \langle (\tau_\ell)_*(T_n), \alpha\wedge (\Pi_{p+1})_\diamond(\Phi^\sharp)\rangle+\langle (\tau_\ell)_*(T_n), (\tilde\tau^*_\ell-\id)(\alpha\wedge (\Pi_{p+1})_\diamond(\Phi^\sharp))\rangle\\
&+&\langle (\tau_\ell)_*( T_n), \tilde\tau^*_\ell(\alpha\wedge (\Pi_{p+1})_\diamond(\Phi-\Phi^\sharp))\rangle.
\end{eqnarray*}
Applying Lemma 
\ref{L:pushforward_Pi_p} to  $\Phi^\sharp$ yields that
 $$ |\langle (\tau_\ell)_*(T_n), \alpha\wedge (\Pi_{p+1})_\diamond(\Phi^\sharp)\rangle|\leq  c \sum_{0\leq j\leq l,\  0\leq q\leq k-l-p} \int  (\tau_\ell)_*(T_n) \wedge  \pi^*(\omega^j)\wedge \alpha_\ver^q\wedge \beta_\ver^{k-p-j-q}.$$
By Proposition  \ref{P:existence-finite} and  \ref{P:existence}, the  RHS is  uniformly bounded independent of $ n.$
So is  $ |\langle (\tau_\ell)_*(T_n), \alpha\wedge (\Pi_{p+1})_\diamond(\Phi^\sharp)\rangle|.$

Set $\Psi':= (\tilde\tau^*_\ell-\id)(\alpha\wedge (\Pi_{p+1})_\diamond(\Phi^\sharp)).$
Applying Lemma \ref{L:tilde-tau-dzeta_j} and Lemma \ref{L:difference} yields that
\begin{equation*}
 \langle (\tau_\ell)_*(T_n), \Psi'\rangle= \langle (\tau_\ell)_*(T_n), (\Psi')^\sharp\rangle \leq  c \sum_{0\leq j\leq l,\  0\leq q\leq k-l-p} \int  (\tau_\ell)_*(T_n) \wedge  \pi^*(\omega^j)\wedge \alpha_\ver^q\wedge \beta_\ver^{k-p-j-q}
\end{equation*}
By Proposition  \ref{P:existence-finite} and  \ref{P:existence}, the  RHS is  uniformly bounded independent of $ n.$
So is  $\big|\langle (\tau_\ell)_*(T_n), (\tilde\tau^*_\ell-\id)(\alpha\wedge (\Pi_{p+1})_\diamond(\Phi^\sharp))\rangle\big|.$

Set $\Psi:=[\tilde\tau^*_\ell(\alpha\wedge (\Pi_{p+1})_\diamond(\Phi-\Phi^\sharp))]^\sharp.$
Applying Lemma \ref{L:pushforward_Pi_p-bis}  yields that
\begin{multline*}
 \langle (\tau_\ell)_*( T_n), \tilde\tau^*_\ell(\alpha\wedge (\Pi_{p+1})_\diamond(\Phi-\Phi^\sharp))\rangle = \langle (\tau_\ell)_*(T_n), \Psi\rangle \\
 \leq  c r\sum_{0\leq j\leq l,\  0\leq q\leq k-l-p} \int  (\tau_\ell)_*(T_n) \wedge  \pi^*(\omega^j)\wedge \alpha_\ver^q\wedge \beta_\ver^{k-p-j-q}
\end{multline*}
By Proposition  \ref{P:existence-finite} and  \ref{P:existence}, the  RHS is  uniformly bounded by $cr$  independent of $ n.$
So is $\big|\langle (\tau_\ell)_*( T_n), \tilde\tau^*_\ell(\alpha\wedge (\Pi_{p+1})_\diamond(\Phi-\Phi^\sharp))\rangle\big|.$

Putting together the above  three estimates, we get  $| \langle (\Pi_{p+1})^\diamond(\tau_*T_{N_n}\wedge \alpha), \Phi\rangle|\leq  c\|\Phi\|_{\Cc^0},$
for a constant $c>0$ independent of $\Phi.$  This proves assertion (1).

\noindent{\bf Proof  of assertion (3).}  Writing
 \begin{equation*}
\widehat T_\ell:=\lim_{n\to\infty}  \big(  \Pi^\diamond_{p+1}((\tau_\ell)_*T_{N_n}\wedge \hat\alpha')  \big)_\bullet -c_1\lim_{n\to\infty}  \big(  \Pi^\diamond_{p+1}(\tau_*T_{N_n}\wedge \pi^*\omega)  \big),
\end{equation*}
we  see that  $\widehat T_\ell$ is  the difference of two positive closed currents. Hence, assertion (3) follows.

\noindent{\bf Proof  of assertion (4).} 
 Fix $1\leq \ell\leq\ell_0.$
Let  $\Phi$ be a continuous  test form of dimension $2p+2$ on $\X_{p+1}$ which is compactly supported on $\Pi_{p+1}^{-1}(\U_\ell).$
Let $0<r\leq\bfr.$
Write
\begin{eqnarray*}
\langle \widehat T-  \widehat S_\ell ,\Phi \rangle_{\Tube(B,r)}&=&\lim_{n\to\infty} \langle (\Pi_{p+1})^\diamond\big((\tau_*( T_n)  -(\tau_\ell)_* ( T_n)) \wedge \alpha\big)  , \Phi\rangle_{\Tube(B,r)}\\
&=&  \langle (\tau_\ell)_*(T_n),(\tilde\tau^*_\ell-\id)( \alpha\wedge (\Pi_{p+1})_\diamond(\Phi))\rangle_{\Tube(B,r)}=I_r.
\end{eqnarray*}
Applying Lemma \ref{L:pushforward_Pi_p-bis} to $I_r$   we infer    that
\begin{equation*}
 I_r \leq  c r^{1\over 2}\sum_{0\leq j\leq l,\  0\leq q\leq k-l-p} \int  (\tau_\ell)_*(T_n) \wedge  \pi^*(\omega^j)\wedge \alpha_\ver^q\wedge \beta_\ver^{k-p-j-q}
\end{equation*} 
By Proposition  \ref{P:existence-finite} and  \ref{P:existence}, $I_r$  uniformly bounded by $cr$  independent of $ n.$
So   $\lim_{r\to 0} \langle  \widehat T-  \widehat T_\ell  ,\Phi \rangle_{\Tube(B,r)}=0.$  This proves assertion (4).
\endproof

The following central result of the section provides a geometric characterization of the  top Lelong number in the case of strongly admissible maps. It  should be  compared with Theorem \ref{T:charac-Lelong-for-psh-hol}.

\begin{theorem}\label{T:charac-Lelong-for-psh}
 Let  $\widetilde  T$ and $\widetilde S$ be the currents defined by Proposition  \ref{P:widetilde-T-S}.
There exist two  functions $f,g\in L^1_\loc(B)$   such that
$$
\ind_{\Pi_p^{-1}(B)}=(f\circ \Pi_p)[ \Pi_p^{-1}(B)]\qquad\text{and}\qquad \ind_{\Pi_{p+1}^{-1}(B)}=(g\circ \Pi_p)[ \Pi_{p+1}^{-1}(B)].
$$
Moreover, both   function  $f$ and $f+g$ are  non-negative and  $f$ is    plurisubharmonic on $B$  and  $f+g$ is the difference of two  plurisubharmonic functions on $B$ and 
\begin{equation*}
 \nu(T,B,\tau)=\int_B (f+g)\omega^l.
\end{equation*}

\end{theorem}
\proof 

Summing up, we have shown that $\widehat T$ and $\widetilde S$  are both $\C$-normal.
By  Proposition 
\ref{P:Cut-Off}, there  exist  non-negative functions $\hat f,$ $\tilde g\in L^1_\loc(\Pi_{p+1}^{-1}(B))$ 
such that
\begin{equation}\label{e:charac-Lelong-for-psh(2)}
 \ind_{\Pi_{p+1}^{-1}(B)}\widehat T=\hat f [\Pi_{p+1}^{-1}(B)]\qquad\text{and}\qquad \ind_{\Pi_{p+1}^{-1}(B)}\widetilde  S=\tilde g [\Pi_{p+1}^{-1}(B)].
\end{equation}
Moreover,  by  Theorem  \ref{T:Bassanelli-psh} $\hat f$ and $\hat f+\hat g$ are  positive plurisubharmonic. Therefore,  they are constant on fibers.
So there are   functions $ f_0,$ $ g\in L^1_\loc(B)$  such that
\begin{equation}\label{e:hatf-f_tildeg-g-bis}\hat f=f_0\circ \Pi_{p+1}\qquad \text{and}\qquad   \tilde  g=g\circ \Pi_{p+1}\qquad \text{on}\qquad\Pi_{p+1}^{-1}(B).
\end{equation}
By  Lemma \ref{L:T(1)-S(0)},  we have 
\begin{multline*}
 \nu(T,B,\tau)=\lim_{N\to\infty}\nu(T^{(1)}+S^{(0)},B,r_N,\tau)\\
 =\lim_{N\to \infty}\lim_{n\to\infty} {1\over  r_N^{2(k-l-p-1)}}\int_{\Tube(B,r_N)} \big(\tau_*T_{n}\wedge\alpha-(\log\varphi)\tau_*(\ddc T_n) \big)\wedge\beta^{k-l-p-1}\wedge \pi^*(\omega^l).
\end{multline*}
By  Lemma 
\ref{L:T_n-change-alpha-to-beta-hol}, the last line  is  equal to
\begin{equation*}
\lim_{N\to \infty}\lim_{n\to\infty} \int_{\Tube(B,r_N)} \big(\tau_*T_n\wedge\alpha-(\log\varphi)\tau_*(\ddc T_n) \big)\wedge\alpha^{k-l-p-1}\wedge \pi^*(\omega^l).
\end{equation*}
Since $\alpha^{k-l-p-1}\wedge \pi^*(\omega^l)=\alpha_\ver^{k-l-p-1}\wedge \pi^*(\omega^l),$ the last line  is  equal  to 
\begin{equation*}
\lim_{N\to \infty}\lim_{n\to\infty} \int_{\Tube(B,r_N)} \big(\tau_*T_n\wedge\alpha-(\log\varphi)\tau_*(\ddc T_n) \big)\wedge\alpha_\ver^{k-l-p-1}\wedge \pi^*(\omega^l)
\end{equation*}
 By assertion (1), this is  equal to
\begin{equation*}
\lim_{N\to \infty}  \int_{\Pi_{p+1}^{-1}(\Tube(B,r_N))} (\widehat T+\widetilde S )\wedge\Pi_{p+1}^\diamond\big (\alpha_\ver^{k-l-p-1}\wedge \pi^*(\omega^l)\big).
\end{equation*}
By  Lemma \ref{L:Upsilon_j} and equality \eqref{e:hatf-f_tildeg-g-bis},  the last  expression is  equal to
\begin{equation*}
 \lim_{N\to \infty}  \int_{\Pi_{p+1}^{-1}(\Tube(B,r_N))} (\widehat T+\widetilde S )\wedge\Pr\nolimits_{p+1}^\diamond(\Upsilon^{(p+1)(k-l-p-1)})\wedge \Pi_{p+1}^\diamond\big(\pi^*(\omega^l)\big)
 =\int_B (f_0+g)\omega^l,
\end{equation*}
where the  equality  follows from \eqref{e:charac-Lelong-for-psh(2)} and  Proposition \ref{P:Cut-Off}.

On the other  hand,
by Proposition \ref{P:Cut-Off-hol} for $j=l$ and hence $\hat j=p$, there is a function  $f\in L^1_\loc(B)$   such that 
$$\ind_{\Pi_p^{-1}(B)}\widetilde T=(f\circ \Pi_p)[ \Pi_p^{-1}(B)]$$
and that for every  $\Cc^2$-piecewise smooth subdomain  $D\subset B,$ 
\begin{equation*}
 \lim_{N\to \infty}\lim_{n\to\infty} \int_{\Tube(D,r)} \tau_*T_n\wedge\alpha^{k-l-p}\wedge \pi^*(\omega^l)=\|\widetilde T\|(\Pi^{-1}_p(D))
 =\int_D f\omega^l.
\end{equation*}
Observe that the expression on the LHS is  also equal to $\|\widehat T\|(\Pi^{-1}_{p+1}(D))=\int_D f_0\omega^l.$
So   $\int_D (f_0-f)\omega^l=0.$ Since this equality holds for every  $\Cc^2$-piecewise smooth subdomain  $D\subset B,$ we  infer that $f_0=f.$
The proof of assertion (2) is thereby completed.
\endproof

\section{The top  Lelong number is totally intrinsic}
\label{S:Top-Lelong-strong-intrinsic}

We  keep the notation  introduced in   Sections  \ref{S:Intro} and  the  Standing Hypothesis  introduced in  Subsection \ref{SS:Global-setting}.
We are in the position to state the second collection of main results.
The first  two results are devoted to  positive  closed currents  for  strongly  admissible maps and for holomorphic  admissible maps.
\begin{theorem}\label{T:tot-intrinsic-closed}
Suppose that  one of  the  following two conditions is  fulfilled:
\begin{enumerate}
 \item  $T$ is a positive closed current  in the class $\CL^{2,2}_p(B)$ and $\ddc \omega^j=0$ for $1\leq j\leq \upm-1;$
 \item  $T$ is a positive closed current  in the class $\CL^{1,1}_p(B)$ and $\omega$ is  K\"ahler.
\end{enumerate}
Then the top Lelong number of $T$ along $B$   is  totally  intrinsic, that is,
$\nu(T,B,\tau,h)$ is independent of the choice of  a  strongly admissible  map $\tau$ and a Hermitian metric $h$ on $\E.$ 
\end{theorem}

\begin{theorem}\label{T:tot-intrinsic-closed-hol}
Suppose that  one of  the  following two conditions is  fulfilled:
\begin{enumerate}
 \item  $T$ is a positive closed current  in the class $\CL^{2}_p(B)$ and $\ddc \omega^j=0$ for $1\leq j\leq \upm-1;$
 \item  $T$ is a positive closed current  in the class $\CL^{1}_p(B)$ and $\omega$ is  K\"ahler.
\end{enumerate}
Then the top Lelong number of $T$ along $B$   is  totally  intrinsic, that is,
$\nu(T,B,\tau,h)$ is independent of the choice of  a  holomorphic admissible  map $\tau$ and a Hermitian metric $h$ on $\E.$ 
\end{theorem}

The next two results deal with  positive  pluriharmonic currents  for  strongly  admissible maps and for holomorphic  admissible maps.

\begin{theorem}\label{T:tot-intrinsic-ph}
Assume that $\omega$ is  K\"ahler and
 $T$ is a positive pluriharmonic current  in the class $\PH^{2,2}_p(B).$
Then the Lelong number of $T$ along $B$ is  totally  intrinsic, that is,
$\nu(T,B,\tau,h)$ is independent of the choice of  a  strongly  admissible  map $\tau$ and a Hermitian metric $h$ on $\E.$ 
\end{theorem}

\begin{theorem}\label{T:tot-intrinsic-ph-hol}
Assume that $\omega$ is  K\"ahler and
 $T$ is a positive pluriharmonic current  in the class $\PH^{2}_p(B).$
Then the Lelong number of $T$ along $B$ is  totally  intrinsic, that is,
$\nu(T,B,\tau,h)$ is independent of the choice of  a  holomorphic  admissible  map $\tau$ and a Hermitian metric $h$ on $\E.$ 
\end{theorem}

The last two results of the  section  discuss  positive  plurisubharmonic currents  for  strongly  admissible maps and for holomorphic  admissible maps.
\begin{theorem}\label{T:tot-intrinsic-psh}
Assume that $\omega$ is  K\"ahler and
 $T$ is a positive plurisubharmonic current  in the class $\SH^{3,3}_p(B).$
Then the Lelong number of $T$ along $B$ is  totally  intrinsic, that is,
$\nu(T,B,\tau,h)$ is independent of the choice of  a strongly  admissible  map $\tau$ and a Hermitian metric $h$ on $\E.$ 
\end{theorem}

\begin{theorem}\label{T:tot-intrinsic-psh-hol}
Assume that $\omega$ is  K\"ahler and
 $T$ is a positive plurisubharmonic current  in the class $\SH^{3}_p(B).$
Then the Lelong number of $T$ along $B$ is  totally  intrinsic, that is,
$\nu(T,B,\tau,h)$ is independent of the choice of  a holomorphic admissible  map $\tau$ and a Hermitian metric $h$ on $\E.$ 
\end{theorem}

Let  $h$ and $h'$ be  two  Hermitian metrics  on $\E.$  Let $\varphi$ (resp. $\varphi'$) be  the      function given  by \ref{e:varphi} 
corresponding  to the metric $\|\cdot\|:=h$  (resp.  $\|\cdot\|:=h'$).  Fix  a strongly  admissible  map $\tau.$
Since the Lelong numbers are intrinsic,  
we only need  to show that
\begin{equation}\label{e:tot-intrinsic-psh(1)}
 \nu(T,B,\tau,h)=\nu(T,B,\tau,h').
\end{equation}
We may assume  without loss of generality that $p<k-l.$
Let $(T_n)_{n=1}^\infty$ be a  sequence of approximating forms  for  $T$ as  an element of  the class $\SH^{3,3}_p(B).$
 We may assume  without loss of generality that  $(T_n)_{n=1}^\infty\subset\widetilde \SH^{3,3}_p(\bfU,\bfW).$
 Fix a  small open  neighborhood $\U'$ of $\overline B$ in $\E$   such that 
  $\U'\Subset  \U.$
By Theorem  \ref{T:charac-Lelong-for-psh-hol},
 we can assume that the following  currents  are well-defined:
 \begin{eqnarray*}
  \widetilde T&:=&\lim_{n\to\infty} \Pi^*_p  (\tau_*T_{n})\qquad\text{on}\qquad \Pi_p^{-1}(\U')\\
  \widetilde S&:=&\lim_{n\to\infty} \Pi^\diamond_{p+1}  \big(-(\log{\varphi}) (\tau_*(\ddc T_{n})\big)_\bullet\qquad\text{on}\qquad \Pi_{p+1}^{-1}(\U').
 \end{eqnarray*}
and 
\begin{equation}\label{e:tot-intrinsic-psh(2)}
 \nu(T,B,\tau,h)=c_{\widetilde T}+c_{\widetilde S},
\end{equation}
where
\begin{equation*}
 c_{\widetilde T}:=\|\widetilde T\|(\Pi^{-1}_p(B)) 
 \qquad\text{and}\qquad
 c_{\widetilde S}:=\|\widetilde S\|(\Pi^{-1}_p(B)).
\end{equation*}

Let $\pi:\ \GL_{k-l}(\E,\C)\to  V$ be the canonical  holomorphic  projection whose fiber over $x\in V$ is $\GL(\E_x,\C),$ the general linear  group of degree  $k-l$ over $\C.$
 
\begin{lemma}\label{L:F}
 There exists a   smooth map $F:\  \pi^{-1}(\overline B)  \to  \pi^{-1}(\overline B)$  
 such that 
 \begin{enumerate}
 \item   for every $x\in\overline B,$  $F$ sends $\E_x$ onto $\E_x$ and  $F|_{\E_x}$  is  $\C$-linear, in other words,
 $F$  is  a  section over $\overline B$  of    the  projection  $\pi:\ \GL_{k-l}(\E,\C)\to V;$
 \item 
 $
 \|y\|_{h'}=\|F(y)\|_h  $ for $y\in \E,$  in other words,  $\varphi\circ F=\varphi'$ on $\E.$
 \end{enumerate}
\end{lemma}
\proof
Fix a  point $x\in B$ and  consider the  inner product  given by $\varphi|_{\E_x}.$  Then the matrix $A'(x)$ of the inner product $\varphi'|_{\E_x}$ can be written as  $U(x)\circ F(x),$ where $F(x)\in\GL(\E_x,\C)$ is a diagonal matrix   and $U(x)\in\U(k-l)$  is a unitary
matrix. Hence,   we infer that $\varphi\circ F=\varphi'$ on $\E_x.$ 

Since $\varphi$ and $\varphi'$  are smooth,  we obtain   the decomposition $A'(x)=U(x)\circ F(x)$ for all $x\in\overline B.$  Therefore,  we can construct a smooth map $F:\  \pi^{-1}(\Omega)  \to  \pi^{-1}(\Omega)$ satisfying the conclusion of the lemma.
\endproof

\begin{lemma}\label{L:F-bis}  \begin{enumerate} \item The map $F$  defined by Lemma \ref{L:F} induces a diffeomorphic map $F_p:\  \X_p\to \X_p$  such that
\begin{equation*}
\Pi_p\circ F_p=F\circ \Pi_p\qquad\text{and}\qquad \Pi_{p}\circ F^{-1}_{p}=F^{-1}\circ \Pi_{p} 
\end{equation*}
and  that $F_p(\Pi_p^{-1}(B))= \Pi_p^{-1}(B).$
\item  Let $S$ be  a  current  of bidegree $(p,p)$ on $\X_p$ which enjoys the cut-off  property through $\Pi_p^{-1}(B).$
Suppose  that $\ind_{\Pi_p^{-1}(B)} S=f[\Pi_p^{-1}(B)],$ where $f$ is a function  on $\Pi_p^{-1}(B).$
By Lemma \ref{L:Cut-Off} applied to $F_p$, $(F_p)_*S$ also  enjoys the cut-off  property through $\Pi_p^{-1}(B).$
Then 
$$  \int_{\Pi_p^{-1}(B)}\ind_{\Pi_p^{-1}(B)} S= \int_{\Pi_p^{-1}(B)}\ind_{\Pi_p^{-1}(B)} (F_p)_*S .$$
 \end{enumerate}
\end{lemma}
\proof

To prove assertion (1), 
fix  a point $x\in V.$ Choose   a  system of   coordinates $z=(z_1,\ldots,z_{k-l})$ on $\E_x$  so that the  hyperplane $\{z_1=0\}$ is invariant  by the  $\C$-linear map  $F(x).$
Add to the  coordinates $z$ 
the  coordinates $w=(w_1,\ldots,w_k)$ so that
$(z,w)$ is a local coordinate  around $x.$
Let  $H$ be an element of $\G_p(\E_x).$ 
We may assume without loss of generality that  $H_0:=H\cap\{z_1=0\}$ is  a linear subspace of dimension $p-1.$
So $H_0$ defines an element in  $\G_{p-1}(\E_x).$
We may assume without loss of generality that
$$
H_0:=  \left\lbrace z_1=\cdots= z_{k-l-p+1}=0\right\rbrace.
$$
For $z=(z_1,\ldots,z_{k-l}),$  write  $z^{('p)}=(z_1,\ldots,z_{k-l-p+1})\in\C^{k-l-p+1}.$
 Recall  from \eqref{e:X_p,H_0} and \eqref{e:blow-up-H_0} that    $\X_{p,H_0}$ defined in  is the closure of  $\X'_{p,H_0}$ in $\C^{k-l-p+1}\times \P^{k-l-p},$  and 
\begin{equation*} 
 \X'_{p,H_0}\simeq  \{(z^{('p)},[z^{('p)}]):\  z^{('p)}\in\C^{k-l-p+1}\setminus\{0\}\}\quad\text{and}\quad
\Pi_{p,H_0}(z,H)=z^{('p)}.
\end{equation*}
Set   $H_1:=F(x)(H_0).$ This  is  a linear subspace of dimension $p-1$ of the  hyperplane $\{z_1=0\}.$ 
Consider the map $F_p:\  \X_p\setminus \Pi_p^{-1}(V)\to \X_p\setminus \Pi_p^{-1}(V)$  given by
$$F_p((z^{('p)},[z^{('p)}])):=(F(x)z^{('p)},[F(x)z^{('p)}])\in  \X'_{p,H_1}.$$
We extends the map continuously  through  $\Pi_p^{-1}(V)$ in order to obtain a  continuous map  on $\X_p.$   
Using this explicit  formula, we can check that this  map  satisfies the conclusion of assertion (1).

 It follows from assertion (1) that $\ind_{\Pi_p^{-1}(B)} (F_p)_*S= (F_p)_*\big( \ind_{\Pi_p^{-1}(B)}S\big).$
 Hence, since $F_p|_{\Pi_p^{-1}(B)}$ is  diffeomorphic, we infer that 
\begin{eqnarray*}  \int_{\Pi_p^{-1}(B)}\ind_{\Pi_p^{-1}(B)} S&=&\int_{\Pi_p^{-1}(B)} f\cdot \Pi_p^*(\omega^l)\wedge \Upsilon_p^{k-l-p}\\&=& \int_{\Pi_p^{-1}(B)} (F_p)_*\big( f\cdot \Pi_p^*(\omega^l)\wedge \Upsilon_p^{k-l-p}\big)  =\int_{\Pi_p^{-1}(B)}\ind_{\Pi_p^{-1}(B)} (F_p)_*S .
 \end{eqnarray*}
 This proves assertion (2).

\endproof
Let $T'_n$ be currents  on $\bfU$  defined  by $\tau_*T'_n:=F_*(\tau_*T_n).$  Similarly,  let   $T'$ be currents  on $\bfU$  defined  by  $\tau_*T'=F_*(\tau_*T).$
Define 
\begin{eqnarray*}
  \widetilde T'&:=&\lim_{n\to\infty} \Pi^*_p  (\tau_*T'_{n})\qquad\text{on}\qquad \Pi_p^{-1}(\U')\\
  \widetilde S'&:=&\lim_{n\to\infty} \Pi^\diamond_{p+1}  \big(-(\log{\varphi})\tau_*(\ddc T'_{n})\big)_\bullet\qquad\text{on}\qquad \Pi_{p+1}^{-1}(\U').
 \end{eqnarray*}
We  deduce from Lemma \ref{L:F-bis} that 
\begin{equation}\label{e:tot-intrinsic-psh(3)}
 \nu(T,B,\tau,h')=\nu(T',B,\tau,h)=c_{\widetilde T'}+c_{\widetilde S'},
\end{equation}
where
\begin{equation*}
 c_{\widetilde T'}:=\|\widetilde T'\|(\Pi^{-1}_p(B))  
 \qquad\text{and}\qquad 
 c_{\widetilde S'}:=\|\widetilde S'\|(\Pi^{-1}_p(B)).
\end{equation*}
We infer  from  \eqref{e:tot-intrinsic-psh(2)} and \eqref{e:tot-intrinsic-psh(3)} that  in order to prove \eqref{e:tot-intrinsic-psh(1)},
it suffices to show that
$c_{\widetilde T}=c_{\widetilde T'}$ and $c_{\widetilde S}=c_{\widetilde S'}.$

By Lemma \ref{L:F-bis},  $F$ induces a  map $F_p:\  \X_p\to \X_p$ and  $F_{p+1}:\  \X_{p+1}\to \X_{p+1}$  such that
\begin{eqnarray*}
\Pi_p\circ F_p&=&F\circ \Pi_p\qquad\text{and}\qquad \Pi_{p+1}\circ F_{p+1}=F\circ \Pi_{p+1},\\
\Pi_p\circ F^{-1}_p&=&F^{-1}\circ \Pi_p\qquad\text{and}\qquad \Pi_{p+1}\circ F^{-1}_{p+1}=F^{-1}\circ \Pi_{p+1}.
\end{eqnarray*}
Therefore,  we infer that
$$
(F_p)_*( \Pi_p^*(\tau_*T_n))=\Pi_p^*(\tau_* T'_n).
$$
So $(F_p)_*\widetilde T=\widetilde T'.$ Hence, $c_{\widetilde T}=c_{\widetilde T'}.$

On the  other hand, 
we have
$$
(-\log{\varphi\circ F})  \tau_*(\ddc T_n)=F^*((-\log{\varphi})\tau_* (\ddc T'_n)).
$$
Hence,
$$
F_{p+1}^*(\widetilde S')=\lim_{n\to\infty} \big[ \Pi^\diamond_{p+1} ((-\log{\varphi\circ F})\tau_*(\ddc  T_n))   \big]_\bullet.
$$

Now  for a suitable   constant $c>1$  we have  $c^{-1}\varphi(y)\leq\varphi\circ F(y)\leq c\varphi(y)$ for $y\in \E.$
Thus,
\begin{multline*}
\Big(\log c \cdot( \tau_*(\ddc T_n)-(\log \varphi) \tau_*(\ddc T_n)\Big)^\sharp\geq\Big( -(\log \varphi\circ F)\tau_*(\ddc T_n)\Big)^\sharp\\
\geq \Big(-\log c \cdot( \tau_*(\ddc T_n))-(\log \varphi ) \tau_*(\ddc T_n)\Big)^\sharp.
\end{multline*}
This implies that 
$$
\Big(\log c \big(\Pi_{p+1}^*(  \tau_*(\ddc T))\big)+\widetilde S\Big)^\sharp\geq \big (F_{p+1}^*\widetilde S'\big)^\sharp\geq \Big(-\log c \big(\Pi_{p+1}^*( \tau_*(\ddc T))\big)+\widetilde S\Big)^\sharp.
$$
On the other hand, by Theorem \ref{T:vanishing-Lelong}, $\nu_\top(\ddc T,B,\tau)=0.$
Consequently, we  infer from  the  last two  estimates that   $c_{\widetilde S}=c_{ F^*_{p+1}\widetilde S'}.$
By Lemma \ref{L:F-bis}, we have   $c_{\widetilde S}=c_{ F^*_{p+1}\widetilde S'}.$
Thus,  $c_{\widetilde S}=c_{\widetilde S'}.$
This completes the proof.


\section{Proof of the main general theorems and  concluding remarks}
\label{S:Proofs}

\subsection{Proofs of the main general results}\label{SS:Proofs-Main-General-Results}
 Recall  that  $X$ is a complex manifold of dimension $k$
  and $V\subset X$ is a submanifold of dimension $1\leq l<k.$
  Fix $0\leq p\leq k$ and  define $\lowm$ and $\upm$  by  \eqref{e:m}.   
The  vector bundle $\E$ (that is, the normal  bundle  to $V$ in $X$) is  endowed with a Hermitian  metric $h.$
$V$ is  endowed with a Hermitian  metric $\omega.$ Let $B$ be   a piecewise $\Cc^2$-smooth open subset  of $V$ and that  there exists a strongly   admissible map for $B.$

 \proof[Proof of Theorem \ref{T:main_1} {\rm (Tangent Theorem I) }]
  Let $X,$ $V,$ $B$  be as  above  and  suppose that $(V,\omega)$ is  K\"ahler.
  Let $T$ be  a   positive plurisubharmonic  $(p,p)$-current  on a neighborhood of $\overline B$ in $X$ such that  $T=T^+-T^-$  for some $T^\pm\in\SH_p^{3,3}( B).$
  
  Assertion (1) is proved in  Theorem \ref{T:Lelong-psh} (1).
  
   Assertion (2) is proved in  Theorem \ref{T:Lelong-psh} (5).
   
    Assertion (3) is proved in  Theorem \ref{T:Lelong-psh} (3).
    
     The non-negativity of the top Lelong number $\nu_\upm(T,B,\omega,h)$ stated in assertion (4) is proved in  Theorem \ref{T:Lelong-psh} (4).
By Theorem \ref{T:tot-intrinsic-psh}, the top Lelong number  is totally intrinsic, that is, it does not depend neither on $h$ nor on $\omega.$
By Theorem \ref{T:charac-Lelong-for-psh}, the top Lelong number  has a geometric meaning in the sense of  in the sense of Siu and Alessandrini--Bassanelli (see Theorem \ref{T:AB-2}). This  completes the proof of assertion (4).

 Assertion (5) is proved in  Theorem \ref{T:top-Lelong-psh}.

 Assertion (6) is proved in  Theorem \ref{T:conic-psh}.

 In assertion (7) $T$ is a   positive pluriharmonic  $(p,p)$-current  on a neighborhood of $\overline B$ in $X$ such that  $T=T^+-T^-$  for some $T^\pm\in\PH_p^{2,2}( B).$ 
 Then this assertion 
 follows  by combining  Theorem  \ref{T:Lelong-psh} (6) and  Theorem \ref{T:conic-ph}.

\endproof


 \proof[Proof of Theorem \ref{T:main_2} {\rm (Tangent Theorem II) }]
  Let $X,$ $V,$ $B$  be as  above  and  suppose that  $\ddc \omega^j=0$  on $V$ for $1\leq j\leq  \upm-1.$
  Let $T$ be  a   positive closed  $(p,p)$-current  on a neighborhood of $\overline B$ in $X$ such that  $T=T^+-T^-$  for some $T^\pm\in\CL_p^{2,2}( B).$
  
  Assertion (1) is proved in  Theorem \ref{T:Lelong-closed} (2).
  
   Assertion (2) is proved in  Theorem \ref{T:Lelong-closed} (5).
   
    Assertion (3) is proved in  Theorem \ref{T:Lelong-closed} (4).
    
     The non-negativity of the top Lelong number $\nu_\upm(T,B,\omega,h)$ stated in assertion (4) is proved in  Theorem \ref{T:Lelong-closed} (6).
By Theorem \ref{T:tot-intrinsic-closed}, the top Lelong number  is totally intrinsic, that is, it does not depend neither on $h$ nor on $\omega.$
By Theorem \ref{T:charac-closed-all-degrees} (1), the top Lelong number  has a geometric meaning in the sense of  in the sense of Siu and Alessandrini--Bassanelli (see Theorem \ref{T:AB-2}). This  completes the proof of assertion (4).

 Assertion (5) is proved in  Theorem \ref{T:Lelong-closed-all-degrees} and  Theorem  \ref{T:charac-closed-top-hol}.

 Assertion (6) is proved in  Remark \ref{R:conic-closed}.

 We come  to the proof of   assertion (7). 
  If  instead of the above  assumption on $\omega$ and  $T,$  we assume  that  the form $\omega$  is  K\"ahler and $T$ is a   positive closed  $(p,p)$-current  on a neighborhood of $\overline B$ in $X$ such that  $T=T^+-T^-$  for some $T^\pm\in\CL_p^{1,1}( B),$  then all the above  assertions still hold  by Theorem  \ref{T:Lelong-closed-Kaehler}  and   Theorem \ref{T:conic-closed}.  
 If  moreover $\tau$ is  holomorphic and  $T=T^+-T^-$  for some $T^\pm\in\CL_p^{1}( B),$ then the above four assertions (1)--(4) still hold for $j=\upm$  by Theorem  \ref{T:top-Lelong-closed} and  by Theorem \ref{T:charac-closed-top-hol}.

\endproof

 We consider the  special but very important   case  where $\supp( T)\cap V$ is compact in $V.$
 In this case
 we can  choose  any piecewise smooth open  neighborhood $B$  of  $\supp( T)\cap V$ in $V$
 and define following  \eqref{e:def_nu_j_T,V}:
 \begin{equation*}
 \nu_j(T,V,\omega,h):=\nu_j(T,B,\omega,h).
 \end{equation*}
 Using \eqref{e:Lelong-numbers}  and   the inclusion    $\supp( T)\cap V\Subset B,$   we see easily that   this definition is independent of the choice of  such a $B.$ 
 
 \proof[Proof of Theorem \ref{T:main_1'} {\rm (Tangent Theorem I') }]
 Using the above discussion, we may without loss of generality fix a piecewise smooth open  neighborhood $B$  of  $\supp( T)\cap V$ in $V.$
 Then the theorem   follows from Theorem \ref{T:main_1}.  
 \endproof

 \proof[Proof of Theorem \ref{T:main_2'} {\rm (Tangent Theorem II') }]
 Using the above discussion, we may without loss of generality fix a piecewise smooth open  neighborhood $B$  of  $\supp( T)\cap V$ in $V.$
 Then the theorem   follows from Theorem \ref{T:main_2}.  
 \endproof
 
 In the remainder of this subsection 
   $X$ is supposed to be   K\"ahler. Consequently,  we are able to apply
 Theorem \ref{T:approx-admiss}. Note that a  proof of this  theorem  will be given  in  Appendix \ref{A:approximation}.
 
 \proof[Proof of Corollary  \ref{C:main_1} {\rm(Tangent Corollary I) }]
 Let $T$ and $T^\pm$  be  positive  plurisubharmonic $(p,p)$-currents on a neighborhood of $\overline B$ in $X$ such that $T=T^+-T^-$
 and  that  assumptions (i)--(ii) are fulfilled.
 Applying Theorem \ref{T:approx-admiss} to $T^\pm$ and  for $m=m'=3$
 yields that $T^\pm\in\SH^{3,3}_p(B).$
 Then assertions (1)--(6) follow  from assertions (1)--(6) of Theorem \ref{T:main_1}.
 
 To prove  assertion (7), $T$ and $T^\pm$ are now  positive  pluriharmonic $(p,p)$-currents on a neighborhood of $\overline B$ in $X$ such that $T=T^+-T^-$
 and  that two $\bullet$  assumptions  are fulfilled.
 Applying Theorem \ref{T:approx-admiss} to $T^\pm$ and  for $m=m'=2$
 yields that $T^\pm\in\SH^{2,2}_p(B).$
 Then assertion (7) follow  from  of Theorem \ref{T:main_1} (7).
 \endproof

 \proof[Proof of Corollary  \ref{C:main_2} {\rm(Tangent Corollary II) }]
 Let $T$ and $T^\pm$  be  positive  closed  $(p,p)$-currents on a neighborhood of $\overline B$ in $X$ such that $T=T^+-T^-$
 and  that   $T^\pm$ are of class $\Cc^2$ in a neighborhood of $\partial B$ in $X.$
 Suppose that  $\ddc \omega^j=0$ on $V$  for $1\leq j\leq  \upm-1.$
 Applying Theorem \ref{T:approx-admiss} to $T^\pm$ and  for $m=m'=2$
 yields that $T^\pm\in\CL^{2,2}_p(B).$
 Then assertions (1)--(6) follow  from assertions (1)--(6) of Theorem \ref{T:main_2}.
 
 To prove  assertion (7),  $\omega$ is  now  a  K\"ahler form on $V$ and   $T^\pm$ are of class $\Cc^1$ in a neighborhood of $\partial B$ in $X.$
 Applying Theorem \ref{T:approx-admiss} to $T^\pm$ and  for $m=m'=1$
 yields that $T^\pm\in\CL^{1,1}_p(B).$
 Then assertion (7) follow  from  of Theorem \ref{T:main_2} (7).
 \endproof

 This is a  consequence of   Appendix \ref{A:admissible_maps}  and Appendix \ref{A:approximation}.

 \proof[Proof of Theorem
 \ref{T:approx-admiss}]
 
Assertion (1) follows from Theorem  \ref{T:strong-admissible-maps} in Appendix \ref{A:admissible_maps}.

Assertion (2) is a consequence of Theorem  \ref{T:approximation}  Appendix \ref{A:approximation}.
  \endproof

\subsection{Dependence of the generalized Lelong numbers on the metrics}

We keep the  hypothesis and notation  of Theorem \ref{T:main_1} (resp.  Theorem \ref{T:main_2}).
Let $h$ be a Hermitian  metric  on $\E.$ 
Consider the function $\varphi_h:\ \E_{\pi^{-1}(V_0)}\to\R^+$ defined by
$$\varphi_h(y):=\|y\|^2_h\qquad\text{for}\qquad  y\in\pi^{-1}(V_0)\subset \E.  $$
Consider also  the  $(1,1)$-closed  smooth  form $\beta_h:=\ddc \varphi_h$ on $\pi^{-1}(V_0)\subset \E.  $  
Fix  a  constant $c_h>0$   such that \begin{equation}\label{e:hat-beta-bis}
  \hat\beta_h:=c\varphi_h\cdot  \pi^*\omega+\beta_h
 \end{equation}  is  positive  on $\pi^{-1}(\overline B)$ and is strictly positive on $\pi^{-1}(\overline B)\setminus \overline B.$
Recall  the following mass  indicators already defined in  \eqref{e:other-mass-indicator}:
\begin{equation}\label{e:other-mass-indicator-bis} \hat\nu_j(T,B,r,\tau,\omega,h):={1\over  r^{2(k-p-j)}}\int\limits_{\Tube(B,r)}
 \tau_* T\wedge (\beta_h+c_hr^2\pi^*\omega)^{k-p-j}\wedge  \pi^*\omega^j.
 \end{equation}
 Recall  from  Propositions \ref{P:comparison-Mc} and \ref{P:comparison-Mc-bis} that 
 we have, for $\lowm \leq j\leq \upm,$   $\hat\nu_j(T,B,\tau,\omega,h)\geq 0$ and 
 \begin{equation}\label{e:hat-nu_j} \hat\nu_j(T,B,\tau,\omega,h)=\sum\limits_{q=0}^{k-p-j}  {k-p-j\choose q} c^q_h\nu_{j+q}(T,B,\tau,\omega,h).
 \end{equation}
Here $\nu_j(T,B,\tau,\omega,h)=0$ for $j>\upm.$
\begin{theorem}
\label{T:Lelong-comparison}
Let $X,$ $V$ be as in  Theorem \ref{T:main_1}. Let $B$  be a  piecewice  $\Cc^2$-smooth open subset of $V$ which admits  a  strongly admissible map $\tau.$
Let    $\omega$ and $\omega'$  be  two Hermitian  forms  on $V$  which satisfy the assumption of  Theorem  \ref{T:main_1}  (resp. of Theorem \ref  {T:main_2}).
Let    $h$ and $h'$  be  two Hermitian  forms  on $V$  which satisfy the assumption of  Theorem  \ref{T:main_1}  (resp. of Theorem \ref  {T:main_2}).
Then there is a constant  $c>1$  with  the following  property.
For every  positive  plurisubharmonic (resp.  positive closed) $(p,p)$-current $T$ on a neighborhood of $\overline B$ in  $X$ such  that
$T=T^+-T^-$ for some $T^\pm \in \SH^{3,3}_p(B)$  (resp.   for some  $T^\pm \in \PH^{2,2}_p(B)$, resp.  for some $T^\pm \in \SH^{1,1}_p(B)$),   we  have 
$$
0\leq  c^{-1}\hat\nu(T,B,\omega',h')\leq \hat \nu_j(T,B,\omega,h)\leq c\hat\nu(T,B,\omega',h')\qquad\text{for}\qquad  \lowm\leq j\leq \upm.
$$
 \end{theorem}

\begin{remark}\rm 
Theorem  \ref{T:Lelong-comparison} means  that  the dependence of the Lelong numbers on the Hermitian  form $\omega$ on $V$ and on the  metric $h$ on $\E$  is not so important.
So  in  the remainder of the article  we often   omit the form $\omega,$ and  write  $\nu_j(T,V)$  (resp.  $\nu_j(T,B)$) instead of   $\nu_j(T,V,\omega,h)$ (resp.
 $\nu_j(T,B,\omega,h)$).  
\end{remark}
\proof
Observe that there is a constant $c>1$ such that for $0<r\leq\bfr,$ it holds that
\begin{equation}\label{e:comparison-h-vs-h'}
c^{-1}(\beta_h+c_hr^2\pi^*\omega)\leq (\beta_{h'}+c_{h'}r^2\pi^*\omega)\leq  c(\beta_h+c_hr^2\pi^*\omega) \qquad\text{on}\qquad \Tube(B,r).
\end{equation}
Following formula \eqref{e:other-mass-indicators-bis}, introduce for $\lowm\leq j\leq \upm,$
\begin{equation*}
 \Mc^\hash_j(T,r,h):={1\over r^{2(k-p-j)}}\int
 T^\hash_r \wedge (\beta_h+c_hr^2\pi^*\omega)^{k-p-j}\wedge  \pi^*\omega^j.
\end{equation*}
Using this, we infer from inequality \eqref{e:comparison-h-vs-h'}
that  there is a constant $c>1$ independent of $T$ such that
\begin{equation*}
 0\leq c^{-1}\Mc^\hash_j(T,r,h)\leq \Mc^\hash_j(T,r,h')\leq c \Mc^\hash_j(T,r,h) .
\end{equation*}
On the other hand, by   Proposition  \ref{P:comparison-Mc} and \eqref{e:hat-nu_j}, we  know that
$$
\lim_{r\to 0+} \Mc^\hash_j(T,r,h)=\hat\nu(T,B,\tau,h)\qquad\text{and}\qquad \lim_{r\to 0+} \Mc^\hash_j(T,r,h')=\hat\nu(T,B,\tau,h').
$$
This, combined  with the previous inequality, implies the result.
\endproof

\subsection{The classical case of a single point}
Finally,  we  discuss   tangent theorems  for the classical  case of a single  point. 
Let $X$  be a complex  manifold  of dimension $k$  
  and  $x$ a point of $X.$  In the classical case where $V=B:=$   the single point $\{x\}$, $\E$ is  just one  fiber $\C^k,$ that is,   $\E\simeq \{x\}\times \C^k\simeq \C^k.$  We can imagine  that there is  only one  form
$\omega$ which is  just the Dirac mass  at $x.$ A Hermitian  metric  $h$ on $\E$
  is  identified  canonically with a constant Hermitian  form on $\C^k.$  
  Let $z$ be a local  chart near $x$ such that $\{x\}=\{z=0\}.$ By Definition \ref{D:admissible-maps}, a strongly  admissible map along $\{x\}$
  is a $\Cc^2$-diffeomorphism  $\tau$ from an open neighborhood $U$ of $\{x\}$ into an open neighborhood of $0$ in $\C^k$
  such that for $z\in U,$ 
  $
  \tau(z)=z+zAz^T+O(\|z\|^3),
  $ 
  where $A$ is a $k\times k$-matrix  with complex  entries.
  
  Theorems  \ref{T:main_1} and \ref{T:main_2}  have  the 
  following  version in the context of the  single point  $x.$
 
 \begin{theorem}\label{T:main_1-single-point} {\rm (Tangent Theorem for a single point) }
  Let $T$ be  a  positive plurisubharmonic $(p,p)$-current  on a neighborhood of $x$ in $X.$ 
  Then,  the following  assertions  hold.
  \begin{enumerate}
  \item The following limit  exists and is  finite
  $$    \nu_0(T,\{x\},h):=\lim\limits_{r\to 0+}\nu_0(T,\{x\},r,\tau ,h)                              $$
   for  all strongly  admissible maps  $\tau$     for $\{x\}$ and for all  Hermitian  metric  $h$ on $\E.$   

  \item  The following  equality holds
  $$\lim\limits_{r\to 0+}\kappa_0(T,\{x\},r,\tau,h)=\nu_0(T,\{x\})$$ 
  for all strongly  admissible maps $\tau$  for $\{x\}$ and for all Hermitian  metrics $h$ on $\E.$
  \item  The  real number  $\nu_0(T,\{x\})$ is  nonnegative and is  totally intrinsic, i.e.  it is independent of  the choice  of   $\tau$ and $h.$ 
  Moreover, it has  a  geometric meaning in the sense of Siu (see Theorem \ref{T:AB-2}).

   \item 
There  exists   tangent currents to $T$  along $\{x\},$   and all tangent  currents $T_\infty$  are  positive plurisubharmonic   on $ \E.$

\item  If morever,  $T$ is  pluriharmonic, then every  tangent  current $T_\infty$ is also   $V$-conic pluriharmonic  on $ \E.$

\item  If morever,  $T$ is  closed, then every  tangent  current $T_\infty$ is also   $V$-conic  closed on $ \E.$

  \end{enumerate}
 \end{theorem}
 
 \proof
 Using local  regularization  we  see that  $T\in\SH_p^{2}( \{x\},\comp).$
 
  When  $T$ is  positive pluriharmonic, we  see that  $T\in\PH_p^{2}( \{x\},\comp).$ When  $T$ is  positive closed, we  see that  $T\in\CL_p^{1}( \{x\},\comp).$

  Here is   the main point. We apply the classical  Lelong--Jensen formula for a ball in $\C^k$ 
  (see \cite{BlelDemaillyMouzali,Demailly93,Demailly}) instead of the our Lelong--Jensen formulas  for tubes in a vector bundle developed in Section \ref{S:Lelong-Jensen}. 
  
  More conceretely,  this classical formula  states that   given    a real   $\Cc^2$-smooth form $S$ of dimension $2q$ defined  in a ball $\B(x,\bfr)$ in $\C^k,$
  then   all $r_1,r_2\in  (0,\bfr)$ with  $r_1<r_2,$  we have that  
\begin{equation}\label{e:Lelong-Jensen-single-point}\begin{split}
 {1\over  r_2^{2q}} \int_{\B(x,r_2)} S\wedge \beta^q- {1\over  r_1^{2q}} \int_{\B(x,r_1)} S\wedge \beta^q
 =    \int_{\B(x,r_1,r_2)} S\wedge \alpha^q\\
 +  \int_{r_1}^{r_2} \big( {1\over t^{2q}}-{1\over r_2^{2q}}  \big)2tdt\int_{\B(x,t)} \ddc S\wedge \beta^{q-1} 
 +  \big( {1\over r_1^{2q}}-{1\over r_2^{2q}}  \big) \int_{0}^{r_1}2tdt\int_{\B(x,t)} \ddc S\wedge \beta^{q-1}.
 \end{split}
\end{equation} 
  Here, $\B(x,r)$ denotes  the  ball with center $x$ and radius $r$ and  $\B(x,r_1,r_2)$ denotes  the corona $\{ y\in\C^k:\ r_1<\|y-x\|<r_2\}.$
  Observe that in comparison with the general formula \eqref{e:Lelong-Jensen},
  formula \eqref{e:Lelong-Jensen-single-point} does not have a vertical boundary term.

  We will apply  formula  \eqref{e:Lelong-Jensen-single-point}  to $S:=\tau_*T.$
  The main difference  in comparison with 
  Theorems  \ref{T:main_1} and \ref{T:main_2} is that we do not have  an $(1,1)$-positive form $\omega$ living on $V=\{x\}.$
  So the technique developed in the proof of Theorem \ref{T:conic-psh} does not work here.
  That is why when $T$ is only positive plurisubharmonic, $T_\infty$ is positive plurisubharmonic, but in general  it is neither  pluriharmonic nor conic. We leave the details of the proof to the interested reader.
 \endproof
   
\begin{remark}\label{R:BlelDemaillyMouzali}
 \rm   When $T$ is   positive closed and $\tau=\id,$  Theorem \ref{T:main_1-single-point} is basically proved by 
 Blel--Demailly--Mouzali in \cite{BlelDemaillyMouzali}.
\end{remark}

\begin{appendix}

\section{Construction of a strongly admissible map} \label{A:admissible_maps}

The main purpose  of this  section is to prove  the first  part of Theorem  \ref{T:approx-admiss}

\begin{theorem}\label{T:strong-admissible-maps}
  Let $X$ be a complex  K\"ahler manifold of dimension $k.$ Let  $V\subset X$ be  a submanifold
  of dimension $l$. 
   Then there  exists a  strongly  admissible map  for $V.$
 \end{theorem}

 Let $\omega$ be  a  K\"ahler  form on $X.$ This induces a Hermitian metric on the tangent bundle $\Tan(X)$ of $X.$
For each point $x \in V$ denote by $N_x$
the orthogonal complement of the tangent space $\Tan_x(V)$ to $V$ at $x$ in the tangent space
$\Tan_x(X)$ to $X$ at $x,$ with respect to the considered metric. 
The union of $N_x$ for $x \in V$
can be identified with the normal bundle $\E$ to $V$ in $X,$ but this identification is not a holomorphic map
in general. We construct the map $\tau^{-1}$ from a neighbourhood of the zero section $V$  in $\E$ to
a neighbourhood of $V$ in $X$ in the following way: for $y \in N_x$ close enough to
$x,$ $\tau^{-1} (y)$ is the image of $y$ by the exponential map  at $x$, which is defined on a neighborhood of $x=0$ in  $\Tan_x(X).$ We can check that $\tau$ is well-defined on an oepn neighborhood $U$  of $V$ in $X$ and is smooth
admissible with $d\tau (x) = \id$ for $x \in V $  (see \cite[Lemma 4.2]{DinhSibony18b}).

We  follow  the  proof  of Dinh-Sibony \cite[Proposition 3.8]{DinhSibony18b}  who treat  the case where  $\dim X=2$ and $\dim V=1.$

\begin{proposition}\label{P:admissible_maps}
  In every local
chart $y=(z,w)$  near $V\cap U$ with  $V\cap U=\{z=0\}$,
 we have   
 \begin{eqnarray*}
 \tau_v(z,w)&= &  z+  z A z^T  +O(\|z\|^3),\\
 \tau_h(z,w)&= &  w + Bz 
 +O(\|z\|^2).
 \end{eqnarray*}
 Here, $A$ is a $(k-l)\times(k-l)$-matrix and $B$ is  a $l\times(k-l)$-matrix whose entries are $\Cc^2$-smooth functions in $w,$ $z^T$ is the transpose of $z,$ and  we  write $$\tau(y)=(\tau_{(1)}(y),\ldots,\tau_{(k-l)}(y),\tau_{(k-l+1)}(y),  \ldots \tau_{(k)}(y))=(\tau_v(y),\tau_h(y))\in\C^{k-l}\times\C^l.$$
\end{proposition}
\proof
 Observe 
that the  identity for $\tau$  is equivalent to the similar identity for $
\tilde \tau := \tau^{-1} .$ We
will prove the last one. Since $d\tilde\tau (z  , w )$ is the identity when $z = 0,$ we have
$\tilde\tau(z,w) =(z,w +a(w)z) +O(\|z\|^2  ).$ So if we write $$\tilde\tau=(\tilde\tau_v,\tilde\tau_h)\in\C^{k-l}\times\C^l$$
in coordinates $(z,w),$  we
only have to check that 
\begin{equation}\label{tau_v}
\tilde\tau_v (z,w)= z+ \sum_{p,q=1}^{k-l} O(1)z_pz_{q}+O(\|z\|^3). 
\end{equation}
This property means
there are no terms with $\bar z_p \bar z_q,$  $\bar z_p z_q$ in the Taylor expansion of $\tilde\tau_v$ in $z ,$  $\bar z$ with
functions in $w$  as coefficients. So it is enough to check it on each complex plane
$\{w \} \times \C^{k-l}.$ Recall that in the local coordinates $(z  , w )$ as above, we identify this
complex plane with the fiber of $\E$ over $(0 , w).$ We will need to make some
changes of coordinates. So we first check that the property does not depend
on our choice of coordinates.

Now consider another system of local holomorphic coordinates $(z' , w' )$
such that $z' = 0$ on $V .$ We can write $w' = H (z  , w )$ and $z' = \alpha(w )z  +
\sum_{p,q=1}^{k-l}h_{pq}(z  , w )z_p z_{q}  ,$   where $H$   and $h_{pq}$ are
$(k-l)\times 1$ matrix whose entries are holomorphic functions, and  $\alpha$ is a $(k-l)\times (k-l)$ matrix whose entries are holomorphic functions. 
For $b' = H (0, b),$
the two complex planes $\C^{k-l}\times \{b\}$ for the coordinates $(z , w )$ and $\C^{k-l}\times \{b'\}$ for the
coordinates $(z' , w' )$ are both identified with the same fiber of $\E.$ The linear
map connecting them is $(z,b ) \mapsto (\alpha(b)z , b').$ We will keep the notation 
$\tilde\tau=(\tilde\tau_v,\tilde\tau_h)$
 for the map $\tilde\tau$  in coordinates $(z , w )$ and use 
$\tilde\tau' = (\tilde\tau'_v,\tilde\tau'_h )$ for the same
map in coordinates $(z' , w').$  With these notations, the point 
$\tilde\tau ( \alpha(b)^{-1}(a)a',b )$ 
 in coordinates $(z , w )$     and the  point  $\tilde\tau'(a',b')$ in coordinates $(z' , w' )$ represent
the same point of $X\times X .$ It follows that
\begin{equation*}
\tilde\tau'_v (a',b')=   \alpha(b)\tilde\tau_v ( \alpha(b)^{-1}a',b ) + \sum_{p,q=1}^{k-l}
 h_{pq}\big( \tilde\tau ( \alpha(b)^{-1}a',b )\big)\tilde\tau_p ( \alpha(b)^{-1}a',b )\tilde\tau_{q}( \alpha(b)^{-1}a',b ).
\end{equation*}
We see that if 
$\tilde\tau_v (a, b) = a+ \sum_{p,q=1}^{k-l} O(1)a_pa_q   + O(\|a\|^3 )$ then 
$\tilde\tau'_v (a', b')$ satisfies a
similar property.

In the rest of the proof, we show \eqref{tau_v}. Without
loss of generality, we will only check the property for $w= 0$ and $z = t\zeta$
with $t \in \R^+$ and $|\zeta| = 1.$  
In a neighbourhood of $0,$ we can write $\tilde \omega= \ddc u$ with $u$ a smooth strictly
psh function. Subtracting from $u$ a pluriharmonic function, we can assume the
existence of a positive definite $k\times  k$-matrix $(c_{ij} )$ such that
$$u(z,w) = (z , w) (c_{ij})\begin{pmatrix}
           \bar z \\
           \bar w 
         \end{pmatrix}+ O(\|(z,w)\|^3 ).$$
We will make changes of coordinates keeping the property $V= \{z = 0\}.$
With a linear change of coordinates $(z , w ) \mapsto (\alpha z, \beta\begin{pmatrix}
            z \\
            w 
         \end{pmatrix}),$
where    $\alpha$ is a $(k-l)\times  (k-l)$-matrix  and  $\beta$  is a $l\times  k$-matrix,       we can
assume that
$$ u(z,w) = \|z\|^2 + \|w \|^2 + O(\|(z,w)\|^3 ).$$
Then, using a change of coordinates of type $$z_p\mapsto z_p+  z_q A_{pq}(z,w),$$  $ A_{pq}(z,w)$ being a linear form in $z$ and $w,$
         we can
assume that the coefficients of  all monomials  in the last
$O(\|(z,w)\|^3 )$ which  can be factored by $z_p\bar z_p$  vanish. Note  that since $u$ is real, when we eliminate the coefficient of
a monomial, the coefficient of its complex conjugate is also eliminated. Next,
using a change of coordinates of type  $w_j\mapsto w_j+\text{quadratic form in $z$ and $w$} ,$ we can
assume that   the coefficients of  all monomials $w_j\bar z_p\bar z_q,$   $w_j\bar z_p\bar w_r,$  $w_j\bar w_r\bar w_s,$  in the last  expression
$O(\|(z,w)\|^3 )$    and their conjugates vanish.
It follows that there remain only monomials $z_p^3,$ $w_j^3,$ $w_jw_r\bar z_p$ and their conjugates, that is,
\begin{eqnarray*}
\tilde \omega &=& \sum_{p=1}^{k-l}idz_p \wedge d\bar z_p+  \sum_{j=1}^l idw_j \wedge d\bar w_j\\
&+&\sum_{1\leq p\leq k-l,\ 1\leq j\leq l} O(|w|) idz_p \wedge d\bar w_j +
\sum_{ 1\leq p\leq k-l,\ 1\leq j\leq l} O(|w|) idw_j \wedge d\bar z_p+
O(\|(z,w)\|^2 ).
\end{eqnarray*}

For the rest of the proof, we use real coordinates $\bfx = (x^1 , \ldots, x^{2k} )$ such
that $z_1 = x^1 + i x^2,$ $\ldots,$ $z_{k-l} = x^{2k-2l-1} + i x^{2k-2l},$ and  $w_1=x^{2k-2l+1}+ix^{2k-2l+2},\ldots,$  $w_l=x^{2k-1}+ix^{2k} .$
Denote by $v = (v^1 ,\ldots,v^{2k}) $ the
unit tangent vector to $X$ at $0$ corresponding to $(\zeta,0) \in \C^{k-l}\times\C^l,$ i.e. $v^{2k-2l+1} = \cdots=v^{2k} = 0$ 
and $\eta_1=v^1+iv^2,$ $\ldots,$  $\eta_{k-l}=v^{2k-2l-1}+iv^{2k-2l}.$ So
$\tilde \tau ( t\eta,0)$ is equal to $\exp(t\eta),$ where $\exp$ denotes the
exponential map from the tangent space to $X$ at $0.$
If we write $\tilde\tau ( t\zeta,0) =
(x^1 (t), \ldots, x^{2k} (t)),$ then $x^j (t)$ satisfy the geodesic equations
$$\ddot x^j = -\Gamma^j_{pq} \dot x^p  \dot x^q
\qquad\text{and}\qquad
\dot x^j (0) = v^j ,
$$
where $\Gamma^j_{pq}$   for $1\leq j,p,q\leq 2k$ are the Christoffel symbols associated with the considered K\"ahler
metric.

We will show in the present setting that $\tilde\tau_v ( t\zeta,0) = t\zeta + O(t^3 )$ and we
already know that $
\tilde \tau_v (t\zeta,0) = t\zeta + O(t^2).$ This is equivalent to checking that
$\ddot x^j (0) = 0$ for $1\leq j \leq 2k-2l.$ Note that the property implies that there is no term of
order $2$ in the Taylor expansion of 
$\tilde\tau_{1 2} ( z,0 )$ in the latest system of coordinates.

According to the discussion at the beginning of the proof, the terms with $z_pz_q$  ($1\leq p,q\leq k-l$)  
may appear when we come back to the original coordinates. Since $v^j = 0$
for $2k-2l+1\leq j \leq 2k,$ we only need to show that $\Gamma^j_{ pq} (0) = 0$ for $j, p,q \in\{2k-2l+1, \ldots,2k\}.$ Let
$g = (g_{jp} )$ be the Riemannian metric associated with $\omega.$ The above description
of $\omega$ implies that $g_{jp} = \delta_{ jp} +O(\sum_{q=2k-2l+1}^{2k}|x^{q} |+\|x\|^2 )$ for all $j, p,$ where $\delta_{jp} = 1$ if
$j = p$ and $0$ otherwise. The coefficients of the inverse $(g^{jp} )$ of the matrix $(g_{jp} )
$ satisfy a similar property. Recall that the Christoffel symbols are given by
$$
\Gamma^j_{ pq}={1\over 2}g^{jm}\Big( {\partial g_{mp}\over  \partial x_q} +{\partial g_{mq}\over  \partial x_p}-{\partial g_{pq}\over  \partial x_m} \Big).
$$
It is now easy to check that $\Gamma^j_{ pq} (0) = 0$ for $j, p, q \in \{2k-2l+1,\ldots, 2k\}.$ The proposition
follows.
\endproof

\section{Approximations of currents} \label{A:approximation}

The main purpose  of this  section is to prove  the second part of Theorem  \ref{T:approx-admiss}

\begin{theorem}\label{T:approximation}
  Let $X$ be a  K\"ahler manifold of dimension $k.$ Let  $V\subset X$ be  a submanifold
  of dimension $l$ and  $B\subset  V$  a relatively compact piecewisely $\Cc^2$-smooth open subset.
  Let $m,m'\in\N$ with $m\geq m'.$  
  Let $T$ be a   positive plurisubharmonic (resp.  positive  pluriharmonic, resp.  positive closed) $(p,p)$-current   on $X$ which satisfies the following conditions (i)--(ii):
   \begin{itemize}
   \item[(i)] $T$ is  of class $\Cc^{m'}$ near  $\partial B;$
   \item[(ii)] There is   a relatively compact open subset $\Omega$ of $X$ with  $B\Subset \Omega$
   and $dT$ is of class $\Cc^0$ near $\partial\Omega.$
   \end{itemize}
  Then  $T$ can be written  in an open neighborhood of $\overline B$ in $X$ as $T=T^+-T^-$  for some $T^\pm  \in  \SH^{m,m'}_p(B)$ (resp. $T^\pm  \in  \PH^{m,m'}_p(B),$ $T^\pm  \in  \CL^{m,m'}_p(B)$).
 \end{theorem}

  We  will adapt  Dinh-Sibony's construction of regularizing kernel \cite{DinhSibony04} to our present context of  
  open  K\"ahler manifolds. Note that ther  construction  was  initially used for compact  K\"ahler manifolds.
 Let $\Delta$ be    the  diagonal  of $X\times X.$
 Let $\pi:\ \widetilde {X\times X}\to X\times X$ be the blow-up of $X\times X$  along $\Delta.$   Following  Blanchard \cite{Blanchard},  $\widetilde {X\times X}$ is a K\"ahler manifold.   Set $\widetilde \Delta:=\pi^{-1}(\Delta).$
 Since $[\widetilde \Delta]$ is  a positive closed $(1,1)$-current, there exist a quasi plurisubharmonic  function  
 $\phi$  and  a smooth closed $(1,1)$-form $\Phi'$ such that $\ddc\phi= [\widetilde \Delta]-\Theta'.$  Note that  $\phi$ is   smooth out  of  $\widetilde\Delta.$
  Let $\chi:\ \R\cup\{-\infty\}\to\R$ be an increasing convex smooth   function such that  $\chi(t)=0$ for $t\in [-\infty,-1],$
  $\chi(t)=t$ for $t\in [1,\infty)$ and $0\leq \chi'\leq 1.$ Define, for $n\in\N,$ $\chi_n(t):=\chi(t+n)-n$ and $\phi_n:=\chi_n\circ\phi.$
  So $\phi_n=\phi$   outside a tubular neighborhood of $\widetilde\Delta$  with radius of order $e^{ -n+1}$ and $\phi_n=0$ inside a
  tubular neighborhood of $\widetilde\Delta$  with radius of order $e^{ -n-1}.$ Moreover,
  the  functions $\phi_n$ are smooth decreasing  to $\phi,$ and we have
  \begin{equation}\label{e:ddc_phi_n}
  \begin{split}
  \ddc\phi_n&= (\chi_n''\circ\phi) d\phi\wedge \dc \phi+ (\chi'_n\circ\phi)\ddc\phi \\
  &\geq  (\chi'_n\circ\phi)\ddc\phi=-(\chi'_n\circ\phi)\Theta'\geq -\Theta,
  \end{split}
  \end{equation}
  where  we choose  the  smooth positive closed form $\Theta$ big enough such that $\Theta-\Theta'$ is positive.
  Define, for $n\in\N,$ the  positive closed smooth $(1,1)$-form on $\widetilde{X\times X}:$ 
  $$\Theta^+_n:=\ddc\phi_n+\Theta\quad\text{and}\quad \Theta^-_n:=\Theta-\Theta'.$$
  So  $\Theta^+_n-\Theta^-_n\to  [\widetilde\Delta].$
  Let $\gamma$ be a  closed smooth $(k-1,k-1)$-form   on $\widetilde{X\times X}$ which is  strictly positive  on a nonempty open subset
  of $\widetilde\Delta.$ 
  Then  $\pi_*(\gamma\wedge[\widetilde\Delta])$ is a nonzero positive closed $(k,k)$-current  on $X\times X$ supported on $\Delta.$ So, 
  it is a multiple of $[\Delta].$
  We choose  $\gamma$ so that  $\pi_*(\gamma\wedge[\widetilde\Delta])=[\Delta].$  Define, for $n\in\N,$ 
  \begin{equation}\label{e:Phi_n}
  \begin{split}
  \widetilde K^\pm_n&:=\gamma\wedge \Theta^\pm_n\quad\text{and}\quad K^\pm_n:=\pi_*(\widetilde K^\pm_n),\\
   \widetilde K_n&:=\widetilde K^+_n-\widetilde K^-_n- \gamma\wedge [\widetilde \Delta]\quad\text{and}\quad K_n:=\pi_*(\widetilde K_n).
   \end{split}
  \end{equation}
Observe  that $K^\pm_n$ are  positive  closed  $(k,k)$-forms with coefficients in $L^1_\loc$ which are smooth out of $\Delta.$
Note also that $K^+_n-K^-_n\to [\Delta]$  weakly, so $K_n$   tends to $0$ weakly. 
 Define
 \begin{equation}\label{e:K-pm_n}
 T^\pm_n(x):=\int_{y\in X} K^\pm_n(x,y)\wedge T(y)\quad\text{and}\quad T_n(x):=\int_{y\in X} K_n(x,y)\wedge T(y).
  \end{equation}
In other words, we have
\begin{equation}\label{e:K-pm_n-bis}
T^\pm_n=(\pi_1)_*(K^\pm_n\wedge \pi_2^*(T))\quad\text{and}\quad T_n=(\pi_1)_*(K_n\wedge \pi_2^*(T)),
\end{equation}
 where $\pi_i$ denotes the  canonical projections of $X\times X$ onto  its factors.
 
 We recall  from  \cite{DinhSibony04}  a classical lemma.
 Let $\B$ be the unit ball in $\R^N,$ and  let $K(x,y)$ be a function with compact support in $\B\times \B.$  
 Consider the  function  $K^*:\B\times\B\to\C$ defined by   $K^*(x,z):=K(x,x+z)$ for $(x,z)\in \B\times \B.$
  By the change of variable $(x,z)\mapsto(x,y:=x+z),$ we get the following identity 
\begin{equation}\label{e:K-K-star}
Pf(x)=\int_{z\in\B} K^*(x,z)f(x-z)\quad\text{for every  smooth test form $f.$}
\end{equation}
 Assume that one of the  following two conditions is  satisfied:
 \begin{itemize}
  \item [$\bullet$]
For $(x,y)\in \B\times \B,$ there is a constant $c>0$ such that
 \begin{equation}\label{e:kernel}
  |K(x,y)|\leq  c|x-y|^{2-N}\quad\text{and}\quad |\nabla_x K(x,y)|\leq  c|x-y|^{1-N} \quad\text{and}\quad |\nabla^2_x K(x,y)|\leq  c|x-y|^{-N}.
 \end{equation}
Here, $x=(x_1,\ldots,x_N)$  are coordinates of $\R^N$ and $\nabla_x$ is the derivative with respect to one of the variables $x_1,\ldots,x_N.$
\item For $(x,y)\in \B\times \B,$ there is a constant $c>0$ such that
 \begin{equation}\label{e:kernel-bis}
  |K(x,y)|\leq  c|y|^{2-N}\quad\text{and}\quad |\nabla_x K(x,y)|\leq  c|y|^{1-N} \quad\text{and}\quad |\nabla^2_x K(x,y)|\leq  c|y|^{-N}.
 \end{equation}
\end{itemize}

Let $\mathcal M$ be the set of Radon measures on $\B.$  The integral operator $P=P_K$ associated to the integral  kernel $K$ is  defined  
on $\mathcal M$   by
\begin{equation}\label{e:P_K}
P\mu(x):=\int_{y\in\B} K(x,y)d\mu(y).
\end{equation}
 Note that $P\mu$ is  supported in $\B.$
 
 For every $n\in\N,$ let $\Lambda^n(\B)$ be the space of all functions $f:\ \B\to\C$ such that
 $f$ is $n$-times diffentiable and its $n$-th derivative  $f^{(n)}\in L^\infty(\B).$ Note that $\Lambda^0=L^\infty$ and $\Cc^n\subset \Lambda^n\subset\Cc^{n-1}.$
 
 The  following  result is needed.
 \begin{proposition} \label{P:Lp-Lq}  The following assertions hold:
 \begin{enumerate}  
 \item    If $K$ satisfies either \eqref{e:kernel} or \eqref{e:kernel-bis}, then for every  $1< \delta_0<{N\over N-2},$ $P$ maps continuously
 $\mathcal M$ into $L^\delta.$  It also maps  continuously  $L^p$ into $L^q,$ where $q=\infty$ if   $p^{-1}+
 \delta_0^{-1}\leq 1$
 and $p^{-1}+\delta^{-1}= 1+q^{-1}$ otherwise.
 \item   If $K$ satisfies either \eqref{e:kernel} or \eqref{e:kernel-bis}, then $P$ maps continuously $L^\infty$ into $\Lambda^1.$
 \item If there is  a  $m\in\N$  such that $K^*$  satisfies
 \begin{eqnarray*}
 |\nabla^j_x K^*(x,z)|&\leq & c|z|^{2-N} \qquad\text{for}\qquad   0\leq j\leq m,\\
 |\nabla^{m+1}_x K^*(x,z)|&\leq & c|z|^{1-N}\qquad\text{and}\qquad|\nabla^{m+2}_x K^*(x,z)|\leq  c|z|^{-N},
 \end{eqnarray*}
   then  $P$ maps continuously
 $\Lambda^m$  into  $\Lambda^{m+1}.$
 \end{enumerate}
 \end{proposition}
 \proof
 \noindent {\bf Proof of assertion (1).} 
It  follows  from  Young's inequality, see  \cite[Theorem 0.3.1]{Sogge}.

 \noindent {\bf Proof of assertion (2) when  $K$ satisfies condition \eqref{e:kernel}.}
 Fix $x^0\in \B.$  
 Fix a test form $f\in L^\infty(\B).$ Suppose without loss of generality that $\|f\|_\infty\leq 1.$  For $x,z\in\B,$ write
 $
 (Pf)(z)=(P_{1,x}f)(z)+(P_{2,x}f)(z),  $ where
 $$
 (P_{1,x}f)(z):=\int_{y\in\B:\ \|y-x^0\|<2\|x-x^0\|} K(z,y)f(y)\quad\text{and}\quad (P_{2,x}f)(z):=\int_{y\in\B:\ \|y-x^0\|\geq 2\|x-x^0\|} K(z,y)f(y).
 $$
 Since  $ \|y-x\|<3\|x-x^0\|$ for $ y\in\B$ with $ \|y-x^0\|<2\|x-x^0\|,$ we infer that 
 \begin{multline*}
  \int_{ \|y-x^0\|<2\|x-x^0\|} | K (x,y)- K (x^0,y)||f(y)|\lesssim\int_{ \|y-x\|<3\|x-x^0\|} | K (x,y)||f(y)|\\ +
   \int_{ \|y-x^0\|<2\|x-x^0\|} | K (x_0,y)||f(y)|.
 \end{multline*}
 Using the first inequality in 
\eqref{e:kernel} and the assumption  $\|f\|_\infty\leq 1,$ we  see that the each term on the RHS of the last line is $\lesssim  \|x-x^0\|^2.$
So 
\begin{equation}\label{e:lim_P_1,x}
\lim_{x\to x^0} {(P_{1,x}f)(x)-(P_{1,x}f)(x_0)\over \|x-x^0\|}=\lim_{x\to x^0} O( \|x-x^0\|)= 0. 
\end{equation}
 
 Next,  we show that  
 \begin{equation}\label{e:lim_P_2,x}
\lim_{x\to x^0} {(P_{2,x}f)(x)-(P_{2,x}f)(x_0)\over \|x-x^0\|}=  \int_{y\in\B} {\partial  K\over \partial x} (x_0,y)f(y). 
\end{equation}
 For $x\in \B$ close  to $x^0,$ we can find $\xi\in[x,x^0]$ such that
 $$
 K(x,y)-K(x^0,y)= {\partial  K\over \partial x} (\xi,y)(x-x^0).
 $$
 Observe that there is $\theta\in[\xi,x^0]$ such that 
 $$
 {\partial  K\over \partial x} (\xi,y)-{\partial  K\over \partial x} (x_0,y)={\partial^2  K\over \partial x^2} (\theta,y)(\xi-x^0).
 $$
 So $\theta\in[x,x^0].$
 Since  for $
 y\in\B$ with  $\|y-x^0\|\geq 2\|x-x^0\|$ we have  $\|\theta-y\|\geq  \|x-x^0\|,$ it follows from the third inequality  in 
\eqref{e:kernel} that 
 $$
\big |{\partial^2  K\over \partial x^2} (\theta,y)(\xi-x^0)\big|\lesssim    {\|x-x^0\|\over \|x_0-y\|^m}\lesssim  {\|x-x^0\|^{1\over 2}\over \|x_0-y\|^{N-{1\over  2}}}.
 $$
 This, combined with the two previous  equalities and the definition of $P_{2,x}f,$ implies that  
\begin{equation*}
 {(P_{2,x}f)(x)-(P_{2,x}f)(x_0)\over \|x-x^0\|}-\int_{y\in\B} {\partial  K\over \partial x} (x_0,y)f(y)
\leq \int_{y\in\B:\ \|y-x^0\|\geq 2\|x-x^0\|} {\|x-x^0\|^{1\over 2}\over \|x_0-y\|^{N-{1\over  2}}}    f(y).
\end{equation*}
Since the RHS is of order $\|x-x^0\|^{1\over 2},$  we get  \eqref{e:lim_P_2,x}.
The result follows from combining  \eqref{e:lim_P_1,x} and \eqref{e:lim_P_2,x}.


 \noindent {\bf Proof of assertion (2) when  $K$ satisfies condition \eqref{e:kernel-bis}.}
 Fix $x^0\in \B.$  
 Fix a test form $f\in L^\infty(\B).$ Suppose without loss of generality that $\|f\|_\infty\leq 1.$  For $x,z\in\B,$ write
 $
 (Pf)(z)=(P_{1,x}f)(z)+(P_{2,x}f)(z),  $ where
 $$
 (P_{1,x}f)(z):=\int_{y\in\B:\ \|y\|<2\|x-x^0\|} K(z,y)f(y)\quad\text{and}\quad (P_{2,x}f)(z):=\int_{y\in\B:\ \|y\|\geq 2\|x-x^0\|} K(z,y)f(y).
 $$
 We infer that 
 \begin{multline*}
  \int_{ \|y\|<2\|x-x^0\|} | K (x,y)- K (x^0,y)||f(y)|\lesssim\int_{ \|y\|<2\|x-x^0\|} | K (x,y)||f(y)|\\ +
   \int_{ \|y\|<2\|x-x^0\|} | K (x_0,y)||f(y)|.
 \end{multline*}
 Using the first inequality in 
\eqref{e:kernel-bis} and the assumption  $\|f\|_\infty\leq 1,$ we  see that the each term on the RHS of the last line is $\lesssim  \|x-x^0\|^2.$
So 
\begin{equation}\label{e:lim_P_1,x-bis}
\lim_{x\to x^0} {(P_{1,x}f)(x)-(P_{1,x}f)(x_0)\over \|x-x^0\|}=\lim_{x\to x^0} O( \|x-x^0\|)= 0. 
\end{equation}
 
 Next,  we show that  
 \begin{equation}\label{e:lim_P_2,x-bis}
\lim_{x\to x^0} {(P_{2,x}f)(x)-(P_{2,x}f)(x_0)\over \|x-x^0\|}=  \int_{y\in\B} {\partial  K\over \partial x} (x_0,y)f(y). 
\end{equation}
 For $x\in \B$ close  to $x^0,$ we can find $\xi\in[x,x^0]$ such that
 $$
 K(x,y)-K(x^0,y)= {\partial  K\over \partial x} (\xi,y)(x-x^0).
 $$
 Observe that there is $\theta\in[\xi,x^0]$ such that 
 $$
 {\partial  K\over \partial x} (\xi,y)-{\partial  K\over \partial x} (x_0,y)={\partial^2  K\over \partial x^2} (\theta,y)(\xi-x^0).
 $$
 So $\theta\in[x,x^0].$
  It follows from the third inequality  in 
\eqref{e:kernel-bis} that 
 $$
\big |{\partial^2  K\over \partial x^2} (\theta,y)(\xi-x^0)\big|\lesssim    {\|x-x^0\|\over \|y\|^N}\lesssim  {\|x-x^0\|^{1\over 2}\over \|y\|^{N-{1\over  2}}}.
 $$
 This, combined with the two previous  equalities and the definition of $P_{2,x}f,$ implies that  
\begin{equation*}
 {(P_{2,x}f)(x)-(P_{2,x}f)(x_0)\over \|x-x^0\|}-\int_{y\in\B} {\partial  K\over \partial x} (x_0,y)f(y)
\leq \int_{y\in\B:\ \|y\|\geq 2\|x-x^0\|} {\|x-x^0\|^{1\over 2}\over \|y\|^{N-{1\over  2}}}    f(y).
\end{equation*}
Since the RHS is of order $\|x-x^0\|^{1\over 2},$  we get  \eqref{e:lim_P_2,x-bis}.
The result follows from combining  \eqref{e:lim_P_1,x-bis} and \eqref{e:lim_P_2,x-bis}.

 \noindent {\bf Proof of assertion (3).}
Since $\nabla_x^{n} f(x-z)=f^{(n)}(x-z),$ applying     Leibnitz's rule yields  that
 \begin{equation}\label{e:Leibnitz}
 (\nabla_x)^{j} \big(K^*(x,z)f(x-z)\big)=\sum_{n=0}^j   {j\choose n} (\nabla_x)^{n} \big(K^*(x,z)f^{(j-n)}(x-z)\big).
 \end{equation}
 We deduce from  this  and  from the first inequality   that there is a constant $c>0$  such that for $f\in\Cc^m(\B)$ with $\|f\|_{\Cc^m}\leq 1$ and for $(x,z)\in\B\times \B,$
$$  (\nabla_x)^{j} \big(K^*(x,z)f(x-z)\big)\leq  c|z|^{2-N} \qquad\text{for}\qquad   0\leq j\leq m.$$
We will prove  by  induction on $m$  that
\begin{equation}\label{e:induction-nabla}
  (\nabla_x)^{m} (Pf)(x)=\int_{z\in \B}   (\nabla_x)^{m} \big(K^*(x,z)f(x-z)\big).
\end{equation}
Formula \eqref{e:induction-nabla} is true for $m=0$ by \eqref{e:K-K-star}.

Suppose that  \eqref{e:induction-nabla} is true for $m.$  We need to prove it for $m+1.$
Let $L:\B\times\B\to\C$ be a kernel  such that $L^*(x,z)=  (\nabla_x)^{m} \big(K^*(x,z)f(x-z)\big).$
By \eqref{e:Leibnitz}   we infer from the   assumption of assertion (3) that
$L$ satisfies \eqref{e:kernel}. Hence,   writing
\begin{equation*}
  (\nabla_x)^{m} (Pf)(x)=\int_{y\in \B}  L(x,y)  \ind  dy,
\end{equation*}
and applying   assertion (2) to the RHS with the function $\ind\equiv 1,$  it follows that  the function on the LHS is  in $\Lambda^1$. This proves assertion (3).
 \endproof
 The following lemma  shows that  the coefficients of $K^\pm_n,$ $K_n$  satisfy inequality of type  
 \eqref{e:kernel}. Let $(x,y)=(x_1,\ldots,x_k,y_1,\ldots,y_k),$ $|x_j|<3,$ $|y_j|<3,$ be local  holomorphic  coordinates of a chart of
  $\bfU\times \bfU$ such that  $\Delta\cap (\bfU\times \bfU)=\{y=0\}$ in  that chart.  For $n\in \N$ let $W_n:=\{(x,y):\ |y|<e^{-n}\}.$
 
 \begin{lemma}\label{L:Hpm_n}
  For $n\in\N,$ let $H^\pm_n$ (resp. $H_n$) be  a coefficient of $K^\pm_n$ (resp.  $K_n$)  in these coordinates. Then:
  
\begin{enumerate} \item $H^\pm_n$ is of the form 
$A^\pm_n+ B^\pm_n  dy_k+  C^\pm_n  d\bar y_k+D^\pm_n  dy_k\wedge d\bar y_k.$  Here,  $A^\pm_n,$ $B^\pm_n,$ $C^\pm_n$ and $D^\pm_n$  are of the form
  $$  \sum_{I,J} f_{I,J}(x,{y_1\over y_k},\ldots ,{y_{k-1}\over y_k},y_k) \bigwedge_{i\in I} d \big({y_i\over y_k}\big)\bigwedge_{j\in J} d\big({\bar y_j\over \bar y_k} \big), $$
  where the  $f_{I,J}$'s  are smooth functions.
  \item There exists a constant $c>0,$ independent of  $n,$  such that
  \begin{eqnarray*}
   |H^+_n(x,y)|&\leq& c|\phi_n(y)||y|^{2-2k}\quad\text{and}\quad |\nabla_x H^+_n(x,y)|\leq c|\phi_n(y)||y|^{1-2k},\\
   |\nabla^2_x H^+_n(x,y)|&\leq& c|\phi_n(y)||y|^{-2k}
  \end{eqnarray*}
for $|x_j|\leq 1,$ $|y_j|\leq 1$ and $y\not=0.$
\item There exists a constant $c>0,$ independent of  $n,$  such that
  \begin{eqnarray*}
   |H^-_n(x,y)|&\leq& c|y|^{2-2k}\quad\text{and}\quad |\nabla_x H^-_n(x,y)|\leq c|y|^{1-2k},\\
   |\nabla^2_x H^-_n(x,y)|&\leq& c|y|^{-2k}
  \end{eqnarray*}
for $|x_j|\leq 1,$ $|y_j|\leq 1$ and $y\not=0.$

\item $H_n$ is of the form 
$A_n+ B_n  dy_k+  C_n  d\bar y_k+D_n  dy_k\wedge d\bar y_k.$  Here,  $A_n,$ $B_n,$ $C_n$ and $D_n$  are of the form
  $$  \sum_{I,J} f_{I,J}(x,{y_1\over y_k},\ldots ,{y_{k-1}\over y_k},y_k) \bigwedge_{i\in I} d \big({y_i\over y_k}\big)\bigwedge_{j\in J} d\big({\bar y_j\over \bar y_k} \big), $$
  where the  $f_{I,J}$'s  are  functions supported in $W_n.$ Moreover, there exists a constant $c>0,$ independent of  $n,$  such that
  \begin{eqnarray*}
   |H_n(x,y)|&\leq& c\ind_{W_n} |\log{|y|}| |y|^{2-2k}\quad\text{and}\quad |\nabla_x H_n(x,y)|\leq c\ind_{W_n} |y|^{1-2k},\\
   |\nabla^2_x H_n(x,y)|&\leq& c\ind_{W_n}|y|^{-2k}
  \end{eqnarray*}
for $|x_j|\leq 1,$ $|y_j|\leq 1$ and $y\not=0.$
\end{enumerate}
 \end{lemma}
\proof
 Since  the proof is  not difficult, we leave it   to the interested reader.   
\endproof

 \begin{proposition}\label{P:uniform-cv-K_n} Let $m\in\N.$
 \begin{enumerate}
 \item  The integral operator $P_{K^-_n}$ associated to the integral  kernel $K^-_n$ given in 
formula \eqref{e:P_K} is a bounded operator from $\Lambda^m$ into $\Lambda^{m+1}.$ Moreover,
its norm is $\leq c ,$ where $c$ is a constant independent of $n.$
 \item
  The integral operator $P_{K_n}$ associated to the integral  kernel $K_n$ given in 
formula \eqref{e:P_K} is a bounded operator from $\Lambda^m$ into $\Lambda^{m+1}.$ Moreover,
its norm is $\leq c ne^{-2n},$ where $c$ is a constant independent of $n.$
\end{enumerate}
 \end{proposition}
\proof
Assertion (1) follows from  Lemma \ref{L:Hpm_n}(3).

Assertion (2) follows from  Lemma \ref{L:Hpm_n}(4).

\endproof
Consider the projection $\Pi:= \pi_2\circ \pi:\ \widetilde{X\times X}\to X.$ We prove that $\Pi$ is a submersion and 
$\Pi|_{\overline \Delta}$ is also  a submersion from $\widetilde \Delta$ onto $X.$ To  this end pick charts $U\Subset V'\subset X$
that  we identify with open sets in $\C^k.$  We may suppose that  $U$ is  small enough and $0\in U.$ We can, using the  change of coordinates
$(z,w)\mapsto (z-w,w)$ on $V'\times U,$ reduce to the  product situation $V\times U,$ $U\Subset V\subset \C^k,$ and $\Delta$
is identified  to  $\{0\}\times U.$ Hence,  the  blow-up along   $\{0\}\times U$ is  also a product. So $\Pi^*$ of a  current is  just integration on fibers.

 \begin{lemma}\label{L:CV-T-pm_n}
 Suppose that one  of the following  condition is  fulfilled:
  \begin{enumerate}
   \item  $T$ is a positive closed current in $\CL^{1,1}_p(B);$
    \item $T$ is a positive  pluriharmonic current in $\PH^{2,2}_p(B).$
    \item   $T$ is a positive  plurisubharmonic current in $\SH^{3,3}_p(B).$
  \end{enumerate}
  Then,  the currents $T^+_n-T^-_n$ converge  weakly to $T$ as $n\to\infty.$
 \end{lemma}
\proof[Proof  of assertion (1)]
The potential of $\widetilde \Delta$ is integrable with respect to  $\Pi^*(T)$ since its singularity is like $\log\dist(z,\widetilde\Delta)$ and this  function has bounded integral on fibers of $\Pi.$ In particular,  $[\widetilde \Delta]\wedge \Pi^*(T)$ is  well-defined and is  equal to  $(\Pi|_{\widetilde\Delta})^*(T),$ and   $[\widetilde \Delta]$ has no mass
for  $\Pi^*(T)$ nor for  $\widetilde K^\pm_n\wedge \Pi^*(T).$ We then have
\begin{equation}\label{e:K_n-wedge-pi_2T}
 K_n^\pm\wedge \pi_2^*(T)=\pi_*(\widetilde K^\pm_n\wedge \Pi^*(T))
\end{equation}
since the formula is  valid out of $\Delta$ and there is no mass on $\Delta.$ The potentials of $\widetilde K^+_n$ are decreasing and the currents  $\widetilde K^-_n$ are independent of $n,$ hence
\begin{equation}\label{e:cv-widetilde-K_n}
 \widetilde K^+_n\wedge \Pi^*(T)-\widetilde K^-_n\wedge \Pi^*(T)\to \gamma \wedge [\widetilde \Delta]\wedge \Pi^*(T)=
 \gamma \wedge (\Pi|_{\widetilde\Delta})^*(T).
\end{equation}
Since $\pi|_{\widetilde \Delta}$ is a submersion onto $\Delta,$ we have  $(\Pi|_{\widetilde\Delta})^*(T)=
(\pi|_{\widetilde \Delta})^*(\pi_2|_\Delta)^*(T).$  Hence,
\begin{equation*}
 \pi_*\big(\gamma\wedge (\Pi|_{\widetilde\Delta})^*(T)\big)=(\pi_2|_\Delta)^*T.
\end{equation*}
This, and \eqref{e:K_n-wedge-pi_2T} and \eqref{e:cv-widetilde-K_n} imply that
\begin{equation*}
 K^+_n\wedge \pi_2^*(T)- K^-_n\wedge \pi_2^*(T)\to   (\pi_2|_\Delta)^*(T).
\end{equation*}
Taking the direct image  under $\pi_1$  gives $T^+_n-T^-_n\to T.$
\endproof

\proof[Proof  of assertion (2) and (3)]  We only need to prove the following  analog of \eqref{e:cv-widetilde-K_n}:
\begin{equation}\label{e:cv-widetilde-K_n-analog}
 (\ddc\phi_n+\Theta')\wedge \Pi^*(T)\to (\Pi|_{\widetilde\Delta})^*(T).
\end{equation}
The problem is local. Define $S:= (\Pi|_{\widetilde\Delta})^*(T).$ We choose as in  Lemmas \ref{L:Hpm_n} and   \ref{L:CV-T-pm_n} local holomorphic coordinates $(x_1,\ldots, x_{2k})$ of an open set $\widetilde \bfU$ of $\pi^{-1}(\bfU\times \bfU)$ in $\widetilde{X\times X},$ $|x_j|<1,$ so that in $\widetilde\bfU$ 
\begin{itemize}
 \item [$\bullet$]$\widetilde \Delta=\{x_{2k}=0\};$ hence $\psi:=\phi-\log{|x_{2k}|}$ is  smooth  and $\ddc\psi=-\Theta';$
 \item[$\bullet$] $\Pi(x_1,\ldots,x_{2k})=(x_1,\ldots,x_k).$
\end{itemize}
Define $\sigma(x_1,\ldots,x_{2k})=(x_1,\ldots,x_{2k-1}).$ Since $\Pi=\Pi|_{\widetilde \Delta}\circ\sigma,$
we have $\Pi^*(T)=\sigma^*(S)$ in $\widetilde \bfU.$

Observe that $(\ddc\phi_n+\Theta')\wedge \sigma^*(S)$ is  supported in $\{  \phi<-n+2\}$ and by \eqref{e:ddc_phi_n},
\begin{equation*}
 (\ddc\phi_n+\Theta')\wedge \sigma^*(S)\geq (1-\chi'_n\circ\phi)\Theta'\wedge\sigma^*(S).
\end{equation*}
The definition of $\chi_n$ implies that  the measures $(1-\chi'_n\circ\phi)\Theta'\wedge\sigma^*(S)$ tend to $0.$
Hence,  every limit value of $(\ddc\phi_n+\Theta')\wedge \sigma^*(S)$ is a positive current supported in $\widetilde\Delta.$ On the other hand,  since  $S$ is  plurisubharmonic  and $\ddc\phi_n+\Theta'$ is positive closed,
we see that   $(\ddc\phi_n+\Theta')\wedge \sigma^*(S)$ is plurisubharmonic. Hence, every limit value of $(\ddc\phi_n+\Theta')\wedge \sigma^*(S)$ is a positive plurisubharmonic current supported in $\widetilde\Delta.$
By   Theorem \ref{T:Bassanelli} of Bassanelli, it is a current on $\widetilde \Delta.$ Hence, the proof of  
\eqref{e:cv-widetilde-K_n-analog} is reduced to that of 
\begin{equation*}
 \int_{\widetilde\bfU}\Psi(x_{2k})(\ddc\phi_n+\Theta')\wedge \sigma^*(\Phi\wedge S)\to  \int_{\widetilde\Delta}\Phi\wedge S
\end{equation*}
for every test $(2k-p-1,2k-p-1)$-form  $\Phi$ with compact support in $\widetilde \Delta\cap\widetilde\bfU$
and for every function $\Psi(x_{2k})$ supported in $\{|x_{2k}|<1\},$ such that $\Psi(0)=1.$
Observe that since $\sigma^*(\Phi\wedge S)$ is proportional to $dx_1\wedge d\overline x_1\wedge \cdots\wedge dx_{2k-1}\wedge d\overline x_{2k-1},$   only the component of $\ddc\phi_n+\Theta'$ with respect to $dx_{2k}\wedge d\overline x_{2k}$ is relevant.  When $(x_1,\ldots,x_{2k-1})$ is fixed, we have
$$
\int_{x_{2k}} \Phi(\ddc_{x_{2k}}\phi_n+\Theta')\to 1
$$
since $\ddc_{x_{2k}}\phi_n+\Theta'$ converges to the Dirac mass $\delta_0$ and $\Phi(0)=1.$ The last integral is  uniformly bounded  with respect  to  $n$  and $x_1,\ldots,x_{2k-1}$ because the measures $\ddc_{x_{2k}}\phi_n+\Theta'$ on compact subsets of $\{|x_{2k}|<1,\ x_1,\ldots, x_{2k-1}\,\text{fixed}\}$ are uniformly bounded. This completes the proof.
\endproof

 \begin{lemma}\label{L:finiteness-T-pm_n}
  We keep the  assumption of Lemma  \ref{L:CV-T-pm_n}. Then  the masses  of $(T_n^\pm)$  are uniformly bounded in an open neighborhood of $\overline B$ in  $X.$
 \end{lemma}
\proof
  Let  $\Omega,$ $\Omega'$ be relatively compact open subsets  of $X$ with  $B\Subset \Omega\Subset \Omega'.$
  Set $\widetilde \Omega:=\Pi^{-1}(\Omega)\subset\widetilde{X\times X}.$
  When  the current $T$ is  not closed, we just take $\Omega$ as in assumption (ii) and choose $\Omega'$ slightly bigger than $\Omega$
  such that $T$ is of class $\Cc^1$ near $\partial\Omega'.$
We  will prove that there is a constant $c>0$ independent of $n$ and $T$  such that 
\begin{equation}\label{e:T-pm_n-finite-mass}
 \|\widetilde K^\pm_n\wedge  \Pi^*(T)\|_{\widetilde\Omega}\leq c\|T\|_{\Omega'}.
\end{equation}
Write
\begin{equation*} 
\widetilde K^\pm_n\wedge  \Pi^*(T)=\big(\ddc\phi_n\wedge \gamma\wedge  \Pi^*(T)\big)+\big(\Theta\wedge \gamma\wedge  \Pi^*(T)\big).
\end{equation*}
Since $\Theta$ and $\gamma$ are smooth forms on  $\widetilde{X\times X},$ we have $\|\Theta\wedge \gamma\wedge  \Pi^*(T)\|_{\widetilde\Omega}\lesssim 1.$

Using integration by part formula  (see \cite[Formula III.3.1, p. 144]{Demailly}), we have
\begin{eqnarray*}\int_{\widetilde\Omega}\ddc\phi_n\wedge \gamma\wedge  \Pi^*(T)& =&\int_{\widetilde\Omega}\phi_n\wedge \gamma\wedge  \Pi^*(\ddc T)\\
 &+&\int_{\partial\widetilde\Omega}\dc\phi_n\wedge \gamma\wedge  \Pi^*(T)-\int_{\partial\widetilde\Omega}\phi_n\wedge \gamma\wedge  \Pi^*(\dc T)\equiv I_1+I_2+I_3.
\end{eqnarray*}

Observe that for there is a constant $c>0$ independent of $x\in X$ such that
$$\int_{\Pi^{-1}(x)} \phi_n(\bfx) d(\bfx)<c,$$
where $d\bfx$ is the Lebesgue measure on the $k$-dimensional complex manifold $\Pi^{-1}(x).$
This  follows from the fact that $\int_{0}^1|\log t|dt <\infty.$ Hence, 
$$I_1\leq c\|\Pi^*(\ddc T)\|_{\widetilde\Omega}\lesssim  \|\ddc T\|_{\Omega'}.$$
Using the property  of $\phi_n$ and $T$ near $\partial B,$ we see that the $\Cc^1$-norms of them are uniformly bounded independent of $n.$
Hence,  $I_2\lesssim  \|T\|_{\Omega'}$ and  $I_3\lesssim  \|T\|_{\Omega'}.$ This completes the proof.
\endproof

 \begin{lemma}\label{L:CL-PH-SH}
  If $T$ is  positive closed (resp.  positive  pluriharmonic, resp. positive  plurisubharmonic), then  so are the currents $T^\pm_n$ for $n\in\N.$  
 \end{lemma}

 \proof
 It  follows from
\eqref{e:K-pm_n},  \eqref{e:K-pm_n-bis} and \eqref{e:K_n-wedge-pi_2T}.
 \endproof

 \begin{lemma}\label{L:control-boundary}
  The  $\Cc^{m'}$-norms of $T^\pm_n$  are uniformly bounded  in a  neighborhood  of $\partial B$ in $X.$ 
 \end{lemma}
\proof
Since  $T$ is of class $ \Cc^{m'}$ near $\partial B,$ we may find an open neighborhood $\bfW'$ of $\partial B$ in $X$
such that   $T|_{\bfW'}$ belons to the  class $ \Cc^{m'}.$
Consider $n\geq  m.$ Since $m'\leq m,$ we  get  $n\geq m'.$
By Proposition \ref{P:uniform-cv-K_n} (1), we get that
$\| T^-_n \|_{\Cc^{m'}(\bfW)}\leq  c \| T \|_{\Cc^{m'}(\bfW')} $ for a constant $c>0$ independent of $n.$
On the other hand,  by Proposition \ref{P:uniform-cv-K_n} (1), we get that
$$\|T^+_n- T^-_n -T\|_{\Cc^{m'}(\bfW)}\leq  c \| T \|_{\Cc^{m'}(\bfW')} $$
for a constant $c>0$ independent of $n.$ 
Hence,  $\| T^\pm_n \|_{\Cc^{m'}(\bfW)}\leq  c \| T \|_{\Cc^{m'}(\bfW')} .$ 
\endproof

 \proof[End of the proof of Theorem \ref{T:approximation}]
 It is  divided into three steps.
 
 \noindent {\bf Step 1.} We show first that   we can  choose  in Theorem \ref{T:approximation}  forms $T^\pm_n$ with $L^1(\bfU)$-coefficients.  Define  $T^\pm_n$ as in  \eqref{e:K-pm_n} and \eqref{e:K-pm_n-bis}. We use partitions of unity of $\bfU  $
 and of $\bfU\times \bfU$ in order to reduce the problem to the case of $\R^m.$
 Following Lemma \ref{L:Hpm_n} and Proposition \ref{P:Lp-Lq}, the forms $T^\pm_n$ have $L^1(\bfU)$-coefficients.
 Lemmas \ref{L:CV-T-pm_n}  and \ref{L:finiteness-T-pm_n}   implies that  $T^+_n-T^-_n\to T$ and $\|T^\pm_n\|_\bfU\leq c\|T\|_\Fc$ for a constant $c$
 independent of $T.$
 
 \noindent {\bf Step 2.} We can now assume that $T$ is a form with $L^1$ coefficients.  Define $T^\pm_n$ as in 
  \eqref{e:K-pm_n} and \eqref{e:K-pm_n-bis}.  Lemmas \ref{L:Hpm_n} and Proposition \ref{P:Lp-Lq} imply that  the forms $T^\pm_n$ have $L^{1+\delta}(\bfU)$-coefficients. We also have   $T^+_n-T^-_n\to T$ and $\|T^\pm_n\|_\bfU\leq c\|T\|_\Fc$ for a constant $c$
 independent of $T.$ Hence, we can assume that $T$ is a form with $L^{1+\delta}$ coefficients. We repeat this process $N$ times with $N\geq \delta^{-1}.$  Lemma \ref{L:Hpm_n} and
 Proposition \ref{P:Lp-Lq} and  Lemmas \ref{L:CV-T-pm_n}  and \ref{L:finiteness-T-pm_n} allows us to reduce the problem to the case  where $T$ is a form with $L^\infty$-coefficients.  If we repeat this process one more times, we can assume that $T$ is a $\Cc^1$-form.   If we repeat this process $m$ more times, we can assume that $T$ is a $\Cc^m$-form. 
 \endproof

 \end{appendix}

\end{document}